\documentclass{article}
\usepackage{amsmath}
\usepackage{amssymb}
\usepackage{amsthm}
\usepackage{cite}
\usepackage{euscript}
\usepackage[arrow,curve,matrix,arc,2cell]{xy}
\UseAllTwocells
\usepackage{hyperref}
\DeclareFontFamily{U}{rsfs}{} \DeclareFontShape{U}{rsfs}{n}{it}{<->
rsfs10}{} \DeclareSymbolFont{mscr}{U}{rsfs}{n}{it}
\DeclareSymbolFontAlphabet{\scr}{mscr}
\def\mathscr{\scr}
\begin{document}
\def\e#1\e{\begin{equation}#1\end{equation}}
\def\ea#1\ea{\begin{align}#1\end{align}}
\def\eq#1{{\rm(\ref{#1})}}
\theoremstyle{plain}
\newtheorem{thm}{Theorem}[section]
\newtheorem{lem}[thm]{Lemma}
\newtheorem{prop}[thm]{Proposition}
\newtheorem{cor}[thm]{Corollary}
\theoremstyle{definition}
\newtheorem{dfn}[thm]{Definition}
\newtheorem{ex}[thm]{Example}
\newtheorem{rem}[thm]{Remark}
\newtheorem{ax}[thm]{Axiom}
\newtheorem{ass}[thm]{Assumption}
\newtheorem{property}[thm]{Property}
\numberwithin{figure}{section}
\numberwithin{equation}{section}
\def\dim{\mathop{\rm dim}\nolimits}
\def\codim{\mathop{\rm codim}\nolimits}
\def\vdim{\mathop{\rm vdim}\nolimits}
\def\depth{\mathop{\rm depth}\nolimits}
\def\sign{\mathop{\rm sign}\nolimits}
\def\Im{\mathop{\rm Im}\nolimits}
\def\Ker{\mathop{\rm Ker}}
\def\Coker{\mathop{\rm Coker}}
\def\GL{\mathop{\rm GL}}
\def\Sh{\mathop{\rm Sh}}
\def\vol{\mathop{\rm vol}}
\def\Stab{\mathop{\rm Stab}\nolimits}
\def\supp{\mathop{\rm supp}}
\def\rank{\mathop{\rm rank}\nolimits}
\def\Hom{\mathop{\rm Hom}\nolimits}
\def\id{{\mathop{\rm id}\nolimits}}
\def\Id{{\mathop{\rm Id}\nolimits}}
\def\Pd{\mathop{\rm Pd}\nolimits}
\def\rsi{{\rm si}}
\def\ssi{{\rm ssi}}
\def\ev{{\rm ev}}
\def\top{{\rm top}}
\def\sing{{\rm sing}}
\def\virt{{\rm virt}}
\def\coh{{\rm coh}}
\def\hom{{\rm hom}}
\def\cs{{\rm cs}}
\def\lf{{\rm lf}}
\def\eff{{\rm eff}}
\def\dR{{\rm dR}}
\def\Mh{{\rm Mh}}
\def\QMh{{\Q{\rm Mh}}}
\def\dRMh{{\rm dRMh}}
\def\Mc{{\rm Mc}}
\def\QMc{{\Q{\rm Mc}}}
\def\dRMc{{\rm dRMc}}
\def\Obj{{\rm Obj}}
\def\MC{{\mathop{\cal{MC}}\nolimits}}
\def\uMC{{\mathop{\smash{\ul{\cal{MC}\!}\,}}\nolimits}}
\def\dga{{\mathop{\bf dga}\nolimits}}
\def\cdga{{\mathop{\bf cdga}\nolimits}}
\def\Rmod{\mathop{R\text{\rm -mod}}}
\def\Top{{\mathop{\bf Top}\nolimits}}
\def\Topem{{\mathop{\bf Top_{em}}\nolimits}}
\def\CSch{{\mathop{\bf C^{\bs\iy}Sch}}}
\def\Man{{\mathop{\bf Man}}}
\def\Manpr{{\mathop{\bf Man_{pr}}}}
\def\Manb{{\mathop{\bf Man^b}}}
\def\Manc{{\mathop{\bf Man^c}}}
\def\tManc{{\mathop{\bf \ti Man^c}}}
\def\Mangc{{\mathop{\bf Man^{gc}}}}
\def\Mancst{{\mathop{\bf Man^c_{st}}}}
\def\Mancwe{{\mathop{\bf Man^c_{we}}}}
\def\Mancsi{{\mathop{\bf Man^c_{si}}}}
\def\Mancis{{\mathop{\bf Man^c_{is}}}}
\def\Mancin{{\mathop{\bf Man^c_{in}}}}
\def\Mangcin{{\mathop{\bf Man^{gc}_{in}}}}
\def\Mangcsi{{\mathop{\bf Man^{gc}_{si}}}}
\def\cManc{{\mathop{\bf\check{M}an^c}}}
\def\dMan{{\mathop{\bf dMan}}}
\def\dOrb{{\mathop{\bf dOrb}}}
\def\Orb{{\mathop{\bf Orb}}}
\def\Orbb{{\mathop{\bf Orb^b}}}
\def\Orbc{{\mathop{\bf Orb^c}}}
\def\tOrbc{{\mathop{\bf \ti Orb^c}}}
\def\tOrbeffc{{\mathop{\bf \ti Orb^c_{eff}}}}
\def\Orbeff{{\mathop{\bf Orb_{eff}}}}
\def\Orbeffb{{\mathop{\bf Orb^b_{eff}}}}
\def\Orbeffc{{\mathop{\bf Orb^c_{eff}}}}
\def\OrbSa{{\mathop{\bf Orb_{Sa}}}}
\def\go{{\rm go}}
\def\cla{{\rm cla}}
\def\ul{\underline}
\def\bs{\boldsymbol}
\def\ge{\geqslant}
\def\le{\leqslant\nobreak}
\def\O{{\mathcal O}}
\def\H{{\mathbin{\mathbb H}}}
\def\K{{\mathbin{\mathbb K}}}
\def\R{{\mathbin{\mathbb R}}}
\def\Z{{\mathbin{\mathbb Z}}}
\def\Q{{\mathbin{\mathbb Q}}}
\def\N{{\mathbin{\mathbb N}}}
\def\C{{\mathbin{\mathbb C}}}
\def\CP{{\mathbin{\mathbb{CP}}}}
\def\RP{{\mathbin{\mathbb{RP}}}}
\def\fC{{\mathbin{\mathfrak C}\kern.05em}}
\def\fD{{\mathbin{\mathfrak D}}}
\def\fE{{\mathbin{\mathfrak E}}}
\def\fF{{\mathbin{\mathfrak F}}}
\def\cA{{\mathbin{\cal A}}}
\def\cB{{\mathbin{\cal B}}}
\def\cC{{\mathbin{\cal C}}}
\def\cD{{\mathbin{\cal D}}}
\def\cE{{\mathbin{\cal E}}}
\def\cF{{\mathbin{\cal F}}}
\def\cG{{\mathbin{\cal G}}}
\def\cH{{\mathbin{\cal H}}}
\def\cI{{\mathbin{\cal I}}}
\def\cJ{{\mathbin{\cal J}}}
\def\cK{{\mathbin{\cal K}}}
\def\cL{{\mathbin{\cal L}}}
\def\cM{{\mathbin{\cal M}}}
\def\cN{{\mathbin{\cal N}}}
\def\cP{{\mathbin{\cal P}}}
\def\cQ{{\mathbin{\cal Q}}}
\def\cR{{\mathbin{\cal R}}}
\def\cS{{\mathbin{\cal S}}}
\def\cT{{\mathbin{\cal T}\kern -0.1em}}
\def\cW{{\mathbin{\cal W}}}
\def\cX{{\cal X}}
\def\cY{{\cal Y}}
\def\cZ{{\cal Z}}
\def\ucC{{\mathbin{\ul{{\cal C}\!}\,}}}
\def\ucD{{\mathbin{\ul{{\cal D}\!}\,}}}
\def\ucE{{\mathbin{\ul{{\cal E}\!}\,}}}
\def\ucF{{\mathbin{\ul{{\cal F}\!\!}\,\,}}}
\def\oM{{\mathbin{\smash{\,\,\overline{\!\!\mathcal M\!}\,}}}}
\def\cV{{\cal V}}
\def\cW{{\cal W}}
\def\bA{{\bs A}}
\def\bB{{\bs B}}
\def\bC{{\bs C}}
\def\bD{{\bs D}}
\def\bE{{\bs E}}
\def\bF{{\bs F}}
\def\bG{{\bs G}}
\def\bH{{\bs H}}
\def\bM{{\bs M}}
\def\bN{{\bs N}}
\def\bO{{\bs O}}
\def\bP{{\bs P}}
\def\bQ{{\bs Q}}
\def\bS{{\bs S}}
\def\bT{{\bs T}}
\def\bU{{\bs U}}
\def\bV{{\bs V}}
\def\bW{{\bs W}\kern -0.1em}
\def\bX{{\bs X}}
\def\bY{{\bs Y}\kern -0.1em}
\def\bZ{{\bs Z}}
\def\uf{{\underline{f\!}\,}{}}
\def\ug{{\underline{g\!}\,}{}}
\def\uh{{\underline{h\!}\,}{}}
\def\ui{{\underline{i\kern -0.07em}\kern 0.07em}{}}
\def\uid{{\underline{\id\kern -0.1em}\kern 0.1em}}
\def\uphi{{\underline{\phi\!}\,}{}}
\def\al{\alpha}
\def\be{\beta}
\def\ga{\gamma}
\def\de{\delta}
\def\io{\iota}
\def\ep{\epsilon}
\def\la{\lambda}
\def\ka{\kappa}
\def\th{\theta}
\def\ze{\zeta}
\def\up{\upsilon}
\def\vp{\varphi}
\def\si{\sigma}
\def\om{\omega}
\def\De{\Delta}
\def\La{\Lambda}
\def\Om{\Omega}
\def\Ga{\Gamma}
\def\Si{\Sigma}
\def\Th{\Theta}
\def\pd{\partial}
\def\ts{\textstyle}
\def\st{\scriptstyle}
\def\sst{\scriptscriptstyle}
\def\w{\wedge}
\def\sm{\setminus}
\def\lt{\ltimes}
\def\bu{\bullet}
\def\sh{\sharp}
\def\op{\oplus}
\def\od{\odot}
\def\op{\oplus}
\def\ot{\otimes}
\def\ov{\overline}
\def\bigop{\bigoplus}
\def\bigot{\bigotimes}
\def\iy{\infty}
\def\es{\emptyset}
\def\ra{\rightarrow}
\def\rra{\rightrightarrows}
\def\Ra{\Rightarrow}
\def\Longra{\Longrightarrow}
\def\ab{\allowbreak}
\def\longra{\longrightarrow}
\def\hookra{\hookrightarrow}
\def\dashra{\dashrightarrow}
\def\ha{{\ts\frac{1}{2}}}
\def\t{\times}
\def\ci{\circ}
\def\ti{\tilde}
\def\d{{\rm d}}
\def\md#1{\vert #1 \vert}
\def\bmd#1{\big\vert #1 \big\vert}
\def\an#1{\langle #1 \rangle}
\def\ban#1{\bigl\langle #1 \bigr\rangle}
\title{Some new homology and cohomology theories \\ of manifolds and orbifolds}
\author{Dominic Joyce}
\date{Preliminary version, September 2015}
\maketitle

\begin{abstract} 
For each manifold or effective orbifold $Y$ and commutative ring $R$, we define a new homology theory $MH_*(Y;R)$, {\it M-homology}, and a new cohomology theory $MH^*(Y;R)$, {\it M-cohomology}. For $MH_*(Y;R)$ the chain complex $\bigl(MC_*(Y;R),\pd\bigr)$ is generated by quadruples $[V,n,s,t]$ satisfying relations, where $V$ is an oriented manifold with corners, $n\in\N$, and $s:V\ra\R^n$, $t:V\ra Y$ are smooth with $s$ proper near 0 in~$\R^n$. 

We show that $MH_*(Y;R),MH^*(Y;R)$ satisfy the Eilenberg--Steenrod axioms, and so are canonically isomorphic to conventional (co)homology. The usual operations on (co)homology --- pushforwards $f_*$, pullbacks $f^*$, fundamental classes $[Y]$ for compact oriented $Y$, cup, cap and cross products $\cup,\cap,\t$ --- are all defined and well-behaved at the (co)chain level. Chains $MC_*(Y;R)$ form flabby cosheaves on $Y$, and cochains $MC^*(Y;R)$ form soft sheaves on $Y$, so they have good gluing properties.

We also define {\it compactly-supported M-cohomology\/} $MH^*_\cs(Y;R)$, {\it locally finite M-homology\/} $MH_*^\lf(Y;R)$ (a kind of Borel--Moore homology), and two variations on the entire theory, {\it rational M-(co)homology\/} and {\it de Rham M-(co)homology}. All of these are canonically isomorphic to the corresponding type of conventional (co)homology.

The reason for doing this is that our M-(co)homology theories are very well behaved at the (co)chain level, and will be better than other (co)homology theories for some purposes, particularly in problems involving transversality. In a sequel \cite{Joyc9} we will construct virtual classes and virtual chains for Kuranishi spaces in M-(co)homology, with a view to applications of M-(co)homology in areas of Symplectic Geometry involving moduli spaces of $J$-holomorphic curves.
\end{abstract}

\setcounter{tocdepth}{2}
\tableofcontents

\section{Introduction}
\label{kh1}

This book defines some new homology and cohomology theories of manifolds $Y$, or more generally effective orbifolds $Y$. There are three families of theories:
\begin{itemize}
\setlength{\itemsep}{0pt}
\setlength{\parsep}{0pt}
\item (Integral) M-homology and M-cohomology $MH_*(Y;R),MH^*(Y;R),\ldots,$ which are defined over an arbitrary commutative ring $R$, such as $R=\Z$.  
\item Rational M-homology and M-cohomology $MH_*^\Q(Y;R),MH^*_\Q(Y;R),\ldots,$ which are defined over a $\Q$-algebra $R$, such as $R=\Q,\R$ or~$\C$.

Rational M-(co)homology has better symmetry properties than integral M-(co)homology, for example, the cup product $\cup$ is supercommutative on rational M-cochains $MC_*^\Q(Y;R)$, but not on M-cochains $MC_*(Y;R)$.
\item De Rham M-homology and M-cohomology $MH_*^\dR(Y;\R),\ldots$ over $R=\R$, which is a hybrid of rational M-(co)homology and de Rham cohomology.
\end{itemize}
The `M-' stands for `Manifold'. For each we define four kinds of (co)homology:
\begin{itemize}
\setlength{\itemsep}{0pt}
\setlength{\parsep}{0pt}
\item Homology $MH_*(Y;R),MH_*^\Q(Y;R)$ and $MH_*^\dR(Y;\R)$.
\item Cohomology $MH^*(Y;R),MH^*_\Q(Y;R)$ and $MH^*_\dR(Y;\R)$.
\item Compactly-supported cohomology $MH^*_\cs(Y;R),MH^*_{\cs,\Q}(Y;R)$ and \hfil\break $MH^*_{\cs,\dR}(Y;\R)$.
\item Locally finite homology $MH_*^\lf(Y;R),MH_*^{\lf,\Q}(Y;R)$ and $MH_*^{\lf,\dR}(Y;\R)$, also called homology with closed supports, or Borel--Moore homology.
\end{itemize}

Each of our (co)homology theories $MH_*(Y;R),MH^*(Y;R),\ldots$ is defined as the (co)homology of a (co)chain complex $\bigl(MC_*(Y;R),\pd\bigr),\bigl(MC^*(Y;R),\d\bigr),\ldots,$ where the (co)chain spaces $MC_k(Y;R),MC^k(Y;R),\ldots$ are $R$-modules defined by generators and relations. For integral or rational M-homology, the generators of $MC_k(Y;R),MC_k^\Q(Y;R)$ are isomorphism classes $[V,n,s,t]$ of quadruples $(V,n,s,t)$, where $V$ is an oriented manifold with corners, $n\in\N$ with $\dim V=n+k$, and $s:V\ra\R^n$, $t:V\ra Y$ are smooth maps with $s$ proper near 0 in $\R^n$. The boundary operator $\pd:MC_k(Y;R)\ra MC_{k-1}(Y;R)$ acts by $\pd:[V,n,s,t]\mapsto [\pd V,n,s\vert_{\pd V},t\vert_{\pd V}]$. For de Rham M-homology, the generators $[V,n,s,t,\om]$ include a de Rham form $\om$ on $V$, with~$\dim V=n+k+\deg\om$.

All of these (co)homology theories satisfy the Eilenberg--Steenrod axioms, and so are canonically isomorphic to conventional (co)homology theories such as singular homology, \v Cech cohomology, or de Rham cohomology of manifolds. So, the interest in these theories is not in the (co)homology groups themselves, which are already well understood, but in special properties of the theories at the cochain level, which make them especially well suited for certain applications.
 
Section \ref{kh12} will discuss applications of our (co)homology theories in Symplectic Geometry, which were the author's motivation for developing them. Our theories should also have interesting applications in other areas. They are particularly suitable for questions involving transversality. The author is also planning projects applying M-(co)homology to build chain-level models for String Topology \cite{ChSu}, in Singularity Theory, and to perverse sheaves of vanishing cycles on complex manifolds. Readers are invited to find their own applications.

In the rest of the introduction, \S\ref{kh11} explains M-homology $MH_*(Y;R)$ for manifolds, \S\ref{kh12} discusses applications in Symplectic Geometry, and \S\ref{kh13} summarizes related work. Chapter \ref{kh2} gives background on homology and cohomology. 

Chapter \ref{kh3} discusses manifolds with corners, which are vital to our project as our (co)chains $[V,n,s,t]$ or $[V,n,s,t,\om]$ include a manifold with corners $V$ with smooth maps $s:V\ra\R^n$, $t:V\ra Y$. There are many ways to define categories of manifolds with corners, and our constructions are insensitive to the details. 

For future applications, it may be useful for instance to allow $V$ to be a `pseudomanifold' or `stratified manifold' with singularities in codimension 2, or to allow $s,t$ to be only piecewise smooth. So, after explaining the `usual' category of manifolds with corners $\Manc$ in \S\ref{kh31}, we give a list of minimal assumptions in \S\ref{kh33} on a category $\tManc$ of `manifolds with corners', such that our M-(co)homology theories work with (co)chains $[V,n,s,t],[V,n,s,t,\om]$ with $V$ an object in $\tManc$, and $s:V\ra\R^n$, $t:V\ra Y$ morphisms in~$\tManc$.

Chapter \ref{kh4} defines and studies M-(co)homology $MH_*(Y;R)$, $MH^*(Y;R)$, $MH^*_\cs(Y;R)$, $MH_*^\lf(Y;R)$ for manifolds $Y$, the most basic version of our theory. We show that they are canonically isomorphic to conventional (co)homology, and have functorial pushforwards $f_*$, pullbacks $f^*$, and products $\cup,\cap,\t$ defined at the (co)chain level, and have good (co)sheaf-theoretic properties.

Chapter \ref{kh5} explains variations on M-(co)homology of manifolds: rational M-(co)homology, de Rham M-(co)homology, M-(co)homology for effective orbifolds $Y$, and extending M-(co)homology to a bivariant theory in the sense of Fulton and MacPherson~\cite{FuMa}. Chapters \ref{kh6}--\ref{kh8} give longer proofs of results deferred from Chapters \ref{kh2}, \ref{kh4} and \ref{kh5}, respectively.

\subsection{The definition of M-homology}
\label{kh11}

We will eventually define 24 different (co)homology theories, parametrized by the product of sets
\begin{gather*}
\text{$\bigl\{$integral, rational, de Rham$\bigr\}\t\bigl\{$homology, cohomology, compactly-supported}\\ 
\text{cohomology, locally finite homology$\bigr\}\t\bigl\{$manifolds, effective orbifolds$\bigr\}.$}
\end{gather*}
The number of options largely explains the length of this book. Any one of the (co)homology theories can be defined quite succinctly, and for any given application, you probably need only one.

We now explain the easiest theory, integral M-homology $MH_*(Y;R)$ of manifolds $Y$ from \S\ref{kh41}, but all the theories have the same flavour. We simplify slightly by working in the `usual' category of manifolds with corners $\Manc$ from \S\ref{kh31}, rather than a category $\tManc$ satisfying the assumptions in~\S\ref{kh33}.

\begin{dfn} Let $Y$ be a manifold, and $R$ a commutative ring. Consider quadruples $(V,n,s,t)$, where $V$ is an oriented manifold with corners, and $n=0,1,\ldots,$ and $s:V\ra\R^n$ is a smooth map which is proper over an open neighbourhood of 0 in $\R^n$, and $t:V\ra Y$ is a smooth map.

Define an equivalence relation $\sim$ on such quadruples by $(V,n,s,t)\sim(V',\ab n',\ab s',\ab t')$ if $n=n'$, and there exists an orientation-preserving diffeomorphism $f:V\ra V'$ with $s=s'\ci f$ and $t=t'\ci f$. Write $[V,n,s,t]$ for the $\sim$-equivalence class of $(V,n,s,t)$. We call $[V,n,s,t]$ a {\it generator}.

For each $k\in\Z$, define the {\it M-chains\/} $MC_k(Y;R)$ to be the $R$-module generated by such $[V,n,s,t]$ with $\dim V=n+k$, subject to the relations:
\begin{itemize}
\setlength{\itemsep}{0pt}
\setlength{\parsep}{0pt}
\item[(i)] For each generator $[V,n,s,t]$ and each $i=0,\ldots,n$ we have
\begin{equation*}
[V,n,s,t]=(-1)^{n-i}[V\t\R,n+1,s',t\ci\pi_V]\quad\text{in $MC_k(Y;R)$,}
\end{equation*}
where writing $s=(s_1,\ldots,s_n):V\ra\R^n$ with $s_j:V\ra\R$ for $j=1,\ldots,n$ and $\pi_V:V\t\R\ra V$, $\pi_\R:V\t\R\ra\R$ for the projections, then 
\begin{equation*}
s'=(s_1\ci\pi_V,\ldots,s_i\ci\pi_V,\pi_\R,s_{i+1}\ci\pi_V,\ldots,s_n\ci\pi_V):V\t\R\longra\R^{n+1}.
\end{equation*}
\item[(ii)] Let $I$ be a finite indexing set, $a_i\in R$ for $i\in I$, and $[V_i,n,s_i,t_i]$, $i\in I$ be generators for $MC_k(Y;R)$, all with the same $n$. Suppose there exists an open neighbourhood $X$ of $0$ in $\R^n$, such that $s_i:V_i\ra\R^n$ is proper over $X$ for all $i\in I$, and the following condition holds:
\begin{itemize}
\setlength{\itemsep}{0pt}
\setlength{\parsep}{0pt}
\item[$(*)$] Suppose $(x,y)\in X\t Y$, such that for all $i\in I$ and $v\in V_i$ with $(s_i,t_i)(v)=(x,y)$, we have that $v\in V_i^\ci$ and 
\begin{equation*}
T_v(s_i,t_i):T_vV_i\longra T_xX\op T_yY
\end{equation*}
is injective. We require that for all oriented $(n+k)$-planes $P\subseteq T_xX\op T_yY=\R^n\op T_yY$, we have
\e
\begin{split}
&\sum_{\begin{subarray}{l} i\in I,\; v\in V_i^\ci:(s_i,t_i)(v)=(x,y),\;  T_v(s_i,t_i)[T_vV_i]=P \\ \text{$T_v(s_i,t_i):T_vV_i\,{\buildrel\cong\over\longra}\,P$ is orientation-preserving}\end{subarray}\!\!\!\!\!\!\!} a_i=\\
&\sum_{\begin{subarray}{l} i\in I,\; v\in V_i^\ci:(s_i,t_i)(v)=(x,y),\; T_v(s_i,t_i)[T_vV_i]=P \\ \text{$T_v(s_i,t_i):T_vV_i\,{\buildrel\cong\over\longra}\,P$ is orientation-reversing}\end{subarray}\!\!\!\!\!\!\!} a_i\qquad\text{in $R$,}
\end{split}
\label{kh1eq1}
\e
where the injectivity and properness assumptions imply that there are only finitely many nonzero terms in \eq{kh1eq1}.
\end{itemize}
Then
\begin{equation*}
\sum_{i\in I}a_i\,[V_i,n,s_i,t_i]=0\qquad\text{in $MC_k(Y;R)$.}
\end{equation*}
\end{itemize}

Define $\pd:MC_k(Y;R)\ra MC_{k-1}(Y;R)$ to be the unique $R$-linear map with
\e
\pd[V,n,s,t]=[\pd V,n,s\ci i_V,t\ci i_V].
\label{kh1eq2}
\e
Then $\pd\ci\pd=0:MC_k(Y;R)\ra MC_{k-2}(Y;R)$ for all $k$. Define the {\it M-homology groups\/} (or {\it integral M-homology groups\/}) $MH_*(Y;R)$ to be the homology of the chain complex $\bigl(MC_*(Y;R),\pd\bigr)$. That is, for $k\in\Z$ we define $R$-modules
\begin{equation*}
MH_k(Y;R)=\frac{\ts \Ker\bigl(\pd: MC_k(Y;R)\longra MC_{k-1}(Y;R)\bigr)}{\ts \Im\bigl(\pd:MC_{k+1}(Y;R)\longra MC_k(Y;R)\bigr)}\,.
\end{equation*}

If $Y$ is compact and oriented with $\dim Y=m$, define the {\it fundamental cycle\/} $[Y]=[Y,0,0,\id_Y]\in MC_m(Y;R)$. We have $\pd[Y]=0$ as $\pd Y=\es$, so passing to homology gives the {\it fundamental class\/}~$[[Y]]\in MH_m(Y;R)$.

Now let $f:Y_1\ra Y_2$ be a smooth map of manifolds. Define the {\it pushforward\/} $f_*:MC_k(Y_1;R)\ra MC_k(Y_2;R)$ for $k\in\Z$ to be the unique $R$-linear map defined on generators $[V,n,s,t]$ of $MC_k(Y_1;R)$ by
\e
f_*[V,n,s,t]=[V,n,s,f\ci t].
\label{kh1eq3}
\e
Equations \eq{kh1eq2} and \eq{kh1eq3} give $f_*\ci\pd=\pd\ci f_*:MC_k(Y_1;R)\ra MC_{k-1}(Y_2;R)$. So the $f_*$ induce pushforwards $f_*:MH_k(Y_1;R)\ra MH_k(Y_2;R)$ on homology.
\label{kh1def1}
\end{dfn}

\begin{rem}{\bf(a)} This is a complete definition of M-homology $MH_*(Y;R)$, though we have omitted the proofs that $\pd,f_*$ are well defined. In \S\ref{kh41} we also define {\it relative M-homology\/} $MH_*(Y,Z;R)$ for open $Z\subseteq Y$. The most subtle part of the definition is relation (ii). It turns out that (ii) implies that $MC_k(Y;R)=0$ for $k>\dim Y$, but in general $MC_k(Y;R)<0$ for all~$k<0$.
\smallskip

\noindent{\bf(b)} The definitions of {\it M-cohomology\/} $MH^*(Y;R)$ in \S\ref{kh42}, {\it compactly-supported M-cohomology\/} $MH^*_\cs(Y;R)$ in \S\ref{kh43}, and {\it locally finite M-homology\/} $MH_*^\lf(Y;R)$ in \S\ref{kh44}, are similar to Definition \ref{kh1def1}. The main differences are:
\begin{itemize}
\setlength{\itemsep}{0pt}
\setlength{\parsep}{0pt}
\item[(i)] For $MH^*(Y;R)$, in generators $[V,n,s,t]$ in $MC^k(Y;R)$ we take $V$ to be unoriented with $\dim V+k=\dim Y+n$, and $t:V\ra Y$ to be a cooriented submersion, and require $(s,t):V\ra\R^n\t Y$ to be proper near $\{0\}\t Y$ in~$\R^n\t Y$.

The M-cochains $MC^k(Y;R)$ are a kind of sheaf-theoretic completion of $R$-modules of `M-precochains' $\cP MC^k(Y;R)$.

For $f:Y_1\ra Y_2$ smooth, the {\it pullback\/} $f^*:MC^k(Y_2;R)\ra MC^k(Y_1;R)$ is 
\e
f^*[V,n,s,t]=\bigl[V\t_{t,Y_2,f}Y_1,n,s\ci\pi_V,\pi_{Y_1}\bigr],
\label{kh1eq4}
\e
where the fibre product $V\t_{t,Y_2,f}Y_1$ in $\Manc$ exists as $t$ is a submersion.
\item[(ii)] For $MH^*_\cs(Y;R)$, in generators $[V,n,s,t]$ in $MC^k_\cs(Y;R)$ we take $V$ to be unoriented with $\dim V+k=\dim Y+n$, $t:V\ra Y$ to be a cooriented submersion, and require $s:V\ra\R^n$ to be proper near 0 in $\R^n$. Pullbacks $f^*:MC^k_\cs(Y_2;R)\ra MC^k_\cs(Y_1;R)$ are as in \eq{kh1eq4}, but only for $f$ proper.

We may identify $MC^k_\cs(Y;R)$ with the submodule of compactly-supported cochains in $MC^k(Y;R)$.
\item[(iii)] For $MH_*^\lf(Y;R)$, in generators $[V,n,s,t]$ in $MC_k^\lf(Y;R)$ we take $V$ to be oriented with $\dim V=n+k$, and require $(s,t):V\ra\R^n\t Y$ to be proper near $\{0\}\t Y$ in~$\R^n\t Y$.

The locally finite M-chains $MC_k^\lf(Y;R)$ are sheaf-theoretic completions of $R$-modules of `locally finite M-prechains' $\cP MC_k^\lf(Y;R)$.

We may identify $MC_k(Y;R)$ with the submodule of compactly-supported chains in $MC_k^\lf(Y;R)$.
\end{itemize}

\noindent{\bf(c)} In M-cohomology there is a well behaved {\it cup product}, on M-cochains
\begin{equation*}
\cup:MC^k(Y;R)\t MC^l(Y;R)\longra MC^{k+l}(Y;R),
\end{equation*}
defined on generators $[V,n,s,t]\in MC^k(Y;R)$, $[V',n',s',t']\in MC^l(Y;R)$ by
\begin{align*}
\bigl[V,n&,(s_1,\ldots,s_n),t\bigr]\cup\bigl[V',n',(s'_1,\ldots,s_{n'}'),t'\bigr]=(-1)^{ln}\bigl[V\t_{t,Y,t'}V',\\
&n+n',(s_1\ci\pi_V,\ldots,s_n\ci\pi_V,s_1'\ci\pi_{V'},\ldots,s_{n'}'\ci\pi_{V'}),t\ci\pi_V\bigr].
\end{align*}
It is associative but not supercommutative on cochains, although for rational and de Rham M-cohomology the cup products are supercommutative on cochains. Passing to cohomology, $\cup$ is identified with the usual cup product by the canonical isomorphism $MH^*(Y;R)\cong H^*(Y;R)$. Similarly, there is a {\it cap product\/} $\cap$ relating M-cohomology and M-homology, defined at the (co)chain level, and identified with the usual cap product on (co)homology.
\label{kh1rem1}
\end{rem}

We will relate M-homology to (smooth) singular homology.

\begin{dfn} For all $k\ge 0$, define the $k$-{\it simplex\/} to be
\begin{equation*}
\De_k=\bigl\{(x_0,\ldots,x_k)\in\R^{k+1}:x_i\ge 0,\;\>
x_0+\ldots+x_k=1\bigr\}.
\end{equation*}
Define maps $F_j^k:\De_{k-1}\ra\De_k$ for $k>0$ and $j=0,\ldots,k$ by
\begin{equation*}
F_j^k:(x_0,\ldots,x_{k-1})\longmapsto(x_0,\ldots,\ab x_{j-1},\ab 0,\ab x_j,\ldots,x_{k-1}).
\end{equation*}

Let $Y$ be a topological space, and $R$ a commutative ring. Write $C_k^\rsi(Y;R)$ for the free $R$-module spanned by {\it singular\/ $k$-simplices\/} $\si$ in $Y$, which are continuous maps $\si:\De_k\ra Y$. The boundary operator $\pd:C_k^\rsi(Y;R)\ra C_{k-1}^\rsi(Y;R)$ is
\e
\pd:\ts\sum_{i\in I}\rho_i\,\si_i\longmapsto \ts\sum_{i\in
I}\sum_{j=0}^k(-1)^j\rho_i\,(\si_i\ci F_j^k).
\label{kh1eq5}
\e
Then {\it singular homology\/} $H_*^\rsi(Y;R)$ is the homology of $\bigl(C_*^\rsi(Y;R),\pd\bigr)$.

Let $f:Y_1\ra Y_2$ be a continuous map of topological spaces. Define the {\it pushforward\/} $f_*:C_k^\rsi(Y_1;R)\ra C_k^\rsi(Y_2;R)$ by
\e
f_*:\ts\sum_{i\in I}\rho_i\,\si_i\mapsto \ts\sum_{i\in I}\rho_i\,(f\ci\si_i).
\label{kh1eq6}
\e
Then $f_*\ci\pd=\pd\ci f_*$, so $f_*$ induces morphisms $f_*:H_k^\rsi(Y_1;R)\ra H_k^\rsi(Y_2;R)$. 

Now let $Y$ be a manifold. Write $C_k^\ssi(Y;R)$ for the free $R$-module spanned by {\it smooth singular\/ $k$-simplices\/} $\si$ in $Y$, which are smooth maps $\si:\De_k\ra Y$, regarding $\De_k$ as a manifold with corners. Since smooth maps are continuous we have $C_k^\ssi(Y;R)\subseteq C_k^\rsi(Y;R)$. Define $\pd:C_k^\ssi(Y;R)\ra C_{k-1}^\ssi(Y;R)$ as in \eq{kh1eq5}. {\it Smooth singular homology\/} $H_*^\ssi(Y;R)$ is the homology of~$\bigl(C_*^\ssi(Y;R),\pd\bigr)$.

If $f:Y_1\ra Y_2$ is a smooth map of manifolds, we define the {\it pushforward\/} $f_*:C_k^\ssi(Y_1;R)\ra C_k^\ssi(Y_2;R)$ as in~\eq{kh1eq6}. 

The inclusion $F_\ssi^\rsi:\bigl(C_*^\ssi(Y;R),\pd\bigr)\hookra \bigl(C_*^\rsi(Y;R),\pd\bigr)$ induces morphisms on homology $F_\ssi^\rsi:H_*^\ssi(Y;R)\ra H_*^\rsi(Y;R)$. Bredon \cite[\S V.5 \& \S V.9]{Bred1} shows that these are isomorphisms, so smooth singular homology is canonically isomorphic to  singular homology of manifolds.

Now define $F_\ssi^\Mh:C_k^\ssi(Y;R)\ra MC_k(Y;R)$ to be the unique $R$-linear morphism acting on generators $\si$ by
\e
F_\ssi^\Mh:\si\longmapsto[\De_k,0,0,\si],
\label{kh1eq7}
\e
so that $V=\De_k$, $n=0$, $s=0:V\ra\R^0$ and $t=\si:V\ra Y$. We have $\pd\ci F_\ssi^\Mh=F_\ssi^\Mh\ci\pd:C_k^\ssi(Y;R)\ra MC_{k-1}(Y;R)$. Thus $F_\ssi^\Mh$ induces morphisms $F_\ssi^\Mh:H_k^\ssi(Y;R)\ra MH_k(Y;R)$ on homology. From \eq{kh1eq3}, \eq{kh1eq6} and \eq{kh1eq7}  we see that $F_\ssi^\Mh$ commutes with pushforwards $f_*:C_k^\ssi(Y_1;R)\ra C_k^\ssi(Y_2;R)$ and $f_*:MC_k(Y_1;R)\ra MC_k(Y_2;R)$ for $f:Y_1\ra Y_2$ a smooth map.
\label{kh1def2}
\end{dfn}

The next theorem summarizes the main results of~\S\ref{kh41}.

\begin{thm}{\bf(a)} M-homology is a homology theory of manifolds. There are canonical isomorphisms $MH_k(Y;R)\cong H_k(Y;R)$ for all\/ $Y,k,$ preserving the data $f_*$ and isomorphisms $MH_0(*;R)\cong R\cong H_0(*;R),$ where $H_*(-;R)$ is any other homology theory of manifolds over $R,$ such as singular homology\/ $H_*^\rsi(-;R)$. For smooth singular homology $H_*^\ssi(-;R),$ these canonical isomorphisms are the morphisms $F_\ssi^\Mh:H_k^\ssi(Y;R)\ra MH_k(Y;R)$ from Definition\/~{\rm\ref{kh1def2}}.

\smallskip

\noindent{\bf(b)} Let\/ $Y$ be a manifold and\/ $k\in\Z,$ write $\uMC_k(Y;R)(U)=MC_k(U;R)$ for open $U\subseteq Y,$ and define $\si_{TU}:\uMC_k(Y;R)(T)\ra\uMC_k(Y;R)(U)$ for open $T\subseteq U\subseteq Y$ by $\si_{TU}=i_*,$ for $i:T\hookra U$ the inclusion. Then $\uMC_k(Y;R)$ is a \begin{bfseries}flabby cosheaf of $R$-modules on\end{bfseries} $Y,$ in the sense of\/~{\rm\S\ref{kh25}}.
\label{kh1thm1}
\end{thm}

\begin{rem}{\bf(i)} We prove Theorem \ref{kh1thm1}(a) by verifying that M-homology satisfies the Eilenberg--Steenrod axioms for homology. This is far from obvious, in particular, the proof of the Dimension Axiom in \S\ref{kh74}, which says that $MH_0(*;R)\cong R$ and $MH_k(*;R)=0$ for $k\ne 0$, is long and difficult.
\smallskip

\noindent{\bf(ii)} {\it Cosheaves\/} are dual to sheaves, which are more well known. Cosheaves are to homology as sheaves are to cohomology, so the analogue of Theorem \ref{kh1thm1}(b) for M-cohomology in \S\ref{kh42} says that M-cochains form a soft sheaf of $R$-modules $\MC^k(Y;R)$ on $Y$, with $\MC^k(Y;R)(U)=MC^k(U;R)$ for open~$U\subseteq Y$.

For comparison, singular chains $C_k^\rsi(-;R)$ do not form a cosheaf on $Y$, but de Rham cochains $C^k_\dR(-;\R)$ do form a soft sheaf on $Y$.

For M-(co)chains to be (co)sheaves on $Y$ means that M-(co)chains are local objects on $Y$, which can (roughly) be glued together from their values on the sets of an open cover $\{U_i:i\in I\}$ of $Y$. It allows us to use powerful techniques from algebraic geometry and homological algebra, on stalks of sheaves, derived categories of sheaves, sheaf cohomology, and so on.

It also means that M-(co)chains $\al$ have a {\it support\/} $\supp\al$, a closed subset of $Y$ (compact for M-homology), which can for instance be a (singular) submanifold of $Y$. For a generator $[V,n,s,t]$ of $MC_k(Y;R)$ we have
\begin{equation*}
\supp[V,n,s,t]\subseteq t[s^{-1}(0)]\subseteq Y.
\end{equation*}
In applications it may be useful to consider $R$-submodules $MC_k(Y;R)_S\subseteq MC_k(Y;R)$ of $\al\in MC_k(Y;R)$ with $\supp\al\subseteq S\subseteq Y$, for subsets~$S\subseteq Y$.

\label{kh1rem2}
\end{rem}

\subsection{Future applications in Symplectic Geometry}
\label{kh12}

The applications of M-(co)homology the author is most interested in are to areas of Symplectic Geometry involving `counting' moduli spaces of $J$-holomorphic curves, including Gromov--Witten invariants \cite{FuOn,HWZ2}, quantum cohomology \cite{McSa}, Lagrangian Floer cohomology \cite{FOOO1}, Fukaya categories \cite{Seid}, Symplectic Field Theory \cite{EGH}, and contact homology~\cite{EES}.

In each of these fields, one studies moduli spaces $\oM$ of $J$-holomorphic curves $u:\Si\ra M$ in a symplectic manifold $(S,\om)$, where $J$ is an almost complex structure on $S$ compatible with $\om$, and the Riemann surface $\Si$ may have boundary, or nodal singularities, or marked points. The moduli spaces come with `evaluation maps' ${\rm ev}:\oM\ra Y$ to a manifold~$Y$. 

Given ${\rm ev}:\oM\ra Y$, we choose a homology theory (or cohomology theory) $H_*(Y;R)$, defined by a chain complex (or cochain complex) $\bigl(C_*(Y;R),\pd\bigr)$. Then the aim is to construct a `virtual class' $[[\oM]]_\virt$ in $H_k(Y;R)$ if $\oM$ is `without boundary', or a `virtual chain' $[\oM]_\virt$ in $C_k(Y;R)$ otherwise, where $k\in\Z$ is the `virtual dimension' of $\oM$. Geometric relationships between moduli spaces should map to algebraic relationships between virtual classes or chains. One then does homological algebra with the virtual chains to build the theory.

For example, in the Lagrangian Floer cohomology of Fukaya, Oh, Ohta and Ono \cite{FOOO1}, one studies moduli spaces $\oM_k(\be)$ of $J$-holomorphic discs $u:D\ra S$ in $S$ with boundary $u(\pd D)$ in a Lagrangian $L\subset S$, with homology class $u_*([D])=\be\in H_2(M,L;\Z)$, and with $k>0$ boundary marked points. Then $\oM_k(\be)$ has evaluation maps $\ev_i:\oM_k(\be)\ra L$ for $i=1,\ldots,k$. We have a geometrical boundary equation of the form
\e
\pd\oM_k(\be)\cong \coprod_{\be=\ga+\de}\coprod_{k=i+j}\oM_{i+1}(\ga)\t_{\ev_{i+1},L,\ev_{j+1}}\oM_{j+1}(\de).
\label{kh1eq8}
\e

We want virtual chains $[\oM_k(\be)]_\virt$ for $\ev_1\t\cdots\t\ev_k:\oM_k(\be)\ra L^k$ in $C_*(L^k;R)$ such that \eq{kh1eq8} translates into an algebraic equation of the form
\e
\pd[\oM_k(\be)]_\virt=\sum_{\be=\ga+\de}\sum_{k=i+j}[\oM_{i+1}(\ga)]_\virt\bu[\oM_{j+1}(\de)]_\virt,
\label{kh1eq9}
\e
where $\bu:C_*(L^{i+1};R)\t C_*(L^{j+1};R)\ra C_*(L^k;R)$ should be the composition of an external product $\boxtimes:C_*(L^{i+1};R)\t C_*(L^{j+1};R)\ra C_*(L^{k+2};R)$ and contraction $[\De_L]\cdot :C_*(L^{k+2};R)\ra C_*(L^k;R)$ with the diagonal~$\De_L\subset L\t L$.

To define the virtual chain $[\oM]_\virt$, one must first put a suitable geometric structure on the moduli space $\oM$. In the general case, there are two main candidates for this geometric structure in the literature: the {\it Kuranishi spaces\/} of Fukaya, Oh, Ohta and Ono \cite{FOOO1,FOOO3,FuOn}, and the {\it polyfolds\/} of Hofer, Wysocki and Zehnder \cite{HWZ1,HWZ2,HWZ3}. Polyfolds contain more information than Kuranishi spaces, and there should be a forgetful functor from polyfolds to Kuranishi spaces,~\cite{Yang}.

Kuranishi spaces are connected to the author's theory of {\it derived differential geometry} \cite{Joyc4,Joyc5,Joyc6,Joyc7}. The author claims that Kuranishi spaces are really {\it derived orbifolds}, and has given a new definition of Kuranishi spaces \cite{Joyc7} which form a 2-category $\bf Kur$ equivalent to the 2-category $\dOrb$ of d-orbifolds from~\cite{Joyc4,Joyc5,Joyc6}.

In current approaches to Lagrangian Floer cohomology, the (co)homology theory used for virtual chains is either (smooth) singular homology $H_*^\ssi(L;\Q)$ \cite{FOOO1,FOOO2}, or de Rham cohomology $H^*_\dR(L;\R)$ \cite{FOOO3}. (Pardon \cite{Pard} also gives an interesting topological virtual chain construction using the dual $\check H^*(L;\Q)^*$ of \v Cech cohomology.) Singular homology has severe technical disadvantages: one cannot define the product $\bu$ in \eq{kh1eq9}, for example. De Rham cohomology also has disadvantages, for instance it works only over~$R=\R$.
 
Our new (co)homology theories have been designed specially for defining virtual chains of moduli spaces of $J$-holomorphic curves, and doing homological algebra with them. These virtual chains will be constructed in the sequel \cite{Joyc9}. Here, we explain why M-(co)homology is well adapted for this task.

Kuranishi spaces are built out of Kuranishi neighbourhoods, \cite[Def.~4.1]{FOOO2}.

\begin{dfn} Let $X$ be a topological space. A {\it Kuranishi neighbourhood\/} on $X$ is a quintuple $(V_i,E_i,\Ga_i,s_i,\psi_i)$ such that:
\begin{itemize}
\setlength{\itemsep}{0pt}
\setlength{\parsep}{0pt}
\item[(a)] $V_i$ is a smooth manifold, which may have boundary or corners.
\item[(b)] $E_i$ is a finite-dimensional real vector space.
\item[(c)] $\Ga_i$ is a finite group with a smooth, effective action on $V_i$, and a linear representation on $E_i$.
\item[(d)] $s_i:V_i\ra E_i$ is a $\Ga_i$-equivariant smooth map.
\item[(e)] $\psi_i$ is a homeomorphism from $s_i^{-1}(0)/\Ga_i$ to an open subset $\Im\psi_i$ in $X$.
\end{itemize}

\label{kh1def3}
\end{dfn}

A {\it Kuranishi space\/} $\bX=(X,\cK)$ in the sense of Fukaya et al.~\cite{FOOO1,FOOO2,FOOO3} is a topological space $X$ with a family $\bigl\{(V_i,E_i,\Ga_i,s_i,\psi_i):i\in I\bigr\}$ of Kuranishi neighbourhoods on $X$ with $\dim V_i-\rank E_i=\vdim\bX\in\Z$ for $i\in I$, and a family of `coordinate changes' $\Phi_{ij}$ from $(V_i,E_i,\Ga_i,s_i,\psi_i)$ to $(V_j,E_j,\Ga_j,s_j,\psi_j)$ over $\Im\psi_i\cap\Im\psi_j\subseteq X$ for $i,j\in I$. They define {\it smooth maps\/} $\bs f:\bX\ra Y$ from a Kuranishi space $\bX$ to a manifold $Y$ \cite[Def.~4.6]{FOOO2}, giving a smooth map $f_i:V_i\ra Y$ for $i\in I$ compatible with the $\Phi_{ij}$. They call $\bs f$ {\it weakly submersive\/} if the $f_i:V_i\ra Y$ are submersions. 

Let $\bX$ be a compact, oriented Kuranishi space with $\vdim\bX=k$, $Y$ be an oriented manifold with $\dim Y=m$, and $\bs f:\bX\ra Y$ be smooth. Ignoring properness conditions on $s,t$ in generators $[V,n,s,t]$, we will sketch the construction in \cite{Joyc9} of a virtual chain $[\bX]_\virt$ for $\bs f:\bX\ra Y$ in M-(co)homology:
\begin{itemize}
\setlength{\itemsep}{0pt}
\setlength{\parsep}{0pt}
\item[(i)] For $i\in I$, first suppose $\Ga_i=\{1\}$. Choose an isomorphism $E_i\cong\R^{n_i}$, either as real vector spaces, or more generally an isomorphism of vector bundles $V_i\t E_i\cong V_i\t\R^{n_i}$ over $V_i$. Then $s_i$ maps $V\ra\R^{n_i}$, so $[V_i,n_i,s_i,f_i]$ is a chain in $MC_k(Y;\Z)$.

If $\bs f$ is weakly submersive, we can instead regard $[V_i,n_i,s_i,f_i]$ as a cochain in~$MC^{m-k}_\cs(Y;\Z)$.
\item[(ii)] If $\Ga_i\ne\{1\}$, suppose we can choose an isomorphism $E_i\cong\R^{n_i}$ as above such that the action of $\Ga_i$ on $E_i$ is identified with an action of $\Ga_i$ on $\R^{n_i}$ by permutation of the coordinates $(x_1,\ldots,x_{n_i})$ in $\R^{n_i}$. Then the chain $[V_i,n_i,s_i,f_i]$ in $MC_k^\Q(Y;\Q)$ is $\Ga_i$-equivariant, and it is natural to associate the chain $\frac{1}{\md{\Ga_i}}[V_i,n_i,s_i,f_i]$ in $MC_k^\Q(Y;\Q)$ to $(V_i,E_i,\Ga_i,s_i,\psi_i)$.

If $\bs f$ is weakly submersive , we can instead regard $\frac{1}{\md{\Ga_i}}[V_i,n_i,s_i,f_i]$ as a cochain in $MC^{m-k}_{\cs,\Q}(Y;\Q)$.
\item[(iii)] Let $\Phi_{ij}$ be a coordinate change in the Kuranishi structure on $\bX$. If the isomorphisms $E_i\cong\R^{n_i}$, $E_j\cong\R^{n_j}$ in (i),(ii) are compatible under $\Phi_{ij}$ in a certain (strong) sense, then Definition \ref{kh1def1}(i) implies that the (co)chains $\frac{1}{\md{\Ga_i}}[V_i,n_i,s_i,f_i],\frac{1}{\md{\Ga_j}}[V_j,n_j,s_j,f_j]$ are equal over~$\Im\psi_i\cap\Im\psi_j$.
\item[(iv)] Suppose we choose equivariant isomorphisms $E_i\cong\R^{n_i}$ as in (i),(ii) for all $i\in I$, compatible with coordinate changes $\Phi_{ij}$ as in (iii) for all $i,j\in I$. Then the (co)sheaf property of M-(co)chains in Theorem \ref{kh1thm1}(b) implies that we can glue the (co)chains $\frac{1}{\md{\Ga_i}}[V_i,n_i,s_i,f_i]$ for $i\in I$ to build a global virtual chain $[\bX]_\virt$ for $\bs f:\bX\ra Y$.

This $[\bX]_\virt$ may be defined in $MC_k(Y;\Z)$ if $\Ga_i=\{1\}$ for all $i\in I$, and in $MC_k^\Q(Y;\Q)$ otherwise. If $\bs f$ is weakly submersive, we may define $[\bX]_\virt$ in $MC^{m-k}_\cs(Y;\Z)$ if $\Ga_i=\{1\}$ for all $i\in I$, and in $MC^{m-k}_{\cs,\Q}(Y;\Q)$ otherwise.
\end{itemize}

So, to define a virtual chain for $\bs f:\bX\ra Y$, we must choose compatible trivializations for the obstruction bundles $E_i\ra V_i$ in the Kuranishi structure on $\bX$. It is not difficult to do this using known techniques for constructing `good coordinate systems'. It is also easy to choose virtual chains $[\bX]_\virt$ with boundary compatibilities, such as in \eq{kh1eq9}. Note that unlike current techniques \cite{FOOO1,FOOO2,FOOO3}, we do not perturb $\bX$ or $s_i:V_i\ra E_i$, and~$\supp[\bX]_\virt\subseteq \bs f(\bX)\subseteq Y$. 

\subsection{Related work}
\label{kh13}

In 2007 I wrote two preprints \cite{Joyc1,Joyc2} on `Kuranishi homology and cohomology'.  These defined (co)homology theories $KH_*(Y;R),KH^*(Y;R)$ of manifolds $Y$ in which the (co)chains were generated by triples $[\bX,\bs f,\bs G]$ of a Kuranishi space $\bX$ (roughly in the sense of \cite{FuOn,FOOO1}), a smooth map $\bs f:\bX\ra Y$, and some extra data $\bs G$. As in \S\ref{kh12}, the idea was to make defining virtual chains in Symplectic Geometry easier, by using $KH_*(Y;R)$ instead of singular homology.

I did not publish \cite{Joyc1,Joyc2}, as I was unhappy with the definition of Kuranishi space I was using, but at the time I could not solve the problems, so I shelved the project. In 2014 I finally found a theory of Kuranishi spaces I was satisfied with \cite{Joyc7}, so I decided to revise \cite{Joyc1,Joyc2} with the new definition. The result was this book and the sequel \cite{Joyc9}. I chose to split the project into two parts, with (co)homology theories based on manifolds with corners, and then virtual chains for Kuranishi spaces in these (co)homology theories, as I hoped the first part (this book) might also appeal to people not interested in Kuranishi spaces.

Several authors have defined bordism theories or homology theories $\ti H_*(Y;R)$ of manifolds $Y$ (or other topological spaces) in which the chains involve smooth maps $t:V\ra Y$ from a manifold $V$ (plus extra data), as for our chains $[V,n,s,t]$. Oversimplifying a bit, Jakob \cite{Jako} defines `geometric homology groups' $\ti H_k^{\rm Ja}(Y;R)$ as the set of bordism classes $[V,\al,t]$ of triples $(V,\al,t)$ for $V$ a compact, oriented manifold, $\al\in H^*(V;\Z)$ with $\dim V=k+\deg\al$, and $t:V\ra Y$ continuous, and shows that~$\ti H_k^{\rm Ja}(Y;R)\cong H_k(Y;R)$.

McDuff and Salamon \cite[\S 7.1]{McSa} define {\it pseudocycles}, which are smooth maps of manifolds $t:V\ra Y$, with $V$ generally non-compact, such that the `boundary' of $t(V)$ in $Y$ has `dimension' $\le\dim V-2$, and use them to define virtual chains for some $J$-holomorphic curve moduli spaces. Kahn \cite{Kahn} and Zinger \cite{Zing} both define homology groups $\ti H_k^{\rm Ka}(Y;R),\ti H_k^{\rm Zi}(Y;R)$ of manifolds $Y$ as sets of bordism classes $[V,t]$ of pseudocycles $t:V\ra Y$ with $\dim V=k$, with $Y$ compact in \cite{Kahn}, and show that~$\ti H_k^{\rm Ka}(Y;R)\cong H_k(Y;R)\cong\ti H_k^{\rm Zi}(Y;R)$.

Kreck \cite{Krec} defines homology groups $\ti H_k^{\rm Kr}(Y;R)$ of a manifold $Y$ as sets of bordism classes $[V,t]$ of smooth maps $t:V\ra Y$ where $V$ is a compact, oriented `stratifold', a singular space with a stratification into smooth manifolds, and shows that~$\ti H_k^{\rm Kr}(Y;R)\cong H_k(Y;R)$.

Motivated by the problem of defining virtual chains in Floer theories as in \S\ref{kh12}, Lipyanskiy \cite{Lipy} defines homology groups $\ti H_k^{\rm Li}(Y;R)$ of a manifold $Y$, in which the chains $[V,t]$ are a compact, oriented manifold with corners $V$ with $\dim V=k$ and a smooth map $t:V\ra Y$, and shows that $\ti H_k^{\rm Li}(Y;R)\cong H_k(Y;R)$. He also defines cohomology groups $\ti H^k_{\rm Li}(Y;R)$ with cochains $[V,t]$ a manifold with corners $V$ with $\dim V=\dim Y-k$ and a proper, cooriented smooth map $t:V\ra Y$, and shows that $\ti H^k_{\rm Li}(Y;R)\cong H^k(Y;R)$. These are similar to our $MH_*(Y;R),MH^*(Y;R)$ but restricting to $n=0$ in (co)chains $[V,n,s,t]$, and are the closest to our project the author knows in the literature.

Irie \cite{Irie} studies homology theories of `differentiable spaces' $Y$, which include smooth manifolds, and also some classes of infinite-dimensional spaces such as loop spaces $\cL M$. He defines {\it smooth singular homology\/} $H^\ssi_*(Y;R)$, which is as in Definition \ref{kh1def2} in the manifold case, and {\it de Rham homology\/} $H^\dR_*(Y;\R)$, with a morphism $H^\ssi_*(Y;\R)\ra H^\dR_*(Y;\R)$, and shows this is an isomorphism for $Y$ a manifold. His {\it de Rham chains\/} $[V,t,\om]$ are an oriented manifold $V$, a smooth map $t:V\ra Y$, and a compactly-supported exterior form $\om$ on $V$. These are similar to our generators $[V,n,s,t,\om]$ in $MC_*^\dR(V;\R)$ in \S\ref{kh52} when~$n=0$.

Irie was motivated by an earlier attempt by Fukaya \cite[\S 6]{Fuka} to define a homology theory of loop spaces $\cL M$ using a complex $\bigl(C_*^{\rm adR}(\cL M;\R),\pd\bigr)$ of `approximate de Rham chains' $[V,t,\om]$ of an oriented manifold with boundary $V$, a smooth map $t:V\ra\cL M$, and a compactly-supported exterior form $\om$ on $V$, although unfortunately the homology of $\bigl(C_*^{\rm adR}(\cL M;\R),\pd\bigr)$ is not $H_*(\cL M;\R)$, invalidating the results of \cite{Fuka}. The author also hopes in future to extend M-homology to loop spaces, and so complete the programme of~\cite{Fuka}.
\medskip

\noindent{\it Acknowledgements.} I would like to thank Lino Amorim, Kenji Fukaya, Helmut Hofer, Dusa McDuff, and Jack Waldron for helpful conversations. This research was supported by EPSRC grants EP/H035303/1 and EP/J016950/1.

\section{Background on homology and cohomology}
\label{kh2}

Homology and cohomology were developed in the $19^{\rm th}$ century, and are now a mature subject, with a huge literature, and many excellent textbooks, including Bott and Tu \cite{BoTu}, Bredon \cite{Bred1}, Eilenberg and Steenrod \cite{EiSt}, Hatcher \cite{Hatc}, MacLane \cite{MacL}, Massey \cite{Mass2}, Maunder \cite{Maun}, Munkres \cite{Munk}, Spanier \cite{Span}, and tom Dieck~\cite{toDi}.

We present the material in a non-standard way designed for our applications in \S\ref{kh4}--\S\ref{kh5}. Sections \ref{kh21}--\ref{kh22} characterize homology and cohomology of manifolds, just as $R$-modules, using the Eilenberg--Steenrod axioms. Sections \ref{kh23}--\ref{kh24} explain compactly-supported cohomology, and locally finite homology, and \S\ref{kh25} expresses all of these in terms of sheaf cohomology. Sections \ref{kh26}--\ref{kh210} discuss other topics, including products $\cup,\cap,\t$ on (co)homology, and Poincar\'e duality.

We postpone proofs of five theorems until \S\ref{kh6}. Theorems \ref{kh2thm1} and \ref{kh2thm2} are very similar to known results, in particular Eilenberg and Steenrod \cite[Th.~10.1]{EiSt} and Kreck and Singhof \cite[Prop.~10]{KrSi}, but we give proofs in \S\ref{kh61} as these theorems are crucial for establishing the isomorphisms between our new (co)homology theories and conventional (co)homology. The material on strong presheaves at the end of \S\ref{kh25} is new, and Theorem \ref{kh2thm4} on their properties will be proved in \S\ref{kh62}. Sections \ref{kh63} and \ref{kh64} prove Theorems \ref{kh2thm6} and \ref{kh2thm7}, on orbifolds.

\subsection{Homology}
\label{kh21}

A {\it homology theory\/} over a commutative ring $R$ associates to any topological space $Y$ (possibly satisfying some conditions) a sequence of $R$-modules $H_0(Y;R),H_1(Y;R),H_2(Y;R),\ab\ldots,$ with $H_0(*;R)=R$ for $Y=*$ the point, and to any continuous map $f:Y_1\ra Y_2$ (possibly satisfying some conditions) a sequence of $R$-module morphisms $f_*:H_k(Y_1;R)\ra H_k(Y_2;R)$ for $k=0,1,\ldots,$ satisfying a package of properties we explain below. Actually we define $H_k(Y;R)$ for all $k\in\Z$, but it turns out that $H_k(Y;R)=0$ for~$k<0$.

There are several different homology theories, the most popular being singular homology \cite{Bred1,EiSt,Maun,Munk,Span}. There are theorems which say that over nice topological spaces (e.g.~finite CW-complexes, or compact manifolds), these theories are all canonically isomorphic, so often one does not distinguish them.

Often it is better to work with the {\it relative homology groups\/} $H_*(Y,Z;R)$ of a pair $(Y,Z)$ of a topological space $Y$ and a subspace $Z\subseteq Y$. Writing $i:Z\hookra Y$ for the inclusion, these lie in the long exact sequence
\begin{equation*}
\xymatrix@C=15pt{ \cdots \ar[r] & H_k(Z;R) \ar[r]^{i_*} & H_k(Y;R) \ar[r] & H_k(Y,Z;R) \ar[r]^\pd & H_{k-1}(Z;R) \ar[r] & \cdots. }
\end{equation*}

We will be concerned only with homology $H_*(Y;R)$ of smooth manifolds $Y$, and pushforwards $f_*:H_*(Y_1;R)\ra H_*(Y_2;R)$ for smooth maps $f:Y_1\ra Y_2$. So we begin by stating axioms for homology theories of smooth manifolds.

\begin{ax} (Eilenberg--Steenrod axioms for homology of manifolds.) Write $\cC$ for the category whose objects are pairs $(Y,Z)$ with $Y$ a smooth manifold (without boundary) and $Z\subseteq Y$ an open set, and whose morphisms $f:(Y_1,Z_1)\ra (Y_2,Z_2)$ are smooth maps $f:Y_1\ra Y_2$ with~$f(Z_1)\subseteq Z_2$.

Fix a commutative ring $R$. A {\it homology theory over\/} $R$ assigns the data:
\begin{itemize}
\setlength{\itemsep}{0pt}
\setlength{\parsep}{0pt}
\item[(a)] For each object $(Y,Z)$ in $\cC$ and each $k\in\Z$, an $R$-module $H_k(Y,Z;R)$. For brevity we write~$H_k(Y;R)=H_k(Y,\es;R)$.
\item[(b)] For each object $(Y,Z)$ in $\cC$ and each $k\in\Z$, an $R$-module morphism $\pd:H_k(Y,Z;R)\ra H_{k-1}(Z;R)$, called the {\it connecting morphism}.
\item[(c)] For each morphism $f:(Y_1,Z_1)\ra (Y_2,Z_2)$ in $\cC$ and each $k\in\Z$, an $R$-module morphism $f_*:H_k(Y_1,Z_1;R)\ra H_k(Y_2,Z_2;R)$.
\end{itemize}
All this data must satisfy the axioms:
\begin{itemize}
\setlength{\itemsep}{0pt}
\setlength{\parsep}{0pt}
\item[(i)] (Functoriality.) For each $k\in\Z$, mapping $(Y,Z)$ to $H_k(Y,Z;R)$ and $f:(Y_1,Z_1)\ra (Y_2,Z_2)$ to $f_*:H_k(Y_1,Z_1;R)\ra H_k(Y_2,Z_2;R)$ gives a functor $\cC\ra\Rmod$, where $\Rmod$ is the abelian category of $R$-modules.
\item[(ii)] (Exactness.) For each $(Y,Z)\in\cC$, write $i:(Z,\es)\ra(Y,\es)$ and $j:(Y,\es)\ra (Y,Z)$ for the morphisms in $\cC$ induced by the inclusion $Z\hookra Y$ and the identity $\id_Y:Y\ra Y$. Then the following is exact in $\Rmod$:
\end{itemize}
\e
\xymatrix@C=12pt{ \cdots \ar[r] & H_k(Z;R) \ar[r]^(0.47){i_*} & H_k(Y;R) \ar[r]^(0.45){j_*} & H_k(Y,Z;R) \ar[r]^(0.48)\pd & H_{k-1}(Z;R) \ar[r] & \cdots. }
\label{kh2eq1}
\e
\begin{itemize}
\setlength{\itemsep}{0pt}
\setlength{\parsep}{0pt}
\item[(iii)] (Functoriality of $\pd$.) For each $f:(Y_1,Z_1)\ra (Y_2,Z_2)$ in $\cC$ and $k\in\Z$, the following commutes:
\begin{equation*}
\xymatrix@C=90pt@R=14pt{ *+[r]{H_k(Y_1,Z_1;R)} \ar[d]^{f_*} \ar[r]_\pd & *+[l]{H_{k-1}(Z_1;R)} \ar[d]_{(f\vert_{Z_1})_*} \\ *+[r]{H_k(Y_2,Z_2;R)} \ar[r]^\pd & *+[l]{H_{k-1}(Z_2;R).\!\!} }
\end{equation*}
\item[(iv)] (Homotopy.) Suppose $(Y_1,Z_1),(Y_2,Z_2)\in\cC$, and $g:Y_1\t[0,1]\ra Y_2$ is smooth with $g(Z_1\t[0,1])\subseteq Z_2$. Define $f,f':Y_1\ra Y_2$ by $f(y)=g(y,0)$ and $f'(y)=g(y,1)$ for $y\in Y_1$. Then $f_*=f'_*:H_k(Y_1,Z_1;R)\ra H_k(Y_2,Z_2;R)$ for all $k\in\Z$.
\item[(v)] (Excision.) Suppose $(Y,Z)\in\cC$, and $S\subseteq Z$ is closed in $Y$, so the inclusion $j:Y\sm S\hookra Y$ is a morphism $j:(Y\sm S,Z\sm S)\ra (Y,Z)$ in $\cC$. Then $j_*:H_k(Y\sm S,Z\sm S;R)\ra H_k(Y,Z;R)$ is an isomorphism for all~$k\in\Z$.
\item[(vi)] (Additivity Axiom.) Suppose $(Y,Z)\in\cC$ with $Y=\coprod_{a\in A}Y_a$ for $A$ a countable indexing set and each $Y_a$ open and closed in $Y$, and set $Z_a=Z\cap Y_a$. Then for all $k\in\Z$ we have canonical isomorphisms
\e
H_k(Y,Z;R)\cong \ts\bigop_{a\in A}H_k(Y_a,Z_a;R),
\label{kh2eq2}
\e
compatible with the morphisms $(i_a)_*:H_k(Y_a,Z_a;R)\ra H_k(Y,Z;R)$ for $a\in A$ induced by the inclusion $i_a:(Y_a,Z_a)\hookra(Y,Z)$.
\item[(vii)] (Dimension Axiom.) When $Y=*$ is the point we have $H_0(*;R)\cong R$ and $H_k(*;R)=0$ for all~$k\ne 0$.
\end{itemize}
\label{kh2ax1}
\end{ax}

\begin{rem}{\bf(a)} Parts (i)--(v) and (vii) of Axiom \ref{kh2ax1} were introduced by Eilenberg and Steenrod \cite[\S I.3]{EiSt}, and are known as the {\it Eilenberg--Steenrod Axioms}. They are usually applied to pairs $(Y,Z)$ with $Y$ a topological space and $Z\subseteq Y$ a subspace, and morphisms $f:(Y_1,Z_1)\ra (Y_2,Z_2)$ with $f:Y_1\ra Y_2$ continuous with $f(Z_1)\subseteq Z_2$, but Eilenberg and Steenrod start with an {\it admissible category for homology theory\/} \cite[p.~5]{EiSt}, which include our category~$\cC$.

The Additivity Axiom (vi) was added by Milnor \cite{Miln}, and says that the homology functors must commute with certain infinite colimits. We have included the assumption that $A$ is countable as our definition of manifolds $Y$ requires $Y$ to be a second countable topological space, so $Y$ can be the disjoint union of at most countably many nonempty open subsets~$Y_a$.
\smallskip

\noindent{\bf(b)} If $\{M_a:a\in A\}$ is a family of $R$-modules, we define the {\it direct sum\/} $\bigop_{a\in A}M_a$, used in \eq{kh2eq2}, and {\it direct product\/} $\prod_{a\in A}M_a$ as $R$-modules by
\ea
\bigop_{a\in A}M_a&=\bigl\{(r_a)_{a\in A}:r_a\in R_a,\;\> \text{$r_a\ne 0$ for only finitely many $a\in A$}\bigr\},
\label{kh2eq3}\\
\prod_{a\in A}M_a&=\bigl\{(r_a)_{a\in A}:r_a\in R_a,\;\> a\in A\bigr\}.
\label{kh2eq4}
\ea

\noindent{\bf(c)} Homology groups $H_*(Y,Z;R)$ are usually defined as the homology of some chain complex of $R$-modules $\bigl(C_*(Y,Z;R),\pd\bigr)$. But these chain complexes may not satisfy analogues of (i)--(vii), nice properties appear on passing to homology.

\smallskip

\noindent{\bf(d)} Homology is a homotopy invariant. If $Y\simeq Y'$ are homotopic topological spaces, there are canonical isomorphisms $H_k(Y;R)\cong H_k(Y';R)$ for all~$k,R$.
\smallskip

\noindent{\bf(e)} {\it Generalized homology theories}, such as bordism, are required to satisfy conditions (i)--(v) and possibly (vi), but not the Dimension Axiom~(vii). 
\smallskip

\noindent{\bf(f)} (Fundamental classes.) Suppose $Y$ is a compact, oriented manifold of dimension $m$. Then there is a natural {\it fundamental class\/} $[[Y]]$ in $H_m(Y;R)$, as in \cite[\S VI.8]{Bred1}, \cite[\S 3.3]{Hatc}, \cite[\S 55]{Munk}, and \cite[\S 16.4]{toDi}. If $R$ has characteristic 2 (e.g.\ $R=\Z_2$) then $[[Y]]$ is well-defined even if $Y$ is not oriented.
\smallskip

\noindent{\bf(g)} Following Eilenberg and Steenrod \cite[\S I.3]{EiSt}, we give axioms for $H_k(Y,Z;R)$, $H_k(Y;R)$ for all $k\in\Z$. In fact, as in \cite[p.~13]{EiSt}, every homology theory satisfying Axiom \ref{kh2ax1} has $H_k(Y,Z;R)=H_k(Y;R)=0$ for all $Y,Z$ and $k<0$, so many authors define homology as functors $H_k(-;R)$ for $k\in\N$ only, rather than~$k\in\Z$.

For our new homology theories $MH_*(-;R),\ldots$ in \S\ref{kh4}--\S\ref{kh5}, it will be natural to define $MH_k(Y,Z;R),MH_k(Y;R)$ for $k\in\Z$ rather than $k\in\N$, and it will not be at all obvious that $MH_k(Y,Z;R)=MH_k(Y;R)=0$ for $k<0$, but this will follow from Theorem \ref{kh2thm1}, once we verify Axiom~\ref{kh2ax1}.
\label{kh2rem1}
\end{rem}

The next theorem will be proved in \S\ref{kh61}, using results of Eilenberg and Steenrod \cite{EiSt}, Milnor \cite{Miln}, and Kreck and Singhof~\cite{KrSi}. 

\begin{thm} Any two homology theories $H_*(-;R),\ti H_*(-;R)$ satisfying Axiom\/ {\rm\ref{kh2ax1}} over the same commutative ring\/ $R$ are canonically isomorphic. That is, there exist\/ $R$-module isomorphisms $I_{Y,Z}:H_*(Y,Z;R)\ra\ti H_*(Y,Z;R)$ for all\/ $(Y,Z)\in\cC$ commuting with the given morphisms\/ $f_*,\pd$ and isomorphisms\/ $H_0(*;R)\cong R\cong\ti H_0(*;R),$ and any other assignment of morphisms\/ $J_{Y,Z}:H_*(Y,Z;R)\ra\ti H_*(Y,Z;R)$ for all\/ $(Y,Z)\in\cC$ commuting with the\/ $f_*,\pd$ and\/ $H_0(*;R)\cong R\cong\ti H_0(*;R)$ have\/ $J_{Y,Z}=I_{Y,Z}$ for all\/~$(Y,Z)$. 
\label{kh2thm1}
\end{thm}

Having defined our new homology theories of manifolds, we will verify they satisfy Axiom \ref{kh2ax1}(i)--(vii), and then Theorem \ref{kh2thm1} implies that they are canonically isomorphic to conventional homology groups such as singular homology.

We now define {\it singular homology\/} $H_*^\rsi(Y;R)$ for general topological spaces $Y$ (including manifolds), as in \cite[\S IV.1]{Bred1}, \cite[\S VII]{EiSt}, \cite[\S 4]{Maun}, \cite[\S 4]{Munk}, \cite[\S 4]{Span}, and {\it smooth singular homology\/} $H_*^\ssi(Y;R)$ for manifolds $Y$, as in Bredon~\cite[\S V.5]{Bred1}.

\begin{ex} For all $k\ge 0$, define the $k$-{\it simplex\/} to be
\e
\De_k=\bigl\{(x_0,\ldots,x_k)\in\R^{k+1}:x_i\ge 0,\;\>
x_0+\ldots+x_k=1\bigr\}.
\label{kh2eq5}
\e
It is a compact manifold with corners, of dimension $k$. Define an
orientation on $\De_k$ such that $(\d x_1\w\cdots\w\d x_k)\vert_{\De_k}$ is a positive volume form. We may write $\pd\De_k=\coprod_{j=0}^k \pd_j\De_k$, where $\pd_j\De_k$ is the boundary face of $\De_k$ upon which $x_j\equiv 0$. There is a unique diffeomorphism $F_j^k:\De_{k-1}\ra\pd_j\De_k$ for $j=0,\ldots,k$ with $i_{\De_k}\ci F_j^k:(x_0,\ldots,x_{k-1})\mapsto(x_0,\ldots,x_{j-1},0,x_j,\ldots,x_{k-1})$. Note that $\De_{k-1}$ and $\pd_j\De_k\subset\pd\De_k$ are oriented, and $F_j^k$ multiplies orientations by~$(-1)^j$.

Let $Y$ be a topological space, and $R$ a commutative ring. Write $C_k^\rsi(Y;R)$ for the free $R$-module spanned by {\it singular\/ $k$-simplices\/} $\si$ in $Y$, which are continuous maps $\si:\De_k\ra Y$. Elements of $C_k^\rsi(Y;R)$ will be written $\sum_{i\in I}\rho_i\,\si_i$, for $I$ a finite indexing set, $\rho_i\in R$ and $\si_i:\De_k\ra Y$ continuous for each $i\in I$. As in \cite[\S IV.1]{Bred1}, the boundary operator $\pd:C_k^\rsi(Y;R)\ra C_{k-1}^\rsi(Y;R)$ is
\e
\pd:\ts\sum_{i\in I}\rho_i\,\si_i\longmapsto \ts\sum_{i\in
I}\sum_{j=0}^k(-1)^j\rho_i\,(\si_i\ci i_{\De_k}\ci F_j^k).
\label{kh2eq6}
\e
Then $\pd^2=0$, and $H_*^\rsi(Y;R)$ is the homology of $\bigl(C_*^\rsi(Y;R),\pd\bigr)$, that is, 
\begin{equation*}
H_k^\rsi(Y;R)=\frac{\ts \Ker\bigl(\pd: C_k^\rsi(Y;R)\longra C_{k-1}^\rsi(Y;R)\bigr)}{\ts \Im\bigl(\pd:C_{k+1}^\rsi(Y;R)\longra C_k^\rsi(Y;R)\bigr)}\,.
\end{equation*}

If $Z\!\subseteq\! Y$ then $C_k^\rsi(Z;R)\!\subseteq\! C_k^\rsi(Y;R)$. Set $C_k^\rsi(Y,Z;R)\!=\!C_k^\rsi(Y;R)/C_k^\rsi(Z;R)$. Then $\pd$ descends to $\pd:C_k^\rsi(Y,Z;R)\ra C_{k-1}^\rsi(Y,Z;R)$, and {\it relative singular homology\/} $H_*^\rsi(Y,Z;R)$ is the homology of $\bigl(C_*^\rsi(Y,Z;R),\pd\bigr)$. We have a short exact sequence of chain complexes
\e
\xymatrix@C=15pt{ 0 \ar[r] & \bigl(C_*^\rsi(Z;R),\pd\bigr) \ar[r] & \bigl(C_*^\rsi(Y;R),\pd\bigr) \ar[r] & \bigl(C_*^\rsi(Y,Z;R),\pd\bigr) \ar[r] & 0. }
\label{kh2eq7}
\e
In the usual way this induces a long exact sequence \eq{kh2eq1} on (relative) singular homology, and defines the morphisms~$\pd:H_k^\rsi(Y,Z;R)\ra H_{k-1}^\rsi(Z;R)$.

Let $f:Y_1\ra Y_2$ be a continuous map of topological spaces. Define $f_*:C_k^\rsi(Y_1;R)\ra C_k^\rsi(Y_2;R)$ by
\e
f_*:\ts\sum_{i\in I}\rho_i\,\si_i\mapsto \ts\sum_{i\in I}\rho_i\,(f\ci\si_i).
\label{kh2eq8}
\e 
Then $f_*\ci\pd=\pd\ci f_*$, so $f_*$ induces morphisms of homology groups $f_*:H_k^\rsi(Y_1;R)\ra H_k^\rsi(Y_2;R)$. If $Z_1\subseteq Y_1$ with $f(Z_1)\subseteq Z_2\subseteq Y_2$ we get morphisms $f_*:H_k^\rsi(Y_1,Z_1;R)\ra H_k^\rsi(Y_2,Z_2;R)$ in the same way.

As in \cite[\S IV]{Bred1}, \cite[\S VII]{EiSt}, \cite[\S 30--\S 31]{Munk}, singular homology satisfies Axiom \ref{kh2ax1} for topological spaces, and therefore also for manifolds. 
\label{kh2ex1}
\end{ex}

\begin{ex} Let $Y$ be a smooth manifold, and $R$ a commutative ring. Use the notation of Example \ref{kh2ex1}. Write $C_k^\ssi(Y;R)$ for the free $R$-module spanned by {\it smooth singular\/ $k$-simplices\/} $\si$ in $Y$, which are smooth maps $\si:\De_k\ra Y$, regarding $\De_k$ as a manifold with corners. Since smooth maps are continuous we have $C_k^\ssi(Y;R)\subseteq C_k^\rsi(Y;R)$. Define $\pd:C_k^\ssi(Y;R)\ra C_{k-1}^\ssi(Y;R)$ as in \eq{kh2eq6}. {\it Smooth singular homology\/} $H_*^\ssi(Y;R)$ is the homology of~$\bigl(C_*^\ssi(Y;R),\pd\bigr)$.

If $Z\subseteq Y$ is open then $Z$ is a manifold and $C_k^\ssi(Z;R)\subseteq C_k^\ssi(Y;R)$. Write $C_k^\ssi(Y,Z;R)=C_k^\ssi(Y;R)/C_k^\ssi(Z;R)$, and define {\it relative smooth singular homology\/} $H_*^\ssi(Y,Z;R)$ to be the homology of $\bigl(C_*^\ssi(Y,Z;R),\pd\bigr)$. Define pushforwards $f_*$ and connecting morphisms $\pd$ as in Example~\ref{kh2ex1}.

The inclusion $\bigl(C_*^\ssi(Y;R),\pd\bigr)\hookra \bigl(C_*^\rsi(Y;R),\pd\bigr)$ induces morphisms on homology $H_*^\ssi(Y;R)\ra H_*^\rsi(Y;R)$, and similarly $H_*^\ssi(Y,Z;R)\ra H_*^\rsi(Y,Z;R)$. Bredon \cite[\S V.5 \& \S V.9]{Bred1} shows that these are isomorphisms, so (relative) smooth singular homology is canonically isomorphic to (relative) singular homology of manifolds. Therefore smooth singular homology satisfies Axiom~\ref{kh2ax1}.

Suppose $Y$ is compact and oriented of dimension $m$, so that as in Remark \ref{kh2rem1}(f) it has a fundamental class $[[Y]]\in H_m^\ssi(Y;R)$. We can construct this at the chain level as follows. Choose a triangulation of $Y$ which cuts it into finitely many smooth $m$-simplices, which we write as $\si_i:\De_m\ra Y$ for $i$ in $I$, a finite indexing set. Here each $\si_i$ embeds $\De_m$ into $Y$ as a submanifold with corners. Write $\ep_i=1$ if $\si_i$ is orientation-preserving, and $\ep_i=-1$ otherwise. Then $\pd\sum_{i\in I}\ep_i\,\si_i=0$, since each boundary face of each $m$-simplex $\si_i(\De_m)$ is cancelled by a boundary face of a neighbouring $m$-simplex in the triangulation, and the fundamental class is~$[[Y]]=\bigl[\sum_{i\in I}\ep_i\,\si_i\bigr]$.
\label{kh2ex2}
\end{ex}

Smooth singular homology is important in relating singular homology to de Rham cohomology of manifolds, as in Bredon \cite[\S V.5 \& \S V.9]{Bred1}, and will also be used below to map to our new homology theories. The next example discusses {\it barycentric subdivision\/} in (smooth) singular homology, following \cite[\S II.6]{EiSt}, \cite[\S IV.17]{Bred1},\cite[\S 15]{Munk}. It will be used in \S\ref{kh25} to define the (smooth) singular cosheaves. 

\begin{ex} Let $\De_k$ be the $k$-simplex in \eq{kh2eq5}. The {\it barycentre\/} of a simplex in an affine space is its `centre of gravity', the average of its vertices. For $\De_k$ in \eq{kh2eq5}, the barycentre is $\bigl(\frac{1}{k+1},\ldots,\frac{1}{k+1}\bigr)$. Similarly, every $l$-dimensional boundary face of $\De_k$ for $0\le l\le k$ is an $l$-simplex with its own barycentre, a point of $\De_k$. There are $2^{k+1}-1$ such boundary faces, giving $2^{k+1}-1$ points in $\De_k$. These are the vertices of a triangulation of $\De_k$ into $(k+1)!$ smaller affine $k$-simplices $\De_k^1,\De_k^2,\ldots,\De_k^{(k+1)!}$, called the {\it barycentric subdivision\/} of $\De_k$. 

Each $\De_k^j$ has volume $\frac{1}{(k+1)!}\cdot\vol(\De_k)$, and the $k+1$ vertices of $\De_k^j$ are the barycentres of one boundary face of $\De_k$ of each dimension $l=0,1,\ldots,k$. There is a unique affine-linear map $B_k^j:\De_k\ra\De_k$ for each $j=1,\ldots,(k+1)!$ such that $\De_k^j=B_k^j(\De_k)$, and writing $v_l=(0,\ldots,0,1,0,\ldots,0)$ for the vertex of $\De_k$ with coordinate $x_l=1$, then $B_k^j(v_l)$ is the barycentre of an $l$-dimensional boundary face of $\De_k$ for $l=0,\ldots,k$. Define $\ep_k^j=1$ if $B_k^j:\De_k\ra\De_k$ is orientation-preserving, and $\ep_k^j=-1$ if $B_k^j$ is orientation-reversing.

We can iterate barycentric subdivision, dividing each $k$-simplex $\De_k^i\subset\De_k$ into $(k+1)!$ smaller simplices, and so on. In this way we can cut $\De_k$ into many arbitrarily small sub-simplices. Figure \ref{kh2fig1} illustrates the first two barycentric subdivisions of the 2-simplex (triangle)~$\De_2$.

\begin{figure}[htb]
\centerline{$\splinetolerance{.8pt}
\begin{xy}
0;<.6mm,0mm>:
,(-95,0)*{\bu}
,(-65,60)*{\bu}
,(-35,0)*{\bu}
,(-95,0);(-65,60)**\crv{}
,(-35,0);(-65,60)**\crv{}
,(-95,0);(-35,0)**\crv{(-60,0)}
,(-30,0)*{\bu}
,(0,60)*{\bu}
,(30,0)*{\bu}
,(0,0)*{\bu}
,(0,20)*{\bu}
,(-15,30)*{\bu}
,(15,30)*{\bu}
,(-30,0);(0,60)**\crv{}
,(30,0);(0,60)**\crv{}
,(-30,0);(30,0)**\crv{(0,0)}
,(-30,0);(0,20)**\crv{}
,(30,0);(0,20)**\crv{}
,(0,0);(0,20)**\crv{(0,10)}
,(-15,30);(0,20)**\crv{}
,(15,30);(0,20)**\crv{}
,(0,60);(0,20)**\crv{(0,30)}
,(35,0)*{\bu}
,(65,60)*{\bu}
,(95,0)*{\bu}
,(65,0)*{\bu}
,(65,20)*{\bu}
,(50,30)*{\bu}
,(80,30)*{\bu}
,(50,10)*{\bu}
,(80,10)*{\bu}
,(65,10)*{\bu}
,(57.5,25)*{\bu}
,(72.5,25)*{\bu}
,(65,40)*{\bu}
,(42.5,15)*{\bu}
,(57.5,45)*{\bu}
,(87.5,15)*{\bu}
,(72.5,45)*{\bu}
,(80,0)*{\bu}
,(50,0)*{\bu}
,(95,0);(65,60)**\crv{}
,(35,0);(65,60)**\crv{}
,(95,0);(35,0)**\crv{(60,0)}
,(55,6.667)*{\bu}
,(75,6.667)*{\bu}
,(50,16.667)*{\bu}
,(80,16.667)*{\bu}
,(60,36.667)*{\bu}
,(70,36.667)*{\bu}
,(35,0);(65,20)**\crv{}
,(95,0);(65,20)**\crv{}
,(65,0);(65,20)**\crv{(65,10)}
,(50,30);(65,20)**\crv{}
,(80,30);(65,20)**\crv{}
,(65,60);(65,20)**\crv{(65,30)}
,(35,0);(65,10)**\crv{}
,(50,0);(65,20)**\crv{}
,(95,0);(65,10)**\crv{}
,(80,0);(65,20)**\crv{}
,(65,0);(50,10)**\crv{}
,(65,0);(80,10)**\crv{}
,(50,10);(50,30)**\crv{(50,20)}
,(80,10);(80,30)**\crv{(80,20)}
,(42.5,15);(65,20)**\crv{}
,(87.5,15);(65,20)**\crv{}
,(57.5,25);(35,0)**\crv{}
,(72.5,25);(95,0)**\crv{}
,(57.5,25);(65,60)**\crv{}
,(72.5,25);(65,60)**\crv{}
,(50,30);(65,40)**\crv{}
,(80,30);(65,40)**\crv{}
,(57.5,45);(65,20)**\crv{}
,(72.5,45);(65,20)**\crv{}
\end{xy}$}
\caption{The 2-simplex $\De_2$ and its first two barycentric subdivisions}
\label{kh2fig1}
\end{figure}
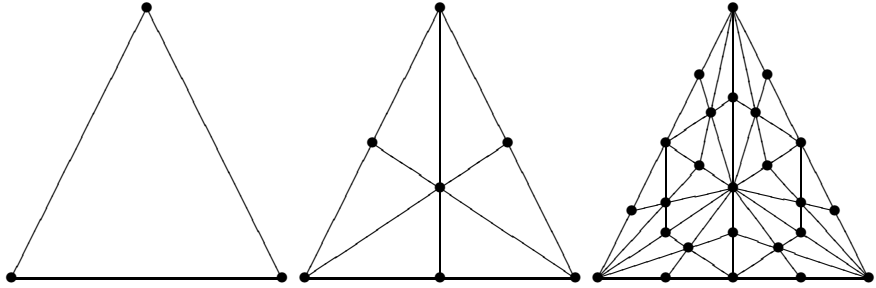

Barycentric subdivision is used in (smooth) singular homology as follows. Let $Y$ be a topological space and $R$ a commutative ring, and define $C_k^\rsi(Y;R)$ as in Example \ref{kh2ex1}. Define an $R$-linear map $B:C_k^\rsi(Y;R)\ra C_k^\rsi(Y;R)$ by
\e
B:\ts\sum_{i\in I}\rho_i\,\si_i\longmapsto \ts\sum_{i\in
I}\sum_{j=1}^{(k+1)!}\ep_k^j\rho_i\,(\si_i\ci B_k^j).
\label{kh2eq9}
\e
We call $B$ the {\it barycentric subdivision morphism}. It has the important property that $B\ci\pd=\pd\ci B:C_k^\rsi(Y;R)\ra C_{k-1}^\rsi(Y;R)$, since the boundary faces of $\De_k^i$ for $i=1,\ldots,(k+1)!$ can be divided into two kinds, (A) those coming from the barycentric subdivision of a boundary face of $\De_k$, and (B) those coming from the interior $\De_k^\ci$, which cancel out in pairs in \eq{kh2eq9}.

Thus, $B$ induces morphisms $B_*:H_k^\rsi(Y;R)\ra H_k^\rsi(Y;R)$ for all $k=0,1,\ldots.$ By considering a subdivision of $\De_k\t[0,1]$ into $(k+1)$-simplices, one can show that these $B_*$ are the identity maps. Also if $f:Y_1\ra Y_2$ is a continuous map of topological spaces then~$B\ci f_*=f_*\ci B:C_k^\rsi(Y_1;R)\ra C_k^\rsi(Y_2;R)$.

In the same way, if $Y$ is a manifold, $R$ a commutative ring, and $C_k^\ssi(Y;R)$ is as in Example \ref{kh2ex2}, we define $B:C_k^\ssi(Y;R)\ra C_k^\ssi(Y;R)$, and then $B\ci\pd=\pd\ci B:C_k^\ssi(Y;R)\ra C_{k-1}^\ssi(Y;R)$, and the induced maps on smooth singular homology $H_k^\ssi(Y;R)$ are~$B_*=\id:H_k^\ssi(Y;R)\ra H_k^\ssi(Y;R)$.

Barycentric subdivision is used to prove that $H_*^\rsi(-;R),H_*^\ssi(-;R)$ satisfy the excision axiom, Axiom \ref{kh2ax1}(v). The main idea is that if $Z\subseteq Y$ is open, and $S\subseteq Z$ is closed in $Y$, then given a chain $\sum_{i\in I}\rho_i\,\si_i$ in $C_k^\rsi(Y;R)$, we would like to write it as the sum of chains in $C_k^\rsi(Y\sm S;R)$ and $C_k^\rsi(Z;R)$. In general this is impossible, as $\rho_i(\De_k)\subseteq Y$ for $i\in I$ need not be contained in $Y\sm S$ or $Z$. However, if we instead consider $B^n\bigl(\sum_{i\in I}\rho_i\,\si_i\bigr)$ for $n\gg 0$, this is the sum of many very small $k$-simplices in $Y$, each of which lies in either $Y\sm Z$ or $Z$ or both, so $B^n\bigl(\sum_{i\in I}\rho_i\,\si_i\bigr)$ lies in~$C_k^\rsi(Y\sm S;R)+C_k^\rsi(Z;R)\subseteq C_k^\rsi(Y;R)$.
\label{kh2ex3}
\end{ex}

\subsection{Cohomology}
\label{kh22}

Cohomology is dual to homology. A {\it cohomology theory\/} over a commutative ring $R$ associates to any topological space $Y$ (possibly satisfying conditions) a sequence of $R$-modules $H^0(Y;R),H^1(Y;R),\ldots,$ with $H^0(*;R)=R$, and to any continuous map $f:Y_1\ra Y_2$ a sequence of morphisms $f^*:H^k(Y_2;R)\ra H^k(Y_1;R)$ for $k=0,1,\ldots,$ satisfying a package of properties we explain below. As for homology, cohomology of topological spaces is homotopy invariant.

Some examples of cohomology theories are singular cohomology \cite[\S VII.2]{EiSt}, \cite[\S 44]{Munk}, \cite[\S 5.4]{Span}, \v Cech cohomology \cite[\S IX]{EiSt}, \cite[\S 73]{Munk}, Alexander--Spanier cohomology \cite[\S 6]{Span}, de Rham cohomology of manifolds \cite[\S V]{Bred1}, and sheaf cohomology $H^*(Y,R_Y)$ as in \S\ref{kh25}. This section concerns cohomology groups only as $R$-modules; the ring structure on cohomology will be discussed in \S\ref{kh26}. Here is the analogue of Axiom~\ref{kh2ax1}.
 
\begin{ax}[Eilenberg--Steenrod axioms for cohomology of manifolds] Let $\cC$ be as in Axiom \ref{kh2ax1}, and write $\cC^{\rm op}$ for the opposite category. Fix a commutative ring $R$. A {\it cohomology theory over\/} $R$ assigns the data:
\begin{itemize}
\setlength{\itemsep}{0pt}
\setlength{\parsep}{0pt}
\item[(a)] For each object $(Y,Z)$ in $\cC$ and each $k\in\Z$, an $R$-module $H^k(Y,Z;R)$. For brevity we write~$H^k(Y;R)=H^k(Y,\es;R)$.
\item[(b)] For each object $(Y,Z)$ in $\cC$ and each $k\in\Z$, an $R$-module morphism $\d:H^k(Z;R)\ra H^{k+1}(Y,Z;R)$, called the {\it connecting morphism}.
\item[(c)] For each morphism $f:(Y_1,Z_1)\ra (Y_2,Z_2)$ in $\cC$ and each $k\in\Z$, an $R$-module morphism $f^*:H^k(Y_2,Z_2;R)\ra H^k(Y_1,Z_1;R)$.
\end{itemize}
All this data must satisfy the axioms:
\begin{itemize}
\setlength{\itemsep}{0pt}
\setlength{\parsep}{0pt}
\item[(i)] (Functoriality.) For each $k\in\Z$, mapping $(Y,Z)$ to $H^k(Y,Z;R)$ and $f:(Y_1,Z_1)\ra (Y_2,Z_2)$ to $f^*:H^k(Y_2,Z_2;R)\ra H^k(Y_1,Z_1;R)$ gives a functor $\cC^{\rm op}\ra\Rmod$, where $\Rmod$ is the abelian category of $R$-modules.
\item[(ii)] (Exactness.) For each $(Y,Z)\in\cC$, write $i:(Z,\es)\ra(Y,\es)$ and $j:(Y,\es)\ra (Y,Z)$ for the morphisms in $\cC$ induced by the inclusion $Z\hookra Y$ and the identity $\id_Y:Y\ra Y$. Then the following is exact in $\Rmod$:
\end{itemize}
\e
\xymatrix@C=11pt{ \cdots \ar[r] & H^{k-1}(Z;R) \ar[r]^(0.48)\d & H^k(Y,Z;R) \ar[r]^(0.52){j^*} & H^k(Y;R) \ar[r]^(0.47){i^*} & H^k(Z;R)    \ar[r] & \cdots. }
\label{kh2eq10}
\e
\begin{itemize}
\setlength{\itemsep}{0pt}
\setlength{\parsep}{0pt}
\item[(iii)] (Functoriality of $\d$.) For each $f:(Y_1,Z_1)\ra (Y_2,Z_2)$ in $\cC$ and $k\in\Z$, the following commutes:
\begin{equation*}
\xymatrix@C=90pt@R=14pt{ *+[r]{H^{k-1}(Z_2;R)} \ar[d]^{f\vert_{Z_1}^*} \ar[r]_\d & *+[l]{H^k(Y_2,Z_2;R)} \ar[d]_{f^*} \\ *+[r]{H^{k-1}(Z_1;R)} \ar[r]^\d & *+[l]{H^k(Y_1,Z_1;R).\!\!} }
\end{equation*}
\item[(iv)] (Homotopy.) Suppose $(Y_1,Z_1),(Y_2,Z_2)\in\cC$, and $g:Y_1\t[0,1]\ra Y_2$ is smooth with $g(Z_1\t[0,1])\subseteq Z_2$. Define $f,f':Y_1\ra Y_2$ by $f(y)=g(y,0)$ and $f'(y)=g(y,1)$ for $y\in Y_1$. Then $f^*=f^{\prime *}:H^k(Y_2,Z_2;R)\ra H^k(Y_1,Z_1;R)$ for all $k\in\Z$.
\item[(v)] (Excision.) Suppose $(Y,Z)\in\cC$, and $S\subseteq Z$ is closed in $Y$, so the inclusion $j:Y\sm S\hookra Y$ is a morphism $j:(Y\sm S,Z\sm S)\ra (Y,Z)$ in $\cC$. Then $j^*:H^k(Y,Z;R)\ra H^k(Y\sm S,Z\sm S;R)$ is an isomorphism for all~$k\in\Z$.
\item[(vi)] (Additivity Axiom.) Suppose $(Y,Z)\in\cC$ with $Y=\coprod_{a\in A}Y_a$ for $A$ a countable indexing set and each $Y_a$ open and closed in $Y$, and set $Z_a=Z\cap Y_a$. Then for all $k\in\Z$ we have canonical isomorphisms
\e
H^k(Y,Z;R)\cong \ts\prod_{a\in A}H^k(Y_a,Z_a;R),
\label{kh2eq11}
\e
compatible with the morphisms $i_a^*:H^k(Y,Z;R)\ra H^k(Y_a,Z_a;R)$ for $a\in A$ induced by the inclusion $i_a:(Y_a,Z_a)\hookra(Y,Z)$.
\item[(vii)] (Dimension Axiom.) When $Y=*$ is the point we have $H^0(*;R)\cong R$ and $H^k(*;R)=0$ for all~$k\ne 0$.
\end{itemize}
\label{kh2ax2}
\end{ax}

Note that direct sums in \eq{kh2eq2} are replaced by direct products in \eq{kh2eq11}, as in \eq{kh2eq3}--\eq{kh2eq4}. Here is the analogue of Theorem \ref{kh2thm1}, also proved in~\S\ref{kh61}.

\begin{thm} Any two cohomology theories $H^*(Y,Z;R),\ti H^*(Y,Z;R)$ satisfying Axiom\/ {\rm\ref{kh2ax2}} over the same commutative ring\/ $R$ are canonically isomorphic.

\label{kh2thm2}
\end{thm}

As in Remark \ref{kh2rem1}(g), any cohomology theory $H^*(-;R)$ satisfying Axiom \ref{kh2ax2} has $H^k(Y,Z;R)=H^k(Y;R)=0$ for all $Y,Z$ and $k<0$, so many authors define cohomology as functors $H^k(-;R)$ for $k\in\N$ rather than~$k\in\Z$.

\begin{ex} Let $\K$ be a field, and suppose we are given a homology theory over $R=\K$ in the sense of Axiom \ref{kh2ax1}, with homology groups $H_k(Y,Z;\K)$. For all $(Y,Z)\in\cC$ and $k\ge 0$, define $H^k(Y,Z;\K)=H_k(Y,Z;\K)^*$ to be the dual $\K$-vector space. Let $\d:H^{k-1}(Z;\K)\ra H^k(Y,Z;\K)$ be the dual linear map to $\pd:H_k(Y,Z;\K)\ra H_{k-1}(Z;\K)$, and for $f:(Y_1,Z_1)\ra(Y_2,Z_2)$, let $f^*:H^k(Y_2,Z_2;\K)\ra H^k(Y_1,Z_1;\K)$ be dual to~$f_*:H_k(Y_1,Z_1;\K)\ra H_k(Y_2,Z_2;\K)$. Then Axiom \ref{kh2ax1}(i)--(vii) imply Axiom \ref{kh2ax2}(i)--(vii) by duality. So the dual of a homology theory over a field $\K$ is a cohomology theory over~$\K$.

Note that for the Additivity Axioms (vi), this is subtle. If $V_a$, $a\in A$ is a family of $\K$-vector spaces, using \eq{kh2eq3}--\eq{kh2eq4} we define a canonical isomorphism
\begin{equation*}
\bigl[\ts \bigop_{a\in A}V_a\bigr]{}^*\cong \prod_{a\in A}(V_a^*).
\end{equation*}
But as we can have $V\not\cong (V^*)^*$ for infinite-dimensional $\K$-vector spaces $V$, in general for $A$ infinite we may have
\begin{equation*}
\bigl[\ts \prod_{a\in A}V_a\bigr]{}^*\not\cong \bigop_{a\in A}(V_a^*).
\end{equation*}
Thus, the dual of a cohomology theory over $\K$ is {\it not\/} a homology theory, and $H_k(Y,Z;\K)\not\cong H^k(Y,Z;\K)^*$ when $H_k(Y,Z;\K)$ is infinite-dimensional.
\label{kh2ex4}
\end{ex}

Following \cite[\S VII.2]{EiSt}, \cite[\S 44]{Munk}, \cite[\S 5.4]{Span} we define {\it singular cohomology\/} of topological spaces, which is dual to singular homology in Example~\ref{kh2ex1}.

\begin{ex} Fix a commutative ring $R$. Let $Y$ be a topological space and $Z\subseteq Y$ a subspace. Then Example \ref{kh2ex1} defines the singular chain complex $\bigl(C_*^\rsi(Y,Z;\Z),\pd\bigr)$ over $\Z$. For each $k\ge 0$, define an $R$-module $C^k_\rsi(Y,Z;R)=\Hom_\Z\bigl(C_k^\rsi(Y,Z;\Z),R\bigr)$ of {\it singular cochains}. 

For $k\ge 0$ define $\d:C^k_\rsi(Y,Z;R)\ra C^{k+1}_\rsi(Y,Z;R)$ by $\d\al=\al\ci\pd$, for $\al:C_k^\rsi(Y,Z;\Z)\ra R$ and $\pd:C_{k+1}^\rsi(Y,Z;\Z)\ra C_k^\rsi(Y,Z;\Z)$ as in Example \ref{kh2ex1}. Then $\d^2=0$ as $\pd^2=0$, so $\bigl(C^*_\rsi(Y,Z;R),\d\bigr)$ is a cochain complex of $R$-modules. Define the {\it singular cohomology\/} $H^k_\rsi(Y,Z;R)$ to be the $k^{\rm th}$ cohomology group of $\bigl(C^*_\rsi(Y,Z;R),\d\bigr)$. Starting with \eq{kh2eq7} for $R=\Z$ and applying $\Hom_\Z(-,R)$ gives a short exact sequence of cochain complexes
\begin{equation*}
\xymatrix@C=15pt{ 0 \ar[r] & \bigl(C^*_\rsi(Y,Z;R),\d\bigr)  \ar[r] & \bigl(C^*_\rsi(Y;R),\d\bigr) \ar[r] &  \ar[r] \bigl(C^*_\rsi(Z;R),\d\bigr) & 0. }
\end{equation*}
In the usual way this induces a long exact sequence \eq{kh2eq10} on (relative) singular cohomology, and defines the morphisms~$\d:H^k_\rsi(Y,Z;R)\ra H^{k+1}_\rsi(Z;R)$.

Suppose $f:Y_1\ra Y_2$ is a continuous map of topological spaces, and $Z_1\subseteq Y_1$ with $f(Z_1)\subseteq Z_2\subseteq Y_2$. Then Example \ref{kh2ex1} defines $f_*:\bigl(C_*^\rsi(Y_1,Z_1;\Z),\pd\bigr)\ra \bigl(C_*^\rsi(Y_2,Z_2;\Z),\pd\bigr)$. Applying $\Hom_\Z(-,R)$ gives a morphism of cochain complexes $f^*:\bigl(C^*_\rsi(Y_2,Z_2;R),\d\bigr)\ra\bigl(C^*_\rsi(Y_1,Z_1;R),\d\bigr)$. Let $f^*:H^k_\rsi(Y_2,Z_2;R)\ra H^k_\rsi(Y_1,Z_1;R)$ be the induced morphism on $k^{\rm th}$ cohomology groups.

As in \cite[\S VII]{EiSt}, \cite[\S 44]{Munk}, \cite[\S 5.4]{Span}, singular cohomology satisfies Axiom \ref{kh2ax2} for topological spaces, and therefore also for manifolds.

Note that although singular chains in $C_k^\rsi(Y;R)$ are simple objects -- finite $R$-linear combinations of continuous maps $\si:\De_k\ra Y$ -- singular cochains in $C^k_\rsi(Y;R)$ are $\Z$-linear maps $C_k^\rsi(Y;\Z)\ra R$, where $C_k^\rsi(Y;\Z)$ is generally huge, so singular cochains are difficult to write down and work with explicitly.
\label{kh2ex5}
\end{ex}

\begin{ex} Fix a commutative ring $R$. Let $Y$ be a smooth manifold and $Z\subseteq Y$ an open set. Then Example \ref{kh2ex2} defines the smooth singular chain complex $\bigl(C_*^\ssi(Y,Z;\Z),\pd\bigr)$ over $\Z$, a subcomplex of $\bigl(C_*^\rsi(Y,Z;\Z),\pd\bigr)$ in Example \ref{kh2ex1}. For each $k\ge 0$, define an $R$-module $C^k_\ssi(Y,Z;R)=\Hom_\Z\bigl(C_k^\ssi(Y,Z;\Z),R\bigr)$ of {\it smooth singular cochains}. Define $\d:C^k_\ssi(Y,Z;R)\ra C^{k+1}_\ssi(Y,Z;R)$ as in Example \ref{kh2ex5}, and let the {\it smooth singular cohomology\/} $H^k_\ssi(Y,Z;R)$ be the $k^{\rm th}$ cohomology group of $\bigl(C^*_\ssi(Y,Z;R),\d\bigr)$. Define $\d:H^k_\ssi(Y,Z;R)\ra H^{k+1}_\ssi(Z;R)$ and $f^*:H^k_\ssi(Y_2,Z_2;R)\ra H^k_\ssi(Y_1,Z_1;R)$ for $f:(Y_1,Z_1)\ra(Y_2,Z_2)$ a morphism in $\cC$ as in Example \ref{kh2ex5}. 

Example \ref{kh2ex2} gives an inclusion $\bigl(C_*^\ssi(Y,Z;\Z),\pd\bigr)\hookra\bigl(C_*^\rsi(Y,Z;\Z),\pd\bigr)$ of chain complexes inducing isomorphisms $H_k^\ssi(Y,Z;\Z)\ra H_k^\rsi(Y,Z;\Z)$ on homology. Applying $\Hom(-,R)$ gives a surjection $\bigl(C^*_\rsi(Y,Z;R),\d\bigr)\twoheadrightarrow\bigl(C^*_\ssi(Y,Z;R),\d\bigr)$ inducing isomorphisms $H^k_\rsi(Y,Z;R)\ra H^k_\ssi(Y,Z;R)$ on cohomology. These isomorphisms commute with $\d,f^*$. Thus, since singular cohomology satisfies Axiom \ref{kh2ex2}, smooth singular cohomology also satisfies Axiom~\ref{kh2ex2}.
\label{kh2ex6}
\end{ex}

We discuss de Rham cohomology, as in Bott and Tu \cite[\S I]{BoTu} and Bredon~\cite[\S V]{Bred1}.

\begin{ex} Work over $R=\R$. Let $Y$ be a smooth manifold. Define the {\it de Rham cochains\/} $C^k_\dR(Y;\R)=C^\iy(\La^kT^*Y)$, the smooth $k$-forms on $Y$, for $k=0,1,2,\ldots,$ where $C^k_\dR(Y;\R)=0$ for $k>\dim Y$. Define $\d:C^k_\dR(Y;\R)\ra C^{k+1}_\dR(Y;\R)$ to be the usual exterior derivative on $k$-forms. Then $\bigl(C^*_\dR(Y;\R),\d\bigr)$ is a cochain complex over $\R$. Define the {\it de Rham cohomology group\/} $H^k_\dR(Y;\R)$ to be the $k^{\rm th}$ cohomology group of~$\bigl(C^*_\dR(Y;\R),\d\bigr)$.

Let $Z\subseteq Y$ be open. As in \cite[p.~78-79]{BoTu}, define the {\it relative de Rham cochains\/}
\e
C^k_\dR(Y,Z;\R)=C^\iy(\La^kT^*Y)\op C^\iy(\La^{k-1}T^*Z),
\label{kh2eq12}
\e
and define $\d:C^k_\dR(Y,Z;\R)\ra C^{k+1}_\dR(Y,Z;\R)$ by $\d(\al,\be)=(\d\al,\al\vert_Z-\d\be)$. Then $\bigl(C^*_\dR(Y,Z;\R),\d\bigr)$ is a cochain complex over $\R$. Define $H^k_\dR(Y,Z;\R)$ to be the $k^{\rm th}$ cohomology group of $\bigl(C^*_\dR(Y,Z;\R),\d\bigr)$. Then \cite[Prop.~6.49]{BoTu} says that \eq{kh2eq10} is exact, where $\d:H^{k-1}_\dR(Z;\R)\ra H^k_\dR(Y,Z;\R)$ maps $\d:[\be]\mapsto[0,\be]$.

If $f:(Y_1,Z_1)\ra (Y_2,Z_2)$ is a morphism in $\cC$, define $f^*:C^k_\dR(Y_2,Z_2;\R)\ra C^k_\dR(Y_1,Z_1;\R)$ by $f^*:(\al,\be)\mapsto(f^*(\al),f^*(\be))$. Then $f^*\ci\d=\d\ci f^*$, so $f^*$ induces morphisms of cohomology groups $f^*:H^k_\dR(Y_2,Z_2;\R)\ra H^k_\dR(Y_1,Z_1;\R)$.

One can show as in \cite[\S V]{Bred1} that de Rham cohomology satisfies Axiom \ref{kh2ax2},  and so Theorem \ref{kh2thm2} implies that de Rham cohomology is canonically isomorphic to (smooth) singular cohomology over $\R$ in Examples \ref{kh2ex5} and \ref{kh2ex6}. Thus, for each manifold $Y$ we have functorial isomorphisms
\e
H^k_\dR(Y;\R)\cong H^k_\ssi(Y;\R).
\label{kh2eq13}
\e

We can realize these at the level of cochain complexes: define
\e
\Pi:C^k_\dR(Y;\R)\longra C^k_\ssi(Y;\R)\;\>\text{by}\;\> \Pi(\om):\sum_{i\in I}a_i\si_i\longmapsto \sum_{i\in I}a_i\int_{\De_k}\si_i^*(\om)
\label{kh2eq14}
\e
for all $\om\in C^k_\dR(Y;\R)=C^\iy(\La^kT^*Y)$ and $\sum_{i\in I}a_i\si_i\in C_k^\ssi(Y;\Z)$, so that $I$ is finite, $a_i\in\Z$ and $\si_i:\De_k\ra Y$ is smooth for $i\in I$. Clearly, $\Pi(\om)$ is a $\Z$-linear map $C_k^\ssi(Y;\Z)\ra\R$, so that $\Pi(\om)\in  C^k_\ssi(Y;\R)$. Then $\Pi\ci\d=\d\ci\Pi$, so that $\Pi:\bigl(C^*_\dR(Y;\R),\d\bigr)\ra\bigl(C^*_\ssi(Y;\R),\d\bigr)$ is a morphism of cochain complexes, and induces $\Pi_*:H^*_\dR(Y;\R)\ra H^*_\ssi(Y;\R)$, which are isomorphisms as in \eq{kh2eq13}. This is known as {\it de Rham's Theorem}, and is proved directly in~\cite[\S V.9]{Bred1}.
\label{kh2ex7}
\end{ex}

\subsection{Compactly-supported cohomology}
\label{kh23}

{\it Compactly-supported cohomology\/} $H^*_\cs(Y;R)$ (or {\it cohomology of the second kind\/}) is a variation on cohomology $H^*(Y;R)$. Some compactly-supported cohomology theories are compactly-supported singular cohomology \cite[\S 3.3]{Hatc}, \cite[p.~323]{Span}, compactly-supported Alexander--Spanier cohomology \cite[\S 1]{Mass2}, \cite[p.~320]{Span}, com\-pac\-tly-supported de Rham cohomology of manifolds \cite[\S 1.1]{BoTu}, and compactly-supported sheaf cohomology $H^*_\cs(Y,R_Y)$ as in~\S\ref{kh25}.

For (nice) Hausdorff topological spaces $Y$ we can define compactly-supported cohomology using relative cohomology, by 
\e
H^k_\cs(Y;R)\cong H^k(Y\amalg\{\iy\},\{\iy\};R),
\label{kh2eq15}
\e
where $Y\amalg\{\iy\}$ is the {\it one-point compactification\/} of $Y$, with open sets $U$ for $U\subseteq Y$ open and $U\amalg\{\iy\}$ for $U\subseteq Y$ with $Y\sm U$ compact. An alternative to \eq{kh2eq15} is the more complicated formula 
\e
H^k_\cs(Y;R)\cong \underrightarrow{\lim}_{\,\text{$Z\subseteq Y$: $Y\sm Z$ is compact}\,} H^k(Y,Z;R).
\label{kh2eq16}
\e

Here $\underrightarrow{\lim}$ is a {\it direct limit\/} in the category $\Rmod$. That is, for each open set $Z\subseteq Y$ with $Y\sm Z$ compact there is a natural $R$-module morphism
\e
I_{Y,Z}:H^k(Y,Z;R)\longra H^k_\cs(Y;R),
\label{kh2eq17}
\e
and if $Z_2\subseteq Z_1\subseteq Y$ with $Y\sm Z_1,Y\sm Z_2$ compact then 
the following commutes:
\begin{equation*}
\xymatrix@C=50pt@R=5pt{
H^k(Y,Z_1;R) \ar[rr]_{I_{Y,Z_1}} \ar[dr]_{\id_Y^*} && H^k_\cs(Y;R), \\
& H^k(Y,Z_2;R) \ar[ur]_{I_{Y,Z_2}} }
\end{equation*}
where $\id_Y^*:H^k(Y,Z_1;R)\ra H^k(Y,Z_2;R)$ comes from $\id_Y:(Y,Z_2)\ra(Y,Z_1)$ in $\cC$ by Axiom \ref{kh2ax2}(c), and $H^k_\cs(Y;R)$ is universal with these properties.

For compactly-supported cohomology of manifolds $Y$, equation \eq{kh2eq16} may be more useful than \eq{kh2eq15}, since $Y,Z$ in \eq{kh2eq16} are manifolds, but $Y\amalg\{\iy\}$ in \eq{kh2eq15} is usually not a manifold.

Later we will define some new cohomology theories of manifolds $MH^*(-;R),\ab MH^*_\Q(-;R),MH^*_\dR(-;\R)$, and corresponding compactly-supported cohomology theories of manifolds $MH^*_\cs(-;R), MH^*_{\cs,\Q}(-;R),MH^*_{\cs,\dR}(-;\R)$. We will first prove $MH^*(-;R),\ldots$ satisfy Axiom \ref{kh2ax2}, and so are canonically isomorphic to conventional cohomology of manifolds by Theorem \ref{kh2thm2}. Then we will use facts about sheaf cohomology in \S\ref{kh25} to deduce that $MH^*_\cs(-;R),\ldots$ are canonically isomorphic to conventional compactly-supported cohomology of manifolds.

Because of this, we will not need an axiomatic characterization of compactly-supported cohomology (see Petkova \cite{Petk} and Skljarenko \cite{Sklj1} for results of this type), or even a formal definition. For simplicity, we will also not discuss relative compactly-supported cohomology, as we will not need it.

Here are some useful properties of compactly-supported cohomology of manifolds. It is not intended to be a complete definition, or an axiomatization.

\begin{property} Let $R$ be a commutative ring, and fix a cohomology theory of manifolds $H^*(-;R)$ as in \S\ref{kh22}. A {\it compactly-supported cohomology theory of manifolds\/}  $H^*_\cs(-;R)$ gives an $R$-module $H^k_\cs(Y;R)$ for each manifold $Y$ and $k\in\N$, in a direct limit with morphisms $I_{Y,Z}$ as in \eq{kh2eq16}--\eq{kh2eq17}. Furthermore:
\begin{itemize}
\setlength{\itemsep}{0pt}
\setlength{\parsep}{0pt}
\item[(a)] There are natural $R$-module morphisms $\Pi:H^k_\cs(Y;R)\ra H^k(Y;R)$ for all $Y,k$. If $Y$ is compact then $\Pi:H^k_\cs(Y;R)\ra H^k(Y;R)$ is an isomorphism.

In terms of the direct limit \eq{kh2eq16}, $\Pi$ is the unique morphism such that the following commutes for all $Z\subseteq Y$ with $Y\sm Z$ compact:
\begin{equation*}
\xymatrix@C=50pt@R=5pt{
H^k(Y,Z;R) \ar[dr]_{I_{Y,Z}} \ar[rr]_{\id_Y^*} && H^k(Y;R). \\
& H^k_\cs(Y;R) \ar[ur]_{\Pi} }
\end{equation*}
\item[(b)] Let $f:Y_1\ra Y_2$ be a proper, smooth map of manifolds, and $k\in\N$. Then there is a unique $R$-module morphism $f^*:H^k_\cs(Y_2;R)\ra H^k_\cs(Y_1;R)$ such that the following commutes for all $Z_2\subseteq Y_2$ with $Y_2\sm Z_2$ compact:
\begin{equation*}
\xymatrix@C=130pt@R=13pt{
*+[r]{H^k(Y_2,Z_2;R)} \ar[d]^{f^*} \ar[r]_(0.55){I_{Y_2,Z_2}} & *+[l]{H^k_\cs(Y_2;R)} \ar[d]_{f^*} \\
*+[r]{H^k(Y_1,f^{-1}(Z_2);R)} \ar[r]^(0.55){I_{Y_1,f^{-1}(Z_2)}} & *+[l]{H^k_\cs(Y_1;R).\!\!} }
\end{equation*}
Here we need $f$ proper to ensure that $Y_1\sm f^{-1}(Z_2)$ is compact. The following also commutes:
\begin{equation*}
\xymatrix@C=130pt@R=13pt{
*+[r]{H^k_\cs(Y_2;R)} \ar[d]^{f^*} \ar[r]_\Pi & *+[l]{H^k(Y_2;R)} \ar[d]_{f^*} \\
*+[r]{H^k_\cs(Y_1;R)} \ar[r]^\Pi & *+[l]{H^k(Y_1;R).\!\!} }
\end{equation*}
If $f:Y_1\ra Y_2$, $g:Y_2\ra Y_3$ are proper, smooth maps then $(g\ci f)^*=f^*\ci g^*:H^k_\cs(Y_3;R)\ra H^k_\cs(Y_1;R)$. That is, compactly-supported cohomology is a {\it contravariant\/} functor under proper, smooth maps $f:Y_1\ra Y_2$.
\item[(c)] Let $Y$ be a manifold, $U\subseteq Y$ an open set, and $k\in\N$. Write $i:U\hookra Y$ for the inclusion. Then there is a unique $R$-module morphism $i_*:H^k_\cs(U;R)\ra H^k_\cs(Y;R)$ such that the following commutes for all $Z\subseteq Y$ with $Y\sm Z$ compact and $Y=U\cup Z$, writing $Z'=U\cap Z$ and $S=Y\sm U$:
\begin{equation*}
\xymatrix@C=140pt@R=15pt{
*+[r]{H^k(Y\sm S,Z\sm S;R)=H^k(U,Z';R)} \ar[d]^{(e^*)^{-1}} \ar[r]_(0.73){I_{U,Z'}} & *+[l]{H^k_\cs(U;R)} \ar[d]_{i_*} \\
*+[r]{H^k(Y,Z;R)} \ar[r]^(0.55){I_{Y,Z}} & *+[l]{H^k_\cs(Y;R).\!\!} }
\end{equation*}
Here $e^*:H^k(Y,Z;R)\ra H^k(Y\sm S,\ab Z\sm S;R)$ is an isomorphism by Axiom \ref{kh2ax2}(v) (excision), with $e:Y\sm S\hookra Y$ the inclusion.

If $j:U'\hookra U$ is an inclusion of open sets then $(i\ci j)_*=\ab i_*\ci j_*:H^k_\cs(U';R)\ra H^k_\cs(Y;R)$. That is, compactly-supported cohomology is a {\it covariant\/} functor under inclusions of open sets~$i:U\hookra Y$.  
\item[(d)] (Homotopy.) Unlike (co)homology, compactly-supported cohomology is not homotopy invariant (although it is invariant under proper homotopies). For example, $H^n_\cs(\R^n;R)=R$ and $H^k_\cs(\R^n;R)=0$ for $k\ne n$, though the homotopy type of $\R^n$ is independent of~$n$.
\item[(e)] (Additivity Axiom.) Suppose $Y$ is a manifold with $Y=\coprod_{a\in A}Y_a$ for $A$ a countable indexing set and each $Y_a$ open and closed in $Y$. Then for all $k\ge 0$ we have canonical isomorphisms
\begin{equation*}
H^k_\cs(Y;R)\cong \ts\bigop_{a\in A}H^k_\cs(Y_a;R),
\end{equation*}
compatible with the morphisms $(i_a)_*:H^k_\cs(Y_a;R)\ra H^k_\cs(Y;R)$ for $a\in A$ induced by the inclusion $i_a:Y_a\hookra Y$.
\item[(f)] (Dimension Axiom.) When $Y=*$ is the point we have $H^0_\cs(*;R)\cong R$ and $H^k_\cs(*;R)=0$ for all~$k\ne 0$.
\end{itemize}
Any two such compactly-supported cohomology theories $H^*_\cs(-;R),\ti H_\cs^*(-;R)$ over the same ring $R$ are canonically isomorphic, as in Theorem \ref{kh2thm2}. We can also take $H^k_\cs(Y;R)$ to be defined for $k\in\Z$ rather than $k\in\N$, but then $H^k_\cs(Y;R)=0$ for~$k<0$.
\label{kh2pr1}
\end{property}

\begin{rem}{\bf(a)} We can combine the two kinds of functoriality in Property \ref{kh2pr1}(b),(c) as follows. Write $\Manpr$ for the category with objects smooth manifolds $Y$, and with morphisms $\uf:Y_1\ra Y_2$ pairs $\uf=(Y_1',f)$ with $Y_1'\subseteq Y_1$ an open set, and $f:Y_1'\ra Y_2$ a smooth proper map. Composition of morphisms $\uf=(Y_1',f):Y_1\ra Y_2$ and $\ug=(Y_2',g):Y_2\ra Y_3$ is $\ug\ci\uf:=\bigl(f^{-1}(Y_2'),g\ci f\vert_{f^{-1}(Y_2')}\bigr):Y_1\ra Y_3$. Identities are $\uid_Y=(Y,\id_Y):Y\ra Y$.

For each $\uf:Y_1\ra Y_2$ in $\Manpr$, define $\uf^*:H^k_\cs(Y_2;R)\ra H^k_\cs(Y_1;R)$ by $\uf^*=i_*\ci f^*$, where $f^*:H^k_\cs(Y_2;R)\ra H^k_\cs(Y_1';R)$ and $i_*:H^k_\cs(Y_1';R)\ra H^k_\cs(Y_1;R)$ are as in Property \ref{kh2pr1}(b),(c). Then mapping $Y\mapsto H^k_\cs(Y;R)$ and $\uf\mapsto\uf^*$ defines a functor $(\Manpr)^{\bf op}\ra\Rmod$, similar to Axiom~\ref{kh2ax2}(i).
\smallskip

\noindent{\bf(b)} We interpret the category $\Manpr$ in terms of the one-point compactifications $Y\amalg\{\iy\}$ of manifolds $Y$, as used in \eq{kh2eq15}. Given a morphism $\uf=(Y_1',f):Y_1\ra Y_2$ in $\Manpr$, define a map $\ti f:Y_1\amalg\{\iy\}\ra Y_2\amalg\{\iy\}$ between the one-point compactifications of $Y_1,Y_2$ by $\ti f(y)=f(y)$ if $y\in Y_1'\subseteq Y_1$, and $\ti f(y)=\iy$ if $y\in Y_1\sm Y_1'$, and $\ti f(\iy)=\iy$. Then $Y_1'$ open and $f$ proper imply that $\ti f$ is continuous. This establishes a functorial 1-1 correspondence between morphisms $\uf:Y_1\ra Y_2$ in $\Manpr$, and continuous maps $\ti f:Y_1\amalg\{\iy\}\ra Y_2\amalg\{\iy\}$ with $\ti f(\iy)=\iy$ which are smooth on~$\ti f^{-1}(Y_2)\subseteq Y_1$. 
\label{kh2rem2}
\end{rem}

Following Hatcher \cite[p.~233]{Hatc} and Spanier \cite[p.~323]{Span} we define compactly-supported singular cohomology:

\begin{ex} Let $Y$ be a smooth manifold, $R$ a commutative ring and $k\ge 0$. Example \ref{kh2ex1} defined $C_k^\rsi(Y;\Z)$, and Example \ref{kh2ex5} defined $C^k_\rsi(Y;R)=\Hom_\Z\bigl(C_k^\rsi(Y;\Z),R\bigr)$. Define the {\it compactly-supported singular cochains\/} to be
\begin{align*}
C^k_{\cs,\rsi}(Y;R)=\bigl\{\al\in C^k_\rsi(Y;R):\text{there exists a compact $K\subseteq Y$}&\\
\text{such that if $\ga\in C_k^\rsi(Y\sm K;\Z)\subseteq C_k^\rsi(Y;\Z)$ then $\al(\ga)=0$}&\bigr\}.
\end{align*}
It is an $R$-submodule of $C^k_\rsi(Y;R)$. Define $\d:C^k_{\cs,\rsi}(Y;R)\ra C^{k+1}_{\cs,\rsi}(Y;R)$ to be the restriction of $\d:C^k_\rsi(Y;R)\ra C^{k+1}_\rsi(Y;R)$, so that $\bigl(C^*_{\cs,\rsi}(Y;R),\d\bigr)$ is a subcomplex of $\bigl(C^*_\rsi(Y;R),\d\bigr)$. Define the {\it compactly-supported singular cohomology\/} $H^k_{\cs,\rsi}(Y;R)$ to be the $k^{\rm th}$ cohomology group of~$\bigl(C^*_{\cs,\rsi}(Y;R),\d\bigr)$.
\label{kh2ex8}
\end{ex}

\begin{ex} Define {\it compactly-supported smooth singular cochains\/} and {\it cohomology\/} $\bigl(C^*_{\cs,\ssi}(Y;R),\d\bigr),H^*_{\cs,\ssi}(Y;R)$ as in Example \ref{kh2ex8}, but using the smooth versions $C_k^\ssi(Y;\Z),C^k_\ssi(Y;R)$ from Examples \ref{kh2ex2} and~\ref{kh2ex6}.
\label{kh2ex9}
\end{ex}

We also define compactly-supported de Rham cohomology, as in~\cite[\S 1.1]{BoTu}. 

\begin{ex} Work over $R=\R$. Let $Y$ be a smooth manifold, and $k\ge 0$. Call a smooth $k$-form $\al$ on $Y$ {\it compactly-supported\/} if the (closed) support $\supp\al$ of $\al$ is a compact subset of $Y$. Equivalently, $\al$ is compactly-supported if there exists an open set $Z\subseteq Y$ with $\al\vert_Z=0$ and $Y\sm Z$ compact.

Define the {\it compactly-supported de Rham cochains\/} $C^k_{\cs,\dR}(Y;\R)$ to be the vector subspace of compactly-supported $k$-forms $\al$ in $C^\iy(\La^kT^*Y)$. Define $\d:C^k_{\cs,\dR}(Y;\R)\ra C^{k+1}_{\cs,\dR}(Y;\R)$ to be the usual exterior derivative on $k$-forms. Then $\bigl(C^*_{\cs,\dR}(Y;\R),\d\bigr)$ is a cochain complex over $\R$, a subcomplex of $\bigl(C^*_\dR(Y;\R),\d\bigr)$ in Example \ref{kh2ex7}. Define the {\it compactly-supported de Rham cohomology group\/} $H^k_{\cs,\dR}(Y;\R)$ to be the $k^{\rm th}$ cohomology group of~$\bigl(C^*_{\cs,\dR}(Y;\R),\d\bigr)$.

The morphism $\Pi:\bigl(C^*_\dR(Y;\R),\d\bigr)\ra\bigl(C^*_\ssi(Y;\R),\d\bigr)$ in \eq{kh2eq14} maps $\Pi:\bigl(C^*_{\cs,\dR}(Y;\R),\d\bigr)\ra\bigl(C^*_{\cs,\ssi}(Y;\R),\d\bigr)$, and induces the natural isomorphisms $\Pi_*:H^*_{\cs,\dR}(Y;\R)\ra H^*_{\cs,\ssi}(Y;\R)$.
\label{kh2ex10}
\end{ex}

\subsection{Locally finite homology}
\label{kh24}

We can also ask whether there is a variation on homology in \S\ref{kh21}, which is analogous to compactly-supported cohomology in \S\ref{kh23} as a variation on cohomology in \S\ref{kh22}. Ordinary homology is already compactly-supported, though ordinary cohomology is not. There is such a variation, called {\it locally finite homology\/} \cite{Geog,HuRa}, or {\it homology with closed supports\/} \cite{Dimc}, or {\it homology of the second kind\/} \cite{Petk,Sklj1}, or {\it Borel--Moore homology\/} \cite{BoMo,Iver}. We recommend Hughes and Ranicki \cite[\S 3]{HuRa} for an introduction. We will write locally finite homology as~$H_*^\lf(Y;R)$.

Some examples of locally finite homology theories are locally finite singular homology \cite[\S 3]{HuRa}, Borel--Moore homology \cite{BoMo}, Massey's version of Alexander--Spanier homology \cite[\S I]{Mass2}, and hypercohomology $\H^{-*}(Y,\om_Y)$ of the dualizing complex $\om_Y$, as in \S\ref{kh25}. Axioms for locally finite homology are given by Eilenberg and Steenrod \cite[\S X.7]{EiSt}, Petkova \cite{Petk}, and Skljarenko~\cite{Sklj1}.

As for Property \ref{kh2pr1}, here are some useful properties of locally finite homology of manifolds. This is not a definition or an axiomatization.

\begin{property} Let $R$ be a commutative ring, and fix a homology theory of manifolds $H_*(-;R)$ as in \S\ref{kh22}. A {\it locally finite homology theory of manifolds\/}  $H_*^\lf(-;R)$ gives an $R$-module $H_k^\lf(Y;R)$ for each manifold $Y$ and $k\in\N$. Furthermore:
\begin{itemize}
\setlength{\itemsep}{0pt}
\setlength{\parsep}{0pt}
\item[(a)] There are natural $R$-module morphisms $\Pi:H_k(Y;R)\ra H_k^\lf(Y;R)$ for all $Y,k$. If $Y$ is compact then $\Pi:H_k(Y;R)\ra H_k^\lf(Y;R)$ is an isomorphism.
\item[(b)] Let $f:Y_1\ra Y_2$ be a proper, smooth map of manifolds, and $k\in\N$. Then there is a natural $R$-module morphism $f_*:H_k^\lf(Y_1;R)\ra H_k^\lf(Y_2;R)$. The following also commutes:
\begin{equation*}
\xymatrix@C=130pt@R=13pt{
*+[r]{H_k(Y_1;R)} \ar[d]^{f_*} \ar[r]_\Pi & *+[l]{H_k^\lf(Y_1;R)} \ar[d]_{f_*} \\
*+[r]{H_k(Y_2;R)} \ar[r]^\Pi & *+[l]{H_k^\lf(Y_2;R).\!\!} }
\end{equation*}
If $f:Y_1\ra Y_2$, $g:Y_2\ra Y_3$ are proper, smooth maps then $(g\ci f)_*=g_*\ci f_*:H_k^\lf(Y_1;R)\ra H_k^\lf(Y_3;R)$. That is, locally finite homology is a {\it covariant\/} functor under proper, smooth maps $f:Y_1\ra Y_2$.
\item[(c)] Let $Y$ be a manifold, $U\subseteq Y$ an open set, and $k\in\N$. Write $i:U\hookra Y$ for the inclusion. Then there is a natural $R$-module morphism $i^*:H_k^\lf(Y;R)\ra H_k^\lf(U;R)$. If $j:U'\hookra U$ is another inclusion of open sets then $(i\ci j)^*=\ab j^*\ci i^*:H_k^\lf(Y;R)\ra H_k^\lf(U';R)$. That is, locally finite homology is a {\it contravariant\/} functor under open inclusions~$i:U\hookra Y$.

As in Remark \ref{kh2rem2}(a), we can combine these two forms of functoriality in (b),(c) to define a functor $\Manpr\ra\Rmod$ mapping~$Y\mapsto H_k^\lf(Y;R)$.  
\item[(d)] (Homotopy.) Locally finite homology is not homotopy invariant (although it is invariant under proper homotopies). For example, $H_n^\lf(\R^n;R)=R$ and $H_k^\lf(\R^n;R)=0$ for $k\ne n$, though the homotopy type of $\R^n$ is independent of~$n$.
\item[(e)] (Additivity Axiom.) Suppose $Y$ is a manifold with $Y=\coprod_{a\in A}Y_a$ for $A$ a countable indexing set and each $Y_a$ open and closed in $Y$. Then for all $k\ge 0$ we have canonical isomorphisms
\begin{equation*}
H_k^\lf(Y;R)\cong \ts\prod_{a\in A}H_k^\lf(Y_a;R),
\end{equation*}
compatible with the morphisms $(i_a)^*:H_k^\lf(Y;R)\ra H_k^\lf(Y_a;R)$ for $a\in A$ induced by the inclusion $i_a:Y_a\hookra Y$.
\item[(f)] (Dimension Axiom.) When $Y=*$ is the point we have $H_0^\lf(*;R)\cong R$ and $H_k^\lf(*;R)=0$ for all~$k\ne 0$.
\item[(g)] (Field coefficients.) When $R$ is a field $\K$, we have canonical isomorphisms
\begin{equation*}
H_k^\lf(Y;\K)\cong H^k_\cs(Y;\K)^*
\end{equation*}
for all $Y,k$, as for the duality of (co)homology over $\K$ in Example \ref{kh2ex4}, though note that homology and cohomology are exchanged. Under these isomorphisms, $\Pi,f_*,i^*$ above are dual to $\Pi,f^*,i_*$ in Property~\ref{kh2pr1}.

Also in this case we may write $H_k^\lf(Y;\K)$ as an inverse limit of relative homology groups, in a similar way to \eq{kh2eq16}, by 
\e
H_k^\lf(Y;\K) \cong \mathop{\underleftarrow{\lim}}\limits_{Z\subseteq Y:\text{ $Y\sm Z$ compact}\!\!\!\!\!\!\!\!\!\!\!\!\!\!\!\!\!\!\!\!\!\!\!\!\!\!\!\!\!\!\!\!\!\!\!\!\!\!\!\!}H_k(Y,Z;\K).
\label{kh2eq18}
\e
However, the analogue of \eq{kh2eq18} for general commutative rings $R$ is false.
\item[(h)] As in \cite[\S 7]{HuRa}, if $Y$ is {\it forward tame} (a condition of being `nice at infinity', which holds if $Y$ is the interior of a compact manifold with boundary $\bar Y$), then in a similar way to \eq{kh2eq15} we have 
\e
H_k^\lf(Y;R)\cong H_k(Y\amalg\{\iy\},\{\iy\};R),
\label{kh2eq19}
\e
where $Y\amalg\{\iy\}$ is the one-point compactification of $Y$. However, \eq{kh2eq19} fails for general $Y$, such as the 0-manifold $Y=\N$, \cite[Ex.~3.18]{HuRa}.
\item[(i)] (Fundamental classes.) Let $Y$ be an oriented manifold of dimension $m$, not necessarily compact. Then there is a natural {\it fundamental class\/} $[[Y]]$ in $H_m^\lf(Y;R)$. If $R$ has characteristic 2 (e.g.\ $R=\Z_2$) then $[[Y]]$ exists even if $Y$ is not oriented. If $Y$ is compact then $H_m^\lf(Y;R)=H_m(Y;R)$, so $[[Y]]\in H_m(Y;R)$ as in Remark~\ref{kh2rem1}(f).
\end{itemize}
Any two such locally finite homology theories of manifolds $H_*^\lf(-;R),\ab\ti H_*^\lf(-;R)$ over the same ring $R$ are canonically isomorphic, as in Theorem \ref{kh2thm1}. We can also define $H_k^\lf(Y;R)$ for $k\!\in\!\Z$ rather than $k\!\in\!\N$, but then $H_k^\lf(Y;R)\!=\!0$ for~$k<0$.
\label{kh2pr2}
\end{property}

The next two examples define {\it locally finite singular homology\/} $H_*^{\lf,\rsi}(Y;R)$, following Hughes and Ranicki \cite[Def.~3.1]{HuRa} and generalizing $H_*^\rsi(Y;R)$ in Example \ref{kh2ex1}, and {\it locally finite smooth singular homology\/} $H_*^{\lf,\ssi}(Y;R)$, generalizing $H_*^\ssi(Y;R)$ in Example~\ref{kh2ex2}.

\begin{ex} Let $Y$ be a topological space, $R$ a commutative ring, and $k\ge 0$. Consider formal sums $\sum_{i\in I}\rho_i\,\si_i$, where $I$ is an indexing set, $\rho_i\in R$, and $\si_i:\De_k\ra Y$ is a continuous map for $i\in I$. We call $\sum_{i\in I}\rho_i\,\si_i$ {\it locally finite\/} if any $y\in Y$ has an open neighbourhood $U\subseteq Y$ such that $U\cap\si_i(\De_k)\ne\es$ for only finitely many $i\in I$. Define the {\it locally finite singular chains\/} $C_k^{\lf,\rsi}(Y;R)$ to be the $R$-module of locally finite sums $\sum_{i\in I}\rho_i\,\si_i$, with the obvious notions of equality of sums, addition, and scalar multiplication. 

Note that $C_k^\rsi(Y;R)$ in Example \ref{kh2ex1} is the $R$-submodule of $\sum_{i\in I}\rho_i\,\si_i$ in $C_k^{\lf,\rsi}(Y;R)$ with $I$ finite. Write $\Pi:C_k^\rsi(Y;R)\hookra C_k^{\lf,\rsi}(Y;R)$ for the inclusion. As in \cite[Prop.~3.16(i)]{HuRa}, if $Y$ is locally compact we may write $C_k^{\lf,\rsi}(Y;R)$ as the inverse limit
\begin{equation*}
C_k^{\lf,\rsi}(Y;R)=\mathop{\underleftarrow{\lim}\,}\nolimits_{\text{$Z:Z\subseteq Y$, $Y\sm Z$ is compact}}C_k^\rsi(Y,Z;R).
\end{equation*}

Define $\pd:C_k^{\lf,\rsi}(Y;R)\ra C_{k-1}^{\lf,\rsi}(Y;R)$ as in \eq{kh2eq6}. Then $\pd^2=0$, and {\it locally finite singular homology\/} $H_*^{\lf,\rsi}(Y;R)$ is the homology of $\bigl(C_*^{\lf,\rsi}(Y;R),\pd\bigr)$, that~is, 
\begin{equation*}
H_k^{\lf,\rsi}(Y;R)=\frac{\ts \Ker\bigl(\pd: C_k^{\lf,\rsi}(Y;R)\longra C_{k-1}^{\lf,\rsi}(Y;R)\bigr)}{\ts \Im\bigl(\pd:C_{k+1}^{\lf,\rsi}(Y;R)\longra C_k^{\lf,\rsi}(Y;R)\bigr)}\,.
\end{equation*}
Then $\Pi\ci\pd=\pd\ci\Pi:C_k^\rsi(Y;R)\ra C_{k-1}^{\lf,\rsi}(Y;R)$, so these $\Pi$ induce morphisms $\Pi:H_k^\rsi(Y;R)\hookra H_k^{\lf,\rsi}(Y;R)$, as in Property~\ref{kh2pr2}(a).

Let $f:Y_1\ra Y_2$ be a proper map of topological spaces. Define the {\it pushforward\/} $f_*:C_k^{\lf,\rsi}(Y_1;R)\ra C_k^{\lf,\rsi}(Y_2;R)$ by \eq{kh2eq8}, where $f$ proper implies that the r.h.s.\ of \eq{kh2eq8} is a locally finite sum. Then $f_*\ci\pd=\pd\ci f_*$, so the $f_*$ induce morphisms $f_*:H_k^{\lf,\rsi}(Y_1;R)\hookra H_k^{\lf,\rsi}(Y_2;R)$, as in Property \ref{kh2pr2}(b). Pushforwards are covariantly functorial on both $C_*^{\lf,\rsi}(Y_i;R)$ and $H_*^{\lf,\rsi}(Y_i;R)$.

If $i:U\!\hookra\! Y$ is an open inclusion, the pullback $i^*\!:\!H_k^{\lf,\rsi}(Y;R)\!\ra\! H_k^{\lf,\rsi}(U;R)$ in Property \ref{kh2pr2}(c) does not have a nice expression at the chain level.
\label{kh2ex11}
\end{ex}

\begin{ex} Let $Y$ be a manifold and $R$ a commutative ring. We follow Example \ref{kh2ex11}, but using {\it smooth\/} singular chains $\si:\De_k\ra Y$ as in Example \ref{kh2ex2}. So we define the {\it locally finite smooth singular chains\/} $C_k^{\lf,\ssi}(Y;R)$ to be the $R$-module of locally finite sums $\sum_{i\in I}\rho_i\,\si_i$ with $\si_i:\De_k\ra Y$ smooth. Define $\pd:C_k^{\lf,\ssi}(Y;R)\ra C_{k-1}^{\lf,\ssi}(Y;R)$ as in \eq{kh2eq6}, and define {\it locally finite smooth singular homology\/} $H_*^{\lf,\ssi}(Y;R)$ to be the homology of $\bigl(C_*^{\lf,\ssi}(Y;R),\pd\bigr)$.

The $R$-submodule of $\sum_{i\in I}\rho_i\,\si_i$ in $C_k^{\lf,\ssi}(Y;R)$ with $I$ finite is $C_k^\ssi(Y;R)$ in Example \ref{kh2ex2}. Write $\Pi:C_k^\ssi(Y;R)\hookra C_k^{\lf,\ssi}(Y;R)$ for the inclusions. They induce morphisms $\Pi:H_k^\ssi(Y;R)\hookra H_k^{\lf,\ssi}(Y;R)$. We also have
\begin{equation*}
C_k^{\lf,\ssi}(Y;R)=\mathop{\underleftarrow{\lim}\,}\nolimits_{\text{$Z:Z\subseteq Y$, $Y\sm Z$ is compact}}C_k^\ssi(Y,Z;R).
\end{equation*}

Let $f:Y_1\ra Y_2$ be a proper smooth map of manifolds. Define the {\it pushforward\/} $f_*:C_k^{\lf,\ssi}(Y_1;R)\ra C_k^{\lf,\ssi}(Y_2;R)$ by \eq{kh2eq8}. The $f_*$ induce morphisms $f_*:H_k^{\lf,\ssi}(Y_1;R)\hookra H_k^{\lf,\ssi}(Y_2;R)$, as in Property \ref{kh2pr2}(b). Pushforwards are covariantly functorial on both $C_*^{\lf,\ssi}(Y_i;R)$ and~$H_*^{\lf,\ssi}(Y_i;R)$.

Again, if $i:U\hookra Y$ is an inclusion of an open set, the pullback $i^*:H_k^{\lf,\ssi}(Y;R)\ra H_k^{\lf,\ssi}(U;R)$ in Property \ref{kh2pr2}(c) does not have a nice expression at the chain level. But we will explain pullbacks $i^*$ at the chain level for locally finite sheaf smooth singular homology $\hat H_k^{\lf,\ssi}(Y;R)$ in Example \ref{kh2ex17} below.

The inclusion $\bigl(C_*^{\lf,\ssi}(Y;R),\pd\bigr)\hookra \bigl(C_*^{\lf,\rsi}(Y;R),\pd\bigr)$ induces morphisms on homology $H_*^{\lf,\ssi}(Y;R)\ra H_*^{\lf,\rsi}(Y;R)$. The argument of Bredon \cite[\S V.5 \& \S V.9]{Bred1} shows that these are isomorphisms.

Suppose $Y$ is  oriented of dimension $m$, so that as in Property \ref{kh2pr2}(i) it has a fundamental class $[[Y]]\in H_m^{\lf,\ssi}(Y;R)$. As in Example \ref{kh2ex2} for $Y$ compact, we can construct this at the chain level as follows. 

Choose a locally finite triangulation of $Y$ which cuts it into smooth $m$-simplices, which we write as $\si_i:\De_m\ra Y$ for $i$ in $I$, an indexing set. Here each $\si_i$ embeds $\De_m$ into $Y$ as a submanifold with corners. Write $\ep_i=1$ if $\si_i$ is orientation-preserving, and $\ep_i=-1$ otherwise. Then $\sum_{i\in I}\ep_i\,\si_i$ is a locally finite sum and lies in $C_m^\lf(Y;R)$, and $\pd\sum_{i\in I}\ep_i\,\si_i=0$, since each boundary face of each $m$-simplex $\si_i(\De_m)$ is cancelled by a boundary face of a neighbouring $m$-simplex in the triangulation. The fundamental class is~$[[Y]]=\bigl[\sum_{i\in I}\ep_i\,\si_i\bigr]$.
\label{kh2ex12}
\end{ex}

\subsection{(Co)homology via sheaf cohomology}
\label{kh25}

One can also define and study both the cohomology and the homology of topological spaces using {\it sheaf cohomology}. Sheaves are most commonly used in algebraic geometry, as in Hartshorne \cite[\S II--\S III]{Hart} for instance, but they are also important in topology and algebraic topology. Some books on sheaves and sheaf cohomology from the point of view we need are Bredon \cite{Bred2}, Dimca \cite{Dimc}, Gelfand and Manin \cite{GeMa}, Godement \cite{Gode}, Iversen \cite{Iver}, Kashiwara and Schapira \cite{KaSc}, and Strooker \cite{Stro}. We begin with the basic definitions of sheaf theory:

\begin{dfn} Let $Y$ be a topological space, and $R$ a commutative ring. A {\it presheaf of
$R$-modules\/} $\cE$ on $Y$ consists of the data of an $R$-module $\cE(U)$ for every open set $U\subseteq Y$, and a morphism of $R$-modules $\rho_{UV}:\cE(U)\ra\cE(V)$ called the {\it restriction map\/} for every inclusion $V\subseteq U\subseteq Y$ of open sets, satisfying
\begin{itemize}
\setlength{\itemsep}{0pt}
\setlength{\parsep}{0pt}
\item[(i)] $\cE(\es)=0$;
\item[(ii)] $\rho_{UU}=\id_{\cE(U)}:\cE(U)\ra\cE(U)$ for all open
$U\subseteq Y$; and
\item[(iii)] $\rho_{UW}=\rho_{VW}\ci\rho_{UV}:\cE(U)\ra\cE(W)$ for all
open~$W\subseteq V\subseteq U\subseteq Y$.
\end{itemize}
We often write $s\vert_V$ rather than $\rho_{UV}(s)$, for $s\in\cE(U)$.

A presheaf of $R$-modules $\cE$ on $Y$ is called a {\it sheaf\/} if it also satisfies
\begin{itemize}
\setlength{\itemsep}{0pt}
\setlength{\parsep}{0pt}
\item[(iv)] If $U\subseteq Y$ is open, $\{V_i:i\in I\}$ is an open
cover of $U$, and $s\in\cE(U)$ has $\rho_{UV_i}(s)=0$ in
$\cE(V_i)$ for all $i\in I$, then $s=0$ in $\cE(U)$; and
\item[(v)] If $U\subseteq Y$ is open, $\{V_i:i\in I\}$ is an open cover of
$U$, and we are given elements $s_i\in\cE(V_i)$ for all $i\in I$
such that $\rho_{V_i(V_i\cap V_j)}(s_i)=\rho_{V_j(V_i\cap
V_j)}(s_j)$ in $\cE(V_i\cap V_j)$ for all $i,j\in I$, then there
exists $s\in\cE(U)$ with $\rho_{UV_i}(s)=s_i$ for all $i\in I$.
This $s$ is unique by~(iv).
\end{itemize}

Suppose $\cE,\cF$ are presheaves or sheaves of $R$-modules on
$Y$. A {\it morphism\/} $\phi:\cE\ra\cF$ consists of a morphism of
$R$-modules $\phi(U):\cE(U)\ra\cF(U)$ for all open $U\subseteq
Y$, such that the following diagram commutes for all open
$V\subseteq U\subseteq Y$
\begin{equation*}
\xymatrix@C=75pt@R=13pt{
\cE(U) \ar[r]_{\phi(U)} \ar[d]_{\rho_{UV}} & \cF(U)
\ar[d]_{\rho_{UV}'} \\ \cE(V) \ar[r]^{\phi(V)} & \cF(V),\!{} }
\end{equation*}
where $\rho_{UV}$ is the restriction map for $\cE$, and $\rho_{UV}'$
the restriction map for~$\cF$.
\label{kh2def1}
\end{dfn}

\begin{dfn} Let $\cE$ be a presheaf of $R$-modules on $Y$. For each $y\in Y$, the {\it stalk\/} $\cE_y$ is the direct limit of the $R$-modules $\cE(U)$ for all $y\in U\subseteq Y$, via the restriction maps $\rho_{UV}$. It is an $R$-module. A morphism $\phi:\cE\ra\cF$ induces morphisms $\phi_y:\cE_y\ra\cF_y$ for all $y\in Y$. If $\cE,\cF$ are sheaves then $\phi$ is an isomorphism if and only if $\phi_y$ is an isomorphism for all~$y\in Y$.
\label{kh2def2}
\end{dfn}

Sheaves of $R$-modules on $Y$ form an abelian category $\Sh(Y;R)$. Thus we have (category-theoretic) notions of when a morphism $\phi:\cE\ra\cF$ in $\Sh(Y;R)$ is injective or surjective, and when a sequence $\cE\ra\cF\ra\cG$ in $\Sh(Y;R)$ is exact.

\begin{dfn} Let $\cE$ be a presheaf of $R$-modules on $Y$. A {\it sheafification\/} of $\cE$ is a sheaf of $R$-modules $\hat\cE$ on $Y$ and a morphism $\pi:\cE\ra\hat\cE$, such that whenever $\cF$ is a sheaf of $R$-modules on $Y$ and $\phi:\cE\ra\cF$ is a morphism, there is a unique morphism $\hat\phi:\hat\cE\ra\cF$ with $\phi=\hat\phi\ci\pi$. Sheafifications always exist, and are unique up to canonical isomorphism. If $\pi:\cE\ra\hat\cE$ is a sheafification then the induced morphisms on stalks $\pi_y:\cE_y\ra\hat\cE_y$ are isomorphisms for all~$y\in Y$. 
\label{kh2def3}
\end{dfn}

\begin{dfn} Let $f:Y\ra Z$ be a continuous map of topological spaces, and $\cE$ a sheaf of $R$-modules on $Y$. Define the {\it pushforward\/} ({\it direct image\/}) sheaf
$f_*(\cE)$ on $Z$ by $\bigl(f_*(\cE)\bigr)(U)=\cE\bigl(f^{-1}(U)\bigr)$ for all
open $U\subseteq V$, with restriction maps $\rho'_{UV}=\rho_{f^{-1}(U)f^{-1}(V)}:\bigl(f_*(\cE)\bigr)(U)\ra\bigl(f_*(\cE)\bigr)(V)$ for all open $V\subseteq U\subseteq Z$.
Then $f_*(\cE)$ is a sheaf of $R$-modules on~$Z$.

If $\phi:\cE\ra\cF$ is a morphism in $\Sh(Y;R)$ we define $f_*(\phi):f_*(\cE)\ra f_*(\cF)$ by $\bigl(f_*(\phi)\bigr)(u)=\phi\bigl(f^{-1}(U)\bigr)$ for all open $U\subseteq Z$. Then $f_*(\phi)$ is a morphism in $\Sh(Z;R)$, and $f_*$ is a functor $\Sh(Y;R)\ra\Sh(Z;R)$. It is a left exact functor between abelian categories, but in general is not exact. For continuous maps $e:X\ra Y$, $f:Y\ra Z$ we have~$(f\ci e)_*=f_*\ci e_*$.
\label{kh2def4}
\end{dfn}

\begin{dfn} Let $f:Y\ra Z$ be a continuous map of topological spaces, and $\cE$ a sheaf of $R$-modules on $Z$. Define a presheaf of $R$-modules ${\cal P}f^{-1}(\cE)$ on $Y$ by $\bigl({\cal P}f^{-1}(\cE)\bigr)(U)=\underrightarrow{\lim}\,_{A\supseteq f(U)}\cE(A)$, where the direct limit is taken over all open $A\subseteq Z$ containing $f(U)$, using the restriction maps $\rho_{AB}$ in $\cE$. For open $V\subseteq U\subseteq Y$, define $\rho_{UV}': \bigl({\cal P}f^{-1}(\cE)\bigr)(U)\ra \bigl({\cal P}f^{-1}(\cE)\bigr)(V)$ as the direct limit of the morphisms $\rho_{AB}$ in $\cE$ for open $B\subseteq A\subseteq Z$ with $f(U)\subseteq A$ and $f(V)\subseteq B$. Then we define the {\it pullback\/} ({\it inverse image\/}) $f^{-1}(\cE)$ to be the sheafification of the presheaf~${\cal P}f^{-1}(\cE)$.

If $\phi:\cE\ra\cF$ is a morphism in $\Sh(Z;R)$, one can define a pullback morphism
$f^{-1}(\phi):f^{-1}(\cE)\ra f^{-1}(\cF)$. Then $f^{-1}:\Sh(Z;R)\ra\Sh(Y;R)$ is an exact functor between abelian categories, which is left adjoint to $f_*:\Sh(Y;R)\ra\Sh(Z;R)$. That is, there are natural bijections
\e
\Hom_Y\bigl(f^{-1}(\cE),\cF\bigr)=\Hom_Z\bigl(\cE,f_*(\cF)\bigr)
\label{kh2eq20}
\e
for all $\cE\in\Sh(Z;R)$ and $\cF\in\Sh(Y;R)$, with functorial properties. For continuous maps $e:X\ra Y$, $f:Y\ra Z$ we have~$(f\ci e)^{-1}=e^{-1}\ci f^{-1}$.
\label{kh2def5}
\end{dfn}

\begin{dfn} Let $\cE$ be a presheaf of $R$-modules on $Y$, and suppose $U\subseteq Y$ is open and $s\in\cE(U)$. We define the {\it support\/} $\supp s$ to be
\begin{equation*}
\supp s=\bigl\{y\in U:\text{$s_y\ne 0$ in the stalk $\cE_y$}\bigr\},
\end{equation*}
where $s_y$ is the germ of $s$ at $y$. Then $\supp s$ is a closed subset of $U$.

We call the section $s$ {\it compactly-supported\/} if there exists a compact subset $K\subseteq U$ with $\rho_{U(U\sm K)}(s)=0$. Then $\supp s\subseteq K$, so $\supp S$ is compact. We write $\cE_\cs(U)$ for the $R$-submodule of compactly-supported $s$ in~$\cE(U)$.
\label{kh2def6}
\end{dfn}

\begin{rem} In Definition \ref{kh2def6}, if $\cE$ is a sheaf then $\rho_{U(U\sm\supp s)}(s)=0$, so $s$ is compactly-supported if and only if $\supp s$ is compact. However, if $\cE$ is only a presheaf then we can have $\rho_{U(U\sm\supp s)}(s)\ne 0$, and $\supp s$ can be compact without $s$ being compactly-supported, which may cause confusion.
\label{kh2rem3}
\end{rem}

\begin{ex} Let $Y$ be a topological space, and $R$ a commutative ring. For open $U\subseteq Y$, write $R_Y(U)$ for the $R$-module of locally constant functions $s:U\ra R$. For open $V\!\subseteq\! U\!\subseteq\! Y$, define $\rho_{UV}:R_Y(U)\ra R_Y(V)$ by $\rho_{UV}(s)=s\vert_V$. Then $R_Y$ is a sheaf of $R$-modules on $Y$, called the {\it constant sheaf}. The stalk $R_{Y,y}$ of $R_Y$ at each $y\in Y$ is~$R_{Y,y}\cong R$.

Let $f:Y\ra Z$ be a continuous map of topological spaces. Definition \ref{kh2def6} gives a presheaf ${\cal P}f^{-1}(R_Z)$ on $Y$ by $\bigl({\cal P}f^{-1}(R_Z)\bigr)(U)=\lim_{A\supseteq f(U)}R_Z(A)$. We have natural morphisms $R_Z(A)\ra R_Y(U)$ mapping $s\mapsto s\ci f$ for open $f(U)\subseteq A\subseteq Z$, and these induce morphisms $\bigl({\cal P}f^{-1}(R_Z)\bigr)(U)\ra R_Y(U)$, giving a morphism of presheaves ${\cal P}f^{-1}(R_Z)\ra R_Y$, which factors through a unique morphism $f^\sh:f^{-1}(R_Z)\ra R_Y$ from the sheafification $f^{-1}(R_Z)$. 

On the stalks at $y\in Y$, $f^\sh_y:f^{-1}(R_Z)_y\ra R_{Y,y}$ is just $\id:R\ra R$, so $f^\sh:f^{-1}(R_Z)\ra R_Y$ is an isomorphism. Under the bijection \eq{kh2eq20}, $f^\sh$ corresponds to a natural morphism $f_\sh:R_Z\ra f_*(R_Y)$, which is generally not an isomorphism.

\label{kh2ex13}
\end{ex}

Next we discuss sheaf cohomology.

\begin{dfn} Let $Y$ be a topological space, $R$ a commutative ring, and $\cE$ a sheaf of $R$-modules on $Y$. Then the {\it sheaf cohomology groups\/} $H^k(Y,\cE)$ for $k=0,1,\ldots$ and the {\it compactly-supported sheaf cohomology groups\/} $H^k_\cs(Y,\cE)$ for $k=0,1,\ldots,$ are $R$-modules. There are several equivalent ways to define them. 

One method, following \cite{Dimc,GeMa,Hart,Iver,KaSc,Stro}, is to first define the {\it global sections functor\/} $\Ga_Y:\Sh(Y;R)\ra\Rmod$ and the {\it compactly-supported global sections functor\/} $\Ga_{\cs,Y}:\Sh(Y;R)\ra\Rmod$ by $\Ga_Y(\cE)=\cE(Y)$, $\Ga_{\cs,Y}(\cE)=\cE_\cs(Y)$ on objects $\cE$, and $\Ga_Y(\phi)=\phi(Y)$, $\Ga_{\cs,Y}(\phi)=\phi(Y)\vert_{\cE_\cs(Y)}$ on morphisms $\phi:\cE\ra\cF$. Then $\Ga_Y,\Ga_{\cs,Y}:\Sh(Y;R)\ra\Rmod$ are left exact functors of abelian categories, and so have right derived functors $R\Ga_Y^k,R\Ga_{\cs,Y}^k:\Sh(Y;R)\ra\Rmod$ for $k=0,1,\ldots$ (which we will not define). Then we write $H^k(Y,\cE)=R\Ga_Y^k(\cE)$ and~$H^k_\cs(Y,\cE)=R\Ga_{\cs,Y}^k(\cE)$.

Alternatively, \cite{Bred2}, \cite[\S II.4]{Gode} define sheaf cohomology using Godement's can\-on\-ic\-al resolution. For nice topological spaces (including manifolds), one can also define sheaf cohomology using \v Cech cohomology \cite[\S 10]{BoTu}, \cite[\S II.5]{Gode}.
\label{kh2def7}
\end{dfn}

The relation to cohomology of topological spaces, as in \S\ref{kh22}--\S\ref{kh23}, is that as in Bredon \cite[\S III]{Bred2}, for sufficiently nice topological spaces $Y$ (including manifolds), with $R_Y$ the constant sheaf from Example \ref{kh2ex13} we have canonical isomorphisms
\e
H^k(Y;R)\cong H^k(Y,R_Y)\qquad\text{and}\qquad H^k_\cs(Y;R)\cong H^k_\cs(Y,R_Y).
\label{kh2eq21}
\e
Pullbacks $f^*:H^k(Z;R)\ra H^k(Y;R)$ and $f^*:H^k_\cs(Z;R)\ra H^k_\cs(Y;R)$ for (proper) continuous $f:Y\ra Z$ can be written using $f^\sh:f^{-1}(R_Z)\ra R_Y$ in Example \ref{kh2ex13}, by~$H^k(Z,R_Z)\ra H^k(Y,f^{-1}(R_Z))\,{\buildrel f^\sh_*\over\longra}\, H^k(Y,R_Y)$. 

Often it is useful to work not with individual sheaves $\cE$, but with (bounded below) complexes $\cE^\bu=(\cdots\,{\buildrel\d\over\longra}\,\cE^k\,{\buildrel\d\over\longra}\,\cE^{k+1}\,{\buildrel\d\over\longra}\,\cdots)$ of sheaves of $R$-modules on $Y$ in the derived category $D^+(Y;R):=D^+\Sh(Y;R)$. Then one can define {\it hypercohomology groups\/} $\H^k(Y,\cE^\bu)$ and {\it compactly-supported hypercohomology groups\/} $\H^k_\cs(Y,\cE^\bu)$ for $k\in\Z$ \cite{Dimc}, which reduce to sheaf cohomology $H^k(Y,\cE^0)$, $H^k_\cs(Y,\cE^0)$ if $\cE^i=0$ for $i\ne 0$. Quasi-isomorphisms of complexes $\cE^\bu\simeq\cF^\bu$ induce isomorphisms $\H^k(Y,\cE^\bu)\cong\H^k(Y,\cF^\bu)$ and~$\H^k_\cs(Y,\cE^\bu)\cong\H^k_\cs(Y,\cF^\bu)$.

For computing sheaf (hyper)cohomology groups, one often uses {\it acyclic resolutions\/} or {\it soft resolutions}, \cite[\S II.9]{Bred2}, \cite[\S I.5]{GeMa}, \cite[\S II.3.7]{Gode}, \cite[\S 4.12]{Stro}.

\begin{dfn} Let $Y$ be a paracompact, locally compact, Hausdorff topological space (e.g. a manifold), $R$ a commutative ring, and $\cE$ a sheaf of $R$-modules on $Y$. Then:
\begin{itemize}
\setlength{\itemsep}{0pt}
\setlength{\parsep}{0pt}
\item[(a)] We call $\cE$ {\it fine\/} if given any open cover $\{U_i:i\in I\}$ of $Y$, there exists a family of morphisms $\phi_i:\cE\ra\cE$ in $\Sh(Y;R)$ for $i\in I$ such that for each $y\in Y$, the restrictions to the stalks $\phi_{i,y}:\cE_y\ra\cE_y$ for $i\in I$ are nonzero for only finitely many $i\in I$ and only if $y\in U_i$, and~$\sum_{i\in I}\phi_{i,y}=\id:\cE_y\ra\cE_y$.
\item[(b)] If $S\subseteq Y$ is closed, we write
\e
\cE_{\rm cl}(S)=\underrightarrow{\lim}_{\,\text{$U:S\subseteq U\subseteq Y$, $U$ is open in $Y$}\,\,} \cE(U),
\label{kh2eq22}
\e
the direct limit of $\cE(U)$ over open subsets $U\subseteq Y$ containing $S$, using the $\rho_{UV}:\cE(U)\ra\cE(V)$ for $S\subseteq U\subseteq V\subseteq Y$ with $U,V$ open. There are restriction maps $\rho_{US}:\cE(U)\ra\cE_{\rm cl}(S)$ for all open $U\subseteq Y$ with $S\subseteq U$.

We call $\cE$ {\it soft\/} if $\rho_{YS}:\cE(Y)\ra\cE_{\rm cl}(S)$ is surjective for all closed $S\subseteq Y$.
\item[(c)] We call $\cE$ {\it c-soft\/} if $\rho_{YS}:\cE(Y)\ra\cE_{\rm cl}(S)$ is surjective for all compact $S\subseteq Y$, where $\cE_{\rm cl}(S)$ is defined as in \eq{kh2eq22}.
\item[(d)] We call $\cE$ {\it acyclic\/} if $H^k(Y,\cE)=0$ for all $k>0$.
\end{itemize}
Then $\cE$ fine implies $\cE$ soft, implies $\cE$ c-soft and acyclic, under our assumptions on $Y$. If also $Y$ is countable at infinity (the union of countably many compact subsets), which holds for manifolds, then c-soft is equivalent to soft \cite[Ex.~II.6]{KaSc}. Roughly speaking, $\cE$ is a fine sheaf if it has `partitions of unity'. If $\cA$ is a fine sheaf of $R$-algebras on $Y$, and $\cE$ is a sheaf of modules over $\cA$, then $\cE$ is fine.

Later we will also consider presheaves $\cE$ which are soft or c-soft, defined as in (b),(c). Again, soft implies c-soft.
\label{kh2def8}
\end{dfn}

Let $\cE^\bu=(\cdots\,{\buildrel\d\over\longra}\,\cE^j\,{\buildrel\d\over\longra}\,\cE^{j+1}\,{\buildrel\d\over\longra}\,\cdots)$ lie in $D^+(Y;R)$ with $\cE^j$ acyclic (e.g. with $\cE^j$ soft) for all $j\in\Z$. Then there is a canonical isomorphism
\e
\H^k(Y,\cE^\bu)\cong H^k\bigl(\cdots\,{\buildrel\d\over\longra}\,\cE^j(Y)\,{\buildrel\d\over\longra}\,\cE^{j+1}(Y)\,{\buildrel\d\over\longra}\,\cdots\bigr).
\label{kh2eq23}
\e
Similarly, if the $\cE^j$ are c-soft rather than acyclic (e.g. with $\cE^j$ soft), there is a canonical isomorphism
\e
\H_\cs^k(Y,\cE^\bu)\cong H^k\bigl(\cdots\,{\buildrel\d\over\longra}\,\cE_\cs^j(Y)\,{\buildrel\d\over\longra}\,\cE_\cs^{j+1}(Y)\,{\buildrel\d\over\longra}\,\cdots\bigr).
\label{kh2eq24}
\e
Now suppose $\cE$ is a sheaf of $R$-modules on $Y$, and we can find an exact sequence 
\e
\xymatrix@C=20pt{ 0 \ar[r] & \cE \ar[r]^i & \cF^0 \ar[r]^\d & \cF^1 \ar[r]^\d & \cF^2 \ar[r]^\d & \cdots }
\label{kh2eq25}
\e
in $\Sh(Y;R)$ with $\cF^j$ soft for all $j$. Then $i:\cE\ra\cF^\bu$ is a quasi-isomorphism to the {\it soft resolution\/} $\cF^\bu=(0\ra\cF^0\,{\buildrel\d\over\longra}\,\cF^1\,{\buildrel\d\over\longra}\,\cdots)$, so \eq{kh2eq23}--\eq{kh2eq24} give
\ea
H^k(Y,\cE)&\cong H^k\bigl(0\longra \cF^0(Y)\,{\buildrel\d\over\longra}\,\cF^1(Y)\,{\buildrel\d\over\longra}\,\cF^2(Y)\,{\buildrel\d\over\longra}\,\cdots\bigr),
\label{kh2eq26}\\
H_\cs^k(Y,\cE)&\cong H^k\bigl(0\longra \cF^0_\cs(Y)\,{\buildrel\d\over\longra}\,\cF^1_\cs(Y)\,{\buildrel\d\over\longra}\,\cF^2_\cs(Y)\,{\buildrel\d\over\longra}\,\cdots\bigr).
\label{kh2eq27}
\ea

We illustrate this using de Rham cohomology, as in Examples \ref{kh2ex7} and~\ref{kh2ex10}.
 
\begin{ex} Let $Y$ be a smooth manifold. Write $\Om^k_Y$ for the sheaf of smooth sections of $\La^kT^*Y$, as a sheaf of $\R$-vector spaces on $Y$, so that $\Om^k_Y(U)=C^\iy(\La^kT^*U)$ for open $U\subseteq Y$. Write $\d:\Om^k_Y\ra\Om^{k+1}_Y$ for the exterior derivative, and $i:\R_Y\ra\Om^0_Y$ for the inclusion of the locally constant functions $Y\ra\R$ into the smooth functions. Then as in \eq{kh2eq25} we have a complex in $\Sh(Y;\R)$
\begin{equation*}
\xymatrix@C=20pt{ 0 \ar[r] & \R_Y \ar[r]^i & \Om^0_Y \ar[r]^\d & \Om^1_Y \ar[r]^\d & \Om^2_Y \ar[r]^\d & \cdots, }
\end{equation*}
which is exact by the Poincar\'e Lemma. Also, because partitions of unity exist in smooth functions on $Y$, the sheaves $\Om^k_Y$ are fine, and hence soft. So equations \eq{kh2eq21} and \eq{kh2eq26}--\eq{kh2eq27} give
\begin{align*}
H^k(Y;\R)&\cong H^k\bigl(0\ra C^\iy(\La^0T^*Y)\,{\buildrel\d\over\longra}\,C^\iy(\La^1T^*Y)\,{\buildrel\d\over\longra}\,\cdots\bigr)=:H^k_\dR(Y;\R),\\
H^k_\cs(Y;\R)&\cong H^k\bigl(0\ra C^\iy_\cs(\La^0T^*Y)\,{\buildrel\d\over\longra}\,C^\iy_\cs(\La^1T^*Y)\,{\buildrel\d\over\longra}\,\cdots\bigr)=:H^k_{\cs,\dR}(Y;\R).
\end{align*}

\label{kh2ex14}
\end{ex}

Perhaps surprisingly, once can also write homology of topological spaces in terms of sheaf cohomology. As in \cite[\S 3.2]{Dimc}, \cite[\S III.8]{GeMa}, \cite[\S V--\S VI]{Iver}, \cite[\S III.3.1]{KaSc}, for nice topological spaces $Y$ and provided $R$ is a noetherian ring, one can define a {\it dualizing complex\/} $\om_Y\in D^+(Y;R)$ with
\e
H_k(Y;R)\cong \H^{-k}_\cs(Y,\om_Y)\qquad\text{and}\qquad H_k^\lf(Y;R)\cong \H^{-k}(Y,\om_Y).
\label{kh2eq28}
\e
For manifolds $Y$ we can identify the dualizing complex:

\begin{dfn} Let $Y$ be a manifold of dimension $m$. Write $\pi:P\ra Y$ for the principal $\Z_2$-bundle of orientations on $Y$, and $\si:P\ra P$ for the free $\Z_2$-action. For each open $U\subseteq Y$, define $O_Y(U)$ to be the $R$-module of locally constant functions $s:\pi^{-1}(U)\ra R$ with $s\ci\si\vert_{\pi^{-1}(U)}U=-s$, and for open $V\subseteq U\subseteq Y$ define $\rho_{UV}:O_Y(U)\ra O_Y(V)$ to map $s\mapsto s\vert_{\pi^{-1}(V)}$. Then $O_Y$ is a sheaf of $R$-modules on $Y$ called the {\it orientation sheaf}.

If $U\subseteq Y$ is open, write $U=\coprod_{i\in I}U_i$ for the decomposition of $U$ into connected components $U_i$. Then elements $\al$ of $O_Y(U)$ may equivalently be written as formal sums $\al=\sum_{i\in J}a_i\, o_{U_i}$, where $a_i\in R$ and $o_{U_i}$ is an orientation on $U_i$ for $i\in J\subseteq I$, with the convention that $a_i\, o_{U_i}=(-a_i)(-o_{U_i})$, with $-o_{U_i}$ the opposite orientation to $o_{U_i}$. If $U_i$ for $i\in I$ is not orientable then we must have $i\notin J$, or set $a_i=0$, in which case we do not need to choose~$o_{U_i}$.  

Then the dualizing complex $\om_Y$ of $Y$ satisfies $\om_Y\simeq O_Y[m]$, so that
\e
H_k(Y;R)\cong H^{m-k}_\cs(Y,O_Y)\qquad\text{and}\qquad H_k^\lf(Y;R)\cong H^{m-k}(Y,O_Y).
\label{kh2eq29}
\e
If $Y$ is oriented then $O_Y\cong R_Y$, giving Poincar\'e duality isomorphisms as in~\S\ref{kh28}
\e
\begin{split}
H_k(Y;R)&\cong H^{m-k}_\cs(Y,R_Y)\cong H^{m-k}_\cs(Y;R),\\
H_k^\lf(Y;R)&\cong H^{m-k}(Y,R_Y)\cong H^{m-k}(Y;R).
\end{split}
\label{kh2eq30}
\e
The fundamental class $[[Y]]\in H_m^\lf(Y;R)$ from Property \ref{kh2pr2}(i) is identified with $1_Y\in H^0(Y;R)$ by the second isomorphism.

The theory of dualizing complexes in \cite{Dimc,GeMa,Iver,KaSc} requires the base commutative ring $R$ to be noetherian. However, the definition of $O_Y$ and the isomorphisms \eq{kh2eq29}--\eq{kh2eq30} work for arbitrary~$R$.
\label{kh2def9}
\end{dfn}

Another way to explain compactly-supported cohomology is in terms of {\it cosheaves}, which are less well-known than sheaves. Our treatment is based on Bredon~\cite[\S V.1]{Bred2}.

\begin{dfn} Let $Y$ be a paracompact, locally compact, Hausdorff topological space, and $R$ a commutative ring. A {\it precosheaf of $R$-modules\/} $\ucE$ on $Y$ consists of the data of an $R$-module $\ucE(U)$ for every open set $U\subseteq Y$, and a morphism of $R$-modules $\si_{VU}:\ucE(V)\ra\ucE(U)$ called the {\it inclusion map\/} for every inclusion $V\subseteq U\subseteq Y$ of open sets, satisfying
\begin{itemize}
\setlength{\itemsep}{0pt}
\setlength{\parsep}{0pt}
\item[(i)] $\ucE(\es)=0$;
\item[(ii)] $\si_{UU}=\id_{\ucE(U)}:\ucE(U)\ra\ucE(U)$ for all open
$U\subseteq Y$; and
\item[(iii)] $\si_{WU}=\si_{VU}\ci\si_{WV}:\ucE(W)\ra\ucE(U)$ for all
open~$W\subseteq V\subseteq U\subseteq Y$.
\end{itemize}

A precosheaf of $R$-modules $\ucE$ on $Y$ is called a {\it cosheaf\/} if it also satisfies
\begin{itemize}
\setlength{\itemsep}{0pt}
\setlength{\parsep}{0pt}
\item[(iv)] If $U,V\subseteq Y$ are open, the following sequence is exact in $\Rmod$:
\end{itemize}
\e
\xymatrix@C=12pt{ \ucE(U\!\cap\! V) \ar[rrrr]^(0.48){\si_{(U\cap V)U}\op -\si_{(U\cap V)V}} &&&& \ucE(U)\!\op\! \ucE(V) \ar[rrrr]^(0.52){\si_{U(U\cup V)}\op \si_{V(U\cup V)}} &&&& \ucE(U\!\cup\! V) \ar[r] & 0. }\!\!{}
\label{kh2eq31}
\e
\begin{itemize}
\setlength{\itemsep}{0pt}
\setlength{\parsep}{0pt}
\item[(v)] Suppose $U_1\subseteq U_2\subseteq \cdots\subseteq Y$ are open with $U=\bigcup_{i=1}^\iy U_i$. Then we have an isomorphism with the direct limit
\begin{equation*}
\ucE(U)\cong \underrightarrow{\lim}\,_{a=1}^\iy\, \ucE(U_a),
\end{equation*}
compatible with $\si_{U_aU}:\ucE(U_a)\!\ra\!\ucE(U)$ and $\si_{U_aU_b}:\ucE(U_a)\!\ra\!\ucE(U_b)$ for~$a\!\le\! b$.
\end{itemize}

Suppose $\ucE,\ucF$ are precosheaves or cosheaves of $R$-modules on
$Y$. A {\it morphism\/} $\uphi:\ucE\ra\ucF$ consists of a morphism of
$R$-modules $\uphi(U):\ucE(U)\ra\ucF(U)$ for all open $U\subseteq
Y$, such that the following diagram commutes for all open
$V\subseteq U\subseteq Y$
\begin{equation*}
\xymatrix@C=90pt@R=13pt{
*+[r]{\ucE(V)} \ar[r]_{\uphi(V)} \ar[d]^{\si_{VU}} & *+[l]{\ucF(V)}
\ar[d]_{\si_{VU}'} \\ *+[r]{\ucE(U)} \ar[r]^{\uphi(U)} & *+[l]{\ucF(U),\!{}} }
\end{equation*}
where $\si_{VU}$ is the inclusion map for $\ucE$, and $\si_{VU}'$
the inclusion map for~$\ucF$.

A cosheaf $\ucE$ is called {\it flabby\/} if $\si_{VU}:\ucE(V)\ra\ucE(U)$ is an injective morphism of $R$-modules for all open~$V\subseteq U\subseteq Y$.

If $\ucE$ is a flabby cosheaf, $U\subseteq Y$ is open, and $\al\in\ucE(U)$, define the {\it support\/} $\supp\al$ of $\al$ to be the set of points $x\in U$ such that there does not exist open $V\subseteq U\sm\{x\}$ with $\al\in\Im\si_{VU}$. Then $\supp\al$ is a compact subset of~$U$. 

If $\ucE$ is a flabby cosheaf, $\si_{UY}:\ucE(U)\ra\{\al\in\ucE(Y):\supp\al\subseteq U\}$ is an isomorphism for all open $U\subseteq Y$. Thus, knowing $\ucE$ is equivalent to knowing the $R$-module $\ucE(Y)$ and the compact subsets $\supp\al\subseteq Y$ for each~$\al\in\ucE(Y)$.
\label{kh2def10}
\end{dfn}

Bredon \cite[\S V.1]{Bred2} relates c-soft sheaves and flabby cosheaves on~$Y$:

\begin{thm} Let\/ $Y$ be a paracompact, locally compact, Hausdorff topological space, and\/ $R$ a commutative ring. Then
\begin{itemize}
\setlength{\itemsep}{0pt}
\setlength{\parsep}{0pt}
\item[{\bf(a)}] Let $\cE$ be a c-soft sheaf of\/ $R$-modules on $Y$. Define $\ucE(U)=\cE_\cs(U)$ for open $U\subseteq Y$. For all open $V\subseteq U\subseteq Y,$ define $\si_{VU}:\ucE(V)\ra\ucE(U)$ such that\/ $\si_{VU}(\al)\in \ucE(U)\subseteq\cE(U)$ is the unique element with $\rho_{UV}\ci\si_{VU}(\al)=\al$ in $\cE(V)$ and\/ $\rho_{U(U\sm\supp\al)}\ci\si_{VU}(\al)=0$ in $\cE(U\sm\supp\al),$ for all\/~$\al\in\ucE(V)$. 

Then $\ucE$ is a flabby cosheaf of\/ $R$-modules on $Y$.
\item[{\bf(b)}] Let\/ $\ucE$ be a flabby cosheaf of\/ $R$-modules on $Y$. For open $U\subseteq Y,$ define $\cP\cE(U)=\ucE(Y)/\{\al\in\ucE(Y):\supp(\al)\cap U=\es\}$. For open $V\subseteq U\subseteq Y,$ define $\rho_{UV}:\cP\cE(U)\ra\cP\cE(V)$ to be the quotient morphism of\/ $\id:\ucE(Y)\ra\ucE(Y)$. Then $\cP\cE$ is a presheaf of $R$-modules on $Y$. Write $\cE$ for the sheafification of\/ $\cP\cE$. Then $\cE$ is a c-soft sheaf of\/ $R$-modules on $Y$.
For all open $U\subseteq Y,$ define $i_{\ucE,U}:\ucE(U)\ra \cE(U)$ to be the composition
\begin{equation*}
\xymatrix@C=9.5pt{ \ucE(U) \ar[rr]^{\si_{UY}} && \ucE(Y) \ar[r] &
\frac{\ts\ucE(Y)}{\ts\{\al\!\in\!\ucE(Y):\supp(\al)\!\cap\! U\!=\!\es\}} =\cP\cE(U) \ar[rr]^(0.77){\text{sheafify}} && \cE(U). }
\end{equation*}
Then $i_{\ucE,U}$ is an isomorphism $\ucE(U)\ra\cE_\cs(U)\subseteq\cE(U)$. It preserves supports, that is, $\supp \al=\supp i_{\ucE,U}(\al)$ for all\/~$\al\in\ucE(U)$.

The stalks $\cE_y$ of\/ $\cE$ for $y\in Y$ have canonical isomorphisms
\e
\cE_y\cong \ucE(Y)/\si_{(Y\sm\{y\})Y}\bigl[\ucE(Y\sm\{y\})\bigr]
\label{kh2eq32}
\e
compatible with the composition $\ucE(Y)=\cP\cE(Y)\,{\buildrel\text{sheafify}\over\longra}\, \cE(Y)\,{\buildrel\text{stalk}\over\longra}\,\cE_y$.
\item[{\bf(c)}] The constructions of\/ {\bf(a),\bf(b)} are inverse, up to canonical isomorphism. 

Suppose $\phi:\cE\ra\cF$ is a morphism of c-soft sheaves of\/ $R$-modules on $Y,$ and define flabby cosheaves of\/ $R$-modules $\ucE,\ucF$ from $\cE,\cF$ as in {\bf(a)}. Define $\uphi(U)=\phi(U)\vert_{\ucE(U)}:\ucE(U)\ra\ucF(U)$ for all open $U\subseteq Y$. Then $\uphi:\ucE\ra\ucF$ is a morphism of cosheaves, and this gives a functorial\/ $1$-$1$ correspondence between morphisms $\phi:\cE\ra\cF$ and morphisms $\uphi:\ucE\ra\ucF$.

Thus, we have constructed an equivalence of categories between the category of c-soft sheaves of\/ $R$-modules on $Y,$ and the category of flabby cosheaves of\/ $R$-modules on~$Y$.
\end{itemize}
\label{kh2thm3}
\end{thm}

When we compute the compactly-supported cohomology $H^k_\cs(Y,\cE)$ of a sheaf $\cE$ using a c-soft resolution $\cF^\bu=(\cdots\,{\buildrel\d\over\longra}\,\cF^j\,{\buildrel\d\over\longra}\,\cF^{j+1}\,{\buildrel\d\over\longra}\,\cdots)$, as in \eq{kh2eq24} and \eq{kh2eq27}, we use $\cF^j_\cs(Y)$, which is just the global sections $\ucF^j(Y)$ of the flabby cosheaf $\ucF^j$ corresponding to $\cF^j$. So for defining compactly-supported cohomology $H^k_\cs(Y;R)\cong H^k_\cs(Y,R_Y)$ or homology $H_k(Y;R)\cong H^{\dim Y-k}_\cs(Y,O_Y)$ of manifolds, it is more natural to use cosheaves than sheaves.

Following Bredon \cite[\S V.1.3, \S V.1.18, \S VI.12 \& Th.~V.12.14]{Bred2} and Skljarenko \cite{Sklj2}, we define a modified version of singular homology, in which the chains form flabby cosheaves.

\begin{ex} Let $Y$ be a paracompact, locally compact, Hausdorff topological space, and $R$ a commutative ring. Then for $k=0,1,\ldots,$ Example \ref{kh2ex1} defined the singular chains $C_k^\rsi(Y;R)$. Observe that mapping $U\mapsto C_k^\rsi(U;R)$ for all $U\subseteq Y$ open defines a precosheaf of $R$-modules on $Y$, in the sense of Definition \ref{kh2def10}. However, it is not a cosheaf, since if $U,V\subseteq Y$ are open, then
\e
{}\!\!\xymatrix@C=7.5pt{ C_k^\rsi(U\!\cap\! V;R) \ar[rrr]^(0.42){{\rm inc}_*\op -{\rm inc}_*} &&& C_k^\rsi(U;R)\!\op\! C_k^\rsi(V;R) \ar[rrr]^(0.55){{\rm inc}_*\op {\rm inc}_*} &&& C_k^\rsi(U\!\cup\! V;R) \ar[r] & 0 }\!\!{}
\label{kh2eq33}
\e
from \eq{kh2eq31} need not be exact at the third term, so Definition \ref{kh2def10}(iv) fails. This is because there are continuous maps $\si:\De_k\ra U\cup V$ with $\si(\De_k)\not\subseteq U$ and $\si(\De_k)\not\subseteq V$, so $\si\in C_k^\rsi(U\cup V)$ does not lie in the image of~$C_k^\rsi(U)\op C_k^\rsi(V)$.

Define an $R$-module $\hat C_k^\rsi(Y;R)$ of {\it cosheaf singular chains\/} to be the direct limit of the directed system
\begin{equation*}
\xymatrix@C=18pt{C_k^\rsi(Y;R) \ar[r]^B & C_k^\rsi(Y;R) \ar[r]^B & C_k^\rsi(Y;R) \ar[r]^(0.6)B & \cdots, }
\end{equation*}
where $B$ is the barycentric subdivision morphism from Example \ref{kh2ex3}. Then we have a commutative diagram of $R$-modules with a universal property
\e
\begin{gathered}
\xymatrix@C=18pt@R=13pt{C_k^\rsi(Y;R) \ar[r]^B \ar@/_1.1pc/[drrrr]_(0.15){\Pi_0} & C_k^\rsi(Y;R) \ar[r]^B \ar@/_.4pc/[drrr]_(0.23){\Pi_1} & C_k^\rsi(Y;R) \ar[r]^(0.6)B \ar[drr]_(0.3){\Pi_2} & \cdots \\
&&&& \hat C_k^\rsi(Y;R). }
\end{gathered}
\label{kh2eq34}
\e
Since $B\ci\pd=\pd\ci B:C_k^\rsi(Y;R)\ra C_{k-1}^\rsi(Y;R)$ for $k>0$ as in Example \ref{kh2ex3}, there is a unique morphism $\pd:\hat C_k^\rsi(Y;R)\ra\hat C_{k-1}^\rsi(Y;R)$ with $\Pi_j\ci\pd=\pd\ci\Pi_j:C_k^\rsi(Y;R)\ra\hat C_{k-1}^\rsi(Y;R)$ for all~$j=0,1,\ldots.$ 

We have $\pd\ci\pd=0:\hat C_k^\rsi(Y;R)\ra\hat C_{k-2}^\rsi(Y;R)$ as this holds on $C_*^\rsi(Y;R)$. Hence $\bigl(\hat C_*^\rsi(Y;R),\pd\bigr)$ is a chain complex. Define the {\it cosheaf singular homology\/} $\hat H_k^\rsi(Y;R)$ to be the $k^{\rm th}$ homology group of this complex. As $\Pi_j\ci\pd=\pd\ci\Pi_j:C_k^\rsi(Y;R)\ra\hat C_{k-1}^\rsi(Y;R)$, the $\Pi_j$ induce morphisms $(\Pi_j)_*:H_k^\rsi(Y;R)\ra\hat H_k^\rsi(Y;R)$. As in Bredon \cite[\S V.1.3 \& \S VI.12]{Bred2}, these $(\Pi_j)_*$ are isomorphisms, and are independent of~$j=0,1,\ldots.$

If $f:Y_1\ra Y_2$ is a continuous map of topological spaces then since $B\ci f_*=f_*\ci B:C_k^\rsi(Y_1;R)\ra C_k^\rsi(Y_2;R)$ as in Example \ref{kh2ex3}, there is a unique pushforward $f_*:\hat C_k^\rsi(Y_1;R)\ra\hat C_k^\rsi(Y_2;R)$ with $\Pi_j\ci\pd=\pd\ci\Pi_j:C_k^\rsi(Y;R)\ra\hat C_{k-1}^\rsi(Y;R)$ for all $j=0,1,\ldots.$ These $f_*$ are functorial, and induce morphisms $f_*:\hat H_k^\rsi(Y_1;R)\ra\hat H_k^\rsi(Y_2;R)$ on homology, which are identified with the usual pushforwards $f_*:H_k^\rsi(Y_1;R)\ra H_k^\rsi(Y_2;R)$ by the isomorphisms~$(\Pi_j)_*:H_k^\rsi(Y_a;R)\ra\hat H_k^\rsi(Y_a;R)$. 

If $Z\subseteq Y$ is open with inclusion $i:Z\hookra Y$ we define the {\it relative cosheaf singular chains\/} $\hat C_k^\rsi(Y,Z;R)=\hat C_k^\rsi(Y;R)/i_*\bigl(\hat C_k^\rsi(Z;R)\bigr)$. Then $\pd:\hat C_k^\rsi(Y;R)\ra\hat C_{k-1}^\rsi(Y;R)$ induces $\pd:\hat C_k^\rsi(Y,Z;R)\ra\hat C_{k-1}^\rsi(Y,Z;R)$ with $\pd\ci\pd=0$. Define the {\it relative cosheaf singular homology\/} $\hat H_k^\rsi(Y,Z;R)$ to be the $k^{\rm th}$ homology group of the chain complex $\bigl(\hat C_*^\rsi(Y,Z;R),\pd\bigr)$. The morphisms $\Pi_j:C_k^\rsi(Y;R)\ra\hat C_k^\rsi(Y;R)$ induce $\Pi_j:C_k^\rsi(Y,Z;R)\ra\hat C_k^\rsi(Y,Z;R)$ with $\Pi_j\ci\pd=\pd\ci\Pi_j$, so they descend to isomorphisms $(\Pi_j)_*:H_k^\rsi(Y,Z;R)\ra\hat H_k^\rsi(Y,Z;R)$, which are independent of~$j$.

Now, for each open $U\subseteq Y$, define $\hat\ucC{}_k^\rsi(Y;R)(U)=\hat C_k^\rsi(U;R)$, and for all open $V\subseteq U\subseteq Y$, define $\si_{VU}:\hat\ucC{}_k^\rsi(Y;R)(V)\ra\hat\ucC{}_k^\rsi(Y;R)(U)$ by $\si_{VU}=i_*:\hat C_k^\rsi(V;R)\ra\hat C_k^\rsi(U;R)$, for $i:V\hookra U$ the inclusion. Then as in Bredon \cite[\S V.1.3 \& \S VI.12]{Bred2}, this defines a flabby cosheaf $\hat\ucC{}_k^\rsi(Y;R)$ of $R$-modules on $Y$, which we call the $k^{\rm th}$ {\it singular cosheaf}. In particular, the analogue of \eq{kh2eq33} for $\hat C_k(-;R)$ is exact at the third term, since given any chain $\al\in C_k^\rsi(U\cup V;R)$, we may write $B^n(\al)$ for $n\gg 0$ as the sum of chains in $C_k^\rsi(U;R)$ and $C_k^\rsi(V;R)$, as in the final part of Example~\ref{kh2ex3}. 

The morphisms $\pd:\hat C_k^\rsi(U;R)\ra\hat C_{k-1}^\rsi(U;R)$ for open $U\subseteq Y$ define a morphism of cosheaves $\pd:\hat\ucC{}_k^\rsi(Y;R)\ra \hat\ucC{}_{k-1}^\rsi(Y;R)$, with $\pd\ci\pd=0$. This gives a complex of flabby cosheaves on~$Y$:
\e
\hat\ucC{}_{-\bu}^\rsi(Y;R)\!=\!\bigl(\xymatrix@C=14pt{\cdots \ar[r]^(0.35)\pd  & \hat\ucC{}_2^\rsi(Y;R) \ar[r]^\pd & \hat\ucC{}_1^\rsi(Y;R)\ar[r]^\pd & \hat\ucC{}_0^\rsi(Y;R) \ar[r]^(0.65)\pd & 0\bigr), }
\label{kh2eq35}
\e
where we put $\hat\ucC{}_k^\rsi(Y;R)$ in degree~$-k$.
\label{kh2ex15}
\end{ex}

When $Y$ is a manifold, we can do the same using smooth singular homology:

\begin{ex} Let $Y$ be a manifold, and $R$ a commutative ring. Then we can repeat the whole of Example \ref{kh2ex15} using smooth singular chains $C_k^\ssi(Y;R)$ from Example \ref{kh2ex2} instead of singular chains $C_k^\rsi(Y;R)$ from Example \ref{kh2ex1}. In this way we define $R$-modules $\hat C_k^\ssi(Y;R)$ of {\it cosheaf smooth singular chains\/} in a diagram 
\e
\begin{gathered}
\xymatrix@C=18pt@R=13pt{C_k^\ssi(Y;R) \ar[r]^B \ar@/_1.1pc/[drrrr]_(0.15){\Pi_0} & C_k^\ssi(Y;R) \ar[r]^B \ar@/_.4pc/[drrr]_(0.23){\Pi_1} & C_k^\ssi(Y;R) \ar[r]^(0.6)B \ar[drr]_(0.3){\Pi_2} & \cdots \\
&&&& \hat C_k^\ssi(Y;R), }
\end{gathered}
\label{kh2eq36}
\e
as for \eq{kh2eq34}, and morphisms $\pd:\hat C_k^\ssi(Y;R)\ra\hat C_{k-1}^\ssi(Y;R)$ with $\pd^2=0$, and we define the {\it cosheaf smooth singular homology\/} $\hat H_k^\ssi(Y;R)$ to be the $k^{\rm th}$ homology group of $\bigl(\hat C_*^\ssi(Y;R),\pd\bigr)$. The morphisms $\Pi_j:C_k^\ssi(Y;R)\ra\hat C_k^\ssi(Y;R)$ induce isomorphisms $(\Pi_j)_*:H_k^\ssi(Y;R)\ra\hat H_k^\ssi(Y;R)$ which are independent of~$j$. 

If $f:Y_1\ra Y_2$ is a smooth map of manifolds then since $B\ci f_*=f_*\ci B:C_k^\ssi(Y_1;R)\ra C_k^\ssi(Y_2;R)$, there are functorial pushforwards $f_*:\hat C_k^\ssi(Y_1;R)\ra\hat C_k^\ssi(Y_2;R)$ with $\Pi_j\ci\pd=\pd\ci\Pi_j:C_k^\ssi(Y;R)\ra\hat C_{k-1}^\ssi(Y;R)$ for all $j$, which induce morphisms $f_*:\hat H_k^\ssi(Y_1;R)\ra\hat H_k^\ssi(Y_2;R)$ on homology.

We define a flabby cosheaf $\hat\ucC{}_k^\ssi(Y;R)$ of $R$-modules on $Y$ called the $k^{\rm th}$ {\it smooth singular cosheaf}, with $\hat\ucC{}_k^\ssi(Y;R)(U)=\hat C_k^\ssi(U;R)$ for each open $U\subseteq Y$. We have morphisms $\pd:\hat\ucC{}_k^\ssi(Y;R)\ra\hat\ucC{}_{k-1}^\ssi(Y;R)$ with $\pd\ci\pd=0$. So as in \eq{kh2eq35} we have a complex of flabby cosheaves on~$Y$:
\begin{equation*}
\hat\ucC{}_{-\bu}^\ssi(Y;R)\!=\!\bigl(\xymatrix@C=14pt{\cdots \ar[r]^(0.35)\pd  & \hat\ucC{}_2^\ssi(Y;R) \ar[r]^\pd & \hat\ucC{}_1^\ssi(Y;R)\ar[r]^\pd & \hat\ucC{}_0^\ssi(Y;R) \ar[r]^(0.65)\pd & 0\bigr), }
\end{equation*}
where we put $\hat\ucC{}_k^\ssi(Y;R)$ in degree~$-k$.
\label{kh2ex16}
\end{ex}

Applying Theorem \ref{kh2thm3}, we can transform the flabby cosheaves $\hat\ucC{}_k^\ssi(Y;R)$ in Examples \ref{kh2ex15} and \ref{kh2ex16} into (c-)soft sheaves. 

\begin{ex} Let $Y$ be a manifold and $R$ a commutative ring, and use the notation of Example \ref{kh2ex15}. For each $k=0,1,\ldots,$ define the $k^{\rm th}$ {\it singular sheaf\/} $\hat\cC{}^{\lf,\rsi}_k(Y;R)$ to be the c-soft sheaf of $R$-modules on $Y$ associated to the flabby cosheaf $\hat\ucC{}_k^\rsi(Y;R)$ by Theorem \ref{kh2thm3}(b). As $Y$ is a manifold, $\hat\ucC{}_k^\rsi(Y;R)$ is a soft sheaf. By Theorem \ref{kh2thm3}(c), the morphisms $\pd:\hat\ucC{}_k^\rsi(Y;R)\ra\hat\ucC{}_{k-1}^\rsi(Y;R)$ lift to $\pd:\hat\cC{}^{\lf,\rsi}_k(Y;R)\ra\hat\cC{}^{\lf,\rsi}_{k-1}(Y;R)$ with $\pd\ci\pd=0$. So corresponding to \eq{kh2eq35} we have a complex of soft sheaves on $Y$, with $\hat\cC{}^{\lf,\rsi}_k(Y;R)$ in degree $-k$:
\e
\hat\cC{}^{\lf,\rsi}_{-\bu}(Y;R)\!=\!\bigl(\xymatrix@C=10pt{\cdots \ar[r]^(0.3)\pd  & \hat\cC{}^{\lf,\rsi}_2(Y;R) \ar[r]^\pd & \hat\cC{}^{\lf,\rsi}_1(Y;R) \ar[r]^\pd & \hat\cC{}^{\lf,\rsi}_0(Y;R) \ar[r]^(0.7)\pd & 0\bigr). }
\label{kh2eq37}
\e

The stalk $\hat\cC{}^{\lf,\rsi}_k(Y;R)_y$ is $\hat C_k^\rsi(Y;R)/\hat C_k^\rsi(Y\sm\{y\};R)=\hat C_k^\rsi(Y,Y\sm\{y\};R)$ at $y\in Y$ by \eq{kh2eq32}, so the cohomology in degree $-k$ of the complex of stalks at $y$ from \eq{kh2eq37} is $\hat H_k^\rsi(Y,Y\sm\{y\};R)\cong H_k(Y,Y\sm\{y\};R)$, which as $Y$ is a manifold is 0 if $k\ne\dim Y$, and the $R$-module of orientations on $T_yY$ if $k=\dim Y$. Therefore \eq{kh2eq37} is equivalent in $D(Y;R)$ to the dualizing complex $\om_Y\simeq O_Y[\dim Y]$. 

Define the {\it locally finite sheaf singular chains\/} $\hat C_k^{\lf,\rsi}(Y;R)$ to be the global sections $\hat\cC{}^{\lf,\rsi}_k(Y;R)(Y)$. Then the $\pd:\hat\cC{}^{\lf,\rsi}_k(Y;R)\ra\hat\cC{}^{\lf,\rsi}_{k-1}(Y;R)$ induce morphisms $\pd:\hat C_k^{\lf,\rsi}(Y;R)\ra\hat C_{k-1}^{\lf,\rsi}(Y;R)$ making $\bigl(\hat C_*^{\lf,\rsi}(Y;R),\pd\bigr)$ into a chain complex. Define the {\it locally finite sheaf singular homology\/} $\hat H_k^{\lf,\rsi}(Y;R)$ to be the $k^{\rm th}$ homology group of $\bigl(\hat C_*^{\lf,\rsi}(Y;R),\pd\bigr)$. 

Composing $\Pi_0:C_k^\rsi(Y;R)\ra\hat C_k^\rsi(Y;R)$ with the inclusion $\hat C_k^\rsi(Y;R)\hookra\hat C_k^{\lf,\rsi}(Y;R)$ coming from realizing sections of $\hat\ucC{}_k^\rsi(Y;R)$ as sections of the associated sheaf $\hat\cC{}^{\lf,\rsi}_k(Y;R)$ gives a morphism $C_k^\rsi(Y;R)\ra\hat C_k^{\lf,\rsi}(Y;R)$. Now Example \ref{kh2ex12} defined the locally finite singular chains $C_k^{\lf,\rsi}(Y;R)$, which are locally finite sums of elements in~$C_k^\rsi(Y;R)$. 

Since $\hat C_k^{\lf,\rsi}(Y;R)$ is the global sections of a sheaf, locally finite sums make sense in $\hat C_k^{\lf,\rsi}(Y;R)$. Thus, this morphism $C_k^\rsi(Y;R)\ra\hat C_k^{\lf,\rsi}(Y;R)$ extends naturally to a morphism $\Pi_0^\lf:C_k^{\lf,\rsi}(Y;R)\ra \hat C_k^{\lf,\rsi}(Y;R)$. These satisfy $\Pi_0^\lf\ci\pd=\pd\ci\Pi_0^\lf:C_k^{\lf,\rsi}(Y;R)\ra \hat C_{k-1}^{\lf,\rsi}(Y;R)$, and so induce morphisms $(\Pi_0^\lf)_*:H_k^{\lf,\rsi}(Y;R)\ra\hat H_k^{\lf,\rsi}(Y;R)$, which are isomorphisms. Thus, $\hat H_k^{\lf,\rsi}(Y;R)$ is another version of locally finite homology.

There is a natural equivalence in the derived category $D(Y;R)$
\e
\hat\cC{}^{\lf,\rsi}_{-\bu}(Y;R)\simeq \om_Y.
\label{kh2eq38}
\e
This is proved by Arabia \cite[Th.~1.8.6(a)]{Arab} when $Y$ is a `pseudovariety' (a topological space stratified by topological manifolds, which include orbifolds), see also Borel \cite[\S V.7.2]{Bore}, Bredon \cite[Th.~V.12.14]{Bred2}, and Kashiwara and Schapira \cite[Th.~9.2.10]{KaSc}. Taking the hypercohomology of \eq{kh2eq38} in degree $-k$ and using softness of the $\hat\cC{}^{\lf,\rsi}_i(Y;R)$ gives the canonical isomorphism from~\eq{kh2eq28}:
\begin{equation*}
\hat H_k^{\lf,\rsi}(Y;R)\cong 
\H^{-k}(Y,\hat\cC{}^{\lf,\rsi}_{-\bu}(Y;R))\cong \H^{-k}(Y,\om_Y).
\end{equation*}

If $i:U\hookra Y$ is the inclusion of an open set, then as in Property \ref{kh2pr2}(c) there is a natural morphism $i^*:H_k^\lf(Y;R)\ra H_k^\lf(U;R)$ on locally finite homology. We noted in Example \ref{kh2ex12} that $i^*$ does not have a nice expression on the chain level for locally finite singular homology $H_k^{\lf,\rsi}(Y;R)$. 

However, for the sheafified version $\hat H_k^{\lf,\rsi}(Y;R)$, there is a nice expression. As $\hat\cC{}^{\lf,\rsi}_k(Y;R)\vert_U=\hat\cC{}^{\lf,\rsi}_k(U;R)$, the data $\si_{YU}$ in $\hat\cC{}^{\lf,\rsi}_k(Y;R)$ is a morphism $i^*=\si_{YU}:\hat C_k^{\lf,\rsi}(Y;R)\ra\hat C_k^{\lf,\rsi}(U;R)$. These satisfy $i^*\ci\pd=\pd\ci i^*:\hat C_k^{\lf,\rsi}(Y;R)\ra\hat C_{k-1}^{\lf,\rsi}(U;R)$, since $\pd:\hat\cC{}^{\lf,\rsi}_k(Y;R)\ra\hat\cC{}^{\lf,\rsi}_{k-1}(Y;R)$ is a sheaf morphism, and so induce the canonical morphisms $i^*:\hat H_k^{\lf,\rsi}(Y;R)\ra\hat H_k^{\lf,\rsi}(U;R)$, as required.

We can repeat all of the above for the smooth singular cosheaves $\hat\ucC{}_k^\ssi(Y;R)$ from Example \ref{kh2ex16} rather than $\hat\ucC{}_k^\rsi(Y;R)$ from Example \ref{kh2ex15}. This defines a complex $\bigl(\hat\cC{}^{\lf,\ssi}_*(Y;R),\pd\bigr)$ of {\it smooth singular sheaves\/} on $Y$, with global sections the {\it locally finite sheaf smooth singular chains\/} $\bigl(\hat C_*^{\lf,\ssi}(Y;R),\pd\bigr)$ on $Y$, which has homology the {\it locally finite sheaf smooth singular homology\/} $\hat H_k^{\lf,\ssi}(Y;R)$ of $Y$, which is canonically isomorphic to locally finite homology. As for \eq{kh2eq38} we have a natural equivalence in the derived category~$D(Y;R)$:
\e
\hat\cC{}^{\lf,\ssi}_{-\bu}(Y;R)\simeq \om_Y.
\label{kh2eq39}
\e

\label{kh2ex17}
\end{ex}

The next definition and theorem are new, so far as the author knows. We will define a class of {\it strong\/} presheaves, and discuss their properties.

\begin{dfn} Let $Y$ be a topological space, and $R$ a commutative ring. We say that a presheaf of $R$-modules $\cE$ on $Y$ is {\it strong\/} if for all open $U,V\subseteq Y$, the following sequence is exact in $\Rmod$:
\e
\xymatrix@C=12pt{ 0 \ar[r] & \cE(U\!\cup\! V) \ar[rrrr]^(0.48){\rho_{(U\cup V)U}\op \rho_{(U\cup V)V}} &&&& \cE(U)\!\op\! \cE(V) \ar[rrrr]^(0.53){\rho_{U(U\cap V)}\op -\rho_{V(U\cap V)}} &&&& \ucE(U\!\cap\! V). }\!\!{}
\label{kh2eq40}
\e

\label{kh2def11}
\end{dfn}

The next theorem will be proved in \S\ref{kh62}.

\begin{thm} Let\/ $Y$ be a paracompact, locally compact, Hausdorff topological space, and\/ $R$ a commutative ring. Suppose $\cE$ is a strong presheaf of\/ $R$-modules on $Y,$ with sheafification $\pi:\cE\ra\hat\cE$. Then:
\begin{itemize}
\setlength{\itemsep}{0pt}
\setlength{\parsep}{0pt}
\item[{\bf(a)}] The presheaf\/ $\cE$ satisfies the sheaf conditions Definition\/ {\rm\ref{kh2def1}(iv),(v)} for the open cover $\{V_i:i\in I\}$ whenever $I$ is a finite set. Conversely, any presheaf\/ $\cE'$ satisfying Definition\/ {\rm\ref{kh2def1}(iv),(v)} whenever $I$ is finite, is strong.
\item[{\bf(b)}] If\/ $V\subseteq U\subseteq Y$ are open then using the notation of Definition\/ {\rm\ref{kh2def6},} $\rho_{UV}:\cE(U)\ra\cE(V)$ restricts to an isomorphism
\e
\rho_{UV}\vert_{\cdots}:\bigl\{\al\in \cE_\cs(U):\supp\al\subseteq V\subseteq U\bigr\}\,{\buildrel\cong\over\longra}\,\cE_\cs(V).
\label{kh2eq41}
\e
\item[{\bf(c)}] Suppose $V\subseteq U\subseteq Y$ are open, and the closure $\bar V$ of\/ $V$ in $U$ is compact. Then there is a canonical\/ $R$-module morphism\/ $\tau_{UV}:\hat\cE(U)\ra\cE(V)$ making the following diagram commute:
\e
\begin{gathered}
\xymatrix@C=100pt@R=15pt{ *+[r]{\cE(U)} \ar[r]_(0.3){\rho_{UV}} \ar[d]^{\pi(U)} & *+[l]{\cE(V)} \ar[d]_{\pi(V)} \\ 
 *+[r]{\hat\cE(U)} \ar[r]^(0.6){\hat\rho_{UV}} \ar[ur]^(0.3){\tau_{UV}} & *+[l]{\hat\cE(V).\!} }
\end{gathered}
\label{kh2eq42}
\e
For each\/ $\hat s\in\hat\cE(U),$ we can characterize $\tau_{UV}(\hat s)\in\cE(V)$ uniquely by the following property: if\/ $y\in\bar V,$ and\/ $U_y$ is an open neighbourhood of\/ $y$ in $U,$ and\/ $s_y\in\cE(U_y)$ with\/ $\pi(U_y)(s_y)=\rho_{UU_y}(\hat s),$ then there exists an open neighbourhood\/ $U_y'$ of\/ $y$ in $U_y$ such that\/ $\rho_{U_y(U'_y\cap V)}(s_y)=\rho_{V(U'_y\cap V)}\ci\tau_{UV}(\hat s)$ in\/~$\cE(U'_y\cap V)$.

If also\/ $W\subseteq V$ is open and the closure $\bar W$ of\/ $W$ in $U$ is compact then the following diagram commutes:
\e
\begin{gathered}
\xymatrix@C=130pt@R=13pt{ *+[r]{\hat\cE(U)} \ar[d]_{\tau_{UV}} \ar[dr]^(0.4){\tau_{UW}} \\
*+[r]{\cE(V)} \ar[r]^(0.4){\rho_{VW}} & *+[l]{\cE(W).} }
\end{gathered}
\label{kh2eq43}
\e

\item[{\bf(d)}] For any open $U\subseteq Y,$ the map $\pi(U):\cE(U)\ra\hat\cE(U)$ restricts to an isomorphism $\pi_\cs(U):\cE_\cs(U)\ra\hat\cE_\cs(U)$. Hence if\/ $U$ is compact then $\pi(U):\cE(U)\ra\hat\cE(U)$ is an isomorphism.
\item[{\bf(e)}] For all open $U\subseteq Y,$ there is a canonical isomorphism
\e
\hat\cE(U)\cong\mathop{\underleftarrow{\lim}\,}\nolimits_{\text{$V:V\subseteq U$ open, $\bar V$ is compact}}\cE(V).
\label{kh2eq44}
\e
Here the inverse limit is over open $V\subseteq U$ for which the closure $\bar V$ of $V$ in $U$ is compact under the morphisms $\rho_{V_1V_2}:\cE(V_1)\ra\cE(V_2)$ for $V_2\subseteq V_1\subseteq U,$ and the projection $\hat\cE(U)\ra\cE(V)$ from \eq{kh2eq44} is $\tau_{UV}$ from part\/~{\bf(c)}.

\item[{\bf(f)}] If\/ $\cE$ is c-soft then $\hat\cE$ is c-soft, in the sense of Definition\/~{\rm\ref{kh2def8}}.

Then as in Theorem\/ {\rm\ref{kh2thm3}(a)} we may define a natural flabby cosheaf\/ $\ucE$ on $Y$ with\/ $\ucE(U)=\cE_\cs(U)$ for $U\subseteq Y$ open, and\/ $\si_{VU}:\ucE(V)\ra\ucE(U)$ for $V\subseteq U\subseteq Y$ open is the inverse of the isomorphism {\rm\eq{kh2eq41},} and\/ $\hat\cE$ is the c-soft sheaf associated to $\ucE$ in Theorem\/~{\rm\ref{kh2thm3}(b)}.
\end{itemize}
\label{kh2thm4}
\end{thm}

In (c) it is essential that the characterizing property of $\tau_{UV}(\hat s)$ holds for all $y\in\bar V$ with $\bar V$ compact, not just for all $y\in V$. Since $\pi(U),\pi(V)$ in \eq{kh2eq42} need not be injective or surjective, the fact that \eq{kh2eq42} commutes does not itself determine $\tau_{UV}$ uniquely. The moral of the theorem is that strong presheaves are quite close to being sheaves, and the sheafification $\hat\cE$ of a strong presheaf $\cE$ is quite close to $\cE$, for example, $\cE,\hat\cE$ agree on compact open sets $U\subseteq Y$, and on compactly-supported sections. Equation \eq{kh2eq44} gives a useful alternative expression for the sheafification of a strong presheaf.

\subsection{Products on (co)homology}
\label{kh26}

In \S\ref{kh21}--\S\ref{kh25} we considered $H_k(Y;R),H^k(Y;R),H^k_\cs(Y;R),H_k^\lf(Y;R)$ just as $R$-modules. In fact, homology and cohomology have many interesting additional algebraic structures. We now discuss {\it cup}, {\it cap\/} and {\it cross products\/}~$\cup,\cap,\t$.

\subsubsection{The different products and their properties}
\label{kh261}

Let $Y$ be a topological space, and $R$ a commutative ring. Then as in \cite[\S VI.4]{Bred1}, \cite[\S 3.2]{Hatc}, \cite[\S 48]{Munk}, for all $k,l\ge 0$ we have  $R$-bilinear {\it cup products\/}
\e
\begin{split}
\cup&:H^k(Y;R)\t H^l(Y;R)\longra H^{k+l}(Y;R),\\
\cup&:H^k_\cs(Y;R)\t H^l(Y;R)\longra H^{k+l}_\cs(Y;R),\\
\cup&:H^k(Y;R)\t H^l_\cs(Y;R)\longra H^{k+l}_\cs(Y;R),\\
\cup&:H^k_\cs(Y;R)\t H^l_\cs(Y;R)\longra H_\cs^{k+l}(Y;R).
\end{split}
\label{kh2eq45}
\e
As in \cite[\S VI.5]{Bred1}, \cite[\S 3.3]{Hatc}, \cite[\S 66]{Munk}, we also have $R$-bilinear {\it cap products\/}
\e
\begin{split}
\cap&:H^k(Y;R)\t H_l(Y;R)\longra H_{l-k}(Y;R),\\
\cap&:H^k_\cs(Y;R)\t H_l(Y;R)\longra H_{l-k}(Y;R),\\
\cap&:H^k(Y;R)\t H_l^\lf(Y;R)\longra H_{l-k}^\lf(Y;R),\\
\cap&:H^k_\cs(Y;R)\t H_l^\lf(Y;R)\longra H_{l-k}(Y;R).
\end{split}
\label{kh2eq46}
\e

Writing $H_k^?(Y;R)$ to mean either $H_k(Y;R)$ or $H_k^\lf(Y;R)$, and 
$H^k_?(Y;R)$ to mean either $H^k(Y;R)$ or $H^k_\cs(Y;R)$, for $\al\in H^k_?(Y;R)$, $\be\in H^l_?(Y;R)$, $\ga\in H^m_?(Y;R)$, $\de\in H_n^?(Y;R)$, these satisfy the identities
\ea
\al\cup\be&=(-1)^{kl}\be\cup\al, 
\label{kh2eq47}\\
\al\cup(\be\cup\ga)&=(\al\cup\be)\cup\ga, 
\label{kh2eq48}\\
\al\cap(\be\cap\de)&=(\al\cup\be)\cap\de.
\label{kh2eq49}
\ea
There is a natural {\it identity element\/} $[1_Y]\in H^0(Y;R)$ with
\ea
[1_Y]\cup\al=\al\cup [1_Y]&=\al, 
\label{kh2eq50}\\ 
[1_Y]\cap\de&=\de.
\label{kh2eq51}
\ea
In fact $[1_Y]=\pi^*(1)$, where $\pi:Y\ra *$ is the projection and $1\in H_0(*;R)\cong R$. Also $\cup,\cap$ are compatible with the projections $\Pi:H^*_\cs(Y;R)\ra H^*(Y;R)$ and $\Pi:H_*(Y;R)\ra H_*^\lf(Y;R)$, so that $\Pi(\al\cup\be)=\Pi(\al)\cup\Pi(\be)$, and so on. Thus, $\cup,[1_Y]$ make $H^*(Y;R)$ into a unital supercommutative graded $R$-algebra, and $\cup$ makes $H^*_\cs(Y;R)$ into a non-unital supercommutative graded $R$-algebra, and $\cap$ makes $H_*(Y;R),H_*^\lf(Y;R)$ into graded modules over~$H^*(Y;R),H^*_\cs(Y;R)$.

The {\it Kronecker products\/} \cite[\S VI.3]{Bred1} are the compositions
\begin{equation*}
\xymatrix@C=28pt@R=11pt{ 
H^k(Y;R)\t H_k(Y;R) \ar[r]_(0.6)\cap & H_0(Y;R) \ar[r]_{\pi_*} & H_0(*;R) \ar[r]_(0.58)\cong & R, \\ 
H^k_\cs(Y;R)\t H_k^\lf(Y;R) \ar[r]^(0.6)\cap & H_0(Y;R) \ar[r]^{\pi_*} & H_0(*;R) \ar[r]^(0.58)\cong & R, }
\end{equation*}
where $\pi:Y\ra *$ is the projection. If $R=\K$ is a field, these are perfect pairings, inducing isomorphisms $H^k(Y;\K)\cong H_k(Y;\K)^*$, $H_k^\lf(Y;\K)\cong H^k_\cs(Y;\K)^*$, as in Example \ref{kh2ex4} and Property~\ref{kh2pr2}(g).

If $Y_1,Y_2$ are topological spaces (e.g. manifolds), as in \cite[\S VI.1, \S VI.3]{Bred1}, \cite[\S 3.2, \S 3.B]{Hatc}, \cite[\S 59--\S 60]{Munk}, there are also $R$-bilinear {\it cross products\/}
\ea
\begin{split}
\t&:H^k(Y_1;R)\t H^l(Y_2;R)\longra H^{k+l}(Y_1\t Y_2;R),\\
\t&:H^k_\cs(Y_1;R)\t H^l_\cs(Y_2;R)\longra H^{k+l}_\cs(Y_1\t Y_2;R),
\end{split}
\label{kh2eq52}\\
\begin{split}
\t&:H_k(Y_1;R)\t H_l(Y_2;R)\longra H_{k+l}(Y_1\t Y_2;R),\\
\t&:H_k^\lf(Y_1;R)\t H_l^\lf(Y_2;R)\longra H_{k+l}^\lf(Y_1\t Y_2;R).
\end{split}
\label{kh2eq53}
\ea
For $\al\in H^k_?(Y_1;R)$, $\be\in H^l_?(Y_2;R)$, $\ga\in H^m_?(Y_3;R)$, these satisfy the identities
\begin{equation*}
\al\t(\be\t\ga)=(\al\t\be)\t\ga,\;\> \al\t\be=(-1)^{kl}\be\t\al,\;\> [1_{Y_1}]\t [1_{Y_2}]=[1_{Y_1\t Y_2}]
\end{equation*}
in cohomology, where for the second equation we identify $H^{k+l}_?(Y_1\t Y_2;R)=H^{l+k}_?(Y_2\t Y_1;R)$ using the homeomorphism $Y_1\t Y_2\cong Y_2\t Y_1$. Similarly, for $\de\in H_k^?(Y_1;R)$, $\ep\in H_l^?(Y_2;R)$, $\ze\in H_m^?(Y_3;R)$, in homology we have
\begin{equation*}
\de\t(\ep\t\ze)=(\de\t\ep)\t\ze,\qquad \de\t\ep=(-1)^{kl}\ep\t\de. 
\end{equation*}
Also cross products are compatible with the projections $\Pi:H^*_\cs(Y;R)\ra H^*(Y;R)$ and $\Pi:H_*(Y;R)\ra H_*^\lf(Y;R)$.

The cross product on cohomology may be expressed in terms of the cup product, and vice versa, since if $\al,\be\in H^*_?(Y;R)$ then
\begin{equation*}
\al\cup\be=\De_Y^*(\al\t\be),
\end{equation*}
where $\De_Y:Y\ra Y\t Y$ is the diagonal map, and if $\ga_i\in H^*_?(Y_i;R)$ then
\e
\ga_1\t\ga_2=\pi_{Y_1}^*(\ga_1)\cup\pi_{Y_2}^*(\ga_2),
\label{kh2eq54}
\e
where $\pi_{Y_1}:Y_1\t Y_2\ra Y_1$, $\pi_{Y_2}:Y_1\t Y_2\ra Y_2$ are the projections.

All of $\cup,\cap,\t$ and $[1_Y]$ have the obvious functoriality under (proper) pullbacks $f^*$ and pushforwards $f_*$. So, for example, if $f:Y_1\ra Y_2$ is a continuous map of topological spaces and $\al,\be\in H^*(Y_2;R)$, $\de\in H_*(Y_1;R)$ then
\ea
f^*(\al\cup\be)&=f^*(\al)\cup f^*(\be),
\label{kh2eq55}\\
f^*([1_{Y_2}])&=[1_{Y_1}],
\label{kh2eq56}\\
f_*(f^*(\al)\cap\de)&=\al\cap f_*(\de).
\label{kh2eq57}
\ea
All of $\cup,\cap,\t$ extend to relative (co)homology $H^*_?(Y,Z;R),H_*^?(Y,Z;R)$, but for simplicity we will not discuss the relative versions.

The K\"unneth Theorem \cite[Th.~VI.1.6]{Bred1}, \cite[Th.~3B.6]{Hatc}, \cite[Th.~59.3]{Munk} says that on homology, if $R$ is a principal ideal domain, we have functorial exact sequences for all $Y_1,Y_2,m$:
\e
\xymatrix@C=12.3pt{
0 \ar[r] & \bigop\limits_{\begin{subarray}{l} k+l\\ =m \end{subarray}} {\begin{aligned}[t] &H_k(Y_1;R)\ot_R\\[-4pt] & \; H_l(Y_2;R)\end{aligned}} \ar[r]^(0.52)\t & H_m(Y_1\!\t\! Y_2;R) \ar[r] & \bigop\limits_{\begin{subarray}{l} k+l= \\ m-1 \end{subarray}} {\begin{aligned}[t] & H_k(Y_1;R)*{} \\[-4pt] & \; H_l(Y_2;R) \end{aligned}} \ar[r] & 0. }
\label{kh2eq58}
\e
Here $A*B$ is the {\it torsion product\/} of abelian groups $A,B$, \cite[\S V.6]{Bred1}, \cite[\S 3.A]{Hatc}, \cite[\S 54]{Munk}. If $A,B$ are vector spaces over a field $\K$ then $A*B=0$, so when $R=\K$, equation \eq{kh2eq58} yields isomorphisms
\e
\t:\ts\bigop_{k+l=m} H_k(Y_1;\K)\ot_\K H_l(Y_2;\K)\,{\buildrel\cong\over\longra}\, H_m(Y_1\t Y_2;\K).
\label{kh2eq59}
\e
The analogues of \eq{kh2eq58}--\eq{kh2eq59} also hold for $H_*^\lf(-;R)$.

\subsubsection{Defining cup, cap and cross products in sheaf cohomology}
\label{kh262}

Our preferred approach to defining the products $\cup,\cap,\t$ in \S\ref{kh4}--\S\ref{kh5} will be via sheaf cohomology. At the sheaf level, this is very simple: if $Y$ is a topological space there is a natural isomorphism of sheaves of $R$-modules on $Y$:
\e
I_\cup:R_Y\ot_R R_Y\,{\buildrel\cong\over\longra}\, R_Y,
\label{kh2eq60}
\e
which on stalks $R_{Y,y}=R$ at $y\in Y$ acts by the isomorphism $R\ot_R R\,{\buildrel\cong\over\longra}\, R$ induced by multiplication in $R$. The cup products on $H^*(Y;R)=H^*(Y,R_Y)$ and $H^*_\cs(Y;R)=H^*_\cs(Y,R_Y)$ are then induced by the compositions
\begin{equation*}
\xymatrix@C=13pt@R=13pt{ H^k(Y,R_Y)\!\ot_R\! H^l(Y,R_Y) \ar[r] & H^{k+l}(Y,R_Y\!\ot_R\! R_Y) \ar[rrr]_(0.55){H^{k+l}(I_\cup)} &&& H^{k+l}(Y,R_Y), \\
H^k_\cs(Y,R_Y)\!\ot_R\! H^l_\cs(Y,R_Y) \ar[r] & H^{k+l}_\cs(Y,R_Y\!\ot_R\! R_Y) \ar[rrr]^(0.55){H^{k+l}_\cs(I_\cup)} &&& H^{k+l}_\cs(Y,R_Y). }
\end{equation*}

Similarly, for sufficiently nice topological spaces $Y$ and rings $R$ the cap product is induced by a natural quasi-isomorphism of complexes of $R$-modules on~$Y$:
\e
I_\cap:R_Y\ot_R \om_Y\,{\buildrel\simeq\over\longra}\, \om_Y.
\label{kh2eq61}
\e
If $Y$ is a manifold of dimension $m$, so that $\om_Y\simeq O_Y[m]$, this reduces to an isomorphism of sheaves of $R$-modules on $Y$:
\e
I_\cap:R_Y\ot_R O_Y\,{\buildrel\cong\over\longra}\, O_Y,
\label{kh2eq62}
\e
which on stalks at $y\in Y$ acts by the isomorphism $R\ot_R O_{Y,y}\ra O_{Y,y}$ induced by multiplication by $R$. Then the cap product between $H^k(Y;R)=H^k(Y,R_Y)$ and $H_l(Y;R)=H^{m-l}_\cs(Y,O_Y)$ is induced by the composition
\begin{equation*}
\xymatrix@C=7.2pt{ H^k(Y,R_Y)\!\ot_R\! H^{m-l}_\cs(Y,O_Y) \ar[r] & H^{k+m-l}_\cs(Y,R_Y\!\ot_R\! O_Y) \ar[rr]^(0.55){H^*(I_\cap)} && H^{m-(l-k)}_\cs(Y,O_Y). }
\end{equation*}

If $Y_1,Y_2$ are topological spaces then cross products $\t$ on cohomology \eq{kh2eq52} are induced by a natural isomorphism of sheaves of $R$-modules on $Y_1\t Y_2$:
\begin{equation*}
I_\t^\coh:R_{Y_1}\boxtimes_R R_{Y_2}\,{\buildrel\cong\over\longra}\, R_{Y_1\t Y_2}
\end{equation*}
acting on stalks $R_{Y_1,y_1}=R$, $R_{Y_2,y_2}=R$, $R_{Y_1\t Y_2,(y_1,y_2)}=R$ at $(y_1,y_2)\in Y_1\t Y_2$ by the isomorphism $R\ot_R R\ra R$ induced by multiplication in~$R$. 

Similarly, for nice $Y_1,Y_2,R,$ cross products $\t$ on homology \eq{kh2eq53} are induced by a natural quasi-isomorphism of complexes of $R$-modules on $Y_1\t Y_2$:
\begin{equation*}
I_\t^\hom:\om_{Y_1}\boxtimes_R \om_{Y_2}\,{\buildrel\cong\over\longra}\, \om_{Y_1\t Y_2}.
\end{equation*}
If $Y_1,Y_2$ are manifolds of dimensions $m,n$, so that $\om_{Y_1}\simeq O_{Y_1}[m]$, $\om_{Y_2}\simeq O_{Y_2}[n]$, $\om_{Y_1\t Y_2}\simeq O_{Y_1\t Y_2}[m+n]$, this reduces to an isomorphism of sheaves
\begin{equation*}
I_\t^\hom:O_{Y_1}\boxtimes_R O_{Y_2}\,{\buildrel\cong\over\longra}\, O_{Y_1\t Y_2}.
\end{equation*}

As explained in \S\ref{kh25}, one often computes sheaf cohomology using soft resolutions. Suppose $Y$ is a topological space, and we have chosen a soft resolution 
\e
\xymatrix@C=20pt{ 0 \ar[r] & R_Y \ar[r]^i & \cF^0 \ar[r]^\d & \cF^1 \ar[r]^\d  & \cF^2 \ar[r]^(0.45)\d  & \cdots }
\label{kh2eq63}
\e
of $R_Y$, giving isomorphisms $H^k(Y;R)\cong H^k\bigl(0\longra \cF^0(Y)\,{\buildrel\d\over\longra}\,\cF^1(Y)\,{\buildrel\d\over\longra}\,\cdots\bigr)$. Suppose also that we can choose morphisms $\psi_{k,l}:\cF^k\ot_R\cF^l\ra\cF^{k+l}$ of sheaves of $R$-modules on $Y$ for all $k,l\ge 0$, satisfying the conditions, for $I_\cup$ as in~\eq{kh2eq60}:
\e
\begin{split}
i\ci I_\cup&=\psi_{0,0}\ci (i\ot i):R_Y\ot_R R_Y\longra \cF^0,\\
\d\ci\psi_{k,l}&=\psi_{k+1,l}\ci(\d\ot\id_{\cF^l})+(-1)^k\psi_{k,l+1}\ci(\id_{\cF^k}\ot\d):\\
&\qquad\cF^k\ot_R\cF^l\longra \cF^{k+l+1},\qquad \text{for all $k,l\ge 0$.}
\end{split}
\label{kh2eq64}
\e
That is, we have a commutative diagram of sheaves of $R$-modules on $Y$:
\begin{equation*}
\raisebox{67pt}{\text{\begin{small}$\displaystyle\xymatrix@C=9pt@R=7pt{ 0 \ar[r] & R_Y\!\ot_R\! R_Y \ar[r]^{i\ot i} \ar[d]_{I_\cup} & \cF^0\!\ot_R\!\cF^0 \ar[rr]^(0.45){\begin{pmatrix} \st \d\ot\id \\ \st \id\ot\d \end{pmatrix}} \ar[d]_{\psi_{0,0}} && {\raisebox{5.6pt}{$\begin{subarray}{l}\ts \cF^1\!\ot_R\!\cF^0\!\op{} \\ \ts \cF^0\!\ot_R\!\cF^1\end{subarray}$}} \ar[d]_{\begin{pmatrix} \st\psi_{1,0}\!\!\!\! & \st\psi_{0,1} \end{pmatrix}}
\ar[rrrr]^(0.48){\begin{pmatrix} \st \!\d\ot\id \!\!\!\!\! & \st 0 \\ \st \!-\id\ot\d \!\!\!\!\! & \st \d\ot\id \\ \st 0 \!\!\!\!\! & \st \!\id\ot\d \end{pmatrix}}   &&&& {\raisebox{10.3pt}{${\begin{subarray}{l}\ts \cF^2\!\ot_R\!\cF^0\!\op{} \\ \ts \cF^1\!\ot_R\!\cF^1\!\op{} \\ \ts \cF^0\!\ot_R\!\cF^2\end{subarray}}$}} \ar[r] \ar[d]_{\begin{pmatrix} \st\psi_{2,0}\!\!\!\!\! & \st\psi_{1,1}\!\!\!\!\! & \st\psi_{0,2} \end{pmatrix}}  & \cdots 
\\
0 \ar[r] & R_Y \ar[r]^i & \cF^0 \ar[rr]^\d && \cF^1 \ar[rrrr]^\d  &&&& \cF^2 \ar[r]^(0.45)\d  & \cdots,\!\! }$\end{small}}}
\end{equation*}\vskip -38pt
\noindent where the rows are exact, and are soft resolutions of $R_Y\ot_RR_Y,R_Y$.

Hence $(\psi_{k,l})_{k,l\ge 0}$ defines a morphism of complexes $\cF^\bu\ot_R\cF^\bu\ra \cF^\bu$ equivalent to $I_\cup:R_Y\ot_RR_Y\ra R_Y$, and taking global sections, or compactly-supported global sections, gives a formula for the cup product $\cup$ at the cochain level. So, under the isomorphisms $H^k(Y;R)\cong H^k\bigl(\cF^*(Y),\d\bigr)$, $H^k_\cs(Y;R)\cong H^k\bigl(\cF_\cs^*(Y),\d\bigr)$, the cup products in \eq{kh2eq45} are given by
\e
[\al]\cup[\be]=[\psi_{k,l}(Y)(\al\ot\be)]
\label{kh2eq65}
\e
for $\al\in\cF^k(Y)$, $\be\in\cF^l(Y)$ with $\d\al=\d\be=0$, so that $\al\ot\be\in(\cF^k\ot_R\cF^l)(Y)$, and $\psi_{k,l}(Y)(\al\ot\be)\in\cF^{k+l}(Y)$ with $\d\bigl[\psi_{k,l}(Y)(\al\ot\be)\bigr]=0$ by \eq{kh2eq64}.

The procedure for computing the other products in \eq{kh2eq45}, \eq{kh2eq46}, \eq{kh2eq52}, \eq{kh2eq53} at the (co)chain level is the same. Here is an example:

\begin{ex} Let $Y$ be a manifold, and consider the de Rham cohomology $H^*_\dR(Y;\R)$ from Example \ref{kh2ex7}, and the compactly-supported de Rham cohomology $H^*_{\cs,\dR}(Y;\R)$ from Example \ref{kh2ex10}, which were explained using sheaves in Example \ref{kh2ex14}. Writing $\Om^k_Y$ for the sheaf of smooth sections of $\La^kT^*Y$, as in Example \ref{kh2ex14}, define morphisms $\psi_{k,l}:\Om^k_Y\ot_\R\Om^l_Y\ra \Om^{k+l}_Y$ by $\psi_{k,l}(U):\al\ot\be\mapsto \al\w\be$, where $U\subseteq Y$ is open, $\al\in\Om^k_Y(U)$, $\be\in\Om^l_Y(U)$ are $k$- and $l$-forms on $U$, and $\al\w\be$ is the wedge product of exterior forms on $U$.

Since $k$- and $l$-forms $\al,\be$ satisfy $\d(\al\!\w\!\be)\!=\!(\d\al)\!\w\!\be\!+\!(-1)^k\al\!\w\!(\d\be)$, equation \eq{kh2eq64} holds, so these morphisms $\psi_{k,l}$ allow us to compute the cup product on $H^*_\dR(Y;\R), H^*_{\cs,\dR}(Y;\R)$ at the cochain level. That is, we have $[\al]\cup[\be]\!=\![\al\w\be]$, where $\w:C_\dR^k(Y;\R)\!\t\! C_\dR^l(Y;\R)\!\ra\! C_\dR^{k+l}(Y;\R)$ or $\w:C_{\cs,\dR}^k(Y;\R)\!\t\! C_{\cs,\dR}^l(Y;\R)\!\ra\! C_{\cs,\dR}^{k+l}(Y;\R)$ is the wedge product of exterior forms on~$Y$.

Note that $\cup=\w:C_\dR^k(Y;\R)\t C_\dR^l(Y;\R)\ra C_\dR^{k+l}(Y;\R)$ satisfies all of identities \eq{kh2eq47}, \eq{kh2eq48}, \eq{kh2eq50}, \eq{kh2eq55}, \eq{kh2eq56} at the cochain level. 
\label{kh2ex18}
\end{ex}

\subsubsection{Cup products in singular cohomology}
\label{kh263}

Following \cite[\S VI.4]{Bred1}, \cite[\S 3.2]{Hatc}, \cite[\S 48]{Munk}, we define the cup product on singular cohomology $H^*_\rsi(Y;R)$ and smooth singular cohomology~$H^*_\ssi(Y;R)$.

\begin{ex} Let $Y$ be a topological space and $R$ a commutative ring, and use the notation of Example \ref{kh2ex5}. For all $k,l\ge 0$, define an $R$-bilinear map $\cup:C^k_\rsi(Y;R)\t C^l_\rsi(Y;R)\ra C^{k+l}_\rsi(Y;R)=\Hom_\Z\bigl(C_{k+l}^\rsi(Y;\Z),R\bigr)$ by
\e
\al\cup\be:\ts\sum_{i\in I}\rho_i\,\si_i\longmapsto \sum_{i\in I}\rho_i\cdot\al(\si_i\ci G_k^{k+l})\cdot \be(\si_i\ci H_l^{k+l}).
\label{kh2eq66}
\e

Here $\sum_{i\in I}\rho_i\,\si_i\in C_{k+l}^\rsi(Y;\Z)$ so that $I$ is a finite indexing set, $\rho_i\in\Z$, and $\si_i:\De_{k+l}\ra Y$ is continuous for $i\in I$. Also $G_k^{k+l}:\De_k\ra\De_{k+l}$ and $H_l^{k+l}:\De_l\ra\De_{k+l}$ are defined by $G_k^{k+l}:(x_0,\ldots,x_k)\mapsto(x_0,\ldots,x_k,0,\ldots,0)$ and $H_l^{k+l}:(x_0,\ldots,x_l)\mapsto(0,\ldots,0,x_0,\ldots,x_l)$, for $\De_k,\De_l,\De_{k+l}$ as in \eq{kh2eq5}. Thus $\si_i\ci G_k^{k+l}:\De_k\ra Y$ is continuous, so that $\si_i\ci G_k^{k+l}\in C_k^\rsi(Y;\Z)$, and $\al\in C^k_\rsi(Y;R)=\Hom_\Z\bigl(C_k^\rsi(Y;\Z),R\bigr)$, so that $\al(\si_i\ci G_k^{k+l})\in R$. Similarly $\be(\si_i\ci H_l^{k+l})\in R$, so the r.h.s.\ of \eq{kh2eq66} is a finite sum in $R$, and $\al\cup\be$ lies in $\Hom_\Z\bigl(C_{k+l}^\rsi(Y;\Z),R\bigr)= C^{k+l}_\rsi(Y;R)$ as required.

Define $1_Y\in C^0_\rsi(Y;R)$ by $1_Y:\sum_{i\in I}\rho_i\,\si_i\mapsto\sum_{i\in I}\pi(\rho_i)$, where $\pi:\Z\ra R$ is the natural ring morphism. Then $\d\, 1_Y=0$, so $[1_Y]\in H^0_\rsi(Y;R)$.

As in \cite[\S VI.4]{Bred1}, \cite[\S 3.2]{Hatc}, \cite[\S 48]{Munk}, these $\cup,1_Y$ satisfy identities \eq{kh2eq48}, \eq{kh2eq50}, \eq{kh2eq55}, and \eq{kh2eq56} at the cochain level. They also satisfy
\e
\d(\al\cup\be)=(\d\al)\cup\be+(-1)^k\al\cup(\d\be)
\label{kh2eq67}
\e
for all $\al\in C^k_\rsi(Y;R)$ and $\be\in C^l_\rsi(Y;R)$. This implies that $\cup$ descends to cohomology, giving $\cup:H^k_\rsi(Y;R)\t H^l_\rsi(Y;R)\ra H^{k+l}_\rsi(Y;R)$.

However, $\cup$ {\it does not satisfy \eq{kh2eq47} on cochains\/} $C^*_\rsi(Y;R)$, even though it does satisfy \eq{kh2eq47} on cohomology $H^*_\rsi(Y;R)$. That is, $\bigl(C^*_\rsi(Y;R),\d,\cup,1_Y\bigr)$ is a unital differential graded $R$-algebra, but it may not be supercommutative. When $Y$ is a manifold, the same construction defines the cup product on smooth singular cochains $C^*_\ssi(Y;R)$ and cohomology $H^*_\ssi(Y;R)$ in Example~\ref{kh2ex5}.
\label{kh2ex19}
\end{ex}

\begin{rem}{\bf(i)} As in Bredon \cite[\S VI.4]{Bred1}, there are actually many different ways to define a functorial cup product $\cup$ on $C^*_\rsi(Y;R)$, depending on a choice of `diagonal approximation'. Equation \eq{kh2eq66} is just one possibility, coming from the `Alexander--Whitney diagonal approximation'. It has the advantages of being a simple formula, and being associative at the cochain level.
\smallskip 

\noindent{\bf(ii)} In Example \ref{kh2ex18} the cup product $\cup=\w$ is supercommutative on cochains $\bigl(C^*_\dR(Y;\R),\d\bigr)$. However, in Example \ref{kh2ex19}, for singular cohomology over general $R$, the cup product $\cup$ is not supercommutative on cochains~$\bigl(C^*_\rsi(Y;R),\d\bigr)$.

For some $R$, such as $\Z$ and $\Z_2$, {\it it is not possible\/} to define a cochain complex $\bigl(C^*(Y;R),\d\bigr)$ computing $H^*(Y;R)$ which has a supercommutative cup product $\cup$ defined on cochains. This because of {\it Steenrod squares\/} \cite[\S VI.15--\S VI.16]{Bred1}, which are cohomology operations defined using the failure of $\cup$ to be supercommutative at the cochain level, and are nontrivial when $R=\Z$ or~$\Z_2$.
\label{kh2rem4}
\end{rem}

\subsubsection{Axiomatic characterization of cup products}
\label{kh264}

Since we characterized homology and cohomology as $R$-modules axiomatically in Axioms \ref{kh2ax1} and \ref{kh2ax2}, it is natural to ask whether we can also characterize the products $\cup,\cap,\t$ uniquely by some axioms. Kreck and Singhof \cite[\S 6]{KrSi} prove:

\begin{prop} Let\/ $R$ be $\Z$ or $\Z_n,$ and suppose that for all manifolds $Y$ we are given $R$-bilinear maps $\bu:H^k(Y;R)\t H^l(Y;R)\ra H^{k+l}(Y;R)$ such that
\begin{itemize}
\setlength{\itemsep}{0pt}
\setlength{\parsep}{0pt}
\item[{\bf(i)}] $\al\bu(\be\bu\ga)=(\al\bu\be)\bu\ga$ and\/ $1_Y\bu \al=\al\bu 1_Y=\al$ for all\/ $\al,\be,\ga\in H^*(Y;R);$
\item[{\bf(ii)}] If\/ $f:Y_1\ra Y_2$ is a smooth map of manifolds then $f^*(\al\bu\be)=f^*(\al)\bu f^*(\be)$ for all\/ $\al,\be\in H^*(Y_2;R);$ and
\item[{\bf(iii)}] If\/ $\al\in H^m(\cS^m;R)$ and\/ $\be\in H^n(\cS^n;R)$ for $m,n\ge 0$ then $\pi_{\cS^m}^*(\al)\bu\pi_{\cS^n}^*(\be)=\pi_{\cS^m}^*(\al)\cup\pi_{\cS^n}^*(\be)$ in $H^{m+n}(\cS^m\t \cS^n;R)$.
\end{itemize}
Then $\bu$ equals the usual cup product.
\label{kh2prop2}
\end{prop}

If $R$ is any $\Q$-algebra such as $\Q,\R$ or $\C$ then $H^*(Y;R)\cong H^*(Y,\Z)\ot_\Z\nobreak R$. Thus, $\cup$ on $H^*(Y;\Z)$ also determines $\cup$ on $H^*(Y;R)$ for manifolds $Y$. The cross product $\t$ on cohomology is then given in terms of $\cup$ by \eq{kh2eq54}. For manifolds (at least in the compact oriented case), Poincar\'e duality in \S\ref{kh28} identifies $\cap$ with $\cup$, and the cross product on homology with the cross product on cohomology. So the cup product on manifolds determines all of~$\cup,\cap,\t$.

See also Massey \cite[\S 7.8]{Mass2}, who proves that cup and cross products $\cup,\t$ on compactly-supported cohomology $H^*_\cs(-;R)$ of locally-compact topological spaces, for general $R$, are characterized uniquely by some lists of axioms.

\subsection{Differential graded algebras and homotopy theory}
\label{kh27}

We saw in \S\ref{kh21}--\S\ref{kh26} that given a topological space $Y$ (e.g.\ a manifold), the cohomology $H^*(Y;R)$ with the cup product $\cup$ and identity $[1_Y]\in H^0(Y;R)$ is a supercommutative graded $R$-algebra. Usually $H^*(Y;R)$ is the cohomology of a cochain complex $\bigl(C^*(Y;R),\d\bigr)$. Under good conditions $\cup$ on $H^*(Y;R)$ is induced by an associative graded product $\cup$ on $C^*(Y;R)$, compatible with $\d$ as in \eq{kh2eq67}, with an identity $1_Y\in C^0(Y;R)$ inducing $[1_Y]$ in $H^0(Y;R)$. For instance, this works for singular cohomology with the Alexander--Whitney cup product, and for de Rham cohomology of manifolds.

As in Remark \ref{kh2rem4}(ii), whether $\cup$ is supercommutative on cochains $C^*(Y;R)$ depends on the ring $R$ and the cohomology theory. For de Rham cohomology $\cup$ is supercommutative in Example \ref{kh2ex18}, but for singular cohomology $\cup$ is not supercommutative in Example \ref{kh2ex19}, and for some rings such as $R=\Z,\Z_2$ it is not possible to make $\cup$ supercommutative on $C^*(Y;R)$ in any cohomology theory. Thus, $\bigl(C^*(Y;R),\d,\cup,1_Y\bigr)$ is a {\it differential graded\/ $R$-algebra} ({\it dga\/}), which may not be supercommutative, though it is supercommutative on cohomology.

There is subtle but important homotopy-invariant information about the topological space $Y$ which is remembered by the dga $\bigl(C^*(Y;R),\ab\d,\ab\cup,\ab 1_Y\bigr)$, but forgotten by the graded $R$-algebra $\bigl(H^*(Y;R),\cup,1_Y\bigr)$. Some of this information can be used to define additional algebraic structures on cohomology, such as Steenrod squares \cite[\S VI.15]{Bred1}, \cite[\S 4.L]{Hatc} and Massey products \cite{Mass1}, which can be used to distinguish topological spaces with isomorphic cohomology.

However, a more modern approach in homotopy theory, as in F\'elix, Halperin and Thomas \cite{FHT1} for instance, is to work with the dga $\bigl(C^*(Y;R),\ab\d,\ab\cup,\ab 1_Y\bigr)$ up to equivalence in a suitable $\iy$-category $\dga_R^\iy$ of dgas, rather than working with cohomology $H^*(Y;R)$ plus additional structures up to isomorphism. When $R$ is a $\Q$-algebra, one can choose $\bigl(C^*(Y;R),\ab\d,\ab\cup,\ab 1_Y\bigr)$ to be supercommutative, that is, a {\it commutative differential graded\/ $R$-algebra} ({\it cdga\/}), regarded as an object up to equivalence in a suitable $\iy$-category $\cdga_R^\iy$ of cdgas.

For example, in the subject of {\it rational homotopy theory\/} \cite[\S 19]{BoTu}, \cite{FHT2}, to any simply-connected topological space $Y$ one functorially associates another topological space $Y_\Q$ called the {\it rationalization\/} of $Y$, informally by killing all torsion in the homotopy groups $\pi_*(Y)$. One then studies $Y_\Q$ up to homotopy. It turns out that this is equivalent to studying certain cdgas $A^\bu$ over $\Q$ associated to $Y$ up to equivalence in $\cdga_\Q^\iy$. Such cdgas $A^\bu$ include cochain models $\bigl(C^*(Y;\Q),\ab\d,\ab\cup,\ab 1_Y\bigr)$ for $H^*(Y;\Q)$ in suitable cohomology theories. One popular model is Sullivan's cdga $A_{PL}(Y;\Q)$ of {\it polynomial differential forms on $Y$ with coefficients in\/} $\Q$, as in F\'elix, Halperin and Thomas \cite[\S 10(c)]{FHT2} and Sullivan~\cite{Sull}.

For our new cohomology theories of manifolds $MH^*(Y;R),MH^*_\Q(Y;R),\ldots$ introduced in \S\ref{kh4}--\S\ref{kh5}, it may in future be useful to know not only that there are canonical isomorphisms $MH^*(Y;\R)\cong H^*(Y;R)$, but also that the cochain dgas or cdgas lie in the expected equivalence classes of (c)dgas, so that we can use them to compute Steenrod squares, Massey products, rational homotopy types, and so on. To prove this, it is enough to show that, as in \S\ref{kh262}, the cochain (c)dgas $\bigl(MC^*(Y;R),\ab\d,\ab\cup,\ab 1_Y\bigr),\ldots$ of our cohomology theories are constructed using sheaf cohomology from a soft resolution $\cF^\bu$ of $R_Y$ as in \eq{kh2eq63}, with identity $1_Y=i(1)\in\cF^0(Y)$, and cup product $\cup$ defined using sheaf morphisms $\psi_{k,l}:\cF^k\ot_R\cF^l\ra\cF^{k+l}$ satisfying \eq{kh2eq64} and associativity
\e
\psi_{j+k,l}\!\ci\!(\psi_{j,k}\!\ot\!\id_{\cF^l})\!=\!\psi_{j,k+l}\!\ci\!(\id_{\cF^j}\!\ot\!\psi_{k,l}):\cF^j\!\ot_R\!\cF^k\!\ot_R\!\cF^l\!\longra\!\cF^{j+k+l}.
\label{kh2eq68}
\e
For example, all this holds for the sheaf cohomology presentation of de Rham cohomology in Example~\ref{kh2ex14}.

\subsection{(Co)homology of manifolds, and Poincar\'e duality}
\label{kh28}

Almost all the theory of \S\ref{kh21}--\S\ref{kh27} works for general topological spaces, although we often explained it only for manifolds. But the material of this section works only for manifolds, or for topological spaces which are very like manifolds.

Let $Y$ be an oriented manifold, of dimension $m$, not necessarily compact, and $R$ a commutative ring. Then there exists a natural {\it fundamental class\/} $[[Y]]$ in $H_m^\lf(Y;R)$, in locally finite homology. To define $[[Y]]$ in locally finite (smooth) singular homology $H_m^{\lf,\rsi}(Y;R)$ or $H_m^{\lf,\ssi}(Y;R)$ from Examples \ref{kh2ex11} and \ref{kh2ex12}, we choose a (smooth) locally finite triangulation of $Y$ by $m$-simplices $\De_m$.

Define {\it Poincar\'e duality morphisms\/} using the cap products in \eq{kh2eq46}:
\e
\begin{gathered}
\Pd:H^k_\cs(Y;R)\longra H_{m-k}(Y;R),\quad
\Pd:H^k(Y;R)\longra H_{m-k}^\lf(Y;R),\\
\text{by}\qquad \Pd:\al\longmapsto \al\cap[[Y]].
\end{gathered}
\label{kh2eq69}
\e
The {\it Poincar\'e duality theorem\/} says that these are isomorphisms, as in \cite[\S VI]{Bred1}, \cite[\S 3.3]{Hatc}, \cite[Prop.~5.10(i) \& Rem.~15.7(ii)]{HuRa}, \cite[\S 11.3]{Mass2}, \cite[Ch.~8]{Munk}, and~\cite[\S 18]{toDi}.

If $Y$ is compact then $H_*^\lf(Y;R)=H_*(Y;R)$ and $H^*_\cs(Y;R)=H^*(Y;R)$, so $[[Y]]\in H_m(Y;R)$, and we have isomorphisms $\Pd:H^k(Y;R)\ra H_{m-k}(Y;R)$. Poincar\'e duality is often stated just in this case, or for compact manifolds with boundary, since few authors discuss locally finite homology~$H_*^\lf(Y;R)$.

Poincar\'e duality and the cup product on $H^*_\cs(Y;R)$ induce an associative, supercommutative {\it intersection product\/} $\bu:H_k(Y;R)\t H_l(Y;R)\ra H_{k+l-m}(Y;R)$ on homology $H_*(Y;R)$. This can be described geometrically at the chain level in singular homology, as in~\cite[\S VI.11]{Bred1} and~\cite[\S IV]{Lefs}.

Poincar\'e duality implies that $H_k(Y;R)\!=\!H^k(Y;R)\!=\!H^k_\cs(Y;R)\!=\!\ab H_k^\lf(Y;R)\ab=\!0$ for all $k>m=\dim Y$. This also holds for non-oriented manifolds~$Y$.

Let $Y,Z$ be oriented manifolds of dimensions $m,n$, and $f:Y\ra Z$ be a smooth map. Define morphisms $f^!:H^*_\cs(Y;R)\ra H^{*-m+n}_\cs(Z;R)$ and $f_!:H_*^\lf(Z;R)\ra H_{*-n+m}^\lf(Y;R)$ by the commutative diagrams
\begin{equation*}
\xymatrix@!0@C=21pt@R=30pt{
*+[r]{H^k_\cs(Y;R)} \ar[rrrrrrr]_(0.45)\Pd \ar[d]^{f^!} &&&&&&& *+[l]{H_{m-k}(Y;R)} \ar[d]_{f_*} & *+[r]{H_k^\lf(Z;R)} \ar[rrrrrrr]_(0.45){\Pd^{-1}} \ar[d]^{f_!} &&&&&&& *+[l]{H^{n-k}(Z;R)} \ar[d]_{f^*} \\
*+[r]{H^{k-m+n}_\cs(Z;R)} &&&&&&& *+[l]{H_{m-k}(Z;R),\!\!} \ar[lllllll]_(0.42){\Pd^{-1}}  & *+[r]{H_{k+m-n}^\lf(Y;R)} &&&&&&& *+[l]{H^{n-k}(Y;R),\!\!} \ar[lllllll]_(0.42)\Pd }\!\!\!{}
\end{equation*}
using the fact that $\Pd$ in \eq{kh2eq69} are invertible. We have $f_!([[Z]])=[[Y]]$ in $H_m(Y;R)$. Similarly, if $f$ is proper, define $f^!:H^*(Y;R)\ra H^{*-m+n}(Z;R)$ and $f_!:H_*(Z;R)\ra H_{*-n+m}(Y;R)$ by the commutative diagrams
\e
\begin{gathered}
\xymatrix@!0@C=19.5pt@R=30pt{
*+[r]{H^k(Y;R)} \ar[rrrrrrr]_(0.45)\Pd \ar[d]^{f^!} &&&&&&& *+[l]{H_{m-k}^\lf(Y;R)} \ar[d]_{f_*} & *+[r]{H_k(Z;R)} \ar[rrrrrrr]_(0.45){\Pd^{-1}} \ar[d]^{f_!} &&&&&&& *+[l]{H^{n-k}_\cs(Z;R)} \ar[d]_{f^*} \\
*+[r]{H^{k-m+n}(Z;R)} &&&&&&& *+[l]{H_{m-k}^\lf(Z;R),\!\!} \ar[lllllll]_(0.42){\Pd^{-1}}  & *+[r]{H_{k+m-n}(Y;R)} &&&&&&& *+[l]{H^{n-k}_\cs(Y;R).\!\!} \ar[lllllll]_(0.42)\Pd }\!\!\!{}
\end{gathered}
\label{kh2eq70}
\e
These morphisms $f^!,f_!$ are known as {\it wrong way maps\/} (or {\it shriek maps}, or {\it umkehr maps}, or {\it Gysin maps\/}) as they have the opposite functoriality that one expects, covariant on cohomology and contravariant on homology.

In fact it is not necessary for $Y,Z$ to be oriented to define $f^!,f_!$, we only need a coorientation (relative orientation) for $f:Y\ra Z$, as in Definition~\ref{kh3def6}.

Here are some compatibility conditions between morphisms of type $f^!,f_!$ and $f^*,f_*$. Suppose we have a Cartesian square in $\Man$:
\begin{equation*}
\xymatrix@C=100pt@R=15pt{ *+[r]{W} \ar[r]_f \ar[d]^e & *+[l]{Y} \ar[d]_h \\ 
*+[r]{X} \ar[r]^g & *+[l]{Z,\!} }
\end{equation*}
with $g,h$ transverse, and set $\dim W=k$, $\dim X=l$, $\dim Y=m$, $\dim Z=n$, so that $k=l+m-n$ by transversality. Suppose $h$ is cooriented, which implies that $e$ is cooriented. Then the following commute:
\e
\begin{gathered}
\xymatrix@!0@C=19.5pt@R=30pt{
*+[r]{H^j_\cs(W;R)}  \ar[d]^{e^!} &&&&&&& *+[l]{H^j_\cs(Y;R)} \ar[lllllll]^(0.45){f^*} \ar[d]_{h^!} & *+[r]{H_{j+k-l}^\lf(W;R)} \ar[rrrrrrr]_(0.48){f_*}  &&&&&&& *+[l]{H_{j+m-n}^\lf(Y;R)}  \\
*+[r]{H^{j-k+l}_\cs(X;R)} &&&&&&& *+[l]{H^{j-m+n}_\cs(Z;R),\!\!} \ar[lllllll]_(0.48){g^*}  & *+[r]{H_j^\lf(X;R)}  \ar[rrrrrrr]^(0.5){g_*} \ar[u]_{e_!} &&&&&&& *+[l]{H_j^\lf(Z;R).\!\!} \ar[u]^{h_!} }\!\!\!{}
\end{gathered}
\label{kh2eq71}
\e
If also $h$ is proper, which implies that $e$ is proper, then the analogues of \eq{kh2eq71} for $H^*(-;R),H_*(-;R)$ also commute. 

\subsection{Introduction to orbifolds, and their (co)homology}
\label{kh29}

Orbifolds are generalizations of manifolds, which are locally modelled on $\R^m/G$, for $G$ a finite group acting linearly on $\R^m$. They were introduced by Satake \cite{Sata}, who called them `V-manifolds'. Later they were studied by Thurston \cite[Ch.~13]{Thur} who gave them the name `orbifold'. An orbifold $X$ is called {\it effective\/} if it is locally modelled on $\R^m/G$ for $G$ acting effectively on $\R^m$, that is, the morphism $G\ra\GL(m,\R)$ is injective, so we can regard $G$ as a finite subgroup of $\GL(m,\R)$. Many authors, including Satake \cite{Sata}, consider only effective orbifolds.

It turns out that Satake's original definition \cite{Sata} of a category of orbifolds $\OrbSa$ is in some ways badly behaved differential-geometrically. For example, the pullback $f^*(E)$ of an orbifold vector bundle $E\ra Y$ by a morphism $f:X\ra Y$ in $\OrbSa$ cannot always be defined. Behrend and Xu \cite{BeXu}, Chen and Ruan \cite{ChRu1}, Lerman \cite{Lerm}, Metzler \cite{Metz}, Moerdijk and Pronk \cite{Moer,MoPr,Pron}, the author \cite{Joyc4,Joyc7}, and others have given alternative definitions of categories or 2-categories of orbifolds, not equivalent to $\OrbSa$, with better differential-geometric behaviour. We discuss these other definitions in Remark~\ref{kh2rem5}(c).

In \S\ref{kh4} and \S\ref{kh51}--\S\ref{kh52} we construct new (co)homology theories $MH_*(Y;R),\ab \ldots,MH^*_\dR(Y;\R)$ for manifolds $Y$. Sections \ref{kh35} and \ref{kh53} will explain how to extend all of this to effective orbifolds $Y$. To do this, we have chosen to work in a (rather crude) category of effective orbifolds $\Orbeff$, defined in \S\ref{kh291}, which is very similar to Satake's $\OrbSa$, rather than any of the more sophisticated versions of \cite{BeXu,ChRu1,Joyc4,Joyc7,Lerm,Metz,Moer,MoPr,Pron}. There are three main reasons for this:
\begin{itemize}
\setlength{\itemsep}{0pt}
\setlength{\parsep}{0pt}
\item[(a)] The differential-geometric properties of manifolds $\Man$ we need for our constructions also work in $\Orbeff$, including in particular submersions and transverse fibre products, as in~\S\ref{kh292}.

In general, transverse fibre products of orbifolds are subtle, and only work properly when orbifolds are defined as a 2-category. But we avoid all these problems by adopting a more restrictive definition of `submersions' and `transverse smooth maps' in $\Orbeff$, which are strong enough to ensure that transverse fibre products exist in $\Orbeff$ in the sense of category theory, but still weak enough to be used in our theory.
\item[(b)] As in \S\ref{kh293}, any effective orbifold $X$ is homotopic in $\Orbeff$ to a manifold $Y$, which is false for more sophisticated definitions. This makes it easy to extend results on (co)homology of manifolds, such as characterization by Eilenberg--Steenrod axioms in \S\ref{kh21}--\S\ref{kh22}, to orbifolds.
\item[(c)] If $\Orb$ is any of the more modern (2-)categories of (not necessarily effective) orbifolds in the literature, then there is a truncation functor $F:\Orb\ra\Orbeff$, since the other definitions factor through $\Orbeff$ by forgetting additional structure. Thus, our results on (co)homology of orbifolds in $\Orbeff$ pull back immediately to~$\Orb$.
\end{itemize}

We begin in \S\ref{kh291} by defining a category of effective orbifolds $\Orbeff$, and discussing alternative definitions. Section \ref{kh292} explains some differential geometry in $\Orbeff$, including submersions and transverse fibre products. Then \S\ref{kh293} extends \S\ref{kh21}--\S\ref{kh28} on (co)homology of manifolds to orbifolds. 

\subsubsection{A category of effective orbifolds $\Orbeff$}
\label{kh291}

We define a category $\Orbeff$ of effective orbifolds.

\begin{dfn} Let $X$ be a second countable Hausdorff topological space. An {\it $m$-dimensional orbifold chart on\/} $X$ is a triple $(U,G,\phi)$, where $U\subseteq\R^m$ is an open set, and $G\subseteq\GL(m,\R)$ is a finite subgroup preserving $U$, and $\phi:U\ra X$ is a continuous map with image $\Im\phi\subseteq X$ an open set in $X$, such that $\phi\ci\ga=\phi$ for all $\ga\in G$, so that $\phi$ factors through a map $U/G\ra\Im\phi$, which we require to be a homeomorphism.

Let $(U,G,\phi),(V,H,\psi)$ be $m$-dimensional orbifold charts on $X$. We call $(U,G,\phi)$ and $(V,H,\psi)$ {\it compatible\/} if for all $u\in \phi^{-1}(\Im\psi)\subseteq U$, there exists an open neighbourhood $U'$ of $u$ in $U$ and an \'etale map of manifolds $\xi:U'\ra V$ with~$\phi\vert_{U'}=\psi\ci\xi:U'\ra X$.

An $m$-{\it dimensional orbifold atlas\/} for $X$ is a set $\cA=\{(U_a,G_a,\phi_a):a\in A\}$ of pairwise compatible $m$-dimensional orbifold charts on $X$ with $X=\bigcup_{a\in A}\Im\phi_a$. We call such an atlas {\it maximal\/} if it is not a proper subset of any other orbifold atlas. Any orbifold atlas $\{(U_a,G_a,\phi_a):a\in A\}$ is contained in a unique maximal orbifold atlas, the set of all orbifold charts $(U,G,\phi)$ on $X$ which are compatible with $(U_a,G_a,\phi_a)$ for all~$a\in A$.

An $m$-{\it dimensional effective orbifold\/} $(X,\cA)$ is a second
countable Hausdorff topological space $X$ equipped with a maximal
$m$-dimensional orbifold atlas $\cA$. Usually we refer to $X$ as the orbifold, leaving the atlas implicit, and by a {\it chart\/ $(U,G,\phi)$ on\/} $X$, we mean an element of the maximal atlas.

Now let $X,Y$ be effective orbifolds, and $f:X\ra Y$ be a continuous map. We call $f$ {\it smooth\/} if for all $x\in X$ with $f(x)=y\in Y$, there exist charts $(U,G,\phi)$ on $X$ with $x\in\Im\phi$ and $(V,H,\psi)$ on $Y$ with $y\in\Im\psi$, a smooth map of manifolds $F:U\ra V$, and a group morphism $\rho:G\ra H$, such that $F$ is equivariant under $\rho$, and $f\ci\phi=\psi\ci F:U\ra Y$. Smooth maps include identities $\id_X:X\ra X$ and are closed under composition, so they make effective orbifolds into a category, which we write as~$\Orbeff$.

There is an obvious full and faithful functor $F_\Man^\Orbeff:\Man\ra\Orbeff$ acting on objects by $F_\Man^\Orbeff:(X,\cA)\mapsto(X,\cA')$ and on morphisms by $F_\Man^\Orbeff:f\mapsto f$, where $\cA=\{(U_a,\phi_a):a\in A\}$ is the maximal atlas on $X$ in the usual sense of manifolds, and $\cA'$ is the unique maximal orbifold atlas on $X$ containing $\{(U_a,\{1\},\phi_a):a\in A\}$. Thus we can regard $\Man\subset\Orbeff$ as a full subcategory, and manifolds as examples of effective orbifolds.
\label{kh2def12}
\end{dfn}

\begin{rem}{\bf(a) (Noneffective orbifolds.)} Effectiveness of orbifolds $X$ is built into Definition \ref{kh2def12} by taking orbifold charts to be $(U,G,\phi)$ with $U\subseteq\R^m$ open and $G\subseteq\GL(m,\R)$, rather than taking $G$ to be a finite group acting possibly non-effectively on~$U\subseteq\R^m$.

In other definitions of orbifolds $X$ allowing $X$ to be non-effective, there is a natural morphism $\pi:X\ra X^\eff$, for $X^\eff$ an effective orbifold with the same topological space as $X$, where if $X$ is locally modelled on $\R^m/G$ for $G$ a finite group and $\rho:G\ra\GL(m,\R)$ a morphism, then $X^\eff$ is locally modelled on $\R^m/\rho(G)$. Since $X,X^\eff$ have isomorphic homology and cohomology, in studying (co)homology of orbifolds we lose little by restricting to effective orbifolds.
\smallskip 

\noindent{\bf(b) (Comparison with Satake's definition.)} Definition \ref{kh2def12} differs from Satake's category \cite{Sata} of effective orbifolds $\OrbSa$ in two ways:
\begin{itemize}
\setlength{\itemsep}{0pt}
\setlength{\parsep}{0pt}
\item[(i)] In defining orbifold charts $(U,G,\phi)$, Satake requires $G$ to have fixed points in codimension at least 2, that is, Satake's orbifolds cannot have orbifold strata of codimension 1, but we allow this. An example of an orbifold in $\Orbeff$ but not in $\OrbSa$ is $\R/\{\pm 1\}$. Topologically this is $[0,\iy)$, so orbifolds with codimension 1 orbifold strata are like orbifolds with corners.
\item[(ii)] In defining smooth maps $f:X\ra Y$, Satake requires that for each $x\in X$ with $f(x)=y\in Y$ there should exist charts $(U,G,\phi)$ on $X$ with $x\in\Im\phi$ and $(V,H,\psi)$ with $y\in\Im\psi$ and a smooth map of manifolds $F:U\ra V$ with $f\ci\phi=\psi\ci F:U\ra Y$, as we do. But he does not require, as we do, that $F$ should be equivariant under a group morphism~$\rho:G\ra H$.

There exist pathological examples of morphisms $f:X\ra Y$ in $\OrbSa$ which are not morphisms in $\Orbeff$, though they are not easy to construct. We exclude them as they would cause problems with the existence of transverse fibre products in \S\ref{kh292}. All the (2-)categories of orbifolds in \cite{ALR,BeXu,ChRu1,Joyc4,Joyc7,Lerm,Metz,Moer,MoPr,Pron} satisfy our extra equivariance condition.
\end{itemize}

\noindent{\bf(c) (Other definitions of orbifolds.)} Other than Satake's \cite{Sata}, there are two main definitions of ordinary categories of orbifolds in the literature:
\begin{itemize}
\setlength{\itemsep}{0pt}
\setlength{\parsep}{0pt}
\item[(A)] Chen and Ruan \cite[\S 4]{ChRu1} define orbifolds $X$ in a similar way to \cite{Sata}, but using germs of orbifold charts $(U_p,G_p,\phi_p)$ for $p\in X$. Their morphisms $f:X\ra Y$ are called `good maps', giving a category $\Orb_{\rm CR}$. 
\item[(B)] Moerdijk and Pronk \cite{Moer,MoPr} define a category of orbifolds $\Orb_{\rm MP}$ as {\it proper \'etale Lie groupoids\/} in $\Man$. Their definition of smooth map $f:X\ra Y$, called {\it strong maps\/} \cite[\S 5]{MoPr} is complicated: it is an equivalence class of diagrams $\smash{X\,{\buildrel\phi\over \longleftarrow}\,X'\,{\buildrel\psi\over\longra}\,Y}$, where $X'$ is a third orbifold, and $\phi,\psi$ are morphisms of groupoids with $\phi$ an equivalence (loosely, a diffeomorphism).
\end{itemize}
A book on orbifolds in the sense of \cite{ChRu1,Moer,MoPr} is Adem, Leida and Ruan~\cite{ALR}.

There are also four main ways to define 2-categories of orbifolds: 
\begin{itemize}
\setlength{\itemsep}{0pt}
\setlength{\parsep}{0pt}
\item[(i)] Pronk \cite{Pron} defines a strict 2-category $\bf LieGpd$ of Lie groupoids in $\Man$ as in (B), with the obvious 1-morphisms of groupoids, and localizes by a class of weak equivalences $\cW$ to get a weak 2-category~$\mathfrak{Orb}_{\rm Pr}={\bf LieGpd}[\cW^{-1}]$.

Lerman \cite[\S 3.3]{Lerm} defines a weak 2-category $\mathfrak{Orb}_{\rm Le}$ of Lie groupoids in $\Man$ as in (B), without localizing, using `Hilsum--Skandalis morphisms'.
\item[(ii)] Behrend and Xu \cite[\S 2]{BeXu}, Lerman \cite[\S 4]{Lerm} and Metzler \cite[\S 3.5]{Metz} define a strict 2-category of orbifolds $\mathfrak{Orb}_{\rm ManSta}$ as a class of Deligne--Mumford stacks on the site $(\Man,{\cal J}_\Man)$ of manifolds.
\item[(iii)] The author \cite{Joyc4} defines a strict 2-category of orbifolds $\mathfrak{Orb}_{C^\iy{\rm Sta}}$ as a class of Deligne--Mumford stacks on the site $(\CSch,{\cal J}_\CSch)$ of $C^\iy$-schemes.
\item[(iv)] The author \cite[\S 4.5]{Joyc7} defines a weak 2-category of orbifolds $\mathfrak{Orb}_{\rm Kur}$ as examples of Kuranishi spaces.
\end{itemize}
It is known (loc.~cit.) that (i)--(iv) are equivalent as weak 2-categories, and their homotopy categories are equivalent to~(B).

All of these (2-)categories of orbifolds admit natural truncation functors to $\Orbeff$ in \S\ref{kh291}. Furthermore, isomorphism/equivalence classes of effective orbifolds in each of these (2-)categories are naturally in 1-1 correspondence with isomorphism classes of objects in $\Orbeff$.

Morally speaking, effective orbifolds in each of these (2-)categories are the same, just as objects, and non-effective orbifolds $X$ can be transformed to effective orbifolds $X^\eff$ as in {\bf(a)}. Also, morally, (1-)morphisms in each of these (2-)categories are continuous maps $f:X\ra Y$ of the topological spaces, plus some extra data. The truncation functor to $\Orbeff$ forgets this extra data, remembering only the continuous map.
\smallskip

\noindent{\bf(d) (Orbifolds with boundary and corners.)} As for manifolds with boundary and corners in \S\ref{kh31}, we can modify Definition \ref{kh2def12} to define categories $\Orbeffb$ of {\it effective orbifolds with boundary\/} and $\Orbeffc$ of {\it effective orbifolds with corners}, by also allowing orbifold charts $(U,G,\phi)$ with $U\subseteq [0,\iy)\t\R^{m-1}$ open (for the boundary case), or $U\subseteq [0,\iy)^k\t\R^{m-k}$ open (for the corners case).

To define smooth maps $f:X\ra Y$ in $\Orbeffb,\Orbeffc$, the map $F:U\ra V$ in Definition \ref{kh2def12} should be required to be a smooth map between $U\subseteq[0,\iy)^k\t\R^{m-k}$ and $V\subseteq[0,\iy)^l\t\R^{n-l}$, in the sense of Definition \ref{kh3def1}(b) below.
\label{kh2rem5}
\end{rem}

\subsubsection{Differential geometry in $\Orbeff$}
\label{kh292}

In the category of manifolds $\Man$, a morphism $f:X\ra Y$ is a {\it submersion\/} if $T_xf:T_xX\ra T_yY$ is surjective for all $x\in X$ with $f(x)=y\in Y$, and morphisms $g:X\ra Z$, $h:Y\ra Z$ are {\it transverse\/} if $T_xg\op T_yh:T_xX\op T_yY\ra T_zZ$ is surjective for all $x\in X$, $y\in Y$ with $g(x)=h(y)=z\in Z$. Then $g:X\ra Y$ a submersion implies $g,h$ are transverse for any $h:Y\ra Z$, and $g,h$ transverse implies that a fibre product $X\t_{g,Z,h}Y$ exists in the category~$\Man$.

In our new (co)homology theories for manifolds in \S\ref{kh4} and \S\ref{kh51}--\S\ref{kh52}, we will make heavy use of these facts, and their extension to manifolds with corners. So to extend \S\ref{kh4}--\S\ref{kh52} to orbifolds in \S\ref{kh53}, we need good notions of submersions and transverse morphisms in $\Orbeff$, such that transverse fibre products exist.

The definitions of submersion and transverse morphisms above extend immediately to orbifolds. But with this natural definition, the existence of transverse fibre products turns out to be a subtle question. One can prove that some transverse fibre products of orbifolds do not exist in {\it any\/} ordinary category of orbifolds $\Orb$. Here we mean fibre products in the sense of category theory, satisfying a universal property in $\Orb$. But in the 2-categories of orbifolds $\mathfrak{Orb}$ in Definition \ref{kh2rem5}(c)(i)--(iv), all transverse fibre products exist, in the sense of 2-category theory, characterized by a universal property in $\mathfrak{Orb}$ involving 2-morphisms. This is an important reason for making orbifolds into a 2-category.

Working with a 2-category of orbifolds $\mathfrak{Orb}$ rather than $\Orbeff$ would cause us other problems, for example, Theorem \ref{kh2thm6} below would fail. So we are going to cheat, and adopt stricter definitions of `submersions' and `transverse morphisms' in $\Orbeff$ imposing conditions on orbifold groups as well as on tangent spaces, and with these definitions, transverse fibre products exist in~$\Orbeff$.

The cost of this is that there are fewer submersions and transverse morphisms, for example,  if $V$ is a manifold and $G$ a finite group acting effectively but not freely on $V$ then the natural projection $\pi:V\ra[V/G]$ is not a submersion in our sense. Also for a general orbifold $X$, the projection $\pi:TX\ra X$ is not a submersion. But all the properties of submersions and transverse maps we need in $\Man$, given in Assumption \ref{kh3ass5} below, also hold in~$\Orbeff$.

We first explain orbifold groups and tangent spaces of effective orbifolds.

\begin{dfn} Let $X\in\Orbeff$ be an effective orbifold, and $x\in X$. Then one can define  the {\it orbifold group\/} (or {\it isotropy group\/}) $G_xX$, a finite group, and the {\it tangent space\/} $T_xX$ of $X$ at $x$, a real vector space of dimension $\dim X$ with an effective representation of $G_xX$. Explicitly, if $(U,G,\phi)$ is an orbifold chart on $X$ and $u\in U$ with $\phi(u)=x$ then
\e
G_xX\cong\Stab_G(u)=\bigl\{\ga\in G:\ga\cdot u=u\bigr\}\quad\text{and}\quad T_xX\cong T_uU.
\label{kh2eq72}
\e

In fact $G_xX,T_xX$ above are natural up to isomorphism, but not up to canonical isomorphism; the pair $(G_xX,T_xX)$ is natural up to conjugation by an element of $G_xX$, and to define $G_xX,T_xX$ unambiguously for all $X,x$ we have to make arbitrary choices, e.g.\ choosing $u\in\phi^{-1}(x)\subseteq U$ so \eq{kh2eq72} makes sense.

Now let $f:X\ra Y$ be a morphism in $\Orbeff$, and $x\in X$ with $f(x)=y\in Y$. Then we can define a linear map $T_xf:T_xX\ra T_yY$ called the {\it derivative of\/ $f$ at\/} $x$. Explicitly, as in Definition \ref{kh2def12} there exist charts $(U,G,\phi)$ on $X$ with $x\in\Im\phi$ and $(V,H,\psi)$ on $Y$ with $y\in\Im\psi$, a smooth map of manifolds $F:U\ra V$, and a group morphism $\rho:G\ra H$, such that $F$ is equivariant under $\rho$, and $f\ci\phi=\psi\ci F:U\ra Y$. Then choosing $u\in U$ with $\phi(u)=x$ and $v=F(u)\in V$ with $\psi(v)=y$, we define $T_xf$ by the commutative diagram
\begin{equation*}
\xymatrix@C=100pt@R=15pt{ *+[r]{T_uU} \ar[d]^\cong \ar[r]_{T_uF} & *+[l]{T_vV} \ar[d]_\cong \\
*+[r]{T_xX} \ar[r]^{T_xf} & *+[l]{T_yY,\!} }
\end{equation*}
where the vertical isomorphisms come from~\eq{kh2eq72}.

In fact $T_xf$ above is unique only up to composition with the action of some $\ga\in G_yY$ on $T_yY$, since $F$ is not unique, and to define $T_xf$ unambiguously for all $f,X,Y,x,y$ we have to make arbitrary choices.

Also there exists a group morphism $G_xf:G_xX\ra G_yY$ such that $T_xf:T_xX\ra T_yY$ is equivariant under $G_xf$, which in the situation above may be defined by the commutative diagram
\begin{equation*}
\xymatrix@C=100pt@R=15pt{ *+[r]{\Stab_G(u)} \ar[d]^\cong \ar[r]_{\rho\vert_{\Stab_G(u)}} & *+[l]{\Stab_H(v)} \ar[d]_\cong \\
*+[r]{G_xX} \ar[r]^{G_xf} & *+[l]{G_yY,\!} }
\end{equation*}
where the vertical isomorphisms come from \eq{kh2eq72}. However, {\it we stress that $G_xf$ may not be unique, not even up to conjugation by some\/} $\ga\in G_yY$. This is because in the situation above with $x,y,(U,G,\phi),(V,H,\psi)$ and $F:U\ra V$ fixed, the group morphism $\rho:G\ra H$ may not be uniquely determined, if $F$ maps into the fixed locus of a nontrivial subgroup of $H$ in~$V$.

For the (2-)categories of orbifolds discussed in Remark \ref{kh2rem5}(c), the morphisms $G_xf:G_xX\ra G_yY$ are canonical up to conjugation by $\ga\in G_yY$, as this is part of the data in (1-)morphisms beyond the continuous map. But in $\Orbeff$ this is not true, which will be important in the proof of Theorem~\ref{kh2thm6}.
\label{kh2def13}
\end{dfn}

We define `submersions' and `transverse morphisms' in~$\Orbeff$.

\begin{dfn} We call a morphism $f:X\ra Y$ in $\Orbeff$ a {\it submersion\/} if whenever $x\in X$ with $f(x)=y\in Y$, then
\begin{itemize}
\setlength{\itemsep}{0pt}
\setlength{\parsep}{0pt}
\item[(i)] $T_xf:T_xX\ra T_yY$ is surjective; and
\item[(ii)] The morphism $G_xf:G_xX\ra G_yY$ in Definition \ref{kh2def13} is surjective.
\end{itemize}
As above, $T_xf,G_xf$ depend on arbitrary choices, but one can show that (i),(ii) together are independent of these choices.

Let $g:X\ra Z$ and $h:Y\ra Z$ be morphisms in $\Orbeff$. We call $g,h$ {\it transverse\/} if whenever $x\in X$, $y\in Y$ with $g(x)=h(y)=z\in Z$, then
\begin{itemize}
\setlength{\itemsep}{0pt}
\setlength{\parsep}{0pt}
\item[(i$)'$] $T_xg\op T_yh:T_xX\op T_yY\ra T_zZ$ is surjective; and
\item[(ii$)'$] $G_zZ=G_xg(G_xX)\cdot G_yh(G_yY)$, that is, every $\ep\in G_zZ$ may be written as $\ep=G_x(\ga)G_yh(\de)$ for some $\ga\in G_xX$ and $\de\in G_yY$.
\end{itemize}
Again, $T_xg,T_yh,G_xg,G_yh$ depend on arbitrary choices, but one can show that (i$)'$,(ii$)'$ together are independent of these choices. If $g:X\ra Z$ and $h:Y\ra Z$ are morphisms in $\Orbeff$ with $g$ a submersion then $g,h$ are transverse.
\label{kh2def14}
\end{dfn}

It is then not difficult to prove:

\begin{thm} Let\/ $g:X\ra Z$ and\/ $h:Y\ra Z$ be transverse morphisms in $\Orbeff$. Then a fibre product\/ $W=X\t_{g,Z,h}Y$ exists in $\Orbeff,$ in the sense of category theory. It has $\dim W=\dim X+\dim Y-\dim Z,$ and topological space
\e
W\cong\bigl\{(x,y)\in X\t Y:g(x)=h(y)\bigr\},
\label{kh2eq73}
\e
as a subset of\/ $X\t Y$ with the subspace topology.
\label{kh2thm5}
\end{thm}

\begin{rem} For comparison, if we work in a 2-category of orbifolds $\mathfrak{Orb}$ as in Remark \ref{kh2rem5}(d)(i)--(iv), then the picture above is modified as follows:
\begin{itemize}
\setlength{\itemsep}{0pt}
\setlength{\parsep}{0pt}
\item We should define a 1-morphism $f:\cX\ra\cY$ in $\mathfrak{Orb}$ to be a {\it submersion\/} if $T_xf:T_x\cX\ra T_y\cY$ is surjective for all $x\in\cX$ with $f(x)=y\in\cY$.
\item We should define 1-morphisms $g:\cX\ra\cZ$ and $h:\cY\ra\cZ$ in $\mathfrak{Orb}$ to be {\it transverse\/} if for all $x\in\cX$ and $y\in\cY$ with $g(x)=h(y)=z\in\cZ$ and all $\ep\in G_z\cZ$, the following linear map is surjective:
\begin{equation*}
T_xg\op (\ep\cdot T_yh):T_x\cX\op T_y\cY\longra T_z\cZ.
\end{equation*}
\item If $g,h$ are transverse then a fibre product $\cW=\cX\t_{g,\cZ,h}\cY$ exists in the 2-category $\mathfrak{Orb}$, satisfying a universal property involving 2-morphisms, with $\dim\cW=\dim\cX+\dim\cY-\dim\cZ$. The set of points of $\cW$ is
\e
\begin{split}
\cW\cong\bigl\{(x,y,C):x\in\cX,\;\>y\in\cY,\;\> g(x)=h(y)=z\in\cZ,&\\ 
C\in G_xg(G_x\cX)\backslash G_z\cZ/G_yh(G_y\cY)&\bigr\},
\end{split}
\label{kh2eq74}
\e
where $G_xg:G_x\cX\ra G_z\cZ$, $G_yh:G_y\cY\ra G_z\cZ$ are the natural morphisms of orbifold groups, which are defined in $\mathfrak{Orb}$ but not in~$\Orbeff$.
\end{itemize}

So, in the 2-category $\mathfrak{Orb}$ we can use weaker notions of `submersion' and `transverse', corresponding roughly to assuming (i),(i$)'$ but not (ii),(ii$)'$. We need the richer 2-category structure on $\mathfrak{Orb}$ to define transverse fibre products.

Now in Definition \ref{kh2def14} and Theorem \ref{kh2thm5} we consider transverse fibre products in the (rather crude) ordinary category of orbifolds $\Orbeff$, rather than the (much superior) 2-category $\mathfrak{Orb}$. We also want the fibre products in $\Orbeff$ and $\mathfrak{Orb}$ to coincide, since we regard the fibre product in $\mathfrak{Orb}$ as `correct'.

If a fibre product $W=X\t_{g,Z,h}Y$ exists in $\Orbeff$, applying the universal property to maps from a point $*$ shows that $W$ must be given as a set by \eq{kh2eq73}. So, the fibre products in $\Orbeff$ and $\mathfrak{Orb}$ can coincide only if \eq{kh2eq73} and \eq{kh2eq74} agree, which holds if $G_z\cZ=G_xg(G_x\cX)\cdot G_yh(G_y\cY)$ for all $x\in\cX$, $y\in\cY$ with $g(x)=h(y)=z\in\cZ$. This is the reason for Definition~\ref{kh2def14}(ii$)'$.
\label{kh2rem6}
\end{rem}

The next theorem, proved in \S\ref{kh63}, will be important in extending results on (co)homology of manifolds to orbifolds in~\S\ref{kh293}.

\begin{thm} Every effective orbifold\/ $X$ is homotopic in $\Orbeff$ to a manifold\/ $Y$ in $\Man\subset\Orbeff$. That is, there exist smooth maps $f:X\ra Y,$ $g:Y\ra X,$ $F:[0,1]\t X\ra X$ and\/ $G:[0,1]\t Y\ra Y$ with\/ $F(0,x)=g\ci f(x),$ $F(1,x)=x,$ $G(0,y)=f\ci g(y)$ and\/ $G(1,y)=y$ for all\/ $x\in X$ and\/~$y\in Y$.

Here we call\/ $F:[0,1]\t X\ra X$ \begin{bfseries}smooth\end{bfseries} if there exists an open neighbourhood\/ $U$ of\/ $[0,1]\t X$ in $\R\t X$ and a smooth map (morphism in $\Orbeff$) $F':U\ra X$ with $F'\vert_{[0,1]\t X}=F,$ and similarly for~$G:[0,1]\t Y\ra Y$.
\label{kh2thm6}
\end{thm}

\begin{rem} In Theorem \ref{kh2thm6}, it is {\it essential\/} that we are working in the rather crude category of orbifolds $\Orbeff$ from \S\ref{kh291}, in which morphisms $f:X\ra Y$ are continuous maps satisfying conditions, rather than any of the more complicated (2-)categories of orbifolds in \cite{ALR,BeXu,ChRu1,Joyc4,Joyc7,Lerm,Metz,Moer,MoPr,Pron} discussed in Remark \ref{kh2rem5}(c), in which (1-)morphisms are continuous maps with extra data.

In any of these other (2-)categories, the last part of the proof in \S\ref{kh63} would fail, and we would not be able to choose the (1-)morphism $F:[0,1]\t X\ra X$ with $F\vert_{\{0\}\t X}=g\ci f$ and $F\vert_{\{1\}\t X}=\id_X$. As in Definition \ref{kh2def13}, if $f:X\ra Y$ is a (1-)morphism of orbifolds, and $x\in X$ with $f(x)=y\in Y$, then there exists a morphism of orbifold groups $G_xf:G_xX\ra G_yY$. In our category $\Orbeff$, these $G_xf$ are {\it not canonical}, and there can be several very different choices. In the orbifolds of \cite{ALR,BeXu,ChRu1,Joyc4,Joyc7,Lerm,Metz,Moer,MoPr,Pron}, the morphisms $G_xf:G_xX\ra G_yY$ are canonical up to conjugation by $\ga\in G_yY$, and the $G_xf$ are part of the extra data in (1-)morphisms $f:X\ra Y$ in these (2-)categories.

These morphisms $G_xf$ have some continuity properties, which imply that for a (1-)morphism $F:[0,1]\t X\ra X$ in any of these more sophisticated \hbox{(2-)categories} of orbifolds, the morphism $G_{(t,x)}F:G_xX\ra G_{F(t,x)}X$ must have kernel $\Ker G_{(t,x)}F$ independent of $t\in[0,1]$, for fixed $x\in X$. But when $t=0$ as $F\vert_{\{0\}\t X}=g\ci f$ we have $G_{(0,x)}F=G_{f(x)}g\ci G_xf=1$, since $G_{f(x)}Y=\{1\}$, and when $t=1$ as $F\vert_{\{1\}\t X}=\id_X$ we have~$G_{(1,x)}F=\id_{G_xX}$. 

Hence $\Ker G_{(0,x)}F=G_xX$ and $\Ker G_{(1,x)}F=\{1\}$, so $\Ker G_{(t,x)}F$ cannot be independent of $t\in[0,1]$ unless $G_xX=\{1\}$. Thus, Theorem \ref{kh2thm6} fails in all the (2-)categories of \cite{ALR,BeXu,ChRu1,Joyc4,Joyc7,Lerm,Metz,Moer,MoPr,Pron}, for any orbifold $X$ which is not a manifold. The argument works in $\Orbeff$ because of the nonuniqueness of morphisms $G_xf:G_xX\ra G_yY$ there.
\label{kh2rem7}
\end{rem}

\subsubsection{Homology and cohomology of effective orbifolds}
\label{kh293}

Next we explain how to explain the material of \S\ref{kh21}--\S\ref{kh28} on (co)homology of manifolds, to orbifolds. First we should ask what we mean by the homology or cohomology of an orbifold $Y$. There are several possible answers:
\begin{itemize}
\setlength{\itemsep}{0pt}
\setlength{\parsep}{0pt}
\item[(a)] The most obvious is the (co)homology $H_*(Y_\top;R), H^*(Y_\top;R)$ (e.g. singular (co)homology) of the underlying topological space $Y_\top$ of $Y$.
\item[(b)] As in Moerdijk \cite[\S 4]{Moer}, to each orbifold $Y$ one can associate a {\it classifying space\/} $Y_\cla$, which is a topological space with a continuous map $\pi:Y_\cla\ra Y_\top$ whose fibre $\pi^{-1}(y)$ over each $y\in Y_\top$ is homotopic to $B(G_yY)$, the classifying space for the orbifold group $G_yY$. Here $Y_\cla$ is canonical only up to homotopy, but this means $H_*(Y_\cla;R),H^*(Y_\cla;R)$ are natural up to isomorphism, so we can consider these to be the (co)homology of~$Y$.

Behrend \cite{Behr} defines versions of singular (co)homology of orbifolds $Y$ considered as proper \'etale Lie groupoids as in \cite{Moer,MoPr,Pron}, which compute $H_*(Y_\cla;R),H^*(Y_\cla;R)$. When $Y$ is a global quotient $[V/G]$ we can take $Y_\cla=(V\t EG)/G$, and then $H_*(Y_\cla;R),H^*(Y_\cla;R)$ are the equivariant (co)homology $H_*^G(V;R),\ab H^*_G(V;R)$, a desirable property.
\item[(c)] Takeuchi and Yokoyama \cite{TaYo1,TaYo2,TaYo3,Yoko} {\it t-singular homology\/} $t\text{-}H_*(Y;R)$ and {\it ws-singular cohomology\/} $ws\text{-}H^*(Y;R)$ of orbifolds $Y$, with the nice properties that $t\text{-}H_1(Y;\Z)$ is the abelianization of the orbifold fundamental group $\pi_1^{\rm orb}(Y)$, and Poincar\'e duality holds. 
\item[(d)] If $Y$ is an orbifold with an almost complex structure, Chen and Ruan \cite{ChRu2} define the {\it orbifold cohomology\/} $H^*_{\sst\rm CR}(Y;R)$, which appears to be natural in Gromov--Witten theory and String Theory of orbifolds. In general $H^*_{\sst\rm CR}(Y;R)$ is different from $H^*(Y_\top;R),H^*(Y_\cla;R)$ and~$ws\text{-}H^*(Y;R)$.

Orbifold cohomology is related to the cohomology $H^*(Y_{\rm in};R)$ of the {\it inertia orbifold\/} $Y_{\rm in}$ of $Y$, where points of $Y_{\rm in}$ are pairs $(y,C)$ with $y\in Y_\top$ and $C$ a conjugacy class in~$G_yY$.
\end{itemize}

We are interested in option (a), and from now on $H_*(Y;R),H^*(Y;R),\ldots$ will mean $H_*(Y_\top;R),H^*(Y_\top;R),\ldots.$ But if $R$ is a $\Q$-algebra, then (b),(c) agree with (a) as in \cite[Prop.~36]{Behr}, \cite{TaYo3,Yoko}. We generalize the approach to (co)homology of manifolds using Eilenberg--Steenrod axioms in \S\ref{kh21}--\S\ref{kh22} to orbifolds.

\begin{dfn} Fix a commutative ring $R$. A {\it homology theory of effective orbifolds $H_*(-;R)$ over\/} $R$ satisfies Axiom \ref{kh2ax1}, but with the category of effective orbifolds $\Orbeff$ in \S\ref{kh291} substituted for the category of manifolds $\Man$ throughout. That is, we are given the data: 
\begin{itemize}
\setlength{\itemsep}{0pt}
\setlength{\parsep}{0pt}
\item[(a)] For each effective orbifold $Y$ and open suborbifold $Z\subseteq Y$, an $R$-module $H_k(Y,Z;R)$ for all $k\in\Z$, where we write $H_k(Y;R)=H_k(Y,\es;R)$; 
\item[(b)] For $Y,Z$ as above, $R$-module morphisms $\pd:H_k(Y,Z;R)\ra H_{k-1}(Z;R)$ for $k\in\Z$, called {\it connecting morphisms\/}; and
\item[(c)] For each morphism $f:Y_1\ra Y_2$ in $\Orbeff$ and open suborbifolds $Z_1\subseteq Y_1$, $Z_2\subseteq Y_2$ with $f(Z_1)\subseteq Z_2$, {\it pushforwards\/} $f_*:H_k(Y_1,Z_1;R)\ra H_k(Y_2,Z_2;R)$ for $k\in\Z$;
\end{itemize}
and all this data should satisfy Axiom \ref{kh2ax1}(i)--(vii) with $\Orbeff$ in place of $\Man$. 

Similarly, a {\it cohomology theory of effective orbifolds $H^*(-;R)$ over\/} $R$ satisfies Axiom \ref{kh2ax2}, but with $\Orbeff$ substituted for $\Man$ throughout.
\label{kh2def15}
\end{dfn}

Here is the orbifold analogue of Theorems \ref{kh2thm1} and \ref{kh2thm2}, which will be proved in \S\ref{kh64} using Theorem~\ref{kh2thm6}.

\begin{thm}{\bf(a)} Any two homology theories of effective orbifolds $H_*(-;R),\ab \ti H_*(-;R)$ over the same commutative ring\/ $R$ are canonically isomorphic. That is, there exist\/ $R$-module isomorphisms $I_{Y,Z}:H_*(Y,Z;R)\ra\ti H_*(Y,Z;R)$ for all effective orbifolds $Y$ and open $Z\subseteq Y$ commuting with the given morphisms\/ $f_*,\pd$ and isomorphisms\/ $H_0(*;R)\cong R\cong\ti H_0(*;R),$ and any other assignment of morphisms\/ $J_{Y,Z}:H_*(Y,Z;R)\ra\ti H_*(Y,Z;R)$ for all\/ $(Y,Z)$ commuting with the\/ $f_*,\pd$ and\/ $H_0(*;R)\cong R\cong\ti H_0(*;R)$ have\/ $J_{Y,Z}=I_{Y,Z}$ for all\/~$(Y,Z)$. 
\smallskip

\noindent{\bf(b)} Any two cohomology theories of effective orbifolds $H^*(-;R),\ab \ti H^*(-;R)$ over the same commutative ring\/ $R$ are canonically isomorphic.
\label{kh2thm7}
\end{thm}

\begin{ex}{\bf(i)} Singular homology $H_*^\rsi(-;R)$ is a homology theory of effective orbifolds, as in Example \ref{kh2ex1}, as it is a topological homology theory.
\smallskip

\noindent{\bf(ii)} Smooth singular homology $H_*^\ssi(-;R)$ in Example \ref{kh2ex2} also extends to orbifolds without change, giving a homology theory of effective orbifolds.
\smallskip

\noindent{\bf(iii)} If $Y$ is an effective orbifold, then one can define $k$-forms $C^\iy(\La^kT^*Y)$ on $Y$ as for manifolds, and they have the same properties under de Rham differentials, pullbacks, and so on. Concretely, if $Y=[V/G]$ is the quotient of a smooth manifold $V$ by a finite group $G$, with projection $\pi:V\ra [V/G]=Y$, then $\pi^*:C^\iy(\La^kT^*Y)\ra C^\iy(\La^kT^*V)$ is an isomorphism $C^\iy(\La^kT^*Y)\cong C^\iy(\La^kT^*V)^G$ with the vector subspace of $G$-invariant $k$-forms on $V$. 

We can define de Rham cohomology $H^*_\dR(-;\R)$ for effective orbifolds as in Example \ref{kh2ex7}, and it is a cohomology theory of effective orbifolds.
\smallskip

\noindent{\bf(iv)} Sheaf cohomology $H^*(Y;R_Y)$ from \S\ref{kh25} is a cohomology theory of orbifolds,  as it is a cohomology theory of (nice) topological spaces.
\label{kh2ex20}
\end{ex}

Almost all the rest of \S\ref{kh21}--\S\ref{kh27} is facts about (co)homology of general topological spaces, and so extends to (co)homology of effective orbifolds immediately. The only other parts of \S\ref{kh21}--\S\ref{kh28} which need attention are the special facts about (co)homology of manifolds, rather than topological spaces, namely:
\begin{itemize}
\setlength{\itemsep}{0pt}
\setlength{\parsep}{0pt}
\item[(a)] If $Y$ is an oriented manifold of dimension $m$, as in Remark \ref{kh2rem1}(f) and Property \ref{kh2pr2}(i) there is a natural {\it fundamental class\/} $[[Y]]$ in $H_m(Y;R)$ (for $Y$ compact) or $H_m^\lf(Y;R)$. If $R$ has characteristic 2 (e.g.\ $R=\Z_2$) then $[[Y]]$ is well-defined even if $Y$ is not oriented.
\item[(b)] As in \S\ref{kh25}, for nice topological spaces $Y$ and noetherian $R$ we have a {\it dualizing complex\/} $\om_Y$ over $R$ in the derived category $D^+(Y;R)$ of $R$-modules on $Y$, and as in \eq{kh2eq28} we have natural isomorphisms
\e
H_k(Y;R)\cong \H^{-k}_\cs(Y,\om_Y)\quad\text{and}\quad H_k^\lf(Y;R)\cong \H^{-k}(Y,\om_Y).
\label{kh2eq75}
\e

If $Y$ is a manifold of dimension $m$, then as in Definition \ref{kh2def9} we have an equivalence $\om_Y\simeq O_Y[m]$ in $D^+(Y;R)$, where $O_Y$ is the orientation sheaf of $Y$ over $R$, a locally constant sheaf of $R$-modules associated to the principal $\Z_2$-bundle $P_Y\ra Y$ of orientations on~$Y$.  

Thus as in \eq{kh2eq29} we have isomorphisms
\e
H_k(Y;R)\cong H^{m-k}_\cs(Y,O_Y)\quad\text{and}\quad H_k^\lf(Y;R)\cong H^{m-k}(Y,O_Y).
\label{kh2eq76}
\e
\item[(c)] All of \S\ref{kh28} on Poincar\'e duality and wrong way maps $f^!,f_!$.\end{itemize}

The next three properties explain how (a)--(c) extend to orbifolds.

\begin{property}{\bf(a)} If $Y$ is an oriented effective orbifold, then $Y$ has a natural {\it fundamental class\/} $[[Y]]$ in $H_m(Y;R)$ (for $Y$ compact) or $H_m^\lf(Y;R)$, for any commutative ring $R$, just as for manifolds.
\smallskip

\noindent{\bf(b)} When considering (co)homological properties of unoriented orbifolds, it is helpful to divide into three cases:
\begin{itemize}
\setlength{\itemsep}{0pt}
\setlength{\parsep}{0pt}
\item[(i)] Orbifolds $Y$ which are {\it locally orientable}, that is, which are locally modelled on $\R^m/G$ for $G\subset\GL_+(m,\R)$ orientation-preserving.
\item[(ii)] Orbifolds $Y$ which are not locally orientable, but which do not have orbifold singularities in real codimension 1.
\item[(iii)] Orbifolds $Y$ which have orbifold singularities in real codimension 1 (this implies $Y$ is not locally orientable). Satake \cite{Sata} excludes this case.
\end{itemize}
Any oriented orbifold is automatically locally orientable. If $R$ has characteristic 2 (e.g.\ $R=\Z_2$) then there is a natural fundamental class $[[Y]]$ in $H_m(Y;R)$ (for $Y$ compact) or $H_m^\lf(Y;R)$ in cases (i),(ii), but not in case~(iii).
\smallskip

\noindent{\bf(c)} Suppose $V$ is an oriented manifold of dimension $m$ and $G$ a finite group acting effectively on $Y$ preserving orientations, and set $Y=[V/G]$, so that $Y$ is an oriented effective orbifold. Write $\pi:V\ra [V/G]=Y$ for the projection, which is proper. Let $R$ be a $\Q$-algebra. Then the fundamental class $[[Y]]\in H_m^\lf(Y;R)$ of $Y$, as an orbifold, is given in terms of the fundamental class $[[V]]\in H_m^\lf(V;R)$ of $V$, as a manifold, by
\begin{equation*}
[[Y]]=\frac{1}{\md{G}}\cdot\pi_*([[V]]).
\end{equation*}

\noindent{\bf(d) (Noneffective orbifolds.)} One can define noneffective orbifolds $Y$, although they are not included as objects in our category $\Orbeff$. As in Remark \ref{kh2rem5}(a), every noneffective orbifold $Y$ has an associated effective orbifold $Y^\eff$, with a projection $\pi:Y\ra Y^\eff$ which is a homeomorphism of topological spaces. Hence $H_*(Y;R)\cong H_*(Y^\eff;R)$, $H^*(Y;R)\cong H^*(Y^\eff;R)$, so for the purposes of homology and cohomology one can often restrict to effective orbifolds.

One exception to this is fundamental classes. Suppose $Y$ is a connected, oriented, noneffective orbifold of dimension $m$. Then $\pi:Y\ra Y^\eff$ is a `gerbe' with fibre $[*/G]$ for some finite group $G$. We should restrict to $R$ a $\Q$-algebra, and the natural fundamental class for $Y$ is then
\begin{equation*}
[[Y]]=\frac{1}{\md{G}}\cdot [[Y^\eff]]\in H_m^\lf(Y^\eff;R)\cong H_m^\lf(Y;R).
\end{equation*}

\label{kh2pr3}
\end{property}

\begin{property}{\bf(a)} Property \ref{kh2pr3} defined {\it locally orientable\/} effective orbifolds $Y$. If $Y$ is locally orientable then it has a principal $\Z_2$-bundle $P_Y\ra Y$ of orientations, and the orientation sheaf $O_Y$ exists as the locally constant sheaf of $R$-modules on $Y$ associated to $P_Y$, for any ring $R$.

If $Y$ is not locally orientable then the principal $\Z_2$-bundle $P_Y\ra Y$ is not defined. The sheaf $O_Y$ is still defined, but it is not a locally constant.
\smallskip

\noindent{\bf(b)} For any effective orbifold $Y$ and noetherian ring $R$, the dualizing complex $\om_Y$ is well defined, and equation \eq{kh2eq75} holds.
\smallskip

\noindent{\bf(c)} If $Y$ is an effective orbifold of dimension $m$ and $R$ is a $\Q$-algebra, then the equivalence $\om_Y\simeq O_Y[m]$ holds in $D^+(Y;R)$ as for manifolds, and thus equation \eq{kh2eq76} holds. If $R$ is not a $\Q$-algebra then in general we have $\om_Y\not\simeq O_Y[m]$, and \eq{kh2eq76} may be false.
\label{kh2pr4}
\end{property}

\begin{property}{\bf(a)} Let $Y$ be an oriented, effective orbifold of dimension $m$. Then as in \eq{kh2eq69} and with $[[Y]]$ as in Property \ref{kh2pr3}(a), for any commutative ring $R$ we may define {\it Poincar\'e duality morphisms\/}
\begin{gather}
\Pd:H^k_\cs(Y;R)\ra H_{m-k}(Y;R),\quad
\Pd:H^k(Y;R)\ra H_{m-k}^\lf(Y;R),
\nonumber\\
\text{by}\qquad \Pd:\al\longmapsto \al\cap[[Y]].
\label{kh2eq77}
\end{gather}

If $R$ is a $\Q$-algebra, these morphisms $\Pd$ are isomorphisms (that is, the Poincar\'e duality theorem holds), as for manifolds. But if $R$ is not a $\Q$-algebra, examples show that they may not be isomorphisms.
\smallskip

\noindent{\bf(b)} The wrong way maps $f^!,f_!$ on (co)homology of manifolds were defined in \S\ref{kh28} using the inverse $\Pd^{-1}$ of $\Pd$ in \eq{kh2eq77}. So by {\bf(a)}, if $R$ is a $\Q$-algebra then we may define wrong way maps $f^!,f_!$ for (proper) cooriented morphisms $f:Y\ra Z$ in $\Orbeff$, and they have the same properties as for manifolds. But if $R$ is not a $\Q$-algebra, in general we cannot define $f^!,f_!$ for orbifolds.
\label{kh2pr5}
\end{property}

One moral of Properties \ref{kh2pr3}--\ref{kh2pr5} is that (co)homology of effective orbifolds $Y$ over a $\Q$-algebra $R$ behaves exactly like (co)homology of manifolds over $R$. But when $R$ is not a $\Q$-algebra, some nice properties of (co)homology of manifolds do not extend to orbifolds.

We discuss different kinds of singular homology for orbifolds:

\begin{ex} Let $Y$ be an effective orbifold. Then we can define:
\begin{itemize}
\setlength{\itemsep}{0pt}
\setlength{\parsep}{0pt}
\item[(i)] The {\it singular homology\/} $H_*^\rsi(Y;R)$, as in Example~\ref{kh2ex1}.
\item[(ii)] The {\it smooth singular homology\/} $H_*^\ssi(Y;R)$, as in Example~\ref{kh2ex2}.
\item[(iii)] The {\it cosheaf singular homology\/} $\hat H_k^\rsi(Y;R)$, as in Example~\ref{kh2ex15}. 
\item[(iv)] The {\it cosheaf smooth singular homology\/} $\hat H_k^\ssi(Y;R)$, as in Example~\ref{kh2ex16}. 
\item[(v)] The {\it locally finite sheaf singular homology\/} $\hat H_k^{\lf,\rsi}(Y;R)$, as in Example~\ref{kh2ex17}. 
\item[(vi)] The {\it locally finite sheaf smooth singular homology\/} $\hat H_k^{\lf,\ssi}(Y;R)$, as in Example~\ref{kh2ex17}. 
\end{itemize}
Here (i)--(iv) are homology theories of effective orbifolds, with pushforwards $f_*$ by morphisms $f:Y_1\ra Y_2$ in $\Orbeff$, and (v),(vi) are locally finite homology theories of effective orbifolds, with pushforwards $f_*$ for proper~$f$.

Now (iii),(iv) come with complexes $\hat\ucC{}_\bu^\rsi(Y;R),\hat\ucC{}_\bu^\ssi(Y;R)$ of flabby cosheaves of $R$-modules on $Y$, and (v),(vi) with complexes $\hat\cC{}^{\lf,\rsi}_{\bu}(Y;R),\hat\cC{}^{\lf,\ssi}_{\bu}(Y;R)$ of soft sheaves of $R$-modules on $Y$. As in \eq{kh2eq38}--\eq{kh2eq39}, since the proof of \eq{kh2eq38} in \cite[Th.~1.8.6(a)]{Arab} works for orbifolds, we have natural equivalences in~$D(Y;R)$:
\e
\hat\cC{}^{\lf,\rsi}_{-\bu}(Y;R)\simeq \om_Y\simeq \hat\cC{}^{\lf,\ssi}_{-\bu}(Y;R),
\label{kh2eq78}
\e
where $\om_Y$ is the dualizing complex of~$Y$.

Note that \eq{kh2eq78} holds for any effective orbifold $Y$ and noetherian ring $R$. it does not depend on the equivalence $\om_Y\simeq O_Y[m]$, which as in Property \ref{kh2pr4}(c) needs $R$ to be a $\Q$-algebra.
\label{kh2ex21}
\end{ex}

\subsection{\texorpdfstring{Fulton and MacPherson's bivariant theories}{Fulton and MacPherson\textquoteright s bivariant theories}}
\label{kh210}

For manifolds, because of Poincar\'e duality in \S\ref{kh28}, homology and
cohomology can be identified into a single theory with {\it products\/} (cup, cap and intersection products), {\it pushforwards\/} (covariant functoriality), and {\it pullbacks}
(contravariant functoriality). Fulton and MacPherson \cite{FuMa} defined {\it bivariant theories}, a way of naturally integrating pairs of a homology theory and a cohomology theory into a single larger theory, which works for singular spaces.

We now summarize the notion of a bivariant theory, following \cite[\S 1.1]{FuMa}. Let $\cC$ be a category. A bivariant theory from $\cC$ to abelian groups (or graded vector spaces, etc.) associates a group $B(g:Y\ra Z)$ (or graded vector space, etc.) to each morphism $g:Y\ra Z$ in $\cC$. It also has three operations:
\begin{itemize}
\setlength{\itemsep}{0pt}
\setlength{\parsep}{0pt}
\item {\bf Products.} If $f:X\ra Y$ and $g:Y\ra Z$ are morphisms in
$\cC$, and $\al\in B(f:X\ra Y)$, $\be\in B(g:Y\ra Z)$, we can form
the product $\al\cdot\be$ in~$B(g\ci f:X\ra Z)$.
\item {\bf Pushforwards.} For a certain class of morphisms $f:X\ra
Y$ in $\cC$ called {\it confined morphisms}, if $g:Y\ra Z$ is a
morphism in $\cC$ and $\al\in B(g\ci f:X\ra Z)$ we can form the
{\it pushforward\/} $f_*(\al)\in B(g:Y\ra Z)$.
\item {\bf Pullbacks.} For a certain class of commutative squares in
$\cC$ called {\it independent squares}, if
\e
\begin{gathered}
\xymatrix@C=70pt@R=13pt{
*+[r]{Y'} \ar[r]_{g'} \ar[d]^{h'} & *+[l]{Z'} \ar[d]_h \\
*+[r]{Y} \ar[r]^g & *+[l]{Z,\!} }
\end{gathered}
\label{kh2eq79}
\e
is an independent square then for each $\al\in B(g:Y\ra Z)$, we can
form the {\it pullback\/} $h^*(\al)\in B(g':Y'\ra Z')$.
\end{itemize}
All this data must satisfy a number of compatibility axioms.

The following example, taken from \cite[\S 3.1]{FuMa} modified as in
\cite[\S 3.3.2]{FuMa}, combines homology and compactly-supported
cohomology into a bivariant theory.

\begin{ex} Let $\Topem$ be the category of topological spaces $Y$ which admit an embedding as a closed subspace of $\R^n$ for some $n\gg 0$ (in particular, $Y$ must be Hausdorff and paracompact), and morphisms are continuous maps. Let $R$ be a commutative ring. Suppose $g:Y\ra Z$ is a morphism in $\Topem$ with $Y$ compact. Choose an embedding $\phi:Y\ra\R^n$ for some $n\gg 0$, and define $H^k(g:Y\ra Z;R)=H^{k+n}\bigl(Z\t\R^n,(Z\t\R^n)\sm(g\t\phi)(Y);R\bigr)$, using relative cohomology. Then $H^*(g:Y\ra Z;R)$ is a graded $R$-module, which turns out to be independent of the choice of $\phi$ up to canonical isomorphism. For general morphisms in $\Topem$, we define $H^*(g:Y\ra Z;R)$ to be the direct limit $\varinjlim H^*(g\ci p:K\ra Z;R)$ over all $p:K\ra Y$ in $\Topem$ with $K$ compact.

If $Y$ is a topological space in $\Topem$, we can associate two natural morphisms in $\Topem$ to $Y$: the projection to a point $\pi:Y\ra *$, and the identity map $\id_Y:Y\ra Y$. By standard algebraic topology one can show that
\e
\begin{split}
H^k(\pi:Y\ra *;R)&\cong H_{-k}(Y;R)\\
\text{and}\qquad H^k(\id_Y:Y\ra Y;R)&\cong H^k_\cs(Y;R).
\end{split}
\label{kh2eq80}
\e
Thus, the bivariant theory $H^*$ specializes to homology and
compactly-supported cohomology over~$R$.

{\it Products\/} $\cdot:H^k(f:X\ra Y)\t H^l(g:Y\ra Z)\ra H^{k+l}(g\ci f:X\ra Z)$ are defined as in \cite[\S 3.1.7]{FuMa}, where
\begin{gather}
\cdot:H^k(\id_Y:Y\ra Y;R)\t H^l(\id_Y:Y\ra Y;R)\longra H^{k+l}(\id_Y:Y\ra Y;R)
\nonumber\\
\text{is}\qquad
\cup:H^k_\cs(Y;R)\t H^l_\cs(Y;R)\longra H^{k+l}_\cs(Y;R),
\label{kh2eq81}\\
\cdot:H^k(\id_Y:Y\ra Y;R)\t H^l(\pi:Y\ra *;R)\longra H^{k+l}(\pi:Y\ra *;R)
\nonumber\\
\text{is}\qquad \cap:H^k_\cs(Y;R)\t H_{-l}(Y;R)\longra H_{-k-l}(Y;R),
\label{kh2eq82}
\end{gather}
so that products generalize cup and cap products $\cup,\cap$.

All morphisms in $\Topem$ are confined, and {\it pushforwards\/} are defined as in \cite[\S 3.1.8]{FuMa}. For $g:Y\ra Z$ a morphism in $\Topem$, the pushforward
\e
\begin{split}
g_*:H^k(\pi:Y\ra *;R)\longra H^k(\pi:Z\ra *;R)\\
\text{is}\qquad g_*:H_{-k}(Y;R)\longra H_{-k}(Z;R).
\end{split}
\label{kh2eq83}
\e

We define {\it independent squares\/} to be squares \eq{kh2eq79} which are Cartesian in $\Topem$, so that $Y'\cong Y\t_{g,Z,h}Z'$, and with $h$ proper, so that $h'$ is also proper. {\it Pullbacks\/} are defined as in \cite[\S 3.1.6]{FuMa}. If $h:Y'\ra Y$ is a proper morphism then
\e
\begin{split}
h^*:H^k(\id_Y:Y\ra Z;R)\longra H^k(\id_{Y'}:Y'\ra Y';R)\\
\text{is}\qquad h^*:H_\cs^k(Y;R)\ra H_\cs^k(Y';R).
\end{split}
\label{kh2eq84}
\e

\label{kh2ex22}
\end{ex}

Let $\Man$ be the category of manifolds, with morphisms smooth maps. Applying the natural functor $\Man\ra\Topem$ taking a manifold to its underlying topological space, the bivariant theory $H^*$ of Example \ref{kh2ex22} induces a bivariant theory on $\Man$. However, because of Poincar\'e duality the picture simplifies, and we can write $H^*$ on $\Man$ more directly. The following is an equivalent way of defining the bivariant theory of Example \ref{kh2ex22} on the category~$\Man$.

\begin{ex} Let $\Man$ be the category of manifolds, with morphisms smooth maps, and $R$ be a commutative ring. Suppose $g:Y\ra Z$ is a morphism in $\Man$, with $\dim Y=l$ and $\dim
Z=m$. Define
\e
\begin{split}
H^k(g:Y\ra Z;R)&=H^{k+l-m}_\cs(Y,O_Y\ot_R g^*(O_Z))\\
&\cong H_{n-k}(Y;g^*(O_Z)),
\end{split}
\label{kh2eq85}
\e
where as in Definition \ref{kh2def9}, $O_Y,O_Z$ are the {\it orientation sheaves\/} of $Y,Z$, which are locally constant sheaves of $R$-modules on $Y,Z$, so that $O_Y\ot_R g^*(O_Z)$ is a locally constant sheaf on $Y$ and $H^*_\cs(Y,O_Y\ot_R g^*(O_Z))$ is its compactly-supported sheaf cohomology, which we can also interpret as $H_*(Y;g^*(O_Z))$, the homology of $Y$ twisted by the locally constant sheaf~$g^*(O_Z)$.

If $g:Y\ra Z$ is $\pi:Y\ra *$ then $O_Z$ and $g^*(O_Z)$ are
trivial and $n=0$, so $H^k(\pi:Y\ra *;R)\cong H_{-k}(Y;R)$. If
$g$ is $\id_Y:Y\ra Y$ then $g^*(O_Z)\cong O_Y$ and $O_Y\ot_R
g^*(O_Z)=R_Y$, so that $H^k(\id_Y:Y\ra Y;R)=H^k_\cs(Y;R)$. So \eq{kh2eq80} holds.

To define products, suppose $f:X\ra Y$ and $g:Y\ra Z$ are smooth with $\dim X=m$, $\dim Y=n$ and $\dim Z=p$. Define 
\begin{equation*}
\cdot:H^k(f:X\ra Y)\t H^l(g:Y\ra Z)\ra H^{k+l}(g\ci f:X\ra Z)
\end{equation*}
to be the natural product of sheaf cohomology groups
\begin{gather*}
\cdot:H^{k+m-n}_\cs\bigl(X,O_X\ot_R f^*(O_Y)\bigr)\t H^{l+n-p}_\cs\bigl(Y,O_Y\ot_R g^*(O_Z)\bigr)\\
\longra H^{k+l+m-p}_\cs\bigl(X,O_X\ot_R(g\ci f)^*(O_Z)\bigr)\\
\text{mapping}\quad (\al,\be)\longmapsto \al\cdot\be=\al\ot_R f^*(\be),
\end{gather*}
where $f^*(\be)\in H^{l+n-p}_\cs\bigl(X,f^*(O_Y)\ot_R (g\ci f)^*(O_Z)\bigr)$, and we use the isomorphism of sheaves on $X$
\begin{equation*}
[O_X\ot_R f^*(O_Y)]\ot_R[f^*(O_Y)\ot_R (g\ci f)^*(O_Z)]\cong O_X\ot_R(g\ci f)^*(O_Z),
\end{equation*}
since $O_Y\ot_RO_Y\cong R_Y$, as $O_Y$ is $R_Y$ twisted by a principal $\Z_2$-bundle.

All morphisms in $\Man$ are confined. {\it Pushforwards\/}
\begin{gather*}
f_*:H^k(g\ci f:X\ra Z;R)\longra H^k(g:Y\ra Z;R)\\
\text{are pushforwards}\quad f_*:H_{p-k}(X;f^*\ci g^*(O_Z))\longra
H_{p-k}(Y;g^*(O_Z))
\end{gather*}
in twisted homology. {\it Independent squares\/} are defined to be transverse Cartesian squares \eq{kh2eq79} in $\Man$ with $h$ proper, and hence $h'$ proper. Then with $\dim Y=m$, $\dim Z=n$, $\dim Y'=m'$, $\dim Z'=n'$, the Cartesian property gives $(h')^*\bigl(O_Y\ot_Rg^*(O_Z)\bigr)\cong O_{Y'}\ot_R(g')^*(O_{Z'})$ and $m-n=m'-n'$. So we define the {\it pullback\/} $h^*:H^k(g:Y\ra Z;R)\ra H^k(g':Y'\ra Z';R)$ to be
\begin{align*}
(h')^*:\,&H^{k+m-n}_\cs(Y,O_Y\ot_R g^*(O_Z))\longra\\
&H^{k+m-n}_\cs(Y',(h')^*(O_Y\ot_R g^*(O_Z)))\cong\\
&H^{k+m'-n'}_\cs(Y',O_{Y'}\ot_R (g')^*(O_{Z'})),
\end{align*}
using proper pullbacks on compactly-supported sheaf cohomology.
\label{kh2ex23}
\end{ex}

For {\it oriented\/} manifolds $Y,Z$, with $O_Y,O_Z$ trivial, equation \eq{kh2eq85} becomes
\begin{equation*}
H^k(g:Y\ra Z;R)\cong H^{k+m-n}_\cs(Y;R)\cong H_{n-k}(Y;R),
\end{equation*}
and we lose the dependence on $g,Z$, apart from the dimension $n=\dim Z$. Thus, for oriented manifolds and smooth maps, the bivariant theories of Examples \ref{kh2ex22} and \ref{kh2ex23} reduce to homology, compactly-supported cohomology, and Poincar\'e duality. That is, {\it bivariant (co)homology theories really give us nothing new for manifolds}, they are of interest {\it only for singular spaces}.

As in \cite[\S 3.1]{FuMa}, we can also combine cohomology $H^*(-;R)$ and locally finite homology $H_*^\lf(-;R)$ into bivariant theories on $\Topem$ and~$\Man$.

\section{Manifolds with corners}
\label{kh3}

In \S\ref{kh4}--\S\ref{kh5} we will define new (co)homology theories $MH_*(Y;R),MH^*(Y;R),\ldots$ of manifolds with boundaries $Y$, in which the (co)chains involve quadruples $(V,n,s,t)$ with $V$ a `manifold with corners' and $s:V\ra\R^n$, $t:V\ra Y$ smooth. In \S\ref{kh31} we explain the usual definition of a category $\Manc$ of manifolds with corners, based on the author \cite{Joyc3,Joyc6,Joyc8} and Melrose \cite{Melr1,Melr2,Melr3}, and some facts from the theory we will need, principally concerning boundaries $\pd X$, the existence of transverse fibre products $X\t_{g,Z,h}Y$ in $\Manc$ for $Z\in\Man$, and orientations.

In fact our definitions of $MH_*(Y;R),MH^*(Y;R),\ldots$ also work with `manifolds with corners' $V$ replaced by other classes of geometric objects satisfying a list of basic properties -- we can let $V$ have singularities or exotic corners, allow $s:V\ra\R^n$, $t:V\ra Y$ to be only piecewise-smooth, etc. 

There are important applications in forming virtual cycles for moduli spaces of $J$-holomorphic curves $\oM$ in Symplectic Geometry, discussed in \S\ref{kh12} and \cite{Joyc9}, in which it will be an advantage to allow such singularities and non-smoothness, since the `Kuranishi neighbourhoods' $(V,E,\Ga,s,\psi)$ on $\oM$ naturally have singularities and non-smoothness over points in $\oM$ representing nodal curves, and one can only make $V,s$ smooth by doing a lot of extra work.

Because of this, in \S\ref{kh4}--\S\ref{kh5} we will suppose we are given a category of `manifolds with corners' $\tManc$ satisfying a list of assumptions we give in \S\ref{kh33}, which are all we need to develop $MH_*(Y;R),MH^*(Y;R),\ldots.$ We have tried to make these assumptions fairly weak, to make our theory widely applicable for forming virtual cycles in singular situations. 

\subsection{The usual definition of manifolds with corners}
\label{kh31}

We now define our favourite category $\Manc$ of manifolds with corners. The relation of our definitions to others in the literature is explained in Remark~\ref{kh3rem1}.

\begin{dfn} Use the notation $\R^m_k=[0,\iy)^k\t\R^{m-k}$
for $0\le k\le m$, and write points of $\R^m_k$ as $u=(u_1,\ldots,u_m)$ for $u_1,\ldots,u_k\in[0,\iy)$, $u_{k+1},\ldots,u_m\in\R$. Let $U\subseteq\R^m_k$ and $V\subseteq \R^n_l$ be open, and $f=(f_1,\ldots,f_n):U\ra V$ be a continuous map, so that $f_j=f_j(u_1,\ldots,u_m)$ maps $U\ra[0,\iy)$ for $j=1,\ldots,l$ and $U\ra\R$ for $j=l+1,\ldots,n$. Then we say that:
\begin{itemize}
\setlength{\itemsep}{0pt}
\setlength{\parsep}{0pt}
\item[(a)] $f$ is {\it weakly smooth\/} if all derivatives $\frac{\pd^{a_1+\cdots+a_m}}{\pd u_1^{a_1}\cdots\pd u_m^{a_m}}f_j(u_1,\ldots,u_m):U\ra\R$ exist and are continuous in for all $j=1,\ldots,m$ and $a_1,\ldots,a_m\ge 0$, including one-sided derivatives where $u_i=0$ for $i=1,\ldots,k$.

By Seeley's Extension Theorem, this is equivalent to requiring $f_j$ to extend to a smooth function $f_j':U'\ra\R$ on open neighbourhood $U'$ of $U$ in~$\R^m$.
\item[(b)] $f$ is {\it smooth\/} if it is weakly smooth and every $u=(u_1,\ldots,u_m)\in U$ has an open neighbourhood $\ti U$ in $U$ such that for each $j=1,\ldots,l$, either:
\begin{itemize}
\setlength{\itemsep}{0pt}
\setlength{\parsep}{0pt}
\item[(i)] we may uniquely write $f_j(\ti u_1,\ldots,\ti u_m)=F_j(\ti u_1,\ldots,\ti u_m)\cdot\ti u_1^{a_{1,j}}\cdots\ti u_k^{a_{k,j}}$ for all $(\ti u_1,\ldots,\ti u_m)\in\ti U$, where $F_j:\ti U\ra(0,\iy)$ is weakly smooth and $a_{1,j},\ldots,a_{k,j}\in\N=\{0,1,2,\ldots\}$, with $a_{i,j}=0$ if $u_i\ne 0$; or 
\item[(ii)] $f_j\vert_{\smash{\ti U}}=0$.
\end{itemize}
\item[(c)] $f$ is {\it strongly smooth\/} if it is smooth, and in case (b)(i), for each $j=1,\ldots,l$ we have $a_{i,j}=1$ for at most one $i=1,\ldots,k$, and $a_{i,j}=0$ otherwise. 
\item[(d)] $f$ is {\it interior\/} if it is smooth, and case (b)(ii) does not occur.
\item[(e)] $f$ is a {\it diffeomorphism\/} if it is a bijection, and both $f:U\ra V$ and $f^{-1}:V\ra U$ are weakly smooth (which implies smooth and strongly smooth).
\end{itemize}
\label{kh3def1}
\end{dfn}

\begin{dfn} Let $X$ be a second countable Hausdorff topological space. An {\it $m$-dimensional chart on\/} $X$ is a pair $(U,\phi)$, where
$U\subseteq\R^m_k$ is open for some $0\le k\le m$, and $\phi:U\ra X$ is a
homeomorphism with an open set~$\phi(U)\subseteq X$.

Let $(U,\phi),(V,\psi)$ be $m$-dimensional charts on $X$. We call $(U,\phi)$ and $(V,\psi)$ {\it compatible\/} if $\psi^{-1}\ci\phi:\phi^{-1}\bigl(\phi(U)\cap\psi(V)\bigr)\ra
\psi^{-1}\bigl(\phi(U)\cap\psi(V)\bigr)$ is a diffeomorphism between open subsets of $\R^m_k,\R^m_l$, in the sense of Definition~\ref{kh3def1}(e).

An $m$-{\it dimensional atlas\/} for $X$ is a system
$\{(U_a,\phi_a):a\in A\}$ of pairwise compatible $m$-dimensional
charts on $X$ with $X=\bigcup_{a\in A}\phi_a(U_a)$. We call such an
atlas {\it maximal\/} if it is not a proper subset of any other
atlas. Any atlas $\{(U_a,\phi_a):a\in A\}$ is contained in a unique
maximal atlas, the set of all charts $(U,\phi)$ on $X$
which are compatible with $(U_a,\phi_a)$ for all~$a\in A$.

An $m$-{\it dimensional manifold with corners\/} is a second
countable Hausdorff topological space $X$ equipped with a maximal
$m$-dimensional atlas. Usually we refer to $X$ as the manifold,
leaving the atlas implicit, and by a {\it chart\/ $(U,\phi)$ on\/}
$X$, we mean an element of the maximal atlas.

Now let $X,Y$ be manifolds with corners of dimensions $m,n$, and $f:X\ra Y$ a continuous map. We call $f$ {\it weakly smooth}, or {\it smooth}, or {\it strongly smooth}, or {\it interior}, if whenever $(U,\phi),(V,\psi)$ are charts on $X,Y$ with $U\subseteq\R^m_k$, $V\subseteq\R^n_l$ open, then
\begin{equation*}
\psi^{-1}\ci f\ci\phi:(f\ci\phi)^{-1}(\psi(V))\longra V
\end{equation*}
is weakly smooth, or smooth, or strongly smooth, or interior, respectively, as maps between open subsets of $\R^m_k,\R^n_l$ in the sense of Definition \ref{kh3def1}. 

These three classes of (a) weakly smooth, (b) smooth, and (c) strongly smooth maps of manifolds with corners, all contain identities and are closed under composition, so each makes manifolds with corners into a category, which we write as $\Mancwe,\Manc,\Mancst$, respectively.
\label{kh3def2}
\end{dfn}

\begin{rem}{\bf(a)} The differences between $\Mancwe,\Manc,\Mancst$ will not concern us much in this paper, as we will be principally working with maps $f:V\ra Y$ with $V$ a manifold with corners and $Y$ a manifold without boundary, and for such $f$ weakly smooth = smooth = strongly smooth.

However, the differences would matter if we wanted to extend our theory to define $MH_*(Y;R),MH^*(Y;R),\ldots$ for $Y$ a manifold with corners.
\smallskip

\noindent{\bf(b)} Some references on manifolds with corners are Cerf \cite{Cerf}, Douady \cite{Doua}, Kottke and Melrose \cite{KoMe}, Margalef-Roig and Outerelo Dominguez \cite{MaOu}, Melrose \cite{Melr1,Melr2,Melr3}, and the author \cite{Joyc3}, \cite[\S 5]{Joyc7}. Just as objects, without considering morphisms, most authors define manifolds with corners as in Definition \ref{kh3def2}. However, Melrose \cite{KoMe,Melr1,Melr2,Melr3} and authors who follow him impose an extra condition: in Definition \ref{kh3def4} we will define the boundary $\pd X$ of a manifold with corners $X$, with an immersion $i_X:\pd X\ra X$. Melrose requires that $i_X\vert_C:C\ra X$ should be injective for each connected component $C$ of $\pd X$ (such $X$ are sometimes called {\it manifolds with faces\/}).

There is no general agreement in the literature on how to define smooth maps, or morphisms, of manifolds with corners: 
\begin{itemize}
\setlength{\itemsep}{0pt}
\setlength{\parsep}{0pt}
\item[(i)] Our notion of `smooth map' in Definitions \ref{kh3def1} and \ref{kh3def2} is due to Melrose \cite[\S 1.12]{Melr3}, \cite[\S 1]{Melr1}, \cite[\S 1]{KoMe}, who calls them {\it b-maps}. 
\item[(ii)] The author \cite{Joyc3} defined and studied `strongly smooth maps' above (which were just called `smooth maps' in \cite{Joyc3}). 
\item[(iii)] Most other authors, such as Cerf \cite[\S I.1.2]{Cerf}, define smooth maps of manifolds with corners to be weakly smooth maps, in our notation.
\end{itemize}
\label{kh3rem1}
\end{rem}

\begin{dfn} Let $U\subseteq\R^m_k$ be open. For each $u=(u_1,\ldots,u_m)$ in $U$, define the {\it depth\/} $\depth_Uu$ of $u$ in $U$ to be the number of $u_1,\ldots,u_k$ which are zero. That is, $\depth_Uu$ is the number of boundary faces of $U$ containing~$u$.

Let $X$ be an $m$-manifold with corners. For $x\in X$, choose a
chart $(U,\phi)$ on the manifold $X$ with $\phi(u)=x$ for $u\in U$,
and define the {\it depth\/} $\depth_Xx$ of $x$ in $X$ by
$\depth_Xx=\depth_Uu$. This is independent of the choice of
$(U,\phi)$. For each $l=0,\ldots,m$, define the {\it depth\/ $l$
stratum\/} of $X$ to be
\begin{equation*}
S^l(X)=\bigl\{x\in X:\depth_Xx=l\bigr\}.
\end{equation*}
Then $X=\coprod_{l=0}^mS^l(X)$ and $\overline{S^l(X)}=
\bigcup_{k=l}^m S^k(X)$. The {\it interior\/} of $X$ is
$X^\ci=S^0(X)$. Each $S^l(X)$ has the structure of an
$(m-l)$-manifold without boundary.

\label{kh3def3}
\end{dfn}

\begin{dfn} Let $X$ be a manifold with corners, and $x\in X$. A
{\it local boundary component\/ $\be$ of\/ $X$ at\/} $x$ is a local
choice of connected component of $S^1(X)$ near $x$. That is, for
each sufficiently small open neighbourhood $V$ of $x$ in $X$, $\be$
gives a choice of connected component $W$ of $V\cap S^1(X)$ with
$x\in\overline W$, and any two such choices $V,W$ and $V',W'$ must
be compatible in the sense that~$x\in\overline{(W\cap W')}$.

As a set, define the {\it boundary\/}
\begin{equation*}
\pd X=\bigl\{(x,\be):\text{$x\in X$, $\be$ is a local boundary
component for $X$ at $x$}\bigr\}.
\end{equation*}
Define a map $i_X:\pd X\ra X$ by $i_X:(x,\be)\mapsto x$. Note that
$i_X$ {\it need not be injective}, as $\bmd{i_X^{-1}(x)}=\depth_Xx$
for all~$x\in X$.

If $(U,\phi)$ is a chart on $X$ with $U\subseteq\R^m_k$ open, then for each $i=1,\ldots,k$ we can define a chart $(U_i,\phi_i)$ on $\pd X$ by
\begin{align*}
&U_i=\bigl\{(x_1,\ldots,x_{m-1})\in \R^{m-1}_{k-1}:
(x_1,\ldots,x_{i-1},0,x_i,\ldots,x_{m-1})\in
U\subseteq\R^m_k\bigr\},\\
&\phi_i:(x_1,\ldots,x_{m-1})\longmapsto\bigl(\phi
(x_1,\ldots,x_{i-1}, 0,x_i,\ldots,x_{m-1}),\phi_*(\{x_i=0\})\bigr).
\end{align*}
All such charts on $\pd X$ are compatible, so they induce the structure of a manifold with corners on $\pd X$, with $\dim(\pd X)=\dim X-1$. Also $i_X:\pd X\ra X$ is a proper, smooth map, which is strongly smooth, but not interior. We can iterate boundaries to form $\pd^2X=\pd(\pd X),\pd^3X,\ldots,$ with $\pd^kX$ a manifold with corners of dimension~$\dim X-k$.

We call $X$ a {\it manifold with boundary\/} if $\pd^2X=\es$, and a {\it manifold without boundary\/} if $\pd X=\es$. Write $\Man,\Manb$ for the full subcategories of manifolds without boundary, and manifolds with boundary, in $\Manc$, so that $\Man\subset\Manb\subset\Manc$, and $\Man$ is the usual category of manifolds.
\label{kh3def4}
\end{dfn}

Let $X$ be a manifold with corners. Then there are natural identifications 
\e
\begin{split}
\pd^kX\cong\bigl\{(x,\be_1,\ldots,\be_k):\,&\text{$x\in X,$
$\be_1,\ldots,\be_k$ are distinct}\\
&\text{local boundary components for $X$ at $x$}\bigr\}.
\end{split}
\label{kh3eq1}
\e
Write $S_k$ for the symmetric group on $k$ elements, the group of bijections $\si:\{1,\ldots,k\}\ra\{1,\ldots,k\}$. From \eq{kh3eq1} we see that $\pd^kX$ has a natural, free action of $S_k$, which is by diffeomorphisms, given by
\begin{equation*}
\si:(x,\be_1,\ldots,\be_k)\longmapsto
(x,\be_{\si(1)},\ldots,\be_{\si(k)}).
\end{equation*}

Manifolds with corners $X$ have two notions of tangent bundle with
functorial properties, the ({\it ordinary\/}) {\it tangent bundle\/}
$TX$, the obvious generalization of tangent bundles of manifolds
without boundary, and the {\it b-tangent bundle\/} ${}^bTX$
introduced by Melrose \cite[\S 2.2]{Melr2}, \cite[\S I.10]{Melr3},
\cite[\S 2]{Melr1}. Taking duals gives two notions of cotangent
bundle $T^*X,{}^bT^*X$. If $(U,\phi)$ is a chart on $X$, with $U\subseteq\R^m_k$ open, and $(x_1,\ldots,x_m)$ the coordinates on $U$ with $x_1,\ldots,x_k\in[0,\iy)$ and $x_{k+1},\ldots,x_m\in\R$, then over $\phi(U)$ we have
\begin{itemize}
\setlength{\itemsep}{0pt}
\setlength{\parsep}{0pt}
\item[(i)] $TX$ is the trivial vector bundle with basis of sections $\frac{\pd}{\pd x_1},\ldots,\frac{\pd}{\pd x_m}$;
\item[(ii)] $T^*X$ is the vector bundle with basis $\d x_1,\ldots,\d x_m$;
\item[(iii)] ${}^bTX$ is the vector bundle with basis $x_1\frac{\pd}{\pd x_1},\ldots,x_k\frac{\pd}{\pd x_k},\ab\frac{\pd}{\pd x_{k+1}},\ab\ldots,\frac{\pd}{\pd x_m}$; and
\item[(iv)] ${}^bT^*X$ is the vector bundle with basis $x_1^{-1}\d x_1,\ldots,\ab x_k^{-1}\d x_k,\ab\d x_{k+1},\ldots,\d x_m$.
\end{itemize}
If $f:X\ra Y$ is smooth we have a vector bundle morphism $\d f:TX\ra f^*(TY)$. If $f$ is interior we have a vector bundle morphism ${}^b\d f:{}^bTX\ra f^*({}^bTY)$.

\begin{dfn} An interior map $f:X\ra Y$ in $\Manc$ is called a {\it b-submersion\/} if ${}^b\d f\vert_x:{}^bT_xX\ra{}^bT_yY$ is surjective for all $x\in X$ with~$f(x)=y\in Y$.

Interior maps $g:X\ra Z$, $h:Y\ra Z$ in $\Manc$ are called {\it b-transverse\/} if 
\begin{equation*}
{}^b\d g\vert_x\op {}^b\d h\vert_y:{}^bT_xX\op{}^bT_yY\longra{}^bT_zZ
\end{equation*}
is surjective for all $x\in X$ and $y\in Y$ with $g(x)=h(y)=z$.
\label{kh3def5}
\end{dfn}

Note that if $g:X\ra Z$ is a b-submersion and $h:Y\ra Z$ is interior then $g,h$ are b-transverse. The question of when (b-)transverse fibre products $X\t_{g,Z,h}Y$ exist in $\Manc$ (or $\Mancst,\Mancwe,\ldots$) is complicated; see the author \cite{Joyc3,Joyc8} and Kottke and Melrose \cite{KoMe} for some results, of which the strongest \cite[Th.~4.27]{Joyc8} says that all b-transverse fibre products exist in the category $\Man^{\bf gc}_{\bf in}$ of manifolds with {\it generalized\/} corners and interior maps. However, we only need the case in which $\pd Z=\es$, i.e.\ $Z\in\Man\subset\Manc$, which is much simpler.

\begin{prop} Suppose $g:X\ra Z,$ $h:Y\ra Z$ are b-transverse maps in $\Manc$ with\/ $\pd Z=\es$. Then a fibre product\/ $W=X\t_{g,Z,h}Y$ exists in $\Manc$ with\/ $\dim W=\dim X+\dim Y-\dim Z$. Furthermore, $g\ci i_X,h$ and\/ $g,h\ci i_Y$ are also b-transverse, and there is a natural diffeomorphism
\begin{equation*}
\pd(X\t_{g,Z,h}Y)\cong \bigl(\pd X\t_{g\ci i_X,Z,h}Y\bigr)\amalg \bigl(X\t_{g,Z,h\ci i_Y}\pd Y\bigr).
\end{equation*}

\label{kh3prop1}
\end{prop}

Orientations on manifolds with corners are discussed by the author \cite[\S 7]{Joyc3}, \cite[\S 5.8]{Joyc6} and Fukaya et al.\ \cite[\S 8.2]{FOOO1}.

\begin{dfn} Let $X$ be a manifold with corners with $\dim X=m$. Then $\La^m(T^*X)$ is a real line bundle on $X$. An {\it orientation\/} $o_X$ on $X$ is an equivalence class $[\om]$ of top-dimensional forms $\om\in C^\iy(\La^mT^*X)$ with $\om\vert_x\ne 0$ for all $x\in X$, where two such $\om,\om'$ are equivalent if $\om'=K\cdot\om$ for $K:X\ra(0,\iy)$ smooth. The {\it opposite orientation\/} is $-o_X=[-\om]$. Then we call $(X,o_X)$ an {\it oriented manifold with corners}. Usually we suppress the orientation $o_X$, and just refer to $X$ as an oriented manifold with corners. When $X$ is an oriented manifold with corners, we write $-X$ for $X$ with the opposite orientation.

We can also define orientations using $\La^m({}^bT^*X)$ instead of $\La^mT^*X$. Since $T^*X$ and ${}^bT^*X$ coincide on $X^\ci$, the two notions turn out to be equivalent.

We will also need a notion of relative orientation of a smooth map $f:X\ra Y$. A {\it coorientation\/} $c_f$ for $f$ is a an equivalence class $[\ga]$ of $\ga\in C^\iy\bigl(\La^{\dim X}T^*X\ot f^*(\La^{\dim Y}T^*Y)^*\bigr)$ with $\ga\vert_x\ne 0$ for all $x\in X$, where two such $\ga,\ga'$ are equivalent if $\ga'=K\cdot\ga$ for $K:X\ra(0,\iy)$ smooth. The {\it opposite coorientation\/} is $-c_f=[-\ga]$. Then we call $(f,c_f):X\ra Y$ a {\it cooriented smooth map}. Usually we suppress  $c_f$, and just refer to $f$ as a cooriented smooth map. If $Y$ is oriented then orientations on $X$ are equivalent to coorientations on $f:X\ra Y$. Orientations on $X$ are equivalent to coorientations on $\pi:X\ra *$, for $*$ the point.
\label{kh3def6}
\end{dfn}

If $X$ is an oriented manifold with corners we can define a natural orientation on $\pd X$, and hence on $\pd^2X,\pd^3X,\ldots,\pd^{\dim X}X$, and if $X,Y,Z$ are oriented manifolds with corners and $W=X\t_{g,Z,h}Y$ is a transverse fibre product in $\Manc$ then we can define a natural orientation on $W$. Products $X\t Y$ are fibre products over $Z=*$, so this includes orienting products.

There are similar facts about coorientations, e.g.\ for any $X$ (oriented or not) $i_X:\pd X\ra X$ has a natural coorientation, and if $g:X\ra Y$ is cooriented and $Y$ is oriented then a transverse fibre product $W=X\t_{g,Z,h}Y$ is naturally oriented. To do all this requires a choice of {\it orientation convention}. Ours follow Fukaya et al.\  \cite[Conv.~45.1]{FOOO1}. Different conventions would change the signs in~\eq{kh3eq2}--\eq{kh3eq5}.

If $X$ is an oriented manifold with corners, then as above each $\si\in S_k$ acts on the oriented manifold with corners $\pd^kX$. This multiplies the orientation on $\pd^kX$ by $\mathop{\rm sign}\si$. In particular, it will be important below that the free action of $S_2\cong\Z_2$ on $\pd^2X$ is orientation-reversing.

In Proposition \ref{kh3prop1}, if $X,Y,Z$ are oriented then in oriented manifolds with corners, as in \cite[Prop.~7.4]{Joyc3}, equation \eq{kh3eq2} becomes
\e
\pd(X\t_{g,Z,h}Y)\cong \bigl(\pd X\t_{g\ci i_X,Z,h}Y\bigr)\amalg (-1)^{\dim X+\dim Z}\bigl(X\t_{g,Z,h\ci i_Y}\pd Y\bigr).
\label{kh3eq2}
\e
Here \cite[Prop.~7.5]{Joyc3} are some more identities on orientations:

\begin{prop}{\bf(a)} If\/ $g:X\ra Z,$ $h:Y\ra Z$ are transverse
smooth maps of oriented manifolds with corners then in oriented
manifolds we have
\e
X\t_{g,Z,h}Y\cong(-1)^{(\dim X-\dim Z)(\dim Y-\dim Z)}Y\t_{h,Z,g}X.
\label{kh3eq3}
\e

\noindent{\bf(b)} If\/ $e:V\ra Y,$ $f:W\ra Y,$ $g:W\ra Z,$ $h:X\ra
Z$ are smooth maps of oriented manifolds with corners then in
oriented manifolds we have
\e
V\t_{e,Y,f\ci\pi_W}\bigl(W\t_{g,Z,h}X\bigr)\cong
\bigl(V\t_{e,Y,f}W\bigr)\t_{g\ci\pi_W,Z,h}X,
\label{kh3eq4}
\e
provided all four fibre products are transverse.
\smallskip

\noindent{\bf(c)} If\/ $e:V\ra Y,$ $f:V\ra Z,$ $g:W\ra Y,$ $h:X\ra
Z$ are smooth maps of oriented manifolds with corners then in
oriented manifolds we have
\e
\begin{split}
&V\t_{(e,f),Y\t Z,g\t h}(W\t X)\cong \\
&\quad(-1)^{\dim Z(\dim Y+\dim W)}
(V\t_{e,Y,g}W)\t_{f\ci\pi_V,Z,h}X,
\end{split}
\label{kh3eq5}
\e
provided all three fibre products are transverse.
\label{kh3prop2}
\end{prop}

Equations \eq{kh3eq3}--\eq{kh3eq4} will be important in proving that the cup product on cohomology $MH^*(Y;R)$ is supercommutative and associative.

\subsection{\texorpdfstring{Hausdorff measure and Sard's Theorem}{Hausdorff measure and Sard\textquoteright s Theorem}}
\label{kh32}

In \S\ref{kh33} we will give our list Assumptions \ref{kh3ass1}--\ref{kh3ass7} of properties of a category $\tManc$ of `manifolds with corners' we need to construct our (co)homology theories. The final one, Assumption \ref{kh3ass7}, requires a version of Sard's Theorem to hold for morphisms $f:X\ra Y$ in $\tManc$ with $Y\in\Man$, stated in terms of Hausdorff measure in $Y$, so we explain these here.

\begin{dfn} Let $Y$ be a smooth manifold, $S\subseteq Y$ a subset, and $d\ge 0$. Fix a Riemannian metric $g$ on $Y$. We say that $S$ {\it has zero Hausdorff\/ $d$-measure}, or is $d$-{\it null}, written $\cH^d(S)=0$, if for all $\ep>0$ there exist $y_i\in Y$ and $r_i> 0$ for all $i$ in some finite or countable indexing set $I$, such that $\sum_{i\in I}r_i^d<\ep$ and $S\subseteq \bigcup_{i\in I}B_{r_i}(y_i)$, for $B_{r_i}(y_i)$ the open ball about $y_i$ in $Y$ defined using the metric $g$. Whether $\cH^d(S)=0$ or not is independent of the choice of~$g$.

Some properties of this definition: if $d>\dim Y$ then $\cH^d(S)=0$ for any $S\subseteq Y$. If $Y_1,Y_2$ are manifolds and $f:Y_1\ra Y_2$ is smooth (or more generally $C^1$, or locally Lipschitz) and $S\subseteq Y_1$ with $\cH^d(S)=0$, then $\cH^d(f(S))=0$ in~$Y_2$.

\label{kh3def7}
\end{dfn}

Then Sard's Theorem \cite{Sard} states:

\begin{thm}{\bf(a)} Suppose $X,Y$ are manifolds with\/ $\dim X=m,$ $\dim Y=n$ for $m\le n,$ and\/ $f:X\ra Y$ is a smooth map (or more generally a $C^k$ map for any $k\ge 1$). Define $S\subseteq X$ to be the closed subset of\/ $x\in X$ for which\/ $T_xf:T_xX\ra T_{f(x)}Y$ is not injective, so that\/ $f(S)$ is a subset of\/ $Y,$ which is closed if\/ $f$ is proper. Then $\cH^m(f(S))=0$ in~$Y$.
\smallskip

\noindent{\bf(b)} Suppose $X,Y$ are manifolds with\/ $\dim X=m,$ $\dim Y=n$ for $m\ge n,$ and\/ $f:X\ra Y$ is a smooth map (or more generally a $C^k$ map for any $k\ge m-n+1$). Define $S\subseteq X$ to be the closed subset of\/ $x\in X$ for which\/ $T_xf:T_xX\ra T_{f(x)}Y$ is not surjective. Then $\cH^n(f(S))=0$ in $Y$.
\label{kh3thm1}
\end{thm}

Essentially this says that for $f:X\ra Y$ a smooth map of manifolds, if $\dim X\le\dim Y$ then $f$ is an immersion over $Y\sm S$ for $S\subset Y$ a `negligible' subset (roughly, $\dim S<\dim X\le\dim Y$), and if $\dim X\ge\dim Y$ then $f$ is a submersion over $Y\sm S$ for $S\subset Y$ a `negligible' subset (roughly, $\dim S<\dim Y\le\dim X$). The requirements on $k$ are sharp.

\subsection{\texorpdfstring{Assumptions we need on `manifolds with corners'}{Assumptions we need on \textquoteleft manifolds with corners\textquoteright}}
\label{kh33}

The next seven assumptions give the properties of a category of `manifolds with corners' $\tManc$ we need for both `integral M-(co)homology' in \S\ref{kh4}, and `rational M-(co)homology' in \S\ref{kh51}. Examples of categories $\tManc$ satisfying these assumptions are given in Example~\ref{kh3ex1}.

\begin{ass}{\bf(Category-theoretic properties.)} {\bf(a)} We are given a category $\tManc$. For simplicity, from \S\ref{kh4} onwards, objects $X$ in $\tManc$ will be called {\it manifolds with corners} (although they may in examples not be manifolds, but some kind of singular space), and morphisms $f:X\ra Y$ in $\tManc$ will be called {\it smooth maps} (although they may in examples be non-smooth). 

Isomorphisms in $\tManc$ are called {\it diffeomorphisms}.
\smallskip

\noindent{\bf(b)} There is an object $\es\in\tManc$ called the {\it empty set}, which is an initial object in $\tManc$ (i.e.\ every $X\in\tManc$ has a unique morphism~$\es\ra X$).
\smallskip

\noindent{\bf(c)} There is an object $*\in\tManc$ called the {\it point}, which is a final object in $\tManc$ (i.e.\ every $X\in\tManc$ has a unique morphism~$X\ra *$).
\smallskip

\noindent{\bf(d)} Each object $X$ in $\tManc$ has a {\it dimension\/} $\dim X\in\N=\{0,1,\ldots\}$, except that $\dim\es$ is undefined, or allowed to take any value. We have $\dim *=0$.
\smallskip

\noindent{\bf(e)} {\it Products\/} $X\t Y$ of objects $X,Y\in\tManc$ exist in $\tManc$, in the sense of category theory (fibre products over $*$), with projections $\pi_X:X\t Y\ra X$ and $\pi_Y:X\t Y\ra Y$. They have $\dim (X\t Y)=\dim X+\dim Y$. Hence {\it products\/} $f\t g:W\t X\ra Y\t Z$ of morphisms $f:W\ra Y$, $g:X\ra Z$, and {\it direct products\/} $(f,g):X\ra Y\t Z$ of $f:X\ra Y$, $g:X\ra Z$, exist in $\tManc$.
\smallskip

\noindent{\bf(f)} If $X,Y\in\tManc$ with $\dim X=\dim Y$ there is a {\it disjoint union\/} $X\amalg Y$ in $\tManc$ with inclusion morphisms $\io_X:X\hookra X\amalg Y$, $\io_Y:Y\hookra X\amalg Y$. It is a coproduct in the sense of category theory, with~$\dim(X\amalg Y)=\dim X=\dim Y$.

\label{kh3ass1}
\end{ass}

\begin{ass}{\bf(Underlying topological spaces.)} {\bf(a)} There is a functor $F_\tManc^\Top:\tManc\ra\Top$ from $\tManc$ to the category of topological spaces $\Top$, mapping objects $X\in\tManc$ to the {\it underlying topological space\/} $X_\top:=F_\tManc^\Top(X)$, and morphisms $f:X\ra Y$ to $f_\top:=F_\tManc^\Top(f):X_\top\ra Y_\top$. So we can think of objects $X$ of $\tManc$ as `topological spaces $X_\top$ with extra structure'. From \S\ref{kh4} on, we will often write $X,f$ instead of~$X_\top,f_\top$.
\smallskip

\noindent{\bf(b)} Underlying topological spaces $X_\top$ are Hausdorff and locally compact, and $F_\Manc^\Top(\es)=\es$, and $F_\Manc^\Top(*)$ is a point.
\smallskip

\noindent{\bf(c)} $F_\tManc^\Top$ takes products and disjoint unions in $\tManc$ functorially to products and disjoint unions in $\Top$.
\smallskip

\noindent{\bf(d)} If $X\in\tManc$ and $U'\subseteq X_\top$ is open with inclusion $i':U'\hookra X_\top$, there is a natural object $U$ in $\tManc$ called an {\it open submanifold\/} with $U_\top=U'$ and $\dim U=\dim X$, and an {\it inclusion morphism\/} $i:U\hookra X$ with $i_\top=i'$. If $U'=\es$ then $U=\es$. Inclusion morphisms are functorial under inclusions of open sets $U'\hookra V'\hookra X_\top$. Given a morphism $f:X\ra Y$, we often write $f\vert_U:U\ra Y$ instead of~$f\ci i:U\ra Y$.

If $f:W\ra X$ is a morphism in $\tManc$ with $f_\top(W_\top)\subseteq U_\top\subseteq X_\top$ then $f$ factorizes uniquely as $f=i\ci f'$ for a morphism $f':W\ra U$ in~$\tManc$.

Inclusions $\io_X:X\hookra X\amalg Y$, $\io_Y:Y\hookra X\amalg Y$ are open submanifolds.
\label{kh3ass2}
\end{ass}

In later sections we will generally drop the distinction between $X$ and $X_\top$, and write $x\in X$ rather than $x\in X_\top$, identify open submanifolds $i:U\hookra X$ with open sets $U\subseteq X$, and so on, just as one does for ordinary manifolds in differential geometry.

\begin{ass}{\bf(Boundaries.)} {\bf(a)} Each object $X$ in $\tManc$ has a {\it boundary\/} $\pd X$ in $\tManc$, with $\dim\pd X=\dim X-1$, where $\pd X=\es$ if $\dim X=0$, and a canonical morphism $i_X:\pd X\ra X$.
\smallskip

\noindent{\bf(b)} There is a class of morphisms $f:X\ra Y$ in $\tManc$ called {\it simple\/} morphisms. Simple morphisms include identities, diffeomorphisms, and inclusions of open submanifolds, and are closed under composition.

If $f:X\ra Y$ is simple, there is a canonical simple morphism $\pd f:\pd X\ra\pd Y$ with $f\ci i_X=i_Y\ci\pd f:\pd X\ra Y$, which is functorial in $f$. Thus $X\mapsto\pd X$, $f\mapsto\pd f$ is a functor $\pd:\mathop{\bf\check Man^c_{si}}\ra\mathop{\bf\check Man^c_{si}}$ on the subcategory $\mathop{\bf\check Man^c_{si}}$ of $\cManc$ with only simple morphisms.

If $j:U\hookra X$ is an open submanifold then $\pd j:\pd U\hookra\pd X$ makes $\pd U$ into an open submanifold of $\pd X$, with~$(\pd U)_\top=i_{X,\top}^{-1}(U_\top)$.
\smallskip

\noindent{\bf(c)} For each object $X\in\tManc$ there is a canonical diffeomorphism $\si_X:\pd^2X\ra\pd^2X$ with $\si_X\ci\si_X=\id_{\pd^2X}$ and $i_X\ci i_{\pd X}\ci\si=i_X\ci i_{\pd X}:\pd^2X\ra X$.
\smallskip

\noindent{\bf(d)} If $X,Y\in\cManc$ there is a canonical diffeomorphism
\e
\pd(X\t Y)\cong \bigl(\pd X\t Y\bigr)\amalg\bigl(X\t\pd Y\bigr),
\label{kh3eq6}
\e
identifying $\pi_X\ci i_{X\t Y}$ with $(i_X\ci\pi_{\pd X})\amalg\pi_X$ and $\pi_Y\ci i_{X\t Y}$ with~$\pi_Y\amalg(i_Y\ci\pi_{\pd Y})$.

\smallskip

\noindent{\bf(e)} If $X,Y\in\cManc$ with $\dim X=\dim Y$ there is a canonical diffeomorphism
\begin{equation*}
\pd(X\amalg Y)\cong \pd X\amalg \pd Y,
\end{equation*}
identifying $i_{X\amalg Y}$ with~$i_{\pd X}\amalg i_{\pd Y}$.
\smallskip

\noindent{\bf(f)} $i_{X,\top}:\pd X_\top\ra X_\top$ is a finite, proper map in $\Top$ for any $X\in\tManc$.
\label{kh3ass3}
\end{ass}

\begin{ass}{\bf(Relation with ordinary manifolds.)} {\bf(a)} The category $\Man$ of (ordinary) smooth manifolds without boundary and (ordinary) smooth maps between them is included as a strictly full subcategory~$\Man\subset\tManc$.

Dimensions of objects in $\Man\subset\tManc$ are as usual in $\Man$. Products and disjoint unions in $\tManc$ of $X,Y\in\Man$ agree with those in $\Man$. The empty set $\es$ and point $*$ in Assumption \ref{kh3ass1}(b),(c) lie in $\Man\subset\tManc$.

The underlying topological space functor $\smash{F_\tManc^\Top}$ is as usual on $\Man\subset\Manc$. Open submanifolds in $\Man,\tManc$ agree.
\smallskip

\noindent{\bf(b)} For each object $X\in\tManc$ there is a natural open submanifold $X^\ci\hookra X$ called the {\it interior\/} of $X$, with $X^\ci\in\Man\subset\tManc$. It is the largest open submanifold of $X$ lying in $\Man$, that is, if $i:U\hookra X$ is an open submanifold then $U\in\Man$ if and only if $U_\top\subseteq X_\top^\ci\subseteq X_\top$.

For objects $X,Y\in\tManc$ we have $(X\t Y)^\ci=X^\ci\t Y^\ci$. For objects $X,Y\in\tManc$ with $\dim X=\dim Y$ we have $(X\amalg Y)^\ci=X^\ci\amalg Y^\ci$.
\smallskip

\noindent{\bf(c)} Intervals $I$ in $\R$ such as $[0,1]$, $[0,\iy)$, $(-\iy,0]$ are objects in $\tManc$. For all $k\ge 0$, the $k$-{\it simplex\/} 
\begin{equation*}
\De_k=\bigl\{(x_0,\ldots,x_k)\in\R^{k+1}:x_i\ge 0,\;\>
x_0+\ldots+x_k=1\bigr\}
\end{equation*}
is an object in $\tManc$. Dimensions, and topological spaces, and interiors, of all these are as in $\Manc$ in~\S\ref{kh31}.

The following morphisms in $\Manc$ in \S\ref{kh31} are also morphisms in $\tManc$:
\begin{itemize}
\setlength{\itemsep}{0pt}
\setlength{\parsep}{0pt}
\item[(A)] Smooth maps $f:I\ra Y$ for $I$ an interval in $\R$ and $Y\in\Man$.
\item[(B)] Smooth maps $f:I\t X\ra Y$ for $I$ an interval in $\R$ and $X,Y\in\Man$.
\item[(C)] Smooth maps $f:\De_k\ra Y$ for $Y\in\Man$.
\item[(D)] $i_{\De_k}:\pd\De_k\cong\coprod_{j=0}^k\De_{k-1}\ra\De_k$ for $k>0$.
\end{itemize}
These behave as expected on topological spaces, and under compositions, products, and disjoint unions, with each other and with morphisms in~$\Man$.
\smallskip

\noindent{\bf(d)} Any $X\in\Man\subset\tManc$ has $\pd X=\es$. Boundaries of intervals in $\R$ and $k$-simplices $\De_k$ behave as in $\Manc$ in \S\ref{kh31}, so for example $\pd [0,\iy)=\{0\}$, and $\pd\De_k$ is diffeomorphic to the disjoint union of $k+1$ copies of $\De_{k-1}$ for~$k>0$.
\smallskip

\noindent{\bf(e)} Any $X\in\tManc$ with $\dim X=0$ has $X^\ci=X$, so $X$ lies in $\Man\subset\tManc$, and $X_\top$ is a set with the discrete topology.
\label{kh3ass4}
\end{ass}

\begin{ass}{\bf(Submersions and transverse fibre products.)} {\bf(a)} There is a class of morphisms $f:X\ra Y$ in $\tManc$ called {\it submersions}, with $\dim X\ge\dim Y$. They include identities, diffeomorphisms, inclusions of open submanifolds $i:U\hookra X$, and projections~$\pi_X:X\t Y\ra X$. 

Submersions in $\tManc$ are closed under composition. 

If $f:X\ra Y$ is a morphism in $\Man$ then $f$ is a submersion in $\tManc$ if and only if it is a submersion of manifolds in the usual sense.
\smallskip

\noindent{\bf(b)} If $f:X\ra Y$ is a submersion in $\tManc$ and $Y\in\Man$ then $f\ci i_X:\pd X\ra Y$ is a submersion in $\tManc$.
\smallskip

\noindent{\bf(c)} Let $g:X\ra Z,$ $h:Y\ra Z$ be morphisms in $\tManc$ with $Z\in\Man\subset\tManc$. Then there is a notion of when $g,h$ are {\it transverse}. 

If $g,h$ are transverse then a fibre product $W=X\t_{g,Z,h}Y$ exists in $\tManc$ in the sense of category theory, with $\dim W=\dim X+\dim Y-\dim Z$, in a commutative, Cartesian square in~$\tManc$:
\e
\begin{gathered}
\xymatrix@C=70pt@R=13pt{
*+[r]{W} \ar[r]_{\pi_Y} \ar[d]^{\pi_X} & *+[l]{Y} \ar[d]_h \\
*+[r]{X} \ar[r]^g & *+[l]{Z.\!} }
\end{gathered}
\label{kh3eq7}
\e

The functor $F_\tManc^\Top:\tManc\ra\Top$ should take such transverse fibre products in $\tManc$ to fibre products in $\Top$, so that we have a homeomorphism
\e
W_\top\cong\bigl\{(x,y)\in X_\top\t Y_\top:g_\top(x)=h_\top(y)\bigr\}.
\label{kh3eq8}
\e
There should be a natural diffeomorphism
\e
\pd(X\t_{g,Z,h}Y)\cong \bigl(\pd X\t_{g\ci i_X,Z,h}Y\bigr)\amalg \bigl(X\t_{g,Z,h\ci i_Y}\pd Y\bigr),
\label{kh3eq9}
\e
where $g\ci i_X,h$ and $g,h\ci i_Y$ are also transverse.

If $g$ is a submersion (or $h$ is a submersion) then $g,h$ are transverse, and $\pi_Y:W\ra Y$ is a submersion (or $\pi_X:W\ra X$ is a submersion, respectively).
\smallskip

\noindent{\bf(d)} If $X,Y\in\Man$ in {\bf(c)} then transversality of $g,h$ in $\tManc$ is the usual notion of transverse smooth maps in $\Man$, and $W=X\t_{g,Z,h}Y$ lies in $\Man\subset\tManc$ and is the usual transverse fibre product in differential geometry.

More generally, for $X,Y\in\tManc$ in {\bf(c)} we have
\e
(X\t_{g,Z,h}Y)^\ci=X^\ci\t_{g\vert_{X^\ci},Z,h\vert_{Y^\ci}}Y^\ci,
\label{kh3eq10}
\e
where the r.h.s.\ of \eq{kh3eq10} is the usual transverse fibre product of manifolds.
\smallskip

\noindent{\bf(e)} Let $g:X\ra Z,$ $h:Y\ra Z$ be morphisms in $\tManc$ with $Z\in\Man\subset\tManc$. Then $g,h$ transverse implies that $g\ci i:U\ra Z$ and $h\ci j:V\ra Z$ are transverse for any open submanifolds $i:U\hookra X$ and $j:V\hookra Y$.

Conversely, if for all $x\in X_\top$, $y\in Y_\top$ with $g_\top(x)=h_\top(y)$ in $Z_\top$, there exist open submanifolds $i:U\hookra X$, $j:V\hookra Y$ with $x\in U_\top$, $y\in V_\top$ such that $g\ci i$ and $h\ci j$ are transverse, then $g,h$ are transverse.
\smallskip

\noindent{\bf(f)} Let $(g_1,g_2):X\ra Z_1\t Z_2$, $(h_1,h_2):Y\ra Z_1\t Z_2$ be morphisms in $\tManc$ with $Z_1,Z_2\in\Man$, and suppose $g_1:X\ra Z_1$ and $h_2:Y\ra Z_2$ are submersions. Then $(g_1,g_2)$ and $(h_1,h_2)$ are transverse.
\label{kh3ass5}
\end{ass}

\begin{ass}{\bf(Orientations and coorientations.)} {\bf(a)} For each object $X$ in $\tManc$, there is a notion of {\it orientation\/} $o_X$ on $X$, an additional geometric structure on $X$. Every orientation $o_X$ has an {\it opposite orientation\/} $-o_X$, with $-(-o_X)=o_X$. We call the pair $(X,o_X)$ an {\it oriented manifold with corners}. Often we suppress $o_X$ and call $X$ an oriented manifold with corners, and then we write $-X$ for $(X,-o_X)$ with the opposite orientation.
\smallskip

\noindent{\bf(b)} For each morphism $f:X\ra Y$ in $\tManc$, there is a notion of {\it coorientation\/} $c_f$, an additional geometric structure on $f$. Every coorientation $c_f$ has an {\it opposite coorientation\/} $-c_f$, with $-(-c_f)=c_f$. We call the pair $(f,c_f)$ a {\it cooriented morphism\/} in $\tManc$. Often we suppress $c_f$, and call $f$ a cooriented morphism.
\smallskip

\noindent{\bf(c)} Suppose $f:X\ra Y$ is a morphism in $\tManc$, and $o_Y$ an orientation on $Y$. Then there is a natural 1-1 correspondence between orientations $o_X$ on $X$ and coorientations $c_f$ for $f$. We write these as operations on coorientations by $o_X=c_f\cdot o_Y$ and $c_f=o_X/o_Y$. We have $(-c_f)\cdot o_Y=c_f\cdot(-o_Y)=-(c_f\cdot o_Y)$ and~$(-o_X)/o_Y=o_X/(-o_Y)=-(o_X/o_Y)$.
\smallskip

\noindent{\bf(d)} If $f:X\ra Y$, $g:Y\ra Z$ are morphisms in $\tManc$ and $c_f,c_g$ are coorientations for $f,g$, there is a natural coorientation $c_g\ci c_f$ for $g\ci f$. This composition of coorientations is associative, with $(-c_g)\ci c_f=c_g\ci(-c_f)=-(c_g\ci c_f)$. If $X,Y,Z$ have orientations $o_X,o_Y,o_Z$ with $c_f=o_X/o_Y$, $c_g=o_Y/o_Z$ then~$c_g\ci c_f=o_X/o_Z$.

\smallskip

\noindent{\bf(e)} Identity morphisms, diffeomorphisms, and inclusions of open submanifolds in $\tManc$ all have natural coorientations,  functorial under composition.

So, for example, if $f:X\ra Y$ is a diffeomorphism it has a natural coorientation $c_f$, and setting $o_X=c_f\cdot o_Y$ gives a 1-1 correspondence between orientations $o_X$ on $X$ and orientations $o_Y$ on $Y$, thought of as {\it pullback\/} $o_X=f^*(o_Y)$ or {\it pushforward\/}~$o_Y=f_*(o_X)$. 

Similarly, if $i:U\hookra X$ is an open submanifold it has a natural coorientation $c_i$, so orientations $o_X$ on $X$ map to orientations $o_U=c_i\cdot o_X$ on $U$, thought of as {\it restriction\/}~$o_U=o_X\vert_U$.

All the operations and properties of (co)orientations in {\bf(a)}--{\bf(n)} behave functorially under pullback / pushforward by diffeomorphisms, and restriction to open submanifolds.
\smallskip

\noindent{\bf(f)} If $(X,o_X),(Y,o_Y)$ are oriented objects in $\tManc$ there is a natural orientation $o_X\t o_Y$ on $X\t Y$, with $(-o_X)\t o_Y=o_X\t(-o_Y)=-(o_X\t o_Y)$. Under the identification $X\t Y\cong Y\t X$ we have $o_X\t o_Y=(-1)^{\dim X\dim Y}o_Y\t o_X$. If $(Z,o_Z)$ is an oriented object then $(o_X\t o_Y)\t o_Z=o_X\t(o_Y\t o_Z)$ on~$X\t Y\t Z$.

This also gives coorientations $c_{\pi_X}=(o_X\t o_Y)/o_X$ on $\pi_X:X\t Y\ra X$ and $c_{\pi_Y}=(o_X\t o_Y)/o_Y$ on $\pi_Y:X\t Y\ra Y$. Here $c_{\pi_X}$ depends only on $o_Y$, and is defined even if $X$ is not oriented, and changing the sign of $o_Y$ changes that of $c_{\pi_X}$. Similarly, $c_{\pi_Y}$ depends only on $o_X$, and is defined even if $Y$ is not oriented, and changing the sign of $o_X$ changes that of $c_{\pi_Y}$. 
\smallskip

\noindent{\bf(g)} If $(X,o_X),(Y,o_Y)$ are oriented objects in $\tManc$ with $\dim X=\dim Y$ then there is a natural orientation $o_X\amalg o_Y$ on $X\amalg Y$. We have open submanifolds $\io_X:X\hookra X\amalg Y$, $\io_Y:Y\hookra X\amalg Y$, and $(o_X\amalg o_Y)\vert_X=o_X$, $(o_X\amalg o_Y)\vert_Y=o_Y$.

\smallskip

\noindent{\bf(h)} For $X\in\tManc$, there is a natural coorientation $c_{i_X}$ for $i_X:\pd X\ra X$. Hence any orientation $o_X$ on $X$ induces an orientation $o_{\pd X}=c_{i_X}\cdot o_X$ on~$\pd X$.

If $(X,o_X),(Y,o_Y)$ are oriented objects in $\tManc$, equation \eq{kh3eq6} becomes
\begin{equation*}
\pd(X\t Y)\cong \bigl(\pd X\t Y\bigr)\amalg(-1)^{\dim X}\bigl(X\t\pd Y\bigr).
\end{equation*}

\noindent{\bf(i)} For any $X\in\tManc$, Assumption \ref{kh3ass3}(c) gives a diffeomorphism $\si_X:\pd^2X\ra\pd^2X$ with $i_X\ci i_{\pd X}\ci\si=i_X\ci i_{\pd X}:\pd^2X\ra X$. We have natural coorientations $c_{\si_X},c_{i_X},c_{i_{\pd X}}$ for $\si_X,i_X,i_{\pd X}$ by {\bf(e)},{\bf(h)}. So $c_{i_X}\ci c_{i_{\pd X}}\ci c_{\si_X}$ and $c_{i_X}\ci c_{i_{\pd X}}$ are coorientations for $i_X\ci i_{\pd X}\ci\si=i_X\ci i_{\pd X}$. We have~$c_{i_X}\ci c_{i_{\pd X}}\ci c_{\si_X}=-c_{i_X}\ci c_{i_{\pd X}}$.

Thus, if $X$ is oriented, so that $\pd X$ and $\pd^2X$ are also oriented by {\bf(h)}, then $\si:\pd^2X\ra\pd^2X$ is an orientation-reversing diffeomorphism.
\smallskip

\noindent{\bf(j)} For objects $X\in\Man\subset\cManc$, and for intervals $I$ in $\R$ and $k$-simplices $\De_k$ as objects of $\cManc$, orientations are naturally identified with orientations in $\Man$ and $\Manc$ in the usual sense in differential geometry.

For smooth maps $f:X\ra Y$ in $\Man\subset\tManc$, and the smooth maps in $\Manc$ in Assumption \ref{kh3ass4}(c)(A)--(D), coorientations in $\tManc$ are naturally identified with coorientations in $\Man,\Manc$ as usual in differential geometry.

These identifications of (co)orientations in $\Man,\Manc$ and $\tManc$ preserve all the operations and properties in {\bf(a)}--{\bf(i)} above.

In particular, for any $X\in\tManc$, an orientation $o_X$ on $X$ restricts to an orientation $o_X\vert_{X^\ci}$ on the (ordinary) manifold $X^\ci\in\Man$, in the usual sense of differential geometry. Also, if $X\in\tManc$ with $\dim X=0$, so that $X\in\Man$ and $X_\top$ is a set with the discrete topology as in Assumption \ref{kh3ass4}(e), then orientations on $X$ are identified with maps~$X_\top\ra\{\pm 1\}$.
\smallskip

\noindent{\bf(k)} On $\R^n$, and on intervals $I$ in $\R$, and on the simplex $\De_k$, we have {\it standard orientations}, such that $\d x_1\w\cdots\w\d x_n$ on $\R^n$, and $\d x$ on $I$, and $\d x_1\w\cdots\w\d  x_k$ on $\De_k$, are positive top-degree forms. So, for example, if the interval $[a,b]$ for $a<b$ has its standard orientation, then $\pd[a,b]=-\{a\}\amalg\{b\}$ in oriented 0-manifolds.

\smallskip

\noindent{\bf(l)} Let $g:X\ra Z,$ $h:Y\ra Z$ be transverse morphisms in $\tManc$ with $Z\in\Man\subset\tManc$, and $W=X\t_{g,Z,h}Y$ be the fibre product in $\tManc$ in Assumption \ref{kh3ass5}(c). Suppose $o_X,o_Y,o_Z$ are orientations on $X,Y,Z$, giving coorientations $c_g=o_X/o_Z$ on $g$ and $c_h=o_Y/o_Z$ on $h$. 

Then there is a natural orientation $o_W$ on $W$, giving natural coorientations $c_{\pi_X}=o_W/o_X$ on $\pi_X:W\ra X$ and $c_{\pi_Y}=o_W/o_Y$ on $\pi_Y:W\ra Y$.

Each of $o_W,c_{\pi_X},c_{\pi_Y}$ only depends on particular subsets of $o_X,o_Y,o_Z,c_g,c_h$, and is still well-defined if the rest of $o_X,o_Y,o_Z,c_g,c_h$ are not chosen, or do not exist. So, $o_W$ depends only on $o_X,o_Y,o_Z$, and also only on $o_X,c_h$, and also only on $o_Y,c_g$. Also $c_{\pi_X}$ depends only on $c_h$, and $c_{\pi_Y}$ depends only on $c_g$.

In each case, changing the sign of an element of the subset on which $o_W,c_{\pi_X}$ or $c_{\pi_Y}$ depends, changes the sign of $o_W,c_{\pi_X}$ or~$c_{\pi_Y}$.

If $X,Y\in\Man$, so that $W=X\t_{g,Z,h}Y$ is the usual transverse fibre product in $\Man$ by Assumption \ref{kh3ass5}(d), and (co)orientations on $W,X,Y,Z,g,h,\pi_X,\pi_Y$ are as usual in $\Man$ by {\bf(j)}, then $o_W,c_{\pi_X},c_{\pi_Y}$ are defined from $o_X,o_Y,o_Z,c_g,c_h$ as usual in differential geometry.
\smallskip

\noindent{\bf(m)} Equations \eq{kh3eq2}--\eq{kh3eq5} hold for oriented objects in $\tManc$, using the orientations on fibre products in {\bf(l)}, provided each fibre product in the formulae is transverse over an object in $\Man$, as in Assumption \ref{kh3ass5}(c). Corresponding identities hold for coorientations, and for combinations of orientations and coorientations.
\smallskip

\noindent{\bf(n)} Suppose $(X,o_X)$ is an oriented object in $\tManc$ with $\dim X=1$ and $X_\top$ compact. Then $\pd X$ is a compact 0-manifold, so $(\pd X)_\top$ is a finite set, and by part {\bf(h)} it has an orientation $o_{\pd X}=c_{i_X}\cdot o_X$, which as in part {\bf(j)} may be interpreted as a map $o_{\pd X}:(\pd X)_\top\ra\{\pm 1\}$. Then we have
\e
\sum_{x'\in (\pd X)_\top}o_{\pd X}(x')=0\qquad\text{in $\Z$.}
\label{kh3eq11}
\e

\label{kh3ass6}
\end{ass}

\begin{ass}{\bf(Sard's Theorem type conditions.)} {\bf(a)} Suppose $f:X\ra Y$ is a morphism in $\tManc$ with $Y\in\Man$ and $\dim X=m$, $\dim Y=n$ for $m\le n$. Then $\cH^n(f_\top(X_\top\sm X_\top^\ci))=0$ in~$Y$.
\smallskip

\noindent{\bf(b)} Suppose $f:X\ra Y$ is a morphism in $\tManc$ with $Y\in\Man$ and $\dim X=m$, $\dim Y=n$ for $m\ge n$. Then there is a subset $S\subseteq Y_\top$ with $\cH^n(S)=0$ in $Y$, such that if $y\in Y_\top\sm S$ then there is an open submanifold $i:U\hookra X$ as in Assumption \ref{kh3ass2}(d), such that $U_\top$ is an open neighbourhood of $f_\top^{-1}(y)$ in $X_\top$ and $f\ci i:U\ra X$ is a submersion, as in Assumption~\ref{kh3ass5}(a).
\label{kh3ass7}
\end{ass}

\begin{rem}{\bf(i)} We can modify the requirement in Assumption \ref{kh3ass4}(a) that $\tManc$ contains manifolds $\Man$ as a full subcategory. Sections \ref{kh35} and \ref{kh53} give a version of our theory in which $\Man$ is replaced by effective orbifolds~$\Orbeff$.

It would also be possible to produce a version of the theory not requiring morphisms $f:X\ra Y$ in $\tManc$ to give smooth maps between the interiors $X^\ci,Y^\ci$, for example, a version based on $C^k$ maps for~$k\ge 1$. 

Note however that Assumption \ref{kh3ass7}(b) is based on Theorem \ref{kh3thm1}(b), which requires $f:X\ra Y$ to be $C^k$ for $k\ge\dim X-\dim Y+1$. Taking $\dim X\gg 0$, we see that Assumption \ref{kh3ass7}(b) fails if $\tManc$ contains all $C^k$ maps between manifolds for any fixed $k\ge 1$. So we would need a substitute for Assumption~\ref{kh3ass7}(b).

\smallskip

\noindent{\bf(ii)} Defining the product orientations $o_X\t o_Y$ in Assumption \ref{kh3ass6}(f), and the coorientations $c_{i_X}$ in Assumption \ref{kh3ass6}(h), and $o_W,c_{\pi_X},c_{\pi_Y}$ in Assumption \ref{kh3ass6}(l), all require a choice of {\it orientation convention}. Ours follow Fukaya, Oh, Ohta and Ono \cite[Conv.~45.1]{FOOO1}. Different conventions would change the signs in various places, by factors depending on dimensions.
\smallskip

\noindent{\bf(iii)} Several of the assumptions above are crucial to our theory, but are used in only a few places, so it is easy to miss their importance. In particular, in the material of \S\ref{kh4} and the proofs in~\S\ref{kh7}:
\begin{itemize}
\setlength{\itemsep}{0pt}
\setlength{\parsep}{0pt}
\item Assumptions \ref{kh3ass3}(c) and \ref{kh3ass6}(i) are only used to prove $\pd\ci\pd=0$ in $MC_*(Y;R)$ in \S\ref{kh41} and $\d\ci\d=0$ in $MC^*(Y;R)$ in~\S\ref{kh42}.
\item Assumption \ref{kh3ass6}(n) is only used in the proofs of Proposition \ref{kh4prop1} and Theorem \ref{kh4thm2} in \S\ref{kh71} and~\S\ref{kh74}.
\item Assumption \ref{kh3ass7}(a) is only used in the proof of Theorem \ref{kh4thm2} in~\S\ref{kh74}.
\item Assumption \ref{kh3ass7}(b) is only used in in the proofs of Proposition \ref{kh4prop1} and Theorem \ref{kh4thm2} in \S\ref{kh71} and~\S\ref{kh74}.
\end{itemize}
\label{kh3rem2}
\end{rem}

\begin{ex} Here are some categories satisfying Assumptions \ref{kh3ass1}--\ref{kh3ass7}:
\begin{itemize}
\setlength{\itemsep}{0pt}
\setlength{\parsep}{0pt}
\item[(a)] Manifolds with corners and smooth maps $\Manc$ from~\S\ref{kh31}.
\item[(b)] Manifolds with corners and strongly smooth maps $\Mancst$ from~\S\ref{kh31}.
\item[(c)] Manifolds with corners and weakly smooth maps $\Mancwe$ from~\S\ref{kh31}.
\item[(d)] The category $\Mangc$ of manifolds with generalized corners from~\cite{Joyc8}.
\item[(e)] Fix $k\ge 1$. Define categories $\Man^{\bf c}_k$ and $\Man^{\bf c}_{{\bf st},k}$  to have objects manifolds with corners $X,Y$, and morphisms continuous maps $f:X\ra Y$ which are $C^k$, or strongly $C^k$, respectively, in the obvious $C^k$ generalization of smooth and strongly smooth maps from \S\ref{kh31}, such that in addition $f\vert_{S^l(X)}:S^l(X)\ra Y$ is $C^\iy$ for each $l=0,\ldots,\dim X$, using the depth stratification $X=\coprod_{l=0}^{\dim X}S^l(X)$ from~\S\ref{kh31}. 

That is, morphisms $f:X\ra Y$ are smooth on each stratum $S^l(X)$, but only $C^k$ over the transitions between different strata $S^l(X),S^m(X)$. We do not define a weakly smooth version $\Man^{\bf c}_{{\bf we},k}$, since weakly $C^k$ maps $f:X\ra Y$ need not preserve the stratifications $X=\coprod_{l\ge 0}S^l(X)$, $Y=\coprod_{m\ge 0}S^m(Y)$, and so the condition $f\vert_{S^l(X)}$ smooth would not be preserved under composition.

Actually, it would be better to define objects $X$ in $\Man^{\bf c}_k,\Man^{\bf c}_{{\bf st},k}$ not as manifolds with corners, but as being covered by an atlas $\{(U_a,\phi_a):a\in A\}$ as in \S\ref{kh31} in which the transition functions $\psi^{-1}\ci\phi$ are required to be $C^k$ diffeomorphisms which are smooth on each stratum.

We do not allow $k=0$, as this causes problems with transverse fibre products in Assumption \ref{kh3ass5}(c), and with Assumption~\ref{kh3ass7}(b).
\end{itemize}
\label{kh3ex1}
\end{ex}

As in Remark \ref{kh3rem3}(a) below, the categories of effective orbifolds with corners in Example \ref{kh3ex2} also satisfy Assumptions~\ref{kh3ass1}--\ref{kh3ass7}.

\subsection{Assumptions for de Rham M-(co)homology}
\label{kh34}

Section \ref{kh52} defines variants of M-(co)homology called {\it de Rham M-(co)homology\/} $MH_*^\dR(Y;\R),MH^*_\dR(Y;\R)$, which include exterior forms as in ordinary de Rham cohomology. To do this we will need the following additional assumption about exterior forms on `manifolds with corners'. 

Suppose Assumptions \ref{kh3ass1}--\ref{kh3ass7} hold for the category~$\tManc$.

\begin{ass}{\bf(a)} Let $X$ be an object in $\tManc$. Then we are given real vector spaces $\Om^k(X)$ for $k=0,1,\ldots,$ with $\Om^k(X)=0$ for $k>\dim X$. Elements $\al$ of $\Om^k(X)$ are called $k$-{\it forms on\/} $X$. The {\it degree\/} of $\al\in\Om^k(X)$ is~$\deg\al=k$.

There are $\R$-linear {\it exterior derivatives\/} $\d:\Om^k(X)\ra\Om^{k+1}(X)$ for $k=0,1,\ldots,$ with $\d\ci\d=0$.

There are $\R$-bilinear {\it exterior products\/} $\w:\Om^k(X)\t\Om^l(X)\ra\Om^{k+l}(X)$ with $\al\w\be=(-1)^{kl}\be\w\al$ and $(\al\w\be)\w\ga=\al\w(\be\w\ga)$ and $\d(\al\w\be)=(\d\al)\w\be+(-1)^k\al\w\d\be$ for all $\al\in\Om^k(X)$, $\be\in\Om^l(X)$, and $\ga\in\Om^m(X)$.

There is a natural element $1_X\in\Om^0(X)$ with $\d 1_X=0$ and $1_X\w\al=\al\w 1_X=\al$ for all $\al\in\Om^k(X)$.
\smallskip

\noindent{\bf(b)} Let\/ $f:X\ra Y$ be a morphism in $\tManc$. Then there are $\R$-linear {\it pullback maps\/} $f^*:\Om^k(Y)\ra\Om^k(X)$ for all $k=0,1,\ldots,$ with $\d\ci f^*=f^*\ci\d:\Om^k(Y)\ra\Om^{k+1}(X)$, and $f^*(\al\w\be)=f^*(\al)\w f^*(\be)$ for all $\al,\be\in\Om^*(Y)$, and $f^*(1_Y)=1_X$. These are functorial, i.e.\ $(g\ci f)^*=f^*\ci g^*$ and $\id_X^*=\id_{\Om^k(X)}$. When $i:U\hookra X$ is inclusion of an open submanifold in $\tManc$ and $\al\in\Om^k(X)$, we often write $\al\vert_U$ instead of $i^*(\al)$ in $\Om^k(U)$.
\smallskip

\noindent{\bf(c)} Each $\al\in\Om^k(X)$ has a ({\it closed\/}) {\it support\/} $\supp\al$, a closed subset of $X_\top$, with the property that $\al\vert_U=0$ if and only if $U_\top\cap\supp\al=\es$ for all open submanifolds $U\subseteq X$. 

We call $\al$ {\it compactly-supported\/} if $\supp\al$ is compact. We have $\supp(\d\al)\subseteq\supp\al$, and $\supp(\al\w\be)\subseteq\supp\al\cap\supp\be$, and~$\supp(f^*(\al))\subseteq f_\top^{-1}(\supp\al)$.
\smallskip

\noindent{\bf(d)} If $X\in\Man\subset\tManc$ then $\Om^k(X)=C^\iy(\La^kT^*X)$, the usual vector space of smooth $k$-forms on $X$, and $\d:\Om^k(X)\ra\Om^{k+1}(X)$ is the usual de Rham exterior derivative on $k$-forms, and $\w$ is the usual exterior product, and $1_X$ is the constant function $x\mapsto 1$ in $\Om^0(X)=C^\iy(X)$. Thus, for any $X\in\tManc$ and $\al\in\Om^k(X)$, the restriction $\al\vert_{X^\ci}$ is an ordinary $k$-form on the manifold~$X^\ci$.

\smallskip

\noindent{\bf(e)} If $X\in\tManc$ with $\dim X=n$ and $\al\in\Om^n(X)$ is compactly-supported, then $\al\vert_{X^\ci}$ is an integrable (i.e. $L^1$) $n$-form on the $n$-manifold $X^\ci$, even though $\al\vert_{X^\ci}$ will not be compactly-supported if $\supp\al\cap X^\ci_\top$ is not compact.

Thus, if $X$ is oriented, so that $X^\ci$ is oriented, then $\int_{X^\ci}\al\vert_{X^\ci}$ exists in $\R$, as the integral of an integrable $n$-form on an oriented $n$-manifold. We write
\begin{equation*}
\int_X\al:=\int_{X^\ci}\al\vert_{X^\ci}\in\R.
\end{equation*}
This is a definition of $\int_X\al$, motivated by supposing $X\sm X^\ci$ has `measure zero'.

\smallskip

\noindent{\bf(f) (Pushforwards.)} Let $f:X\ra Y$ be a cooriented submersion in $\tManc$ with $Y\in\Man$, and $\dim X=m$, $\dim Y=n$ with $m\ge n$. Suppose $\al\in\Om^k(X)$ for $k\ge m-n$ with $f_\top\vert_{\supp\al}:\supp\al\ra Y_\top$ proper. Then there is a natural {\it pushforward\/} $f_*(\al)\in\Om^{k-m+n}(Y)$, which is characterized by the property that if $U\subseteq Y$ is open and oriented, and $\eta\in\Om^{m-k}(U)$ with $\supp\eta$ compact, then
\e
\int_U\eta\w f_*(\al)=\int_{f^{-1}(U)}f^*(\eta)\w\al,
\label{kh3eq12}
\e
where $f^{-1}(U)\subseteq X$ is an open submanifold with orientation determined by combining the coorientation on $f\vert_{f^{-1}(U)}:f^{-1}(U)\ra U$ and the orientation on $U$, and $\eta\w f_*(\al),f^*(\eta)\w\al$ are both compactly-supported as $\supp\eta$ is compact and $f_\top\vert_{\supp\al}:\supp\al\ra Y_\top$ is proper.
\smallskip

\noindent{\bf(g) (Stokes' Theorem.)} Suppose $X\in\tManc$ is oriented with $\dim X=n$, and $\al\in\Om^{n-1}(X)$ is compactly-supported. This implies that $\pd X$ is oriented with $\dim\pd X=n-1$, and $i_X^*(\al)\in\Om^{n-1}(\pd X)$, $\d\al\in\Om^n(X)$ are compactly-supported, as $\supp [i_X^*(\al)]\subseteq i_{X,\top}^{-1}(\supp\al)$ and $\supp\d\al\subseteq\supp\al$ by {\bf(c)}. Then
\begin{equation*}
\int_X\d\al=\int_{\pd X}i_X^*(\al).
\end{equation*}

\label{kh3ass8}
\end{ass}

Assumption \ref{kh3ass8} holds for all of Example \ref{kh3ex1}(a)--(e). Note that Assumption \ref{kh3ass8} implies Assumption \ref{kh3ass6}(n), since if $(X,o_X)$ is a compact oriented object in $\tManc$ with $\dim X=1$ then
\begin{equation*}
\sum_{x'\in (\pd X)_\top}o_{\pd X}(x')=\int_{\pd X}1_{\pd X}=\int_{\pd X}i_X^*(1_X)=\int_X\d 1_X=\int_X0=0.
\end{equation*}

Let $f:X\ra Y$ and $\al\in\Om^k(X)$ be as in Assumption \ref{kh3ass8}(f). Then one can show using Assumption \ref{kh3ass8}(f),(g) that
\begin{equation*}
\d(f_*(\al))=f_*(\d\al)+(-1)^{\dim X-k}(f\ci i_X)_*(i_X^*(\al)).
\end{equation*}

\subsection{Assumptions for the extension to orbifolds}
\label{kh35}

Section \ref{kh29} discussed orbifolds, and defined a category $\Orbeff$ of effective orbifolds in which we will work. In \S\ref{kh53} we will extend the M-(co)homology theories $MH_*(Y;R),MH^*(Y;R),\ldots,MH^*_\dR(Y;\R)$ of \S\ref{kh4}--\S\ref{kh52} from manifolds $Y$ to effective orbifolds. To do this, we now state assumptions on a category $\tOrbeffc$ of `effective orbifolds with corners', which is to replace the category $\tManc$ of `manifolds with corners' in Assumptions \ref{kh3ass1}--\ref{kh3ass7} and \ref{kh3ass8}. Note that we want $\tOrbeffc$ to be a category, not a 2-category. It should contain the category $\Orbeff$ of \S\ref{kh291}, just as $\tManc$ should contain $\Man$, and have `submersions' and `transverse fibre products' generalizing those in $\Orbeff$ discussed in \S\ref{kh292}.

\begin{ass}{\bf(a)} In Assumption \ref{kh3ass1}(a), in place of the category $\tManc$, we are given a category $\tOrbeffc$. Objects $X$ in $\tOrbeffc$ will be called {\it effective orbifolds with corners\/} (although they may in examples not be orbifolds, but some kind of singular space), and morphisms $f:X\ra Y$ in $\tOrbeffc$ will be called {\it smooth maps} (although they may in examples be non-smooth). Isomorphisms in $\tOrbeffc$ will be called {\it diffeomorphisms}.

Throughout Assumptions \ref{kh3ass1}--\ref{kh3ass7} and \ref{kh3ass8}, we replace $\tManc$ by $\tOrbeffc$ and `manifolds' (with boundary or corners) by `effective orbifolds' (with boundary or corners). The rest of Assumptions \ref{kh3ass1}--\ref{kh3ass3} extend to orbifolds with only trivial changes.
\smallskip

\noindent{\bf(b)} In Assumption \ref{kh3ass4}(a), the category of effective orbifolds $\Orbeff$ from \S\ref{kh291} is included as a full subcategory~$\Orbeff\subset\tOrbeffc$.

In Assumption \ref{kh3ass4}(b), for each object $X\in\tOrbeffc$ there is a natural open suborbifold $X^\ci\hookra X$ called the {\it interior\/} of $X$, with $X^\ci\in\Orbeff\subset\tOrbeffc$. 

Throughout Assumptions \ref{kh3ass4}--\ref{kh3ass7} and \ref{kh3ass8}, objects $X\in\Man$ should be replaced by objects~$X\in\Orbeff$.

As in Assumption \ref{kh3ass4}(e), we require that if $X\in\tOrbeffc$ with $\dim X=0$ then $X\in\Orbeff$. But a 0-dimensional effective orbifold is a manifold, so $X\in\Man\subset\Orbeff\subset\tOrbeffc$, and $X_\top$ is a set with the discrete topology.

The rest of Assumption \ref{kh3ass4} extends to orbifolds with only trivial changes.
\smallskip

\noindent{\bf(c)} Assumption \ref{kh3ass5} extends to orbifolds with only trivial changes, but as this is surprising, we discuss it anyway.

In Assumption \ref{kh3ass5}(a), if $f:X\ra Y$ is a morphism in $\Orbeff$ then $f$ is a submersion in $\tOrbeffc$ if and only if it is a submersion in $\Orbeff$, in the sense of Definition~\ref{kh2def14}.

In Assumption \ref{kh3ass5}(c), we are given a notion of when morphisms $g:X\ra Z,$ $h:Y\ra Z$ in $\tOrbeffc$ with $Z\in\Orbeff\subset\tOrbeffc$ are {\it transverse}.

We require that if $g,h$ are transverse then a fibre product $W=X\t_{g,Z,h}Y$ exists in $\tOrbeffc$ in the sense of category theory, with $\dim W=\dim X+\dim Y-\dim Z$, in a Cartesian square~\eq{kh3eq7}.

The functor $F_\tOrbeffc^\Top:\tOrbeffc\ra\Top$ should take such transverse fibre products in $\tOrbeffc$ to fibre products in $\Top$, so that \eq{kh3eq8} holds.

In Assumption \ref{kh3ass5}(d), if $X,Y,Z\in\Orbeff\subset\tOrbeffc$ then transversality of $g,h$ in $\tOrbeffc$ should be as for $\Orbeff$ in Definition \ref{kh2def14}, and the fibre product $X\t_{g,Z,h}Y$ in $\tOrbeffc$ should lie in $\Orbeff\subset\tOrbeffc$, and be as in Theorem~\ref{kh2thm5}.

Now, much of this may be worrying to those familiar with fibre products of orbifolds. With the {\it natural\/} definition of transverse morphisms of orbifolds:
\begin{itemize}
\setlength{\itemsep}{0pt}
\setlength{\parsep}{0pt}
\item In any ordinary category of orbifolds $\Orb$, some transverse fibre products $X\t_{g,Z,h}Y$ cannot exist, in the sense of category theory.
\item To make all transverse fibre products exist, we have to work in a 2-category of orbifolds $\mathfrak{Orb}$, with fibre products in the 2-category sense.
\item The (2-)functor $F_{\mathfrak{Orb}}^\Top:\mathfrak{Orb}\ra\Top$ does not take some transverse fibre products in $\mathfrak{Orb}$ to fibre products in~$\Top$.
\end{itemize}

The explanation, as in \S\ref{kh292}, is that we do not use the natural definition of transverse morphisms of effective orbifolds. Instead we use more restrictive notions of submersions and transverse morphisms in Definition \ref{kh2def14}, including extra conditions on actions of morphisms on orbifold groups. Such restricted transverse fibre products do exist in the ordinary category $\Orbeff$, and map to fibre products in $\Top$. In $\tOrbeffc$, we should include the same extra conditions in the definitions of submersions and transverse morphisms.

Just within $\Orbeff$, the notions of submersions and transverse morphisms in Definition \ref{kh2def14} do satisfy Assumption \ref{kh3ass5}, so there is no contradiction.
\smallskip

\noindent{\bf(d)} Assumption \ref{kh3ass6} extends to orbifolds with only trivial changes.

For Assumption \ref{kh3ass6}(n), we note as in {\bf(b)} that a 0-dimensional effective orbifold is a manifold, so $\pd X$ is a manifold, and \eq{kh3eq11} does not need to be modified to take orbifold groups into account, as one might guess.
\smallskip

\noindent{\bf(e)} Assumption \ref{kh3ass7} extends to orbifolds with only trivial changes, interpreting Hausdorff measure in orbifolds in the obvious way.
\smallskip

\noindent{\bf(f)} Assumption \ref{kh3ass8} extends to orbifolds with only trivial changes, provided we interpret `exterior forms on orbifolds' and `integration over orbifolds' correctly.

To make these interpretations clear, consider a global quotient orbifold $X=[V/G]$, for $V$ a manifold and $G$ a finite group acting effectively on $V$ by diffeomorphisms. Then in the usual way we have $k$-forms $\Om^k(V)=C^\iy(\La^kT^*V)$ on $V$ for $k=0,1,\ldots$ and exterior derivative $\d:\Om^k(V)\ra\Om^{k+1}(V)$. 

The action of $G$ on $V$ induces a linear action of $G$ on $\Om^k(V)$ by pullbacks. The projection $\pi:V\ra X=[V/G]$ induces a pullback morphism $\pi^*:\Om^k(X)\ra\Om^k(V)$, which is an isomorphism
\begin{equation*}
\pi^*:\Om^k(X)\,{\buildrel\cong\over\longra}\,\Om^k(V)^G
\end{equation*}
with the $\R$-vector subspace $\Om^k(V)^G\subseteq\Om^k(V)$ of $G$-invariant $k$-forms on $V$. So we can identify $k$-forms on $X=[V/G]$ with $G$-invariant $k$-forms on $V$.

Suppose $\dim V=\dim X=m$, and $V$ has a $G$-invariant orientation, which descends to an orientation on $X=[V/G]$. Let $\om$ be a compactly-supported $m$-form on $X$, so that $\pi^*(\om)$ is a $G$-invariant $m$-form on $V$. Then integration over $X=[V/G]$ may be defined by
\begin{equation*}
\int_X\om=\frac{1}{\md{G}}\int_V\pi^*(\om).
\end{equation*}

\label{kh3ass9}
\end{ass}

To summarize Assumption \ref{kh3ass9}: in Assumptions \ref{kh3ass1}--\ref{kh3ass7} and \ref{kh3ass8} we replace `manifolds' by `effective orbifolds' throughout, replace $\Man$ by the category $\Orbeff$ in \S\ref{kh291}, and use the definitions of submersions and transverse morphisms in $\Orbeff$ in \S\ref{kh292}, and no other important changes are necessary.

\begin{ex} As in Remark \ref{kh2rem5}(d), by combining Definition \ref{kh2def5} for the category $\Orbeff$ of effective orbifolds, and Definitions \ref{kh3def1}--\ref{kh3def2} for the category $\Manc$ of manifolds with corners, we may define a category $\Orbeffc$ of {\it effective orbifolds with corners}, which satisfies Assumption~\ref{kh3ass9}.

More generally, to each of the categories $\Manc,\Mancst,\ab\Mancwe,\ab\Mangc,\ab\Man^{\bf c}_k,\ab\Man^{\bf c}_{k,{\bf st}}$ of `manifolds with corners' in Example \ref{kh3ex1} satisfying Assumptions \ref{kh3ass1}--\ref{kh3ass7} and \ref{kh3ass8}, there is a corresponding category of `effective orbifolds with corners' $\Orbeffc,\Orb^{\bf c}_{\bf eff,st},\Orb^{\bf c}_{\bf eff,we},\Orb^{\bf gc}_{\bf eff},\Orb^{\bf c}_{{\bf eff},k},\Orb^{\bf c}_{{\bf eff},{\bf st},k}$ which satisfies Assumption~\ref{kh3ass9}.
\label{kh3ex2}
\end{ex}

\begin{rem}{\bf(a)} Any category $\tOrbeffc$ which satisfies Assumption \ref{kh3ass9}, also satisfies Assumptions \ref{kh3ass1}--\ref{kh3ass7} and \ref{kh3ass8} (setting $\tManc=\tOrbeffc$), since all the claims in Assumption \ref{kh3ass9} about $\Orbeff\subset\tOrbeffc$ can be restricted to $\Man\subset\Orbeff\subset\tOrbeffc$, where they become equivalent to Assumptions \ref{kh3ass1}--\ref{kh3ass7} and \ref{kh3ass8}. Thus, we can use the categories in Example \ref{kh3ex2} as the starting point for (co)homology theories of manifolds, as in \S\ref{kh4} and~\S\ref{kh51}--\S\ref{kh52}.
\smallskip

\noindent{\bf(b)} In Remark \ref{kh2rem5}(c) we discussed more sophisticated, and better behaved, definitions of (2-)categories of orbifolds than our $\Orbeff$. There is little written on orbifolds with boundary or with corners; the only foundational work the author knows is \cite[\S 1.12, \S 8.5--\S 8.9]{Joyc6}, which uses the $C^\iy$-stack setting. However, several of the approaches to orbifolds discussed in Remark \ref{kh2rem5}(c) extend automatically to orbifolds with boundary or corners, including those of Moerdijk and Pronk \cite{Moer,MoPr,Pron} and the author~\cite[\S 4.5]{Joyc7}. 

If we define a (2-)category of (effective) orbifolds with corners $\Orbc$ in this way, it will not satisfy Assumption \ref{kh3ass9}, because the (1-)morphisms $f:X\ra Y$ in $\Orbc$ will be continuous maps plus extra data, but $\tOrbeffc$ in Assumption \ref{kh3ass9} is required to contain $\Orbeff$, in which morphisms are continuous maps without extra data. However, any such (2-)category should admit a truncation functor, forgetting the extra data, to a category $\tOrbeffc$ satisfying Assumption~\ref{kh3ass9}.
\label{kh3rem3}
\end{rem}

\section{Integral M-homology and M-cohomology}
\label{kh4}

We now define and study our main new (co)homology theories of a smooth manifold $Y$ and a commutative ring $R$: {\it M-homology\/} $MH_*(Y;R)$, {\it M-cohomology\/} $MH^*(Y;R)$, {\it compactly-supported M-cohomology\/} $MH^*_\cs(Y;R)$, and {\it locally finite M-homology\/} $MH_*^\lf(Y;R)$. Here the M- is short for `Manifold', as the (co)chains $[V,n,s,t]$ used to define the (co)homology theories involve a `manifold with corners' $V$ in $\tManc$ with `smooth maps' $s:V\ra\R^n$ and $t:V\ra Y$.

Throughout we take $\tManc$ to be a category of `manifolds with corners' satisfying Assumptions \ref{kh3ass1}--\ref{kh3ass7} of \S\ref{kh33}. We can choose $\tManc$ to be conventional manifolds with corners $\Manc$ in \S\ref{kh31}, but there are also other possibilities involving singular manifolds and non-smooth maps.

In \S\ref{kh5} we will describe several variations on these, including {\it rational M-homology\/} $MH_*^\Q(Y;R)$ and {\it rational M-cohomology\/} $MH^*_\Q(Y;R)$, which require $R$ to be a $\Q$-algebra, but have slightly better properties at the (co)chain level (e.g.\ the cup product is supercommutative on $MC^*_\Q(Y;R)$). When we wish to distinguish the theories of this section from these others we will call them {\it integral M-homology\/} and {\it integral M-cohomology}, since they work for~$R=\Z$.

Most proofs are postponed to~\S\ref{kh7}. 

\subsection{\texorpdfstring{M-homology $MH_*(Y;R)$}{M-homology}}
\label{kh41}

For the whole of \S\ref{kh4}, fix a commutative ring $R$, and a category $\tManc$ satisfying Assumptions \ref{kh3ass1}--\ref{kh3ass7} of \S\ref{kh33}. For simplicity, objects $X$ in $\tManc$ will be called {\it manifolds with corners\/} (although they may in examples not be manifolds, but some kind of singular space), and morphisms $f:X\ra Y$ in $\tManc$ will be called {\it smooth maps\/} (although they may in examples be non-smooth). We use the notation of \S\ref{kh33} throughout.

When we say that $Y$ is a {\it manifold\/} (rather than a manifold with corners), we mean an ordinary smooth manifold without boundary. When we say that $f:Y_1\ra Y_2$ is a {\it smooth map of manifolds\/} (rather than of manifolds with corners, or just a smooth map), we mean an ordinary smooth map of smooth manifolds without boundary. Since $\Man\subset\tManc$ is a full subcategory as in Assumption \ref{kh3ass4}, such $Y$ and $f:Y_1\ra Y_2$ are also objects and morphisms in~$\tManc$. 

\begin{dfn} Let $Y$ be a manifold. Consider quadruples $(V,n,s,t)$, where $V$ is an oriented manifold with corners (i.e.\ a pair $(V,o_V)$ with $V$ an object in $\tManc$ and $o_V$ an orientation on $V$, usually left implicit), and $n=0,1,\ldots,$ and $s:V\ra\R^n$ is a smooth map (morphism in $\tManc$) which is proper over an open neighbourhood of 0 in $\R^n$ (i.e. the continuous map $s:V\ra\R^n$ is proper near 0 in $\R^n$), and $t:V\ra Y$ is a smooth map (morphism in~$\tManc$).

Define an equivalence relation $\sim$ on such quadruples by $(V,n,s,t)\sim(V',\ab n',\ab s',\ab t')$ if $n=n'$, and there exists an orientation-preserving diffeomorphism $f:V\ra V'$ with $s=s'\ci f$ and $t=t'\ci f$. Write $[V,n,s,t]$ for the $\sim$-equivalence class of $(V,n,s,t)$. We call $[V,n,s,t]$ a {\it generator}.

For each $k\in\Z$, define the {\it M-chains\/} $MC_k(Y;R)$ to be the $R$-module generated by such $[V,n,s,t]$ with $\dim V=n+k$, subject to the relations:
\begin{itemize}
\setlength{\itemsep}{0pt}
\setlength{\parsep}{0pt}
\item[(i)] For each generator $[V,n,s,t]$ and each $i=0,\ldots,n$ we have
\e
[V,n,s,t]=(-1)^{n-i}[V\t\R,n+1,s',t\ci\pi_V]\quad\text{in $MC_k(Y;R)$,}
\label{kh4eq1}
\e
where writing $s=(s_1,\ldots,s_n):V\ra\R^n$ with $s_j:V\ra\R$ for $j=1,\ldots,n$ and $\pi_V:V\t\R\ra V$, $\pi_\R:V\t\R\ra\R$ for the projections, then 
\begin{equation*}
s'=(s_1\ci\pi_V,\ldots,s_i\ci\pi_V,\pi_\R,s_{i+1}\ci\pi_V,\ldots,s_n\ci\pi_V):V\t\R\longra\R^{n+1},
\end{equation*}
and $V\t\R$ has the product orientation from Assumption \ref{kh3ass6}(f) of the given orientation on $V$ and the standard orientation on~$\R$.
\item[(ii)] Let $I$ be a finite indexing set, $a_i\in R$ for $i\in I$, and $[V_i,n,s_i,t_i]$, $i\in I$ be generators for $MC_k(Y;R)$, all with the same $n$. Suppose there exists an open neighbourhood $X$ of $0$ in $\R^n$, such that $s_i:V_i\ra\R^n$ is proper over $X$ for all $i\in I$, and the following condition holds:
\begin{itemize}
\setlength{\itemsep}{0pt}
\setlength{\parsep}{0pt}
\item[$(*)$] Suppose $(x,y)\in X\t Y$, such that for all $i\in I$ and $v\in V_i$ with $(s_i,t_i)(v)=(x,y)$, we have that $v\in V_i^\ci$ and 
\begin{equation*}
T_v(s_i,t_i):T_vV_i^\ci\longra T_xX\op T_yY
\end{equation*}
is injective, noting that $(s_i,t_i)\vert_{V_i^\ci}:V_i^\ci\ra X\t Y$ is a smooth map of (ordinary) manifolds. This implies that $(s_i,t_i)\vert_{V_i^\ci}$ is an embedding near $v\in V_i^\ci$. Hence $(s_i,t_i):V_i\ra X\t Y$ is injective near each $v$ in $(s_i,t_i)^{-1}(x,y)$, so $(s_i,t_i)^{-1}(x,y)$ has the discrete topology, and thus is finite as $s_i$ is proper over $X$. Note too that $V_i^\ci$ is an oriented manifold by Assumption \ref{kh3ass6}(j) with $\dim V_i^\ci=n+k$, so $T_vV_i^\ci$ is an oriented vector space of dimension $n+k$. We require that for all oriented $(n+k)$-planes $P\subseteq T_xX\op T_yY=\R^n\op T_yY$, we have
\e
\begin{split}
&\sum_{\begin{subarray}{l} i\in I,\; v\in V_i^\ci:(s_i,t_i)(v)=(x,y),\;  T_v(s_i,t_i)[T_vV_i^\ci]=P \\ \text{$T_v(s_i,t_i):T_vV_i^\ci\,{\buildrel\cong\over\longra}\,P$ is orientation-preserving}\end{subarray}\!\!\!\!\!\!\!} a_i=\\
&\sum_{\begin{subarray}{l} i\in I,\; v\in V_i^\ci:(s_i,t_i)(v)=(x,y),\; T_v(s_i,t_i)[T_vV_i^\ci]=P \\ \text{$T_v(s_i,t_i):T_vV_i^\ci\,{\buildrel\cong\over\longra}\,P$ is orientation-reversing}\end{subarray}\!\!\!\!\!\!\!} a_i\qquad\text{in $R$.}
\end{split}
\label{kh4eq2}
\e
\end{itemize}
Then
\begin{equation*}
\sum_{i\in I}a_i\,[V_i,n,s_i,t_i]=0\qquad\text{in $MC_k(Y;R)$.}
\end{equation*}
\end{itemize}

Define $\pd:MC_k(Y;R)\ra MC_{k-1}(Y;R)$ to be the unique $R$-linear map with
\e
\pd[V,n,s,t]=[\pd V,n,s\ci i_V,t\ci i_V],
\label{kh4eq3}
\e
for all generators $[V,n,s,t]$, where the orientation on $\pd V$ is induced from that of $V$ as in Assumption \ref{kh3ass6}(h). We show this is well defined in Proposition~\ref{kh4prop1}.

For any generator $[V,n,s,t]$ in $MC_k(Y;R)$ we have
\begin{equation*}
\pd\ci\pd[V,n,s,t]=[\pd^2V,n,s\ci i_V\ci i_{\pd V},t\ci i_V\ci i_{\pd V}].
\end{equation*}
Apply (ii) with $[\pd^2V,n,s\ci i_V\ci i_{\pd V},t\ci i_V\ci i_{\pd V}]$ in place of $\sum_{i\in I}a_i\,[V_i,n,s_i,t_i]$. Assumptions \ref{kh3ass3}(c) and \ref{kh3ass6}(i) give an orientation-reversing diffeomorphism $\si_V:\pd^2V\ra\pd^2V$ with $\si_V^2=\id_{\pd^2V}$. In \eq{kh4eq2}, each contribution to the top line from $v''$ in $\pd^2V$ with $T_{v''}[(s,t)\ci i_V\ci i_{\pd V}]$ orientation-preserving, is exactly matched by a contribution to the bottom line from $\si_V(v'')$ in $\pd^2V$ with $T_{\si_V(v'')}[(s,t)\ci i_V\ci i_{\pd V}]$ orientation-reversing. Thus $(*)$ holds for $[\pd^2V,n,s\ci i_V\ci i_{\pd V},t\ci i_V\ci i_{\pd V}]$, so that (ii) gives~$\pd\ci\pd[V,n,s,t]=0$. 

Hence $\pd\ci\pd=0:MC_k(Y;R)\ra MC_{k-2}(Y;R)$ for all $k$, and $\bigl(MC_*(Y;R),\pd\bigr)$ is a chain complex. Define the {\it M-homology groups\/} (or {\it integral M-homology groups\/}) $MH_*(Y;R)$ to be the homology of this chain complex. That is, for $k\in\Z$ we define $R$-modules
\e
MH_k(Y;R)=\frac{\ts \Ker\bigl(\pd: MC_k(Y;R)\longra MC_{k-1}(Y;R)\bigr)}{\ts \Im\bigl(\pd:MC_{k+1}(Y;R)\longra MC_k(Y;R)\bigr)}\,.
\label{kh4eq4}
\e

If $Y$ is compact and oriented with $\dim Y=m$, define the {\it fundamental cycle\/} $[Y]=[Y,0,0,\id_Y]\in MC_m(Y;R)$. Here $V=Y$ has the given orientation, and $s=0:V\ra\R^0$ is proper as $Y$ is compact. We have $\pd[Y]=0$ as $\pd Y=\es$, so passing to homology gives the {\it fundamental class\/}~$[[Y]]\in MH_m(Y;R)$.
\label{kh4def1}
\end{dfn}

\begin{rem} If the base ring $R$ has characteristic 2, that is, $1+1=0$ in $R$, then we can omit orientations $o_V$ on $V$ in generators $[V,n,s,t]$ in Definition \ref{kh4def1}, and allow $V$ to be unoriented. We must impose an extra condition on $\tManc$, an unoriented version of Assumption \ref{kh3ass6}(n), saying that if $X\in\tManc$ is compact with $\dim X=1$, then the number of points in $\pd X$ is zero modulo~2.

Similarly, if $R$ has characteristic 2 then we can omit orientations and coorientations throughout \S\ref{kh4}, though we will not mention this again.
\label{kh4rem1}
\end{rem}

The next proposition will be proved in \S\ref{kh71}.

\begin{prop} $\pd:MC_k(Y;R)\ra MC_{k-1}(Y;R)$ above is well defined.
\label{kh4prop1}
\end{prop}

\begin{rem} The most nontrivial part of Definition \ref{kh4def1} is relation (ii), and readers are advised to study it carefully before proceeding further.

Two easy consequences of (ii) are that if $[V,n,s,t]$ is a generator then
\e
[-V,n,s,t]=-[V,n,s,t],
\label{kh4eq5}
\e
where $-V$ is $V$ with the opposite orientation, and if $[V_1,n,s_1,t_1]$, $[V_2,n,s_2,t_2]$ are generators with the same $n$ then
\e
[V_1\amalg V_2,n,s_1\amalg s_2,t_1\amalg t_2]=[V_1,n,s_1,t_1]+[V_2,n,s_2,t_2].
\label{kh4eq6}
\e
We prove these using (ii) with $\sum_{i\in I}a_i\,[V_i,n,s_i,t_i]=[V,n,s,t]+[-V,n,s,t]$ for \eq{kh4eq5} and 
$[V_1\amalg V_2,n,s_1\amalg s_2,t_1\amalg t_2]-[V_1,n,s_1,t_1]-[V_2,n,s_2,t_2]$ for~\eq{kh4eq6}.
\label{kh4rem2}
\end{rem}

\begin{lem} For any manifold\/ $Y$ we have $MC_k(Y;R)=0$ for\/ $k>\dim Y,$ so that\/ $MH_k(Y;R)=0$ for\/ $k>\dim Y$.
\label{kh4lem1}
\end{lem}

\begin{proof} Suppose $\sum_{i\in I}a_i\,[V_i,n,s_i,t_i]\in MC_k(Y;R)$ for $k>\dim Y$. Then condition $(*)$ in Definition \ref{kh4def1}(ii) is trivial, as there are no $(n+k)$-planes $P\subseteq\R^n\op T_yY$ since $\dim P=n+k>n+\dim Y=\dim(\R^n\op T_yY)$. Thus Definition \ref{kh4def1}(ii) gives $\sum_{i\in I}a_i\,[V_i,n,s_i,t_i]=0$, and~$MC_k(Y;R)=0$.
\end{proof}

\begin{dfn} Let $f:Y_1\ra Y_2$ be a smooth map of manifolds. Define the {\it pushforward\/} $f_*:MC_k(Y_1;R)\ra MC_k(Y_2;R)$ for $k\in\Z$ to be the unique $R$-linear map defined on generators $[V,n,s,t]$ of $MC_k(Y_1;R)$ by
\e
f_*[V,n,s,t]=[V,n,s,f\ci t].
\label{kh4eq7}
\e
Again, to show $f_*$ is well-defined, we must show that it maps relations (i),(ii) in $MC_k(Y_1;R)$ to relations (i),(ii) in $MC_k(Y_2;R)$. For (i) this is obvious. For (ii), suppose $\sum_{i\in I}a_i\,[V_i,n,s_i,t_i]=0$ in $MC_k(Y_1;R)$ by relation (ii) using open $0\in X\subseteq\R^n$. Then equation \eq{kh4eq2} for $\sum_{i\in I}a_i\,[V_i,n,s_i,f\ci t_i]$ in $MC_k(Y_2;R)$ for allowed $(x,y_2)$ in $X\t Y_2$ and $P_2$ an oriented $n+k$-plane in $T_xX\op T_{y_2}Y_2$ follows by summing \eq{kh4eq2} for $\sum_{i\in I}a_i\,[V_i,n,s_i,t_i]$ in $MC_k(Y_1;R)$ over all allowed $(x,y_1)\in X\t Y_1$ with $f(y_1)=y_2$ and all oriented $n+k$-planes $P_1$ in $T_xX\op T_{y_1}Y_1$ with $T_{(x,y_1)}(\id_X\t f)[P_1]=P_2$, where only finitely many terms in the sum can be nonzero. This proves $(*)$ for $\sum_{i\in I}a_i\,[V_i,n,s_i,f\ci t_i]$, with the same $0\in X\subseteq\R^n$, so (ii) gives $\sum_{i\in I}a_i\,[V_i,n,s_i,f\ci t_i]=0$ in $MC_k(Y_2;R)$. Therefore $f_*$ takes relation (ii) to relation (ii), and is well-defined.

Equations \eq{kh4eq3}, \eq{kh4eq7} give $f_*\ci\pd=\pd\ci f_*:MC_k(Y_1;R)\ra MC_{k-1}(Y_2;R)$. So the $f_*$ induce pushforwards $f_*:MH_k(Y_1;R)\ra MH_k(Y_2;R)$ on homology.

If $g:Y_2\ra Y_3$ is another smooth map of manifolds then $(g\ci f)_*=g_*\ci f_*$, on both M-chains $MC_*(Y_i;R)$ and M-homology $MH_*(Y_i;R)$. Also $(\id_Y)_*$ is the identity on both M-chains $MC_*(Y;R)$ and M-homology $MH_*(Y;R)$.
\label{kh4def2}
\end{dfn}

The next theorem is proved in~\S\ref{kh72}.

\begin{thm} Let\/ $Y$ be a manifold and\/ $R$ a commutative ring. Then:
\begin{itemize}
\setlength{\itemsep}{0pt}
\setlength{\parsep}{0pt}
\item[{\bf(a)}] Suppose $T\subseteq U\subseteq Y$ are open, and write $i:T\hookra U$ for the inclusion. Then $i_*:MC_k(T;R)\ra MC_k(U;R)$ is injective for all\/~$k\in\Z$.
\item[{\bf(b)}] Suppose $T,U\subseteq Y$ are open sets. Write $i:T\cap U\hookra T,$ $i':T\cap U\hookra U,$ $j:T\hookra T\cup U,$ $j':U\hookra T\cup U$ for the inclusions. Then for all\/ $k\in\Z$ the following sequence is exact:
\end{itemize}
\e
\xymatrix@C=11pt{ 0 \ar[r] & MC_k(T\!\cap\! U;R) \ar[rr]^{i_*\op -i'_*} && {\begin{subarray}{l} \ts \; MC_k(T;R)\\ \ts\op MC_k(U;R) \end{subarray}} \ar[rr]^(0.45){j_*\op j'_*} && MC_k(T\!\cup\! U;R) \ar[r] & 0. }
\label{kh4eq8}
\e
\begin{itemize}
\setlength{\itemsep}{0pt}
\setlength{\parsep}{0pt}
\item[{\bf(c)}] Suppose $U_1\subseteq U_2\subseteq \cdots\subseteq Y$ are open with\/ $U=\bigcup_{a=1}^\iy U_a$. Then we have an isomorphism with the direct limit for all\/ $k\in\Z$
\begin{equation*}
MC_k(U;R)\cong \underrightarrow{\lim}\,_{a=1}^\iy\, MC_k(U_a;R),
\end{equation*}
compatible with the pushforwards\/ $(i_a)_*:MC_k(U_a;R)\ra MC_k(U;R)$ and\/ $(i_{a,b})_*:MC_k(U_a;R)\ra MC_k(U_b;R),$ where $i_a:U_a\hookra U$ and\/ $i_{a,b}:U_a\hookra U_b$ are the inclusions for~$a\!\le\! b$.
\end{itemize}

Write $\uMC_k(Y;R)(U)=MC_k(U;R)$ for open $U\subseteq Y$ and\/ $k\in\Z,$ and define $\si_{TU}:\uMC_k(Y;R)(T)\ra\uMC_k(Y;R)(U)$ for open $T\subseteq U\subseteq Y$ by $\si_{TU}=i_*,$ for $i:T\hookra U$ the inclusion. Then in the notation of\/ {\rm\S\ref{kh25},} functoriality of\/ $i_*$ in Definition\/ {\rm\ref{kh4def2}} implies that\/ $\uMC_k(Y;R)$ is a precosheaf of\/ $R$-modules on $Y,$ and\/ {\bf(b)\rm,\bf(c)} that\/ $\uMC_k(Y;R)$ is a cosheaf, and\/ {\bf(a)} that\/ $\uMC_k(Y;R)$ is flabby.

Define $\pd(U):\uMC_k(Y;R)(U)\ra\uMC_{k-1}(Y;R)(U)$ for open $U\subseteq Y$ by $\pd(U)=\pd:MC_k(U;R)\ra MC_{k-1}(Y;R)$. Since $\pd\ci i_*=i_*\ci\pd:MC_k(T;R)\ra MC_{k-1}(U;R)$ for open $T\subseteq U\subseteq Y$ this defines a morphism of cosheaves $\pd:\uMC_k(Y;R)\ra \uMC_{k-1}(Y;R),$ with\/ $\pd\ci\pd=0:\uMC_k(Y;R)\ra \uMC_{k-2}(Y;R)$. So $\uMC_\bu(Y;R)=\bigl(\uMC_*(Y;R),\pd\bigr)$ is a complex of flabby cosheaves on~$Y$.
\label{kh4thm1}
\end{thm}

If $\ucE$ is a flabby cosheaf on $Y$, $U\subseteq Y$ is open, and $\al\in\ucE(U)$, Definition \ref{kh2def10} defined the {\it support\/} $\supp\al$, a compact subset of $U$. As $\uMC_k(Y;R)$ is a flabby cosheaf, this defines $\supp\al$ for all $\al\in MC_k(Y;R)=\uMC_k(Y;R)(Y)$. Let $[V,n,s,t]$ be a generator of $MC_k(Y;R)$. Then $s^{-1}(0)$ is compact in $V$ as $s:V\ra\R^n$ is proper near $0\in\R^n$, so $t[s^{-1}(0)]$ is a compact subset of $Y$. If $U$ is any open neighbourhood of $t[s^{-1}(0)]$ in $Y$ with inclusion $i:U\hookra Y$ then using Definition \ref{kh4def1}(ii) and \eq{kh4eq7} we see that
\begin{equation*}
[V,n,s,t]=\bigl[t^{-1}(U),n,s\vert_{t^{-1}(U)},t\vert_{t^{-1}(U)}\bigr]=i_*\bigl[t^{-1}(U),n,s\vert_{t^{-1}(U)},t\vert_{t^{-1}(U)}\bigr].
\end{equation*}
Hence Definition \ref{kh2def10} implies that
\e
\supp[V,n,s,t]\subseteq t[s^{-1}(0)]\subseteq Y.
\label{kh4eq9}
\e

\begin{dfn} Let $Y$ be a manifold, $Z\subseteq Y$ be open, and $i:Z\hookra Y$ the inclusion. Define the {\it relative M-chains\/} $MC_k(Y,Z;R)\!=\!MC_k(Y;R)/i_*\bigl(MC_k(Z;R)\bigr)$ for $k\in\Z$, and write $j_*:MC_k(Y;R)\ra MC_k(Y,Z;R)$ for the projection. Then $\pd:MC_k(Y;R)\ra MC_{k-1}(Y;R)$ descends to $\pd:MC_k(Y,Z;R)\ab\ra MC_{k-1}(Y,Z;R)$ as $i_*\ci\pd=\pd\ci i_*$, with $\pd\ci\pd=0$. Define the {\it relative M-homology groups\/} (or {\it relative integral M-homology groups\/}) $MH_*(Y,Z;R)$ to be the homology of this chain complex, as in Axiom~\ref{kh2ax1}(a).

By Theorem \ref{kh4thm1}(a) we have a short exact sequence of chain complexes
\begin{equation*}
\xymatrix@C=15pt{ 0 \ar[r] & \bigl(MC_*(Z;R),\pd\bigr) \ar[r]^{i_*} & \bigl(MC_*(Y;R),\pd\bigr) \ar[r]^(0.45){j_*} & \bigl(MC_*(Y,Z;R),\pd\bigr) \ar[r] & 0. }
\end{equation*}
In the usual way \cite[Lem.~24.1]{Munk} this induces a long exact sequence 
\begin{equation*}
\xymatrix@C=10.3pt{ \cdots \ar[r] & MH_k(Z;R) \ar[r]^(0.47){i_*} & MH_k(Y;R) \ar[r]^(0.45){j_*} & MH_k(Y,Z;R) \ar[r]^(0.48)\pd & MH_{k-1}(Z;R) \ar[r] & \cdots, }
\end{equation*}
which defines the {\it connecting morphisms\/} $\pd:MH_k(Y,Z;R)\ra MH_{k-1}(Z;R)$ as in Axiom \ref{kh2ax1}(b), and proves Axiom \ref{kh2ax1}(ii).

Let $f:Y_1\ra Y_2$ be a smooth map of manifolds, and $Z_1\subseteq Y_1$, $Z_2\subseteq Y_2$ be open with $f(Z_1)\subseteq Z_2$; for short we will say that $f:(Y_1,Z_2)\ra (Y_2,Z_2)$ is smooth. Then we have a commutative diagram of chain complexes
\begin{equation*}
\xymatrix@C=14pt@R=15pt{ 0 \ar[r] & \bigl(MC_*(Z_1;R),\pd\bigr) \ar[d]^{(f\vert_{Z_1})_*} \ar[r]_{i_*} & \bigl(MC_*(Y_1;R),\pd\bigr) \ar[d]^{f_*} \ar[r]_(0.45){j_*} & \bigl(MC_*(Y_1,Z_1;R),\pd\bigr) \ar@{.>}[d]^{f_*} \ar[r] & 0 \\
0 \ar[r] & \bigl(MC_*(Z_2;R),\pd\bigr) \ar[r]^{i_*} & \bigl(MC_*(Y_2;R),\pd\bigr) \ar[r]^(0.45){j_*} & \bigl(MC_*(Y_2,Z_2;R),\pd\bigr) \ar[r] & 0, }
\end{equation*}
inducing $f_*:MC_k(Y_1,Z_1;R)\ra MC_k(Y_2,Z_2;R)$ and $f_*:MH_k(Y_1,Z_1;R)\ra MH_k(Y_2,Z_2;R)$, as in Axiom \ref{kh2ax1}(c), and a commutative diagram
\begin{small}\begin{equation*}
\xymatrix@C=12pt@R=15pt{ \cdots \ar[r] & MH_k(Z_1;R) \ar[d]_{(f\vert_{Z_1})_*} \ar[r]_(0.47){\raisebox{-6pt}{$\st(i_1)_*$}} & MH_k(Y_1;R) \ar[d]_{f_*} \ar[r]_(0.45){\raisebox{-6pt}{$\st(j_1)_*$}} & MH_k(Y_1,Z_1;R) \ar[d]^{f_*} \ar[r]_(0.48){\raisebox{-6pt}{$\st\pd$}} & MH_{k-1}(Z_1;R) \ar[d]^{(f\vert_{Z_1})_*} \ar[r] & \cdots \\ 
\cdots \ar[r] & MH_k(Z_2;R) \ar[r]^(0.47){\raisebox{3pt}{$\st(i_1)_*$}} & MH_k(Y_2;R) \ar[r]^(0.45){\raisebox{3pt}{$\st(j_2)_*$}} & MH_k(Y_2,Z_2;R) \ar[r]^(0.48){\raisebox{3pt}{$\st\pd$}} & MH_{k-1}(Z_2;R) \ar[r] & \cdots}
\end{equation*}\end{small}\noindent by \cite[Th.~24.2]{Munk}, proving Axiom~\ref{kh2ax1}(iii).

If $g:(Y_2,Z_2)\ra (Y_3,Z_3)$ is another smooth map then $(g\ci f)_*=g_*\ci f_*$ on $MC_*(Y_i;R)$ implies that $(g\ci f)_*=g_*\ci f_*$ on $MC_*(Y_i,Z_i;R)$, and hence on $MH_*(Y_i,Z_i;R)$. Also $(\id_Y)_*=\id$ on $MC_*(Y;R)$ implies that $(\id_Y)_*$ is the identity on $MC_*(Y,Z;R)$, and hence on $MH_*(Y,Z;R)$. This proves Axiom~\ref{kh2ax1}(i).

\label{kh4def3}
\end{dfn}

Our goal is to show that M-homology $MH_*(-;R)$ satisfies Axiom \ref{kh2ax1}, so that it is canonically isomorphic to ordinary homology of manifolds by Theorem \ref{kh2thm1}. So far we have defined all the data in Axiom \ref{kh2ax1}(a)--(c), and verified Axiom \ref{kh2ax1}(i)--(iii). The next proposition, proved in \S\ref{kh73}, gives Axiom \ref{kh2ax1}(iv), the homotopy axiom. The main idea is to define $R$-module morphisms $G:MC_k(Y_1,Z_1;R)\ra MC_{k+1}(Y_2,Z_2;R)$ for $k\in\Z$ by
\begin{equation*}
G:[V,n,s,t]\longmapsto (-1)^{\dim V}\bigl[V\t[0,1],n,s\ci\pi_V,g\ci(t\t\id_{[0,1]})\bigr],
\end{equation*}
and show that $G$ is a chain homotopy between $f_*,f_*':\bigl(MC_*(Y_1,Z_1;R),\pd\bigr)\ra \bigl(MC_*(Y_1,Z_1;R),\pd\bigr)$.

\begin{prop} Suppose $Y_1,Y_2$ are manifolds, $Z_1\subseteq Y_1,$ $Z_2\subseteq Y_2$ are open, and\/ $g:Y_1\t[0,1]\ra Y_2$ is a smooth map in $\Manc$ with\/ $g(Z_1\t[0,1])\subseteq Z_2$. Define $f,f':Y_1\ra Y_2$ by $f(y)=g(y,0)$ and $f'(y)=g(y,1)$ for $y\in Y_1$. Then $f_*=f'_*:MH_k(Y_1,Z_1;R)\ra MH_k(Y_2,Z_2;R)$ for all\/ $k\in\Z$.
\label{kh4prop2}
\end{prop}

Putting $T=Z$ and $U=Y\sm S$ in Theorem \ref{kh4thm1}(b) and using $MC_k(Y,Z;R)=MC_k(Y;R)/i_*(MC_k(Z;R))$, we deduce Axiom \ref{kh2ax1}(v), the excision axiom:

\begin{prop} Suppose $Y$ is a manifold, $Z\subseteq Y$ is open, and\/ 
$S\subseteq Z$ is closed in $Y$. Then $j_*:MC_k(Y\sm S,Z\sm S;R)\ra MC_k(Y,Z;R)$ is an isomorphism for all\/ $k\in\Z,$ with\/ $j:Y\sm S\hookra Y$ the inclusion. Hence $j_*:MH_k(Y\sm S,Z\sm S;R)\ab\ra MH_k(Y,Z;R)$ is also an isomorphism for all\/~$k\in\Z$. 
\label{kh4prop3}
\end{prop}

We prove Axiom \ref{kh2ax1}(vi), the additivity axiom:

\begin{lem} Suppose $Y$ is a manifold with\/ $Y=\coprod_{a\in A}Y_a$ for $A$ a countable indexing set and each\/ $Y_a$ open and closed in $Y$. Let\/ $Z\subseteq Y$ be open, and set\/ $Z_a=Z\cap Y_a$. Then for all\/ $k\ge 0$ we have canonical isomorphisms
\e
MC_k(Y;R)\cong \ts\bigop\limits_{a\in A}MC_k(Y_a;R),\;\> MC_k(Y,Z;R)\cong \ts\bigop\limits_{a\in A}MC_k(Y_a,Z_a;R),
\label{kh4eq10}
\e
compatible with the morphisms $(i_a)_*:MC_k(Y_a;R)\ra MC_k(Y;R)$ for $a\in A$ induced by the inclusions $i_a:Y_a\hookra Y$. Hence we also have
\e
MH_k(Y,Z;R)\cong \ts\bigop\limits_{a\in A}MH_k(Y_a,Z_a;R),
\label{kh4eq11}
\e
compatible with\/ $(i_a)_*:MH_k(Y_a,Z_a;R)\ra MH_k(Y,Z;R)$ for $a\in A$.
\label{kh4lem2}
\end{lem}

\begin{proof} Let $[V,n,s,t]$ be a generator of $MC_k(Y;R)$. Then $s:V\ra\R^n$ is proper over an open neighbourhood $X$ of 0 in $\R^n$. Choose $0\in X'\subset K\subset X\subset\R^n$, with $X'$ open and $K$ compact, and set $V'=s^{-1}(X)$, and $s'=s\vert_{V'}$, $t'=t\vert_{V'}$. Then Definition \ref{kh4def1}(ii) implies that $[V,n,s,t]=[V',n,s',t']$. Also $V'\subseteq s^{-1}(K)\subseteq V$, where $s^{-1}(K)$ is compact as $K$ is and $s$ is proper over~$X\supset K$. 

Now $t\vert_{s^{-1}(K)}:s^{-1}(K)\ra Y=\coprod_{a\in A}Y_a$ can map to only finitely many $Y_a$, as $s^{-1}(K)$ is compact and $t$ is continuous. Hence $s':V'\ra Y=\coprod_{a\in A}Y_a$ maps to only finitely many $Y_a$, as $V'\subseteq s^{-1}(K)$. Thus we may decompose $V'=\coprod_{b\in B}V'_b$, where $B\subseteq A$ is finite and $s'(V'_b)\subseteq Y_b\subseteq Y$. By \eq{kh4eq6} we now have
\begin{equation*}
[V,n,s,t]=[V',n,s',t']=\ts\sum_{b\in B}\bigl[V'_b,n,s'\vert_{V'_b},t'\vert_{V'_b}\bigl],
\end{equation*}
where $\bigl[V'_b,n,s'\vert_{V'_b},t'\vert_{V'_b}\bigl]\in (i_b)_*\bigl(MC_k(Y_b;R)\bigr)$.

Thus the natural map $\bigop_{a\in A}MC_k(Y_a;R)\ra MC_k(Y;R)$ is surjective, since this holds for every generator. It is also injective, as the relations in $MC_k(Y;R)$ can all be written as finite sums of relations in $MC_k(Y_a;R)$ by a similar argument. This proves the first isomorphism of \eq{kh4eq10}. The second follows by taking quotients, and then taking cohomology yields~\eq{kh4eq11}.
\end{proof}

Axiom \ref{kh2ax1}(vii), the dimension axiom, is the next theorem, proved in \S\ref{kh74}. It is crucial to our programme, and the proof is quite complex.

\begin{thm} We have $MH_k(*;R)=0$ for all\/ $0\ne k\in\Z$. There is a canonical isomorphism $MH_0(*;R)\cong R$ identifying the fundamental class $[[*]]\in MH_0(*;R)$ with\/~$1\in R$.
\label{kh4thm2}
\end{thm}

We have now proved all of Axiom \ref{kh2ax1} for M-homology. So Theorem \ref{kh2thm1} gives:

\begin{thm} M-homology is a homology theory of manifolds. There are canonical isomorphisms $MH_k(Y;R)\cong H_k(Y;R),$ $MH_k(Y,Z;R)\cong H_k(Y,Z;R)$ for all\/ $Y,Z,k,$ preserving the data $f_*,\pd$ and isomorphisms $MH_0(*;R)\cong R\cong H_0(*;R),$ where $H_*(-;R)$ is any other homology theory of manifolds over $R,$ such as singular homology\/~$H_*^\rsi(-;R)$.
\label{kh4thm3}
\end{thm}

\begin{rem}{\bf (Change of base ring $R$.)} Suppose $R,R'$ are commutative rings, and $\rho:R\ra R'$ a ring morphism. Then there are natural $\rho$-linear morphisms $\rho_*:MC_k(Y;R)\ra MC_k(Y;R')$ mapping
\begin{equation*}
\rho_*:\ts\sum_{i\in I}a_i\,[V_i,n_i,s_i,t_i]\longmapsto \sum_{i\in I}\rho(a_i)\,[V_i,n_i,s_i,t_i].
\end{equation*}
These commute with differentials $\pd$ and pushforwards $f_*$, and so induce morphisms $\rho_*:MH_k(Y;R)\ra MH_k(Y;R')$. By the natural extension of Theorem \ref{kh2thm1} to change of base ring as in \cite{EiSt}, under the isomorphisms $MH_k(Y;R)\cong H_k(Y;R)$,  $MH_k(Y;R')\cong H_k(Y;R')$ from Theorem \ref{kh4thm3}, these $\rho_*:MH_k(Y;R)\ab\ra MH_k(Y;R')$ correspond to the usual morphisms $\rho_*:H_k(Y;R)\ra H_k(Y;R')$ from change of base ring in homology.
 
In the same way, though we will not mention it again, throughout \S\ref{kh4} there are canonical morphisms $\rho_*$ for change of base ring $\rho:R\ra R'$, which commute with all available structures, and are identified with the usual change of base ring morphisms under the isomorphisms $MH^k(Y;R)\cong H^k(Y;R),\ldots.$
\label{kh4rem3}
\end{rem}

For smooth singular homology $H_*^\ssi(-;R)$ in Example \ref{kh2ex2}, we can realize the isomorphisms $H_k^\ssi(Y;R)\cong MH_k(Y;R)$ from a morphism of chain complexes. 

\begin{ex} Let $Y$ be a smooth manifold, and let $\bigl(C_*^\ssi(Y;R),\pd\bigr)$ be as in Example \ref{kh2ex2}, so that $C_k^\ssi(Y;R)$ is the free $R$-module spanned by smooth maps $\si:\De_k\ra Y$ in $\Manc$. By Assumption \ref{kh3ass4}(c)(C), such $\si$ are also morphisms in $\tManc$. Define $F_\ssi^\Mh:C_k^\ssi(Y;R)\ra MC_k(Y;R)$ to be the unique $R$-linear morphism acting on generators $\si$ by
\e
F_\ssi^\Mh:\si\longmapsto[\De_k,0,0,\si],
\label{kh4eq12}
\e
so that $V=\De_k$, $n=0$, $s=0:V\ra\R^0$ and $t=\si:V\ra Y$. This is well-defined as there are no relations in $C_k^\ssi(Y;R)$. We have
\ea
\pd\ci F_\ssi^\Mh(\si)&=\bigl[\pd\De_k,0,0,\si\ci i_{\De_k}\bigr]=
\ts\sum_{j=0}^k(-1)^j\bigl[\De_{k-1},0,0,\si\ci i_{\De_k}\ci F_j^k\bigr]\nonumber\\
&=F_\ssi^\Mh\bigl(\ts\sum_{j=0}^k(-1)^j\si\ci i_{\De_k}\ci F_j^k\bigr)=F_\ssi^\Mh\ci\pd\si,
\label{kh4eq13}
\ea
using \eq{kh4eq3} and \eq{kh4eq12} in the first step, the decomposition
$\pd\De_k=\coprod_{j=0}^k\pd_j\De_k$ and diffeomorphism $F_j^k:\De_{k-1}\ra\pd_j\De_k$ for $j=0,\ldots,k$ multiplying orientations by $(-1)^j$ from Example \ref{kh2ex1} and \eq{kh4eq5}--\eq{kh4eq6} in the second, \eq{kh4eq12} in the third and \eq{kh2eq6} in the fourth. As this holds for all generators $\si$ of $C_k^\ssi(Y;R)$, we have $\pd\ci F_\ssi^\Mh=F_\ssi^\Mh\ci\pd:C_k^\ssi(Y;R)\ra MC_{k-1}(Y;R)$. Thus $F_\ssi^\Mh$ induces morphisms~$F_\ssi^\Mh:H_k^\ssi(Y;R)\ra MH_k(Y;R)$. 

Comparing \S\ref{kh31} and above we see that $F_\ssi^\Mh$ commutes with pushforwards $f_*:C_k^\ssi(Y_1;R)\ra C_k^\ssi(Y_2;R)$ and $f_*:MC_k(Y_1;R)\ra MC_k(Y_2;R)$. For open $Z\subseteq Y$, applying this to the inclusion $i:Z\hookra Y$, we see that the $F_\ssi^\Mh:C_k^\ssi(Y;R)\ra MC_k(Y;R)$ descend to $F_\ssi^\Mh:C_k^\ssi(Y,Z;R)\ra MC_k(Y,Z;R)$ with $\pd\ci F_\ssi^\Mh=F_\ssi^\Mh\ci\pd$, and hence induce~$F_\ssi^\Mh:H_k^\ssi(Y,Z;R)\ra MH_k(Y,Z;R)$. 

These $F_\ssi^\Mh$ commute with $f_*:C_k^\ssi(Y_1,Z_1;R)\ra C_k^\ssi(Y_2,Z_2;R)$ and $f_*:MC_k(Y_1,Z_1;R)\ra MC_k(Y_2,Z_2;R)$, and so passing to homology, the morphisms $F_\ssi^\Mh:H_k^\ssi(Y,Z;R)\ra MH_k(Y,Z;R)$ commute with pushforwards $f_*:H_k^\ssi(Y_1,Z_1;R)\ra H_k^\ssi(Y_2,Z_2;R)$, $f_*:MH_k(Y_1,Z_1;R)\ra MH_k(Y_2,Z_2;R)$, and with $\pd:H_k^\ssi(Y,Z;R)\ra H_{k-1}^\ssi(Z;R)$, $\pd:MH_k(Y,Z;R)\ra MH_{k-1}(Z;R)$. Since the chosen isomorphisms $H_0^\ssi(*;R)\cong R$, $MH_0(*;R)\cong R$ both identify $[[*]]$ with $1\in R$, we see that $F_\ssi^\Mh:H_0^\ssi(*;R)\ra MH_0(*;R)$ is compatible with the isomorphisms~$H_0^\ssi(*;R)\cong R\cong MH_0(*;R)$. 

Thus by the last part of Theorem \ref{kh2thm1}, the morphisms $F_\ssi^\Mh:H_k^\ssi(Y;R)\ra MH_k(Y;R)$ and $F_\ssi^\Mh:H_k^\ssi(Y,Z;R)\ra MH_k(Y,Z;R)$, defined above via morphisms of chain complexes, are the canonical isomorphisms from Theorem~\ref{kh4thm3}.

\label{kh4ex1}
\end{ex}

\begin{rem} Suppose that $Y$ is a compact, oriented manifold of dimension $m$. Then as in Remark \ref{kh2rem1}(f) we have a fundamental class $[[Y]]\in H_m^\ssi(Y;R)$. As in Example \ref{kh2ex2} we may define this explicitly in smooth singular homology by choosing a triangulation of $Y$ into finitely many smooth $m$-simplices $\si_i:\De_m\ra Y$ for $i\in I$, and setting $\ep_i=1$ if $\si_i$ is orientation-preserving, and $\ep_i=-1$ otherwise. Then $\sum_{i\in I}\ep_i\,\si_i\in C_m^\ssi(Y;R)$ and~$[[Y]]=\bigl[\sum_{i\in I}\ep_i\,\si_i\bigr]\in H_m^\ssi(Y;R)$.

Definition \ref{kh4def1} defined $[Y]=[Y,0,0,\id_Y]\in MC_m(Y;R)$. We have
\e
\begin{split}
[Y]&-F_\ssi^\Mh\bigl(\ts\sum_{i\in I}\ep_i\,\si_i\bigr)=
[Y,0,0,\id_Y]-\ts\sum_{i\in I}\ep_i\,[\De_m,0,0,\si_i]\\
&=[Y,0,0,\id_Y]-\ts\sum_{i\in I}[\si_i(\De_m),0,0,\id_{\si_i(\De_m)}]=0,
\end{split}
\label{kh4eq14}
\e
using \eq{kh4eq12} in the first step, that $\si_i:\De_m\ra\si_i(\De_m)$ is a diffeomorphism multiplying orientations by $\ep_i$ in the second, and Definition \ref{kh4def1}(ii) in the third. Hence $[Y]=F_\ssi^\Mh\bigl(\sum_{i\in I}\ep_i\,\si_i\bigr)$ in $MC_m(Y;R)$, so $[[Y]]=F_\ssi^\Mh\bigl([[Y]]\bigr)$ in~$MH_m(Y;R)$.

This proves that the fundamental class $[[Y]]\in MH_m(Y;R)$ is identified with the usual fundamental class $[[Y]]\in H_m(Y;R)$ by the canonical isomorphism $MH_m(Y;R)\cong H_m(Y;R)$ of Theorem \ref{kh4thm3}. Note too that although the chain $\sum_{i\in I}\ep_i\,\si_i\in C_m^\ssi(Y;R)$ representing $[[Y]]\in H_k^\ssi(Y;R)$ involves an arbitrary choice of triangulation, $[Y]$ is the unique cycle in $MC_m(Y;R)$ representing $[[Y]]\in MH_m(Y;R)$, as $MC_{m+1}(Y;R)=0$ by Lemma~\ref{kh4lem1}.
\label{kh4rem4}
\end{rem}

We can also do the same for cosheaf smooth singular homology $\hat H_k^\ssi(Y;R)$ in Example~\ref{kh2ex16}.

\begin{ex} We continue in the situation of Example \ref{kh4ex1}. Now Example \ref{kh2ex3} defined barycentric subdivision morphisms $B:C_k^\ssi(Y;R)\ra C_k^\ssi(Y;R)$ for $k=0,1,\ldots$ by the formula \eq{kh2eq9}. We claim that
\begin{equation*}
F_\ssi^\Mh=F_\ssi^\Mh\ci B:C_k^\ssi(Y;R)\longra MC_k(Y;R).
\end{equation*}
To see this, note that if $\si$ is a generator of $C_k^\ssi(Y;R)$, so that $\si:\De_k\ra Y$ is a smooth map in $\Manc$, then $B(\si)=\sum_{j=1}^{(k+1)!}\ep_k^j\,(\si\ci B_k^j)$ is the result of triangulating the $k$-simplex $\De_k$ into $(k+1)!$ smaller simplices $B_k^j(\De_k)\subset \De_k$ for $j=1,\ldots,(k+1)!$, and restricting $\si$ to each $B_k^j(\De_k)$. We have
\begin{align*}
F_\ssi^\Mh(\si)&-F_\ssi^\Mh\ci B(\si)=
\bigl[\De_k,0,0,\si\bigr]-\ts\sum_{j=1}^{(k+1)!}\ep_k^j\,\bigl[\De_k,0,0,\si\ci B_k^j\bigr]\\
&=\bigl[\De_k,0,0,\si\bigr]-\ts\sum_{j=1}^{(k+1)!}\bigl[B_k^j(\De_k),0,0,\si\vert_{B_k^j(\De_k)}\bigr]=0,
\end{align*}
where the second step uses that $B_k^j:\De_k\ra B_k^j(\De_k)$ is a diffeomorphism which multiplies orientations by $\ep_k^j$, and the third uses relation Definition~\ref{kh4def1}(ii).

Example \ref{kh2ex16} defined the cosheaf smooth singular chains $\hat C_k^\ssi(Y;R)$, as the direct limit of $C_k^\ssi(Y;R)\,{\buildrel B\over\longra}\,C_k^\ssi(Y;R)\,{\buildrel B\over\longra}\,\cdots.$ Equation \eq{kh4eq16} and properties of direct limits imply that there is a unique morphism $\hat F_\ssi^\Mh:\hat C_k^\ssi(Y;R)\ra MC_k(Y;R)$ such that $\hat F_\ssi^\Mh\ci\Pi_j=F_\ssi^\Mh:C_k^\ssi(Y;R)\ra MC_k(Y;R)$ for all $j=0,1,\ldots,$ where $\Pi_j$ are the morphisms in \eq{kh2eq36} from the direct limit.

Since $\pd\ci F_\ssi^\Mh=F_\ssi^\Mh\ci\pd:C_k^\ssi(Y;R)\ra MC_{k-1}(Y;R)$ and $\pd\ci B=B\ci\pd:C_k^\ssi(Y;R)\ra C_{k-1}^\ssi(Y;R)$, properties of direct limits imply that $\pd\ci\hat F_\ssi^\Mh=\hat F_\ssi^\Mh\ci\pd:\hat C_k^\ssi(Y;R)\ra MC_{k-1}(Y;R)$. Therefore we have induced morphisms $\hat F_\ssi^\Mh:\hat H_k^\ssi(Y;R)\ra MH_k(Y;R)$ on homology. As $\hat F_\ssi^\Mh\ci\Pi_j=F_\ssi^\Mh$, where both $\Pi_j$ and $F_\ssi^\Mh$ induce isomorphisms on homology by Examples \ref{kh2ex3} and \ref{kh4ex1}, these $\hat F_\ssi^\Mh:\hat H_k^\ssi(Y;R)\ra MH_k(Y;R)$ are isomorphisms.

If $f:Y_1\!\ra\! Y_2$ is a smooth map of manifolds then as both $F_\ssi^\Mh:C_k^\ssi(Y_a;R)\!\ra\!MC_k(Y_a;R)$ and $\Pi_j:C_k^\ssi(Y_a;R)\ra\hat C_k^\ssi(Y_a;R)$ commute with pushforwards $f_*$, we have $f_*\ci\hat F_\ssi^\Mh=\hat F_\ssi^\Mh\ci f_*:\hat C_k^\ssi(Y_1;R)\ra MC_k(Y_2;R)$. Thus the isomorphisms $\hat F_\ssi^\Mh:\hat H_k^\ssi(Y_a;R)\ra MH_k(Y_a;R)$ preserve pushforwards~$f_*$.

Example \ref{kh2ex16} defined a flabby cosheaf $\hat\ucC{}_k^\ssi(Y;R)$ of $R$-modules on $Y$ with $\hat\ucC{}_k^\ssi(Y;R)(U)=\hat C_k^\ssi(U;R)$ for open $U\subseteq Y$. Theorem \ref{kh4thm1} defined a flabby cosheaf $\uMC_k(Y;R)$ of $R$-modules on $Y$ with $\uMC_k(Y;R)(U)=MC_k(U;R)$ for open $U\subseteq Y$. Define  
$\hat {\ul{F\!}\,}{}_\ssi^\Mh:\hat\ucC{}_k^\ssi(Y;R)\ra \uMC_k(Y;R)$ by $\hat {\ul{F\!}\,}{}_\ssi^\Mh(U)=\hat F_\ssi^\Mh:\hat C_k^\ssi(U;R)\ra MC_k(U;R)$ for all open $U\subseteq Y$. This is a morphism of cosheaves because $\hat F_\ssi^\Mh\ci i_*=i_*\ci \hat F_\ssi^\Mh:\hat C_k^\ssi(V;R)\ra MC_k(U;R)$ for open $V\subseteq U\subseteq Y$, where $i:V\hookra U$ is the inclusion. We now have a diagram of cosheaves on $Y$
\e
\begin{gathered}
\xymatrix@C=18pt@R=15pt{ \cdots \ar[r]_(0.25)\pd  & \hat\ucC{}_{k+1}^\ssi(Y;R) \ar[r]_\pd \ar[d]^{\hat {\ul{F\!}\,}{}_\ssi^\Mh} & \hat\ucC{}_k^\ssi(Y;R) \ar[r]_\pd \ar[d]^{\hat {\ul{F\!}\,}{}_\ssi^\Mh} & \hat\ucC{}_{k-1}^\ssi(Y;R) \ar[r]_(0.65)\pd \ar[d]^{\hat {\ul{F\!}\,}{}_\ssi^\Mh} &  \cdots \\
\cdots \ar[r]^(0.25)\pd  & \uMC_{k+1}(Y;R) \ar[r]^\pd & \uMC_k(Y;R) \ar[r]^\pd & \uMC_{k-1}(Y;R) \ar[r]^(0.65)\pd &  \cdots,\!\!{} }
\end{gathered}
\label{kh4eq15}
\e
where the rows are complexes, which commutes because $\pd\ci\hat F_\ssi^\Mh=\hat F_\ssi^\Mh\ci\pd:\hat C_k^\ssi(U;R)\ra MC_{k-1}(U;R)$ for $U\subseteq Y$ open. Note that on the top row $\hat\ucC{}_k^\ssi(Y;R)=0$ for $k<0$, but on the bottom row $\uMC_k(Y;R)=0$ if~$k>\dim Y$.

\label{kh4ex2}
\end{ex}

\subsection{\texorpdfstring{M-cohomology $MH^*(Y;R)$}{M-cohomology}}
\label{kh42}

We now discuss {\it M-cohomology\/} $MH^*(Y;R)$, the dual cohomology theory to M-homology in \S\ref{kh41}. As in \S\ref{kh2}, homology is compactly-supported, but cohomology is not (though compactly-supported cohomology is). Because of this, in M-cochains $MC^*(Y;R)$ we need to allow {\it infinite\/} sums $\sum_{i\in I}a_i\,[V_i,n_i,s_i,t_i]$ satisfying a local finiteness condition over $Y$.

As specifying and working with the relations on such infinite sums would be complicated, we proceed in two stages. First we define spaces $\cP MC^k(Y;R)$ by generators and relations, using only finite sums. We show $U\mapsto\cP MC^k(U;R)$ for open $U\subseteq Y$ forms a strong presheaf $\cP\MC^k(Y;R)$ on $Y$. Then we define $\MC^k(Y;R)$ to be the sheafification of $\cP\MC^k(Y;R)$, and $MC^k(Y;R)=\MC^k(Y;R)(Y)$ to be its global sections. By sheaf theory, the required locally finite sums make sense in~$MC^k(Y;R)$.

\begin{dfn} Let $Y$ be a manifold, of dimension $m$. Consider quadruples $(V,n,s,t)$, where $V$ is a manifold with corners (object in $\tManc$), and $n\in\N,$ and $s:V\ra\R^n$ is a smooth map (morphism in $\tManc$), and $t:V\ra Y$ is a cooriented submersion (i.e.\ a pair $(t,c_t)$ of a submersion $t:V\ra Y$ in $\tManc$ and a coorientation $c_t$ for $t$, usually left implicit), such that $(s,t):V\ra\R^n\t Y$ is proper over an open neighbourhood of $\{0\}\t Y$ in~$\R^n\t Y$.

Define an equivalence relation $\sim$ on such quadruples by $(V,n,s,t)\sim(V',\ab n',\ab s',\ab t')$ if $n=n'$, and there exists a diffeomorphism $f:V\ra V'$ with $s=s'\ci f$ and $t=t'\ci f$ such that the coorientations satisfy $c_t=c_{t'}\ci c_f$, where $c_f$ is the natural coorientation on $f$ from Assumption \ref{kh3ass6}(e). Write $[V,n,s,t]$ for the $\sim$-equivalence class of $(V,n,s,t)$. We call $[V,n,s,t]$ a {\it generator}.

For each $k\in\Z$, define the {\it M-precochains\/} $\cP MC^k(Y;R)$ to be the $R$-module generated by such $[V,n,s,t]$ with with $\dim V+k=m+n$, subject to the relations:
\begin{itemize}
\setlength{\itemsep}{0pt}
\setlength{\parsep}{0pt}
\item[(i)] For each generator $[V,n,s,t]$ and each $i=0,\ldots,n$ we have
\e
[V,n,s,t]=(-1)^{n-i}[V\t\R,n+1,s',t\ci\pi_V]\quad\text{in $\cP MC^k(Y;R)$,}
\label{kh4eq16}
\e
where writing $s=(s_1,\ldots,s_n):V\ra\R^n$ with $s_j:V\ra\R$ for $j=1,\ldots,n$ and $\pi_V:V\t\R\ra V$, $\pi_\R:V\t\R\ra\R$ for the projections, then 
\begin{equation*}
s'=(s_1\ci\pi_V,\ldots,s_i\ci\pi_V,\pi_\R,s_{i+1}\ci\pi_V,\ldots,s_n\ci\pi_V):V\t\R\ra\R^{n+1},
\end{equation*}
and $t\ci\pi_V$ has coorientation $c_{t\ci\pi_V}=c_t\ci c_{\pi_V}$, where $c_t$ is the given coorientation on $t:V\ra Y$, and $c_{\pi_V}$ is the coorientation on $\pi_V:V\t\R\ra V$ induced by the standard orientation on $\R$, as in Assumption~\ref{kh3ass6}(d),(f),(k).
\item[(ii)] Let $I$ be a finite indexing set, $a_i\in R$ for $i\in I$, and $[V_i,n,s_i,t_i]$, $i\in I$ be generators for $\cP MC^k(Y;R)$, all with the same $n$. Suppose there exists an open neighbourhood $X$ of $\{0\}\t Y$ in $\R^n\t Y$, such that $(s_i,t_i):V_i\ra\R^n\t Y$ is proper over $X$ for all $i\in I$, and the following condition holds:
\begin{itemize}
\setlength{\itemsep}{0pt}
\setlength{\parsep}{0pt}
\item[$(*)$] Suppose $(x,y)\!\in\! X$, such that for all $i\!\in\! I$ and $v\!\in\! V_i$ with $(s_i,t_i)(v)\!=\!(x,y)$, we have that $v\in V_i^\ci$ and 
\begin{equation*}
T_v(s_i,t_i):T_vV_i^\ci\longra T_x\R^n\op T_yY
\end{equation*}
is injective. This implies that $(s_i,t_i)\vert_{V_i^\ci}$ is an embedding near $v\in V_i^\ci$. Hence $(s_i,t_i):V_i\ra\R^n\t Y$ is injective near each $v$ in $(s_i,t_i)^{-1}(x,y)$, so $(s_i,t_i)^{-1}(x,y)$ has the discrete topology, and thus is finite as $(s_i,t_i)$ is proper over $X$. Note too that $T_vV_i^\ci$ is a vector space of dimension $m+n-k$ and $\d t\vert_v:T_vV_i^\ci\ra T_yY$ is cooriented, since $t\vert_{V_i^\ci}:V_i^\ci\ra Y$ is a cooriented smooth map of  manifolds by Assumption \ref{kh3ass6}(j). We require that for all $(m+n-k)$-planes $P\subseteq T_x\R^n\op T_yY$ with $\pi_{T_yY}:P\ra T_yY$ cooriented, we have
\e
\begin{split}
&\sum_{\begin{subarray}{l} i\in I,\; v\in V_i^\ci:(s_i,t_i)(v)=(x,y),\;  T_v(s_i,t_i)[T_vV_i^\ci]=P \\ \text{$T_v(s_i,t_i):T_vV_i^\ci\,{\buildrel\cong\over\longra}\,P$ is coorientation-preserving}\end{subarray}\!\!\!\!\!\!\!} a_i=\\
&\sum_{\begin{subarray}{l} i\in I,\; v\in V_i^\ci:(s_i,t_i)(v)=(x,y),\; T_v(s_i,t_i)[T_vV_i^\ci]=P \\ \text{$T_v(s_i,t_i):T_vV_i^\ci\,{\buildrel\cong\over\longra}\,P$ is coorientation-reversing}\end{subarray}\!\!\!\!\!\!\!} a_i\qquad\text{in $R$.}
\end{split}
\label{kh4eq17}
\e
\end{itemize}
Then
\begin{equation*}
\sum_{i\in I}a_i\,[V_i,n,s_i,t_i]=0\qquad\text{in $\cP MC^k(Y;R)$.}
\end{equation*}
\end{itemize}
\label{kh4def4}
\end{dfn}

As in Remark \ref{kh4rem2}, we see that if $[V,n,s,t]$ is a generator in $\cP MC^k(Y;R)$ and $c_t$ is the given coorientation for $t$, we have
\e
\bigl[V,n,s,(t,-c_t)\bigr]=-\bigl[V,n,s,(t,c_t)\bigr],
\label{kh4eq18}
\e
and if $[V_1,n,s_1,t_1]$, $[V_2,n,s_2,t_2]$ are generators in $\cP MC^k(Y;R)$ then
\e
[V_1\amalg V_2,n,s_1\amalg s_2,t_1\amalg t_2]=[V_1,n,s_1,t_1]+[V_2,n,s_2,t_2].
\label{kh4eq19}
\e

We define differentials, identities, and pullbacks on the~$\cP MC^k(Y;R)$.

\begin{dfn} For each manifold $Y$ and $k\in\Z$, let $\d:\cP MC^k(Y;R)\ra \cP MC^{k+1}(Y;R)$ be the unique $R$-linear map with
\e
\d[V,n,s,t]=[\pd V,n,s\ci i_V,t\ci i_V],
\label{kh4eq20}
\e
for all generators $[V,n,s,t]$, where the coorientation on $t\ci i_V$ is $c_{t\ci i_V}=c_t\ci c_{i_V}$, for $c_t$ the given coorientation on $t$, and $c_{i_V}$ the coorientation for $i_V:\pd V\ra V$ from Assumption \ref{kh3ass6}(h). We show that $\d$ is well-defined by an almost identical proof to that of Proposition \ref{kh4prop1} in \S\ref{kh71}, the differences being in the properness conditions, and in using coorientations rather than orientations. We show that $\d\ci\d=0$ as for $\pd$ in Definition~\ref{kh4def1}.

Define the {\it identity precocycle\/} $\Id_Y=[Y,0,0,\id_Y]$ in $\cP MC^0(Y;R)$. Here $t=\id_Y:Y\ra Y$ has the coorientation from Assumption \ref{kh3ass6}(e), and $(s,t):Y\ra\R^0\t Y$ is proper. We have $\d\,\Id_Y=0$ as $\pd Y=\es$.

Suppose $f:Y_1\ra Y_2$ is a smooth map of manifolds. For each $k\in\Z$, define the {\it pullback\/} $f^*:\cP MC^k(Y_2;R)\ra\cP MC^k(Y_1;R)$ to be the unique $R$-module morphism acting on generators $[V,n,s,t]$ of $\cP MC^k(Y_2;R)$ by
\e
f^*[V,n,s,t]=[V',n,s',t']:=\bigl[V\t_{t,Y_2,f}Y_1,n,s\ci\pi_V,\pi_{Y_1}\bigr],
\label{kh4eq21}
\e
where $V'=V\t_{t,Y_2,f}Y_1$ is the fibre product in $\tManc$, which exists by Assumption \ref{kh3ass5}(c) as $Y_2$ is a manifold and $t$ a submersion, with projections $\pi_V:V'\ra V$ and $\pi_{Y_1}:V'\ra Y_1$. Here $t'=\pi_{Y_1}:V'\ra Y_1$ is a submersion by Assumption \ref{kh3ass5}(c), and has a coorientation $c_{t'}$ determined by the given coorientation $c_t$ on $t:V\ra Y_2$ by Assumption \ref{kh3ass6}(l). Proposition \ref{kh4prop4} shows $f^*$ is well-defined.

From \eq{kh4eq20} and \eq{kh4eq21} we see that $f^*\ci\d=\d\ci f^*:\cP MC^k(Y_2;R)\ra\cP MC^{k+1}(Y_1;R)$. As $Y_2\t_{\id_{Y_2},Y_2,f}Y_1\cong Y_1$, we see that $f^*(\Id_{Y_2})=\Id_{Y_1}$.
 
If $g:Y_2\ra Y_3$ is another smooth map of manifolds then as
\begin{equation*}
(V\t_{t,Y_3,g}Y_2)\t_{\pi_{Y_2},Y_2,f}Y_1\cong V\t_{t,Y_3,g\ci f}Y_1,
\end{equation*}
we see from \eq{kh4eq21} that $(g\ci f)^*=f^*\ci g^*:\cP MC^k(Y_3;R)\ra\cP MC^k(Y_1;R)$. Also $\id_Y^*$ is the identity. Thus, pullbacks $f^*$ are contravariantly functorial.
\label{kh4def5}
\end{dfn}

The next proposition is proved in~\S\ref{kh75}.

\begin{prop} $f^*:\cP MC^k(Y_2;R)\ra\cP MC^k(Y_1;R)$ above is well-defined.

\label{kh4prop4}
\end{prop}

The following proposition is proved in \S\ref{kh76}. The proof is similar to that of Theorem \ref{kh4thm1} in~\S\ref{kh72}.

\begin{prop} Let\/ $Y$ be a manifold and\/ $R$ a commutative ring. Then:
\begin{itemize}
\setlength{\itemsep}{0pt}
\setlength{\parsep}{0pt}
\item[{\bf(a)}] Suppose $T,U\subseteq Y$ are open sets. Write $i:T\cap U\hookra T,$ $i':T\cap U\hookra U,$ $j:T\hookra T\cup U,$ $j':U\hookra T\cup U$ for the inclusions. Then for all\/ $k\in\Z$ the following sequence is exact:
\end{itemize}
\e
\xymatrix@C=9.5pt{ 0 \ar[r] & \cP MC^k(T\!\cup\! U;R) \ar[rr]^{j^*\op j^{\prime *}} && {\begin{subarray}{l} \ts \; \cP MC^k(T;R)\\ \ts\op\cP MC^k(U;R) \end{subarray}} \ar[rr]^(0.45){i^*\op -i^{\prime *}} && \cP MC^k(T\!\cap\! U;R). }
\label{kh4eq22}
\e
\begin{itemize}
\setlength{\itemsep}{0pt}
\setlength{\parsep}{0pt}
\item[{\bf(b)}] Suppose $K\subseteq Y$ is closed, $U$ is an open neighbourhood of\/ $K$ in $Y,$ and\/ $\al\in\cP MC^k(U;R)$. Then there exists an open neighbourhood\/ $U'$ of\/ $K$ in $U$ and an element\/ $\be\in\cP MC^k(Y;R)$ with\/ $i^*(\al)=j^*(\be),$ where $i:U'\hookra U$ and\/ $j:U'\hookra Y$ are the inclusions.
\end{itemize}
\label{kh4prop5}
\end{prop}

\begin{dfn} Let $Y$ be a manifold and $k\in\Z$. For all open $U\subseteq Y$ define $\cP\MC^k(Y;R)(U)=\cP MC^k(U;R)$, and for all open $U'\subseteq U\subseteq Y$ define $\rho_{UU'}:\cP\MC^k(Y;R)(U)\ra\cP\MC^k(Y;R)(U')$ by $\rho_{UU'}=i^*:\cP MC^k(U;R)\ra\cP MC^k(U';R)$, with $i:U'\hookra U$ the inclusion. Functoriality of pullbacks $f^*$ in Definition \ref{kh4def5} implies that $\cP\MC^k(Y;R)$ is a presheaf of $R$-modules on $Y$. Proposition \ref{kh4prop5}(a) then means that $\cP\MC^k(Y;R)$ is a strong presheaf, and Proposition \ref{kh4prop5}(b) that $\cP\MC^k(Y;R)$ is soft, and hence c-soft, in the sense of Definitions \ref{kh2def8} and~\ref{kh2def11}.

Write $\MC^k(Y;R)$ for the sheafification of $\cP\MC^k(Y;R)$. Then Theorem \ref{kh2thm4}(f) says that $\MC^k(Y;R)$ is a c-soft sheaf of $R$-modules on $Y$, and hence a soft sheaf of $R$-modules, since c-soft sheaves on manifolds are soft. Define the $R$-module of ({\it integral\/}) {\it M-cochains\/} $MC^k(Y;R)$ to be $MC^k(Y;R)=\MC^k(Y;R)(Y)$, the global sections of $\MC^k(Y;R)$. Since $\cP\MC^k(Y;R)\vert_U=\cP\MC^k(U;R)$ for open $U\subseteq Y$, we have $\MC^k(Y;R)\vert_U=\MC^k(U;R)$, and hence the sheaf $\MC^k(Y;R)$ has $\MC^k(Y;R)(U)=MC^k(U;R)$ for all open~$U\subseteq Y$.

As $\MC^k(Y;R)$ is the sheafification of the strong presheaf $\cP\MC^k(Y;R)$, Theorem \ref{kh2thm4} applies. So Theorem \ref{kh2thm4}(e) gives a canonical isomorphism
\e
MC^k(Y;R)\cong\mathop{\underleftarrow{\lim}\,}\nolimits_{\text{$U:U\subseteq Y$ open, $\bar U$ is compact}}\cP MC^k(U;R),
\label{kh4eq23}
\e
where the right hand side is the inverse limit of $\cP MC^k(U;R)$ over all open $U\subseteq Y$ with closure $\bar U$ compact in $Y$. Such $U$ are partially ordered by inclusion, and if $U'\subseteq U\subseteq Y$ are open with $\bar U,\bar U'$ compact and $i:U'\hookra U$ is the inclusion then Definition \ref{kh4def5} defines $i^*:\cP MC^k(U;R)\ra\cP MC^k(U';R)$, which we use to define the inverse limit.

Write $\Pi:\cP MC^k(Y;R)\ra MC^k(Y;R)$ for the natural projection coming from sheafification. We will use the same notation for elements of $\cP MC^k(Y;R)$, such as generators $[V,n,s,t]$, and for their images under $\Pi$ in $MC^k(Y;R)$. So we have the {\it identity cocycle\/} $\Id_Y\in MC^0(Y;R)$, the image of $\Id_Y\in\cP MC^0(Y;R)$. Applying $\Pi$ shows that equations \eq{kh4eq18}--\eq{kh4eq19} hold in~$MC^k(Y;R)$.

If $Y$ is compact then $U=Y$ is allowed in \eq{kh4eq23}, and $\Pi:\cP MC^k(Y;R)\ra MC^k(Y;R)$ is an isomorphism.

Define $\d(U):\cP\MC^k(Y;R)(U)\ra\cP\MC^{k+1}(Y;R)(U)$ for open $U\subseteq Y$ by $\d(U)=\d:\cP MC^k(U;R)\ra\cP MC^{k+1}(U;R)$. If $U'\subseteq U\subseteq Y$ are open and $i:U'\hookra U$ is the inclusion then Definition \ref{kh4def5} gives $\d\ci i^*=i^*\ci\d:\cP MC^k(U;R)\ra \cP MC^{k+1}(U';R)$. Hence $\d:\cP\MC^k(Y;R)\ra\cP\MC^{k+1}(Y;R)$ is a morphism of presheaves, with $\d\ci\d=0$, and under sheafification descends to a morphism of sheaves $\d:\MC^k(Y;R)\ra\MC^{k+1}(Y;R)$ with $\d\ci\d=0$. Thus $\MC^\bu(Y;R)=\bigl(\MC^*(Y;R),\d\bigr)$ is a complex of soft sheaves of $R$-modules on~$Y$.

Define $\d:MC^k(Y;R)\ra MC^{k+1}(Y;R)$ to be the induced action $\d=\d(Y)$ on global sections. The differentials $\d$ on $\cP MC^*(Y;R),MC^*(Y;R)$ commute with projections $\Pi:\cP MC^k(Y;R)\ra MC^k(Y;R)$. Thus $\d$ on generators $[V,n,s,t]\in MC^k(Y;R)$ is again given by \eq{kh4eq20}. Define the {\it M-cohomology groups\/} (or {\it integral M-cohomology groups\/}) $MH^*(Y;R)$ to be the cohomology of the cochain complex $\bigl(MC^*(Y;R),\d\bigr)$. That is, for $k\in\Z$ we define $R$-modules
\begin{equation*}
MH^k(Y;R)=\frac{\ts \Ker\bigl(\d:MC^k(Y;R)\longra MC^{k+1}(Y;R)\bigr)}{\ts \Im\bigl(\d:MC^{k-1}(Y;R)\longra MC^k(Y;R)\bigr)}\,.
\end{equation*}
As $\d\,\Id_Y=0$ in $\cP MC^1(Y;R)$, we have $\d\,\Id_Y=0$ in $MC^1(Y;R)$, so $\Id_Y$ in $MC^0(Y;R)$ defines an {\it identity cohomology class\/}~$[\Id_Y]\in MH^0(Y;R)$.

Let $f:Y_1\ra Y_2$ be a smooth map of manifolds. Then as in \S\ref{kh25} we have a pushforward presheaf $f_*(\cP\MC^k(Y_1;R))$ with $f_*(\cP\MC^k(Y_1;R))(U_2)=\cP\MC^k(Y_1;R)(U_1)=\cP MC^k(U_1;R)$ for open $U_2\subseteq Y_2$ with $U_1=f^{-1}(U_2)\subseteq Y_1$. Define a presheaf morphism $f_\sh:\cP\MC^k(Y_2;R)\ra f_*(\cP\MC^k(Y_1;R))$ on $Y_2$ by
\e
\begin{split}
f_\sh(U_2)=f\vert_{U_1}^*:\cP\MC^k(Y_2;R)(U_2)=\cP MC^k(U_2;R)&\\
\longra f_*(\cP\MC^k(Y_1;R))(U_2)=\cP MC^k(U_1;R)&.
\end{split}
\label{kh4eq24}
\e
Sheafifying induces a morphism $f_\sh:\MC^k(Y_2;R)\ra f_*(\MC^k(Y_1;R))$ on $Y_2$, where $f_*(\MC^k(Y_1;R))(U_2)\!=\!\MC^k(Y_1;R)(U_1)\!=\!MC^k(U_1;R)$ for $U_1\!=\!f^{-1}(U_2)$, so in particular~$f_*(\MC^k(Y_1;R))(Y_2)=MC^k(Y_1;R)$.

Define the {\it pullback\/} $f^*:MC^k(Y_2;R)\ra MC^k(Y_1;R)$ to be the induced morphism $f^*=f_\sh(Y_2)$ on global sections. Then
\begin{align*}
f^*\ci\Pi=\Pi\ci f^*:\cP MC^k(Y_2;R)\longra MC^k(Y_1;R),
\end{align*}
so $f^*$ on generators $[V,n,s,t]\in MC^k(Y_2;R)$ is again given by~\eq{kh4eq21}.

Since $f^{-1}$ is left adjoint to $f_*$, as in \eq{kh2eq20} $f_\sh$ corresponds to a morphism
\e
f^\sh:f^{-1}(\MC^k(Y_2;R))\longra \MC^k(Y_1;R)
\label{kh4eq25}
\e
of sheaves on $Y_1$, such that the following commutes 
\begin{equation*}
\xymatrix@!0@C=45pt@R=35pt{ 
*+[r]{\cP MC^k(U_2;R)} \ar[d]^{f\vert_{U_1}^*} \ar[rrr]_(0.6)\Pi &&&
MC^k(U_2;R) \ar[d]^(0.36){f\vert_{U_1}^*} \ar[rrrr]_(0.28){\raisebox{-11pt}{$\begin{subarray}{l}\text{project from} \\ \text{direct limit}\end{subarray}$}} &&&& *+[l]{\bigl({\cal P}f^{-1}(\MC^k(Y_2;R))\bigr)(U_1)} \ar[d]_{\text{sheafify}} \\
*+[r]{\cP MC^k(U_1;R)} \ar[rrr]^(0.53)\Pi &&&
{\begin{subarray}{l}\ts \;MC^k(U_1;R)={}\\ \ts \MC^k(Y_1;R)(U_1)\end{subarray}}  &&&& *+[l]{f^{-1}(\MC^k(Y_2;R))\bigr)(U_1),\!\!} \ar[llll]_(0.65){f^\sh(U_1)} }
\end{equation*}
for all open $U_1\subseteq Y_1$, $U_2\subseteq Y_2$ with $f(U_1)\subseteq U_2$, where $\cP f^{-1}(\MC^k(Y_2;R))(U_1)$ is defined as the direct limit of such $MC^k(U_2;R)$ in Definition~\ref{kh2def5}.

If $i:U'\hookra U$ is an inclusion of open sets and $\al\in MC^k(U;R)$, we will often write $\al\vert_{U'}$ for $i^*(\al)\in MC^k(U';R)$.

As $f\vert_{U_1}^*\ci\d=\d\ci f\vert_{U_1}^*:\cP MC^k(U_2;R)\ra\cP MC^{k+1}(U_1;R)$ by Definition \ref{kh4def5}, we deduce that $f^*\ci\d=\d\ci f^*:MC^k(Y_2;R)\ra MC^{k+1}(Y_1;R)$. Hence $f^*:\bigl(MC^*(Y_2;R),\d\bigr)\ra \bigl(MC^*(Y_1;R),\d\bigr)$ is a morphism of cochain complexes, and induces pullbacks $f^*:MH^k(Y_2;R)\ra MH^k(Y_1;R)$ on M-cohomology.

Since pullbacks $f^*$ in Definition \ref{kh4def5} are contravariantly functorial on M-pre\-co\-chains $\cP MC^k(U_i;R)$, we deduce that pullbacks $f^*$ are contravariantly functorial on both M-cochains $MC^k(Y_i;R)$ and M-cohomology~$MH^k(Y_i;R)$.
\label{kh4def6}
\end{dfn}

\begin{rem}{\bf(a)} Let $[V,n,s,t]$ be a generator of $MC^k(Y;R)$. Regarding it as a global section of $\MC^k(Y;R)$, Definition \ref{kh2def6} defines the {\it support\/} $\supp[V,n,s,t]$, a closed subset of $Y$. Now $s^{-1}(0)$ is closed in $V$, so $t[s^{-1}(0)]$ is a closed subset of $Y$ as $(s,t)$ is proper near $\{0\}\t Y$. As for \eq{kh4eq9} we see that
\e
\supp[V,n,s,t]\subseteq t[s^{-1}(0)]\subseteq Y.
\label{kh4eq26}
\e

\noindent{\bf(b)} Consider formal sums $\sum_{i\in I}a_i\,[V_i,n_i,s_i,t_i]$ in $MC^k(Y;R)$, where $I$ is a possibly infinite indexing set, $a_i\in R$ and $[V_i,n_i,s_i,t_i]$ is a generator of $MC^k(Y;R)$ for $i\in I$. We call such a sum {\it locally finite\/} if any $y\in Y$ has an open neighbourhood $U$ such that $\supp[V_i,n_i,s_i,t_i]\cap U\ne\es$ for only finitely many $i\in I$. By \eq{kh4eq26}, this holds if $t_i[s_i^{-1}(0)]\cap U\ne\es$ for only finitely many $i\in I$. This in turn holds provided $\coprod_{i\in I}t_i\vert_{s_i^{-1}(0)}:\coprod_{i\in I}s_i^{-1}(0)\ra Y$ is proper.

By properties of sheaves, if $\sum_{i\in I}a_i\,[V_i,n_i,s_i,t_i]$ is locally finite then there exists a unique element $\al$ of $MC^k(Y;R)$ with
\e
\al\vert_U=\ts\sum_{i\in I}a_i\,[V_i,n_i,s_i,t_i]\vert_U\quad\text{in $MC^k(U;R)$}
\label{kh4eq27}
\e
for all open $U\subseteq Y$ for which $\supp[V_i,n_i,s_i,t_i]\cap U\ne\es$ for only finitely many $i\in I$, so that the sum in \eq{kh4eq27} makes sense as there are only finitely many nonzero terms. We will write $\sum_{i\in I}a_i\,[V_i,n_i,s_i,t_i]=\al$ in $MC^k(Y;R)$. As manifolds are second countable, only countably many terms in a locally finite sum can be nonzero, so we can suppose $I$ is countable.

Using the ideas used to prove Proposition \ref{kh4prop5}(b), one can show that every $\al\in MC^k(Y;R)$ can be written as a locally finite sum $\sum_{i\in I}a_i\,[V_i,n_i,s_i,t_i]$, and we can even take the $[V_i,n_i,s_i,t_i]$ to be compactly-supported. So an alternative way to define $MC^k(Y;R)$ would be as the quotient of the $R$-module of locally finite sums $\sum_{i\in I}a_i\,[V_i,n_i,s_i,t_i]$ by some (rather complicated) relations based on Definition~\ref{kh4def4}(i),(ii).
\label{kh4rem5}
\end{rem}

Here is the analogue of Lemma \ref{kh4lem1}, proved in the same way.

\begin{lem} For any manifold\/ $Y$ we have $\cP MC^k(Y;R)=MC^k(Y;R)=0$ for\/ $k<0,$ so that\/ $MH^k(Y;R)=0$ for\/ $k<0$.
\label{kh4lem3}
\end{lem}

\begin{dfn} Let $Y$ be a manifold, $Z\subseteq Y$ open, and $i:Z\hookra Y$ the inclusion. Define the {\it relative M-cochains\/} $MC^k(Y,Z;R)$ for $k\in\Z$ by
\begin{equation*}
MC^k(Y,Z;R)=MC^k(Y;R)\op MC^{k-1}(Z;R),
\end{equation*}
as in \eq{kh2eq12}, and define $\d:MC^k(Y,Z;R)\ra MC^{k+1}(Y,Z;R)$ by $\d(\al,\be)=(\d\al,i^*(\al)-\d\be)$. Then $\d\ci\d=0:MC^k(Y,Z;R)\ra MC^{k+2}(Y,Z;R)$, since $\d\ci i^*=i^*\ci\d:MC^k(Y;R)\ra MC^{k+1}(Z;R)$, and $\bigl(MC^*(Y,Z;R),\d\bigr)$ is a cochain complex over $R$, the cocone of the morphism $i^*:\bigl(MC^*(Y;R),\d\bigr)\ra\bigl(MC^*(Z;R),\d\bigr)$. Define the {\it relative M-cohomology\/} $MH^k(Y,Z;R)$ to be the $k^{\rm th}$ cohomology group of $\bigl(MC^*(Y,Z;R),\d\bigr)$, as in Axiom~\ref{kh2ax2}(a).

We have a short exact sequence of cochain complexes
\e
{}\!\!\xymatrix@C=7pt{ 0 \ar[r] & \bigl(MC^{*-1}(Z;R),-\d\bigr) \ar[r]^{\d} & \bigl(MC^*(Y,Z;R),\d\bigr) \ar[r]^(0.52){j^*} & \bigl(MC^*(Y;R),\d\bigr) \ar[r] & 0, }\!\!\!\!{}
\label{kh4eq28}
\e
where $\d:MC^{k-1}(Z;R)\ra MC^k(Y,Z;R)$ maps $\d:\be\mapsto(0,\be)$ and $j^*:MC^k(Y,Z;R)\ra MC^k(Y;R)$ maps $j^*:(\al,\be)\mapsto\al$. Thus in the usual way  \cite[Lem.~24.1]{Munk} we get a long exact sequence in cohomology
\begin{equation*}
\xymatrix@C=9.6pt{ \cdots \ar[r] & MH^{k-1}(Z;R) \ar[r]^(0.48)\d & MH^k(Y,Z;R) \ar[r]^(0.52){j^*} & MH^k(Y;R) \ar[r]^(0.47){i^*} & MH^k(Z;R)    \ar[r] & \cdots, }
\end{equation*}
where explicit calculation shows the connecting morphisms are pullbacks $i^*:MH^k(Y;R)\ra MH^k(Z;R)$ from Definition \ref{kh4def6}. This defines the morphisms $\d:MH^k(Z;R)\ra MH^{k+1}(Y,Z;R)$ in Axiom \ref{kh2ax2}(b), and proves Axiom~\ref{kh2ax2}(ii).

Let $f:Y_1\ra Y_2$ be a smooth map of manifolds, and $Z_1\subseteq Y_1$, $Z_2\subseteq Y_2$ be open with $f(Z_1)\subseteq Z_2$; for short we will say that $f:(Y_1,Z_2)\ra (Y_2,Z_2)$ is smooth. Define $f^*:MC^k(Y_2,Z_2;R)\ra MC^k(Y_1,Z_1;R)$ for $k\in\Z$ by $f^*(\al,\be)=\bigl(f^*(\al),f\vert_{Z_1}^*(\be)\bigr)$, using $f^*:MC^k(Y_2;R)\ra MC^k(Y_1;R)$ and $f\vert_{Z_1}^*:MC^{k-1}(Z_2;R)\ra MC^{k-1}(Z_1;R)$ from Definition \ref{kh4def6}. Functoriality of pullbacks and compatibility with $\d$ in Definition \ref{kh4def6} implies that $\d\ci f^*=f^*\ci\d:MC^k(Y_2,Z_2;R)\ra MC^{k+1}(Y_1,Z_1;R)$, so these $f^*$ induce morphisms $f^*:MH^k(Y_2,Z_2;R)\ra MH^k(Y_1,Z_1;R)$ on cohomology, as in Axiom \ref{kh2ax2}(c). These $f^*$ are contravariantly functorial on $MC^k(Y_i,Z_i;R)$ and on $MH^k(Y_i,Z_i;R)$, since they are on $MC^k(Y_i;R)$, proving Axiom \ref{kh2ax2}(i). Note that $j^*$ in \eq{kh4eq28} is pullback by $j=\id_Y:(Y,\es)\ra(Y,Z)$, as in Axiom~\ref{kh2ax2}(ii).

We have a commutative diagram of chain complexes
\begin{equation*}
\xymatrix@C=9pt@R=15pt{ 0 \ar[r] & \bigl(MC^{*-1}(Z_2;R),-\d\bigr) \ar[r]_{\d} \ar[d]^{f\vert_{Z_1}^*} & \bigl(MC^*(Y_2,Z_2;R),\d\bigr) \ar[d]^{f^*} \ar[r]_(0.55){j_2^*} & \bigl(MC^*(Y_2;R),\d\bigr) \ar[d]^{f^*} \ar[r] & 0 \\
0 \ar[r] & \bigl(MC^{*-1}(Z_1;R),-\d\bigr) \ar[r]^{\d} & \bigl(MC^*(Y_1,Z_1;R),\d\bigr) \ar[r]^(0.55){j_1^*} & \bigl(MC^*(Y_1;R),\d\bigr) \ar[r] & 0,\!\!{} }
\end{equation*}
which in the usual way \cite[Th.~24.2]{Munk} induces a commutative diagram
\begin{equation*}
\text{\begin{small}$\displaystyle\xymatrix@C=9pt@R=15pt{ \cdots \ar[r] & MH^{k-1}(Z_2;R) \ar[d]^{f\vert_{Z_1}^*} \ar[r]_(0.48)\d & MH^k(Y_2,Z_2;R) \ar[d]^{f^*} \ar[r]_(0.52){j_2^*} & MH^k(Y_2;R) \ar[d]^{f^*} \ar[r]_(0.47){i_2^*} & MH^k(Z_2;R) \ar[d]^{f\vert_{Z_1}^*} \ar[r] & \cdots \\
\cdots \ar[r] & MH^{k-1}(Z_1;R) \ar[r]^(0.48)\d & MH^k(Y_1,Z_1;R) \ar[r]^(0.52){j_1^*} & MH^k(Y_1;R) \ar[r]^(0.47){i_1^*} & MH^k(Z_1;R) \ar[r] & \cdots,\!\!{} }$\end{small}}
\end{equation*}
proving Axiom~\ref{kh2ax2}(iii).
\label{kh4def7}
\end{dfn}

Our goal is to show that M-cohomology $MH^*(-;R)$ satisfies Axiom \ref{kh2ax2}, so that it is canonically isomorphic to ordinary cohomology by Theorem \ref{kh2thm2}. So far we have defined all the data in Axiom \ref{kh2ax2}(a)--(c), and verified Axiom \ref{kh2ax2}(i)--(iii). The next proposition, proved in \S\ref{kh77}, gives Axiom \ref{kh2ax2}(iv), the homotopy axiom. The main idea is to construct $R$-module morphisms $G:MC^k(Y_2,Z_1;R)\ra MC^{k-1}(Y_1,Z_1;R)$ acting on generators by
\begin{equation*}
G:[V,n,s,t]\longmapsto (-1)^{\dim V}\bigl[V\t_{t,Y_2,g}(Y_1\t[0,1]),n,s\ci\pi_V,\pi_{Y_1}\ci\pi_{Y_1\t[0,1]}\bigr],
\end{equation*}
and show $G$ is a cochain homotopy between $f^*,f^{\prime *}:\bigl(MC^*(Y_2,Z_2;R),\d\bigr)\ra\bigl(MC^*(Y_1,Z_1;R),\d\bigr)$.

\begin{prop} Suppose $Y_1,Y_2$ are manifolds, $Z_1\subseteq Y_1,$ $Z_2\subseteq Y_2$ are open, and\/ $g:Y_1\t[0,1]\ra Y_2$ is a smooth map in $\Manc$ with\/ $g(Z_1\t[0,1])\subseteq Z_2$. Define $f,f':Y_1\ra Y_2$ by $f(y)=g(y,0)$ and\/ $f'(y)=g(y,1)$ for $y\in Y_1$. Then $f^*=f^{\prime *}:MH^k(Y_2,Z_2;R)\ra MH^k(Y_1,Z_1;R)$ for all\/~$k\in\Z$.
\label{kh4prop6}
\end{prop}

Next we prove Axiom \ref{kh2ax2}(v), the excision axiom. The proof uses only softness of the sheaves $\MC^k(Y;R)$ in Definition~\ref{kh4def6}.

\begin{prop} Suppose $Y$ is a manifold, $Z\subseteq Y$ is open, and\/ 
$S\subseteq Z$ is closed in $Y$. Then\/ $j^*:MH^k(Y,Z;R)\ra MH^k(Y\sm S,Z\sm S;R)$ is an isomorphism for all\/ $k\in\Z,$ where\/ $j:Y\sm S\hookra Y$ is the inclusion. 
\label{kh4prop7}
\end{prop}

\begin{proof} To prove $j^*:MH^k(Y,Z;R)\ra MH^k(Y\sm S,Z\sm S;R)$ is injective, suppose $\ga\in MH^k(Y,Z;R)$ with $j^*(\ga)=0$, and let $\ga=[(\al,\be)]$ for $(\al,\be)$ in $MC^k(Y,Z;R)$, so that $\al\in MC^k(Y;R)$, $\be\in MC^{k-1}(Z;R)$ with
\e
\d\al=0\quad\text{and}\quad \al\vert_Z=\d\be.
\label{kh4eq29}
\e
Then $j^*(\al,\be)=\bigl(\al\vert_{Y\sm S},\be\vert_{Z\sm S}\bigr)$ is exact in $\bigl(MC^*(Y\sm S,Z\sm S;R),\d\bigr)$, so there exists $(\al',\be')\in MC^{k-1}(Y\sm S,Z\sm S;R)$ with $\d(\al',\be')=j^*(\al,\be)$, that is, $\al'\in MC^{k-1}(Y\sm S;R)$ and $\be'\in MC^{k-2}(Z\sm S;R)$ with
\e
\d\al'=\al\vert_{Y\sm S} \quad\text{and}\quad \al'\vert_{Z\sm S}-\d\be'=\be\vert_{Z\sm S}.
\label{kh4eq30}
\e

As $Y\sm Z$, $S$ are disjoint closed subsets of $Y$, we may choose an open neighbourhood $U$ of $Y\sm Z$ in $Y$ with $\bar U\cap S=\es$, since manifolds are normal topological spaces. As $\MC^{k-2}(Z;R)$ is a soft sheaf by Definition \ref{kh4def6}, and $\be'$ is a section of $\MC^{k-2}(Z;R)$ over the open set $Z\sm S\subseteq Z$ which contains the closed set $\bar U\cap Z$, there exists $\be''\in MC^{k-2}(Z;R)$ such that for some open neighbourhood $V$ of $\bar U\cap Z$ in $Z\sm S$ we have $\be''\vert_V=\be'\vert_V$, so in particular $\be''\vert_{U\cap Z}=\be'\vert_{U\cap Z}$.

Now $\{U,Z\}$ is an open cover of $Y$, and $\al'\vert_U\in MC^{k-1}(U;R)$, $\be+\d\be''\in MC^{k-1}(Z;R)$ satisfy $\bigl(\al'\vert_U\bigr)\vert_{U\cap Z}=\bigl(\be+\d\be''\bigr)\vert_{U\cap Z}$ in $MC^{k-1}(U\cap Z;R)$, by \eq{kh4eq30} and $\be''\vert_{U\cap Z}=\be'\vert_{U\cap Z}$. So by the sheaf property of $\MC^{k-1}(Y;R)$ there exists a unique $\al''\in MC^{k-1}(Y;R)$ with $\al''\vert_U=\al'\vert_U$ and $\al''\vert_Z=\be+\d\be''$. 

We have $\d\al''\vert_U=\d\al'\vert_U=\al\vert_U$ in $MC^k(U;R)$ by \eq{kh4eq30} and $\d\al''\vert_Z=\d(\be+\d\be'')=\d\be=\al\vert_Z$ in $MC^k(Z;R)$ by \eq{kh4eq29}, so $\d\al''=\al$ in $MC^k(Y;R)$ by the sheaf property of $\MC^k(Y;R)$. Therefore $(\al'',\be'')\in\MC^{k-1}(Y,Z;R)$ with $\d(\al'',\be'')=(\al,\be)$, so $\ga=0$, and $j^*:MH^k(Y,Z;R)\ra MH^k(Y\sm S,Z\sm S;R)$ is injective. A similar proof shows it is surjective.
\end{proof}

The next lemma gives Axiom \ref{kh2ax2}(vi), the additivity axiom. The first equation of \eq{kh4eq31} follows from the sheaf property of $\MC^k(Y;R)$ in Definition \ref{kh4def6}, the second equation of \eq{kh4eq31} follows from the first, and this implies~\eq{kh4eq32}.

\begin{lem} Suppose $Y$ is a manifold with\/ $Y=\coprod_{a\in A}Y_a$ for $A$ a countable indexing set and each\/ $Y_a$ open and closed in $Y$. Let\/ $Z\subseteq Y$ be open, and set\/ $Z_a=Z\cap Y_a$. Then for all\/ $k\ge 0$ we have canonical isomorphisms
\e
MC^k(Y;R)\cong \ts\prod\limits_{a\in A}MC^k(Y_a;R),\;\> MC^k(Y,Z;R)\cong \ts\prod\limits_{a\in A}MC^k(Y_a,Z_a;R),
\label{kh4eq31}
\e
compatible with the morphisms $i_a^*:MC^k(Y;R)\ra MC^k(Y_a;R)$ for $a\in A$ induced by the inclusions $i_a:Y_a\hookra Y$. Hence we also have
\e
MH^k(Y,Z;R)\cong \ts\prod\limits_{a\in A}MH^k(Y_a,Z_a;R),
\label{kh4eq32}
\e
compatible with\/ $i_a^*:MH^k(Y,Z;R)\ra MH^k(Y_a,Z_a;R)$ for $a\in A$.
\label{kh4lem4}
\end{lem}

For Axiom \ref{kh2ax2}(vii), the dimension axiom, observe that when $Y$ is the point $*$, comparing Definitions \ref{kh4def1}, \ref{kh4def4} and \ref{kh4def6} gives canonical isomorphisms
\e
MC_k(*;R)\cong\cP MC^{-k}(*;R)\cong MC^{-k}(*;R),
\label{kh4eq33}
\e
and these isomorphisms identify $\pd:MC_k(*;R)\ra MC_{k-1}(*;R)$ with $\d:\cP MC^{-k}(*;R)\ra\cP MC^{-k+1}(*;R)$ with $\d:MC^{-k}(*;R)\ra MC^{-k+1}(*;R)$. Thus we have canonical isomorphisms $MH_k(*;R)\cong MH^{-k}(*;R)$ for $k\in\Z$, so Axiom \ref{kh2ax2}(vii) follows from Theorem \ref{kh4thm2}. We have now proved all of Axiom \ref{kh2ax2} for M-cohomology. So Theorem \ref{kh2thm2} gives:

\begin{thm} M-cohomology is a cohomology theory of manifolds. There are canonical isomorphisms $MH^k(Y;R)\!\cong\! H^k(Y;R),$ $MH^k(Y,Z;R)\!\cong\! H^k(Y,Z;R)$ for all\/ $Y,Z,k,$ preserving the data $f^*,\d$ and isomorphisms $MH^0(*;R)\cong R\cong H^0(*;R),$ where $H^*(-;R)$ is any other cohomology theory of manifolds over $R,$ such as singular cohomology\/ $H^*_\rsi(Y;R)$ or sheaf cohomology\/~$H^*(Y,R_Y)$.
\label{kh4thm4}
\end{thm}

Next we relate M-cohomology to sheaf cohomology in~\S\ref{kh25}.

\begin{dfn} Let $Y$ be a manifold. Then Example \ref{kh2ex13} defines the constant sheaf of $R$-modules $R_Y$ on $Y$, with $R_Y(U)$ the $R$-module of locally constant functions $s:U\ra R$ for open $U\subseteq Y$. Also Definition \ref{kh4def6} defines a complex $\MC^\bu(Y;R)=\bigl(\MC^*(Y;R),\d\bigr)$ of soft sheaves of $R$-modules $\MC^k(Y;R)$ on $Y$, where $\MC^k(Y;R)=0$ for $k<0$ by Lemma~\ref{kh4lem3}.

We will define a sheaf morphism $i_Y:R_Y\ra\MC^0(Y;R)$. For $U\subseteq Y$ open, write $U=\coprod_{a\in A}U_a$ for the decomposition of $U$ into connected components. Then each $U_a$ is open and closed in $U$, so \eq{kh4eq31} gives
\e
\MC^0(Y;R)(U)=MC^0(U;R)\cong \ts\prod_{a\in A}MC^0(U_a;R).
\label{kh4eq34}
\e
Define an $R$-module morphism $i_Y(U):R_Y(U)\ra\MC^0(Y;R)(U)$ by
\e
i_Y(U):s\longmapsto \ts\sum_{a\in A} s_a\,\Id_{U_a},
\label{kh4eq35}
\e
writing $s\vert_{U_a}=s_a\in R$ for each $a\in A$, since $s$ is locally constant and $U_a$ is connected, so that $s_a\,\Id_{U_a}\in MC^0(U_a;R)$, and $\sum_{a\in A} s_a\,\Id_{U_a}\in \prod_{a\in A}MC^0(U_a;R)$, which we identify with $\MC^0(Y;R)(U)$ by \eq{kh4eq34}.

Since $\Id_U\vert_V=\Id_V$ in $MC^0(V;R)$ for open $V\subseteq U\subseteq Y$, we see that the $i_Y(U):R_Y(U)\ra\MC^0(Y;R)(U)$ are compatible with restriction maps $\rho_{UV}:R_Y(U)\ra R_Y(V)$, $\rho_{UV}:\MC^0(Y;R)(U)\ra\MC^0(Y;R)(V)$, so $i_Y$ is a morphism of sheaves of $R$-modules.
\label{kh4def8}
\end{dfn}

\begin{thm} For each manifold\/ $Y,$ the following is an exact sequence of sheaves of\/ $R$-modules on $Y\!:$
\e
\xymatrix@C=15pt{ 0 \ar[r] & R_Y \ar[r]^(0.3){i_Y} & \MC^0(Y;R) \ar[r]^(0.48)\d & \MC^1(Y;R) \ar[r]^(0.48)\d & \MC^2(Y;R) \ar[r]^(0.62)\d & \cdots. }
\label{kh4eq36}
\e
Hence $\MC^\bu(Y;R)=\bigl(\MC^*(Y;R),\d\bigr)$ is a soft resolution of\/~$R_Y$.
\label{kh4thm5}
\end{thm}

\begin{proof} Equation \eq{kh4eq36} is exact if and only if it is exact on stalks at each $y\in Y$, and so by Definition \ref{kh2def2}, if and only if
\e
\xymatrix@C=21pt{ 0 \ar[r] & R\!=\!R_{Y,y} \ar[r]^(0.37){i_{Y,y}} & 
\mathop{\underrightarrow{\lim}}\limits_{\text{$y\in U\subseteq Y$ open}\!\!\!\!\!\!\!\!\!\!\!\!\!\!\!\!\!\!\!\!\!\!\!\!\!\!\!\!\!} MC^0(U;R) \ar[r]^(0.48){\underrightarrow{\lim}\,\d} & 
\mathop{\underrightarrow{\lim}}\limits_{\text{$y\in U\subseteq Y$ open}\!\!\!\!\!\!\!\!\!\!\!\!\!\!\!\!\!\!\!\!\!\!\!\!\!\!\!\!\!} MC^1(U;R) \ar[r]^(0.68){\underrightarrow{\lim}\,\d} & \cdots }
\label{kh4eq37}
\e
is exact in $R$-modules. As direct limits commute with cohomology, we see that
\e
\begin{split}
H^k&\Bigl(\xymatrix@C=6pt{\cdots \ar[rr] && \mathop{\underrightarrow{\lim}}\limits_{\text{$y\in U\subseteq Y$ open}\!\!\!\!\!\!\!\!\!\!\!\!\!\!\!\!\!\!\!\!\!\!\!\!\!\!\!\!\!} MC^i(U;R) \ar[rrr]^(0.48){\underrightarrow{\lim}\,\d} &&& 
\mathop{\underrightarrow{\lim}}\limits_{\text{$y\in U\subseteq Y$ open}\!\!\!\!\!\!\!\!\!\!\!\!\!\!\!\!\!\!\!\!\!\!\!\!\!\!\!\!\!} MC^{i+1}(U;R) \ar[rr] && \cdots }\Bigr)\\
&\cong \mathop{\underrightarrow{\lim}}\limits_{\text{$y\in U\subseteq Y$ open}\!\!\!\!\!\!\!\!\!\!\!\!\!\!\!\!\!\!\!\!\!\!\!\!\!\!\!\!\!} MH^k(U;R).
\end{split}
\label{kh4eq38}
\e

We may take the limit in \eq{kh4eq38} to be over smaller and smaller balls $U\cong \R^{\dim Y}$ about $y$ in $Y$, so that $MH^0(U;R)\cong H^0(\R^{\dim Y};R)\cong R$ and $MH^k(U;R)\ab\cong H^k(\R^{\dim Y};R)=0$ for $k>0$ by Theorem \ref{kh4thm4}, and the last line of \eq{kh4eq38} is $R$ when $k=0$ and zero for $k>0$. Also $i_{Y,y}$ induces a morphism
\e
(i_{Y,y})_*:R=R_{Y,y}\longra \mathop{\underrightarrow{\lim}}\limits_{\text{$y\in U\subseteq Y$ open}\!\!\!\!\!\!\!\!\!\!\!\!\!\!\!\!\!\!\!\!\!\!\!\!\!\!\!\!\!} MH^0(U;R)\cong R.
\label{kh4eq39}
\e
By definition of $i_Y$ we see that $(i_{Y,y})_*$ maps $1\mapsto \underrightarrow{\lim}\,_U [\Id_U]$, where $[\Id_U]\in MH^0(U;R)$ is identified with $1\in R$ under the isomorphism $MH^0(U;R)\cong R$. Hence \eq{kh4eq39} is an isomorphism. Therefore \eq{kh4eq37} is exact for each $y\in Y$, so \eq{kh4eq36} is exact, as we have to prove. The last part of the theorem follows from Definition~\ref{kh4def6}.
\end{proof}

The material of \S\ref{kh25} now gives an alternative proof, via \eq{kh2eq26}, that
\begin{equation*}
H^k(Y;R)\cong H^k(Y;R_Y)\cong H^k\bigl(MC^*(Y;R),\d\bigr)=MH^k(Y;R).
\end{equation*}

If $f:Y_1\ra Y_2$ is a smooth map of manifolds, we can form a diagram of sheaves of $R$-modules on $Y_1$:
\e
\begin{gathered}
\xymatrix@C=17pt@R=15pt{ 0 \ar[r] & f^{-1}(R_{Y_2}) \ar[d]^{f^\sh}_\cong \ar[r]_(0.4){\raisebox{-9pt}{$\st f^{-1}(i_{Y_2})$}} & f^{-1}(\MC^0(Y_2;R)) \ar[d]^{f^\sh} \ar[r]_{\raisebox{-9pt}{$\st f^{-1}(\d)$}} & f^{-1}(\MC^1(Y_2;R)) \ar[d]^{f^\sh} \ar[r]_(0.65){\raisebox{-9pt}{$\st f^{-1}(\d)$}} & \cdots \\
0 \ar[r] & R_{Y_1} \ar[r]^{i_{Y_1}} & \MC^0(Y_1;R) \ar[r]^\d & \MC^1(Y_1;R) \ar[r]^(0.65)\d & \cdots,\!\!{} }\!\!\!\!\!{}
\end{gathered}
\label{kh4eq40}
\e
where $f^\sh:f^{-1}(R_{Y_2})\ra R_{Y_1}$ is the isomorphism defined in Example \ref{kh2ex13}, by pullback of locally constant functions, and $f^\sh:f^{-1}(\MC^k(Y_2;R))\ra \MC^k(Y_1;R)$ is as in \eq{kh4eq25}. Then $f^*(\Id_{Y_2})=\Id_{Y_1}$ in $MC^0(Y_1;R)$ and the definitions of $i_{Y_1},i_{Y_2}$ in Definition \ref{kh4def8} imply that the left hand square of \eq{kh4eq40} commutes, and $\d\ci f\vert_{U_1}^*=f\vert_{U_1}^*\ci\d:MC^k(U_2;R)\ra MC^{k+1}(U_1;R)$ for open $U_1\subseteq Y_1$, $U_2\subseteq Y_2$ with $f(U_1)\subseteq U_2$ imply that the remaining squares of \eq{kh4eq40} commute.

\subsection{\texorpdfstring{Compactly-supported M-cohomology $MH^*_\cs(Y;R)$}{Compactly-supported M-cohomology}}
\label{kh43}

Next we define compactly-supported M-cohomology.

\begin{dfn} Let $Y$ be a manifold and $k\in\Z$. Then \S\ref{kh42} defined an $R$-module $MC^k(Y;R)$ which is the global sections $MC^k(Y;R)=\MC^k(Y;R)(Y)$ of a c-soft sheaf of $R$-modules $\MC^k(Y;R)$. Define the {\it compactly-supported M-cochains\/} $MC^k_\cs(Y;R)\subseteq MC^k(Y;R)$ to be the $R$-submodule of compactly-supported global sections of $\MC^k(Y;R)$, as in Definition~\ref{kh2def6}.

As in Definition \ref{kh4def6}, $\MC^k(Y;R)$ is the sheafification of the strong presheaf $\cP\MC^k(Y;R)$. Write $\cP MC^k_\cs(Y;R)\subseteq \cP MC^k(Y;R)$ for the $R$-submodule of compactly-supported elements in $\cP MC^k(Y;R)$, which as in Definition \ref{kh2def6} means that $\cP MC^k_\cs(Y;R)$ is the subset of $\al\in\cP MC^k(Y;R)$ such that for some compact $K\subseteq Y$ we have $i^*(\al)=0$ in $\cP MC^k(Y\sm K;R)$, where $i:Y\sm K\hookra Y$ is the inclusion. Then Theorem \ref{kh2thm4}(d) says that
\e
\Pi\vert_{\cP MC^k_\cs(Y;R)}:\cP MC^k_\cs(Y;R)\longra MC^k_\cs(Y;R)
\label{kh4eq41}
\e
is an isomorphism. This gives an alternative description of $MC^k_\cs(Y;R)$ which is more explicit, as it does not involve sheafification.

Note that Lemma \ref{kh4lem3} shows that $MC^k_\cs(Y;R)=0$ for~$k<0$.

Starting with the c-soft sheaf $\MC^k(Y;R)$, Theorem \ref{kh2thm3}(a) defines a flabby cosheaf of $R$-modules on $Y$, which we will write as $\uMC^k_\cs(Y;R)$. Then by definition $MC^k_\cs(Y;R)=\uMC^k_\cs(Y;R)(Y)$ is the global sections of $\uMC^k_\cs(Y;R)$, and more generally $\uMC^k_\cs(Y;R)(U)=MC^k_\cs(U;R)$ for all open~$U\subseteq Y$.

The morphisms $\d:MC^k(Y;R)\!\ra\!MC^{k+1}(Y;R)$ from \S\ref{kh42} with $\d\ci\d=0$ restrict to $\d:MC^k_\cs(Y;R)\!\ra\! MC^{k+1}_\cs(Y;R)$ with $\d\ci\d=0$. Hence $\bigl(MC_\cs^*(Y;R),\d\bigr)$ is a cochain complex. Define the ({\it integral\/}) {\it compactly-supported M-cohomology groups\/} $MH^*_\cs(Y;R)$ to be the cohomology of this cochain complex. 

Write $\Pi:MC^k_\cs(Y;R)\hookra MC^k(Y;R)$ for the inclusion maps. Then $\d\ci\Pi=\Pi\ci\d:MC^k_\cs(Y;R)\ra MC^{k+1}(Y;R)$, so they induce morphisms $\Pi:MH^k_\cs(Y;R)\ra MH^k(Y;R)$ on cohomology, as in Property~\ref{kh2pr1}(a).

If $Y$ is compact then $MC^k_\cs(Y;R)=MC^k(Y;R)$, and so $MH^k_\cs(Y;R)=MH^k(Y;R)$ for all $k\in\Z$, with $\Pi:MH^k_\cs(Y;R)\ra MH^k(Y;R)$ the identity.

Let $f:Y_1\ra Y_2$ be a proper smooth map of manifolds, so that \S\ref{kh42} defines pullback morphisms $f^*:MC^k(Y_2;R)\ra MC^k(Y_1;R)$. If $\al\in MC^k_\cs(Y_2;R)\subseteq MC^k(Y_2;R)$ then $\supp\al\subseteq Y_2$ is compact, so $f^{-1}(\supp\al)\subseteq Y_1$ is compact as $f$ is proper. But $\supp(f^*(\al))\subseteq f^{-1}(\supp\al)$, so $f^*(\al)$ is compactly-supported and lies in $MC^k_\cs(Y_1;R)$. Write $f^*:MC^k_\cs(Y_2;R)\ra MC^k_\cs(Y_1;R)$ for the restriction of $f^*:MC^k(Y_2;R)\ra MC^k(Y_1;R)$. Then $\d\ci f^*=f^*\ci\d:MC^k_\cs(Y_2;R)\ra MC^{k+1}_\cs(Y_1;R)$, as this holds for $MC^*(Y_i;R)$, so the $f^*$ induce morphisms $f^*:MH^k_\cs(Y_2;R)\ra MH^k_\cs(Y_1;R)$, as in Property~\ref{kh2pr1}(b).

Proper pullbacks $f^*$ are contravariantly functorial on $MC^*_\cs(Y_i;R)$ and on $MH^*_\cs(Y_i;R)$, since pullbacks $f^*$ are contravariantly functorial on~$MC^*(Y_i;R)$.

If $U\subseteq Y$ is open and $i:U\hookra Y$ is the inclusion, then as for the morphism $\si_{UY}:\uMC^k_\cs(Y;R)(U)\ra\uMC^k_\cs(Y;R)(Y)$ defined in Theorem \ref{kh2thm3}(a), which is injective as $\uMC^k_\cs(Y;R)$ is flabby, there is an injective {\it pushforward\/} $i_*:MC^k_\cs(U;R)\ra MC^k_\cs(Y;R)$, such that if $\al\in MC^k_\cs(U;R)$ then $i_*(\al)\in MC^k_\cs(Y;R)$ is the unique element with $i_*(\al)\vert_U=\al$ and~$i_*(\al)\vert_{Y\sm\supp\al}=0$.

We have $\d\ci i_*=i_*\ci\d:MC^k_\cs(U;R)\ra MC^{k+1}_\cs(Y;R)$, as $\d:\uMC^k_\cs(Y;R)\ra\uMC^{k+1}_\cs(Y;R)$ is a morphism of cosheaves, so the $i_*$ induce morphisms $i_*:MH^k_\cs(U;R)\ra MH^k_\cs(Y;R)$ on cohomology, as in Property~\ref{kh2pr1}(c).

Pushforwards $i_*$ are covariantly functorial on $MC^*_\cs(-;R)$ and~$MH^*_\cs(-;R)$.

\label{kh4def9}
\end{dfn}

Now Theorem \ref{kh4thm5} shows that $0\ra\MC^0(Y;R)\,{\buildrel\d\over\longra}\,\MC^1(Y;R)\,{\buildrel\d\over\longra}\,\cdots$ is a soft resolution of $R_Y$. So by the material of \S\ref{kh25}, in particular \eq{kh2eq27}, we have natural isomorphisms for all $k\in\Z$
\e
\begin{split}
H^k_\cs(Y,R_Y)&\cong H^k\bigl(\cdots\,{\buildrel\d\over\longra}MC_\cs^i(Y;R)\,{\buildrel\d\over\longra}\,MC_\cs^{i+1}(Y;R)\,{\buildrel\d\over\longra}\,\cdots\bigr)\\
&=MH^k_\cs(Y;R).
\end{split}
\label{kh4eq42}
\e

\begin{thm} Compactly-supported M-cohomology is a compactly-supported cohomology theory of manifolds. As in \eq{kh4eq42} there are canonical isomorphisms $MH^k_\cs(Y;R)\cong H^k_\cs(Y;R)$ for all\/ $Y,k,$ preserving the data $\Pi,f^*,i_*$ described in Property\/ {\rm\ref{kh2pr1}(a)--(c)} and the isomorphisms $MH^0_\cs(*;R)\cong R\cong H^0_\cs(*;R),$ where $H^*_\cs(-;R)$ is any other compactly-supported cohomology theory of manifolds over $R,$ such as compactly-supported singular cohomology\/ $H^*_{\cs,\rsi}(Y;R)$ in Example\/ {\rm\ref{kh2ex8},} or compactly-supported sheaf  cohomology\/ $H^*_\cs(Y,R_Y)$ from\/~{\rm\S\ref{kh25}}.
\label{kh4thm6}
\end{thm}

\begin{proof} We have already constructed canonical isomorphisms $MH^k_\cs(Y;R)\cong H^k_\cs(Y,R_Y)$ in \eq{kh4eq42}. But compactly-supported cohomology theories are known to be canonically isomorphic on sufficiently nice topological spaces (such as manifolds), as in the axiomatic characterizations of Petkova \cite{Petk} and Skljarenko \cite{Sklj1}, for instance. Thus $MH^k_\cs(-;R)$ is also canonically isomorphic to other compactly-supported cohomology theories on manifolds, such as compactly-supported singular cohomology.

That the isomorphisms \eq{kh4eq42} identify $\Pi,i_*$ on $MH^k_\cs(Y;R),MH^k(Y;R)$ and on $H^*_\cs(Y,R_Y),H^*(Y,R_Y)$ is immediate from facts about sheaf cohomology. That they identify pullbacks $f^*$ follows from \eq{kh4eq40} commuting for proper smooth $f:Y_1\ra Y_2$, and properties of sheaf cohomology. That they preserve isomorphisms $MH^0_\cs(*;R)\cong R\cong H^0_\cs(*,R_*)$ follows from the definition of $i_Y$ in Definition \ref{kh4def8} when $Y=*$. This completes the proof.
\end{proof}

\begin{dfn} Let $Y$ be a manifold, $k\in\Z$, and $[V,n,s,t]\in MC^k(Y;R)$ be a generator of $MC^k(Y;R)$, in the sense of Definitions \ref{kh4def4} and \ref{kh4def6}. We say that $[V,n,s,t]$ is a {\it compact generator\/} if $s:V\ra\R^n$ is proper over an open neighbourhood of 0 in $\R^n$. This implies that $(s,t):V\ra\R^n\t Y$ is proper over an open neighbourhood of $\{0\}\t Y$ in $\R^n\t Y$, as assumed in Definition~\ref{kh4def4}.

As $s$ is proper near 0 in $\R^n$, $s^{-1}(0)$ is compact, so Remark \ref{kh4rem5}(a) implies that $\supp[V,n,s,t]$ is compact, and $[V,n,s,t]\!\in\! MC^k_\cs(Y;R)\!\subseteq\! MC^k(Y;R)$. 

\label{kh4def10}
\end{dfn}

The next proposition will be proved in \S\ref{kh78}, using the isomorphism~\eq{kh4eq41}.

\begin{prop} Let\/ $Y$ be a manifold and\/ $k\in\Z$. Then as an $R$-module, $MC^k_\cs(Y;R)$ is generated by compact generators $[V,n,s,t],$ subject only to relations Definition\/ {\rm\ref{kh4def4}(i),(ii)} applied to compact generators.
\label{kh4prop8}
\end{prop}

Thus, we could instead have defined $MC^k_\cs(Y;R)$ by generators and relations, using compact generators, in a very similar way to $MC_k(Y;R)$ in \S\ref{kh41}. The morphisms $\d,f^*,i_*$ in Definition \ref{kh4def9} can all be written explicitly in terms of compact generators. For $\d:MC^k_\cs(Y;R)\ra MC^{k+1}_\cs(Y;R)$, for each compact generator $[V,n,s,t]$ in $MC^k_\cs(Y;R)$ we have
\begin{equation*}
\d[V,n,s,t]=[\pd V,n,s\ci i_V,t\ci i_V],
\end{equation*}
as in \eq{kh4eq20}, where $s\ci i_V:\pd V\ra\R^n$ is proper near 0 in $\R^n$ as $s$ is and $i_V$ is proper, so $[\pd V,n,s\ci i_V,t\ci i_V]$ is a compact generator.

For $f:Y_1\ra Y_2$ a proper smooth map, we write $f^*:MC^k_\cs(Y_2;R)\ra MC^k_\cs(Y_1;R)$ for compact generators $[V,n,s,t]$ in $MC^k_\cs(Y_2;R)$ by
\begin{equation*}
f^*[V,n,s,t]=[V',n,s',t']:=\bigl[V\t_{t,Y_2,f}Y_1,n,s\ci\pi_V,\pi_{Y_1}\bigr],
\end{equation*}
as in \eq{kh4eq21}. Then $s'=s\ci\pi_V:V'\ra\R^n$ is proper near 0 in $\R^n$ as $s$ is, and $\pi_V:V'\ra V$ is proper as $f:Y_1\ra Y_2$ is. So $[V',n,s',t']$ is a compact generator.

For $i:U\hookra Y$ an inclusion of open sets in manifolds, we write $i_*:MC^k_\cs(U;R)\!\ra\! MC^k_\cs(Y;R)$ for compact generators $[V,n,s,t]$ in $MC^k_\cs(U;R)$~by
\e
i_*:[V,n,s,t]=[V,n,s,t],
\label{kh4eq43}
\e
regarding $t$ as a cooriented submersion $t:V\ra U$ for $[V,n,s,t]\in MC^k_\cs(U;R)$, but as a cooriented submersion $t:V\ra Y$ for $[V,n,s,t]\in MC^k_\cs(Y;R)$. Note that despite its simplicity, {\it equation\/ \eq{kh4eq43} does not make sense for non-compact generators\/} $[V,n,s,t]$ in $MC^k(U;R)$, even if $\supp[V,n,s,t]$ is compact, since $(s,t):V\ra \R^n\t U$ proper does not imply $(s,t):V\ra \R^n\t Y$ proper, so $[V,n,s,t]$ may not be a generator of~$MC^k(Y;R)$.

\subsection{\texorpdfstring{Locally finite M-homology $MH_*^\lf(Y;R)$}{Locally finite M-homology}}
\label{kh44}

We define locally finite M-homology, in a similar way to M-cohomology in~\S\ref{kh42}.

\begin{dfn} Let $Y$ be a manifold. Consider quadruples $(V,n,s,t)$, where $V$ is an oriented manifold with corners (i.e.\ a pair $(V,o_V)$ with $V$ an object in $\tManc$ and $o_V$ an orientation on $V$, usually left implicit), and $n\in\N$, and $s:V\ra\R^n$, $t:V\ra Y$ are smooth maps (morphisms in $\tManc$), such that $(s,t):V\ra\R^n\t Y$ is proper over an open neighbourhood of $\{0\}\t Y$ in~$\R^n\t Y$.

Define an equivalence relation $\sim$ on such quadruples by $(V,n,s,t)\sim(V',\ab n',\ab s',\ab t')$ if $n=n'$, and there exists an orientation-preserving diffeomorphism $f:V\ra V'$ with $s=s'\ci f$ and $t=t'\ci f$. Write $[V,n,s,t]$ for the $\sim$-equivalence class of $(V,n,s,t)$. We call $[V,n,s,t]$ a {\it generator}. For each $k\in\Z$, define the {\it locally finite M-prechains\/} $\cP MC_k^\lf(Y;R)$ to be the $R$-module generated by such $[V,n,s,t]$ with $\dim V=n+k$, subject to the relations:
\begin{itemize}
\setlength{\itemsep}{0pt}
\setlength{\parsep}{0pt}
\item[(i)] For each generator $[V,n,s,t]$ and each $i=0,\ldots,n$ we have
\begin{equation*}
[V,n,s,t]=(-1)^{n-i}[V\t\R,n+1,s',t\ci\pi_V]\quad\text{in $\cP MC_k^\lf(Y;R)$,}
\end{equation*}
where writing $s=(s_1,\ldots,s_n):V\ra\R^n$ with $s_j:V\ra\R$ for $j=1,\ldots,n$ and $\pi_V:V\t\R\ra V$, $\pi_\R:V\t\R\ra\R$ for the projections, then 
\begin{equation*}
s'=(s_1\ci\pi_V,\ldots,s_i\ci\pi_V,\pi_\R,s_{i+1}\ci\pi_V,\ldots,s_n\ci\pi_V):V\t\R\longra\R^{n+1},
\end{equation*}
and $V\t\R$ has the product orientation from Assumption \ref{kh3ass6}(f) of the given orientation on $V$ and the standard orientation on~$\R$.
\item[(ii)] Let $I$ be a finite indexing set, $a_i\in R$ for $i\in I$, and $[V_i,n,s_i,t_i]$, $i\in I$ be generators for $MC_k(Y;R)$, all with the same $n$. Suppose there exists an open neighbourhood $X$ of $\{0\}\t Y$ in $\R^n\t Y$, such that $(s_i,t_i):V_i\ra\R^n\t Y$ is proper over $X$ for all $i\in I$, and the following condition holds:
\begin{itemize}
\setlength{\itemsep}{0pt}
\setlength{\parsep}{0pt}
\item[$(*)$] Suppose $(x,y)\!\in\! X$, such that for all $i\!\in\! I$ and $v\!\in\! V_i$ with $(s_i,t_i)(v)\!=\!(x,y)$, we have that $v\in V_i^\ci$ and 
\begin{equation*}
T_v(s_i,t_i):T_vV_i^\ci\longra T_x\R^n\op T_yY
\end{equation*}
is injective. This implies that $(s_i,t_i)\vert_{V_i^\ci}$ is an embedding near $v\in V_i^\ci$. Hence $(s_i,t_i):V_i\ra\R^n\t Y$ is injective near each $v$ in $(s_i,t_i)^{-1}(x,y)$, so $(s_i,t_i)^{-1}(x,y)$ has the discrete topology, and thus is finite as $(s_i,t_i)$ is proper over $X$. Note too that $V_i^\ci$ is an oriented manifold by Assumption \ref{kh3ass6}(j) with $\dim V_i^\ci=n+k$, so $T_vV_i^\ci$ is an oriented vector space of dimension $n+k$. We require that for all oriented $(n+k)$-planes $P\subseteq T_xX\op T_yY=\R^n\op T_yY$, we have
\e
\begin{split}
&\sum_{\begin{subarray}{l} i\in I,\; v\in V_i^\ci:(s_i,t_i)(v)=(x,y),\;  T_v(s_i,t_i)[T_vV_i^\ci]=P \\ \text{$T_v(s_i,t_i):T_vV_i^\ci\,{\buildrel\cong\over\longra}\,P$ is orientation-preserving}\end{subarray}\!\!\!\!\!\!\!} a_i=\\
&\sum_{\begin{subarray}{l} i\in I,\; v\in V_i^\ci:(s_i,t_i)(v)=(x,y),\; T_v(s_i,t_i)[T_vV_i^\ci]=P \\ \text{$T_v(s_i,t_i):T_vV_i^\ci\,{\buildrel\cong\over\longra}\,P$ is orientation-reversing}\end{subarray}\!\!\!\!\!\!\!} a_i\qquad\text{in $R$.}
\end{split}
\label{kh4eq44}
\e
\end{itemize}
Then
\begin{equation*}
\sum_{i\in I}a_i\,[V_i,n,s_i,t_i]=0\qquad\text{in $\cP MC_k^\lf(Y;R)$.}
\end{equation*}
\end{itemize}
These are the same as Definition \ref{kh4def1}(i),(ii), except that they have the properness conditions of Definition~\ref{kh4def4}(i),(ii).

Define $\pd:\cP MC_k^\lf(Y;R)\ra\cP MC_{k-1}^\lf(Y;R)$ to be the unique $R$-linear morphism satisfying \eq{kh4eq3} on generators $[V,n,s,t]$. As in Definition \ref{kh4def1}, it is well-defined, with~$\pd\ci\pd=0:\cP MC_k^\lf(Y;R)\ra\cP MC_{k-2}^\lf(Y;R)$.

Let $f:Y_1\ra Y_2$ be a proper smooth map of manifolds. Define $f_*:\cP MC_k^\lf(Y_1;R)\ra\cP MC_k^\lf(Y_2;R)$ to be the unique $R$-linear morphism acting on generators $[V,n,s,t]$ by \eq{kh4eq7}. We need $f$ proper so that $(s,t):V\ra\R^n\t Y_1$ proper near $\{0\}\t Y_1$ in $\R^n\t Y_1$ implies that $(s,f\ci t):V\ra\R^n\t Y_2$ is proper near $\{0\}\t Y_2$ in $\R^n\t Y_2$. Then $f_*$ is well-defined as in Definition \ref{kh4def2}, with~$\pd\ci f_*=f_*\ci\pd:\cP MC_k^\lf(Y_1;R)\ra\cP MC_{k-1}^\lf(Y_2;R)$.

Let $Y$ be a manifold and $U\subseteq Y$ an open set, with $i:U\hookra Y$ the inclusion. Define $i^*:\cP MC_k^\lf(Y;R)\ra\cP MC_k^\lf(U;R)$ to be the unique $R$-linear morphism acting on generators $[V,n,s,t]$ by
\e
i^*:[V,n,s,t]\longmapsto [V',n,s',t']:=\bigl[t^{-1}(U),n,s\vert_{t^{-1}(U)},t\vert_{t^{-1}(U)}\bigr].
\label{kh4eq45}
\e
Then $(s,t):V\ra\R^n\t Y$ proper near $\{0\}\t Y$ in $\R^n\t Y$ implies that $(s\vert_{t^{-1}(U)},t\vert_{t^{-1}(U)}):t^{-1}(U)\ra\R^n\t U$ is proper near $\{0\}\t U$ in $\R^n\t U$, so the r.h.s.\ of \eq{kh4eq45} is a generator of $\cP MC_k^\lf(U;R)$. Clearly $i^*$ takes relations (i),(ii) in $\cP MC_k^\lf(Y;R)$ to relations (i),(ii) in $\cP MC_k^\lf(U;R)$, so is well-defined. Also $\pd\ci i^*=i^*\ci\pd:\cP MC_k^\lf(Y;R)\ra\cP MC_{k-1}^\lf(U;R)$, and if $j:U'\hookra U$ is another open inclusion then~$j^*\ci i^*=(i\ci j)^*:\cP MC_k^\lf(Y;R)\ra\cP MC_k^\lf(U';R)$.
\label{kh4def11}
\end{dfn}

\begin{rem}{\bf(a)} Note the difference between generators in Definitions \ref{kh4def1} and \ref{kh4def11}: above $(s,t):V\ra\R^n\t Y$ is proper near $\{0\}\t Y$ in $\R^n\t Y$, but Definition \ref{kh4def1} requires the stronger condition that $s:V\ra\R^n$ is proper near 0 in $\R^n$. Thus any generator $[V,n,s,t]$ in $MC_k(Y;R)$ in Definition \ref{kh4def1} is also a generator of $\cP MC_k^\lf(Y;R)$ above.

The only difference between relations Definition \ref{kh4def1}(i),(ii) and Definition \ref{kh4def11}(i),(ii) is that in Definition \ref{kh4def1}(ii) equation \eq{kh4eq2} must hold for all suitable $(x,y)\in X\t Y$ for open $0\in X\subseteq\R^n$, but in Definition \ref{kh4def11}(ii) equation \eq{kh4eq44} must hold for all suitable $(x,y)\in X$ for open $\{0\}\t Y\subseteq X\subseteq\R^n\t Y$. Hence Definition \ref{kh4def1}(i),(ii) imply Definition \ref{kh4def11}(i),(ii). Thus there is a unique $R$-linear morphism $\Pi:MC_k(Y;R)\ra\cP MC_k^\lf(Y;R)$ mapping $\Pi:[V,n,s,t]\mapsto[V,n,s,t]$. These $\Pi$ commute with pushforwards $f_*$ for smooth proper~$f:Y_1\ra Y_2$.

Actually, when applied only to generators $[V,n,s,t]$ with $s:V\ra\R^n$ proper near 0 in $\R^n$ as in Definition \ref{kh4def1}, relations Definition \ref{kh4def1}(ii) and \ref{kh4def11}(ii) are equivalent, since near $s_i^{-1}(0)$ the $V_i$ are confined to a compact subset of~$Y$.
\smallskip

\noindent{\bf(b)} As in Remark \ref{kh4rem2}, if $[V,n,s,t]$ is a generator in $\cP MC_k^\lf(Y;R)$ then
\e
[-V,n,s,t]=-[V,n,s,t],
\label{kh4eq46}
\e
and if $[V_1,n,s_1,t_1]$, $[V_2,n,s_2,t_2]$ are generators in $\cP MC_k^\lf(Y;R)$ then
\e
[V_1\amalg V_2,n,s_1\amalg s_2,t_1\amalg t_2]=[V_1,n,s_1,t_1]+[V_2,n,s_2,t_2].
\label{kh4eq47}
\e
As in Lemma \ref{kh4lem1}, we have $\cP MC_k^\lf(Y;R)=0$ for~$k>\dim Y$.
\label{kh4rem6}
\end{rem}

Now the definition of $\cP MC_k^\lf(Y;R)$ above is similar to that of $\cP MC^k(Y;R)$ in \S\ref{kh42}. The use of orientations/coorientations and submersions is different, but the other details including submersion conditions and relations (i),(ii) are the same. Pullbacks $i^*$ for open inclusions $i:U\hookra Y$ also have identical definitions. Because of this, some results and proofs for $\cP MC^k(-;R)$ translate immediately to results and proofs for $\cP MC_k^\lf(-;R)$ with only cosmetic changes. In particular, as for Proposition \ref{kh4prop5} and its proof in \S\ref{kh76} we can show:

\begin{prop} Let\/ $Y$ be a manifold and\/ $R$ a commutative ring. Then:
\begin{itemize}
\setlength{\itemsep}{0pt}
\setlength{\parsep}{0pt}
\item[{\bf(a)}] Suppose $T,U\subseteq Y$ are open sets. Write $i:T\cap U\hookra T,$ $i':T\cap U\hookra U,$ $j:T\hookra T\cup U,$ $j':U\hookra T\cup U$ for the inclusions. Then for all\/ $k\in\Z$ the following sequence is exact:
\end{itemize}
\begin{equation*}
\xymatrix@C=11.5pt{ 0 \ar[r] & \cP MC^\lf_k(T\!\cup\! U;R) \ar[rr]^{j^*\op j^{\prime *}} && {\begin{subarray}{l} \ts \; \cP MC^\lf_k(T;R)\\ \ts\op\cP MC^\lf_k(U;R) \end{subarray}} \ar[rr]^(0.45){i^*\op -i^{\prime *}} && \cP MC^\lf_k(T\!\cap\! U;R). }
\end{equation*}
\begin{itemize}
\setlength{\itemsep}{0pt}
\setlength{\parsep}{0pt}
\item[{\bf(b)}] Suppose $K\subseteq Y$ is closed, $U$ is an open neighbourhood of\/ $K$ in $Y,$ and\/ $\al\in\cP MC^\lf_k(U;R)$. Then there exists an open neighbourhood\/ $U'$ of\/ $K$ in $U$ and an element\/ $\be\in\cP MC^\lf_k(Y;R)$ with\/ $i^*(\al)=j^*(\be),$ where $i:U'\hookra U$ and\/ $j:U'\hookra Y$ are the inclusions.
\end{itemize}
\label{kh4prop9}
\end{prop}

The next definition follows Definition \ref{kh4def6} closely.

\begin{dfn} Let $Y$ be a manifold and $k\in\Z$. For all open $U\subseteq Y$ define $\cP\MC^\lf_k(Y;R)(U)=\cP MC^\lf_k(U;R)$, and for all open $U'\subseteq U\subseteq Y$ define $\rho_{UU'}:\cP\MC^\lf_k(Y;R)(U)\ra\cP\MC^\lf_k(Y;R)(U')$ by $\rho_{UU'}=i^*:\cP MC^\lf_k(U;R)\ra\cP MC^\lf_k(U';R)$, with $i:U'\hookra U$ the inclusion. Functoriality of pullbacks $i^*$ in Definition \ref{kh4def11} implies that $\cP\MC^\lf_k(Y;R)$ is a presheaf of $R$-modules on $Y$. Proposition \ref{kh4prop9}(a) then means that $\cP\MC^\lf_k(Y;R)$ is a strong presheaf, and Proposition \ref{kh4prop9}(b) that $\cP\MC^\lf_k(Y;R)$ is soft, and hence c-soft.

Write $\MC^\lf_k(Y;R)$ for the sheafification of $\cP\MC^\lf_k(Y;R)$. Then Theorem \ref{kh2thm4}(f) says that $\MC^\lf_k(Y;R)$ is a c-soft sheaf of $R$-modules on $Y$, and hence a soft sheaf, since c-soft sheaves on manifolds are soft. Define the $R$-module of ({\it integral\/}) {\it locally finite M-chains\/} $MC^\lf_k(Y;R)$ by $MC^\lf_k(Y;R)=\MC^\lf_k(Y;R)(Y)$, the global sections of $\MC^\lf_k(Y;R)$. Since $\cP\MC^\lf_k(Y;R)\vert_U=\cP\MC^\lf_k(U;R)$ for open $U\subseteq Y$, we have $\MC^\lf_k(Y;R)\vert_U=\MC^\lf_k(U;R)$, and hence the sheaf $\MC^\lf_k(Y;R)$ has $\MC^\lf_k(Y;R)(U)=MC^\lf_k(U;R)$ for all open~$U\subseteq Y$.

As $\MC^\lf_k(Y;R)$ is the sheafification of the strong presheaf $\cP\MC^\lf_k(Y;R)$, Theorem \ref{kh2thm4} applies. So Theorem \ref{kh2thm4}(e) gives a canonical isomorphism
\e
MC^\lf_k(Y;R)\cong\mathop{\underleftarrow{\lim}\,}\nolimits_{\text{$U:U\subseteq Y$ open, $\bar U$ is compact}}\cP MC^\lf_k(U;R),
\label{kh4eq48}
\e
where the right hand side is the inverse limit of $\cP MC^\lf_k(U;R)$ over all open $U\subseteq Y$ with closure $\bar U$ compact in $Y$. Such $U$ are partially ordered by inclusion, and if $U'\subseteq U\subseteq Y$ are open with $\bar U,\bar U'$ compact and $i:U'\hookra U$ is the inclusion then Definition \ref{kh4def11} defines $i^*:\cP MC^\lf_k(U;R)\ra\cP MC^\lf_k(U';R)$, which we use to define the inverse limit.

Write $\Pi:\cP MC^\lf_k(Y;R)\ra MC^\lf_k(Y;R)$ for the natural projection coming from sheafification. We will use the same notation for elements of $\cP MC_k^\lf(Y;R)$, such as generators $[V,n,s,t]$, and for their images under $\Pi$ in $MC_k^\lf(Y;R)$. Applying $\Pi$ shows that equations \eq{kh4eq46}--\eq{kh4eq47} hold in~$MC_k^\lf(Y;R)$.

If $Y$ is compact then $U=Y$ is allowed in \eq{kh4eq48}, and $\Pi:\cP MC^\lf_k(Y;R)\ra MC^\lf_k(Y;R)$ is an isomorphism.

The morphisms $\pd:\cP MC_k^\lf(U;R)\ra\cP MC_{k-1}^\lf(U;R)$ in Definition \ref{kh4def11} for open $U\subseteq Y$ with $\pd\ci\pd=0$ induce presheaf morphisms $\pd:\cP\MC_k^\lf(Y;R)\ra\cP\MC_{k-1}^\lf(Y;R)$ with $\pd\ci\pd=0$, and under sheafification these descend to sheaf morphisms $\pd:\MC_k^\lf(Y;R)\ra\MC_{k-1}^\lf(Y;R)$ with $\pd\ci\pd=0$. So $\MC_\bu^\lf(Y;R)=\bigl(\MC_*^\lf(Y;R),\pd\bigr)$ is a complex of soft sheaves of $R$-modules on~$Y$.

Define $\pd:MC_k^\lf(Y;R)\ra MC_{k-1}^\lf(Y;R)$ to be the induced action $\pd=\pd(Y)$ on global sections. Then $\Pi\ci\pd=\pd\ci\Pi:\cP MC_k^\lf(Y;R)\ra MC_{k-1}^\lf(Y;R)$, so $\pd$ on generators $[V,n,s,t]$ in $MC_k^\lf(Y;R)$ is again given by \eq{kh4eq3}. Define the ({\it integral\/}) {\it locally finite M-homology groups\/} $MH_*^\lf(Y;R)$ to be the homology of~$\bigl(MC_*^\lf(Y;R),\pd\bigr)$. 

Write $\Pi:MC_k(Y;R)\ra MC_k^\lf(Y;R)$ for the composition of
\e
\Pi:MC_k(Y;R)\longra\cP MC_k^\lf(Y;R),\;\> \Pi:\cP MC_k^\lf(Y;R)\longra MC_k^\lf(Y;R).
\label{kh4eq49}
\e
Then $\Pi\ci\pd=\pd\ci\Pi:MC_k(Y;R)\ra MC_{k-1}^\lf(Y;R)$, so we have induced morphisms $\Pi:MH_k(Y;R)\ra MH_k^\lf(Y;R)$, as in Property \ref{kh2pr2}(a). If $Y$ is compact then as above \eq{kh4eq49} are isomorphisms, so $\Pi:MC_k(Y;R)\ra MC_k^\lf(Y;R)$ and $\Pi:MH_k(Y;R)\ra MH_k^\lf(Y;R)$ are isomorphisms.

If $Y$ is oriented with $\dim Y=m$, define the {\it fundamental cycle\/} $[Y]=[Y,0,0,\id_Y]$ in $MC_m^\lf(Y;R)$. Then $\pd[Y]=0$ as $\pd Y=\es$, so taking homology gives the {\it fundamental class\/}~$[[Y]]\in MH_m^\lf(Y;R)$.

Let $f:Y_1\ra Y_2$ be a proper smooth map of manifolds. As in \eq{kh4eq24}, define a presheaf morphism $f_\sh:f_*(\cP\MC^\lf_k(Y_1;R))\ra\cP\MC^\lf_k(Y_2;R)$ on $Y_2$ by
\begin{align*}
f_\sh(U_2)=(f\vert_{U_1})_*:f_*(\cP\MC^\lf_k(Y_1;R))(U_2)=\cP MC^\lf_k(U_1;R)&\\
\longra \cP\MC^\lf_k(Y_2;R)(U_2)=\cP MC^\lf_k(U_2;R)&
\end{align*}
for open $U_2\subseteq Y_2$ with $U_1=f^{-1}(U_2)\subseteq Y_1$. Note that $f$ proper implies that $f\vert_{U_1}:U_1\ra U_2$ is proper, so $(f\vert_{U_1})_*:\cP MC^\lf_k(U_1;R)\ra\cP MC^\lf_k(U_2;R)$ is well defined. Sheafifying induces  $f_\sh:f_*(\MC^\lf_k(Y_1;R))\ra\MC^\lf_k(Y_2;R)$ on $Y_2$, where $f_*(\MC^\lf_k(Y_1;R))(U_2)\!=\!\MC^\lf_k(Y_1;R)(U_1)\!=\!MC^\lf_k(U_1;R)$ for $U_1\!=\!f^{-1}(U_2)$, so in particular~$f_*(\MC^\lf_k(Y_1;R))(Y_2)=MC^\lf_k(Y_1;R)$.

Define the {\it pushforward\/} $f_*:MC^\lf_k(Y_1;R)\ra MC^\lf_k(Y_2;R)$ to be the induced morphism $f_*=f_\sh(Y_2)$ on global sections. Then
\begin{align*}
f_*\ci\Pi=\Pi\ci f_*:\cP MC^\lf_k(Y_1;R)\longra MC^\lf_k(Y_2;R),
\end{align*}
so $f_*$ on generators $[V,n,s,t]\in MC^\lf_k(Y_1;R)$ is again given by~\eq{kh4eq7}.

As in \eq{kh2eq20} and \eq{kh4eq25}, $f_\sh$ corresponds to a morphism of sheaves on~$Y_1$
\begin{equation*}
f^\sh:\MC^\lf_k(Y_1;R)\longra f^{-1}(\MC^\lf_k(Y_2;R)).
\end{equation*}

Since $(f\vert_{U_1})_*\!\ci\!\pd\!=\!\pd\!\ci\!(f\vert_{U_1})_*:\cP MC^\lf_k(U_1;R)\!\ra\!\cP MC^\lf_{k-1}(U_2;R)$ for $U_2\!\subseteq\! Y$ open with $U_1\!=\!f^{-1}(U_2)$ by Definition \ref{kh4def11}, we see that $f_*\!\ci\!\pd\!=\!\pd\!\ci\! f_*:MC^\lf_k(Y_1;R)\!\ra\! MC^\lf_{k-1}(Y_2;R)$. Thus $f_*:\bigl(MC_*^\lf(Y_1;R),\pd\bigr)\!\ra\! \bigl(MC_*^\lf(Y_2;R),\pd\bigr)$ is a morphism of chain complexes, and induces pushforwards $f_*\!:\!MH^\lf_k(Y_1;R)\!\ra\! MH^\lf_k(Y_2;R)$ on locally finite M-homology.

Since pushforwards $f_*$ in Definition \ref{kh4def11} are covariantly functorial, we deduce that pushforwards $f_*$ are covariantly functorial on both locally finite M-chains $MC^\lf_k(Y_i;R)$ and M-homology~$MH^\lf_k(Y_i;R)$.

If $i:U\hookra Y$ is an inclusion of open sets, write $i^*:MC^\lf_k(Y;R)\ra MC^\lf_k(U;R)$ for the restriction map $\rho_{YU}:\MC^\lf_k(Y;R)(Y)\ra\MC^\lf_k(Y;R)(U)$ in the sheaf $\MC^\lf_k(Y;R)$. Then $\Pi\ci i^*=i^*\ci\Pi:\cP MC^\lf_k(Y;R)\ra MC^\lf_k(U;R)$, so $i^*$ acts on generators $[V,n,s,t]\in MC^\lf_k(Y;R)$ as in \eq{kh4eq45}. 

When $\al\in MC_k^\lf(Y;R)$, we may write $\al\vert_U$ for $i^*(\al)\in MC_k^\lf(U;R)$. As 
$\pd:\MC_k^\lf(Y;R)\ra\MC_{k-1}^\lf(Y;R)$ is a sheaf morphism we have $i^*\ci\pd=\pd\ci i^*:MC_k^\lf(Y;R)\ra MC_{k-1}^\lf(U;R)$, so the $i^*$ induce contravariantly functorial {\it pullbacks\/} $i^*:MH_k^\lf(Y;R)\ra MH_k^\lf(U;R)$, as in Property~\ref{kh2pr2}(c).

\label{kh4def12}
\end{dfn}

\begin{rem}{\bf(a)} Let $[V,n,s,t]$ be a generator of $MC_k^\lf(Y;R)$. Regarding it as a global section of $\MC_k^\lf(Y;R)$, Definition \ref{kh2def6} defines the {\it support\/} $\supp[V,n,s,t]$, a closed subset of $Y$. As for \eq{kh4eq9} and \eq{kh4eq26}, we see that
\begin{equation*}
\supp[V,n,s,t]\subseteq t[s^{-1}(0)]\subseteq Y.
\end{equation*}

\noindent{\bf(b)} Consider formal sums $\sum_{i\in I}a_i\,[V_i,n_i,s_i,t_i]$ in $MC_k^\lf(Y;R)$, where $I$ is a possibly infinite indexing set, $a_i\in R$ and $[V_i,n_i,s_i,t_i]$ is a generator of $MC_k^\lf(Y;R)$ for $i\in I$. We call such a sum {\it locally finite\/} if any $y\in Y$ has an open neighbourhood $U$ such that $\supp[V_i,n_i,s_i,t_i]\cap U\ne\es$ for only finitely many $i\in I$. Then as in Remark \ref{kh4rem5}(b) there exists a unique element $\al$ of $MC_k^\lf(Y;R)$ with
\begin{equation*}
\al\vert_U=\ts\sum_{i\in I}a_i\,[V_i,n_i,s_i,t_i]\vert_U\quad\text{in $MC^\lf_k(U;R)$}
\end{equation*}
for all open $U\subseteq Y$ with $\supp[V_i,n_i,s_i,t_i]\cap U\ne\es$ for only finitely many $i\in I$. We will write $\sum_{i\in I}a_i\,[V_i,n_i,s_i,t_i]=\al$ in~$MC_k^\lf(Y;R)$. 
\label{kh4rem7}
\end{rem}

The first part of next proposition is the analogue of Proposition \ref{kh4prop8}, and can be proved as in \S\ref{kh78} with only cosmetic modifications. To see this, note that $MC_k(Y;R)$ in \S\ref{kh41} is the analogue of the module in Proposition \ref{kh4prop8} spanned by compact generators $[V,n,s,t],$ subject to relations Definition \ref{kh4def4}(i),(ii), and $MC_k^\lf(Y;R)_\cs$ is the analogue of $MC^k_\cs(Y;R)$ in \S\ref{kh43}. The second part, the analogue of \eq{kh4eq41}, follows from Theorem \ref{kh2thm4}(d), since as above $\MC^\lf_k(Y;R)$ is the sheafification of the strong presheaf~$\cP\MC^\lf_k(Y;R)$.

\begin{prop} Let\/ $Y$ be a manifold and\/ $k\in\Z$. Then the morphism $\Pi:MC_k(Y;R)\ra MC_k^\lf(Y;R)$ above is injective, with image the $R$-submodule $MC_k^\lf(Y;R)_\cs$ of compactly-supported elements $\al\in MC_k^\lf(Y;R)$.

Also $\Pi:\cP MC_k^\lf(Y;R)\ra MC_k^\lf(Y;R)$ above restricts to an isomorphism
\begin{equation*}
\Pi\vert_{\cdots}: \cP MC_k^\lf(Y;R)_\cs\,{\buildrel\cong\over\longra}\,MC_k^\lf(Y;R)_\cs.
\end{equation*}

\label{kh4prop10}
\end{prop}

Thus we can regard $MC_k(Y;R)$ as an $R$-submodule of~$MC_k^\lf(Y;R)$.

Combining Proposition \ref{kh4prop10} with Theorems \ref{kh2thm3} and \ref{kh4thm1} and Definition \ref{kh4def12} yields the following corollary, where the isomorphism \eq{kh4eq50} follows from \eq{kh2eq32} and the definition of relative M-chains $MC_k(Y,Y\sm\{y\};R)$ in Definition~\ref{kh4def3}. 

\begin{cor} Let\/ $Y$ be a manifold and\/ $k\in\Z$. Then the flabby cosheaf\/ $\uMC_k(Y;R)$ of\/ $R$-modules on $Y$ in Theorem\/ {\rm\ref{kh4thm1}} and the c-soft sheaf\/ $\MC^\lf_k(Y;R)$ of\/ $R$-modules on $Y$ in Definition\/ {\rm\ref{kh4def12}} are canonically related as in Theorem\/ {\rm\ref{kh2thm3}}. The stalks $\MC_k^\lf(Y;R)_y$ for\/ $y\in Y$ have canonical isomorphisms 
\e
\MC_k^\lf(Y;R)_y\cong MC_k(Y,Y\sm\{y\};R)
\label{kh4eq50}
\e
which identify the stalk morphism $\pd_y:\MC_k^\lf(Y;R)_y\ra\MC_{k-1}^\lf(Y;R)_y$ with\/ $\pd:MC_k(Y,Y\sm\{y\};R)\ra MC_{k-1}(Y,Y\sm\{y\};R)$ from\/~{\rm\S\ref{kh41}}.
\label{kh4cor1}
\end{cor}

We relate locally finite M-homology $MH_*^\lf(Y;R)$ at the chain level to locally finite smooth singular homology $H_*^{\lf,\ssi}(Y;R)$ in Example \ref{kh2ex12}, and to locally finite sheaf smooth singular homology $\hat H_k^{\lf,\ssi}(Y;R)$ in Example~\ref{kh2ex17}.

\begin{ex} Let $Y$ be a manifold, and let $\bigl(C_*^{\lf,\ssi}(Y;R),\pd\bigr)$ and $H_*^{\lf,\ssi}(Y;R)$ be as in Example \ref{kh2ex12}, so that elements of $C_k^{\lf,\ssi}(Y;R)$ are locally finite sums $\sum_{i\in I}\rho_i\,\si_i$ with $\rho_i\in R$ and $\si_i:\De_k\ra Y$ a smooth map in $\Manc$ for $i\in I$. As for the morphisms $F_\ssi^\Mh:C_k^\ssi(Y;R)\ra MC_k(Y;R)$ in Example \ref{kh4ex1}, define $R$-module morphisms $F_{\lf,\ssi}^{\lf,\Mh}:C_k^{\lf,\ssi}(Y;R)\ra MC_k^\lf(Y;R)$ for $k=0,1,\ldots$ by
\e
F_{\lf,\ssi}^{\lf,\Mh}:\ts\sum_{i\in I}\rho_i\,\si_i\longmapsto\ts\sum_{i\in I}\rho_i\,[\De_k,0,0,\si_i].
\label{kh4eq51}
\e
The r.h.s.\ of \eq{kh4eq51} is a locally finite sum in $MC_k^\lf(Y;R)$, and so is well-defined as in Remark \ref{kh4rem7}(b). 

The argument of \eq{kh4eq13} shows that $\pd\ci F_{\lf,\ssi}^{\lf,\Mh}=F_{\lf,\ssi}^{\lf,\Mh}\ci\pd:C_k^{\lf,\ssi}(Y;R)\ra MC_{k-1}^\lf(Y;R)$. Thus we have induced morphisms $F_{\lf,\ssi}^{\lf,\Mh}:H_k^{\lf,\ssi}(Y;R)\ra MH_k^\lf(Y;R)$ for $k=0,1,\ldots,$ which Theorem \ref{kh4thm7} will show are isomorphisms.

Suppose $f:Y_1\ra Y_2$ is a proper smooth map of manifolds. Then comparing \eq{kh4eq51} with the definitions of pushforwards $f_*$ on $C_k^{\lf,\ssi}(Y_a;R)$ in \S\ref{kh24} and $MC_k^\lf(Y_a;R)$ above, we see that $f_*\ci F_{\lf,\ssi}^{\lf,\Mh}=F_{\lf,\ssi}^{\lf,\Mh}\ci f_*:C_k^{\lf,\ssi}(Y_1;R)\ra MC_k^\lf(Y_2;R)$, and thus $f_*\ci F_{\lf,\ssi}^{\lf,\Mh}=F_{\lf,\ssi}^{\lf,\Mh}\ci f_*:H_k^{\lf,\ssi}(Y_1;R)\ra MH_k^\lf(Y_2;R)$.

\label{kh4ex3}
\end{ex}

\begin{ex} Let $Y$ be a manifold. Example \ref{kh4ex2} constructed a commutative diagram \eq{kh4eq15} of flabby cosheaves $\hat\ucC{}_k^\ssi(Y;R),\uMC_k(Y;R)$ on $Y$. Applying Theorem \ref{kh2thm3}(b),(c) gives a corresponding diagram of (c-)soft sheaves of $R$-modules on $Y$. By Example \ref{kh2ex17} and Corollary \ref{kh4cor1}, the soft sheaves corresponding to $\hat\ucC{}_k^\ssi(Y;R),\uMC_k(Y;R)$ are $\hat\cC{}^{\lf,\ssi}_k(Y;R),\MC_k^\lf(Y;R)$. Thus the commutative diagram of soft sheaves corresponding to \eq{kh4eq15} under Theorem \ref{kh2thm3} is
\e
\begin{gathered}
\xymatrix@C=18pt@R=15pt{ \cdots \ar[r]_(0.25)\pd  & \hat\cC{}^{\lf,\ssi}_{k+1}(Y;R) \ar[r]_\pd \ar[d]^{\hat {\ul{F\!}\,}{}_{\lf,\ssi}^{\lf,\Mh}} & \hat\cC{}^{\lf,\ssi}_k(Y;R) \ar[r]_\pd \ar[d]^{\hat {\ul{F\!}\,}{}_{\lf,\ssi}^{\lf,\Mh}} & \hat\cC{}^{\lf,\ssi}_{k-1}(Y;R) \ar[r]_(0.65)\pd \ar[d]^{\hat {\ul{F\!}\,}{}_{\lf,\ssi}^{\lf,\Mh}} &  \cdots \\
\cdots \ar[r]^(0.25)\pd  & \MC_{k+1}^\lf(Y;R) \ar[r]^\pd & \MC_k^\lf(Y;R) \ar[r]^\pd & \MC_{k-1}^\lf(Y;R) \ar[r]^(0.65)\pd &  \cdots,\!\!{} }
\end{gathered}
\label{kh4eq52}
\e
writing $\hat {\ul{F\!}\,}{}_{\lf,\ssi}^{\lf,\Mh}$ for the sheaf morphisms corresponding to $\hat {\ul{F\!}\,}{}_\ssi^\Mh$ in Example~\ref{kh4ex2}.

By definition, the global sections over $Y$ of $\hat\cC{}^{\lf,\ssi}_k(Y;R),\ab\MC_k^\lf(Y;R)$ are $\hat C_k^{\lf,\ssi}(Y;R),\ab MC_k^\lf(Y;R)$. Therefore taking global sections in \eq{kh4eq52} gives a commutative diagram of $R$-modules
\e
\begin{gathered}
\xymatrix@C=18pt@R=15pt{ \cdots \ar[r]_(0.25)\pd  & \hat C^{\lf,\ssi}_{k+1}(Y;R) \ar[r]_\pd \ar[d]^{\hat F{}_{\lf,\ssi}^{\lf,\Mh}} & \hat C^{\lf,\ssi}_k(Y;R) \ar[r]_\pd \ar[d]^{\hat F{}_{\lf,\ssi}^{\lf,\Mh}} & \hat C^{\lf,\ssi}_{k-1}(Y;R) \ar[r]_(0.65)\pd \ar[d]^{\hat F{}^{\lf,\Mh}_{\lf,\ssi}} &  \cdots \\
\cdots \ar[r]^(0.25)\pd  & MC_{k+1}^\lf(Y;R) \ar[r]^\pd & MC_k^\lf(Y;R) \ar[r]^\pd & MC_{k-1}^\lf(Y;R) \ar[r]^(0.65)\pd &  \cdots,\!\!{} }
\end{gathered}
\label{kh4eq53}
\e
writing $\hat F{}_{\lf,\ssi}^{\lf,\Mh}=\hat {\ul{F\!}\,}{}_{\lf,\ssi}^{\lf,\Mh}(Y)$. As \eq{kh4eq53} commutes, these $\hat F{}_{\lf,\ssi}^{\lf,\Mh}:\hat C^{\lf,\ssi}_k(Y;R)\ra MC_k^\lf(Y;R)$ induce morphisms $\hat F{}_{\lf,\ssi}^{\lf,\Mh}:\hat H^{\lf,\ssi}_k(Y;R)\ra MH_k^\lf(Y;R)$ on homology, which Theorem \ref{kh4thm7} will show are isomorphisms.

In Example \ref{kh2ex17} we defined morphisms $\Pi_0^\lf:C_k^{\lf,\ssi}(Y;R)\ra \hat C_k^{\lf,\ssi}(Y;R)$, by extending the morphisms $\Pi_0:C_k^\ssi(Y;R)\ra\hat C_k^\ssi(Y;R)$ from Example \ref{kh2ex16} over locally finite sums. Clearly, from the definitions we have
\e
F_{\lf,\ssi}^{\lf,\Mh}=\hat F{}_{\lf,\ssi}^{\lf,\Mh}\ci\Pi_0^\lf:C_k^{\lf,\ssi}(Y;R)\longra MC_k^\lf(Y;R),
\label{kh4eq54}
\e
and so the analogue holds on homology.

Suppose $i:U\hookra Y$ is an inclusion of open sets. Then as in Property \ref{kh2pr2}(c) there is a natural morphism $i^*:H_k^\lf(Y;R)\ra H_k^\lf(U;R)$ on locally finite homology. In Example \ref{kh2ex17} we noted that for $\hat H_k^{\lf,\ssi}(-;R)$ this is given on the chain level by $i^*=\si_{YU}:\hat C_k^{\lf,\ssi}(Y;R)\ra\hat C_k^{\lf,\ssi}(U;R)$, where $\si_{YU}$ is the restriction morphism in the sheaf $\hat\cC{}^{\lf,\ssi}_k(Y;R)$. Since the pullback $i^*:MC_k^\lf(Y;R)\ra MC_k^\lf(U;R)$ in Definition \ref{kh4def12} is also the restriction morphism $\si_{YU}$ in the sheaf $\MC_k^\lf(Y;R)$, as \eq{kh4eq52} commutes we see that
\begin{equation*}
\hat F{}_{\lf,\ssi}^{\lf,\Mh}\ci i^*=i^*\ci \hat F{}_{\lf,\ssi}^{\lf,\Mh}:\hat C_k^{\lf,\ssi}(Y;R)\longra MC_k^\lf(U;R),
\end{equation*}
and so the analogue holds on homology.
\label{kh4ex4}
\end{ex}

We can now prove the main result of this section:

\begin{thm} Locally finite M-homology is a locally finite homology theory of manifolds. That is, there are canonical isomorphisms $H_k^\lf(Y;R)\cong MH_k^\lf(Y;R)$ for all\/ $Y,k,$ preserving the data $\Pi,f_*,i^*$ described in Property\/ {\rm\ref{kh2pr2}(a)--(c)} and the isomorphisms $H_0^\lf(*;R)\cong R\cong MH_0^\lf(*;R),$ where $H_*^\lf(-;R)$ is any other locally finite homology theory of manifolds over~$R$. 

For locally finite smooth singular homology $H_*^{\lf,\ssi}(Y;R)$ in Example\/ {\rm\ref{kh2ex12}} and locally finite sheaf smooth singular homology\/ $\hat H_*^{\lf,\ssi}(Y;R)$ in Example\/ {\rm\ref{kh2ex17},} the canonical isomorphisms $F_{\lf,\ssi}^{\lf,\Mh}:H_k^{\lf,\ssi}(Y;R)\ra MH_k^\lf(Y;R)$ and\/ $\hat F{}_{\lf,\ssi}^{\lf,\Mh}:\hat H^{\lf,\ssi}_k(Y;R)\ra MH_k^\lf(Y;R)$ are as in Examples\/ {\rm\ref{kh4ex3}} and\/~{\rm\ref{kh4ex4}}.
\label{kh4thm7}
\end{thm}

\begin{proof} Let $Y$ be a manifold and $y\in Y$. Taking stalks at $y$ in \eq{kh4eq52} and using Corollary \ref{kh4cor1} and its analogue for $\hat\cC{}^{\lf,\ssi}_*(Y;R)$ gives a commutative diagram
\e
\begin{gathered}
\xymatrix@C=20pt@R=15pt{ \cdots \ar[r]_(0.25)\pd  & \hat C^{\lf,\ssi}_k(Y,Y\sm\{y\};R) \ar[r]_\pd \ar[d]^{\hat {\ul{F\!}\,}{}_{\lf,\ssi}^{\lf,\Mh}\vert_y} & \hat C^{\lf,\ssi}_{k-1}(Y,Y\sm\{y\};R) \ar[r]_(0.65)\pd \ar[d]^{\hat {\ul{F\!}\,}{}_{\lf,\ssi}^{\lf,\Mh}\vert_y} &  \cdots \\
\cdots \ar[r]^(0.25)\pd   &  MC_k(Y,Y\sm\{y\};R) \ar[r]^\pd &  MC_{k-1}(Y,Y\sm\{y\};R) \ar[r]^(0.65)\pd &  \cdots.\!\!{} }
\end{gathered}
\label{kh4eq55}
\e
The rows of \eq{kh4eq55} have homology $\hat H^{\lf,\ssi}_*(Y,Y\sm\{y\};R),MH_*(Y,Y\sm\{y\};R)$, and the columns of \eq{kh4eq55} induce the natural morphisms $\hat H^{\lf,\ssi}_*(Y,Y\sm\{y\};R)\ra MH_*(Y,Y\sm\{y\};R)$, which are isomorphisms by Theorem \ref{kh4thm3}.

Thus, regarding \eq{kh4eq52} as a morphism $\hat {\ul{F\!}\,}{}_{\lf,\ssi}^{\lf,\Mh}:\hat\cC{}^{\lf,\ssi}_\bu(Y;R)\ra \MC_\bu^\lf(Y;R)$ of complexes of soft sheaves of $R$-modules on $Y$, this morphism induces isomorphisms on the homology of stalks at $y\in Y$ for all $y\in Y$, so $\hat {\ul{F\!}\,}{}_{\lf,\ssi}^{\lf,\Mh}$ is a quasi-isomorphism. As the sheaves are soft, this implies that $\hat {\ul{F\!}\,}{}_{\lf,\ssi}^{\lf,\Mh}$ induces isomorphisms on homology of global sections, that is, the morphisms $\hat F{}_{\lf,\ssi}^{\lf,\Mh}:\hat H^{\lf,\ssi}_k(Y;R)\ra MH_k^\lf(Y;R)$ in Example \ref{kh4ex4} are isomorphisms. The analogue of \eq{kh4eq54} on homology now implies that the $F_{\lf,\ssi}^{\lf,\Mh}:H_k^{\lf,\ssi}(Y;R)\ra MH_k^\lf(Y;R)$ in Example \ref{kh4ex3} are also isomorphisms.

But locally finite homology theories are known to be canonically isomorphic on sufficiently nice topological spaces (such as manifolds), as in Petkova \cite{Petk} and Skljarenko \cite{Sklj1}, for instance. Thus $MH_k^\lf(-;R)$ is also canonically isomorphic to other locally finite homology theories on manifolds.

To see that the isomorphisms $H_k^\lf(Y;R)\cong MH_k^\lf(Y;R)$ preserve the data $\Pi$ in Property \ref{kh2pr2}(a), observe that the following diagram commutes
\begin{equation*}
\xymatrix@C=110pt@R=13pt{ *+[r]{H_k^\ssi(Y;R)} \ar[d]^\Pi \ar[r]_{F_\ssi^\Mh} & 
*+[l]{MH_k(Y;R)} \ar[d]_\Pi \\
*+[r]{H_k^{\lf,\ssi}(Y;R)} \ar[r]^{F{}_{\lf,\ssi}^{\lf,\Mh}} & 
*+[l]{MH_k^\lf(Y;R),\!\!}  }
\end{equation*}
since the chain-level analogue commutes, both routes mapping~$\si\mapsto[\De_k,0,0,\si]$.

To see that they preserve pushforwards $f_*$ and pullbacks $i^*$ in Property \ref{kh2pr2}(b),(c), note that $f_*\ci F_{\lf,\ssi}^{\lf,\Mh}=F_{\lf,\ssi}^{\lf,\Mh}\ci f_*:H_k^{\lf,\ssi}(Y_1;R)\ra MH_k^\lf(Y_2;R)$ as in Example \ref{kh4ex3}, and $\hat F{}_{\lf,\ssi}^{\lf,\Mh}\ci i^*=i^*\ci \hat F{}_{\lf,\ssi}^{\lf,\Mh}:\hat H_k^{\lf,\ssi}(Y;R)\ra MH_k^\lf(U;R)$ as in Example \ref{kh4ex4}. This completes the proof.
\end{proof}

Combining the proof of Theorem \ref{kh4thm7} with \eq{kh2eq39} gives a natural equivalence in $D(Y;R)$, where $\om_Y$ is the dualizing complex of~$Y$:
\e
\MC_{-\bu}^\lf(Y;R)\simeq \om_Y.
\label{kh4eq56}
\e

In Definition \ref{kh4def8} we defined a sheaf morphism $i_Y:R_Y\ra\MC^0(Y;R)$ fitting into an exact sequence \eq{kh4eq36} of sheaves on $Y$. The analogue for homology is a sheaf morphism $j_Y:O_Y\ra\MC_m^\lf(Y;R)$, where $O_Y$ is the orientation sheaf of $Y$ from Definition \ref{kh2def9}, and~$m=\dim Y$.

\begin{dfn} Let $Y$ be a manifold of dimension $m$. We will define a morphism $j_Y:O_Y\ra\MC_m^\lf(Y;R)$ of sheaves of $R$-modules on $Y$. As in Definition \ref{kh2def9}, if $U\subseteq Y$ is open and we write $U=\coprod_{i\in I}U_i$ for $U_i$, $i\in I$ the connected components of $U$, then elements $\al$ of $O_Y(U)$ may equivalently be written as formal sums $\al=\sum_{i\in J}a_i\, o_{U_i}$, where $a_i\in R$ and $o_{U_i}$ is an orientation on $U_i$ for $i\in J\subseteq I$. Define $j_Y(U):O_Y(U)\ra \MC_m^\lf(Y;R)(U)=MC_m^\lf(U;R)$ by
\e
j_Y(U):\ts\sum_{i\in J}a_i\, o_{U_i}\longmapsto \ts\sum_{i\in J}a_i\,\bigl[(U_i,o_{U_i}),0,0,\id_{U_i}\bigr],
\label{kh4eq57}
\e
as in \eq{kh4eq35}. Here $\bigl[(U_i,o_{U_i}),0,0,\id_{U_i}\bigr]$ is a generator of $MC_m^\lf(U;R)$, since $(0,\id_{U_i}):U_i\ra\R^0\t U$ is proper, and the r.h.s.\ of \eq{kh4eq57} is a locally finite sum, as in Remark \ref{kh4rem7}(b), and so makes sense in $MC_m^\lf(U;R)$. It is easy to see that $\rho_{UU'}\ci j_Y(U)=j_Y(U')\ci\rho_{UU'}:O_Y(U)\ra \MC_m^\lf(Y;R)(U')$ for all open $U'\subseteq U\subseteq Y$, so this defines a sheaf morphism~$j_Y:O_Y\ra\MC_m^\lf(Y;R)$.
\label{kh4def13}
\end{dfn}

Here is the analogue of Theorem \ref{kh4thm5}. Note that $\MC_k^\lf(Y;R)=0$ for $k>m=\dim Y$, as this holds for $MC_k(Y;R)$ and~$MC_k^\lf(Y;R)$.

\begin{thm} For each manifold\/ $Y$ of dimension $m,$ the following is an exact sequence of sheaves of\/ $R$-modules on $Y\!:$
\e
\xymatrix@C=11pt{ 0 \ar[r] & O_Y \ar[r]^(0.3){j_Y} & \MC_m^\lf(Y;R) \ar[r]^(0.46)\pd & \MC_{m-1}^\lf(Y;R) \ar[r]^(0.48)\pd & \MC_{m-2}^\lf(Y;R) \ar[r]^(0.66)\pd & \cdots. }
\label{kh4eq58}
\e
Hence $\MC_{m-\bu}^\lf(Y;R)=\bigl(\MC_{m-*}(Y;R),\pd\bigr)$ is a soft resolution of\/~$O_Y$.
\label{kh4thm8}
\end{thm}

\begin{proof} Equation \eq{kh4eq58} is exact if and only if it is exact on stalks at each $y\in Y$, and so by the last part of Corollary \ref{kh4cor1}, if and only if
\e
\xymatrix@C=13.5pt{ 0 \ar[r] & O_{Y,y} \ar[r]^(0.25){j_{Y,y}} & 
MC_m(Y,Y\sm\{y\};R) \ar[r]^(0.47){\pd} & 
MC_{m-1}(Y,Y\sm\{y\};R) \ar[r]^(0.77){\pd} & \cdots }
\label{kh4eq59}
\e
is exact. The cohomology of \eq{kh4eq59} at $MC_k(Y,Y\sm\{y\};R)$ for $k<m$ is $MH_k(Y,Y\sm\{y\};R)\cong H_k(Y,Y\sm\{y\};R)$ by Theorem \ref{kh4thm3}. Identifying $Y$ near $y$ with $\R^m$ near 0, we have $H_k(Y,Y\sm\{y\};R)\cong H_k(\R^m,\R^m\sm\{0\};R)=0$ for $k<m$ by excision, so \eq{kh4eq59} is exact at $MC_k(Y,Y\sm\{y\};R)$ for $k<m$.
 
We have $O_{Y,y}\cong O_{\R^m,0}\cong R$. The kernel of the first $\pd$ in \eq{kh4eq59} is 
\begin{equation*}
MH_m(Y,Y\sm\{y\};R)\cong H_m(Y,Y\sm\{y\};R)\cong H_m(\R^m,\R^m\sm\{0\};R)\cong R,
\end{equation*}
since $MC_{m+1}(Y,Y\sm\{y\};R)=0$. Under these isomorphisms, $j_{Y,y}:O_{Y,y}\ra\Ker\pd$ is identified with $\id:R\ra R$. Therefore \eq{kh4eq59} is exact for each $y\in Y$, so \eq{kh4eq58} is exact. The last part of the theorem follows from Definition~\ref{kh4def12}.
\end{proof}

Here is the analogue of Remark~\ref{kh4rem4}:

\begin{rem} Suppose that $Y$ is an oriented manifold of dimension $m$, not necessarily compact. Then as in Property \ref{kh2pr2}(i) we have a fundamental class $[[Y]]$ in $H_m^{\lf,\ssi}(Y;R)$. Example \ref{kh2ex12} defined this explicitly by choosing a locally finite triangulation of $Y$ into smooth $m$-simplices $\si_i:\De_m\ra Y$ for $i\in I$, and setting $\ep_i=1$ if $\si_i$ is orientation-preserving, and $\ep_i=-1$ otherwise. Then $\sum_{i\in I}\ep_i\,\si_i$ is a locally finite sum in $C_m^{\lf,\ssi}(Y;R)$, and $[[Y]]=\bigl[\sum_{i\in I}\ep_i\,\si_i\bigr]$ in~$H_m^{\lf,\ssi}(Y;R)$.

Definition \ref{kh4def12} defined $[Y]\in MC_m^\lf(Y;R)$. Using \eq{kh4eq14}, Example \ref{kh4ex3}, and a limiting argument we find that $[Y]=F_{\lf,\ssi}^{\lf,\Mh}\bigl(\sum_{i\in I}\ep_i\,\si_i\bigr)$ in $MC_m(Y;R)$, so $[[Y]]=F_{\lf,\ssi}^{\lf,\Mh}\bigl([[Y]]\bigr)$ in $MH_m(Y;R)$. Thus the fundamental class $[[Y]]$ in $MH_m^\lf(Y;R)$ is identified with the usual fundamental class $[[Y]]\in H_m^\lf(Y;R)$ by the canonical isomorphism $MH_m^\lf(Y;R)\cong H_m^\lf(Y;R)$ from Theorem \ref{kh4thm7}. Although the chain $\sum_{i\in I}\ep_i\,\si_i\in C_m^{\lf,\ssi}(Y;R)$ representing $[[Y]]\in H_k^{\lf,\ssi}(Y;R)$ involves an arbitrary choice, the fundamental cycle $[Y]\in MC_m^\lf(Y;R)$ is the {\it unique\/} cycle  representing $[[Y]]\in MH_m^\lf(Y;R)$, as~$MC_{m+1}^\lf(Y;R)=0$.

We can also give another proof that $[[Y]]\in MH_m^\lf(Y;R)$ is identified with the usual fundamental class in $[[Y]]\in H_m^\lf(Y;R)$ using sheaf cohomology and Theorem \ref{kh4thm8}, since under the isomorphisms $H_m^\lf(Y;R)\cong H^0(Y,O_Y)\cong H^0(Y,R_Y)$ from \eq{kh2eq29} and $O_Y\cong R_Y$ as $Y$ is oriented, $[[Y]]\in H_m^\lf(Y;R)$ is identified with~$1\in H^0(Y,R_Y)$.

\label{kh4rem8}
\end{rem}

\subsection{Cup products, and cross products on M-cohomology}
\label{kh45}

In \S\ref{kh26} we discussed cup and cross products $\cup,\t$ on cohomology. We now define these on $MH^*(Y;R),MH^*_\cs(Y;R)$, and prove that they are identified with the usual products $\cup,\t$ by the canonical isomorphisms with~$H^*(Y;R),H^*_\cs(Y;R)$. 

\begin{dfn} Let $Y$ be a manifold. Using the notation of \S\ref{kh42}, for all $k,l\in\Z$, define $R$-bilinear morphisms
\e
\cup:\cP MC^k(Y;R)\t \cP MC^l(Y;R)\longra \cP MC^{k+l}(Y;R)
\label{kh4eq60}
\e
on generators $[V,n,s,t]\in \cP MC^k(Y;R)$, $[V',n',s',t']\in \cP MC^l(Y;R)$ by
\e
\begin{split}
[V,n,s&,t]\cup[V',n',s',t']
=(-1)^{ln}[\ti V,\ti n,\ti s,\ti t]:=(-1)^{ln}\bigl[V\t_{t,Y,t'}V',\\
&n+n',(s_1\ci\pi_V,\ldots,s_n\ci\pi_V,s_1'\ci\pi_{V'},\ldots,s_{n'}'\ci\pi_{V'}),t\ci\pi_V\bigr].
\end{split}
\label{kh4eq61}
\e
Writing $c_t,c_{t'}$ for the coorientations on $t,t'$, the coorientation on $\ti t=t\ci\pi_V:\ti V\ra Y$ is $c_{\ti t}=c_t\ci c_{\pi_V}$, where the coorientation $c_{\pi_V}$ on $\pi_V$ is induced from $c_{t'}$, using Assumption \ref{kh3ass6}(d),(l). In Proposition \ref{kh4prop11} below we show that $\cup$ is well-defined. The sign $(-1)^{ln}$ in \eq{kh4eq61} is needed to ensure that $\cup$ takes relation Definition \ref{kh4def4}(i) in $\cP MC^k(Y;R),\cP MC^l(Y;R)$ to relation (i) in~$\cP MC^{k+l}(Y;R)$.

Given generators $[V,n,s,t]\in \cP MC^j(Y;R)$, $[V',n',s',t']\in \cP MC^k(Y;R)$ and $[V'',n'',s'',t'']\in \cP MC^l(Y;R)$, we have
\ea
&\bigl([V,n,s,t]\cup[V',n',s',t']\bigr)\cup[V'',n'',s'',t'']=(-1)^{kn}\bigl[V\t_{t,Y,t'}V',n+n',
\nonumber\\
&\quad(s_1\ci\pi_V,\ldots,s_n\ci\pi_V,s_1'\ci\pi_{V'},\ldots,s_{n'}'\ci\pi_{V'}),t\ci\pi_V\bigr]\cup[V'',n'',s'',t'']
\nonumber\\
&=(-1)^{kn}\cdot (-1)^{l(n+n')}\bigl[(V\t_{t,Y,t'}V')\t_{t'\ci\pi_{V'},Y,t''}V'',n+n'+n'',
\nonumber\\
&\quad(s_1\ci\pi_V,\ldots,s_n\ci\pi_V,s_1'\ci\pi_{V'},\ldots,s_{n'}'\ci\pi_{V'},s_1''\ci\pi_{V''},\ldots,s_{n''}''\ci\pi_{V''}),t\ci\pi_V\bigr]
\nonumber\\
&=(-1)^{(k+l)n}\cdot (-1)^{ln'}\bigl[V\t_{t,Y,t'\ci\pi_{V'}}(V'\t_{t'\ci Y,t''}V'',n+n'+n'',
\nonumber\\
&\quad(s_1\ci\pi_V,\ldots,s_n\ci\pi_V,s_1'\ci\pi_{V'},\ldots,s_{n'}'\ci\pi_{V'},s_1''\ci\pi_{V''},\ldots,s_{n''}''\ci\pi_{V''}),t\ci\pi_V\bigr]
\nonumber\\
&=(-1)^{(k+l)n}[V,n,s,t]\cup\bigl[V'\t_{t',Y,t''}V'',n'+n'',
\nonumber\\
&\quad(s_1'\ci\pi_{V'},\ldots,s_{n'}'\ci\pi_{V'},s_1''\ci\pi_{V''},\ldots,s_{n''}''\ci\pi_{V''}),t'\ci\pi_{V'}\bigr]
\nonumber\\
&=[V,n,s,t]\cup\bigl([V',n',s',t']\cup[V'',n'',s'',t'']\bigr),
\label{kh4eq62}
\ea
using \eq{kh4eq61} in the first, second, fourth and fifth steps, and natural isomorphisms of fibre products in the third. Equation \eq{kh4eq62} implies that
\e
(\al\cup\be)\cup\ga=\al\cup(\be\cup\ga)\qquad\text{in $\cP MC^{j+k+l}(Y;R)$}
\label{kh4eq63}
\e
for all $\al\in\cP MC^j(Y;R)$, $\be\in\cP MC^k(Y;R)$ and $\ga\in\cP MC^l(Y;R)$.

Applying $\d:\cP MC^{k+l}(Y;R)\ra \cP MC^{k+l+1}(Y;R)$ to \eq{kh4eq61},
we see that
\ea
&\d\bigl([V,n,s,t]\cup[V',n',s',t']\bigr)
=(-1)^{ln}\bigl[\pd(V\t_{t,Y,t'}V'),n+n',
\nonumber\\
&\;\>(s_1\ci\pi_V,\ldots,s_n\ci\pi_V,s_1'\ci\pi_{V'},\ldots,s_{n'}'\ci\pi_{V'})\ci i_{\ti V},t\ci\pi_V\ci i_{\ti V}\bigr]
\nonumber\\
&=(-1)^{ln}\bigl[(\pd V)\t_{t\ci i_V,Y,t'}V',n+n',(s_1\ci i_V\ci\pi_{\pd V},\ldots,s_n\ci i_V\ci\pi_{\pd V},
\nonumber\\
&\qquad\qquad\quad s_1'\ci\pi_{V'},\ldots,s_{n'}'\ci\pi_{V'}),t\ci i_V\ci\pi_{\pd V}\bigr]
\nonumber\\
&\;\>+(-1)^{ln}\cdot(-1)^{(m+n-k)-n}\bigl[V\t_{t,Y,t'\ci i_{V'}}(\pd V'),n+n',(s_1\ci \pi_V,\ldots,s_n\ci\pi_V,
\nonumber\\
&\qquad\qquad\qquad\qquad\qquad\qquad\quad s_1'\ci i_{V'}\ci\pi_{\pd V'},\ldots,s_{n'}'\ci i_{V'}\ci\pi_{\pd V'}),t\ci\ci\pi_V\bigr]
\nonumber\\
&=[\pd V,n,s\!\ci\! i_V,t\!\ci\! i_V]\cup[V',n',s',t']+(-1)^k[V,n,s,t]\cup[\pd V',n',s'\!\ci\! i_{V'},t'\!\ci\! i_{V'}]
\nonumber\\
&=\bigl(\d[V,n,s,t]\bigr)\cup[V',n',s',t']+(-1)^k[V,n,s,t]\cup\bigl(\d[V',n',s',t']\bigr),
\label{kh4eq64}
\ea
using \eq{kh4eq20} in the first and fourth steps,  Assumptions \ref{kh3ass5}(c), \ref{kh3ass6}(m) and \eq{kh3eq2}, \eq{kh4eq19} in the second, and \eq{kh4eq61} in the third. As \eq{kh4eq64} holds for all generators $[V,n,s,t],[V',n',s',t']$ of $\cP MC^k(Y;R),\cP MC^l(Y;R)$, we see that
\e
\d(\al\cup\be)=(\d\al)\cup\be+(-1)^k\al\cup(\d\be)
\label{kh4eq65}
\e
for all $\al\in\cP MC^k(Y;R)$ and $\be\in\cP MC^l(Y;R)$.

As in Definition \ref{kh4def5}, the identity is $\Id_Y=[Y,0,0,\id_Y]\in\cP MC^0(Y;R)$. Given a generator $[V,n,s,t]\in \cP MC^k(Y;R)$, from \eq{kh4eq61} we have
\e
\begin{split}
\Id_Y\cup[V,n,s,t]&=\bigl[Y\t_{\id_Y,Y,t}V,0+n,(s_1\ci\pi_V,\ldots,s_n\ci\pi_V),t\ci\pi_V\bigr]\\
&=[V,n,s,t],
\end{split}
\label{kh4eq66}
\e
by natural isomorphisms of fibre products. Similarly $[V,n,s,t]\cup\Id_Y=[V,n,s,t]$. Therefore for all $\al\in\cP MC^k(Y;R)$ we have
\e
\Id_Y\cup \al=\al\cup\Id_Y=\al.
\label{kh4eq67}
\e

Suppose $f:Y_1\ra Y_2$ is a smooth map of manifolds, so that Definition \ref{kh4def5} defines the pullback $f^*:\cP MC^k(Y_2;R)\ra\cP MC^k(Y_1;R)$. If $[V,n,s,t]\in \cP MC^k(Y_2;R)$ and $[V',n',s',t']\in \cP MC^l(Y_2;R)$ we have
\ea
&f^*\bigl([V,n,s,t]\cup[V',n',s',t']\bigr)=f^*\bigl(
(-1)^{ln}\bigl[V\t_{t,Y,t'}V',n+n',
\nonumber\\
&\qquad\qquad (s_1\ci\pi_V,\ldots,s_n\ci\pi_V,s_1'\ci\pi_{V'},\ldots,s_{n'}'\ci\pi_{V'}),t\ci\pi_V\bigr]\bigr)
\nonumber\\
&=(-1)^{ln}\bigl[(V\t_{t,Y,t'}V')\t_{t\ci\pi_V,Y_2,f}Y_1,n+n',\nonumber\\
&\qquad\qquad (s_1\ci\pi_V,\ldots,s_n\ci\pi_V,
s_1'\ci\pi_{V'},\ldots,s_{n'}'\ci\pi_{V'}),t\ci\pi_V\bigr]\bigr)
\nonumber\\
&=(-1)^{ln}\bigl[(V\t_{t,Y_2,f}Y_1)\t_{\pi_{Y_1},Y_1,\pi_{Y_1}}(V'\t_{t',Y_2,f}Y_1),n+n',\nonumber\\
&\qquad\qquad (s_1\ci\pi_V,\ldots,s_n\ci\pi_V,
s_1'\ci\pi_{V'},\ldots,s_{n'}'\ci\pi_{V'}),t\ci\pi_V\bigr]\bigr)
\nonumber\\
&=\bigl[V\t_{t,Y_2,f}Y_1,n,s\ci\pi_V,\pi_{Y_1}\bigr]\cup
\bigl[V'\t_{t',Y_2,f}Y_1,n',s'\ci\pi_V,\pi_{Y_1}\bigr]
\nonumber\\
&=f^*\bigl([V,n,s,t]\bigr)\cup f^*\bigl([V',n',s',t']\bigr),
\label{kh4eq68}
\ea
using \eq{kh4eq61} in the first and fourth steps, \eq{kh4eq21} in the second and fifth, and natural isomorphisms of fibre products in the third. Equation \eq{kh4eq68} implies that for all $\al\in\cP MC^k(Y_2;R)$ and $\be\in\cP MC^l(Y_2;R)$ we have
\e
f^*(\al\cup\be)=f^*(\al)\cup f^*(\be)\qquad\text{in $\cP MC^{k+l}(Y_1;R)$.}
\label{kh4eq69}
\e

\label{kh4def14}
\end{dfn}

The next result will be proved in~\S\ref{kh79}.

\begin{prop} The product\/ $\cup$ in \eq{kh4eq60}--\eq{kh4eq61} is well defined.
\label{kh4prop11}
\end{prop}

\begin{dfn} Let $Y$ be a manifold, and $k,l\in\Z$. Define 
\e
\begin{split}
\cup_{k,l}&:\cP\MC^k(Y;R)\ot_R\cP\MC^l(Y;R)\longra\cP\MC^{k+l}(Y;R)\quad\text{by}\\
\cup_{k,l}(U)=\cup&:\cP MC^k(U;R)\ot_R\cP MC^l(U;R)\longra\cP MC^{k+l}(U;R)
\end{split}
\label{kh4eq70}
\e
for all open $U\subseteq Y$, where we interpret $\cup$ in \eq{kh4eq60} with $U$ in place of $Y$ as an $R$-linear map on $\cP MC^k(U;R)\ot_R\cP MC^l(U;R)$ rather than an $R$-bilinear map on $\cP MC^k(U;R)\t\cP MC^l(U;R)$. As $\cup$ in Definition \ref{kh4def14} commutes with pullbacks $i^*:\cP MC^*(U;R)\ra\cP MC^*(U';R)$ for open $U'\subseteq U\subseteq Y$ with inclusion $i:U'\hookra U$, this defines a morphism $\cup_{k,l}$ of presheaves of $R$-modules on $Y$. Thus passing to sheafifications induces a morphism
\e
\cup_{k,l}:\MC^k(Y;R)\ot_R\MC^l(Y;R)\longra\MC^{k+l}(Y;R).
\label{kh4eq71}
\e

Then taking global sections $\cup=\cup_{k,l}(Y)$ defines $R$-bilinear {\it cup products\/}
\e
\cup:MC^k(Y;R)\t MC^l(Y;R)\longra MC^{k+l}(Y;R).
\label{kh4eq72}
\e
By construction these satisfy $\Pi(\al\cup\be)=\Pi(\al)\cup\Pi(\be)$ in $MC^{k+l}(Y;R)$ for all $\al\in\cP MC^k(Y;R)$ and $\be\in\cP MC^l(Y;R)$, and $i^*(\ga\cup\de)=i^*(\ga)\cup i^*(\de)$ in $MC^{k+l}(U;R)$ whenever $\ga\in MC^k(Y;R)$, $\de\in MC^l(Y;R)$ and $U\subseteq Y$ is open with inclusion $i:U\hookra Y$, and the morphisms $\cup$ in \eq{kh4eq72} are determined uniquely by these properties.

Applying $\Pi$, or by properties of sheafification, we see that equations \eq{kh4eq61}, \eq{kh4eq63}, \eq{kh4eq65}, \eq{kh4eq67} and \eq{kh4eq69} hold in $MC^*(-;R)$ as well as in~$\cP MC^*(-;R)$.

As in \S\ref{kh42}, $MC^k(Y;R)$ is the global sections of a sheaf $\MC(Y;R)$ on $Y$, and so each $\al\in MC^k(Y;R)$ has a support $\supp\al$, a closed subset of $Y$. Section \ref{kh43} defined the compactly-supported M-cochains $MC^k_\cs(Y;R)$ to be the $R$-submodule of $\al\in MC^k(Y;R)$ with $\supp\al$ compact. Since cup products are compatible with restriction to open subsets, we see that if $\al\in MC^k(Y;R)$ and $\be\in MC^l(Y;R)$ then $\supp(\al\cup\be)\subseteq(\supp\al)\cap(\supp\be)$, and so in particular, if either $\supp\al$ or $\supp\be$ is compact, then $\supp(\al\cup\be)$ is compact. Therefore $\cup$ in \eq{kh4eq72} restricts to $R$-bilinear morphisms
\e
\begin{split}
&\cup:MC^k_\cs(Y;R)\t MC^l(Y;R)\longra MC^{k+l}_\cs(Y;R),\\
&\cup:MC^k(Y;R)\t MC^l_\cs(Y;R)\longra MC^{k+l}_\cs(Y;R),\\
&\cup:MC^k_\cs(Y;R)\t MC^l_\cs(Y;R)\longra MC^{k+l}_\cs(Y;R).
\end{split}
\label{kh4eq73}
\e

Equation \eq{kh4eq65} in $MC^*(Y;R)$ implies that $\cup$ descends to M-cohomology. Thus as in \eq{kh2eq45}, from \eq{kh4eq72} and \eq{kh4eq73} we define $R$-bilinear morphisms
\e
\begin{split}
&\cup:MH^k(Y;R)\t MH^l(Y;R)\longra MH^{k+l}(Y;R),\\
&\cup:MH^k_\cs(Y;R)\t MH^l(Y;R)\longra MH^{k+l}_\cs(Y;R),\\
&\cup:MH^k(Y;R)\t MH^l_\cs(Y;R)\longra MH^{k+l}_\cs(Y;R),\\
&\cup:MH^k_\cs(Y;R)\t MH^l_\cs(Y;R)\longra MH^{k+l}_\cs(Y;R),
\end{split}
\label{kh4eq74}
\e
by $[\al]\cup[\be]=[\al\cup\be]$ for $\al\in MC^k_?(Y;R)$, $\be\in MC^l_?(Y;R)$ with $\d\al=\d\be=0$.
\label{kh4def15}
\end{dfn}

\begin{thm} Under the canonical isomorphisms $MH^k(Y;R)\cong H^k(Y;R),$
$MH^k_\cs(Y;R)\cong H^k_\cs(Y;R)$ from Theorems\/ {\rm\ref{kh4thm4}} and\/ {\rm\ref{kh4thm6},} the cup products in \eq{kh4eq74} are identified with the usual cup products \eq{kh2eq45} on ordinary (compactly-supported) cohomology $H^*(Y;R),H^*_\cs(Y;R)$.
\label{kh4thm9}
\end{thm}

\begin{proof} In \S\ref{kh262} we explained how to define the cup product on ordinary cohomology using sheaf cohomology, and a soft resolution $\cF^\bu$ of the constant sheaf $R_Y$ on $Y$, as in equations \eq{kh2eq63}--\eq{kh2eq65}. We will apply this method to the soft resolution $\MC^\bu(Y;R)$ of $R_Y$ given by Theorem~\ref{kh4thm5}.

Define $\cup_{k,l}$ as in \eq{kh4eq71}. As the morphism $i_Y:R_Y\ra\MC^0(Y;R)$ in Definition \ref{kh4def8} maps $1\mapsto\Id_Y$, and \eq{kh4eq67} gives $\Id_Y\cup\Id_Y=\Id_Y$, we have
\begin{equation*}
i_Y\ci I_\cup=\cup_{0,0}\ci (i_Y\ot i_Y):R_Y\ot_R R_Y\longra \MC^0(Y;R),
\end{equation*}
giving the first equation of \eq{kh2eq64}, where $I_\cup:R_Y\ot_R R_Y\ra R_Y$ is as in~\eq{kh2eq60}. Equation \eq{kh4eq65} implies that
\begin{align*}
\d\ci\cup_{k,l}&=\cup_{k+1,l}\ci(\d\ot\id_{\MC^l(Y;R)})+(-1)^k\cup_{k,l+1}\ci(\id_{\MC^k(Y;R)}\ot\d):\\
&\qquad\MC^k(Y;R)\ot_R\MC^l(Y;R)\longra \MC^{k+l+1}(Y;R),
\end{align*}
giving the second equation of \eq{kh2eq64}. Hence \eq{kh2eq65} gives an expression for the cup product $\cup$ under the isomorphism $H^k(Y;R)\cong H^k\bigl(\MC^*(Y;R)(Y),\d\bigr)$. This coincides with the definition $[\al]\cup[\be]=[\al\cup\be]$ of $\cup$ on $MH^*(Y;R)$ in Definition \ref{kh4def15}, under the isomorphism~$H^k(Y;R)\cong MH^k(Y;R)$. 

So the first line of \eq{kh4eq74} is identified with the usual cup product on $H^*(Y;R)$. Using the same argument, but taking compactly-supported sheaf cohomology in two or three of the factors, shows that the last three lines of \eq{kh4eq74} are also identified with the usual cup products on~$H^*(Y;R),H^*_\cs(Y;R)$.
\end{proof}

An alternative proof of Theorem \ref{kh4thm9} for $\cup$ on $MH^*(Y;R)$ with $R=\Z,\Z_n$ or $R$ a $\Q$-algebra would be to verify axioms (i)--(iii) for cup products on $H^*(Y;R)$ in Proposition \ref{kh2prop2}, where (i),(ii) already follow from \eq{kh4eq63}, \eq{kh4eq67} and~\eq{kh4eq69}.

Observe that $\bigl(MC^*(Y;R),\d,\cup,\Id_Y\bigr)$ is a differential graded algebra (dga) over $R$, with graded product $\cup$ which is associative by \eq{kh4eq63}, and $\Id_Y$ which is a strict identity by \eq{kh4eq67}. Now in \S\ref{kh27} we explained that in homotopy theory one often associates a dga (or cdga) $\bigl(C^*(Y;R),\d,\cup,1_Y\bigr)$ over $R$ to a topological space $Y$. This dga has cohomology $H^*(Y;R)$, and is unique up to equivalence in an $\iy$-category $\dga_R^\iy$ of dgas. From the dga we can compute invariants of $Y$ such as Steenrod squares \cite[\S VI.15]{Bred1}, \cite[\S 4.L]{Hatc} and Massey products~\cite{Mass1}.

\begin{thm} For each manifold\/ $Y,$ the dga $\bigl(MC^*(Y;R),\d,\cup,\Id_Y\bigr)$ over $R$ is equivalent in $\dga_R^\iy$ to the `usual' dga over $R$ associated to $Y$ in topology, as represented for instance by the singular cochains $\bigl(C_\rsi^*(Y;R),\d,\cup,1_Y\bigr)$ with the Alexander--Whitney cup product\/~$\cup$.

Therefore topological invariants of\/ $Y$ depending on the dga up to equivalence, such as Steenrod squares and Massey products, may be computed using the dga $\bigl(MC^*(Y;R),\d,\cup,\Id_Y\bigr),$ and will give the correct answers under the canonical isomorphism $MH^*(Y;R)\cong H^*(Y;R)$ from Theorem\/~{\rm\ref{kh4thm4}}.
\label{kh4thm10}
\end{thm}

\begin{proof} We noted in \S\ref{kh27} that to prove a cochain dga $\bigl(C^*(Y;R),\d,\cup,1_Y\bigr)$ for $Y$ is equivalent to the `usual' dga, it is sufficient that $\bigl(C^*(Y;R),\d\bigr)$ should be the global sections of a soft resolution $\cF^\bu$ of $R_Y$ as in \eq{kh2eq63}, with identity $1_Y=i(1)\in\cF^0(Y)$, and cup product $\cup$ defined using sheaf morphisms $\psi_{k,l}:\cF^k\ot_R\cF^l\ra\cF^{k+l}$ satisfying \eq{kh2eq64} and associativity \eq{kh2eq68}. We showed all this in the proof of Theorem \ref{kh4thm9}, except for \eq{kh2eq68}, which follows from~\eq{kh4eq63}.
\end{proof}

\begin{rem} The cup products $\cup$ on $H^*(Y;R),H^*_\cs(Y;R)$, and hence on $MH^*(Y;R),MH^*_\cs(Y;R)$ by Theorem \ref{kh4thm9}, are associative and supercommutative with strict identity $[1_Y]$. At the cochain level, the cup products $\cup$ on $MC^*(Y;R),MC^*_\cs(Y;R)$ are associative by \eq{kh4eq63}, with strict identity $\Id_Y$ in $MC^0(Y;R)$ by \eq{kh4eq67}. However, they are generally {\it not\/} supercommutative, that is, we can have $\al\cup\be\ne (-1)^{kl}\be\cup\al$ for $\al\in MC^k(Y;R)$ and $\be\in MC^l(Y;R)$.

This is necessary: as in Remark \ref{kh2rem4}(ii) and \S\ref{kh27}, for some rings $R$ such as $\Z$ or $\Z_2$, it is not possible to define a cohomology theory $\bigl(C^*(Y;R),\d\bigr)$ computing $H^*(Y;R)$ with a supercommutative cup product $\cup$ defined on cochains, because Steenrod squares \cite[\S VI.15--\S VI.16]{Bred1} are an obstruction to this. If $R$ is a $\Q$-algebra there is no obstruction, and in \S\ref{kh51} we will define rational M-cohomology $MH^*_\Q(Y;R)$, for which $\cup$ is supercommutative on cochains~$MC^*_\Q(Y;R)$.

The reason $\cup$ is not supercommutative on $MC^*(Y;R)$ is in the defining equation \eq{kh4eq61} for $[V,n,s,t]\!\cup\![V',n',s',t']$, we set $\ti s:V\t_{t,Y,t'}V'\!\ra\!\R^{n+n'}$ to be
\begin{equation*}
\ti s=(s_1\ci\pi_V,\ldots,s_n\ci\pi_V,s_1'\ci\pi_{V'},\ldots,s_{n'}'\ci\pi_{V'}).
\end{equation*}
But $[V',n',s',t']\cup[V,n,s,t]$ would involve $\check s:V'\t_{t',Y,t}V\ra\R^{n+n'}$, where
\begin{equation*}
\check s=(s_1'\ci\pi_{V'},\ldots,s_{n'}'\ci\pi_{V'},s_1\ci\pi_V,\ldots,s_n\ci\pi_V).
\end{equation*}
So even after identifying $V\t_{t,Y,t'}V'\cong V'\t_{t',Y,t}V$, $[V,n,s,t]\cup[V',n',s',t']$ and $[V',n',s',t']\cup[V,n,s,t]$ still differ by a permutation of the $n+n'$ coordinates in the target $\R^{n+n'}$ of~$\ti s:V\t_{t,Y,t'}V'\ra\R^{n+n'}$.
\label{kh4rem9}
\end{rem}

Using this idea of permuting the $n$ coordinates in the target $\R^n$ of $s:V\ra\R^n$ in generators $[V,n,s,t]$, we can define a cochain-level involution $\Xi$ of the entire cohomology theory $MH^*(-;R)$, with $\Xi(\al\cup\be)=(-1)^{kl}\Xi(\be)\cup\Xi(\al)$ for all $\al\in MC^k(Y;R)$ and~$\be\in MC^l(Y;R)$.

\begin{dfn} Let $Y$ be a manifold. For all $k\in\Z$, define $R$-linear maps
\e
\begin{split}
&\Xi:\cP MC^k(Y;R)\longra\cP MC^k(Y;R)\qquad\text{by}\\
&\Xi:\bigl[V,n,(s_1,s_2,\ldots,s_n),t\bigr]\!\longmapsto \!(-1)^{n(n-1)/2}\bigl[V,n,(s_n,s_{n-1},\ldots,s_1),t\bigr]
\end{split}
\label{kh4eq75}
\e
for all generators $[V,n,s,t]$ of $\cP MC^k(Y;R)$. This takes relations Definition \ref{kh4def4}(i),(ii) to themselves, where the sign $(-1)^{n(n-1)/2}$ in \eq{kh4eq75} ensures compatibility with the sign $(-1)^{n-i}$ in \eq{kh4eq16}, and so $\Xi$ is well-defined. Clearly $\Xi^2=\id$, so $\Xi:\cP MC^k(Y;R)\ra\cP MC^k(Y;R)$ is an isomorphism.

Then $\Xi$ commutes with differentials $\d:\cP MC^k(Y;R)\ra\cP MC^{k+1}(Y;R)$ and pullbacks $f^*:\cP MC^k(Y_2;R)\ra\cP MC^k(Y_1;R)$ in Definition \ref{kh4def5}, with $\Xi(\Id_Y)=\Id_Y$, so $\Xi$ descends to $\Xi:MC^k(Y;R)\ra MC^k(Y;R)$ in Definition \ref{kh4def6}. In fact, the involutions $\Xi$ apply to and commute with the whole of the theory of $MH^*(-;R),MH^*_\cs(-;R)$ in \S\ref{kh42}--\S\ref{kh43}. Thus they induce morphisms $\Xi_*:MH^k(Y;R)\ra MH^k(Y;R)$ for all $Y,k$, which as they commute with all the structures on cohomology and act as the identity on $MH^0(*;R)\cong R$, must be the identity maps on $MH^k(Y;R)$ by the Eilenberg--Steenrod Theorem~\ref{kh2thm2}.

Comparing \eq{kh4eq61} and \eq{kh4eq75}, we see that if $[V,n,s,t]\in \cP MC^k(Y;R)$ and $[V',n',s',t']\in \cP MC^l(Y;R)$ are generators then
\begin{equation*}
\Xi\bigl([V,n,s,t]\cup[V',n',s',t']\bigr)=(-1)^{kl}\Xi\bigl([V',n',s',t']\bigr)\cup\Xi\bigl([V,n,s,t]\bigr).
\end{equation*}
Therefore for all $\al\in\cP MC^k(Y;R)$ and $\be\in\cP MC^l(Y;R)$ we have
\e
\Xi(\al\cup\be)=(-1)^{kl}\Xi(\be)\cup\Xi(\al).
\label{kh4eq76}
\e
Sheafifying, we see that \eq{kh4eq76} also holds in~$MC^*(Y;R)$. 
 
This involution $\Xi$ can be understood as some kind of {\it homotopy commutative structure\/} on the dga $\bigl(MC^*(Y;R),\d,\cup,\Id_Y\bigr)$: it induces the identity on cohomology, and therefore \eq{kh4eq76} is a cochain-level identity which implies that $\cup$ is supercommutative on cohomology~$MH^*(Y;R)$.

We can also define an involution $\Xi$ at the chain level on (locally finite) M-homology $MH_*(-;R),MH_*^\lf(-;R)$ in \S\ref{kh41} and \S\ref{kh44} in exactly the same way.

\label{kh4def16}
\end{dfn}

Next we define cross products on M-cohomology:

\begin{dfn} Let $Y_1,Y_2$ be manifolds. Define $R$-bilinear maps
\e
\begin{gathered}
\t:MC^k(Y_1;R)\t MC^l(Y_2;R)\longra MC^{k+l}(Y_1\t Y_2;R)\\
\text{by}\qquad \al\t\be=\pi_{Y_1}^*(\al)\cup\pi_{Y_2}^*(\be),
\end{gathered}
\label{kh4eq77}
\e
where $\cup$ on $MC^*(Y_1\t Y_2;R)$ is as in Definition \ref{kh4def17}, and $\pi_{Y_i}:Y_1\t Y_2\ra Y_i$ for $i=1,2$ are the projections. Applied to generators $[V,n,s,t]\in MC^k(Y_1;R)$ and $[V',n',s',t']\in MC^l(Y_2;R)$, using \eq{kh4eq61} to define the cup product and \eq{kh4eq21} to give $\pi_{Y_1}^*([V,n,s,t]),\pi_{Y_2}^*([V',n',s',t'])$, and the natural isomorphism
\begin{equation*}
\bigl((V\t_{t,Y_1,\pi_{Y_1}}(Y_1\t Y_2)\bigr)\t_{\pi_{Y_1\t Y_2},Y_1\t Y_2,\pi_{Y_1\t Y_2}}
\bigl((V'\t_{t',Y_2,\pi_{Y_2}}(Y_1\t Y_2)\bigr)\cong V\t V',
\end{equation*}
we see that
\e
\begin{split}
[V,n,s&,t]\t[V',n',s',t']=(-1)^{ln}\bigl[V\t V',n+n',\\
&(s_1\ci\pi_V,\ldots,s_n\ci\pi_V,s_1'\ci\pi_{V'},\ldots,s_{n'}'\ci\pi_{V'}),t\t t'\bigr].
\end{split}
\label{kh4eq78}
\e
Thus, an alternative way to define $\t$ would be to first define $\t$ on $\cP MC^*(-;R)$ using \eq{kh4eq78}, as for $\cup$ in \eq{kh4eq61}, and then follow the method of Definition \ref{kh4def14}, Proposition \ref{kh4prop11}, and Definition~\ref{kh4def15}.

Clearly $\supp(\al\!\t\!\be)\!\subseteq\!(\supp\al)\!\t\!(\supp\be)$, so if $\al,\be$ are compactly-supported, then so is $\al\t\be$. Therefore $\t$ in \eq{kh4eq77} restricts to $R$-bilinear maps
\begin{equation*}
\t:MC^k_\cs(Y_1;R)\t MC^l_\cs(Y_2;R)\longra MC^{k+l}_\cs(Y_1\t Y_2;R).
\end{equation*}

From $\Id_Y=[Y,0,0,\id_Y]$ and equation \eq{kh4eq78} we see that
\begin{equation*}
\Id_{Y_1}\t\Id_{Y_2}=\Id_{Y_1\t Y_2}.
\end{equation*}

If $g_1:Y_1\ra Z_1$ and $g_2:Y_2\ra Z_2$ are smooth maps of manifolds and $\al\in MC^k(Z_1;R)$, $\be\in MC^l(Z_2;R)$ then
\ea
(g_1&\t g_2)^*(\al\t\be)=(g_1\t g_2)^*\bigl(\pi_{Z_1}^*(\al)\cup\pi_{Z_2}^*(\be)\bigr)
\nonumber\\
&=\bigl((g_1\t g_2)^*\ci\pi_{Z_1}^*(\al)\bigr)\cup \bigl((g_1\t g_2)^*\ci\pi_{Z_2}^*(\be)\bigr)
\nonumber\\
&=(\pi_{Z_1}\ci(g_1\t g_2))^*(\al)\cup(\pi_{Z_2}\ci(g_1\t g_2))^*(\be)
\nonumber\\
&=(g_1\ci\pi_{Y_1})^*(\al)\cup(g_2\ci\pi_{Y_2})^*(\be)=\bigl(\pi_{Y_1}^*\ci g_1^*(\al)\bigr)\cup\bigl(\pi_{Y_2}^*\ci g_2^*(\be)\bigr)
\nonumber\\
&=g_1^*(\al)\t g_2^*(\be),
\label{kh4eq79}
\ea
using \eq{kh4eq77} in the first and sixth steps, \eq{kh4eq69} in the second, and functoriality of pullbacks in the third and fifth. So cross products are compatible with pullbacks.

Equation \eq{kh4eq65} and the compatibility of $\d$ with pullbacks $\pi_{Y_i}^*$ yields
\begin{equation*}
\d(\al\t\be)=(\d\al)\t\be+(-1)^k\al\t(\d\be)
\end{equation*}
for all $\al\in MC^k(Y_1;R)$ and $\be\in MC^l(Y_2;R)$. Therefore cross products descend to cohomology $MH^*(-;R),MH^*_\cs(-;R)$ by $[\al]\t[\be]=[\al\t\be]$, giving products
\e
\begin{split}
\t&:MH^k(Y_1;R)\t MH^l(Y_2;R)\longra MH^{k+l}(Y_1\t Y_2;R),\\
\t&:MH^k_\cs(Y_1;R)\t MH^l_\cs(Y_2;R)\longra MH^{k+l}_\cs(Y_1\t Y_2;R).
\end{split}
\label{kh4eq80}
\e

\label{kh4def17}
\end{dfn}

From equations \eq{kh2eq54} and \eq{kh4eq77} and Theorem \ref{kh4thm9} we deduce:

\begin{cor} Under the canonical isomorphisms $MH^k(Y;R)\cong H^k(Y;R),$
$MH^k_\cs(Y;R)\cong H^k_\cs(Y;R)$ from Theorems\/ {\rm\ref{kh4thm4}} and\/ {\rm\ref{kh4thm6},} the cross products in \eq{kh4eq80} are identified with the usual cross products \eq{kh2eq52} on ordinary (compactly-supported) cohomology $H^*(Y;R),H^*_\cs(Y;R)$.
\label{kh4cor2}
\end{cor}

\subsection{Cap products, and cross products on M-homology}
\label{kh46}

The story for cap products, and for cross products on M-homology, is similar to that for $\cup,\t$ on M-cohomology in \S\ref{kh45}, so we will be brief in places. 

\begin{dfn} Let $Y$ be a manifold, of dimension $m$. For all $k,l\in\Z$, define $R$-bilinear morphisms
\e
\cap:\cP MC^k(Y;R)\t \cP MC_l^\lf(Y;R)\longra \cP MC_{l-k}^\lf(Y;R)
\label{kh4eq81}
\e
on generators $[V,n,s,t]\in \cP MC^k(Y;R)$, $[V',n',s',t']\in \cP MC_l^\lf(Y;R)$ by
\ea
[V,n,s&,t]\cap[V',n',s',t']
=(-1)^{(l+m)n}[\ti V,\ti n,\ti s,\ti t]:=(-1)^{(l+m)n}\bigl[V\t_{t,Y,t'}V',
\nonumber\\
&n+n',(s_1\ci\pi_V,\ldots,s_n\ci\pi_V,s_1'\ci\pi_{V'},\ldots,s_{n'}'\ci\pi_{V'}),t\ci\pi_V\bigr].
\label{kh4eq82}
\ea

Here \eq{kh4eq82} is apparently the same as \eq{kh4eq61}, except for the sign. In fact in \eq{kh4eq61} we have $\dim V'=m+n'-l$, $\dim Y=m$ by Definition \ref{kh4def4}, but in \eq{kh4eq82} we have $\dim V'=n'+l$, $\dim Y=m$ by Definition \ref{kh4def11}, so in both cases the sign is $(-1)^{n(n'+\dim V'+\dim Y)}$, and from this point of view the signs are the same. The signs are needed to ensure \eq{kh4eq61}, \eq{kh4eq82} are compatible with applying the relations Definitions \ref{kh4def4}(i), \ref{kh4def11}(i) to $[V,n,s,t]$ and $[\ti V,\ti n,\ti s,\ti t]$, where the proof of compatibility involves the equation, with our orientation conventions
\begin{equation*}
(V\t\R)\t_{t\ci\pi_V,Y,t'}V'\cong (-1)^{\dim V'+\dim Y}(V\t_{t,Y,t'}V')\t\R,
\end{equation*}
in manifolds with corners either oriented, or cooriented over~$Y$.

Equations \eq{kh4eq61} and \eq{kh4eq82} also differ in the orientations/coorientations, and submersions: in \eq{kh4eq61}, $t:V\ra Y$ and $t':V'\ra Y$ are cooriented submersions and these ensure $\ti V=V\t_{t,Y,t'}V'$ exists with $\ti t:\ti V\ra Y$ a cooriented submersion, but in \eq{kh4eq82}, $t:V\ra Y$ is a cooriented submersion and $V'$ is oriented, which ensure that $\ti V=V\t_{t,Y,t'}V'$ exists and is oriented. 

By a very similar proof to that of Proposition \ref{kh4prop11} in \S\ref{kh79}, modifying signs, orientations/coorientations and submersions as above, we can show $\cap$ in \eq{kh4eq81}--\eq{kh4eq82} is well defined.

As for \eq{kh4eq62}--\eq{kh4eq63}, we show that if $\al\in\cP MC^j(Y;R)$, $\be\in\cP MC^k(Y;R)$ and $\ga\in \cP MC^\lf_l(Y;R)$ then
\e
(\al\cup\be)\cap\ga=\al\cap(\be\cap\ga)\qquad\text{in $\cP MC^\lf_{l-j-k}(Y;R)$.}
\label{kh4eq83}
\e

As for \eq{kh4eq64}--\eq{kh4eq65}, if $\al\in\cP MC^k(Y;R)$ and $\be\in \cP MC^\lf_l(Y;R)$ then
\e
\pd(\al\cap\be)=(\d\al)\cap\be+(-1)^k\al\cap(\pd\be)\qquad\text{in $\cP MC^\lf_{l-k-1}(Y;R)$.}
\label{kh4eq84}
\e

As for \eq{kh4eq67}, if $\al\in \cP MC^\lf_k(Y;R)$ then
\e
\Id_Y\cap \al=\al.
\label{kh4eq85}
\e

In a similar way to \eq{kh4eq68}--\eq{kh4eq69}, if $f:Y_1\ra Y_2$ is a proper smooth map of manifolds, so that Definition \ref{kh4def5} gives $f^*:\cP MC^k(Y_2;R)\ra\cP MC^k(Y_1;R)$ and (using $f$ proper) Definition \ref{kh4def11} gives $f_*:\cP MC^\lf_l(Y_1;R)\ra\cP MC^\lf_l(Y_2;R)$, we show that if $\al\in\cP MC^k(Y_2;R)$ and $\be\in\cP MC^\lf_l(Y_1;R)$ we have
\e
\al\cap f_*(\be)=f_*\bigl(f^*(\al)\cap\be\bigr)\qquad\text{in $\cP MC^\lf_{k-l}(Y_2;R)$.}
\label{kh4eq86}
\e

By Proposition \ref{kh4prop10}, the morphisms $\Pi:MC_*(Y_a;R)\ra \cP MC^\lf_*(Y_a;R)$ embed $MC_*(Y_a;R)$ as the $R$-submodules of compactly-supported sections in $\cP MC^\lf_*(Y_a;R)$ for $a=1,2$. So we can restrict $\cap$ in \eq{kh4eq81} to
\begin{equation*}
\cap:\cP MC^k(Y_a;R)\t MC_l(Y_a;R)\longra MC_{l-k}(Y_a;R).
\end{equation*}
Since $f_*:MC_l(Y_1;R)\ra MC_l(Y_2;R)$ in \S\ref{kh41} is defined without supposing $f$ proper, for $\be\in MC_l(Y;R)$ equation \eq{kh4eq86} holds without assuming $f$ proper.
\label{kh4def18}
\end{dfn}

Here is the analogue of Definition \ref{kh4def15}.

\begin{dfn} Let $Y$ be a manifold, and $k,l\in\Z$. As in \eq{kh4eq70}, define 
\begin{align*}
\cap_{k,l}&:\cP\MC^k(Y;R)\ot_R\cP\MC^\lf_l(Y;R)\longra\cP\MC^\lf_{l-k}(Y;R)\quad\text{by}\\
\cap_{k,l}(U)=\cap&:\cP MC^k(U;R)\ot_R\cP MC^\lf_l(U;R)\longra\cP MC^\lf_{l-k}(U;R)
\end{align*}
for all open $U\subseteq Y$. This defines a morphism $\cap_{k,l}$ of presheaves of $R$-modules on $Y$. Passing to sheafifications induces a morphism
\e
\cap_{k,l}:\MC^k(Y;R)\ot_R\MC^\lf_l(Y;R)\longra\MC^\lf_{l-k}(Y;R).
\label{kh4eq87}
\e

Then taking global sections $\cap=\cap_{k,l}(Y)$ defines $R$-bilinear {\it cap products\/}
\e
\cap:MC^k(Y;R)\t MC_l^\lf(Y;R)\longra MC_{l-k}^\lf(Y;R).
\label{kh4eq88}
\e
These satisfy $\Pi(\al\cap\be)=\Pi(\al)\cap\Pi(\be)$ in $MC^\lf_{l-k}(Y;R)$ for all $\al\in\cP MC^k(Y;R)$ and $\be\in\cP MC^\lf_l(Y;R)$, and $i^*(\ga\cap\de)=i^*(\ga)\cap i^*(\de)$ in $MC^\lf_{l-k}(U;R)$ whenever $\ga\in MC^k(Y;R)$, $\de\in MC^\lf_l(Y;R)$ and $U\subseteq Y$ is open with inclusion $i:U\hookra Y$, and the morphisms $\cap$ in \eq{kh4eq88} are determined uniquely by these properties.

Applying $\Pi$, or by sheafification, we can show that \eq{kh4eq82}--\eq{kh4eq86} hold with $MC^*(-;R),MC^\lf_*(-;R)$ in place of $\cP MC^*(-;R),\cP MC^\lf_*(-;R)$.

Since $MC^k_\cs(Y;R)\subseteq MC^k(Y;R)$ and $MC_l(Y;R)\subseteq MC_l^\lf(Y;R)$ are the $R$-submodules of compactly-supported sections, as for \eq{kh4eq73} we may restrict \eq{kh4eq88} to $R$-bilinear morphisms
\e
\begin{split}
&\cap:MC^k(Y;R)\t MC_l(Y;R)\longra MC_{l-k}(Y;R),\\
&\cap:MC^k_\cs(Y;R)\t MC_l(Y;R)\longra MC_{l-k}(Y;R),\\
&\cap:MC^k_\cs(Y;R)\t MC_l^\lf(Y;R)\longra MC_{l-k}(Y;R).
\end{split}
\label{kh4eq89}
\e
Equations \eq{kh4eq82}--\eq{kh4eq86} hold for these by restriction, where for \eq{kh4eq86} the map $f:Y_1\ra Y_2$ need not be proper for the first line of~\eq{kh4eq89}.

Equation \eq{kh4eq84} implies that $\cap$ descends to (co)homology. So as in \eq{kh2eq46}, from \eq{kh4eq88}--\eq{kh4eq89} we obtain $R$-bilinear cap products
\e
\begin{split}
\cap&:MH^k(Y;R)\t MH_l(Y;R)\longra MH_{l-k}(Y;R),\\
\cap&:MH^k_\cs(Y;R)\t MH_l(Y;R)\longra MH_{l-k}(Y;R),\\
\cap&:MH^k(Y;R)\t MH_l^\lf(Y;R)\longra MH_{l-k}^\lf(Y;R),\\
\cap&:MH^k_\cs(Y;R)\t MH_l^\lf(Y;R)\longra MH_{l-k}(Y;R),
\end{split}
\label{kh4eq90}
\e
by $[\al]\cap[\be]=[\al\cap\be]$ for $\al\in MC^k_?(Y;R)$, $\be\in MC_l^?(Y;R)$ with~$\d\al=\pd\be=0$.
 
\label{kh4def19}
\end{dfn}

\begin{thm} Under the canonical isomorphisms $MH_k(Y;R)\cong H_k(Y;R),\ab\ldots,MH_k^\lf(Y;R)\cong H_k^\lf(Y;R)$ from {\rm\S\ref{kh41}--\S\ref{kh44},} the cap products in \eq{kh4eq90} are identified with the usual cap products \eq{kh2eq46} on ordinary (co)homology.
\label{kh4thm11}
\end{thm}

\begin{proof} The proof follows that of Theorem \ref{kh4thm9}, applying the method of \S\ref{kh262}. We use the soft resolution $\MC^\bu(Y;R)$ of $R_Y$ from Theorem \ref{kh4thm5}, and the soft resolution $\MC_{m-\bu}^\lf(Y;R)$ of $O_Y$ from Theorem \ref{kh4thm8}, where $m=\dim Y$, and the isomorphism $I_\cap:R_Y\ot_R O_Y\ra O_Y$ from \eq{kh2eq62}. We have $\Id_Y\cap[Y]=[Y]$ by \eq{kh4eq85}, where $\Id_Y\in MC^0(Y;R)=\MC^0(Y;R)(Y)$ and $[Y]\in MC_m^\lf(Y;R)=\MC_m^\lf(Y;R)(Y)$. From this we see that
\begin{equation*}
j_Y\ci I_\cap=\cap_{0,0}\ci (i_Y\ot j_Y):R_Y\ot_R O_Y\longra \MC_m^\lf(Y;R),
\end{equation*}
the analogue of the first equation of \eq{kh2eq64}. The analogue of the second follows from \eq{kh4eq84}. Hence the analogue of \eq{kh2eq65} gives an expression for $\cap$ at the (co)chain level in each line of \eq{kh4eq90}, under the isomorphisms
\begin{equation*}
H^k(Y;R)\cong H^k\bigl(\MC^*(Y;R)(Y),\d\bigr)=MH^k(Y;R),
\end{equation*}
and so on, using compactly-supported sheaf cohomology for $MH_*(Y;R)$ and $MH^*_\cs(Y;R)$. This coincides with the definition $[\al]\cap[\be]=[\al\cap\be]$ of $\cap$ on $MH^*(Y;R),\ldots$ in Definition \ref{kh4def19}.
\end{proof}

Next we discuss cross products on M-homology.

\begin{dfn} Let $Y_1,Y_2$ be manifolds, of dimensions $m_1,m_2$. For all $k,l$ in $\Z$, define $R$-bilinear maps
\e
\t:\cP MC^\lf_k(Y_1;R)\t \cP MC^\lf_l(Y_2;R)\longra \cP MC^\lf_{k+l}(Y_1\t Y_2;R)
\label{kh4eq91}
\e
on generators $[V,n,s,t]\in \cP MC^\lf_k(Y_1;R)$, $[V',n',s',t']\in \cP MC^\lf_l(Y_2;R)$ by
\e
\begin{split}
[V,n,s&,t]\t[V',n',s',t']=(-1)^{(l+m_2)n}\bigl[V\t V',n+n',\\
&(s_1\ci\pi_V,\ldots,s_n\ci\pi_V,s_1'\ci\pi_{V'},\ldots,s_{n'}'\ci\pi_{V'}),t\t t'\bigr].
\end{split}
\label{kh4eq92}
\e
Here \eq{kh4eq92} is apparently the same as \eq{kh4eq78}, except for the sign $(-1)^{(l+m_2)n}$, which is related to the sign $(-1)^{ln}$ in \eq{kh4eq78} in the same way the signs $(-1)^{(l+m)n}$ in \eq{kh4eq82} and $(-1)^{lm}$ in \eq{kh4eq61} are related, as explained in Definition~\ref{kh4def18}.

By a very similar proof to that of Proposition \ref{kh4prop11} in \S\ref{kh79}, we can show $\t$ in \eq{kh4eq91}--\eq{kh4eq92} is well defined.

As for \eq{kh4eq62}--\eq{kh4eq63} and \eq{kh4eq83}, we show that if $Y_1,Y_2,Y_3$ are manifolds and $\al\in \cP MC^\lf_j(Y_1;R)$, $\be\in MC^k(Y_2;R)$ and $\ga\in \cP MC^\lf_l(Y_3;R)$ then
\e
(\al\t\be)\t\ga=\al\t(\be\t\ga)\quad\text{in $\cP MC^\lf_{j+k+l}(Y_1\t Y_2\t Y_3;R)$.}
\label{kh4eq93}
\e

As for \eq{kh4eq64}--\eq{kh4eq65} and \eq{kh4eq84}, if $\al\in \cP MC^\lf_k(Y_1;R)$, $\be\in \cP MC^\lf_l(Y_2;R)$ then
\e
\pd(\al\t\be)=(\pd\al)\t\be+(-1)^k\al\t(\pd\be)\quad\text{in $\cP MC^\lf_{k+l-1}(Y_1\t Y_2;R)$.}
\label{kh4eq94}
\e

If $g_1:Y_1\ra Z_1$ and $g_2:Y_2\ra Z_2$ are proper smooth maps of manifolds and $\al\in \cP MC^\lf_k(Y_1;R)$, $\be\in \cP MC^\lf_l(Y_2;R)$ then in a similar way to equations \eq{kh4eq68}, \eq{kh4eq79} and \eq{kh4eq86} we find that
\e
(g_1\t g_2)_*(\al\t\be)=(g_1)_*(\al)\t(g_2)_*(\be)\quad\text{in $\cP MC^\lf_{k+l}(Z_1\t Z_2;R)$.}
\label{kh4eq95}
\e
When $\al\in MC_k(Y_1;R)\subseteq \cP MC^\lf_k(Y_1;R)$ and $\be\in MC_l(Y_2;R)\subseteq \cP MC^\lf_l(Y_2;R)$ then \eq{kh4eq95} holds in $MC_{k+l}(Z_1\t Z_2;R)$ without supposing $g_1,g_2$ proper.

In a similar way to \eq{kh4eq71} and \eq{kh4eq87}, by defining a presheaf morphism and sheafifying we can show that there is a unique sheaf morphism
\begin{equation*}
\t_{k,l}:\MC_k^\lf(Y_1;R)\boxtimes_R \MC_l^\lf(Y_2;R)\longra \MC_{k+l}^\lf(Y_1\t Y_2;R)
\end{equation*}
such that for all open $U_1\subseteq Y_1$ and $U_2\subseteq Y_2$ we have
\begin{align*}
\t_{k,l}(U_1\t U_2)\ci(\Pi\ot\Pi)=\Pi\ci\t:\cP MC^\lf_k(U_1;R)\ot_R \cP MC^\lf_l(U_2;R)& \\\longra MC^\lf_{k+l}(U_1\t U_2;R)&.
\end{align*}
The we define an $R$-bilinear morphism\e
\t:MC_k^\lf(Y_1;R)\t MC_l^\lf(Y_2;R)\longra MC_{k+l}^\lf(Y_1\t Y_2;R)
\label{kh4eq96}
\e
to be induced by the global sections $\t_{kl}(Y_1\t Y_2)$. It restricts to
\e
\t:MC_k(Y_1;R)\t MC_l(Y_2;R)\longra MC_{k+l}(Y_1\t Y_2;R).
\label{kh4eq97}
\e

Equations \eq{kh4eq92}--\eq{kh4eq95} descend to the sheafification, and so hold for cross products \eq{kh4eq96}--\eq{kh4eq97} on $MC_*^\lf(-;R),MC_*(-;R)$, where $g_1,g_2$ need not be proper for \eq{kh4eq95} on $MC_*(-;R)$. Equation \eq{kh4eq94} implies that the products \eq{kh4eq96}--\eq{kh4eq97} descend to homology, giving $R$-bilinear products
\e
\begin{split}
\t&:MH_k(Y_1;R)\t MH_l(Y_2;R)\longra MH_{k+l}(Y_1\t Y_2;R),\\
\t&:MH_k^\lf(Y_1;R)\t MH_l^\lf(Y_2;R)\longra MH_{k+l}^\lf(Y_1\t Y_2;R),
\end{split}
\label{kh4eq98}
\e
by $[\al]\t[\be]=[\al\t\be]$ for $\al\in MC_k^?(Y;R)$, $\be\in MC_l^?(Y;R)$ with~$\pd\al=\pd\be=0$.

\label{kh4def20}
\end{dfn}

Using the method of Theorems \ref{kh4thm9} and \ref{kh4thm11} and Corollary \ref{kh4cor2}, we prove:

\begin{thm} Under the canonical isomorphisms $MH_k(Y;R)\cong H_k(Y;R),$
$MH_k^\lf(Y;R)\cong H_k^\lf(Y;R)$ from Theorems\/ {\rm\ref{kh4thm3}} and\/ {\rm\ref{kh4thm7},} the cross products in \eq{kh4eq98} are identified with the usual cross products \eq{kh2eq53} on ordinary (locally finite) homology\/~$H_*(Y;R),H_*^\lf(Y;R)$.
\label{kh4thm12}
\end{thm}

\subsection{Poincar\'e duality and wrong way maps}
\label{kh47}

In \S\ref{kh28} we explained that if $Y$ is an oriented manifold of dimension $m$, it has a natural fundamental class $[[Y]]\in H_m^\lf(Y;R)$, and as in \eq{kh2eq46}, cap product with $[[Y]]$ induces Poincar\'e duality isomorphisms
\begin{gather*}
\Pd:H^k_\cs(Y;R)\,{\buildrel\cong\over\longra}\, H_{m-k}(Y;R),\quad
\Pd:H^k(Y;R)\,{\buildrel\cong\over\longra}\, H_{m-k}^\lf(Y;R),\\
\text{given by}\qquad \Pd:\al\longmapsto \al\cap[[Y]].
\end{gather*}

As in \S\ref{kh44}, we have a natural fundamental cycle $[Y]=[Y,0,0,\id_Y]$ in $MC_m^\lf(Y;R)$, with homology class the fundamental class $[[Y]]\in MH_m^\lf(Y;R)$. As in \S\ref{kh46}, cap products are also defined on (co)chains in M-(co)homology. Therefore we may define Poincar\'e duality morphisms at the (co)chain level
\begin{gather}
\Pd:MC^k_\cs(Y;R)\longra MC_{m-k}(Y;R),\quad
\Pd:MC^k(Y;R)\longra MC_{m-k}^\lf(Y;R),
\nonumber\\
\text{by}\qquad \Pd:\al\longmapsto \al\cap[Y].
\label{kh4eq99}
\end{gather}
These then induce isomorphisms $\Pd:MH^k_\cs(Y;R)\ra MH_{m-k}(Y;R)$,
$\Pd:MH^k(Y;R)\ra MH_{m-k}^\lf(Y;R)$, which are identified with the usual isomorphisms \eq{kh2eq69} by the identifications $MH_*(Y;R)\cong H_*(Y;R),\ldots$ in~\S\ref{kh41}--\S\ref{kh44}.

Combining the expressions $[Y]=[Y,0,0,\id_Y]$ for the fundamental cycle and \eq{kh4eq82} for the cap product, we see that the Poincar\'e duality morphism $\Pd$ in \eq{kh4eq99} is given on generators $[V,n,s,t]$ of $MC^k_\cs(Y;R)$ or $MC^k(Y;R)$ by
\e
\Pd:[V,n,s,t]\longmapsto [V,n,s,t].
\label{kh4eq100}
\e
Here on the l.h.s.\ of \eq{kh4eq100} for $[V,n,s,t]\in MC^k_?(Y;R)$, $t:V\ra Y$ is a cooriented submersion, but on the r.h.s.\ of \eq{kh4eq100} for $[V,n,s,t]\in MC_{m-k}^?(Y;R)$, $V$ is oriented, with orientation $o_V$ obtained by combining the coorientation $c_t$ on $t:V\ra Y$ and the given orientation $o_Y$ on $Y$, as in Assumption~\ref{kh3ass6}(c).

Also as in \S\ref{kh28}, if $f:Y\ra Z$ is a (perhaps proper) cooriented smooth map of manifolds, then we can define `wrong way maps' $f^!,f_!$, which have the opposite functoriality that one expects, covariant on cohomology and contravariant on homology. We now show that if $f$ is a submersion, we can define $f^!,f_!$ naturally at the (co)chain level on M-(co)homology.

\begin{dfn} Let $Y,Z$ be manifolds of dimensions $m,n$, and $f:Y\ra Z$ be a proper cooriented submersion. Define $f_!:MC_k(Z;R)\ra MC_{k-n+m}(Y;R)$ to be the $R$-linear morphism given on generators $[V,n,s,t]$ of $MC_k(Z;R)$ by
\e
f_![V,n,s,t]=[V',n,s',t']:=\bigl[V\t_{t,Z,f}Y,n,s\ci\pi_V,\pi_Y\bigr],
\label{kh4eq101}
\e
as in \eq{kh4eq21}. Here the fibre product $V'=V\t_{t,Z,f}Y$ exists by Assumption \ref{kh3ass5}(c) as $f$ is a submersion, and we give $V'$ the orientation induced by the orientation on $V$ and the coorientation on $f$, as in Assumption \ref{kh3ass6}(l). Also $s:V\ra\R^n$ is proper near $0\in\R^n$, and $\pi_V:V'\ra V$ is proper as $f$ is, so $s'=s\ci\pi_V:V'\ra\R^n$ is proper near $0\in\R^n$, as required for $[V',n,s',t']$ to be a generator of $MC_{k-n+m}(Y;R)$. An almost identical proof to that of Proposition \ref{kh4prop4} in \S\ref{kh75} shows that $f_!:MC_k(Z;R)\ra MC_{k-n+m}(Y;R)$ is well defined.

From \eq{kh4eq3} and \eq{kh4eq101} we see that
\begin{equation*}
f_!\ci\pd=\pd\ci f_!:MC_k(Z;R)\longra MC_{k-n+m-1}(Y;R).
\end{equation*}  
Thus the $f_!$ descend to morphisms $f_!:MH_k(Z;R)\ra MH_{k-n+m}(Y;R)$.

If $Z$ is oriented, so that combining the orientation on $Z$ with the coorientation on $f$ gives an orientation on $Y$, then comparing equations \eq{kh4eq21}, \eq{kh4eq100} and \eq{kh4eq101} we see the following diagram commutes on generators $[V,n,s,t]$ of $MC^{n-k}_\cs(Z;R)$, and so commutes:
\e
\begin{gathered}
\xymatrix@C=130pt@R=15pt{ 
*+[r]{MC^{n-k}_\cs(Z;R)} \ar[r]_\Pd \ar[d]^{f^*} & *+[l]{MC_k(Z;R)} \ar[d]_{f_!}  \\
*+[r]{MC^{n-k}_\cs(Y;R)} \ar[r]^\Pd  & *+[l]{MC_{k-n+m}(Y;R).\!}  }
\end{gathered}
\label{kh4eq102}
\e
This induces a corresponding commutative diagram on M-homology. 

Comparing this with the second diagram of \eq{kh2eq70}, and noting that the canonical isomorphisms $MH_*(Y;R)\cong H_*(Y;R)$, $MH^*_\cs(Y;R)\cong H^*_\cs(Y;R)$ from \S\ref{kh41} and \S\ref{kh44} identify pullbacks $f^*$ and Poincar\'e duality isomorphisms on M-(co)homology and ordinary (co)homology, we see that the morphisms $f_!:MH_k(Z;R)\ra MH_{k-n+m}(Y;R)$ above are identified with the usual wrong way morphisms $f_!:H_k(Z;R)\ra H_{k-n+m}(Y;R)$ in \S\ref{kh28} by the isomorphisms $MH_*(Y;R)\cong H_*(Y;R)$ and $MH_*(Z;R)\cong H_*(Z;R)$ from \S\ref{kh41}. By an argument involving compactly-supported cohomology twisted by orientation bundles, we can show this is also true without assuming $Z$ oriented.

In a similar way, by the ideas used to construct $f_*,f^*$ in \S\ref{kh41}--\S\ref{kh44}, we can also define morphisms $f_!:MC_k^\lf(Z;R)\ra MC_{k-n+m}^\lf(Y;R)$ without assuming $f$ proper, given on generators $[V,n,s,t]$ in $MC_k^\lf(Z;R)$ by \eq{kh4eq101}, and we can define $f^!:MC^k_\cs(Y;R)\ra MC^{k-m+n}_\cs(Z;R)$ for general cooriented submersions $f:Y\ra Z$ and $f^!:MC^k(Y;R)\ra MC^{k-m+n}(Z;R)$ for proper cooriented submersions $f:Y\ra Z$ acting on $[V,n,s,t]\in MC^k_?(Y;R)$ by
\begin{equation*}
f^![V,n,s,t]=[V,n,s,f\ci t],
\end{equation*}
as in \eq{kh4eq7}, and all of these descend to M-(co)homology, and are identified with the usual wrong-way morphisms $f_!,f^!$ on ordinary (co)homology.
\label{kh4def21}
\end{dfn}

\section{Other forms of M-(co)homology}
\label{kh5}

We now describe several modifications of the theory of integral M-(co)homology in \S\ref{kh4}. Section \ref{kh51} defines {\it rational M-homology and M-cohomology\/} $MH_*^\Q(Y;R),\ab MH^*_\Q(Y;R),\ldots,$ which is defined over a $\Q$-algebra $R$, and has better symmetry properties than integral M-(co)homology, including the fact that the cup product $\cup$ is supercommutative on M-cochains~$MC^*_\Q(Y;R)$.

Section \ref{kh52} discusses {\it de Rham M-homology and M-cohomology\/} $MH_*^\dR(Y;\R),\ab MH^*_\dR(Y;\R),\ldots,$ a combination of M-(co)homology and de Rham cohomology in which the (co)chains $[V,n,s,t,\om]$ include an exterior form $\om$ on $V$. Section \ref{kh53} extends \S\ref{kh4}--\S\ref{kh52} from manifolds $Y$ to effective orbifolds, and \S\ref{kh54} generalizes M-(co)homology to a bivariant theory, in the sense of \cite{FuMa} and~\S\ref{kh210}. 

\subsection{Rational M-homology and M-cohomology}
\label{kh51}

In \S\ref{kh4} we defined integral M-homology and integral M-cohomology $MH_*(Y;R)$, $ MH^*(Y;R)$, $MH^*_\cs(Y;R)$, $MH_*^\lf(Y;R)$ for $Y$ a manifold and $R$ a commutative ring, such as the integers $R=\Z$. As in Remark \ref{kh4rem9}, the cup product $\cup$ on $MC^*(Y;R)$ is associative, but not supercommutative. As in Remark \ref{kh2rem4}(ii), for some rings $R$ such as $\Z,\Z_2$ (though not $\Q$-algebras) we cannot make $\cup$ supercommutative on cochains, as Steenrod squares are an obstruction to this.

We now discuss a variation on \S\ref{kh4}, called {\it rational M-homology\/} $MH_*^\Q(Y;R)$, {\it rational M-cohomology\/} $MH^*_\Q(Y;R)$, {\it compactly-supported rational M-co\-hom\-ol\-ogy\/} $MH^*_{\cs,\Q}(Y;R)$, and {\it locally finite rational M-homology\/} $MH_*^{\lf,\Q}(Y;R)$. There are three main differences with~\S\ref{kh4}:
\begin{itemize}
\setlength{\itemsep}{0pt}
\setlength{\parsep}{0pt}
\item[(a)] The base ring $R$ is now required to be a $\Q$-algebra, e.g. $R=\Q,\R$ or $\C$.
\item[(b)] We include an extra relation (iii) in the definitions of rational M-chains $MC_*^\Q(Y;R)$, M-precochains $\cP MC^*_\Q(Y;R)$ and locally finite M-prechains $\cP MC_k^{\lf,\Q}(Y;R)$, in the analogues of Definitions \ref{kh4def1}, \ref{kh4def4} and~\ref{kh4def11}.
\item[(c)] The cup product $\cup$ on $MC^*_\Q(Y;R),MC^*_{\cs,\Q}(Y;R)$ is supercommutative.
\end{itemize}
We need $R$ to be a $\Q$-algebra as for generators $[V,n,s,t]$ with maps $s:V\ra\R^n$, some proofs involve averaging over the action of the symmetric group $S_n$ permuting the coordinates $(x_1,\ldots,x_n)$ of $\R^n$, so $1/\md{S_n}=1/n!$ must lie in~$R$.

The theory is intended for applications in which it is an advantage that $\cup$ (and similar operations, such as cross products $\t$) are supercommutative.

We now explain how to modify the material of \S\ref{kh4} to define rational M-(co)homology $MH_*^\Q(Y;R)$, $MH^*_\Q(Y;R)$, $MH^*_{\cs,\Q}(Y;R)$, $MH_*^{\lf,\Q}(Y;R)$. The changes are mostly cosmetic, requiring very little additional work. For the whole of \S\ref{kh51}, fix a $\Q$-algebra $R$, and a category $\tManc$ satisfying Assumptions \ref{kh3ass1}--\ref{kh3ass7} of~\S\ref{kh33}. 

\subsubsection{Rational M-homology $MH_*^\Q(Y;R)$}
\label{kh511}

We begin with rational M-homology~$MH_*^\Q(Y;R)$.

\begin{dfn} Let $Y$ be a manifold. As in Definition \ref{kh4def1}, define generators $[V,n,s,t]$ to be $\sim$-equivalence classes of quadruples $(V,n,s,t)$. For each $k\in\Z$, define the {\it rational M-chains\/} $MC_k^\Q(Y;R)$ to be the $R$-module generated by such $[V,n,s,t]$ with $\dim V=n+k$, subject to the relations Definition \ref{kh4def1}(i),(ii) used to define $MC_k(Y;R)$, and the additional relation:
\begin{itemize}
\setlength{\itemsep}{0pt}
\setlength{\parsep}{0pt}
\item[(iii)] For each generator $[V,n,s,t]$ with $s=(s_1,\ldots,s_n):V\ra\R^n$ and each permutation $\si\in S_n$ of $1,2,\ldots,n$ we have
\begin{equation*}
\bigl[V,n,(s_1,s_2,\ldots,s_n),t\bigr]=\sign(\si)\cdot \bigl[V,n,(s_{\si(1)},s_{\si(2)},\ldots,s_{\si(n)}),t\bigr]
\end{equation*}
in $MC_k^\Q(Y;R)$, where $\sign:S_n\ra\{\pm 1\}$ is the usual group morphism, which is characterized by the property that for each $\si\in S_n$, the diffeomorphism $\R^n\ra\R^n$ mapping $(x_1,x_2,\ldots,x_n)\mapsto (x_{\si(1)},x_{\si(2)},\ldots,x_{\si(n)})$ multiplies the orientation on $\R^n$ by~$\sign(\si)$.
\end{itemize}

As in Definition \ref{kh4def1}, define an $R$-linear map $\pd:MC_k^\Q(Y;R)\!\ra\! MC_{k-1}^\Q(Y;R)$ by \eq{kh4eq3}. To show that $\pd$ is well defined, we just have to add to the proof of Proposition \ref{kh4prop1} in \S\ref{kh71} the obvious fact that $\pd$ maps relation (iii) in $MC_k^\Q(Y;R)$ to relation (iii) in~$MC_{k-1}^\Q(Y;R)$.

The proof in Definition \ref{kh4def1} shows that $\pd\ci\pd=0$. Define the {\it rational M-homology groups\/} $MH_*^\Q(Y;R)$ to be the homology of~$\bigl(MC_*^\Q(Y;R),\pd\bigr)$. \label{kh5def1}
\end{dfn}

\begin{dfn} Let $Y$ be a manifold and $k\in\Z$. Then $MC_k^\Q(Y;R)$ in Definition \ref{kh5def1} is the quotient of $MC_k(Y;R)$ in Definition \ref{kh4def1} by relation (iii). Hence there is a natural, surjective $R$-module morphism
\e
\Pi:MC_k(Y;R)\longra MC_k^\Q(Y;R),\qquad \Pi:[V,n,s,t]\longmapsto[V,n,s,t].
\label{kh5eq1}
\e
This commutes with differentials $\pd$ on $MC_*(-;\!R),MC_*^\Q(-;\!R)$, and so induces morphisms on homology
\e
\Pi:MH_k(Y;R)\longra MH_k^\Q(Y;R).
\label{kh5eq2}
\e

Similarly, define an $R$-linear map $\io:MC_k^\Q(Y;R)\longra MC_k(Y;R)$ by
\e
\io\!:\!\bigl[V,n,(s_1,\ldots,s_n),t\bigr]\!\mapsto\!\frac{1}{n!}\!\sum_{\si\in S_n}\sign(\si)\!\cdot\!\bigl[V,n,(s_{\si(1)},s_{\si(2)},\ldots,s_{\si(n)}),t\bigr].
\label{kh5eq3}
\e
The factor $\frac{1}{n!}$ exists in $R$ as $R$ is a $\Q$-algebra. Then $\io$ maps relation Definition \ref{kh4def1}(i) for $n,i$ in $MC_k^\Q(Y;R)$ to the average over $i'=0,\ldots,n$ of relation Definition \ref{kh4def1}(i) for $n,i'$ in $MC_k(Y;R)$. Also $\io$ maps Definition \ref{kh4def1}(ii) in $MC_k^\Q(Y;R)$ to an average over $S_n$ of Definition \ref{kh4def1}(ii) in $MC_k(Y;R)$. And the r.h.s.\ of \eq{kh5eq3} is unchanged under Definition \ref{kh5def1}(iii). Hence $\io$ is well defined.

As $\io$ commutes with differentials $\pd$ on $MC_*^\Q(-;R),MC_*(-;R)$, it induces morphisms on homology
\e
\io:MH_k^\Q(Y;R)\longra MH_k(Y;R).
\label{kh5eq4}
\e
By equations \eq{kh5eq1} and \eq{kh5eq3} and Definition \ref{kh5def1}(iii) we see that
\begin{equation*}
\Pi\ci\io\bigl([V,n,s,t]\bigr)=\frac{1}{n!}\!\sum_{\si\in S_n}\sign(\si)\!\cdot\!\bigl[V,n,(s_{\si(1)},s_{\si(2)},\ldots,s_{\si(n)}),t\bigr]=[V,n,s,t].
\end{equation*}
Hence $\Pi\ci\io=\id:MC_k^\Q(Y;R)\longra MC_k^\Q(Y;R)$, and $\io$ is a right inverse for $\Pi$. This gives a canonical splitting
\e
MC_k(Y;R)\cong MC_k^\Q(Y;R)\op MC_k^\Q(Y;R)^\perp,
\label{kh5eq5}
\e
where $MC_k^\Q(Y;R)^\perp=\Ker\bigl(\Pi:MC_k(Y;R)\ra MC_k^\Q(Y;R)\bigr)$. The splitting \eq{kh5eq5} is preserved by $\pd$, as $\Pi,\io$ commute with $\pd$, so induces a splitting
\e
MH_k(Y;R)\cong MH_k^\Q(Y;R)\op MH_k^\Q(Y;R)^\perp.
\label{kh5eq6}
\e

\label{kh5def2}
\end{dfn}

Theorem \ref{kh4thm2} says that $MH_k(*;R)=0$ for $0\ne k\in\Z$, so that $MH_k^\Q(*;R)=0$ for $0\ne k\in\Z$ by \eq{kh5eq6}, and $MH_0(*;R)\cong R$, where $1\in R$ is identified with the cohomology class $[[*]]$ of $[*]\in MC_0(*;R)$. In the splitting \eq{kh5eq5}, $[*]$ lies in $MC_0^\Q(*;R)\subseteq MC_0(*;R)$, so in \eq{kh5eq6} $[[*]]$ lies in $MH_0^\Q(*;R)\subseteq MH_0(*;R)$, giving $MH_0^\Q(*;R)\cong R$. This proves:

\begin{cor} The analogue of Theorem\/ {\rm\ref{kh4thm2}} holds for rational M-homology.

\label{kh5cor1}
\end{cor}

To show $f_*:MC_k(Y_1;R)\ra MC_k(Y_2;R)$ in Definition \ref{kh4def2} is well defined on $MC_k^\Q(Y_i;R)$ we note that $f_*$ obviously takes relation (iii) in $MC_k^\Q(Y_1;R)$ to relation (iii) in $MC_k^\Q(Y_2;R)$. The rest of \S\ref{kh41}, and the proofs of Propositions \ref{kh4prop1} and Theorem \ref{kh4thm1} in \S\ref{kh71}--\S\ref{kh73}, need no changes in the rational case. Therefore as in Theorem \ref{kh4thm1} the complex of rational M-chains $\bigl(MC_*^\Q(Y;R),\pd\bigr)$ is the global sections of a complex $\bigl(\uMC_*^\Q(Y;R),\pd\bigr)$ of flabby cosheaves of $R$-modules on $Y$. Also, as for Theorem \ref{kh4thm3}, we have:

\begin{thm} For any $\Q$-algebra $R,$ rational M-homology $MH_k^\Q(-;R)$ is a homology theory of manifolds. There are canonical isomorphisms $MH_k^\Q(Y;R)\cong H_k(Y;R),$ $MH_k^\Q(Y,Z;R)\cong H_k(Y,Z;R)$ for all\/ $Y,Z,k,$ preserving the data $f_*,\pd$ and isomorphisms $MH_0^\Q(*;R)\cong R\cong H_0(*;R),$ where $H_*(-;R)$ is any other homology theory of manifolds over $R,$ such as singular homology\/~$H_*^\rsi(-;R)$.
\label{kh5thm1}
\end{thm}

The morphisms $\Pi,\io$ on M-homology in \eq{kh5eq2} and \eq{kh5eq4} commute with the isomorphisms $H_k^\ssi(Y;R)\cong MH_k(Y;R)$, $H_k^\ssi(Y;R)\cong MH_k^\Q(Y;R)$ from Example \ref{kh4ex1}, and so are the canonical isomorphisms~$MH_k(Y;R)\cong MH_k^\Q(Y;R)$.

\subsubsection{Rational M-cohomology $MH^*_\Q(Y;R)$}
\label{kh512}

Next we discuss rational M-cohomology~$MH^*_\Q(Y;R)$.

\begin{dfn} Let $Y$ be a manifold, of dimension $m$. As in Definition \ref{kh4def4}, define generators $[V,n,s,t]$ to be $\sim$-equivalence classes of quadruples $(V,n,s,t)$. For each $k\in\Z$, define the {\it rational M-precochains\/} $\cP MC^k_\Q(Y;R)$ to be the $R$-module generated by such $[V,n,s,t]$ with $\dim V+k=m+n$, subject to the relations Definition \ref{kh4def4}(i),(ii) used to define $\cP MC^k(Y;R)$, and the additional relation Definition \ref{kh5def1}(iii) (though note that the notion of generator $[V,n,s,t]$ we use now is different to that used in Definition~\ref{kh5def1}).

Define {\it differentials\/} $\d:\cP MC^k_\Q(Y;R)\ra \cP MC^{k+1}_\Q(Y;R)$ with $\d\ci\d=0$, the {\it identity precocycle\/} $\Id_Y=[Y,0,0,\id_Y]$ in $\cP MC^0_\Q(Y;R)$, and {\it pullbacks\/} $f^*:\cP MC^k_\Q(Y_2;R)\ra\cP MC^k_\Q(Y_1;R)$ by smooth maps of manifolds $f:Y_1\ra Y_2$ as in Definition \ref{kh4def5}. To show $\d,f^*$ are well-defined, we add to the proofs of Propositions \ref{kh4prop1} and \ref{kh4prop4} in \S\ref{kh71} and \S\ref{kh75} the obvious fact that $\d,f^*$ map relation (iii) in the domain to relation (iii) in the target.

The proof of Proposition \ref{kh4prop5} in \S\ref{kh76} needs no changes in the rational case. As in Definition \ref{kh4def6}, we define a c-soft strong presheaf $\cP\MC^k_\Q(Y;R)$ of $R$-modules on $Y$ by $\cP\MC^k_\Q(Y;R)(U)=\cP MC^k_\Q(U;R)$ for all open $U\subseteq Y$, and $\rho_{UU'}:\cP\MC^k_\Q(Y;R)(U)\ra\cP\MC^k_\Q(Y;R)(U')$ is $\rho_{UU'}=i^*:\cP MC^k_\Q(U;R)\ra\cP MC^k_\Q(U';R)$ for all open $U'\subseteq U\subseteq Y$, with $i:U'\hookra U$ the inclusion. Then we define $\MC^k_\Q(Y;R)$ to be the sheafification of $\cP\MC^k_\Q(Y;R)$, and we define the {\it rational M-cochains\/} $MC^k_\Q(Y;R)$ to be the global sections $\MC^k_\Q(Y;R)(Y)$. Then we have $\MC^k_\Q(Y;R)(U)=MC^k_\Q(U;R)$ for all open $U\subseteq Y$.
As in \eq{kh4eq23}, we may write $MC^k_\Q(Y;R)$ as an inverse limit
\begin{equation*}
MC^k_\Q(Y;R)=\mathop{\underleftarrow{\lim}\,}\nolimits_{\text{$U:U\subseteq Y$ open, $\bar U$ is compact}}\cP MC^k_\Q(U;R).
\end{equation*}

As in Definition \ref{kh4def6}, differentials $\d$, identities $\Id_Y$ and pullbacks $f^*$ descend to $MC^*_\Q(-;R)$ with $\d\ci\d\!=\!0$. Define the {\it rational M-cohomology groups\/} $MH^*_\Q(Y;R)$ to be the cohomology of the cochain complex~$\bigl(MC^*_\Q(Y;R),\d\bigr)$.

As in Definition \ref{kh5def2} we can define $R$-linear morphisms
\begin{equation*}
\Pi:\cP MC^k(Y;R)\longra\cP MC^k_\Q(Y;R),\quad \io:\cP MC^k_\Q(Y;R)\longra \cP MC^k(Y;R)
\end{equation*}
by the same formulae as in \eq{kh5eq1} and \eq{kh5eq3}, and these induce morphisms
\e
\Pi:MC^k(Y;R)\longra MC^k_\Q(Y;R),\quad \io:MC^k_\Q(Y;R)\longra MC^k(Y;R).
\label{kh5eq7}
\e
They satisfy $\Pi\ci\io=\id$, and commute with $\d$ and so descend to M-cohomology.

The rest of \S\ref{kh42}, and the proofs of Propositions \ref{kh4prop4}, \ref{kh4prop5} and \ref{kh4prop6} in \S\ref{kh7}, need no changes in the rational case. Therefore as in Definition \ref{kh4def6} and Theorem \ref{kh4thm5}, the complex of rational M-cochains $\bigl(MC^*_\Q(Y;R),\d\bigr)$ is the global sections of a complex $\bigl(\MC^*_\Q(Y;R),\d\bigr)$ of soft sheaves of $R$-modules on $Y$, which is a soft resolution of $R_Y$. Also, as for Theorem \ref{kh4thm4}, we have:
\label{kh5def3}
\end{dfn}

\begin{thm} For any $\Q$-algebra $R,$ rational M-cohomology is a cohomology theory of manifolds. There are canonical isomorphisms $MH^k_\Q(Y;R)\!\cong\! H^k(Y;R),$ $MH^k_\Q(Y,Z;R)\!\cong\! H^k(Y,Z;R)$ for all\/ $Y,Z,k,$ preserving the data $f^*,\d$ and isomorphisms $MH^0_\Q(*;R)\cong R\cong H^0(*;R),$ where $H^*(-;R)$ is any other cohomology theory of manifolds over $R,$ such as singular cohomology\/ $H^*_\rsi(Y;R)$ or sheaf cohomology\/~$H^*(Y,R_Y)$.
\label{kh5thm2}
\end{thm}

\subsubsection{Compactly-supported rational M-cohomology $MH^*_{\cs,\Q}(Y;R)$}
\label{kh513}

The rational analogue of \S\ref{kh43} requires only cosmetic changes. In Definition \ref{kh4def9} we define the {\it compactly-supported rational M-cochains\/} $MC^k_{\cs,\Q}(Y;R)\subseteq MC^k_\Q(Y;R)$ to be the $R$-submodule of compactly-supported global sections of $\MC^k_\Q(Y;R)$, and then we define {\it compactly-supported rational M-cohomology\/} $MH^*_{\cs,\Q}(Y;R)$ to be the cohomology of $\bigl(MC^*_{\cs,\Q}(Y;R),\d\bigr)$. Then as in Theorem \ref{kh4thm6}, compactly-supported rational M-cohomology is a compactly-supported cohomology theory of manifolds, so there are canonical, functorial isomorphisms $MH^k_{\cs,\Q}(Y;R)\cong H^k_\cs(Y;R)$ for all $Y,k$.

The rational analogue of Proposition \ref{kh4prop8} says that if $Y$ is a manifold and $k\in\Z$, then as an $R$-module, $MC^k_{\cs,\Q}(Y;R)$ is generated by compact generators $[V,n,s,t],$ subject only to relations Definition \ref{kh4def4}(i),(ii) and Definition \ref{kh5def1}(iii) applied to compact generators. Its proof in \S\ref{kh78} requires only cosmetic changes.  

\subsubsection{Locally finite rational M-homology $MH_*^{\lf,\Q}(Y;R)$}
\label{kh514}

Again, the rational analogue of \S\ref{kh44} requires only cosmetic changes. We first define {\it locally finite rational M-prechains\/} $\cP MC_k^{\lf,\Q}(Y;R)$ as in Definition \ref{kh4def11}, but imposing the additional relation Definition \ref{kh5def1}(iii). Proposition \ref{kh4prop9} holds in the rational case. As in Definition \ref{kh4def12}, we define a c-soft strong presheaf $\cP\MC_k^{\lf,\Q}(Y;R)$ of $R$-modules on $Y$ by $\cP\MC_k^{\lf,\Q}(Y;R)(U)=\cP MC_k^{\lf,\Q}(U;R)$ for all open $U\subseteq Y$, and $\rho_{UU'}:\cP\MC_k^{\lf,\Q}(Y;R)(U)\ra\cP\MC_k^{\lf,\Q}(Y;R)(U')$ is $\rho_{UU'}=i^*:\cP MC_k^{\lf,\Q}(U;R)\ra\cP MC_k^{\lf,\Q}(U';R)$ for all open $U'\subseteq U\subseteq Y$, with $i:U'\hookra U$ the inclusion. Then we define $\MC_k^{\lf,\Q}(Y;R)$ to be the sheafification of $\cP\MC_k^{\lf,\Q}(Y;R)$, and we define the {\it locally finite rational M-chains\/} $MC_k^{\lf,\Q}(Y;R)$ to be the global sections $\MC_k^{\lf,\Q}(Y;R)(Y)$. As in \eq{kh4eq48}, we may write $MC_k^{\lf,\Q}(Y;R)$ as an inverse limit
\begin{equation*}
MC_k^{\lf,\Q}(Y;R)\cong\mathop{\underleftarrow{\lim}\,}\nolimits_{\text{$U:U\subseteq Y$ open, $\bar U$ is compact}}\cP MC_k^{\lf,\Q}(U;R).
\end{equation*}

We define the {\it locally finite rational M-homology groups\/} $MH_*^{\lf,\Q}(Y;R)$ to be the homology of $\bigl(MC_*^{\lf,\Q}(Y;R),\pd\bigr)$. As in Theorem \ref{kh4thm7}, locally finite rational M-homology is a locally finite homology theory of manifolds, so there are canonical, functorial isomorphisms $MH_*^{\lf,\Q}(Y;R)\cong H_k^\lf(Y;R)$ for all~$Y,k$.

\subsubsection{Cup, cap and cross products, and Poincar\'e duality}
\label{kh515}

The analogues of \S\ref{kh45}--\S\ref{kh47} and the proof of Proposition \ref{kh4prop11} in \S\ref{kh79} also require only cosmetic changes in the rational case. Thus we define cup products $\cup$, cap products $\cap$, cross products $\t$, and Poincar\'e duality morphisms $\Pd$ on rational M-(co)chains and rational M-(co)homology, which as in Theorems \ref{kh4thm9}, \ref{kh4thm11}, \ref{kh4thm12} and Corollary \ref{kh4cor2} agree with the usual $\cup,\cap,\t,\Pd$ under the canonical isomorphisms $MH_k^\Q(Y;R)\cong H_k(Y;R),\ldots.$

The maps $\Pi:MC_*(Y;R)\ra MC_*^\Q(Y;R)$, $\Pi:MC^*(Y;R)\ra MC^*_\Q(Y;R)$ in \eq{kh5eq1} and \eq{kh5eq7} preserve the products $\cup,\cap,\t$ at the (co)chain level. However, the maps $\io:MC_*^\Q(Y;R)\ra MC_*(Y;R)$, $\io:MC^*_\Q(Y;R)\longra MC^*(Y;R)$ in \eq{kh5eq3} and \eq{kh5eq7} {\it do not preserve\/ $\cup,\cap,\t$ at the (co)chain level}, since the averaging over $S_n$ in \eq{kh5eq3} is not compatible with the definitions of $\cup,\cap,\t$ in~\S\ref{kh45}--\S\ref{kh46}.

The most important new feature is that if $Y$ is a manifold of dimension $m$ and $[V,n,s,t]\in \cP MC^k_\Q(Y;R)$, $[V',n',s',t']\in \cP MC^l_\Q(Y;R)$ are generators, so that $\dim V=m+n-k$, $\dim V'=m+n'-l$, then 
\ea
[V&,n,s,t]\cup[V',n',s',t']=(-1)^{ln}\bigl[V\t_{t,Y,t'}V',n+n',
\nonumber\\
&\qquad (s_1\ci\pi_V,\ldots,s_n\ci\pi_V,s_1'\ci\pi_{V'},\ldots,s_{n'}'\ci\pi_{V'}),t\ci\pi_V\bigr]
\nonumber\\
&=(-1)^{ln}\cdot (-1)^{(n-k)(n'-l)}\cdot (-1)^{nn'}
\bigl[V'\t_{t',Y,t}V,n+n',
\nonumber\\
&\qquad (s_1'\ci\pi_{V'},\ldots,s_{n'}'\ci\pi_{V'},s_1\ci\pi_V,\ldots,s_n\ci\pi_V),t'\ci\pi_{V'}\bigr]
\nonumber\\
&=(-1)^{kl} [V',n',s',t']\cup [V,n,s,t],
\label{kh5eq8}
\ea
where we use \eq{kh4eq61} in the first and third steps, and in the second step we use the natural isomorphism of manifolds with corners cooriented over $Y$
\begin{equation*}
V\t_{t,Y,t'}V'\cong (-1)^{(\dim V-\dim Y)(\dim V'-\dim Y)}\, V'\t_{t',Y,t}V,
\end{equation*}
by equation \eq{kh3eq3} and Assumption \ref{kh3ass6}(m), and we apply relation Definition \ref{kh5def1}(iii) to $[V\t_{t,Y,t'}V',n+n',(s_1\ci\pi_V,\ldots,s_n\ci\pi_V,s_1'\ci\pi_{V'},\ldots,s_{n'}'\ci\pi_{V'}),t\ci\pi_V]$ in $\cP MC^{k+l}(Y;R)$ with permutation $\si\in S_{n+n'}$ mapping
\begin{equation*}
\si:(1,2,\ldots,n,n+1,\ldots,n+n')\longmapsto(n+1,n+2,\ldots,n+n',1,2,\ldots,n),
\end{equation*}
which has $\sign(\si)=(-1)^{nn'}$.

By equation \eq{kh5eq8}, for $\al\in\cP MC^k_\Q(Y;R)$, $\be\in\cP MC^l_\Q(Y;R)$ we have
\begin{equation*}
\al\cup\be=(-1)^{kl}\be\cup\al.
\end{equation*}
This descends to $\al$ in $MC^k_\Q(Y;R)$ or $MC^k_{\cs,\Q}(Y;R)$ and $\be$ in $MC^l_\Q(Y;R)$ or $MC^l_{\cs,\Q}(Y;R)$. So we have:

\begin{prop} For any $\Q$-algebra $R$ and any manifold\/ $Y,$ the cup products $\cup$ on (compactly-supported) rational M-cochains $MC^*_\Q(Y;R),MC^*_{\cs,\Q}(Y;R)$ are supercommutative, at the cochain level.
\label{kh5prop1}
\end{prop}

Similarly, under the natural isomorphisms
\begin{align*}
MC^*_\Q(Y_1\t Y_2;R)&\cong MC^*_\Q(Y_2\t Y_1;R),\\ 
MC_*^\Q(Y_1\t Y_2;R)&\cong MC_*^\Q(Y_2\t Y_1;R), 
\end{align*}
the cross products on rational M-(co)chains are supercommutative:
\begin{align*}
\t:MC^k_\Q(Y_1;R)\t MC^l_\Q(Y_2;R)&\longra MC^{k+l}_\Q(Y_1\t Y_2;R),\\
\t:MC_k^\Q(Y_1;R)\t MC_l^\Q(Y_2;R)&\longra MC_{k+l}^\Q(Y_1\t Y_2;R).
\end{align*}
 
Note the contrast with integral M-cochains: as in Remark \ref{kh4rem9}, the cup product $\cup$ is generally not supercommutative on $MC^*(Y;R),MC^*_\cs(Y;R)$ even if $R$ is a $\Q$-algebra. We can now enhance Theorem \ref{kh4thm10} in \S\ref{kh45} to work with cdgas over $R$ rather than dgas over $R$, yielding:

\begin{thm} For any $\Q$-algebra $R$ and manifold\/ $Y,$ the cdga $\bigl(MC^*_\Q(Y;R),\ab\d,\ab\cup,\ab\Id_Y\bigr)$ over $R$ is equivalent in $\cdga_R^\iy$ to the `usual' cdga over $R$ associated to $Y$ in topology, as represented for instance by Sullivan's cdga $A_{PL}(Y;R)$ of polynomial differential forms on $Y$ with coefficients in $R,$ as in {\rm\cite[\S 10(c)]{FHT2}, \cite{Sull},} or by the de Rham cdga $\bigl(C^*_\dR(Y;\R),\d,\w,1_Y\bigr)$ from Example\/ {\rm\ref{kh2ex7}} when~$R=\R$.

Therefore the rational homotopy type of\/ $Y,$ and topological invariants of\/ $Y$ depending on the cdga up to equivalence, such as Massey products, may be computed using the cdga $\bigl(MC^*_\Q(Y;R),\ab\d,\ab\cup,\ab\Id_Y\bigr),$ and will give the correct answers under the canonical isomorphism $MH^*_\Q(Y;R)\cong H^*(Y;R)$ from Theorem\/~{\rm\ref{kh5thm2}}.
\label{kh5thm3}
\end{thm}

\subsection{De Rham M-homology and M-cohomology}
\label{kh52}

We now define variants of M-(co)homology called {\it de Rham M-(co)homology\/} $MH_*^\dR(Y;\R),MH^*_\dR(Y;\R),\ldots,$ which include exterior forms $\om$ on $V$ in generators $[V,n,s,t,\om]$, and generalize de Rham cohomology in Example \ref{kh2ex7}. This works for homology as well as cohomology. For simplicity we take the base ring to be $R=\R$, although we could extend to $R$ any $\R$-algebra in the obvious way. As $\R$ is a $\Q$-algebra we include the analogue of relation Definition \ref{kh5def1}(iii) in \S\ref{kh51} throughout, which makes cup products supercommutative on cochains.

For the whole of \S\ref{kh52}, fix a category $\tManc$ satisfying Assumptions \ref{kh3ass1}--\ref{kh3ass7} of \S\ref{kh33}, and Assumption \ref{kh3ass8} in \S\ref{kh34}, which defines $k$-forms $\om\in\Om^k(V)$ on manifolds with corners $V\in\tManc$.

\subsubsection{De Rham M-homology $MH_*^\dR(Y;\R)$}
\label{kh521}

Here is the analogue of the first half of Definition~\ref{kh4def1}.

\begin{dfn} Let $Y$ be a manifold. Consider quintuples $(V,n,s,t,\om)$, where $V$ is an oriented manifold with corners (i.e.\ a pair $(V,o_V)$ with $V$ an object in $\tManc$ and $o_V$ an orientation on $V$, usually left implicit), and $n=0,1,\ldots,$ and $s:V\ra\R^n$ is a smooth map (morphism in $\tManc$), and $t:V\ra Y$ is a smooth map, and $\om\in\Om^p(V)$ is a $p$-form on $V$ for $p=0,1,\ldots,$ as in Assumption \ref{kh3ass8}, such that the continuous map $s\vert_{\supp\om}:\supp\om\ra\R^n$ is proper over an open neighbourhood of 0 in~$\R^n$.

Define an equivalence relation $\sim$ on such quintuples by $(V,n,s,t,\om)\sim(V',\ab n',\ab s',\ab t',\ab\om')$ if $n=n'$, and there exists an orientation-preserving diffeomorphism $f:V\ra V'$ with $s=s'\ci f$ and $t=t'\ci f$ and $f^*(\om')=\om$. Write $[V,n,s,t,\om]$ for the $\sim$-equivalence class of $(V,n,s,t,\om)$. We call $[V,n,s,t,\om]$ a {\it generator}.

For each $k\in\Z$, define the {\it de Rham M-chains\/} $MC^\dR_k(Y;\R)$ to be the $\R$-vector space generated by such $[V,n,s,t,\om]$ with $\dim V=n+k+\deg\om$, subject to the relations:
\begin{itemize}
\setlength{\itemsep}{0pt}
\setlength{\parsep}{0pt}
\item[(i)] For each generator $[V,n,s,t,\om]$ and each $i=0,\ldots,n$ we have
\begin{equation*}
[V,n,s,t,\om]=(-1)^{n-i}[V\t\R,n+1,s',t\ci\pi_V,\pi_V^*(\om)]
\end{equation*}
in $MC^\dR_k(Y;\R)$, where writing $s=(s_1,\ldots,s_n):V\ra\R^n$ with $s_j:V\ra\R$ and $\pi_V:V\t\R\ra V$, $\pi_\R:V\t\R\ra\R$ for the projections, then 
\begin{equation*}
s'=(s_1\ci\pi_V,\ldots,s_i\ci\pi_V,\pi_\R,s_{i+1}\ci\pi_V,\ldots,s_n\ci\pi_V):V\t\R\longra\R^{n+1},
\end{equation*}
and $V\t\R$ has the product orientation from Assumption \ref{kh3ass6}(f) of the given orientation on $V$ and the standard orientation on~$\R$.
\item[(ii)] Let $I$ be a finite indexing set, $a_i\in \R$ for $i\in I$, and $[V_i,n,s_i,t_i,\om_i]$, $i\in I$ be generators for $MC^\dR_k(Y;\R)$, all with the same $n$. Suppose there exists an open neighbourhood $X$ of $0$ in $\R^n$, such that $s_i\vert_{\supp\om_i}:\supp\om_i\ra\R^n$ is proper over $X$ for all $i\in I$, and the following condition holds:
\begin{itemize}
\setlength{\itemsep}{0pt}
\setlength{\parsep}{0pt}
\item[$(*)$] Suppose $\eta\in C^\iy\bigl(\La^{n+k}T^*(\R^n\t Y)\bigr)$ is an $(n+k)$-form on $\R^n\t Y$ with $\supp\eta\subseteq K\t Y$ for some compact subset $K\subseteq X\subseteq\R^n$. Then
\e
\sum_{i\in I}a_i\int_{V_i}(s_i,t_i)^*(\eta)\w\om_i=0
\quad\text{in $\R$.}
\label{kh5eq9}
\e
Here $(s_i,t_i)^*(\eta)\w\om_i$ is a form of degree $n+k+\deg\om_i=\dim V_i$ on $V_i$, with $\supp[(s_i,t_i)^*(\eta)\w\om_i]\subseteq(s_i,t_i)^{-1}[\supp\eta]\cap\supp\om_i$ by Assumption \ref{kh3ass8}(c). Since $\pi_{\R^n}(\supp\eta)\subseteq K\subseteq X\subseteq\R^n$ for $K$ compact and $s_i\vert_{\supp\om_i}:\supp\om_i\ra\R^n$ is proper over $X$, we see that $(s_i,t_i)^*(\eta)\w\om_i$ is compactly-supported, so $\int_{V_i}(s_i,t_i)^*(\eta)\w\om_i$ is defined in $\R$ by Assumption \ref{kh3ass8}(e), and \eq{kh5eq9} makes sense.
\end{itemize}
Then
\begin{equation*}
\sum_{i\in I}a_i\,[V_i,n,s_i,t_i,\om_i]=0\qquad\text{in $MC^\dR_k(Y;\R)$.}
\end{equation*}
\item[(iii)] For each generator $[V,n,s,t,\om]$ with $s=(s_1,\ldots,s_n):V\ra\R^n$ and each permutation $\si\in S_n$ of $1,2,\ldots,n$ we have
\begin{equation*}
\bigl[V,n,(s_1,s_2,\ldots,s_n),t,\om\bigr]=\sign(\si)\cdot \bigl[V,n,(s_{\si(1)},s_{\si(2)},\ldots,s_{\si(n)}),t,\om\bigr]
\end{equation*}
in $MC_k^\dR(Y;\R)$, where $\sign:S_n\ra\{\pm 1\}$ is the usual group morphism.
\end{itemize}
\label{kh5def4}
\end{dfn}

\begin{rem}{\bf(a)} Relations (i),(iii) above are just Definition \ref{kh4def1}(i) and Definition \ref{kh5def1}(iii), with forms $\om$ added in a trivial way. Relation (ii) replaces Definition \ref{kh4def1}(ii), but is significantly different. So to generalize \S\ref{kh4} and \S\ref{kh51} to de Rham M-(co)homology, the main task is to rewrite the material using Definitions \ref{kh4def1}(ii), \ref{kh4def4}(ii) and \ref{kh4def11}(ii) in a nontrivial way, in particular the proofs of Propositions \ref{kh4prop1}, \ref{kh4prop4} and \ref{kh4prop11} and Theorem \ref{kh4thm2} in \S\ref{kh71}, \S\ref{kh74}, \S\ref{kh75} and \S\ref{kh79}, using relation Definition \ref{kh5def4}(ii) and its M-cohomology analogue instead.

As an example, in Definition \ref{kh5def5} below we show that $\pd$ on $MC^\dR_*(Y;\R)$ is well defined, the analogue of the proof of Proposition \ref{kh4prop1} in \S\ref{kh71}. The proof is rather easier in the de Rham case.
\smallskip

\noindent{\bf(b)} As in Remark \ref{kh4rem2}, here are some easy consequences of relation (ii). If $[V,n,s,t,\om]$ is a generator in $MC^\dR_k(Y;\R)$ then
\begin{equation*}
[-V,n,s,t,\om]=-[V,n,s,t,\om],
\end{equation*}
where $-V$ is $V$ with the opposite orientation, if $a,b\in\R$ then
\e
a[V,n,s,t,\om]+b[V,n,s,t,\om']=[V,n,s,t,a\,\om+b\,\om'],
\label{kh5eq10}
\e
and if $[V_1,n,s_1,t_1,\om_1]$, $[V_2,n,s_2,t_2,\om_2]$ are generators then
\e
[V_1\amalg V_2,n,s_1\amalg s_2,t_1\amalg t_2,\om_1\amalg\om_2]=[V_1,n,s_1,t_1,\om_1]+[V_2,n,s_2,t_2,\om_2].
\label{kh5eq11}
\e

Note that \eq{kh5eq10} implies that we can write a general element of $MC^\dR_k(Y;\R)$ as $\sum_{i\in I}[V_i,n_i,s_i,t_i,\om_i]$, omitting our usual factors~$a_i\in\R$.
\smallskip

\noindent{\bf(c)} Let $Y$ be oriented of dimension $m$ and $[V,n,s,t,\om]\in MC_k^\dR(Y;\R)$ with $k\le m$. Suppose $X$ is an open neighbourhood of 0 in $\R^n$ such that 
$s\vert_{\supp\om}:\supp\om\ra\R^n$ is proper over $X$ (this is part of the definition of $[V,n,s,t,\om]$) and $(s,t):V\ra\R^n\t Y$ is a submersion over $X\t Y$ (this is an extra assumption). 

Give $\R^n\t Y$ the product orientation from the standard orientation on $\R^n$ and the orientation on $Y$, as in Assumption \ref{kh3ass6}(f),(k). The orientations on $V$ and $\R^n\t Y$ induce a coorientation on $(s,t):V\ra\R^n\t Y$, as in Assumption \ref{kh3ass6}(c). Then $(s,t)\vert_{s^{-1}(X)}:s^{-1}(X)\ra X\t Y$ is a cooriented submersion, and $\om\vert_{s^{-1}(X)}$ is a form on $s^{-1}(X)$ with $(s,t)\vert_{\supp\om\vert_{s^{-1}(X)}}:\supp\om\vert_{s^{-1}(X)}\ra X\t Y$ proper, so Assumption \ref{kh3ass8}(f) defines a pushforward form $\bigl((s,t)\vert_{s^{-1}(X)}\bigr){}_*(\om\vert_{s^{-1}(X)})$ on $X\t Y$ of degree $\dim Y-k\ge 0$. We claim that
\e
[V,n,s,t,\om]=\bigl[X\t Y,n,\pi_X,\pi_Y,\bigl((s,t)\vert_{s^{-1}(X)}\bigr){}_*(\om\vert_{s^{-1}(X)})\bigr].
\label{kh5eq12}
\e

This holds by Definition \ref{kh5def4}(ii), where if $\eta$ is an $(n+k)$-form on $\R^n\t Y$ with $\supp\eta\subseteq K\t Y$ for some compact $K\subseteq X\subseteq\R^n$ then equation \eq{kh3eq12} yields
\begin{equation*}
\int_V(s,t)^*(\eta)\w\om=\int_{X\t Y}\eta\w\bigl((s,t)\vert_{s^{-1}(X)}\bigr){}_*(\om\vert_{s^{-1}(X)}),
\end{equation*}
and this is equivalent to Definition \ref{kh5def4}(ii)$(*)$ for \eq{kh5eq12}.
\smallskip

\noindent{\bf(d)} The signs in many formulae in this section, such as \eq{kh5eq9}, \eq{kh5eq12}, \eq{kh5eq13}, \eq{kh5eq19}, \eq{kh5eq20}, \eq{kh5eq22}, \eq{kh5eq24}, \eq{kh5eq25}, \eq{kh5eq27}, \eq{kh5eq29}, \eq{kh5eq32}, \eq{kh5eq36}, \eq{kh5eq39}, and \eq{kh5eq50}, depend on separate orientation conventions for homology and for cohomology. Our conventions were chosen to make the signs in \eq{kh5eq9}, \eq{kh5eq12}, \eq{kh5eq19}, \eq{kh5eq20}, \eq{kh5eq25}, \eq{kh5eq29} and \eq{kh5eq32} equal to 1. It seems that any convention must give ugly-looking signs somewhere, such as in \eq{kh5eq22}, \eq{kh5eq27}, \eq{kh5eq39} and~\eq{kh5eq50}. 
\label{kh5rem1}
\end{rem}

Following the second half of Definition \ref{kh4def1}, we define $\pd$ on $MC^\dR_*(Y;\R)$, and de Rham M-homology~$MH^\dR_*(Y;\R)$.

\begin{dfn} Define $\pd:MC^\dR_k(Y;\R)\ra MC^\dR_{k-1}(Y;\R)$ to be the $\R$-linear map such that for all generators $[V,n,s,t,\om]$ of $MC^\dR_k(Y;\R)$ we have
\e
\pd[V,n,s,t,\om]=[\pd V,n,s\ci i_V,t\ci i_V,i_V^*(\om)]+(-1)^{n+k}[V,n,s,t,\d\om],
\label{kh5eq13}
\e
where the orientation on $\pd V$ is induced from that of $V$ as in Assumption \ref{kh3ass6}(h). To show this is well defined, we have to show it takes relations (i)--(iii) in $MC^\dR_k(Y;\R)$ to relations (i)--(iii) in $MC^\dR_{k-1}(Y;\R)$. For (i),(iii) this is obvious. 

For (ii), suppose that $\sum_{i\in I}a_i\,[V_i,n,s_i,t_i,\om_i]=0$ in $MC^\dR_k(Y;\R)$ by (ii), where (ii)$(*)$ holds for some open $0\in X\subseteq\R^n$. We will show that
\e
\begin{split}
&\sum_{i\in I}a_i[\pd V_i,n,s_i\ci i_{V_i},t_i\ci i_{V_i},i_{V_i}^*(\om_i)]\\
&+\sum_{i\in I}(-1)^{n+k}a_i\,[V_i,n,s_i,t_i,\d\om_i]=0 \qquad\text{in $MC^\dR_{k-1}(Y;\R)$.}
\end{split}
\label{kh5eq14}
\e
Let $\eta\in C^\iy\bigl(\La^{n+k-1}T^*(\R^n\t Y)\bigr)$ be an $(n+k-1)$-form on $\R^n\t Y$ with $\supp\eta\subseteq K\t Y$ for compact $K\subseteq X\subseteq\R^n$. Then Assumption \ref{kh3ass8}(g) gives
\begin{equation*}
\int_{V_i}\d[(s_i,t_i)^*(\eta)\w\om_i]=\int_{\pd V_i}i_{V_i}^*((s_i,t_i)^*(\eta)\w\om_i).
\end{equation*}
Rewriting this using Assumption \ref{kh3ass8}(a),(b) yields
\e
\begin{split}
\int_{\pd V_i}&(s_i\ci i_{V_i},t_i\ci i_{V_i})^*(\eta)\w i_{V_i}^*(\om_i)+(-1)^{n+k}\int_{V_i}(s_i,t_i)^*(\eta)\w(\d\om_i)\\
&=\int_{V_i}(s_i,t_i)^*(\d\eta)\w\om_i.
\end{split}
\label{kh5eq15}
\e
Multiplying \eq{kh5eq15} by $a_i$ and summing over $i\in I$ yields
\ea
&\sum_{i\in I}a_i\int_{\pd V_i}(s_i\ci i_{V_i},t_i\ci i_{V_i})^*(\eta)\w i_{V_i}^*(\om_i)
+\sum_{i\in I}(-1)^{n+k}a_i\int_{V_i}(s_i,t_i)^*(\eta)\w(\d\om_i)
\nonumber\\
&\quad=\sum_{i\in I}a_i\int_{V_i}(s_i,t_i)^*(\d\eta)\w\om_i.
\label{kh5eq16}
\ea

The r.h.s.\ of \eq{kh5eq16} is zero by Definition \ref{kh5def4}(ii)$(*)$ for $\sum_{i\in I}a_i\,[V_i,n,s_i,t_i,\om_i]$ with $\d\eta$ in place of $\eta$. So the l.h.s.\ of \eq{kh5eq16} is zero, which proves Definition \ref{kh5def4}(ii)$(*)$ for the l.h.s.\ of \eq{kh5eq14}. Thus relation (ii) in $MC^\dR_{k-1}(Y;\R)$ says that \eq{kh5eq14} holds. Hence $\pd$ takes relation (ii) in $MC^\dR_k(Y;\R)$ to relation (ii) in $MC^\dR_{k-1}(Y;\R)$, and is well defined.

For any generator $[V,n,s,t,\om]$ in $MC^\dR_k(Y;\R)$ we have
\ea
\pd&\ci\pd[V,n,s,t,\om]=[\pd^2V,n,s\ci i_V\ci i_{\pd V},t\ci i_V\ci i_{\pd V},i_{\pd V}^*\ci i_V^*(\om)]
\nonumber\\
&+(-1)^{n+k-1}[\pd V,n,s\ci i_V,t\ci i_V,\d\ci i_V^*(\om)]
\nonumber\\
&+(-1)^{n+k}[\pd V,n,s\ci i_V,t\ci i_V,i_V^*(\d\om)]+(-1)^{n+k}(-1)^{n+k-1}[V,n,s,t,\d\ci\d\om]
\nonumber\\
&=[\pd^2V,n,s\ci i_V\ci i_{\pd V},t\ci i_V\ci i_{\pd V},i_{\pd V}^*\ci i_V^*(\om)],
\label{kh5eq17}
\ea
where in the second step, the second and third terms cancel as $\d\ci i_V^*(\om)=i_V^*(\d\om)$ and $\deg(\d\om)=\deg\om+1$, and the fourth term is zero as $\d\ci\d=0$. Apply Definition \ref{kh5def4}(ii) to $[\pd^2V,n,s\ci i_V\ci i_{\pd V},t\ci i_V\ci i_{\pd V},i_{\pd V}^*\ci i_V^*(\om)]$. Equation \eq{kh5eq9} for suitable $\eta\in C^\iy\bigl(\La^{n+k}T^*(\R^n\t Y)\bigr)$ becomes
\e
\int_{\pd^2V}(s\ci i_V\ci i_{\pd V},t\ci i_V\ci i_{\pd V})^*(\eta)\w(i_V\ci i_{\pd V})^*(\om)=0.
\label{kh5eq18}
\e
Assumptions \ref{kh3ass3}(c) and \ref{kh3ass6}(i) give an orientation-reversing diffeomorphism $\si_V:\pd^2V\ra\pd^2V$ with $\si_V^2=\id_{\pd^2V}$ and $i_V\ci i_{\pd V}\ci\si_V=i_V\ci i_{\pd V}$. Using $\si_V^*$ we see that the l.h.s.\ of \eq{kh5eq18} is minus itself, proving \eq{kh5eq18}. So Definition \ref{kh5def4}(ii)$(*)$ holds for 
$\sum_{i\in I}a_i\,[V_i,n,s_i,t_i,\om_i]=[\pd^2V,n,s\ci i_V\ci i_{\pd V},t\ci i_V\ci i_{\pd V},i_{\pd V}^*\ci i_V^*(\om)]$, and (ii) gives $[\pd^2V,n,s\ci i_V\ci i_{\pd V},t\ci i_V\ci i_{\pd V},i_{\pd V}^*\ci i_V^*(\om)]=0$, and thus \eq{kh5eq17} yields $\pd\ci\pd[V,n,s,t,\om]=0$.

Hence $\pd\ci\pd=0:MC^\dR_k(Y;\R)\ra MC^\dR_{k-2}(Y;\R)$, and $\bigl(MC^\dR_*(Y;\R),\pd\bigr)$ is a chain complex. Define the {\it de Rham M-homology groups\/} $MH^\dR_*(Y;\R)$ to be the homology of this chain complex.

If $Y$ is compact and oriented with $\dim Y=m$, define the {\it fundamental cycle\/} $[Y]=[Y,0,0,\id_Y,1_Y]\in MC^\dR_m(Y;\R)$. Here $V=Y$ has the given orientation, and $s=0:V\ra\R^0$ is proper as $Y$ is compact, and $1_Y\in\Om^0(Y)$ is as in Assumption \ref{kh3ass8}(a). We have $\pd[Y]=0$ by \eq{kh5eq13} as $\pd Y=\es$ and $\d 1_Y=0$, so passing to homology gives the {\it fundamental class\/}~$[[Y]]\in MH^\dR_m(Y;\R)$.
\label{kh5def5}
\end{dfn}

Here is the analogue of Lemma \ref{kh4lem1}. It holds since in Definition \ref{kh5def4}(ii)$(*)$, if $k>\dim Y$ then $C^\iy\bigl(\La^{n+k}T^*(\R^n\t Y)\bigr)=0$, so $(*)$ holds trivially, and (ii) says that any element $\sum_{i\in I}a_i\,[V_i,n,s_i,t_i,\om_i]$ in $MC_k^\dR(Y;\R)$ is zero.

\begin{lem} For any manifold\/ $Y$ we have $MC_k^\dR(Y;\R)=0$ for\/ $k>\dim Y,$ so that\/ $MH_k^\dR(Y;\R)=0$ for\/~$k>\dim Y$.
\label{kh5lem1}
\end{lem}

Here is the analogue of Definition \ref{kh4def2}.

\begin{dfn} Let $f:Y_1\ra Y_2$ be a smooth map of manifolds. Define the {\it pushforward\/} $f_*:MC_k^\dR(Y_1;\R)\ra MC_k^\dR(Y_2;\R)$ for $k\in\Z$ to be the unique $\R$-linear map defined on generators $[V,n,s,t,\om]$ of $MC_k(Y_1;\R)$ by
\e
f_*[V,n,s,t,\om]=[V,n,s,f\ci t,\om].
\label{kh5eq19}
\e
To show $f_*$ is well-defined, we must show that it maps relations (i)--(iii) in $MC_k(Y_1;\R)$ to relations (i)--(iii) in $MC_k(Y_2;\R)$. For (i),(iii) this is obvious. For (ii), suppose $\sum_{i\in I}a_i\,[V_i,n,s_i,t_i,\om_i]=0$ in $MC_k(Y_1;\R)$ by relation (ii) using open $0\in X\subseteq\R^n$. Suppose $\eta\in C^\iy\bigl(\La^{n+k}T^*(\R^n\t Y_2)\bigr)$ with $\supp\eta\subseteq K\t Y_2$ for some compact subset $K\subseteq X\subseteq\R^n$. Then $(\id_{\R^n}\t f)^*(\eta)$ lies in $C^\iy\bigl(\La^{n+k}T^*(\R^n\t Y_1)\bigr)$ with $\supp(\id_{\R^n}\t f)^*(\eta)\subseteq K\t Y_1$, so \eq{kh5eq9} with $(\id_{\R^n}\t f)^*(\eta)$ in place of $\eta$ yields
\begin{equation*}
\sum_{i\in I}a_i\int_{V_i}(s_i,t_i)^*\ci(\id_{\R^n}\t f)^*(\eta)\w\om_i=0.
\end{equation*}
But this is \eq{kh5eq9} for $\sum_{i\in I}a_i\,[V_i,n,s_i,f\ci t_i,\om_i]$, so Definition \ref{kh5def4}(ii) implies that $\sum_{i\in I}a_i\,[V_i,n,s_i,f\ci t_i,\om_i]=0$, and $f^*$ maps (ii) to (ii), and is well defined.

Equations \eq{kh5eq13}, \eq{kh5eq19} give $f_*\ci\pd=\pd\ci f_*:MC_k^\dR(Y_1;\R)\ra MC_{k-1}^\dR(Y_2;\R)$. So the $f_*$ induce pushforwards $f_*:MH_k^\dR(Y_1;\R)\ra MH_k^\dR(Y_2;\R)$ on homology. If $g:Y_2\ra Y_3$ is another smooth map of manifolds then $(g\ci f)_*=g_*\ci f_*$, on both de Rham M-chains $MC_*^\dR(Y_i;\R)$ and M-homology $MH_*^\dR(Y_i;\R)$. Also $(\id_Y)_*$ is the identity on both $MC_*^\dR(Y;\R)$ and $MH_*^\dR(Y;\R)$.
\label{kh5def6}
\end{dfn}

The material in \S\ref{kh41} from Theorem \ref{kh4thm1} down to Lemma \ref{kh4lem2}, and the proofs of Theorem \ref{kh4thm1} in \S\ref{kh72} and Proposition \ref{kh4prop2} in \S\ref{kh73}, all require only very straightforward modifications for the de Rham case, and we leave the details to the reader. For instance, in equation \eq{kh7eq10} in \S\ref{kh72} we should instead write
\begin{align*}
&\Pi_{T\cup U}^{T,f-}:MC_k^\dR(T\cup U;\R)\longra MC_k^\dR(T;\R),\\
&\Pi_{T\cup U}^{T,f-}:\bigl[V,n,(s_1,\ldots,s_n),t,\om\bigr]\longmapsto \bigl[t^{-1}(T)\t(-\iy,0],n+1,\\
&\quad (s_1\ci\pi_{t^{-1}(T)},\ldots,s_n\ci\pi_{t^{-1}(T)},f\ci\pi_{t^{-1}(T)}+\pi_{(-\iy,0]}),t\ci\pi_{t^{-1}(T)},\pi_{t^{-1}(T)}^*(\om)\bigr].
\end{align*}

Therefore the complex of de Rham M-chains $\bigl(MC_*^\dR(Y;\R),\pd\bigr)$ is the global sections of a complex $\bigl(\uMC_*^\dR(Y;\R),\pd\bigr)$ of flabby cosheaves of $\R$-vector spaces on $Y$, as in Theorem \ref{kh4thm1}, and we have notions of relative de Rham M-chains $MC_k^\dR(Y,Z;\R)$ and M-homology $MH_k^\dR(Y,Z;\R)$ for open~$Z\subseteq Y$.

Here is the analogue of Theorem \ref{kh4thm2}, proved in~\S\ref{kh81}.

\begin{thm} We have $MH_k^\dR(*;\R)=0$ for all\/ $0\ne k\in\Z$. There is a canonical isomorphism $MH_0^\dR(*;\R)\cong\R$ identifying the fundamental class $[[*]]\in MH_0^\dR(*;\R)$ with\/~$1\in\R$.
\label{kh5thm4}
\end{thm}

As for Theorem \ref{kh4thm3} we deduce:

\begin{thm} De Rham M-homology is a homology theory of manifolds. There are canonical isomorphisms $MH_k^\dR(Y;\R)\cong H_k(Y;\R),$ $MH_k^\dR(Y,Z;\R)\cong H_k(Y,Z;\R)$ for all\/ $Y,Z,k,$ preserving the $f_*,\pd$ and isomorphisms $MH_0^\dR(*;\R)\ab\cong\R\cong H_0(*;\R),$ where $H_*(-;\R)$ is any other homology theory of manifolds over $\R,$ such as singular homology\/~$H_*^\rsi(-;\R)$.
\label{kh5thm5}
\end{thm}

\begin{ex} Let $Y$ be a manifold and $k\in\Z$. Define $\R$-linear maps
\e
\begin{split}
F_\Mh^\dRMh&:MC_k(Y;\R)\longra MC^\dR_k(Y;\R),\\ 
F_\QMh^\dRMh&:MC_k^\Q(Y;\R)\longra MC^\dR_k(Y;\R) \qquad\text{by}\\
F_\Mh^\dRMh,F_\QMh^\dRMh&:[V,n,s,t]\longmapsto [V,n,s,t,1_V]\quad \text{on generators $[V,n,s,t]$.}
\end{split}
\label{kh5eq20}
\e
We will show in Proposition \ref{kh5prop3} below that these are well defined.

Comparing \eq{kh4eq3}, \eq{kh5eq13} and \eq{kh5eq20} and noting that $\d 1_V=0$ we see that $F_\Mh^\dRMh\ci\pd=\pd\ci F_\Mh^\dRMh:MC_k(Y;\R)\ra MC^\dR_{k-1}(Y;\R)$, and similarly for $F_\QMh^\dRMh$, so \eq{kh5eq20} induces morphisms
\e
\begin{split}
F_\Mh^\dRMh&:MH_k(Y;\R)\longra MH^\dR_k(Y;\R),\\ 
F_\QMh^\dRMh&:MH_k^\Q(Y;\R)\longra MH^\dR_k(Y;\R).
\end{split}
\label{kh5eq21}
\e
These morphisms $F_\Mh^\dRMh,F_\QMh^\dRMh$ in \eq{kh5eq21} extend to relative homology, and commute with pushforwards $f_*$ and continuation maps $\pd:MC_k^?(Y,Z;\R)\ra MC_{k-1}^?(Z;\R)$. Since $F_\Mh^\dRMh([[*]])=F_\QMh^\dRMh([[*]])=[[*]]$, we see that $F_\Mh^\dRMh:MH_0(*;\R)\ra MH^\dR_0(*;\R)$, $F_\QMh^\dRMh:MH_0^\Q(*;\R)\ra MH^\dR_0(*;\R)$ are compatible with the isomorphisms $MH_0(*;\R)\cong MH_0^\Q(*;\R)\cong MH_0^\dR(*;\R)\cong\R$. Thus by the last part of Theorem \ref{kh2thm1}, the morphisms \eq{kh5eq21} are the canonical isomorphisms from Theorem~\ref{kh5thm5}.

We can compose $F_\Mh^\dRMh:MC_k(Y;\R)\ra MC^\dR_k(Y;\R)$ in \eq{kh5eq20} with the morphisms $F_\ssi^\Mh:C_k^\ssi(Y;\R)\ra MC_k(Y;\R)$, $\hat F_\ssi^\Mh:\hat C_k^\ssi(Y;\R)\ra MC_k(Y;\R)$ from Examples \ref{kh4ex1} and \ref{kh4ex2} to get explicit morphisms of chain complexes lifting the isomorphisms $H_*^\ssi(Y;\R)\cong MH^\dR_*(Y;\R)$, $\hat H_*^\ssi(Y;\R)\cong MH^\dR_*(Y;\R)$ given by Theorem~\ref{kh5thm5}.
\label{kh5ex1}
\end{ex}

The next proposition will be proved in~\S\ref{kh82}.

\begin{prop} $F_\Mh^\dRMh,F_\QMh^\dRMh$ in \eq{kh5eq20} above are well defined.
\label{kh5prop3}
\end{prop}

\subsubsection{De Rham M-cohomology $MH^*_\dR(Y;\R)$}
\label{kh522}

Here are the de Rham analogues of Definitions \ref{kh4def4} and~\ref{kh4def5}.

\begin{dfn} Let $Y$ be a manifold, of dimension $m$. Consider quintuples $(V,n,s,t,\om)$, where $V$ is a manifold with corners (object in $\tManc$), and $n\in\N,$ and $s:V\ra\R^n$ is a smooth map (morphism in $\tManc$), and $t:V\ra Y$ is a cooriented submersion (i.e.\ a pair $(t,c_t)$ of a submersion $t:V\ra Y$ in $\tManc$ and a coorientation $c_t$ for $t$, usually left implicit), and $\om\in\Om^p(V)$ is a $p$-form on $V$ for $p=0,1,\ldots,$ as in Assumption \ref{kh3ass8}, such that the continuous map $(s,t)\vert_{\supp\om}:\supp\om\ra\R^n\t Y$ is proper over an open neighbourhood of $\{0\}\t Y$ in~$\R^n\t Y$.

Define an equivalence relation $\sim$ on such quintuples by $(V,n,s,t,\om)\sim(V',\ab n',\ab s',\ab t',\ab\om')$ if $n=n'$, and there exists a diffeomorphism $f:V\ra V'$ with $s=s'\ci f$ and $t=t'\ci f$ and $f^*(\om')=\om$  such that the coorientations satisfy $c_t=c_{t'}\ci c_f$, where $c_f$ is the natural coorientation on $f$ from Assumption \ref{kh3ass6}(e). Write $[V,n,s,t,\om]$ for the $\sim$-equivalence class of $(V,n,s,t,\om)$. We call $[V,n,s,t,\om]$ a {\it generator}.

For each $k\in\Z$, define the {\it de Rham M-precochains\/} $\cP MC^k_\dR(Y;\R)$ to be the $\R$-vector space generated by such $[V,n,s,t,\om]$ with with $\dim V+k=\deg\om+m+n$, subject to the relations:
\begin{itemize}
\setlength{\itemsep}{0pt}
\setlength{\parsep}{0pt}
\item[(i)] For each generator $[V,n,s,t,\om]$ and each $i=0,\ldots,n$ we have
\begin{equation*}
[V,n,s,t,\om]=(-1)^{n-i}[V\t\R,n+1,s',t\ci\pi_V,\pi_V^*(\om)]
\end{equation*}
in $\cP MC^k_\dR(Y;\R)$, where writing $s=(s_1,\ldots,s_n):V\ra\R^n$ with $s_j:V\ra\R$ and $\pi_V:V\t\R\ra V$, $\pi_\R:V\t\R\ra\R$ for the projections, then 
\begin{equation*}
s'=(s_1\ci\pi_V,\ldots,s_i\ci\pi_V,\pi_\R,s_{i+1}\ci\pi_V,\ldots,s_n\ci\pi_V):V\t\R\longra\R^{n+1},
\end{equation*}
and $t\ci\pi_V$ has coorientation $c_{t\ci\pi_V}=c_t\ci c_{\pi_V}$, where $c_t$ is the given coorientation on $t:V\ra Y$, and $c_{\pi_V}$ is the coorientation on $\pi_V:V\t\R\ra V$ induced by the standard orientation on $\R$, as in Assumption~\ref{kh3ass6}(d),(f),(k).
\item[(ii)] Let $I$ be a finite indexing set, $a_i\in \R$ for $i\in I$, and $[V_i,n,s_i,t_i,\om_i]$, $i\in I$ be generators for $\cP MC^k_\dR(Y;\R)$, all with the same $n$. Suppose there exists an open neighbourhood $X$ of $\{0\}\t Y$ in $\R^n\t Y$, such that $(s_i,t_i)\vert_{\supp\om_i}:\supp\om_i\ra\R^n\t Y$ is proper over $X$ for all $i\in I$, and the following condition holds:
\begin{itemize}
\setlength{\itemsep}{0pt}
\setlength{\parsep}{0pt}
\item[$(*)$] Suppose $V'$ is an oriented manifold of dimension $l$, $t':V'\ra Y$ is smooth, and $\eta\in C^\iy\bigl(\La^{l+n-k}T^*(\R^n\t V')\bigr)$ with $\supp\eta$ a compact subset of $(\id_{\R^n}\t t')^{-1}(X)\subseteq\R^n\t V'$. Then
\e
\begin{split}
&\sum_{i\in I}(-1)^{l\deg\om_i+\deg\om_i(\deg\om_i-1)/2}a_i\cdot{}\\
&\int_{V_i\t_{t_i,Y,t'}V'}(s_i\ci\pi_{V_i},\pi_{V'})^*(\eta)\w\pi_{V_i}^*(\om_i)=0.
\end{split}
\label{kh5eq22}
\e
Here the fibre product $V_i\t_{t_i,Y,t'}V'$ exists in $\tManc$ by Assumption \ref{kh3ass5}(c)) as $t_i:V_i\ra Y$ is a submersion, with dimension
\begin{align*}
\dim V_i+\dim V'-\dim Y&=(\deg\om_i+m+n-k)+l-m\\
&=\deg\om_i+l+n-k,
\end{align*}
and Assumption \ref{kh3ass6}(l) gives it an orientation by combining the coorientation $c_{t_i}$ on $t_i$ and the orientation on $V'$. By $\supp\eta$ compact in $(\id_{\R^n}\t t')^{-1}(X)$, and $(s_i,t_i)\vert_{\supp\om_i}:\supp\om_i\ra\R^n\t Y$ proper over $X$, and equation \eq{kh3eq8}, we see that $(s_i\ci\pi_{V_i},\pi_{V'})^*(\eta)\w\pi_{V_i}^*(\om_i)$ is compactly-supported in $V_i\t_YV'$, so that \eq{kh5eq22} is well defined by Assumption~\ref{kh3ass8}(e).
\end{itemize}
Then
\begin{equation*}
\sum_{i\in I}a_i\,[V_i,n,s_i,t_i,\om_i]=0\qquad\text{in $\cP MC^k_\dR(Y;\R)$.}
\end{equation*}
\item[(iii)] For each generator $[V,n,s,t,\om]$ with $s=(s_1,\ldots,s_n):V\ra\R^n$ and each permutation $\si\in S_n$ of $1,2,\ldots,n$ we have
\begin{equation*}
\bigl[V,n,(s_1,s_2,\ldots,s_n),t,\om\bigr]=\sign(\si)\cdot \bigl[V,n,(s_{\si(1)},s_{\si(2)},\ldots,s_{\si(n)}),t,\om\bigr]
\end{equation*}
in $\cP MC^k_\dR(Y;\R)$, where $\sign:S_n\ra\{\pm 1\}$ is the usual group morphism.
\end{itemize}

As in Remark \ref{kh5rem1}(b), some consequences of relation (ii) in $\cP MC^k_\dR(Y;\R)$ are equations \eq{kh5eq10}--\eq{kh5eq11} and the analogue of~\eq{kh4eq18}:
\e
\bigl[V,n,s,(t,-c_t),\om\bigr]=-\bigl[V,n,s,(t,c_t),\om\bigr].
\label{kh5eq23}
\e

\label{kh5def7}
\end{dfn}

\begin{dfn} For each manifold $Y$ and $k\in\Z$, as in equation \eq{kh5eq13} let $\d:\cP MC_\dR^k(Y;\R)\ra \cP MC_\dR^{k+1}(Y;\R)$ be the unique $\R$-linear map with
\e
\d[V,n,s,t,\om]=[\pd V,n,s\ci i_V,t\ci i_V,i_V^*(\om)]+(-1)^{k+n+\deg\om}[V,n,s,t,\d\om],
\label{kh5eq24}
\e
for all generators $[V,n,s,t,\om]$, where the coorientation on $t\ci i_V$ is $c_{t\ci i_V}=c_t\ci c_{i_V}$, for $c_t$ the given coorientation on $t$, and $c_{i_V}$ the coorientation for $i_V:\pd V\ra V$ from Assumption \ref{kh3ass6}(h). We show that $\d$ is well-defined with $\d\ci\d=0$ by similar proofs to those for $\pd$ in Definition~\ref{kh5def5}.

Define the {\it identity precocycle\/} $\Id_Y=[Y,0,0,\id_Y,1_Y]$ in $\cP MC_\dR^0(Y;\R)$. Here $t=\id_Y:Y\ra Y$ has the coorientation from Assumption \ref{kh3ass6}(e), and $(s,t):Y\ra\R^0\t Y$ is proper. We have $\d\,\Id_Y=0$ as $\pd Y=\es$ and~$\d 1_Y=0$.

Suppose $f:Y_1\ra Y_2$ is a smooth map of manifolds. For each $k\in\Z$, define the {\it pullback\/} $f^*:\cP MC_\dR^k(Y_2;\R)\ra\cP MC_\dR^k(Y_1;\R)$ to be the unique $\R$-linear map acting on generators $[V,n,s,t,\om]$ of $\cP MC_\dR^k(Y_2;\R)$ by
\e
f^*[V,n,s,t,\om]\!=\![V',n,s',t',\om']\!:=\!\bigl[V\t_{t,Y_2,f}Y_1,n,s\ci\pi_V,\pi_{Y_1},\pi_V^*(\om)\bigr],
\label{kh5eq25}
\e
where $V'=V\t_{t,Y_2,f}Y_1$ is the fibre product in $\tManc$, which exists by Assumption \ref{kh3ass5}(c) as $Y_2$ is a manifold and $t$ a submersion, with projections $\pi_V:V'\ra V$ and $\pi_{Y_1}:V'\ra Y_1$. Here $t'=\pi_{Y_1}:V'\ra Y_1$ is a submersion by Assumption \ref{kh3ass5}(c), and has a coorientation $c_{t'}$ determined by the given coorientation $c_t$ on $t:V\ra Y_2$ by Assumption~\ref{kh3ass6}(l). 

To show that $f^*$ is well defined (the analogue of Proposition \ref{kh4prop4}), we have to show that it maps relations Definition \ref{kh5def7}(i)--(iii) in $\cP MC^k_\dR(Y_2;\R)$ to relations (i)--(iii) in $\cP MC^k_\dR(Y_1;\R)$. For (i) this is obvious, as $(V\t\R)\t_{t\ci\pi_V,Y_2,f}Y_1\ab\cong (V\t_{t,Y_2,f}Y_1)\t\R$, and for (iii) it is also obvious. For (ii), suppose $\sum_{i\in I}a_i\,[V_i,n,s_i,t_i,\om_i]=0$ in $\cP MC^k_\dR(Y_2;\R)$ by relation (ii), using open $\{0\}\t Y_2\subseteq X_2\subseteq\R^n\t Y_2$. Define $X_1=(\id_{\R^n}\t f)^{-1}(X_2)\subseteq\R^n\t Y_1$. Then $X_1$ is an open neighbourhood of $\{0\}\t Y_1$ in~$\R^n\t Y_1$. 

Let $V'$ be an oriented manifold of dimension $l$, $t':V'\ra Y_1$ be smooth, and $\eta\in C^\iy\bigl(\La^{l+n-k}T^*(\R^n\t V')\bigr)$ with $\supp\eta$ compact in $(\id_{\R^n}\t t')^{-1}(X_1)$. Then
\e
\begin{split}
&\sum_{i\in I}(-1)^{l\deg\om_i+\deg\om_i(\deg\om_i-1)/2}a_i\cdot{} \\
&\int_{(V_i\t_{Y_2}Y_1)\t_{Y_1}V'}(s_i\ci\pi_{V_i}\ci\pi_{V_i\t_{Y_2}Y_1},\pi_{V'})^*(\eta)\w(\pi_{V_i\t_{Y_2}Y_1}^*\ci\pi_{V_i}^*(\om_i))\\
&=\sum_{i\in I}(-1)^{l\deg\om_i+\deg\om_i(\deg\om_i-1)/2}a_i\cdot{}\\
&\quad \int_{V_i\t_{t_i,Y_2,f\ci t'}V'}(s_i\ci\pi_{V_i},\pi_{V'})^*(\eta)\w\pi_{V_i}^*(\om_i)=0,
\end{split}
\label{kh5eq26}
\e
using the natural isomorphism $(V_i\t_{Y_2}Y_1)\t_{Y_1}V'\cong V_i\t_{Y_2}V'$ of fibre products in the first step, and \eq{kh5eq22} for $\sum_{i\in I}a_i\,[V_i,n,s_i,t_i,\om_i]=0$ in $\cP MC^k_\dR(Y_2;\R)$ with $f\ci t'$ in place of $t'$ in the second. But \eq{kh5eq26} is condition (ii)$(*)$ for
\begin{equation*}
\sum_{i\in I}a_i\,f^*[V_i,n,s_i,t_i,\om_i]=\sum_{i\in I}a_i\,\bigl[V_i\t_{t_i,Y_2,f}Y_1,n,s_i\ci\pi_V,\pi_{Y_1},\pi_V^*(\om_i)\bigr]=0
\end{equation*}
in $\cP MC^k_\dR(Y_1;\R)$. Hence $f^*$ takes relation (ii) in $\cP MC^k_\dR(Y_2;\R)$ to relation (ii) in $\cP MC^k_\dR(Y_1;\R)$, and $f^*$ is well defined.

From \eq{kh5eq24} and \eq{kh5eq25} we see that $f^*\ci\d=\d\ci f^*:\cP MC_\dR^k(Y_2;\R)\ra\cP MC_\dR^{k+1}(Y_1;\R)$. As $Y_2\t_{\id_{Y_2},Y_2,f}Y_1\cong Y_1$, we also see that~$f^*(\Id_{Y_2})=\Id_{Y_1}$.
 
If $g:Y_2\ra Y_3$ is another smooth map of manifolds then as
\begin{equation*}
(V\t_{t,Y_3,g}Y_2)\t_{\pi_{Y_2},Y_2,f}Y_1\cong V\t_{t,Y_3,g\ci f}Y_1,
\end{equation*}
we see from \eq{kh5eq25} that $(g\ci f)^*=f^*\ci g^*:\cP MC_\dR^k(Y_3;\R)\ra\cP MC_\dR^k(Y_1;\R)$. Also $\id_Y^*$ is the identity. Thus, pullbacks $f^*$ are contravariantly functorial.
\label{kh5def8}
\end{dfn}

The material in \S\ref{kh42} from  Proposition \ref{kh4prop5} down to Lemma \ref{kh4lem4}, and the proofs of Propositions \ref{kh4prop5} and \ref{kh4prop6} in \S\ref{kh76}--\S\ref{kh77}, all require only very straightforward modifications for the de Rham case, and we leave the details to the reader. For instance, in equation \eq{kh7eq73} in \S\ref{kh76} we should instead write
\begin{align*}
&\Pi^{T\cup U}_{T,f-}:\cP MC^k_\dR(T;\R)\longra \cP MC^k_\dR(T\cup U;\R),\\
&\Pi^{T\cup U}_{T,f-}:\bigl[V,n,(s_1,\ldots,s_n),t,\om\bigr]\longmapsto \bigl[V\t(-\iy,0],n+1,\\
&\qquad (s_1\ci\pi_V,\ldots,s_n\ci\pi_V,f\ci\pi_V+\pi_{(-\iy,0]}),t\ci\pi_V,\pi_V^*(\om)\bigr].
\end{align*}

Thus, as in Definition \ref{kh4def6} we define a soft, strong presheaf $\cP\MC^k_\dR(Y;\R)$ of $\R$-vector spaces on $Y$, with $\cP\MC^k_\dR(Y;\R)(U)=\cP MC^k_\dR(U;\R)$ for all open $U\subseteq Y$. We write $\MC^k_\dR(Y;\R)$ for the sheafification of $\cP\MC^k_\dR(Y;\R)$, a soft sheaf of $\R$-vector spaces on $Y$, and we define the {\it de Rham M-cochains\/} $MC^k_\dR(Y;\R)=\MC^k_\dR(Y;\R)(Y)$ to be its global sections. Then $\MC^k_\dR(Y;\R)(U)\!=\!MC^k_\dR(U;\R)$ for all open $U\subseteq Y$, and as in \eq{kh4eq23} we have a canonical isomorphism
\begin{equation*}
MC^k_\dR(Y;\R)\cong\mathop{\underleftarrow{\lim}\,}\nolimits_{\text{$U:U\subseteq Y$ open, $\bar U$ is compact}}\cP MC^k_\dR(U;\R).
\end{equation*}

We write $\Pi:\cP MC^k_\dR(Y;\R)\ra MC^k_\dR(Y;\R)$ for the projection coming from sheafification, and we use the same notation for elements of $\cP MC^k_\dR(Y;\R)$, such as generators $[V,n,s,t,\om]$, and for their images under $\Pi$. So we have the {\it identity cocycle\/} $\Id_Y\in MC^0_\dR(Y;\R)$, the image of $\Id_Y\in\cP MC^0_\dR(Y;\R)$. Equations \eq{kh5eq10}, \eq{kh5eq11} and \eq{kh5eq23} hold for generators in~$MC^k_\dR(Y;\R)$.

Differentials $\d$ on $\cP MC^*_\dR(Y;\R)$ descend to differentials $\d:MC^k_\dR(Y;\R)\ra MC^{k+1}_\dR(Y;\R)$ with $\d\ci\d=0$, given on generators by \eq{kh5eq24}. We define the {\it de Rham M-cohomology groups\/} $MH^*_\dR(Y;\R)$ to be the cohomology of the cochain complex $\bigl(MC^*_\dR(Y;\R),\d\bigr)$. Then $\Id_Y$ in $MC_\dR^0(Y;\R)$ defines the {\it identity cohomology class\/}~$[\Id_Y]\in MH^0_\dR(Y;\R)$.

If $f:Y_1\ra Y_2$ is a smooth map of manifolds, then pullbacks $f^*$ on M-precochains descend to {\it pullbacks\/} $f^*:MC_\dR^k(Y_2;\R)\ra MC_\dR^k(Y_1;\R)$, given on generators by \eq{kh5eq25}. These satisfy $f^*\ci\d=\d\ci f^*:MC^k_\dR(Y_2;\R)\ra MC^{k+1}_\dR(Y_1;\R)$, and so induce {\it pullbacks\/} $f^*:MH^k_\dR(Y_2;\R)\ra MH^k_\dR(Y_1;\R)$ on de Rham M-cohomology. Both pullbacks are contravariantly functorial. As in Definition \ref{kh4def7} we define {\it relative de Rham M-cochains\/} $\bigl(MC^*_\dR(Y,Z;\R),\d\bigr)$ and {\it relative de Rham M-cohomology\/} $MH^*_\dR(Y,Z;\R)$, for $Z\subseteq Y$ open.

Continuing in \S\ref{kh42} from Lemma \ref{kh4lem4}, note that the de Rham analogue of \eq{kh4eq33} is not obvious. The isomorphism $MC_\dR^{-k}(*;\R)\cong\cP MC_\dR^{-k}(*;\R)$ holds, as presheaves on a point $*$ are sheaves, so sheafifying does not change anything.
For the isomorphism $\cP MC^{-k}_\dR(*;\R)\cong MC_k^\dR(*;\R)$, define a morphism
\e
\begin{split}
&\Pd:\cP MC^{-k}_\dR(*;\R)\longra MC_k^\dR(*;\R)\qquad\text{by} \\
&\Pd:[V,n,s,t,\om]\longmapsto (-1)^{\deg\om(\deg\om-1)/2}\,[V,n,s,t,\om]
\end{split}
\label{kh5eq27}
\e
on generators $[V,n,s,t,\om]$ in $\cP MC^{-k}_\dR(*;\R)$. Equation \eq{kh5eq27} is actually the Poincar\'e duality morphism for $Y=*$, as in \S\ref{kh47}, and may be written as $\Pd:\al\mapsto[\al]\cap[*]$ for the cap product $\cap$ defined in \eq{kh5eq50} below.

We claim that $\Pd$ in \eq{kh5eq27} is well defined, and an isomorphism. Up to sign, it gives a 1-1 correspondence between generators of $\cP MC^{-k}_\dR(*;\R)$ and generators of $MC_k^\dR(*;\R)$. Also $\Pd$ maps relations (i),(iii) in $\cP MC^{-k}_\dR(*;\R)$ to relations (i),(iii) in $MC_k^\dR(*;\R)$. But there is apparently a difference between the image under $\Pd$ of Definition \ref{kh5def7}(ii) for $\cP MC^{-k}_\dR(*;\R)$, and
Definition \ref{kh5def4}(ii) for~$MC_k^\dR(*;\R)$.

In fact relation (ii) in $\cP MC^{-k}_\dR(*;\R)$ and $MC_k^\dR(*;\R)$ are equivalent under $\Pd$. For Definition \ref{kh5def7}(ii)$(*)$ with $V'=*$ and $t'=\id_*$ reduces to Definition \ref{kh5def4}(ii)$(*)$, so that Definition \ref{kh5def7}(ii)$(*)$ implies Definition \ref{kh5def4}(ii)$(*)$. But also, when $Y=*$, in the situation of Definition \ref{kh5def7}(ii)$(*)$ we have
\ea
&\sum_{i\in I}(-1)^{l\deg\om_i+\deg\om_i(\deg\om_i-1)/2}a_i\int_{V_i\t V'}(s_i\ci\pi_{V_i},\pi_{V'})^*(\eta)\w\pi_{V_i}^*(\om_i)
\nonumber\\
&=\sum_{i\in I}(-1)^{(n-k)\deg\om_i+\deg\om_i(\deg\om_i-1)/2}a_i\int_{V_i\t V'}\pi_{V_i}^*(\om_i)\w(s_i\ci\pi_{V_i},\pi_{V'})^*(\eta)
\nonumber\\
&=\sum_{i\in I}(-1)^{(n-k)\deg\om_i+\deg\om_i(\deg\om_i-1)/2}a_i\int_{V_i}\om_i\w(\pi_{V_i})_*\ci(s_i\ci\pi_{V_i},\pi_{V'})^*(\eta)
\nonumber\\
&=\sum_{i\in I}(-1)^{\deg\om_i(\deg\om_i-1)/2}a_i\int_{V_i}(\pi_{V_i})_*\ci(s_i\ci\pi_{V_i},\pi_{V'})^*(\eta)\w\om_i,
\label{kh5eq28}
\ea
using \eq{kh3eq12} in the third step. Since the last line of \eq{kh5eq28} is \eq{kh5eq9} with $(\pi_{V_i})_*\ci(s_i\ci\pi_{V_i},\pi_{V'})^*(\eta)$ in place of $\eta$ applied to $\Pd\bigl(\sum_{i\in I}a_i\,[V_i,n,s_i,t_i,\om_i]\bigr)$, Definition \ref{kh5def4}(ii)$(*)$ implies Definition \ref{kh5def7}(ii)$(*)$ under~$\Pd$.

Thus $\Pd$ is a well defined isomorphism, and the de Rham analogue of \eq{kh4eq33} holds. Also $\pd\ci\Pd=\Pd\ci\d$ by equations \eq{kh5eq13}, \eq{kh5eq24} and \eq{kh5eq27}, so we have canonical isomorphisms $MH_k^\dR(*;\R)\cong MH^{-k}_\dR(*;\R)$ for $k\in\Z$, and Axiom \ref{kh2ax2}(vii) for de Rham M-cohomology follows from Theorem \ref{kh5thm4}. As for Theorem \ref{kh4thm4}, from Theorem \ref{kh2thm2} we deduce:

\begin{thm} De Rham M-cohomology is a cohomology theory of manifolds. There are canonical isomorphisms $MH^k_\dR(Y;\R)\!\cong\! H^k(Y;\R),$ $MH^k_\dR(Y,Z;\R)\!\cong\! H^k(Y,Z;\R)$ for all\/ $Y,Z,k,$ preserving the $f^*,\d$ and isomorphisms $MH^0_\dR(*;\R)\ab\cong\R\cong H^0(*;\R),$ where $H^*(-;\R)$ is any other cohomology theory of manifolds over $\R,$ such as singular or de Rham cohomology\/~$H^*_\rsi(Y;\R),H^*_\dR(Y;\R)$.
\label{kh5thm6}
\end{thm}

As in Definition \ref{kh4def8} we define a sheaf morphism $i_Y:\R_Y\ra\MC^0_\dR(Y;\R)$ on each manifold $Y$ by \eq{kh4eq35}, where now the $\Id_{U_a}\in MC^0_\dR(U_a;\R)$ are as above. Then as in Theorem \ref{kh4thm5} we may show that 
\begin{equation*}
\xymatrix@C=20pt{ 0 \ar[r] & \R_Y \ar[r]^(0.3){i_Y} & \MC^0_\dR(Y;\R) \ar[r]^(0.48)\d & \MC^1_\dR(Y;\R) \ar[r]^(0.48)\d & \MC^2_\dR(Y;\R) \ar[r]^(0.62)\d & \cdots }
\end{equation*}
is an exact sequence of sheaves on $Y$, so that $\MC^\bu_\dR(Y;\R)=\bigl(\MC^*_\dR(Y;\R),\d\bigr)$ is a soft resolution of~$\R_Y$.

Here is the analogue of Example~\ref{kh5ex1}.

\begin{ex} Let $Y$ be a manifold and $k\in\Z$. Define $\R$-linear maps
\e
\begin{split}
F_\Mc^\dRMc&:\cP MC^k(Y;\R)\longra\cP MC_\dR^k(Y;\R),\\ 
F_\QMc^\dRMc&:\cP MC^k_\Q(Y;\R)\longra\cP MC_\dR^k(Y;\R) \qquad\text{by}\\
F_\Mc^\dRMc,F_\QMc^\dRMc&:[V,n,s,t]\longmapsto [V,n,s,t,1_V]\quad \text{on generators $[V,n,s,t]$.}
\end{split}
\label{kh5eq29}
\e

We prove these are well defined following the proof of Proposition \ref{kh5prop3} in \S\ref{kh82}, but with one extra step. The issue is to show that $F_\Mc^\dRMc,F_\QMc^\dRMc$ take relation Definition \ref{kh4def4}(ii) in $\cP MC^k_?(Y;\R)$ to Definition \ref{kh5def7}(ii) in $\cP MC_\dR^k(Y;\R)$. Now Definition \ref{kh5def7}(ii) involves a smooth map $t':V'\ra Y$, which does not occur in \S\ref{kh82}. The proof of Proposition \ref{kh4prop4} in \S\ref{kh75} shows that $t^{\prime *}:\cP MC^k_?(Y;R)\ra\cP MC^k_?(V';R)$ maps relation (ii) in $\cP MC^k_?(Y;\R)$ to (ii) in $\cP MC^k_?(V';\R)$. We then use the argument of \S\ref{kh82} to show that this pullback of (ii) in $\cP MC^k_?(V';\R)$ implies Definition \ref{kh5def7}(ii)$(*)$ for this $t':V'\ra Y$, so that $F_\Mc^\dRMc,F_\QMc^\dRMc$ do map relation (ii) to relation (ii), and are well defined.

Comparing \eq{kh4eq21} and \eq{kh5eq25}, we see that $F_\Mc^\dRMc,F_\QMc^\dRMc$ commute with pullbacks $f^*$ on $\cP MC^k(-;\R)$ and $\cP MC_\dR^k(-;\R)$. Therefore they induce morphisms
\begin{align*}
F_\Mc^\dRMc&:\cP\MC^k(Y;\R)\longra\cP\MC_\dR^k(Y;\R),\\ 
F_\QMc^\dRMc&:\cP\MC^k_\Q(Y;\R)\longra\cP\MC_\dR^k(Y;\R)
\end{align*}
of presheaves on $Y$, so sheafifying gives morphisms
\begin{align*}
F_\Mc^\dRMc&:\MC^k(Y;\R)\longra\MC_\dR^k(Y;\R),\\ 
F_\QMc^\dRMc&:\MC^k_\Q(Y;\R)\longra\MC_\dR^k(Y;\R)
\end{align*}
of sheaves on $Y$, and taking global sections gives $\R$-linear maps
\e
\begin{split}
F_\Mc^\dRMc&:MC^k(Y;\R)\longra MC_\dR^k(Y;\R),\\ 
F_\QMc^\dRMc&:MC^k_\Q(Y;\R)\longra MC_\dR^k(Y;\R).
\end{split}
\label{kh5eq30}
\e

Comparing \eq{kh4eq20}, \eq{kh5eq24} and \eq{kh5eq29} and noting that $\deg 1_V=\d 1_V=0$ we see that $F_\Mc^\dRMc\ci\d=\d\ci F_\Mc^\dRMc:\cP MC^k(Y;\R)\ra\cP MC_\dR^{k+1}(Y;\R)$, and similarly for $F_\QMc^\dRMc$. These descend to the sheafifications, so \eq{kh5eq30} induces morphisms
\e
\begin{split}
F_\Mc^\dRMc&:MH^k(Y;\R)\longra MH_\dR^k(Y;\R),\\ 
F_\QMc^\dRMc&:MH^k_\Q(Y;\R)\longra MH_\dR^k(Y;\R).
\end{split}
\label{kh5eq31}
\e
These $F_\Mc^\dRMc,F_\QMc^\dRMc$ in \eq{kh5eq31} extend to relative cohomology, and commute with pullbacks $f^*$ and continuation maps. Also $F_\Mc^\dRMc:MH^0(*;\R)\ra MH_\dR^0(*;\R)$, $F_\QMc^\dRMc:MH^0_\Q(*;\R)\ra MH_\dR^0(*;\R)$ are compatible with the isomorphisms $MH^0(*;\R)\cong MH^0_\Q(*;\R)\cong MH^0_\dR(*;\R)\cong\R$. Thus by Theorem \ref{kh2thm2}, \eq{kh5eq31} are the canonical isomorphisms from Theorem~\ref{kh5thm6}.
\label{kh5ex2}
\end{ex}

\begin{ex} Let $Y$ be a manifold of dimension $m$, and let the de Rham cochains $\bigl(C^*_\dR(Y;\R),\d\bigr)$ and cohomology $H^*_\dR(Y;\R)$ be as in Example \ref{kh2ex7}. For $k=0,1,\ldots,$ define
\e
\begin{split}
F_\dR^\dRMc&:C^k_\dR(Y;\R)\longra MC^k_\dR(Y;\R)\qquad \text{by} \\
F_\dR^\dRMc&:\om\longmapsto [Y,0,0,\id_Y,\om].
\end{split}
\label{kh5eq32}
\e
Here $t=\id_Y:Y\ra Y$ has the coorientation $c_{\id_Y}$ from Assumption \ref{kh3ass6}(e), and $\om\in C^k_\dR(Y;\R)=\Om^k(Y)$ is a $k$-form on $V=Y$, and $(0,\id_Y)\vert_{\supp\om}:\supp\om\ra\R^0\t Y$ is automatically proper, so that $[Y,0,0,\id_Y,\om]$ is a generator of $\cP MC^k_\dR(Y;\R)$ (and hence $MC^k_\dR(Y;\R)$) by Definition \ref{kh5def7}. Thus $F_\dR^\dRMc$ is well defined. Equations \eq{kh5eq10} and \eq{kh5eq24} imply that $F_\dR^\dRMc$ is $\R$-linear, with $F_\dR^\dRMc\ci\d=\d\ci F_\dR^\dRMc:C^k_\dR(Y;\R)\ra MC^{k+1}_\dR(Y;\R)$, so $F_\dR^\dRMc$ induces morphisms
\e
F_\dR^\dRMc:H^k_\dR(Y;\R)\longra MH^k_\dR(Y;\R)
\label{kh5eq33}
\e
on cohomology. These extend to relative cohomology, and commute with pullbacks $f^*$ and continuation maps. Also $F_\dR^\dRMc:H^0_\dR(*;\R)\ra MH_\dR^0(*;\R)$ is compatible with the isomorphisms $H^0_\dR(*;\R)\cong\R\cong MH^0_\dR(*;\R)$. Thus by Theorem \ref{kh2thm2}, \eq{kh5eq33} are the canonical isomorphisms from Theorem~\ref{kh5thm6}.
\label{kh5ex3}
\end{ex}

Examples \ref{kh5ex2} and \ref{kh5ex3} show that we can think of de Rham M-cohomology $MH_\dR^*(-;\R)$ as a common generalization of rational M-cohomology $MH_\Q^*(-;\R)$ and de Rham cohomology $H_\dR^*(-;\R)$.

\subsubsection{Compactly-supported de Rham M-cohomology $MH^*_{\cs,\dR}(Y;\R)$}
\label{kh523}

Section \ref{kh43} extends to de Rham M-cohomology in a straightforward way. As in Definition \ref{kh4def9}, for $Y$ a manifold and $k\in\Z$, we define the {\it compactly-supported de Rham M-cochains\/} $MC^k_{\cs,\dR}(Y;\R)$ to be the $\R$-vector subspace $MC^k_{\cs,\dR}(Y;\R)\!\subseteq\! MC^k_\dR(Y;\R)$ of compactly-supported sections of $\MC^k_\dR(Y;\R)$.

Then $MC^k_{\cs,\dR}(Y;\R)$ is the global sections of a flabby cosheaf $\uMC^k_{\cs,\dR}(Y;\R)$ of $\R$-vector spaces associated to the c-soft sheaf $\MC^k_\dR(Y;\R)$ by Theorem \ref{kh2thm3}. The morphisms $\d:MC^k_\dR(Y;\R)\ra MC^{k+1}_\dR(Y;\R)$ from \S\ref{kh522} restrict to $\d:MC^k_{\cs,\dR}(Y;\R)\ra MC^{k+1}_{\cs,\dR}(Y;\R)$ with $\d\ci\d=0$. Hence $\bigl(MC_{\cs,\dR}^*(Y;\R),\d\bigr)$ is a cochain complex. Define the {\it compactly-supported de Rham M-cohomology groups\/} $MH^*_{\cs,\dR}(Y;\R)$ to be the cohomology of this cochain complex. 

Suppose $f:Y_1\ra Y_2$ is a proper smooth map of manifolds. Then the pullback $f^*:MC^k_\dR(Y_2;\R)\ra MC^k_\dR(Y_1;\R)$ from \S\ref{kh522} restricts to the {\it pullback\/} $f^*:MC^k_{\cs,\dR}(Y_2;\R)\ra MC^k_{\cs,\dR}(Y_1;\R)$. These satisfy $\d\ci f^*=f^*\ci\d$, and so induce $f^*:MH^k_{\cs,\dR}(Y_2;\R)\ra MH^k_{\cs,\dR}(Y_1;\R)$, as in Property~\ref{kh2pr1}(b).

If $U\subseteq Y$ is open and $i:U\hookra Y$ is the inclusion, then as for the morphism $\si_{UY}:\uMC^k_{\cs,\dR}(Y;\R)(U)\ra\uMC^k_{\cs,\dR}(Y;\R)(Y)$ defined in Theorem \ref{kh2thm3}(a), there is an injective {\it pushforward\/} $i_*:MC^k_{\cs,\dR}(U;\R)\ra MC^k_{\cs,\dR}(Y;\R)$, such that if $\al\in MC^k_{\cs,\dR}(U;\R)$ then $i_*(\al)\in MC^k_{\cs,\dR}(Y;\R)$ is the unique element with $i_*(\al)\vert_U=\al$ and $i_*(\al)\vert_{Y\sm\supp\al}=0$. These satisfy $\d\ci i_*=i_*\ci\d$, and so induce $i_*:MH^k_{\cs,\dR}(U;\R)\ra MH^k_{\cs,\dR}(Y;\R)$, as in Property \ref{kh2pr1}(c). Since $\MC^\bu_\dR(Y;\R)$ is a soft resolution of $\R_Y$ by \S\ref{kh522}, as in Theorem \ref{kh4thm6} we deduce:

\begin{thm} Compactly-supported de Rham M-cohomology is a compactly-supported cohomology theory of manifolds. As in \eq{kh4eq42} there are canonical isomorphisms $MH^k_{\cs,\dR}(Y;\R)\cong H^k_\cs(Y;\R)$ for all\/ $Y,k,$ preserving the data $\Pi,f^*,i_*$ in Property\/ {\rm\ref{kh2pr1}(a)--(c)} and the isomorphisms $MH^0_{\cs,\dR}(*;\R)\cong\R\cong H^0_\cs(*;\R),$ where $H^*_\cs(-;\R)$ is any other compactly-supported cohomology theory of manifolds over $\R,$ such as compactly-supported singular cohomology\/ $H^*_{\cs,\rsi}(Y;\R)$ in Example\/ {\rm\ref{kh2ex8},} or compactly-supported de Rham cohomology\/ $H^*_{\cs,\dR}(Y;\R)$ in Example\/~{\rm\ref{kh2ex10}}.
\label{kh5thm7}
\end{thm}

\begin{ex}{\bf(a)} The morphisms $F_\Mc^\dRMc:MC^k(Y;\R)\ra MC_\dR^k(Y;\R)$ from Example \ref{kh5ex2} for $k\in\Z$ restrict to morphisms
\e
F_{\cs,\Mc}^{\cs,\dRMc}:MC^k_\cs(Y;\R)\longra MC_{\cs,\dR}^k(Y;\R).
\label{kh5eq34}
\e
These satisfy $\d\ci F_{\cs,\Mc}^{\cs,\dRMc}=F_{\cs,\Mc}^{\cs,\dRMc}\ci\d$, and so induce morphisms 
\e
F_{\cs,\Mc}^{\cs,\dRMc}:MH^k_\cs(Y;\R)\longra MH_{\cs,\dR}^k(Y;\R).
\label{kh5eq35}
\e
Since \eq{kh5eq34} comes from a morphism $F_\Mc^\dRMc:\MC^\bu(Y;\R)\ra \MC^\bu_\dR(Y;\R)$ of soft resolutions of $\R_Y$ acting as the identity on $\R_Y$, facts about sheaf cohomology imply that \eq{kh5eq35} are the canonical isomorphisms from Theorem~\ref{kh5thm7}.
\smallskip

\noindent{\bf(b)} The morphisms $F_\dR^\dRMc:C^k_\dR(Y;\R)\ra MC^k_\dR(Y;\R)$ from Example \ref{kh5ex3} for $k\in$ restrict to morphisms
\begin{equation*}
F_{\cs,\dR}^{\cs,\dRMc}:C^k_{\cs,\dR}(Y;\R)\longra MC^k_{\cs,\dR}(Y;\R).
\end{equation*}
These induce morphisms 
\begin{equation*}
F_{\cs,\dR}^{\cs,\dRMc}:H^k_{\cs,\dR}(Y;\R)\longra MH^k_{\cs,\dR}(Y;\R),
\end{equation*}
which are the canonical isomorphisms from Theorem~\ref{kh5thm7}.
\label{kh5ex4}
\end{ex}

Let $Y$ be a manifold, $k\in\Z$, and $[V,n,s,t,\om]\in MC^k_\dR(Y;\R)$ be a generator of $MC^k_\dR(Y;\R)$, in the sense of \S\ref{kh522}. We say that $[V,n,s,t,\om]$ is a {\it compact generator\/} if $s\vert_{\supp\om}:\supp\om\ra\R^n$ is proper over an open neighbourhood of 0 in $\R^n$. This implies that $(s,t)\vert_{\supp\om}:\supp\om\ra\R^n\t Y$ is proper over an open neighbourhood of $\{0\}\t Y$ in $\R^n\t Y$, as assumed in Definition \ref{kh5def7}. By a similar proof to that of Proposition \ref{kh4prop8} in \S\ref{kh78}, we can show:

\begin{prop} Let\/ $Y$ be a manifold and\/ $k\in\Z$. Then as an $\R$-vector space, $MC^k_{\cs,\dR}(Y;\R)$ is generated by compact generators $[V,n,s,t,\om],$ subject only to relations Definition\/ {\rm\ref{kh5def7}(i)--(iii)} applied to compact generators.
\label{kh5prop4}
\end{prop}

\subsubsection{Locally finite de Rham M-homology $MH_*^{\lf,\dR}(Y;\R)$}
\label{kh524}

Here is the de Rham analogue of Definition \ref{kh4def11}. It follows Definition \ref{kh5def4}, except for the properness conditions, which match those in Definition~\ref{kh5def7}.

\begin{dfn} Let $Y$ be a manifold. Consider quintuples $(V,n,s,t,\om)$, where $V$ is an oriented manifold with corners (i.e.\ a pair $(V,o_V)$ with $V$ an object in $\tManc$ and $o_V$ an orientation on $V$, usually left implicit), and $n=0,1,\ldots,$ and $s:V\ra\R^n$ is a smooth map (morphism in $\tManc$), and $t:V\ra Y$ is a smooth map, and $\om\in\Om^p(V)$ is a $p$-form on $V$ for $p=0,1,\ldots,$ as in Assumption \ref{kh3ass8}, such that $(s,t)\vert_{\supp\om}:\supp\om\ra\R^n\t Y$ is proper over an open neighbourhood of $\{0\}\t Y$ in~$\R^n\t Y$.

Define an equivalence relation $\sim$ on such quintuples by $(V,n,s,t,\om)\sim(V',\ab n',\ab s',\ab t',\ab\om')$ if $n=n'$, and there exists an orientation-preserving diffeomorphism $f:V\ra V'$ with $s=s'\ci f$ and $t=t'\ci f$ and $f^*(\om')=\om$. Write $[V,n,s,t,\om]$ for the $\sim$-equivalence class of $(V,n,s,t,\om)$. We call $[V,n,s,t,\om]$ a {\it generator}.

For $k\in\Z$, define the {\it locally finite de Rham M-prechains\/} $\cP MC^{\lf,\dR}_k(Y;\R)$ to be the $\R$-vector space generated by such $[V,n,s,t,\om]$ with $\dim V=n+k+\deg\om$, subject to the relations:
\begin{itemize}
\setlength{\itemsep}{0pt}
\setlength{\parsep}{0pt}
\item[(i)] For each generator $[V,n,s,t,\om]$ and each $i=0,\ldots,n$ we have
\begin{equation*}
[V,n,s,t,\om]=(-1)^{n-i}[V\t\R,n+1,s',t\ci\pi_V,\pi_V^*(\om)]
\end{equation*}
in $\cP MC^{\lf,\dR}_k(Y;\R)$, where writing $s=(s_1,\ldots,s_n):V\ra\R^n$ with $s_j:V\ra\R$ and $\pi_V:V\t\R\ra V$, $\pi_\R:V\t\R\ra\R$ for the projections, then 
\begin{equation*}
s'=(s_1\ci\pi_V,\ldots,s_i\ci\pi_V,\pi_\R,s_{i+1}\ci\pi_V,\ldots,s_n\ci\pi_V):V\t\R\longra\R^{n+1},
\end{equation*}
and $V\t\R$ has the product orientation from Assumption \ref{kh3ass6}(f) of the given orientation on $V$ and the standard orientation on~$\R$.
\item[(ii)] Let $I$ be a finite indexing set, $a_i\in \R$ for $i\in I$, and $[V_i,n,s_i,t_i,\om_i]$, $i\in I$ be generators for $\cP MC^{\lf,\dR}_k(Y;\R)$, all with the same $n$. Suppose there exists an open neighbourhood $X$ of $\{0\}\t Y$ in $\R^n\t Y$, such that $(s_i,t_i)\vert_{\supp\om_i}:\supp\om_i\ra\R^n\t Y$ is proper over $X$ for all $i\in I$, and the following condition holds:
\begin{itemize}
\setlength{\itemsep}{0pt}
\setlength{\parsep}{0pt}
\item[$(*)$] Suppose $\eta\in C^\iy\bigl(\La^{n+k}T^*(\R^n\t Y)\bigr)$ is an $(n+k)$-form on $\R^n\t Y$ with $\supp\eta$ a compact subset of $X\subseteq\R^n\t Y$. Then
\e
\sum_{i\in I}a_i\int_{V_i}(s_i,t_i)^*(\eta)\w\om_i=0
\quad\text{in $\R$.}
\label{kh5eq36}
\e
Here $(s_i,t_i)^*(\eta)\w\om_i$ is a form of degree $n+k+\deg\om_i=\dim V_i$ on $V_i$, with $\supp[(s_i,t_i)^*(\eta)\w\om_i]\subseteq(s_i,t_i)^{-1}[\supp\eta]\cap\supp\om_i$ by Assumption \ref{kh3ass8}(c). Since $(s_i,t_i)\vert_{\supp\om_i}:\supp\om_i\ra\R^n\t Y$ is proper over $X$ and $\supp\eta\subseteq X$ is compact, we see that $(s_i,t_i)^*(\eta)\w\om_i$ is compactly-supported, so $\int_{V_i}(s_i,t_i)^*(\eta)\w\om_i$ is defined in $\R$ by Assumption \ref{kh3ass8}(e), and \eq{kh5eq36} makes sense.
\end{itemize}
Then
\begin{equation*}
\sum_{i\in I}a_i\,[V_i,n,s_i,t_i,\om_i]=0\qquad\text{in $\cP MC^{\lf,\dR}_k(Y;\R)$.}
\end{equation*}
\item[(iii)] For each generator $[V,n,s,t,\om]$ with $s=(s_1,\ldots,s_n):V\ra\R^n$ and each permutation $\si\in S_n$ of $1,2,\ldots,n$ we have
\begin{equation*}
\bigl[V,n,(s_1,s_2,\ldots,s_n),t,\om\bigr]=\sign(\si)\cdot \bigl[V,n,(s_{\si(1)},s_{\si(2)},\ldots,s_{\si(n)}),t,\om\bigr]
\end{equation*}
in $\cP MC^{\lf,\dR}_k(Y;\R)$, where $\sign:S_n\ra\{\pm 1\}$ is the usual group morphism.
\end{itemize}

Define $\pd:\cP MC_k^{\lf,\dR}(Y;\R)\ra\cP MC_{k-1}^{\lf,\dR}(Y;\R)$ to be the unique $\R$-linear map satisfying \eq{kh5eq13} on generators $[V,n,s,t,\om]$. As in Definition \ref{kh5def5}, it is well-defined, with~$\pd\ci\pd=0:\cP MC_k^{\lf,\dR}(Y;\R)\ra\cP MC_{k-2}^{\lf,\dR}(Y;\R)$.

Let $f:Y_1\ra Y_2$ be a proper smooth map of manifolds. Define $f_*:\cP MC_k^{\lf,\dR}(Y_1;\R)\ra\cP MC_k^{\lf,\dR}(Y_2;\R)$ to be the unique $\R$-linear morphism acting on generators $[V,n,s,t,\om]$ by \eq{kh5eq19}. We need $f$ proper so that $(s,t)\vert_{\supp\om}:\supp\om\ra\R^n\t Y_1$ proper near $\{0\}\t Y_1$ in $\R^n\t Y_1$ implies that $(s,f\ci t)\vert_{\supp\om}:\supp\om\ra\R^n\t Y_2$ is proper near $\{0\}\t Y_2$ in $\R^n\t Y_2$. Then $f_*$ is well-defined as in Definition \ref{kh5def6}, with~$\pd\ci f_*=f_*\ci\pd:\cP MC_k^{\lf,\dR}(Y_1;\R)\ra\cP MC_{k-1}^{\lf,\dR}(Y_2;\R)$.

Let $Y$ be a manifold and $U\subseteq Y$ an open set, with $i:U\hookra Y$ the inclusion. Define $i^*:\cP MC_k^{\lf,\dR}(Y;\R)\ra\cP MC_k^{\lf,\dR}(U;\R)$ to be the unique $\R$-linear morphism acting on generators $[V,n,s,t,\om]$ by
\e
\begin{split}
i^*:[V,n,s,t,\om]\longmapsto &\,[V',n,s',t',\om']\\
:=&\,\bigl[t^{-1}(U),n,s\vert_{t^{-1}(U)},t\vert_{t^{-1}(U)},\om\vert_{t^{-1}(U)}\bigr].
\end{split}
\label{kh5eq37}
\e
Then $(s,t)\vert_{\supp\om}:\supp\om\ra\R^n\t Y$ proper near $\{0\}\t Y$ in $\R^n\t Y$ implies that $(s',t')\vert_{\supp\om'}:\supp\om'\ra\R^n\t U$ is proper near $\{0\}\t U$ in $\R^n\t U$, so the r.h.s.\ of \eq{kh5eq37} is a generator of $\cP MC_k^{\lf,\dR}(U;\R)$. Clearly $i^*$ takes relations (i)--(iii) in $\cP MC_k^{\lf,\dR}(Y;\R)$ to relations (i)--(iii) in $\cP MC_k^{\lf,\dR}(U;\R)$, so is well-defined. Also $\pd\ci i^*=i^*\ci\pd:\cP MC_k^{\lf,\dR}(Y;\R)\ra\cP MC_{k-1}^{\lf,\dR}(U;\R)$, and if $j:U'\hookra U$ is another open inclusion then~$j^*\ci i^*=(i\ci j)^*:\cP MC_k^{\lf,\dR}(Y;\R)\ra\cP MC_k^{\lf,\dR}(U';\R)$.

If $Y$ is oriented with $\dim Y=m$, we define the {\it fundamental prechain\/} $[Y]=[Y,0,0,\id_Y,1_Y]$ in $\cP MC_m^{\lf,\dR}(Y;\R)$.
\label{kh5def9}
\end{dfn}

All the rest of \S\ref{kh44}, from Remark \ref{kh4rem6} to Remark \ref{kh4rem8}, extends to the de Rham case in a straightforward way, and we leave the details as an exercise. As in Definition \ref{kh4def12} we define a soft, strong presheaf $\cP\MC^{\lf,\dR}_k(Y;\R)$ on $Y$ with $\cP\MC^{\lf,\dR}_k(Y;\R)(U)=\cP MC^{\lf,\dR}_k(U;\R)$ for open $U\subseteq Y$. Write $\MC^{\lf,\dR}_k(Y;\R)$ for the sheafification of $\cP\MC^{\lf,\dR}_k(Y;\R)$, a soft sheaf of $\R$-vector spaces on $Y$, and we define the {\it locally finite de Rham M-chains\/} $MC^{\lf,\dR}_k(Y;\R)$ to be its global sections $\MC^{\lf,\dR}_k(Y;\R)(Y)$. As in \eq{kh4eq48} we have a canonical isomorphism
\begin{equation*}
MC^{\lf,\dR}_k(Y;\R)\cong\mathop{\underleftarrow{\lim}\,}\nolimits_{\text{$U:U\subseteq Y$ open, $\bar U$ is compact}}\cP MC^{\lf,\dR}_k(U;\R).
\end{equation*}

There are natural projections $\Pi:\cP MC^{\lf,\dR}_k(Y;\R)\ra MC^{\lf,\dR}_k(Y;\R)$, and we use the same notation for elements of $\cP MC_k^{\lf,\dR}(Y;\R)$, such as generators $[V,n,s,t,\om]$, and for their images under $\Pi$ in $MC_k^{\lf,\dR}(Y;\R)$. Differentials $\pd$, pushforwards $f_*$, and pullbacks $i^*$ in Definition \ref{kh5def9} descend to $MC^{\lf,\dR}_k(-;\R)$ by sheafification. We define the {\it locally finite de Rham M-homology groups\/} $MH_*^{\lf,\dR}(Y;\R)$ to be the homology of $\bigl(MC_*^{\lf,\dR}(Y;\R),\pd\bigr)$. 

Pushforwards $f_*$ and pullbacks $i^*$ descend to homology. There are injective $\R$-linear maps $\Pi:MC^{\lf,\dR}_k(Y;\R)\ra MC^{\lf,\dR}_k(Y;\R)$ mapping $\Pi:[V,n,s,t,\om]\mapsto[V,n,s,t,\om]$ on generators, which also descend to homology. Eventually, as in Theorem \ref{kh4thm7} we prove:

\begin{thm} Locally finite de Rham M-homology is a locally finite homology theory of manifolds. That is, there are canonical isomorphisms $H_k^\lf(Y;\R)\cong MH_k^{\lf,\dR}(Y;\R)$ for all\/ $Y,k,$ preserving the data $\Pi,f_*,i^*$ described in Property\/ {\rm\ref{kh2pr2}(a)--(c)} and the isomorphisms $H_0^\lf(*;\R)\cong\R\cong MH_0^{\lf,\dR}(*;\R),$ where $H_*^\lf(-;\R)$ is any other locally finite homology theory of manifolds over~$\R$. 
\label{kh5thm8}
\end{thm}

Following the method of Example \ref{kh5ex2}, we may define $\R$-linear maps
\begin{align*}
F_{\lf,\Mh}^{\lf,\dRMh}&:MC_k^\lf(Y;\R)\longra MC^{\lf,\dR}_k(Y;\R),\\ 
F_{\lf,\QMh}^{\lf,\dRMh}&:MC_k^{\lf,\Q}(Y;\R)\longra MC^{\lf,\dR}_k(Y;\R),
\end{align*}
acting on generators $[V,n,s,t]$ of $\cP MC_k^\lf(Y;\R),\cP MC_k^{\lf,\Q}(Y;\R)$ by
\begin{equation*}
F_{\lf,\Mh}^{\lf,\dRMh},F_{\lf,\QMh}^{\lf,\dRMh}:[V,n,s,t]\longmapsto [V,n,s,t,1_V].
\end{equation*}
These commute with the $\pd,f_*,i^*$, and on homology induce 
\begin{align*}
F_{\lf,\Mh}^{\lf,\dRMh}&:MH_k^\lf(Y;\R)\longra MH^{\lf,\dR}_k(Y;\R),\\ 
F_{\lf,\QMh}^{\lf,\dRMh}&:MH_k^{\lf,\Q}(Y;\R)\longra MH^{\lf,\dR}_k(Y;\R),
\end{align*}
which are the canonical isomorphisms from Theorem \ref{kh5thm8}. We may compose $F_{\lf,\Mh}^{\lf,\dRMh}$ with the morphisms $F_{\lf,\ssi}^{\lf,\Mh},\hat F{}_{\lf,\ssi}^{\lf,\Mh}$ from Examples \ref{kh4ex3} and \ref{kh4ex4} to get
\begin{align*}
F_{\lf,\ssi}^{\lf,\dRMh}:C_k^{\lf,\ssi}(Y;\R)\longra MC_k^{\lf,\dR}(Y;\R),\\
\hat F{}_{\lf,\ssi}^{\lf,\dRMh}:\hat C_k^{\lf,\ssi}(Y;\R)\longra MC_k^{\lf,\dR}(Y;\R),
\end{align*}
which also induce the canonical isomorphisms on homology.

\subsubsection{Cup, cap and cross products and Poincar\'e duality}
\label{kh525}

Here is the analogue of Definition \ref{kh4def14}:

\begin{dfn} Let $Y$ be a manifold. For $k,l\in\Z$, define $\R$-bilinear maps
\e
\cup:\cP MC_\dR^k(Y;\R)\t \cP MC_\dR^l(Y;\R)\longra \cP MC_\dR^{k+l}(Y;\R)
\label{kh5eq38}
\e
on generators $[V,n,s,t,\om]\!\in\! \cP MC_\dR^k(Y;\R)$, $[V',n',s',t',\om']\!\in\! \cP MC_\dR^l(Y;\R)$ by
\ea
&[V,n,s,t,\om]\cup[V',n',s',t',\om']
=(-1)^{nl+\deg\om(l+n')+\deg\om\deg\om'}[\ti V,\ti n,\ti s,\ti t,\ti\om]
\nonumber\\
&:=(-1)^{nl+\deg\om(l+n')+\deg\om\deg\om'}\bigl[V\t_{t,Y,t'}V',n+n',(s_1\ci\pi_V,\ldots,s_n\ci\pi_V,
\nonumber\\
&\qquad s_1'\ci\pi_{V'},\ldots,s_{n'}'\ci\pi_{V'}),t\ci\pi_V,\pi_V^*(\om)\w\pi_{V'}^*(\om')\bigr].
\label{kh5eq39}
\ea
Writing $c_t,c_{t'}$ for the coorientations on $t,t'$, the coorientation on $\ti t=t\ci\pi_V:\ti V\ra Y$ is $c_{\ti t}=c_t\ci c_{\pi_V}$, where the coorientation $c_{\pi_V}$ on $\pi_V$ is induced from $c_{t'}$, using Assumption \ref{kh3ass6}(d),(l). In Proposition \ref{kh5prop5} below we show that $\cup$ is well-defined.

Given $[V,n,s,t,\om]\in \cP MC_\dR^j(Y;\R)$, $[V',n',s',t',\om']\in \cP MC_\dR^k(Y;\R)$ and $[V'',n'',s'',t'',\om'']\!\in\! \cP MC_\dR^l(Y;\R)$, as in \eq{kh4eq62}, \eq{kh5eq8}, \eq{kh4eq64} and \eq{kh4eq66} we have
\ea
\begin{split}
&\bigl([V,n,s,t,\om]\cup[V',n',s',t',\om']\bigr)\cup[V'',n'',s'',t'',\om'']\\
&\quad =[V,n,s,t,\om]\cup\bigl([V',n',s',t',\om']\cup[V'',n'',s'',t'',\om'']\bigr),
\end{split}
\label{kh5eq40}\\
&[V,n,s,t,\om]\cup[V',n',s',t',\om']=(-1)^{jk} [V',n',s',t',\om']\cup [V,n,s,t,\om],
\label{kh5eq41}\\
\begin{split}
&\d\bigl([V,n,s,t,\om]\cup[V',n',s',t',\om']\bigr)=\bigl(\d[V,n,s,t,\om]\bigr)\cup[V',n',s',t',\om']\\
&\qquad +(-1)^j[V,n,s,t,\om]\cup\bigl(\d[V',n',s',t',\om']\bigr),
\end{split}
\label{kh5eq42}\\
&\Id_Y\cup[V,n,s,t,\om]=[V,n,s,t,\om]\cup\Id_Y=[V,n,s,t,\om],
\label{kh5eq43}
\ea
using Definition \ref{kh5def7}(iii) to prove \eq{kh5eq41}. Equations \eq{kh5eq40}--\eq{kh5eq43} imply that for all $\al\in\cP MC_\dR^j(Y;\R)$, $\be\in\cP MC_\dR^k(Y;\R)$ and $\ga\in\cP MC_\dR^l(Y;\R)$ we have
\ea
(\al\cup\be)\cup\ga&=\al\cup(\be\cup\ga),
\label{kh5eq44}\\
\al\cup\be&=(-1)^{jk}\be\cup\al,
\label{kh5eq45}\\
\d(\al\cup\be)&=(\d\al)\cup\be+(-1)^j\al\cup(\d\be),
\label{kh5eq46}\\
\Id_Y\cup\al&=\al\cup\Id_Y=\al.
\label{kh5eq47}
\ea

Suppose $f:Y_1\ra Y_2$ is a smooth map of manifolds, so that Definition \ref{kh5def8} defines the pullback $f^*:\cP MC_\dR^k(Y_2;\R)\ra\cP MC_\dR^k(Y_1;\R)$. If $[V,n,s,t,\om]\in \cP MC_\dR^k(Y_2;\R)$ and $[V',n',s',t',\om']\in \cP MC_\dR^l(Y_2;\R)$ then as in \eq{kh4eq68} we have
\begin{equation*}
f^*\bigl([V,n,s,t,\om]\cup[V',n',s',t',\om']\bigr)
=f^*\bigl([V,n,s,t,\om]\bigr)\cup f^*\bigl([V',n',s',t',\om']\bigr).
\end{equation*}
Hence for all $\al\in\cP MC_\dR^k(Y_2;\R)$ and $\be\in\cP MC_\dR^l(Y_2;\R)$ we have
\e
f^*(\al\cup\be)=f^*(\al)\cup f^*(\be).
\label{kh5eq48}
\e

\label{kh5def10}
\end{dfn}

Here is the analogue of Proposition \ref{kh4prop11}. It will be proved in~\S\ref{kh83}.

\begin{prop} The product\/ $\cup$ in \eq{kh5eq38}--\eq{kh5eq39} is well defined.
\label{kh5prop5}
\end{prop}

Here is the analogue of Definition~\ref{kh4def15}:

\begin{dfn} Let $Y$ be a manifold, and $k,l\in\Z$. As in \eq{kh4eq70}, define  
\begin{align*}
\cup_{k,l}&:\cP\MC_\dR^k(Y;\R)\ot_\R\cP\MC_\dR^l(Y;\R)\longra\cP\MC_\dR^{k+l}(Y;\R)\quad\text{by}\\
\cup_{k,l}(U)=\cup&:\cP MC_\dR^k(U;\R)\ot_\R\cP MC_\dR^l(U;\R)\longra\cP MC_\dR^{k+l}(U;\R)
\end{align*}
for all open $U\subseteq Y$. Then $\cup_{k,l}$ is a presheaf morphism, so sheafifying induces
\begin{equation*}
\cup_{k,l}:\MC_\dR^k(Y;\R)\ot_\R\MC_\dR^l(Y;\R)\longra\MC_\dR^{k+l}(Y;\R),
\end{equation*}
and taking global sections $\cup=\cup_{k,l}(Y)$ defines $\R$-bilinear {\it cup products\/}
\begin{equation*}
\cup:MC_\dR^k(Y;\R)\t MC_\dR^l(Y;\R)\longra MC_\dR^{k+l}(Y;\R),
\end{equation*}
which as in \eq{kh4eq73} restrict to
\begin{align*}
&\cup:MC_{\cs,\dR}^k(Y;\R)\t MC_\dR^l(Y;\R)\longra MC_{\cs,\dR}^{k+l}(Y;\R),\\
&\cup:MC_\dR^k(Y;\R)\t MC_{\cs,\dR}^l(Y;\R)\longra MC_{\cs,\dR}^{k+l}(Y;\R),\\
&\cup:MC_{\cs,\dR}^k(Y;\R)\t MC_{\cs,\dR}^l(Y;\R)\longra MC_{\cs,\dR}^{k+l}(Y;\R).
\end{align*}

Equations \eq{kh5eq44}--\eq{kh5eq48} extend to sheafifications, and hold on $MC^*_\dR(-;\R)$. Note in particular that $\cup$ is supercommutative on $MC_\dR^*(Y;\R)$ by \eq{kh5eq45}, as discussed in Remark \ref{kh4rem9}. Equation \eq{kh5eq46} in $MC_\dR^*(Y;\R)$ implies that $\cup$ descends to M-cohomology. Thus as in \eq{kh4eq74}, we define $\R$-bilinear morphisms
\e
\begin{split}
&\cup:MH_\dR^k(Y;\R)\t MH_\dR^l(Y;\R)\longra MH_\dR^{k+l}(Y;\R),\\
&\cup:MH_{\cs,\dR}^k(Y;\R)\t MH_\dR^l(Y;\R)\longra MH_{\cs,\dR}^{k+l}(Y;\R),\\
&\cup:MH_\dR^k(Y;\R)\t MH_{\cs,\dR}^l(Y;\R)\longra MH_{\cs,\dR}^{k+l}(Y;\R),\\
&\cup:MH_{\cs,\dR}^k(Y;\R)\t MH_{\cs,\dR}^l(Y;\R)\longra MH_{\cs,\dR}^{k+l}(Y;\R),
\end{split}
\label{kh5eq49}
\e
by $[\al]\cup[\be]=[\al\cup\be]$ for $\al\in MC_\dR^k(Y;\R)$, $\be\in MC_\dR^l(Y;\R)$ with~$\d\al=\d\be=0$.

\label{kh5def11}
\end{dfn}

Examples \ref{kh5ex2} and \ref{kh5ex3} defined morphisms of cochain complexes
\begin{align*}
F_\QMc^\dRMc&:\bigl(MC^*_\Q(Y;\R),\d\bigr)\longra \bigl(MC_\dR^k(Y;\R),\d\bigr),\\
F_\dR^\dRMc&:\bigl(C^*_\dR(Y;\R),\d\bigr)\longra \bigl(MC^*_\dR(Y;\R),\d\bigr),
\end{align*}
inducing the canonical isomorphisms on cohomology, and Example \ref{kh5ex4} discussed the compactly-supported versions. Comparing equations \eq{kh4eq61}, \eq{kh5eq29}, \eq{kh5eq32}, and \eq{kh5eq39}, we see that $F_\QMc^\dRMc,F_\dR^\dRMc$ commute with cup products, and also map identities $\Id_Y$ to identities $\Id_Y$, so they are morphisms of cdgas over $\R$. Thus as for Theorems \ref{kh4thm9}, \ref{kh4thm10} and \ref{kh5thm3}, we deduce:

\begin{thm} Under the canonical isomorphisms $MH^k_\dR(Y;\R)\cong H^k(Y;\R),$ $MH^k_{\cs,\dR}(Y;\R)\cong H^k_\cs(Y;\R)$ from Theorems\/ {\rm\ref{kh5thm6}} and\/ {\rm\ref{kh5thm7},} the cup products in \eq{kh5eq49} are identified with the usual cup products \eq{kh2eq45} on ordinary (compactly-supported) cohomology $H^*(Y;\R),H^*_\cs(Y;\R)$.
\label{kh5thm9}
\end{thm}

\begin{thm} For any manifold\/ $Y,$ the cdga $\bigl(MC^*_\dR(Y;\R),\d,\cup,\Id_Y\bigr)$ over $\R$ is equivalent in $\cdga_\R^\iy$ to the `usual' cdga over $\R$ associated to $Y$ in topology, as represented by the de Rham cdga $\bigl(C^*_\dR(Y;\R),\d,\w,1_Y\bigr)$ from Example\/~{\rm\ref{kh2ex7},} or by Sullivan's cdga $A_{PL}(Y;\R)$ of polynomial differential forms on\/~$Y$.

Therefore the real homotopy type of\/ $Y,$ and topological invariants of\/ $Y$ depending on the cdga up to equivalence, such as Massey products, may be computed using the cdga $\bigl(MC^*_\dR(Y;\R),\ab\d,\ab\cup,\ab\Id_Y\bigr),$ and will give the correct answers under the canonical isomorphism $MH^*_\dR(Y;\R)\cong H^*(Y;\R)$ from Theorem\/~{\rm\ref{kh5thm6}}.

\label{kh5thm10}
\end{thm}

As in Definition \ref{kh4def17} we define cross products on de Rham M-cochains
\begin{gather*}
\t:MC^k_\dR(Y_1;\R)\t MC^l_\dR(Y_2;\R)\longra MC^{k+l}_\dR(Y_1\t Y_2;\R)\\
\text{by}\qquad \al\t\be=\pi_{Y_1}^*(\al)\cup\pi_{Y_2}^*(\be),
\end{gather*}
and these are compatible with identities, pullbacks and differentials, restrict to compactly-supported cochains $MC^*_{\cs,\dR}(-;\R)$, and descend to cross products on cohomology $MH^*_\dR(-;\R),MH^*_{\cs,\dR}(-;\R)$, which are identified with the usual cross products \eq{kh2eq52} on ordinary (compactly-supported) cohomology $H^*(Y;\R),H^*_\cs(Y;\R)$ under the isomorphisms $MH^*_\dR(Y;\R)\cong H^*(Y;\R)$, $MH^k_{\cs,\dR}(Y;\R)\cong H^k_\cs(Y;\R)$ from Theorems \ref{kh5thm6} and \ref{kh5thm7}, as in Corollary~\ref{kh4cor2}.

The material of \S\ref{kh46} on cap products, and cross products on M-homology, extends to de Rham M-(co)homology in a very similar way to the material above, and we leave the details to the reader. The analogues of equations \eq{kh4eq81}--\eq{kh4eq82} defining cap products on a manifold $Y$ with $\dim Y=m$ are to define
\begin{equation*}
\cap:\cP MC^k_\dR(Y;\R)\t \cP MC_l^{\lf,\dR}(Y;\R)\longra \cP MC_{l-k}^{\lf,\dR}(Y;\R)
\end{equation*}
on $[V,n,s,t,\om]\in\cP MC^k_\dR(Y;\R)$, $[V',n',s',t',\om']\in\cP MC_l^{\lf,\dR}(Y;\R)$ by
\e
\begin{split}
&[V,n,s,t,\om]\cap[V',n',s',t',\om']\\
&=(-1)^{n(l+m)+\deg\om(l+n')+\deg\om(\deg\om-1)/2}[\ti V,\ti n,\ti s,\ti t,\ti\om]\\
&:=(-1)^{n(l+m)+\deg\om(l+n')+\deg\om(\deg\om-1)/2}\bigl[V\t_{t,Y,t'}V',n\!+\!n',\\
&(s_1\ci\pi_V,\ldots,s_n\ci\pi_V,s_1'\ci\pi_{V'},\ldots,s_{n'}'\ci\pi_{V'}),t\ci\pi_V,\pi_V^*(\om)\w\pi_{V'}^*(\om')\bigr].
\end{split}
\label{kh5eq50}
\e

The analogues of equations \eq{kh4eq91}--\eq{kh4eq92} defining cross products on manifolds $Y_1,Y_2$ of dimensions $m_1,m_2$ are to define
\begin{equation*}
\t:\cP MC^{\lf,\dR}_k(Y_1;\R)\t \cP MC^{\lf,\dR}_l(Y_2;\R)\longra \cP MC^{\lf,\dR}_{k+l}(Y_1\t Y_2;\R)
\end{equation*}
on $[V,n,s,t,\om]\in \cP MC^{\lf,\dR}_k(Y_1;\R)$, $[V',n',s',t',\om']\in \cP MC^{\lf,\dR}_l(Y_2;\R)$ by
\begin{align*}
[V,n&,s,t,\om]\t[V',n',s',t',\om']=(-1)^{n(l+m_2)+\deg\om(l+n')}\bigl[V\t V',n+n',\\
&(s_1\ci\pi_V,\ldots,s_n\ci\pi_V,s_1'\ci\pi_{V'},\ldots,s_{n'}'\ci\pi_{V'}),t\t t',\pi_V^*(\om)\w\pi_{V'}^*(\om')\bigr].
\end{align*}
The de Rham analogues of Theorems \ref{kh4thm11} and \ref{kh4thm12} hold.

The material of \S\ref{kh47} on Poincar\'e duality and wrong way maps also all extends easily to the de Rham case, where as in \eq{kh4eq99}, for $Y$ an oriented manifold of dimension $m$ we define Poincar\'e duality morphisms
\begin{gather*}
\Pd:MC^k_{\cs,\dR}(Y;\R)\ra MC_{m-k}^\dR(Y;\R),\;\>
\Pd:MC^k_\dR(Y;\R)\ra MC_{m-k}^{\lf,\dR}(Y;\R),
\nonumber\\
\text{by}\qquad \Pd:\al\longmapsto \al\cap[Y],
\end{gather*}
for $[Y]=[Y,0,0,\id_Y,1_Y]$ in $MC_m^{\lf,\dR}(Y;\R)$ as in~\S\ref{kh524}.

\subsection{M-homology and M-cohomology of effective orbifolds}
\label{kh53}

We now extend the theories of M-(co)homology, rational M-(co)homology, and de Rham M-(co)homology $MH_*(Y;R),\ldots,MH^*_\dR(Y;R)$ in \S\ref{kh4}--\S\ref{kh52}, from $Y$ a manifold, to $Y$ an effective orbifold. Orbifolds were discussed in \S\ref{kh29}, and Assumption \ref{kh3ass9} in \S\ref{kh35} explained the modifications to Assumptions \ref{kh3ass1}--\ref{kh3ass7} and \ref{kh3ass8} needed for the extension to orbifolds.

We have chosen to work in the ordinary category $\Orbeff$ of effective orbifolds defined in \S\ref{kh291}. Here are some comments about this:
\begin{itemize}
\setlength{\itemsep}{0pt}
\setlength{\parsep}{0pt}
\item[(i)] As in Remark \ref{kh2rem5}(a), if $Y$ is a noneffective orbifold, there is a natural effective orbifold $Y^\eff$ with a projection $\pi:Y\ra Y^\eff$ which is a homeomorphism of topological spaces, so the (co)homology of $Y$ and $Y^\eff$ is the same, and we lose little by restricting to effective orbifolds.
\item[(ii)] As in \S\ref{kh291}, there are many ways of defining categories or 2-categories of orbifolds, not all equivalent. Our $\Orbeff$ is quite crude, in that it is a category rather than a 2-category, and its morphisms are continuous maps satisfying  conditions, rather than continuous maps plus extra data.

This means that $\Orbeff$ has some bad differential-geometric behaviour, e.g.\ pullbacks $f^*(E)$ of orbifold vector bundles $E\ra Y$ by morphisms $f:X\ra Y$ in $\Orbeff$ cannot be defined. But to compensate for this, $\Orbeff$ has good topological behaviour, which means it is easy to extend results on (co)homology of manifolds to effective orbifolds, as in~\S\ref{kh293}.

For most other categories or 2-categories of (not necessarily effective) orbifolds $\mathfrak{Orb}$ in the literature, there is a forgetful functor $\Pi:\mathfrak{Orb}\ra\Orbeff$, so that our (co)homology theories pull back from $\Orbeff$ to~$\mathfrak{Orb}$.
\item[(iii)] Because of the bad differential-geometric behaviour in (ii), we adopt restrictive notions of {\it submersions\/} $f:X\ra Y$ and {\it transverse morphisms\/} $g:X\ra Z$, $h:Y\ra Z$ in $\Orbeff$, explained in \S\ref{kh292}, including conditions on how $f,g,h$ act on orbifold groups. These restrictive notions are used in the orbifold version of Assumption \ref{kh3ass5}, and to define M-cohomology.

This means, for example, that projections to quotient orbifolds $\pi:V\ra[V/G]$, and from tangent bundles of orbifolds $\pi:TX\ra X$, do not count as submersions. But it ensures that fibre products $X\t_{g,Z,h}Y$ of transverse morphisms $g,h$ exist in $\Orbeff$, in the sense of category theory, which will be needed to define pullbacks in M-cohomology. 
\end{itemize}

The category $\Orbeff$ and the results of \S\ref{kh29} were designed to ensure that \S\ref{kh4}--\S\ref{kh52} extend from manifolds to effective orbifolds with very little extra work. For the large majority of the material, we simply replace `manifold' by `effective orbifold' throughout, and no further nontrivial changes are required. 

The only things we have to be careful about are those issues, discussed in Properties \ref{kh2pr3}--\ref{kh2pr5} in \S\ref{kh293}, involving special properties of (co)homology of manifolds, rather than topological spaces. These are (a) fundamental classes of manifolds in homology; (b) the expression \eq{kh2eq29} for homology of manifolds in terms of cohomology of the orientation sheaf; and (c)  the material in \S\ref{kh28} on Poincar\'e duality and wrong way maps~$f^!,f_!$.

For the whole of \S\ref{kh53}, fix a 2-category $\tOrbeffc$ of `effective orbifolds with corners' satisfying Assumption \ref{kh3ass9}(a)--(e) in \S\ref{kh35} (the orbifold analogue of Assumptions \ref{kh3ass1}--\ref{kh3ass7}) when we are discussing (rational) M-(co)homology in \S\ref{kh531}--\S\ref{kh537}, and Assumption \ref{kh3ass9}(a)--(f) (the orbifold analogue of Assumptions \ref{kh3ass1}--\ref{kh3ass7} and \ref{kh3ass8}) when we are discussing de Rham M-(co)homology in \S\ref{kh538}. For simplicity, objects $X$ in $\tOrbeffc$ will be called {\it effective orbifolds with corners}, and morphisms $f:X\ra Y$ in $\tOrbeffc$ will be called {\it smooth maps}.

\subsubsection{M-homology $MH_*(Y;R)$ of effective orbifolds $Y$}
\label{kh531}

Fix a commutative ring $R$ throughout. All the definitions, results, and examples in \S\ref{kh41} (though not Remark \ref{kh4rem1}, as we explain in Remark \ref{kh5rem2} below), and the proofs of Propositions \ref{kh4prop1}, \ref{kh4prop2} and Theorems \ref{kh4thm1}, \ref{kh4thm2} in \S\ref{kh71}--\S\ref{kh74}, generalize to effective orbifolds immediately, with no nontrivial changes.

So, in Definition \ref{kh4def1}, we take $Y$ to be an effective orbifold (object in $\Orbeff$), and we define generators $[V,n,s,t]$ to be $\sim$-equivalence classes of quadruples $(V,n,s,t)$, where $V$ is an oriented effective orbifold with corners (i.e.\ a pair $(V,o_V)$ with $V$ an object in $\tOrbeffc$ and $o_V$ an orientation on $V$, usually left implicit), and $n=0,1,\ldots,$ and $s:V\ra\R^n$ is a smooth map (morphism in $\tOrbeffc$) which is proper over an open neighbourhood of 0 in $\R^n$, and $t:V\ra Y$ is a smooth map (morphism in~$\tOrbeffc$). 

We go on to define the complex of {\it M-chains\/} $\bigl(MC_k(Y;R),\pd\bigr)$, and its homology, ({\it integral\/}) {\it M-homology\/} $MH_*(Y;R)$. If $Y$ is a compact, oriented effective orbifold with $\dim Y=m$ we define the {\it fundamental cycle\/} $[Y]=[Y,0,0,\id_Y]$ in $MC_m(Y;R)$, and {\it fundamental class\/}~$[[Y]]\in MH_m(Y;R)$.

In Definition \ref{kh4def2} we define {\it pushforwards\/} $f_*$ on M-chains and M-homology when $f:Y_1\ra Y_2$ is a smooth map of effective orbifolds (morphism in $\Orbeff$). Theorem \ref{kh4thm1} defines a complex $\uMC_\bu(Y;R)$ of flabby cosheaves on~$Y$.

From Definition \ref{kh4def3} through to Theorem \ref{kh4thm2}, we prove that $MH_*(-;R)$ satisfies the analogue of Axiom \ref{kh2ax1} for effective orbifolds, and so is a homology theory of effective orbifolds in the sense of Definition \ref{kh2def15}. Thus Theorem \ref{kh2thm7}(a) implies the analogue of Theorem \ref{kh4thm3}:

\begin{thm} M-homology is a homology theory of effective orbifolds. There are canonical isomorphisms $MH_k(Y;R)\cong H_k(Y;R)$ and\/ $MH_k(Y,Z;R)\cong H_k(Y,Z;R)$ for all effective orbifolds\/ $Y,$ open suborbifolds $Z\subseteq Y$ and\/ $k\in\Z,$ preserving the data $f_*,\pd$ and isomorphisms $MH_0(*;R)\cong R\cong H_0(*;R),$ where $H_*(-;R)$ is any other homology theory of effective orbifolds over $R,$ such as singular homology\/~$H_*^\rsi(-;R)$.
\label{kh5thm11}
\end{thm}

\begin{rem} Only one issue in \S\ref{kh41} needs modification for the orbifold case, which is Remark \ref{kh4rem1}, where we suggested that if $R$ has characteristic 2, e.g.\ $R=\Z_2$, then we could define $MC_*(Y;R)$ and $MH_*(Y;R)$ using generators $[V,n,s,t]$ in which $V$ is unoriented. We noted that this would require an unoriented version of Assumption \ref{kh3ass6}(n), saying that if $X\in\tManc$ is compact with $\dim X=1$, then the number of points in $\pd X$ is zero modulo 2. This would be needed, for example, for the proof in \S\ref{kh71} that $\pd$ is well defined on~$MC_*(Y;R)$.

There is a problem in extending this to orbifolds if objects in $\Orbeff$ and $\tOrbeffc$ are allowed to have orbifold singularities in real codimension 1, as is permitted in Definition \ref{kh2def12}, and discussed in Property \ref{kh2pr3}(b).  For example, $X=\bigl[[-1,1]/\{\pm 1\}\bigr]$ is a compact, unoriented orbifold with boundary, with orbifold singularities in codimension 1, but $\pd X$ is one point, so the above unoriented version of Assumption \ref{kh3ass6}(n) fails.
\label{kh5rem2}
\end{rem}

\subsubsection{M-cohomology $MH^*(Y;R)$ of effective orbifolds $Y$}
\label{kh532}

Fix a commutative ring $R$. All the definitions, results, and examples in \S\ref{kh42}, and the proofs of Propositions \ref{kh4prop4}, \ref{kh4prop5} and \ref{kh4prop6} in \S\ref{kh75}--\S\ref{kh77}, generalize to effective orbifolds immediately, with no nontrivial changes.

So, for $Y$ an effective orbifold (object in $\Orbeff$) and $R$ a commutative ring, we define the complexes of {\it M-precochains\/} $\bigl(\cP MC^*(Y;R),\d\bigr)$ and {\it M-cochains\/} $\bigl(MC^*(Y;R),\d\bigr)$, and ({\it integral\/}) {\it M-cohomology\/} $MH^*(Y;R)$. The generators $[V,n,s,t]$ of $\cP MC^k(Y;R),MC^k(Y;R)$ are $\sim$-equivalence classes of quadruples $(V,n,s,t)$, where $V$ is an effective orbifold with corners (object in $\tOrbeffc$), $s:V\ra\R^n$ a smooth map (morphism in $\tOrbeffc$) for $n=0,1,\ldots,$ and $t:V\ra Y$ a cooriented submersion (submersion in $\tOrbeffc$, with an implicit coorientation $c_t$), with $(s,t):V\ra\R^n\t Y$ proper near $\{0\}\t Y$ in~$\R^n\t Y$.

For $f:Y_1\ra Y_2$ a smooth map of effective orbifolds (morphism in $\Orbeff$) we define {\it pullbacks\/} $f^*$ on M-precochains, M-cochains and M-cohomology. The definition involves fibre products $V\t_{t,Y_2,f}Y_1$, which exist in $\tOrbeffc$ as $t,y$ are transverse since $t$ is a submersion, generalizing the definition of submersions and transverse morphisms in $\Orbeff$ in~\S\ref{kh292}.

From Definition \ref{kh4def7} through to Lemma \ref{kh4lem4}, we prove that $MH^*(-;R)$ satisfies the analogue of Axiom \ref{kh2ax2} for effective orbifolds, and so is a cohomology theory of effective orbifolds in the sense of Definition \ref{kh2def15}. Thus Theorem \ref{kh2thm7}(b) implies the analogue of Theorem \ref{kh4thm4}:

\begin{thm} M-cohomology is a cohomology theory of effective orbifolds. There are canonical isomorphisms $MH^k(Y;R)\!\cong\! H^k(Y;R),$ $MH^k(Y,Z;R)\!\cong\! H^k(Y,Z;R)$ for all effective orbifolds\/ $Y,$ open suborbifolds $Z\subseteq Y$ and\/ $k\in\Z,$ preserving the data $f^*,\d$ and isomorphisms $MH^0(*;R)\cong R\cong H^0(*;R),$ where $H^*(-;R)$ is any other cohomology theory of effective orbifolds over $R,$ such as singular cohomology\/ $H^*_\rsi(Y;R)$ or sheaf cohomology\/~$H^*(Y,R_Y)$.
\label{kh5thm12}
\end{thm}

As in Definition \ref{kh4def6}, for each effective orbifold $Y$ we define a complex $\MC^\bu(Y;R)$ of soft sheaves of $R$-modules on $Y$, and as in Theorem \ref{kh4thm5} we show that $\MC^\bu(Y;R)=\bigl(\MC^*(Y;R),\d\bigr)$ is a soft resolution of~$R_Y$.

\subsubsection{Compactly-supported M-cohomology of effective orbifolds}
\label{kh533}

All of \S\ref{kh43}, and the proof of Proposition \ref{kh4prop8} in \S\ref{kh78}, also generalize to effective orbifolds immediately, with no nontrivial changes. Thus for effective orbifolds $Y\in\Orbeff$ we define ({\it integral\/}) {\it compactly-supported M-cohomology\/} $MH^*_\cs(Y;R)$, with pullbacks $f^*$ by proper morphisms $f:Y_1\ra Y_2$ in $\Orbeff$, and as in Theorem \ref{kh4thm6} we prove:

\begin{thm} Compactly-supported M-cohomology is a compactly-supported cohomology theory of effective orbifolds. There are canonical isomorphisms $MH^k_\cs(Y;R)\cong H^k_\cs(Y;R)$ for all effective orbifolds $Y$ and\/ $k\in\Z,$ preserving the data $\Pi,f^*,i_*$ in Property\/ {\rm\ref{kh2pr1}(a)--(c)} and the isomorphisms $MH^0_\cs(*;R)\cong R\cong H^0_\cs(*;R),$ where $H^*_\cs(-;R)$ is any other compactly-supported cohomology theory of effective orbifolds over $R,$ such as compactly-supported singular cohomology\/ $H^*_{\cs,\rsi}(Y;R),$ or compactly-supported sheaf  cohomology\/ $H^*_\cs(Y,R_Y)$.
\label{kh5thm13}
\end{thm}

\subsubsection{Locally finite M-homology of effective orbifolds}
\label{kh534}

For locally finite M-homology in \S\ref{kh44}, from Definition \ref{kh4def11} through to Definition \ref{kh4def13}, the material extends from manifolds to effective orbifolds with only trivial changes. Thus, for an effective orbifold $Y$ we define the complex of {\it locally finite M-chains\/} $\bigl(MC_*^\lf(Y;R),\pd\bigr)$, and ({\it integral\/}) {\it locally finite M-homology\/} $MH_*^\lf(Y;R)$, with pushforwards $f_*$ by proper morphisms $f:Y_1\ra Y_2$ in~$\Orbeff$. 

We define a complex $\MC_\bu^\lf(Y;R)$ of soft sheaves of $R$-modules on $Y$, with $\MC_k^\lf(Y;R)(U)=MC_k^\lf(U;R)$ for open $U\subseteq Y$, and we show as in Corollary \ref{kh4cor1} that the complex $\uMC_\bu(Y;R)$ of flabby cosheaves on $Y$ from \S\ref{kh531} and the complex $\MC_\bu^\lf(Y;R)$ of soft sheaves on $Y$ above are canonically related as in Theorem~ \ref{kh2thm3}. As for Theorem \ref{kh4thm8} we prove:

\begin{thm} Locally finite M-homology is a locally finite homology theory of effective orbifolds. That is, there are canonical isomorphisms $H_k^\lf(Y;R)\cong MH_k^\lf(Y;R)$ for all effective orbifolds\/ $Y$ and\/ $k\in\Z,$ preserving the data $\Pi,f_*,i^*$ in Property\/ {\rm\ref{kh2pr2}(a)--(c)} and the isomorphisms $H_0^\lf(*;R)\cong R\cong MH_0^\lf(*;R),$ where $H_*^\lf(-;R)$ is any other locally finite homology theory over $R$. 
\label{kh5thm14}
\end{thm}

Since the equivalence $\hat\cC{}^{\lf,\ssi}_{-\bu}(Y;R)\simeq\om_Y$ in \eq{kh2eq39} also holds for orbifolds $Y$ as in \eq{kh2eq78}, as for \eq{kh4eq56} we have a natural equivalence in $D(Y;R)$, where $\om_Y$ is the dualizing complex of~$Y$:
\e
\MC_{-\bu}^\lf(Y;R)\simeq \om_Y.
\label{kh5eq51}
\e

For the last part of \S\ref{kh44}, the morphism $j_Y:O_Y\ra\MC_m^\lf(Y;R)$ in Definition \ref{kh4def13} is well defined for any effective orbifold $Y$ and commutative ring $R$, although as in Property \ref{kh2pr4}(a) $O_Y$ is not locally constant if $Y$ is not locally orientable. But the proof of Theorem \ref{kh4thm8} relies on the fact that $H_k(Y,Y\sm\{y\};R)=0$ for $k<m$ for $Y$ a manifold of dimension $m$ and $y\in Y$. If instead $Y$ is an effective orbifold, this holds only if $R$ is a $\Q$-algebra. Thus we get the following weakened analogue of Theorem~\ref{kh4thm8}:

\begin{thm} For each effective orbifold\/ $Y$ of dimension $m,$ and each\/ $\Q$-algebra $R,$ the following is an exact sequence of sheaves of\/ $R$-modules on $Y\!:$
\begin{equation*}
\xymatrix@C=13pt{ 0 \ar[r] & O_Y \ar[r]^(0.3){j_Y} & \MC_m^\lf(Y;R) \ar[r]^(0.46)\pd & \MC_{m-1}^\lf(Y;R) \ar[r]^(0.48)\pd & \MC_{m-2}^\lf(Y;R) \ar[r]^(0.66)\pd & \cdots. }
\end{equation*}
Hence $\MC_{m-\bu}^\lf(Y;R)=\bigl(\MC_{m-*}(Y;R),\pd\bigr)$ is a soft resolution of\/~$O_Y$.
\label{kh5thm15}
\end{thm}

\subsubsection{Cup, cap and cross products}
\label{kh535}

The material of \S\ref{kh45}--\S\ref{kh46} on cup, cap and cross products on M-(co)homology, and the proof of Proposition \ref{kh4prop11} in \S\ref{kh79}, extend to effective orbifolds with almost no significant changes. 

The one point which needs attention is the proof of Theorem \ref{kh4thm11} in \S\ref{kh46}. This used the facts in \S\ref{kh262} that if $Y$ is a manifold of dimension $m$ then $\om_Y\simeq O_Y[m]$, and the cap product may be defined using the isomorphism \eq{kh2eq62} 
\begin{equation*}
I_\cap:R_Y\ot_R O_Y\,{\buildrel\cong\over\longra}\, O_Y,
\end{equation*}
applied to the soft resolution $\MC^\bu(Y;R)$ of $R_Y$ from Theorem \ref{kh4thm5}, and the soft resolution $\MC_{m-\bu}^\lf(Y;R)$ of $O_Y$ from Theorem \ref{kh4thm8}.

When $Y$ is an effective orbifold of dimension $m$, then if $R$ is not a $\Q$-algebra we may no longer have $\om_Y\simeq O_Y[m]$, as in Property \ref{kh2pr4}(c), and $\MC_{m-\bu}^\lf(Y;R)$ may no longer be a soft resolution of $O_Y$, as in Theorem~\ref{kh5thm15}.

There is an easy solution: as in \S\ref{kh262}, equation \eq{kh2eq62} is the specialization to manifolds of a quasi-isomorphism \eq{kh2eq61}
\begin{equation*}
I_\cap:R_Y\ot_R \om_Y\,{\buildrel\simeq\over\longra}\, \om_Y,
\end{equation*}
which works for sufficiently nice topological spaces, including effective orbifolds, and as in \eq{kh5eq51} we have $\MC_{-\bu}^\lf(Y;R)\simeq\om_Y$. So the proof of Theorem \ref{kh4thm11} extends to orbifolds if we use $\om_Y$ in place of $O_Y[m]$, and \eq{kh2eq61} in place of \eq{kh2eq62}, and \eq{kh5eq51} in place of Theorem \ref{kh4thm8}. We should probably restrict to $R$ a noetherian ring to use results on the dualizing complex~$\om_Y$.

Thus, we may define cup, cap and cross products $\cup,\cap,\t$ on M-(co)homology of effective orbifolds, at the (co)chain level, by the formulae \eq{kh4eq61}, \eq{kh4eq78}, \eq{kh4eq82}, \eq{kh4eq92}, and they induce products $\cup,\cap,\t$ on (co)homology which are identified with the usual cup, cap and cross products by the canonical isomorphisms of M-(co)homology with ordinary (co)homology, as in Theorems \ref{kh4thm9}, \ref{kh4thm11}, \ref{kh4thm12} and Corollary~\ref{kh4cor2}.

As in the manifold case in Remark \ref{kh4rem9}, for integral M-cohomology of effective orbifolds the cup product $\cup$ is associative but not supercommutative on M-cochains. The analogue of Theorem \ref{kh4thm10} holds, saying that for each effective orbifold $Y$ and commutative ring $R$, the dga $\bigl(MC^*(Y;R),\d,\cup,\Id_Y\bigr)$ over $R$ is equivalent in $\dga_R^\iy$ to the `usual' dga over $R$ associated to $Y$ in topology.

\subsubsection{Poincar\'e duality and wrong way maps}
\label{kh536}

The material of \S\ref{kh47} on Poincar\'e duality and wrong way maps $f_!,f^!$ needs modification for the orbifold case. Let $Y$ be an oriented effective orbifold, and $R$ a commutative ring. Then as in \eq{kh4eq99}, we can define morphisms
\begin{gather*}
\Pd:MC^k_\cs(Y;R)\longra MC_{m-k}(Y;R),\quad
\Pd:MC^k(Y;R)\longra MC_{m-k}^\lf(Y;R),
\nonumber\\
\text{by}\qquad \Pd:\al\longmapsto \al\cap[Y].
\end{gather*}
These then induce morphisms
\e
\begin{split}
&\Pd:MH^k_\cs(Y;R)\longra MH_{m-k}(Y;R),\\ 
&\Pd:MH^k(Y;R)\longra MH_{m-k}^\lf(Y;R),
\end{split}
\label{kh5eq52}
\e
which are identified with the usual morphisms \eq{kh2eq69} by the identifications $MH_*(Y;R)\cong H_*(Y;R),\ldots$ in \S\ref{kh531}--\S\ref{kh534}. But as in Property \ref{kh2pr5}(a) in \S\ref{kh293}, in the orbifold case equation \eq{kh5eq52} are isomorphisms if $R$ is a $\Q$-algebra, but may not be otherwise.

Following Definition \ref{kh4def21}, if $Y,Z$ are oriented orbifolds of dimensions $m,n$, $f:Y\ra Z$ is a proper submersion, and $R$ is a commutative ring, we may define morphisms $f_!:MC_k(Z;R)\ra MC_{k-n+m}(Y;R)$ by \eq{kh4eq101}, which pass to morphisms $f_!:MC_k(Z;R)\ra MC_{k-n+m}(Y;R)$ on homology. Equation \eq{kh4eq102} commutes, and so induces a commutative diagram on (co)homology:
\e
\begin{gathered}
\xymatrix@C=130pt@R=15pt{ 
*+[r]{MH^{n-k}_\cs(Z;R)} \ar[r]_\Pd \ar[d]^{f^*} & *+[l]{MH_k(Z;R)} \ar[d]_{f_!}  \\
*+[r]{MH^{n-k}_\cs(Y;R)} \ar[r]^\Pd  & *+[l]{MH_{k-n+m}(Y;R).\!}  }
\end{gathered}
\label{kh5eq53}
\e
If $R$ is a $\Q$-algebra, the rows of \eq{kh5eq53} are isomorphisms, so $f_!=\Pd\ci f^*\ci\Pd^{-1}$ in the usual way. But if $R$ is not a $\Q$-algebra, then $f_!$ may not be determined by $f^*$, and its (co)homological meaning is unclear. The same applies to the morphisms $f_!:MC_k^\lf(Z;R)\ra MC_{k-n+m}^\lf(Y;R)$, $f^!:MC^k_\cs(Y;R)\ra MC^{k-m+n}_\cs(Z;R)$ and $f^!:MC^k(Y;R)\ra MC^{k-m+n}(Z;R)$ in~\S\ref{kh47}.

\subsubsection{Rational M-(co)homology of effective orbifolds}
\label{kh537}

Following \S\ref{kh531}--\S\ref{kh536}, we can extend the theory of rational M-(co)homology $MH_*^\Q(Y;R),MH^*_\Q(Y;R),\ldots$ in \S\ref{kh51} for manifolds $Y$ and $\Q$-algebras $R$ to the case when $Y$ is an effective orbifold. This includes the extra relation Definition \ref{kh5def1}(iii) in the definitions of M-(co)chains $MC_k^\Q(Y;R),\ldots,$ and has the nice property that the cup product $\cup$ is supercommutative at the cochain level. 

This extension to the orbifold case involves no new issues that we have not already discussed for the integral M-(co)homology case in \S\ref{kh531}--\S\ref{kh536}. Actually,  because $R$ is now a $\Q$-algebra, most of the problems discussed in \S\ref{kh531}--\S\ref{kh536} go away, since as in Property \ref{kh2pr4}(c) we have $\om_Y\simeq O_Y[m]$ for $Y$ an effective orbifold of dimension $m$ when $R$ is a $\Q$-algebra.

\subsubsection{De Rham M-(co)homology of effective orbifolds}
\label{kh538}

Following \S\ref{kh531}--\S\ref{kh537}, we can extend the theory of de Rham M-(co)homology $MH_*^\dR(Y;\R),MH^*_\dR(Y;\R),\ldots$ in \S\ref{kh52} for manifolds $Y$ to the case when $Y$ is an effective orbifold. Again, the extension to the orbifold case involves no new issues that we have not already discussed in \S\ref{kh531}--\S\ref{kh536}, and as $\R$ is a $\Q$-algebra, most of the problems discussed in \S\ref{kh531}--\S\ref{kh536} go away.

\subsection{Extending M-(co)homology to a bivariant theory}
\label{kh54}

In \S\ref{kh210} we discussed Fulton and MacPherson's notion of {\it bivariant theories\/} \cite{FuMa}, which are a mixture of homology and cohomology. Example \ref{kh2ex23} defined a bivariant theory $H^*(g:Y\ra Z;R)$ for smooth maps of manifolds $g:Y\ra Z$ extending homology $H_*(Y;R)$ and compactly-supported cohomology $H^*_\cs(Y;R)$. We said that bivariant theories for manifolds do not really give anything new, the power of the theory is its applications to singular spaces.

At the level of (co)homology, for manifolds $Y$ we can extend M-homology $MH_*(Y;R)$ in \S\ref{kh41} and compactly-supported M-cohomology $MH^*_\cs(Y;R)$ in \S\ref{kh43} to a bivariant theory exactly as in Example \ref{kh2ex23}, by using twisted versions of $MH^*_\cs(Y;R),MH_*(Y;R)$ in \eq{kh2eq85}. This would be isomorphic to Example \ref{kh2ex23}, and yield nothing new.

We will now explain how to define a kind of (partial) bivariant theory $MC^*(g:Y\ra Z;R)$ {\it at the level of (co)chains}, which generalizes M-chains $MC_*(Y;R)$ and compactly-supported M-cochains $MC^*_\cs(Y;R)$. It is partial because we define $MC^*(g:Y\ra Z;R)$ only for submersions $g:Y\ra Z$, not for all smooth maps. This is intended for applications in which it is useful to be able to mix homology and cohomology at the (co)chain level. 

For example, in Fukaya--Oh--Ohta--Ono's theory of Lagrangian Floer cohomology \cite{FOOO1}, one studies moduli spaces $\oM_{k+1}(J,\be)$ of $J$-holomorphic maps $u:\Si\ra S$ from a prestable holomorphic disc $\Si$ into a symplectic manifold $(S,\om)$ with boundary $u(\pd\Si)$ in a Lagrangian $L\subset S$, and $k+1$ boundary marked points giving `evaluation maps' ${\rm ev}_i:\oM_{k+1}(J,\be)\ra L$ for $i=1,\ldots,k+1$. We need to define a `virtual (co)chain' $[\oM_{k+1}(J,\be)]_{\rm virt}$ in the (co)homology of~$L^{k+1}$. 

Now $k$ of the marked points are `inputs', associated to homology, and one is an `output', associated to cohomology. So at the (co)homology level, the virtual (co)chain should live in $H_*(L^k;\Q)\ot_\Q H^*(L;\Q)$, which is the bivariant group $H^*(\pi_{k+1}:L^{k+1}\ra L;\Q)$ from Example \ref{kh2ex23}. Thus, to define Lagrangian Floer cohomology in the most functorial way, it may be helpful to have a `bivariant chain complex' $\bigl(C^*(\pi_{k+1}:L^{k+1}\ra L;\Q),\d\bigr)$ with cohomology $H^*(\pi_{k+1}:L^{k+1}\ra L;\Q)$, and to construct $[\oM_{k+1}(J,\be)]_{\rm virt}$ in~$C^*(\pi_{k+1}:L^{k+1}\ra L;\Q)$.

As in \S\ref{kh4}, fix a commutative ring $R$, and a category $\tManc$ satisfying Assumptions \ref{kh3ass1}--\ref{kh3ass7} of~\S\ref{kh33}.

\begin{dfn} Let $g:Y\ra Z$ be a submersion of manifolds, with $\dim Y=l$, $\dim Z=m$. Consider quadruples $(V,n,s,t)$, where $V$ is a manifold with corners (object in $\tManc$), and $n=0,1,\ldots,$ and $s:V\ra\R^n$ and $t:V\ra Y$ are smooth maps (morphisms in $\tManc$), such that $g\ci t:V\ra Z$ is a submersion, we are given a coorientation $c_{g\ci t}$ for $g\ci t$ (which we leave implicit), and $s:V\ra\R^n$ is proper over an open neighbourhood $X$ of $\{0\}$ in~$\R^n$.

Define an equivalence relation $\sim$ on such quadruples by $(V,n,s,t)\sim(V',\ab n',\ab s',\ab t')$ if $n=n'$, and there exists a diffeomorphism $f:V\ra V'$ with $s=s'\ci f$ and $t=t'\ci f$ such that the coorientations satisfy $c_{g\ci t}=c_{g\ci t'}\ci c_f$, where $c_f$ is the natural coorientation on $f$ from Assumption \ref{kh3ass6}(e). Write $[V,n,s,t]$ for the $\sim$-equivalence class of $(V,n,s,t)$. We call $[V,n,s,t]$ a {\it generator}.

For each $k\in\Z$, define the ({\it integral\/}) {\it M-bichains\/} $MC^k(g:Y\ra Z;R)$ to be the $R$-module generated by such $[V,n,s,t]$ with $\dim V+k=m+n$, subject to the relations:
\begin{itemize}
\setlength{\itemsep}{0pt}
\setlength{\parsep}{0pt}
\item[(i)] For each generator $[V,n,s,t]$ and each $i=0,\ldots,n$ we have
\begin{equation*}
[V,n,s,t]=(-1)^{n-i}[V\t\R,n+1,s',t\ci\pi_V]\quad\text{in $MC^k(g:Y\ra Z;R)$,}
\end{equation*}
where writing $s=(s_1,\ldots,s_n):V\ra\R^n$ with $s_j:V\ra\R$ for $j=1,\ldots,n$ and $\pi_V:V\t\R\ra V$, $\pi_\R:V\t\R\ra\R$ for the projections, then 
\begin{equation*}
s'=(s_1\ci\pi_V,\ldots,s_i\ci\pi_V,\pi_\R,s_{i+1}\ci\pi_V,\ldots,s_n\ci\pi_V):V\t\R\ra\R^{n+1},
\end{equation*}
and $g\ci t\ci\pi_V$ has coorientation $c_{g\ci t\ci\pi_V}=c_{g\ci t}\ci c_{\pi_V}$, where $c_{g\ci t}$ is the given coorientation on $g\ci t:V\ra Z$, and $c_{\pi_V}$ is the coorientation on $\pi_V:V\t\R\ra V$ induced by the standard orientation on $\R$, as in Assumption~\ref{kh3ass6}(d),(f),(k).
\item[(ii)] Let $I$ be a finite indexing set, $a_i\in R$ for $i\in I$, and $[V_i,n,s_i,t_i]$, $i\in I$ be generators for $MC^k(g:Y\ra Z;R)$, all with the same $n$. Suppose there exists an open neighbourhood $X$ of $\{0\}$ in $\R^n$, such that $s_i:V_i\ra\R^n$ is proper over $X$ for all $i\in I$, and the following condition holds:
\begin{itemize}
\setlength{\itemsep}{0pt}
\setlength{\parsep}{0pt}
\item[$(*)$] Suppose $(x,y)\in X\t Y$ with $g(y)=z\in Z$, such that for all $i\!\in\! I$ and $v\!\in\! V_i$ with $(s_i,t_i)(v)\!=\!(x,y)$, we have that $v\in V_i^\ci$ and 
\begin{equation*}
T_v(s_i,t_i):T_vV_i^\ci\longra T_x\R^n\op T_yY
\end{equation*}
is injective. This implies that $(s_i,t_i)\vert_{V_i^\ci}$ is an embedding near $v\in V_i^\ci$. Hence $(s_i,t_i):V_i\ra\R^n\t Y$ is injective near each $v$ in $(s_i,t_i)^{-1}(x,y)$, so $(s_i,t_i)^{-1}(x,y)$ has the discrete topology, and thus is finite as $s_i$ is proper over $X$. Note too that $T_vV_i^\ci$ is a vector space of dimension $m+n-k$ and $\d g\vert_y\ci\d t\vert_v:T_vV_i^\ci\ra T_zZ$ is cooriented, since $g\ci t\vert_{V_i^\ci}:V_i^\ci\ra Z$ is a cooriented smooth map of manifolds by Assumption \ref{kh3ass6}(j). We require that for all $(m+n-k)$-planes $P\subseteq T_x\R^n\op T_yY$ with $\d g\vert_y\ci\pi_{T_yY}:P\ra T_zZ$ cooriented, we have
\begin{align*}
&\sum_{\begin{subarray}{l} i\in I,\; v\in V_i:(s_i,t_i)(v)=(x,y),\;  T_v(s_i,t_i)[T_vV_i^\ci]=P \\ \text{$T_v(s_i,t_i):T_vV_i^\ci\,{\buildrel\cong\over\longra}\,P$ is coorientation-preserving}\end{subarray}\!\!\!\!\!\!\!} a_i=\\
&\sum_{\begin{subarray}{l} i\in I,\; v\in V_i:(s_i,t_i)(v)=(x,y),\; T_v(s_i,t_i)[T_vV_i^\ci]=P \\ \text{$T_v(s_i,t_i):T_vV_i^\ci\,{\buildrel\cong\over\longra}\,P$ is coorientation-reversing}\end{subarray}\!\!\!\!\!\!\!} a_i\qquad\text{in $R$.}
\end{align*}
\end{itemize}
Then
\begin{equation*}
\sum_{i\in I}a_i\,[V_i,n,s_i,t_i]=0\qquad\text{in $MC^k(g:Y\ra Z;R)$.}
\end{equation*}
\end{itemize}
\label{kh5def12}
\end{dfn}

\begin{rem}{\bf(a)} When $g:Y\ra Z$ is $\pi:Y\ra *$, coorientations $c_{\pi\ci t}$ for $\pi\ci t:V\ra *$ are equivalent to orientations $o_V$ for $V$, and relations (i),(ii) above reduce to Definition \ref{kh4def1}(i),(ii). Hence from Definition \ref{kh4def1} we see that
\e
MC^k(\pi:Y\ra *;R)\cong MC_{-k}(Y;R).
\label{kh5eq54}
\e

When $g:Y\ra Z$ is $\id_Y:Y\ra Y$, generators $[V,n,s,t]$ in Definition \ref{kh5def12} are `compact generators' of $MC^k(Y;R)$ in the sense of Definition \ref{kh4def10}, and relations (i),(ii) above are Definition \ref{kh4def4}(i),(ii) applied to compact generators. Hence Proposition \ref{kh4prop8} implies that
\e
MC^k(\id_Y:Y\ra Y;R)\cong MC^k_\cs(Y;R).
\label{kh5eq55}
\e

\noindent{\bf(b)} For generators $[V,n,s,t]$ in Definition \ref{kh5def12}, we required $s:V\ra\R^n$ to be proper near 0 in $\R^n$, giving a theory generalizing $MC_*(Y;R)$ and $MC^*_\cs(Y;R)$. There are two alternative properness conditions we could have imposed:
\begin{itemize}
\setlength{\itemsep}{0pt}
\setlength{\parsep}{0pt}
\item[{\bf(i)}] $(s,t):V\ra\R^n\t Y$ should be proper near $\{0\}\t Y$ in $\R^n\t Y$; or
\item[{\bf(ii)}] $(s,g\ci t):V\ra\R^n\t Z$ should be proper near $\{0\}\t Z$ in~$\R^n\t Z$.
\end{itemize}
In both cases, as in \S\ref{kh42} we first define {\it M-prebichains\/} $\cP MC^k(g:Y\ra Z;R)$, and then make $MC^k(g:Y\ra Z;R)$ by sheafifying $\cP MC^k(g\vert_U:U\ra Z;R)$. In case (i) this would yield
\begin{align*}
MC^k_{\rm(i)}(\pi:Y\ra *;R)&\cong MC_{-k}^\lf(Y;R),\\
MC^k_{\rm(i)}(\id_Y:Y\ra Y;R)&\cong MC^k(Y;R),
\end{align*}
and in case (ii) it would yield
\begin{align*}
MC^k_{\rm(ii)}(\pi:Y\ra *;R)&\cong MC_{-k}(Y;R),\\
MC^k_{\rm(ii)}(\id_Y:Y\ra Y;R)&\cong MC^k(Y;R).
\end{align*}
We can extend these to bivariant theories as below, with differences on the properness requirements on morphisms for defining pushforwards and pullbacks. Here (i) corresponds to the bivariant theory defined by Fulton and MacPherson in \cite[\S 3.1]{FuMa}, and (ii) to that defined in~\cite[\S 3.3.1]{FuMa}.
\smallskip

\noindent{\bf(c)} We restrict to $g:Y\ra Z$ a submersion, since otherwise it would not be reasonable to require $g\ci t:V\ra Z$ to be a submersion in the definition of generators $[V,n,s,t]$ in Definition~\ref{kh5def12}.
\label{kh5rem3}
\end{rem}

\begin{dfn} Let $g:Y\ra Z$ be a submersion of manifolds. Define $\d:MC^k(g:Y\ra Z;R)\ra MC^{k+1}(g:Y\ra Z;R)$ to be the $R$-linear map satisfying
\begin{equation*}
\d[V,n,s,t]=[\pd V,n,s\ci i_V,t\ci i_V],
\end{equation*}
for all generators $[V,n,s,t]$, as in \eq{kh4eq3} and \eq{kh4eq20}. This is well-defined as in Proposition \ref{kh4prop1}, and satisfies $\d\ci\d=0$ as in Definitions \ref{kh4def1} and \ref{kh4def5}. Thus we may define the ({\it integral\/}) {\it M-bihomology\/} 
\begin{equation*}
MH^k(g:Y\!\ra\! Z;R)\!=\!\frac{\ts \Ker\bigl(\d:MC^k(g:Y\!\ra\! Z;R)\!\ra\! MC^{k+1}(g:Y\!\ra\! Z;R)\bigr)}{\ts \Im\bigl(\d:MC^{k-1}(g:Y\!\ra\! Z;R)\!\ra\! MC^k(g:Y\!\ra\! Z;R)\bigr)}\,.
\end{equation*}

\label{kh5def13}
\end{dfn}

We define products, pushforwards and pullbacks on M-bichains, as in~\S\ref{kh210}.

\begin{dfn} Let $f:X\ra Y$, $g:Y\ra Z$ be submersions of manifolds, with $\dim X=p$, $\dim Y=q$, and $\dim Z=r$. Define the {\it product\/}
\begin{equation*}
\cdot:MC^k(f:X\ra Y;R)\t MC^l(g:Y\ra Z;R)\longra MC^{k+l}(g\ci f:X\ra Z;R)
\end{equation*}
to be the unique $R$-linear map satisfying
\begin{align*}
[V,n,s&,t]\cdot[V',n',s',t']
\!=\!(-1)^{(l+q+r)n}[\ti V,\ti n,\ti s,\ti t]:=\!(-1)^{(l+q+r)n}\bigl[V\t_{f\ci t,Y,t'}V',\\
&n+n',(s_1\ci\pi_V,\ldots,s_n\ci\pi_V,s_1'\ci\pi_{V'},\ldots,s_{n'}'\ci\pi_{V'}),t\ci\pi_V\bigr],
\end{align*}
for all $[V,n,s,t]\in MC^k(f:X\ra Y;R)$ and $[V',n',s',t']\in MC^l(g:Y\ra Z;R)$. This specializes to the cup product in \eq{kh4eq61} and the cap product in \eq{kh4eq82} under \eq{kh5eq54}--\eq{kh5eq55}, as in \eq{kh2eq81}--\eq{kh2eq82}. A very similar proof to that of Proposition \ref{kh4prop11} in \S\ref{kh79} shows that `$\,\cdot\,$' is well-defined.

Let $f:X\ra Y$ be a smooth map of manifolds and $g:Y\ra Z$ a submersion, such that $g\ci f:X\ra Z$ is a submersion. Define the {\it pushforward\/}
\begin{equation*}
f_*:MC^k(g\ci f:X\ra Z;R)\longra MC^k(g:Y\ra Z;R)
\end{equation*}
to be the unique $R$-linear map acting on generators $[V,n,s,t]$ by
\begin{equation*}
f_*[V,n,s,t]=[V,n,s,f\ci t],
\end{equation*}
where the coorientation for $(g\ci f)\ci t:V\ra Z$ in $[V,n,s,t]\in MC^k(g\ci f:X\ra Z;R)$ is the coorientation for $g\ci(f\ci t):V\ra Z$ in $[V,n,s,f\ci t]\in MC^k(g:Y\ra Z;R)$. By \eq{kh4eq7}, when $g$ is $\pi:Y\ra *$ this specializes to $f_*:MC_{-k}(X;R)\ra MC_{-k}(Y;R)$ under \eq{kh5eq54}, as in \eq{kh2eq83}. This $f_*$ is well defined as for pushforwards $f_*$ on M-chains in Definition~\ref{kh4def2}.

Suppose $g:Y\ra Z$ is a submersion, and $h:Y\ra Z$ is any proper smooth map of manifolds. Then the fibre product $Y'=X\t_{g,Z,h}Z'$ is transverse and exists in $\Man$, and as in \eq{kh2eq79} we have a Cartesian square 
\e
\begin{gathered}
\xymatrix@C=70pt@R=13pt{
*+[r]{Y'} \ar[r]_{g'} \ar[d]^{h'} & *+[l]{Z'} \ar[d]_h \\
*+[r]{Y} \ar[r]^g & *+[l]{Z,\!} }
\end{gathered}
\label{kh5eq56}
\e
with $g'$ a submersion. We call \eq{kh5eq56} an {\it independent square}. Define the {\it pullback\/}
\begin{equation*}
h^*:MC^k(g:Y\ra Z;R)\longra MC^k(g':Y'\ra Z';R)
\end{equation*}
to be the unique $R$-linear map acting on generators $[V,n,s,t]$ by
\e
h^*[V,n,s,t]=[V',n,s',t']:=\bigl[V\t_{t,Y,h'}Y',n,s\ci\pi_V,\pi_{Y'}\bigr],
\label{kh5eq57}
\e
where $V'=V\t_{t,Y,h'}Y'$ is the fibre product in $\tManc$, which exists as
\begin{equation*}
V'=V\t_{t,Y,h'}Y'\cong V\t_{t,Y,\pi_Y}(Y\t_{g,Z,h}Z')\cong V\t_{g\ci t,Z,h}Z',
\end{equation*}
and $g\ci t:V\ra Z$ is a submersion. Then $g'\ci t':V'\ra Z'$ is the projection $\pi_{Z'}:V'\ra Z'$ in the fibre product $V'\cong V\t_{g\ci t,Z,h}Z'$, so $g'\ci t'$ is a submersion as $g\ci t$ is, and the coorientation on $g\ci t:V\ra Z$ determines one on $g'\ci t':V'\ra Z'$ by Assumption \ref{kh3ass6}(l). Also $h$ proper and $s:V\ra\R^n$ proper near $0$ in $\R^n$ imply that $s':V'\ra\R^n$ is proper near $0$. Hence $[V',n,s',t']$ in \eq{kh5eq57} is a generator of~$MC^k(g':Y'\ra Z';R)$.

Equation \eq{kh5eq57} specializes to \eq{kh4eq21} when $g:Y\ra Z$ is $\id_Y:Y\ra Y$, as in \eq{kh2eq84}. We show $h^*$ is well-defined by a very similar proof to Proposition~\ref{kh4prop4}.

Products, pushforwards and pullbacks are all compatible with the differential $\d$ on $MC^*(g:Y\ra Z;R)$, and so descend to products, pushforwards and pullbacks on M-bihomology~$MH^*(g:Y\ra Z;R)$. It is not difficult to show that these operations satisfy the axioms for a
bivariant theory in~\cite[\S 2.2]{FuMa}. 
\label{kh5def14}
\end{dfn}

Much of the theory of \S\ref{kh41}--\S\ref{kh47} can now be extended to M-bihomology. We outline a few important points. If $g:Y\ra Z$ is a submersion and $U\subseteq Y$ is open with inclusion $i:U\hookra Y$ then we have pushforwards
\begin{equation*}
i_*:MC^k(g\vert_U:U\ra Z;R)\longra MC^k(g:Y\ra Z;R).
\end{equation*}
Using these, as in Theorem \ref{kh4thm1} we can define a flabby cosheaf of $R$-modules $\uMC^k(g:Y\ra Z;R)$ on $Y$ with
\begin{equation*}
\uMC^k(g:Y\ra Z;R)(U)=MC^k(g\vert_U:U\ra Z;R)\quad\text{for $U\subseteq Y$ open.}
\end{equation*}
We have $\uMC^k(\pi:Y\ra *;R)\cong\uMC_{-k}(Y;R)$ and $\uMC^k(\id_Y:Y\ra Y;R)\cong\uMC_\cs^k(Y;R)$. The morphisms $\d:MC^k(g\vert_U:U\ra Z;R)\ra MC^{k+1}(g\vert_U:U\ra Z;R)$ induce cosheaf morphisms $\d:\uMC^k(g:Y\ra Z;R)\ra \uMC^{k+1}(g:Y\ra Z;R)$ with $\d\ci\d=0$.

If $\dim Y=l$, $\dim Z=m$ then we have $MC^k(g:Y\ra Z;R)=0$ and $\uMC^k(g:Y\ra Z;R)=0$ for $k<m-l$ as in Lemmas \ref{kh4lem1} and \ref{kh4lem3}. Write $\bigl(\MC^k(g:Y\ra Z;R),\d\bigr)$ for the complex of soft sheaves on $R$-modules on $Y$ corresponding to $\bigl(\uMC^k(g:Y\ra Z;R),\d\bigr)$ under Theorem \ref{kh2thm3}. Then as in Theorems \ref{kh4thm5} and \ref{kh4thm8}, we have an exact sequence of sheaves on~$Y:$
\begin{equation*}
\xymatrix@C=12pt{ 0\! \ar[r] & \!O_Y\!\ot_R\! g^*(O_Z)\! \ar[r]^(0.4){i_Y} & \!\MC^{m-l}(g\!:\!Y\!\ra\! Z;R)\! \ar[r]^(0.48)\d & \!\MC^{m-l+1}(g\!:\!Y\!\ra\! Z;R)\! \ar[r]^(0.77)\d & \cdots. }
\end{equation*}
Thus by \S\ref{kh25}, we may identify M-bihomology with a sheaf cohomology group
\begin{equation*}
MH^k(g:Y\ra Z;R)\cong H^{k+l-m}_\cs(Y,O_Y\ot_R g^*(O_Z)),
\end{equation*}
as in \eq{kh2eq85}. So M-bihomology is canonically isomorphic to the bivariant theory of Example~\ref{kh2ex23}.

We can extend rational M-(co)homology in \S\ref{kh51}, de Rham M-(co)homology in \S\ref{kh52}, and M-(co)homology of effective orbifolds in \S\ref{kh53}, to bivariant theories in a similar way.

\section{Proofs of results in \S\ref{kh2}}
\label{kh6}

\subsection{Proof of Theorems \ref{kh2thm1} and \ref{kh2thm2}}
\label{kh61}

For (co)homology theories defined on topological spaces rather than on manifolds, Eilenberg and Steenrod \cite[Th.~10.1]{EiSt} proved the analogue of Theorems \ref{kh2thm1} and \ref{kh2thm2} for `triangulable' topological spaces (or up to homotopy, finite CW-complexes), which include all compact smooth manifolds, without assuming Axioms \ref{kh2ax1}(vi) and~\ref{kh2ax2}(vi). 

By including (vi), Milnor \cite{Miln} extended the uniqueness theorems to the category of pairs $(Y,Z)$ where $Y,Z$ are topological spaces with the homotopy type of (possibly infinite) CW complexes, which include all smooth manifolds.

Kreck and Singhof \cite[Prop.~10]{KrSi} prove analogues of Theorems \ref{kh2thm1} and \ref{kh2thm2} giving axiomatic characterizations of (co)homology of manifolds, but with two differences: firstly, they use different sets of axioms involving only absolute rather than relative (co)homology, and secondly, because of their applications they consider (co)homology of smooth manifolds $Y$, but pushforwards and pullbacks by continuous maps $f:Y_1\ra Y_2$. Their method is to show that a (co)homology theory defined on smooth manifolds with continuous maps extends uniquely to a homology theory on topological spaces which are countable CW complexes of finite dimension, and then use the arguments of~\cite{EiSt,Miln}.

We will prove Theorems \ref{kh2thm1} and \ref{kh2thm2} using the method of \cite{KrSi}. We start by adapting material about triples $(Y;Z_1,Z_2)$ in \cite[\S 3]{KrSi} to pairs $(Y,Z)$ by taking $Z_1=Y$, $Z_2=Z$. Consider {\it pairs\/} $(Y,Z)$ of a topological space $Y$ and a subset $Z\subseteq Y$. A {\it morphism of pairs\/} $f:(Y_1,Z_1)\ra (Y_2,Z_2)$ is a continuous map $f:Y_1\ra Y_2$ with $f(Z_1)\subseteq Z_2$. We call $f$ a {\it pseudo-equivalence\/} if $f:Y_1\ra Y_2$ and $f\vert_{Z_1}:Z_1\ra Z_2$ are both homotopies of topological spaces. We call $f$ a {\it homotopy equivalence of pairs\/} if there exists a morphism $g:(Y_2,Z_2)\ra (Y_1,Z_1)$ such that $g:Y_2\ra Y_1$ is a homotopy inverse of $f:Y_1\ra Y_2$, and $g\vert_{Z_2}:Z_2\ra Z_1$ is a homotopy inverse of~$f\vert_{Z_1}:Z_1\ra Z_2$.

We will consider two special kinds of pairs: (a) pairs $(Y,Z)$ with $Y$ a manifold (considered as a topological space), and $Z\subseteq Y$ open, and (b) pairs $(Y,Z)$ with $Y$ a countable CW complex of finite dimension and $Z\subseteq Y$ a CW subcomplex. We call type (a) {\it manifold pairs}, and type (b) {\it CW pairs}. 

We call a manifold pair $(Y,Z)$ {\it good\/} if $Y\sm Z$ is a manifold with boundary embedded as a submanifold of $Y$ with $\dim(Y\sm Z)=Y$. Good manifold pairs can also be given the structure of CW pairs.

The next two lemmas follow from Kreck and Singhof~\cite[Prop.s 5 \& 6]{KrSi}.

\begin{lem} Suppose $(Y,Z)$ is a CW pair. Then there exists a good manifold pair $(Y',Z')$ and a homotopy equivalence of pairs\/~$h:(Y,Z)\ra (Y',Z')$.
\label{kh6lem1}
\end{lem}

\begin{lem} Suppose $(Y,Z)$ is a manifold pair. Then there exist a good manifold pair $(Y',Z'),$ a CW pair\/ $(Y'',Z''),$ a pseudo-equivalence $p:(Y',Z')\ra(Y,Z)$ with\/ $p:Y'\ra Y$ smooth, and a homotopy equivalence of pairs~$h':(Y',Z')\ra(Y'',Z'')$.
\label{kh6lem2}
\end{lem}

Kreck and Singhof do not use good manifold pairs, although $(Y',Z')$ good in Lemmas \ref{kh6lem1}--\ref{kh6lem2} follow easily from their construction. Since we want to work with {\it smooth\/} morphisms $f:(Y_1,Z_1)\ra(Y_2,Z_2)$ of manifold pairs (Kreck and Singhof use continuous morphisms), we need a result on replacing continuous morphisms by homotopic smooth morphisms, as in the next easy lemma, which is false if we do not assume $(Y_1,Z_1),(Y_2,Z_2)$ are good.

\begin{lem}{\bf(a)} Suppose $(Y_1,Z_1),(Y_2,Z_2)$ are good manifold pairs, and $f:(Y_1,Z_1)\ra (Y_2,Z_2)$ is a continuous morphism of pairs. Then we can find a smooth morphism of pairs $f':(Y_1,Z_1)\ra (Y_2,Z_2)$ close to $f$ in $C^0,$ and a continuous homotopy $g:(Y_1\t[0,1],Z_1\t[0,1])\ra(Y_2,Z_2)$ from $f$ to\/~$f'$.
\smallskip

\noindent{\bf(b)} Suppose $(Y_1,Z_1),(Y_2,Z_2)$ are good manifold pairs, $f',f'':(Y_1,Z_1)\ra (Y_2,Z_2)$ are smooth morphisms of pairs, and\/ $g:(Y_1\t[0,1],Z_1\t[0,1])\ra(Y_2,Z_2)$ is a continuous homotopy from $f'$ to $f''$. Then we can find a smooth homotopy $g':(Y_1\t[0,1],Z_1\t[0,1])\ra(Y_2,Z_2)$ from $f'$ to $f'',$ close to $g$ in\/~$C^0$.
\label{kh6lem3}
\end{lem}

Broadly following \cite[\S 4]{KrSi}, we prove:

\begin{prop} Suppose $H_*(-;R)$ is a homology theory on manifold pairs $(Y,Z)$ and smooth maps\/ $f:(Y_1,Z_1)\ra(Y_2,Z_2)$ satisfying Axiom\/ {\rm\ref{kh2ax1}}. Then we can construct a homology theory $\ti H_*(-;R)$ on CW pairs $(Y,Z)$ and continuous maps\/ $f:(Y_1,Z_1)\ra(Y_2,Z_2),$ uniquely up to canonical isomorphism, satisfying the analogue of Axiom\/ {\rm\ref{kh2ax1},} with functorial isomorphisms $H_*(Y,Z;R)\cong\ti H_*(Y,Z;R)$ for all manifold pairs $(Y,Z)$ which are also CW pairs.

Conversely, the CW homology theory $\ti H_*(-;R)$ determines the manifold homology theory $H_*(-;R)$ uniquely up to canonical isomorphism.

The analogue of all the above holds for cohomology theories, using Axiom\/~{\rm\ref{kh2ax2}}.

\label{kh6prop1}
\end{prop}

\begin{proof} Let $H_*(-;R)$ be as in the proposition. To define the CW homology theory $\ti H_*(-;R)$, using the Axiom of Choice, for each CW pair $(Y,Z)$ choose a good manifold pair $(Y',Z')$ and a homotopy equivalence of pairs $h:(Y,Z)\ra (Y',Z')$ as in Lemma \ref{kh6lem1}, and choose a homotopy inverse $i:(Y',Z')\ra(Y,Z)$ for $h$. Define~$\ti H^*(Y,Z)=H^*(Y',Z')$. 

Suppose $f:(Y_1,Z_1)\ra (Y_2,Z_2)$ is a morphism of CW pairs, and let $(Y_1',Z_1'),\ab h_1,i_1$, $(Y_2',Z_2'),h_2,i_2$ be the data chosen for $(Y_1,Z_1),(Y_2,Z_2)$. Then $h_2\ci f\ci i_1:(Y_1',Z_1')\ra(Y_2',Z_2')$ is a continuous morphism of good manifold pairs. So Lemma \ref{kh6lem3}(a) gives smooth $f':(Y_1',Z_1')\ra(Y_2',Z_2')$ continuously homotopic to $h_2\ci f\ci i_1$. Define $f_*:\ti H_*(Y_1,Z_1;R)\ra\ti H_*(Y_2,Z_2;R)$ to be $f'_*:H_*(Y_1',Z_1';R)\ra H_*(Y_2',Z_2';R)$. To see that $f_*$ is well-defined, let $f''$ be an alternative choice for $f'$. Then $f',f''$ are smooth morphisms $(Y_1',Z_1')\ra(Y_2',Z_2')$ which are continuously homotopic, so they are smoothly homotopic by Lemma \ref{kh6lem3}(b), and thus $f'_*=f''_*$ by Axiom~\ref{kh2ax1}(iv). 

Now let $(Y,Z)$ be a CW pair, so that $(Y,\es)$ and $(Z,\es)$ are also CW pairs. Let $(Y',Z'),h_1,i_1$ and $(Y'',\es),h_2,i_2$ and $(Z'',\es),h_3,i_3$ be the data chosen for $(Y,Z)$, $(Y,\es)$ and $(Z,\es)$. Choose smooth $f_2:Y''\ra Y'$ and $f_3:Z''\ra Z'$ which are continuously homotopic to $h_1\ci i_2:Y''\ra Y'$ and to $h_1\vert_Z\ci i_3:Z''\ra Z'$. Define $\pd:\ti H_k(Y,Z;R)\ra\ti H_{k-1}(Z;R)$ by the commutative diagram
\begin{equation*}
\xymatrix@C=12pt@R=12pt{ \cdots \ar[r] & \ti H_k(Z;R) \ar@{=}[d] \ar[r]_(0.47){i_*} & \ti H_k(Y;R) \ar@{=}[d] \ar[r]_(0.45){j_*} & \ti H_k(Y,Z;R) \ar@{=}[d] \ar@{.>}[r]_(0.48)\pd & \ar@{=}[d] \ti H_{k-1}(Z;R) \ar[r] & \cdots \\
& H_k(Z'';R) \ar[d]^{(f_3)_*}_\cong  & H_k(Y'';R) \ar[d]^{(f_2)_*}_\cong & H_k(Y',Z';R) \ar@{=}[d] & H_{k-1}(Z'';R) \ar[d]^{(f_3)_*}_\cong \\
\cdots \ar[r] & H_k(Z';R) \ar[r]^(0.47){i_*} & H_k(Y';R) \ar[r]^(0.45){j_*} & H_k(Y',Z';R) \ar[r]^(0.48)\pd & H_{k-1}(Z';R) \ar[r] & \cdots, }\!\!\!\!\!\!\!\!\!\!\!\!\!\!\!{}
\end{equation*}
where the top and bottom lines are \eq{kh2eq1} for $(Y,Z),(Y',Z')$, and the columns are isomorphisms as $f_2,f_3$ are homotopies. This defines all the data in the CW homology theory $\ti H(-;R)$. It is easy to deduce Axiom \ref{kh2ax1}(i)--(vi) for $\ti H$ from Axiom \ref{kh2ax1}(i)--(vi) for $H$, as in \cite[\S 5]{KrSi}. The first part of the proposition follows.

For the second part, to recover $H(-;R)$ from $\ti H(-;R)$, first note that for {\it good\/} manifold pairs $(Y,Z)$ the argument above constructing $\ti H(-;R)$ from $H(-;R)$ is reversible, as every good manifold pair $(Y,Z)$ is homotopic to a CW pair $(Y',Z')$ and vice versa. Thus $\ti H(-;R)$ determines $H(-;R)$ uniquely up to canonical isomorphism on the full subcategory $\cC_\go$ of $\cC$ in Axiom \ref{kh2ax1} with objects good manifold pairs~$(Y,Z)$.

To determine $H(-;R)$ on non-good manifold pairs $(Y,Z)$, note that Lemma \ref{kh6lem2} gives a good manifold pair $(Y',Z')$ and a pseudo-equivalence $p:(Y',Z')\ra(Y,Z)$. Thus we have a commutative diagram with exact rows
\e
\begin{gathered}
\xymatrix@C=12pt@R=15pt{ \cdots \ar[r] & H_k(Z';R) \ar[d]_\cong^{(p\vert_{Z'})_*}\ar[r]_(0.47){i_*} & H_k(Y';R) \ar[d]_\cong^{p_*} \ar[r]_(0.45){j_*} & H_k(Y',Z';R) \ar[d]^{p_*} \ar[r]_(0.48)\pd & H_{k-1}(Z';R) \ar[d]_\cong^{(p\vert_{Z'})_*} \ar[r] & \cdots \\
\cdots \ar[r] & H_k(Z;R) \ar[r]^(0.47){i_*} & H_k(Y;R) \ar[r]^(0.45){j_*} & H_k(Y,Z;R) \ar[r]^(0.48)\pd & H_{k-1}(Z;R) \ar[r] & \cdots. }\!\!\!\!\!\!\!\!\!\!\!\!\!\!\!{}
\end{gathered}
\label{kh6eq1}
\e
As $p$ is a pseudo-equivalence, $p:Y'\ra Y$ and $p\vert_{Z'}:Z'\ra Z$ are homotopies, so the first, second and fourth columns of \eq{kh6eq1} are isomorphisms. The five lemma now implies that $p_*:H_k(Y',Z';R)\ra H_k(Y,Z;R)$ is an isomorphism. So $H_*(Y,Z;R)$ for the non-good pair $(Y,Z)$ is determined up to canonical isomorphism by $H_*(Y',Z';R)$ for the good pair $(Y',Z')$. Using this, we see that $H(-;R)$ is determined up to canonical isomorphism by~$\ti H(-;R)$.

The proof for cohomology is the same, reversing morphisms on~$H^*,\ti H^*$.
\end{proof}

As noted by Kreck and Singhof in \cite[Lem.~11]{KrSi}, Milnor's proof \cite{Miln} of uniqueness of (co)homology theories satisfying the analogues of Axioms \ref{kh2ax1} and \ref{kh2ax2} for CW complexes also works for our `CW pairs' of a countable CW complex $Y$ of finite dimension and a CW subcomplex $Z\subseteq Y$. Thus Theorems \ref{kh2thm1} and \ref{kh2thm2} follow from Proposition~\ref{kh6prop1}.

\subsection{Proof of Theorem \ref{kh2thm4}}
\label{kh62}

Let $Y,R,\cE$ and $\pi:\cE\ra\hat\cE$ be as in Theorem \ref{kh2thm4} throughout.

\subsubsection{Proof of Theorem \ref{kh2thm4}(a)}
\label{kh621}

First consider Definition \ref{kh2def1}(iv),(v) when $I=\{1,2\}$, so that $U\subseteq Y$ is open and $\{V_1,V_2\}$ is an open cover of $U$. Then by \eq{kh2eq40}, as $\cE$ is strong the following is exact:
\e
\xymatrix@C=10pt{ 0 \ar[r] & \cE(U) \ar[rrr]^(0.4){\rho_{UV_1}\op \rho_{UV_2}} &&& \cE(V_1)\!\op\! \cE(V_2) \ar[rrrrrr]^(0.51){\rho_{V_1(V_1\cap V_2)}\op -\rho_{V_2(V_1\cap V_2)}} &&&&&& \ucE(V_1\!\cap\! V_2). }\!\!{}
\label{kh6eq2}
\e
Then Definition \ref{kh2def1}(iv) is equivalent to \eq{kh6eq2} being exact at the second place, and Definition \ref{kh2def1}(v) equivalent to \eq{kh6eq2} being exact at the third place. Hence $\cE$ being strong is equivalent to Definition \ref{kh2def1}(iv),(v) holding for $\cE$ when $\md{I}=2$. This proves the second part of (a), that is, if $\cE'$ is a presheaf and Definition \ref{kh2def1}(iv),(v) hold for $\cE'$ whenever $I$ is finite, then $\cE'$ is strong.

We will prove that Definition \ref{kh2def1}(iv),(v) hold for $\cE$ whenever $I=\{1,\ldots,n\}$. When $n=0,1$ this is trivial, and when $n=2$ it is equivalent to $\cE$ strong. Suppose by induction on $N$ that for some $N=2,3,\ldots,$ Definition \ref{kh2def1}(iv),(v) hold for $\cE$ whenever $I=\{1,\ldots,n\}$ for all $n=0,1,\ldots,N$. Let $U\subseteq Y$ be open and $\{V_1,\ldots,V_{N+1}\}$ be an open cover of~$U$. 

For (iv), suppose $s\in\cE(U)$ with $\rho_{UV_i}(s)=0$ in $\cE(V_i)$ for $i=1,\ldots,N+1$. Let $V_N'=V_N\cup V_{N+1}$ and $s'=\rho_{U\smash{V_N'}}(s)\in \cE(V_N')$. Then $\rho_{\smash{V_N'}V_N}(s')=\rho_{UV_N}(s)=0$ and $\rho_{\smash{V_N'}V_{N+1}}(s')=\rho_{UV_{N+1}}(s)=0$, so by the inductive hypothesis (iv) for $n=2$ we have $s'=0$. Now apply the inductive hypothesis (iv) for $n=N$ to $s$ on $U$ with open cover $\{V_1,\ldots,V_{N-1},V_N'\}$. We have $\rho_{UV_i}(s)=0$ in $\cE(V_i)$ for $i=1,\ldots,N-1$ and $\rho_{U\smash{V_N'}}(s)=s'=0$ in $\cE(V_N')$, so $s=0$, proving (iv) when~$I=\{1,\ldots,N+1\}$.

For (v), suppose $s_i\in\cE(V_i)$ for $i=1,\ldots,N+1$ with $\rho_{V_i(V_i\cap V_j)}(s_i)=\rho_{V_j(V_i\cap V_j)}(s_j)$ in $\cE(V_i\cap V_j)$ for all $i,j=1,\ldots,N+1$. Then the inductive hypothesis (v) for $n=2$ applied to the open cover $\{V_N,V_{N+1}\}$ of $V_{N'}$ gives $s_{N'}\in\cE(V_N')$ with $\rho_{\smash{V_N'}V_N}(s'_N)=s_N$ and $\rho_{\smash{V_N'}V_{N+1}}(s'_N)=s_{N+1}$. Let $i=1,\ldots,N-1$, and consider $\rho_{V_i(V_i\cap \smash{V_N'})}(s_i)$ and $\rho_{\smash{V_N'}(V_i\cap \smash{V_N'})}(s_N')$ in $\cE(V_i\cap V_N')$. We have
\begin{align*}
&\rho_{(V_i\cap \smash{V_N'})(V_i\cap V_N)}\bigl[\rho_{V_i(V_i\cap \smash{V_N'})}(s_i)\bigr]=\rho_{V_i(V_i\cap V_N)}(s_i)=\rho_{V_N(V_i\cap V_N)}(s_N)=\\
&\rho_{V_N(V_i\cap V_N)}\!\ci\!\rho_{\smash{V_N'}V_N}(s'_N)\!=\!\rho_{\smash{V_N'}(V_i\cap V_N)}(s'_N)\!=\!\rho_{(V_i\cap \smash{V_N'})(V_i\cap V_N)}\bigl[\rho_{\smash{V_N'}(V_i\cap \smash{V_N'})}(s_N')\bigr],
\end{align*}
and $\rho_{(V_i\cap \smash{V_N'})(V_i\cap V_{N+1})}\bigl[\rho_{V_i(V_i\cap \smash{V_N'})}(s_i)\bigr]=\rho_{(V_i\cap \smash{V_N'})(V_i\cap V_{N+1})}\bigl[\rho_{\smash{V_N'}(V_i\cap \smash{V_N'})}(s_N')\bigr]$ in the same way. Thus applying the inductive hypothesis (iv) for $n=2$ to the open cover $\{V_i\cap V_N,V_i\cap V_{N+1}\}$ of $V_i\cap V_N'$ shows that~$\rho_{V_i(V_i\cap \smash{V_N'})}(s_i)=\rho_{\smash{V_N'}(V_i\cap \smash{V_N'})}(s_N')$. 

Hence applying the inductive hypothesis (v) for $n=N$ to the open cover $\{V_1,\ldots,V_{N-1},V_N'\}$ of $U$ and sections $s_1,\ldots,s_{N-1},s_N'$ gives $s\in \cE(U)$ with $\rho_{UV_i}(s)=s_i$ for $i=1,\ldots,N-1$ and $\rho_{U\smash{V_N'}}(s)=s_N'$. But then
\begin{equation*}
\rho_{UV_N}(s)=\rho_{\smash{V_N'}V_N}\ci\rho_{U\smash{V_N'}}(s)=\rho_{\smash{V_N'}V_N}(s_N')=s_N,
\end{equation*}
and similarly $\rho_{UV_N}(s)=s_{N+1}$, proving (v) when $I=\{1,\ldots,N+1\}$. This completes the inductive step. Therefore by induction Definition \ref{kh2def1}(iv),(v) hold when $I=\{1,\ldots,n\}$ for all $n=0,1,\ldots,$ so they hold for all finite~$I$.

\subsubsection{Proof of Theorem \ref{kh2thm4}(b)}
\label{kh622}

Let $V\subseteq U\subseteq Y$ be open, and consider the morphism \eq{kh2eq41}. First we show that if $\al\in\cE_\cs(U)$ with $\supp\al\subseteq V$ then $\rho_{UV}(\al)\in\cE_\cs(V)$, that is, $\rho_{UV}(\al)$ is compactly-supported. This is not obvious since, although $\supp[\rho_{UV}(\al)]=\supp\al$ is compact, as in Remark \ref{kh2rem3} this does not imply $\rho_{UV}(\al)$ is compactly-supported unless $\cE$ is a sheaf rather than a presheaf.

As $\al$ is compactly-supported there is compact $K\subseteq U$ with $\rho_{U(U\sm K)}(\al)=0$ in $\cE(U\sm K)$. Then $K\sm V$ is compact in $U$, as it is closed in $K$, and $(K\sm V)\cap\supp\al=0$ as $\supp\al\subseteq V$. Hence by Definition \ref{kh2def6}, for each $y\in K\sm V$ there exists an open neighbourhood $W_y$ of $y$ in $U$ such that $\rho_{UW_y}(s)=0$ in $\cE(W_y)$. Then $\{W_y:y\in K\sm V\}$ is an open cover of $K\sm V\subseteq U$, which is compact, so there exists a finite set $\{y_1,\ldots,y_N\}\subseteq K\sm V$ with~$K\sm V\subseteq W_{y_1}\cup\cdots\cup W_{y_N}$.

Set $K'=K\sm(W_{y_1}\cup\cdots\cup W_{y_N})$. It is compact, as it is closed in $K$, and $K'\subseteq V$ as $K\sm V\subseteq W_{y_1}\cup\cdots\cup W_{y_N}$. Also $\{W_{y_1},\ldots,W_{y_N},U\sm K\}$ is an open cover of $U\sm K'$, and $\rho_{UW_{y_1}}(\al)=\cdots=\rho_{UW_{y_N}}(\al)=\rho_{U(U\sm K)}(\al)=0$. Thus Definition \ref{kh2def1}(iv) for $\cE$ with $I=\{1,\ldots,N+1\}$ applied to $s=\rho_{U(U\sm K')}(\al)$ in $\cE(U\sm K')$, which holds by Theorem \ref{kh2thm4}(a), gives $\rho_{U(U\sm K')}(\al)=0$. Therefore
\begin{equation*}
\rho_{V(V\sm K')}[\rho_{UV}(\al)]=\rho_{U(V\sm K')}(\al)=\rho_{(U\sm K')(V\sm K')}\ci\rho_{U(U\sm K')}(\al)=0,
\end{equation*}
and so $\rho_{UV}(\al)$ is compactly-supported, as $K'\subseteq V$ is compact, and $\rho_{UV}(\al)$ lies in $\cE_\cs(V)$. Thus the morphism \eq{kh2eq41} is well defined.

Suppose $\al\in\cE_\cs(U)$ with $\supp\al\subseteq V\subseteq U$ and $\rho_{UV}(\al)=0$. As above there exists compact $K'\subseteq V$ with $\rho_{U(U\sm K')}(\al)=0$. Exactness of \eq{kh6eq2} at the second term with $V_1=V$, $V_2=U\sm K'$ gives $\al=0$. Hence \eq{kh2eq41} is injective.

Suppose $\be\in\cE_\cs(V)$. Then by definition there exists compact $K'\subseteq V$ with $\rho_{V(V\sm K')}(\be)=0$. By exactness of \eq{kh6eq2} at the third term with $V_1=V$ and $V_2=U\sm K'$, since $\rho_{V(V\sm K')}(\be)=\rho_{(U\sm K')(V\sm K')}(0)=0$, there exists $\al\in\cE(U)$ with $\rho_{UV}(\al)=\be$ and $\rho_{U(U\sm K')}(\al)=0$. Then $\rho_{U(U\sm K')}(\al)=0$ implies that $\al\in\cE_\cs(U)$ with $\supp\al\subseteq K'\subseteq V$. Therefore \eq{kh2eq41} is surjective, so it is an isomorphism. This proves Theorem~\ref{kh2thm4}(b).

\subsubsection{Proof of Theorem \ref{kh2thm4}(c)}
\label{kh623}

Suppose $V\subseteq U\subseteq Y$ are open, and the closure $\bar V$ of $V$ in $U$ is compact, and let $\hat s\in\hat\cE(U)$. We first show that there exists a unique $\tau_{UV}(\hat s)\in\cE(V)$ satisfying the characterizing property in Theorem \ref{kh2thm4}(c) for each~$y\in\bar V$. 

By definition of sheafification, for each $y\in\bar V$ we can choose a small open neighbourhood $U_y$ of $y$ in $U$ and a section $s_y\in\cE(U_y)$  with $\pi(U_y)(s_y)=\hat\rho_{UU_y}(\hat s)$. As $Y$ is locally compact, we can also choose an open neighbourhood $U_y'$ of $y$ in $U_y$ such that the closure $\bar U_y'$ of $U_y'$ in $U_y$ is compact. Then $\{U_y':y\in\bar V\}$ is a family of open sets in $U$ which cover $\bar V\subseteq U$, so as $\bar V$ is compact there exist $y_1,\ldots,y_N$ in $\bar V$ with~$\bar V\subseteq U_{y_1}'\cup\cdots\cup U_{y_N}'$.

For $i,j=1,\ldots,N$, we have $s_{y_i}\in\cE(U_{y_i})$ and $s_{y_j}\in\cE(U_{y_j})$, but we do not know that $\rho_{U_{y_i}(U_{y_i}\cap U_{y_j})}(s_{y_i})=\rho_{U_{y_j}(U_{y_i}\cap U_{y_j})}(s_{y_j})$ in $\cE(U_{y_i}\cap U_{y_j})$. Instead, we only know that each $y\in U_{y_i}\cap U_{y_j}$ has an open neighbourhood $W_y$ of $y$ in $U_{y_i}\cap U_{y_j}$ with $\rho_{U_{y_i}W_y}(s_{y_i})=\rho_{U_{y_j}W_y}(s_{y_j})$ in $\cE(y_j)$. Now $\bar U_{y_i}'\cap\bar U_{y_j}'$ is a subset of $U_{y_i}\cap U_{y_j}$, and is compact as $\bar U_{y_i}',\bar U_{y_j}'$ are and $Y$ is Hausdorff. Then $\bigl\{W_y:y\in U_{y_i}\cap U_{y_j}\bigr\}$ is a family of open subsets of $U_{y_i}\cap U_{y_j}$ which cover $\bar U_{y_i}'\cap\bar U_{y_j}'$, so we can choose $\ti y_1,\ldots,\ti y_{\ti N}$ in $U_{y_i}\cap U_{y_j}$ with~$U_{y_i}'\cap U_{y_j}'\subseteq\bar U_{y_i}'\cap\bar U_{y_j}'\subseteq W_{\ti y_1}\cup\cdots\cup W_{\ti y_{\ti N}}$. 

Apply Theorem \ref{kh2thm4}(a) to the open cover $\bigl\{U_{y_i}'\cap U_{y_j}'\cap W_{\ti y_k}:k=1,\ldots,\ti N\bigr\}$ of $U_{y_i}'\cap U_{y_j}'$, and $\rho_{U_{y_i}(U'_{y_i}\cap U'_{y_j})}(s_{y_i}),\rho_{U_{y_j}(U'_{y_i}\cap U'_{y_j})}(s_{y_j})$ in $\cE(U'_{y_i}\cap U'_{y_j})$. Since
\begin{align*}
&\rho_{(U_{y_i}'\cap U_{y_j}')(U_{y_i}'\cap U_{y_j}'\cap W_{\ti y_k})}\bigl[\rho_{U_{y_i}(U'_{y_i}\cap U'_{y_j})}(s_{y_i})\bigr]=\rho_{U_{y_i}(U_{y_i}'\cap U_{y_j}'\cap W_{\ti y_k})}(s_{y_i})\\
&\quad =\rho_{U_{y_j}(U_{y_i}'\cap U_{y_j}'\cap W_{\ti y_k})}(s_{y_j})=\rho_{(U_{y_i}'\cap U_{y_j}')(U_{y_i}'\cap U_{y_j}'\cap W_{\ti y_k})}\bigl[\rho_{U_{y_j}(U'_{y_i}\cap U'_{y_j})}(s_{y_j})\bigr]
\end{align*}
for $k=1,\ldots,\ti N$, we have $\rho_{U_{y_i}(U'_{y_i}\cap U'_{y_j})}(s_{y_i})=\rho_{U_{y_j}(U'_{y_i}\cap U'_{y_j})}(s_{y_j})$ for all $i,j=1,\ldots,N$ by Definition~\ref{kh2def1}(iv).

Now apply Theorem \ref{kh2thm4}(a) to the open cover $\bigl\{U_{y_1}'\cap V,\ldots,U_{y_N}'\cap V\bigr\}$ of $V$, and the sections $\rho_{U_{y_i}(U_{y_i}'\cap V)}(s_{y_i})\in\cE(U_{y_i}'\cap V)$ for $i=1,\ldots,N$. Since
\begin{align*}
&\rho_{(U_{y_i}'\cap V)(U_{y_i}'\cap U_{y_j}'\cap V)}\bigl[
\rho_{U_{y_i}(U_{y_i}'\cap V)}(s_{y_i})\bigr]=\rho_{U_{y_i}(U_{y_i}'\cap U_{y_j}'\cap V)}(s_{y_i})\\
&=\rho_{U_{y_j}(U_{y_i}'\cap U_{y_j}'\cap V)}(s_{y_j})=
\rho_{(U_{y_j}'\cap V)(U_{y_i}'\cap U_{y_j}'\cap V)}\bigl[\rho_{U_{y_j}(U_{y_j}'\cap V)}(s_{y_j})\bigr]
\end{align*}
for all $i,j=1,\ldots,N$, this says that that there exists a unique element $\tau_{UV}(\hat s)$ in $\cE(V)$ with $\rho_{V(U_{y_i}'\cap V)}\ci\tau_{UV}(\hat s)=\rho_{U_{y_i}(U_{y_i}'\cap V)}(s_{y_i})$ for all~$i=1,\ldots,N$.

We have not yet shown that $\tau_{UV}(\hat s)$ is independent of the choices of $N,\ab y_i,\ab U_{y_i},\ab s_{y_i},\ab U_{y_i}'$ made above, so we do not yet know that $\tau_{UV}(\hat s)$ depends only on $\hat s$, but we take $\tau_{UV}(\hat s)$ and these choices to be fixed for the next part of the argument. Suppose now that as in Theorem \ref{kh2thm4}(c) that $y\in\bar V$, and $U_y$ is an open neighbourhood of $y$ in $U$, and $s_y\in\cE(U_y)$ with $\pi(U_y)(s_y)=\hat\rho_{UU_y}(\hat s)$. Choose any open neighbourhood $U_y'$ of $y$ in $U_y$ such that the closure $\bar U_y'$ of $U_y'$ in $U_y$ is compact. 

We can now apply the argument in which we constructed $\tau_{UV}(\hat s)$ above replacing the data $N,U_{y_1},\ldots,U_{y_N},s_{y_1},\ldots,s_{y_N},U'_{y_1},\ldots,U'_{y_N}$ by $N+1,\ab U_{y_1},\ab\ldots,\ab U_{y_N},\ab U_y,\ab s_{y_1},\ab\ldots,\ab s_{y_N},\ab s_y,\ab U'_{y_1},\ab\ldots,\ab U'_{y_N},\ab U'_y$. This gives $\tau_{UV}(\hat s)'$ in $\cE(V)$ with $\rho_{V(U_{y_i}'\cap V)}\ci\tau_{UV}(\hat s)'=\rho_{U_{y_i}(U_{y_i}'\cap V)}(s_{y_i})$ for $i=1,\ldots,N$ and $\rho_{U_y(U'_y\cap V)}(s_y)=\rho_{V(U'_y\cap V)}\ci\tau_{UV}(\hat s)'$. But uniqueness of $\tau_{UV}(\hat s)$ above implies that $\tau_{UV}(\hat s)'=\tau_{UV}(\hat s)$. Hence we have $\rho_{U_y(U'_y\cap V)}(s_y)=\rho_{V(U'_y\cap V)}\ci\tau_{UV}(\hat s)$ for all such $y,\ab U_y,\ab s_y,\ab U_y'$, so $\tau_{UV}(\hat s)$ above satisfies the characterizing property in Theorem~\ref{kh2thm4}(c).

Next we show that $\tau_{UV}(\hat s)$ is independent of the choices of $N,\ab y_i,\ab U_{y_i},\ab s_{y_i},\ab U_{y_i}'$ in its construction. Suppose $\check\tau_{UV}(\hat s)$ is an alternative outcome, using different choices. Then as $\check\tau_{UV}(\hat s)$ satisfies the characterizing property we have $\rho_{V(U_{y_i}'\cap V)}\ci\check\tau_{UV}(\hat s)=\rho_{U_{y_i}(U_{y_i}'\cap V)}(s_{y_i})$ for all $i=1,\ldots,N$. But $\tau_{UV}(\hat s)$ was unique with this property, so $\check\tau_{UV}(\hat s)=\tau_{UV}(\hat s)$, and $\tau_{UV}(\hat s)$ is independent of choices. Therefore $\tau_{UV}:\hat\cE(U)\ra\cE(V)$ mapping $\hat s\ra\tau_{UV}(\hat s)$ is well defined, and is uniquely characterized by the property in Theorem \ref{kh2thm4}(c), as required. It remains to show that $\tau_{UV}$ is an $R$-module morphism, and that \eq{kh2eq42}--\eq{kh2eq43} commute.

To see that $\tau_{UV}$ is an $R$-module morphism, let $a,b\in R$ and $\hat s,\hat t\in\hat\cE(U)$. Then in the construction of $\tau_{UV}(\hat s)$ above, for each $y\in\bar V$ we can choose $y\in U_y\subseteq U$ and $s_y,t_y\in\cE(U_y)$ such that $\pi(U_y)(s_y)=\hat\rho_{UU_y}(\hat s)$ and $\pi(U_y)(t_y)=\hat\rho_{UU_y}(\hat t)$. Following the argument through for $\tau_{UV}(\hat s),\tau_{UV}(\hat t)$ and $\tau_{UV}(a\hat s+b\hat t)$ simultaneously, uniqueness of $\tau_{UV}(\hat s)$ implies that $\tau_{UV}(a\hat s+b\hat t)=a\tau_{UV}(\hat s)+b\tau_{UV}(\hat t)$, so $\tau_{UV}$ is $R$-linear.

Suppose $V\ne\es$. If $s\in\cE(U)$ with $\pi(U)(s)=\hat s\in\hat\cE(U)$ then in the construction of $\tau_{UV}(\hat s)$ we can take $N=1$, any $y_1\in V\subseteq\bar V$, $U_{y_1}=U$, $s_{y_1}=s$ and $U_{y_1}'=V$, and then $\rho_{V(U_{y_1}'\cap V)}\ci\tau_{UV}(\hat s)=\rho_{U_{y_1}(U_{y_1}'\cap V)}(s_{y_1})$ becomes
\begin{equation*}
\tau_{UV}\ci\pi(U)(s)=\tau_{UV}(\hat s)=\rho_{VV}(\tau_{UV}(\hat s))=\rho_{UV}(s), 
\end{equation*}
so the top left triangle of \eq{kh2eq42} commutes. If $V=\es$ it commutes trivially.

If $\hat s\in\hat\cE(U)$ then from above there is an open cover $\bigl\{U_{y_1}'\cap V,\ldots,U_{y_N}'\cap V\bigr\}$ of $V$ such that for $i=1,\ldots,N$ we have
\begin{align*}
&\hat\rho_{V(U_{y_i}'\cap V)}\bigl[\pi(V)\ci\tau_{UV}(\hat s)\bigr]=\pi(U_{y_i}'\cap V)\ci\rho_{V(U_{y_i}'\cap V)}\ci\tau_{UV}(\hat s)\\
&\quad =\pi(U_{y_i}'\cap V)\ci
\rho_{U_{y_i}(U_{y_i}'\cap V)}(s_{y_i})=\hat\rho_{U_{y_i}(U_{y_i}'\cap V)}\ci\pi(U_{y_i})(s_{y_i})\\
&\quad =\hat\rho_{U_{y_i}(U_{y_i}'\cap V)}\ci\hat\rho_{UU_{y_i}}(\hat s)=\hat\rho_{U(U_{y_i}'\cap V)}(\hat s)=\hat\rho_{V(U_{y_i}'\cap V)}\bigl[\hat\rho_{UV}(\hat s)\bigr],
\end{align*}
using $\rho_{V(U_{y_i}'\cap V)}\ci\tau_{UV}(\hat s)=\rho_{U_{y_i}(U_{y_i}'\cap V)}(s_{y_i})$ from above, and $\pi_{U_i}(s_{y_i})=\hat\rho_{UU_i}(\hat s)$ by choice of $s_{y_i}$, and $\pi:\cE\ra\hat\cE$ a presheaf morphism. The sheaf property of $\hat\cE$ now implies that $\pi(V)\ci\tau_{UV}(\hat s)=\hat\rho_{UV}(\hat s)$, so the bottom right triangle of \eq{kh2eq42} commutes, and \eq{kh2eq42} commutes.

To see that equation \eq{kh2eq43} commutes, compare the characterizing properties of $\tau_{UV}(\hat s)$ and $\tau_{UW}(\hat s)$ at each $y\in\bar W$, and apply $\rho_{VW}$ to~$\tau_{UV}(\hat s)$.

\subsubsection{Proof of Theorem \ref{kh2thm4}(d)}
\label{kh624}

Let $U\subseteq Y$ be open, and suppose $V\subseteq U$ is open and the closure $\bar V$ of $V$ in $U$ is compact. Consider the diagram
\e
\begin{gathered}
\xymatrix@C=180pt@R=15pt{ *+[r]{\bigl\{\al\in \cE_\cs(U):\supp\al\subseteq V\subseteq U\bigr\}} \ar[r]^(0.6){\rho_{UV}\vert_{\cdots}}_(0.6)\cong \ar[d]^{\pi(U)\vert_{\cdots}} & *+[l]{\cE_\cs(V)} \ar[d]_{\pi(V)\vert_{\cdots}} \\ 
 *+[r]{\bigl\{\hat\al\in\hat\cE_\cs(U):\supp\hat\al\subseteq V\subseteq U\bigr\}} \ar[r]^(0.7){\hat\rho_{UV}\vert_{\cdots}}_(0.7)\cong \ar[ur]^(0.4){\tau_{UV}\vert_{\cdots}} & *+[l]{\hat\cE_\cs(V).\!} }
\end{gathered}
\label{kh6eq3}
\e
This is a subdiagram of \eq{kh2eq42} and so commutes by Theorem \ref{kh2thm4}(c), the top line is an isomorphism by Theorem \ref{kh2thm4}(b), the bottom line is an isomorphism as $\hat\cE$ is a sheaf, and the diagonal morphism does map to $\cE_\cs(V)\subseteq\cE(V)$ as \eq{kh2eq42} commutes. Hence the columns $\pi(U)\vert_{\cdots},\pi(V)\vert_{\cdots}$ in \eq{kh6eq3} are also isomorphisms, with inverses $\rho_{UV}\vert_{\cdots}^{-1}\ci\tau_{UV}\vert_{\cdots}$ and $\tau_{UV}\vert_{\cdots}\ci\hat\rho_{UV}\vert_{\cdots}^{-1}$.

As $Y$ is locally compact, for every compact subset $K\subseteq U$ we can find an open neighbourhood $V$ of $K$ in $U$ with the closure $\bar V$ of $V$ in $U$ compact. Applying this with $K=\supp\al$ for $\al\in\cE_\cs(U)$ or $K=\supp\hat\al$ for $\hat\al\in\hat\cE_\cs(U)$, we see that all $\al\in\cE_\cs(U)$ or $\hat\al\in\hat\cE_\cs(U)$ lie in the domain or target of the left hand morphism $\pi(U)\vert_{\cdots}$ in \eq{kh6eq3} for some such $V$. Thus, $\pi(U)\vert_{\cdots}$ an isomorphism for all such $V$ implies that $\pi(U):\cE_\cs(U)\ra\hat\cE_\cs(U)$ is an isomorphism.

\subsubsection{Proof of Theorem \ref{kh2thm4}(e)}
\label{kh625}

Let $U\subseteq Y$ be open. Consider the diagram:
\e
\begin{gathered}
\xymatrix@C=180pt@R=17pt{ *+[r]{\hat\cE(U)} \ar[r]^(0.3){\rm T} \ar[dr]_(0.4){\hat{\rm R}}^(0.4)\cong &
*+[l]{\mathop{\underleftarrow{\lim}\,}\nolimits_{\text{$V:V\subseteq U$ open, $\bar V$ is compact}}\cE(V)} \ar@<-1.3ex>[d]_\Pi \\
& *+[l]{\mathop{\underleftarrow{\lim}\,}\nolimits_{\text{$V:V\subseteq U$ open, $\bar V$ is compact}}\hat\cE(V).\!}  \ar@<.2ex>[u]_\Xi }
\end{gathered}
\label{kh6eq4}
\e
Here the inverse limits are over open $V\subseteq U$ for which the closure $\bar V$ of $V$ in $U$ is compact, defined using the morphisms $\rho_{V_1V_2}:\cE(V_1)\ra\cE(V_2)$ and $\hat\rho_{V_1V_2}:\hat\cE(V_1)\ra\hat\cE(V_2)$ for~$V_2\subseteq V_1\subseteq U$. 

The morphisms $\rm T$ and $\hat{\rm R}$ are induced by the morphisms $\tau_{UV}:\hat\cE(U)\ra\cE(V)$ from Theorem \ref{kh2thm4}(c) and $\hat\rho_{UV}:\hat\cE(U)\ra\hat\cE(V)$ from the (pre)sheaf $\hat\cE$, using the universal properties of inverse limits, and are well defined as $\rho_{V_1V_2}\ci\tau_{UV_1}=\tau_{UV_2}$ by \eq{kh2eq43} and $\hat\rho_{V_1V_2}\ci\hat\rho_{UV_1}=\hat\rho_{UV_2}$. Since $\hat\cE$ is a sheaf and $Y$ is locally compact, so that $U$ is the union of all open $V\subseteq U$ with $\bar V$ compact, we can show that $\hat{\rm R}$ is an isomorphism.

The morphism $\Pi$ in \eq{kh6eq4} is induced by $\pi(V):\cE(V)\ra\hat\cE(V)$ for $V\subseteq U$ with $\bar V$ compact, using the universal properties of inverse limits, and is well defined as $\hat\rho_{V_1V_2}\ci\pi(V_1)=\pi(V_2)\ci\rho_{V_1V_2}$ since $\pi:\cE\ra\hat\cE$ is a presheaf morphism. As $\pi(V)\ci\tau_{UV}=\hat\rho_{UV}$ by \eq{kh2eq42}, where $\pi(V),\tau_{UV},\hat\rho_{UV}$ are used to define $\Pi,{\rm T},\hat{\rm R}$, we have $\Pi\ci{\rm T}=\hat{\rm R}$ in~\eq{kh6eq4}.

To define $\Xi$, note that for any $V\subseteq U$ with $\bar V$ compact, as $Y$ is locally compact we can choose an open neighbourhood $V'$ of $\bar V$ in $U$ whose closure $\bar V'$ in $U$ is compact. Consider the following diagram:
\e
\begin{gathered}
\xymatrix@C=150pt@R=15pt{
*+[r]{\mathop{\underleftarrow{\lim}\,}\nolimits_{\text{$V':V'\subseteq U$ open, $\bar V'$ is compact}}\hat\cE(V')} \ar[r]_(0.8){\hat\pi_{V'}} \ar@{.>}[d]^\Xi &
*+[l]{\hat\cE(V')} \ar[d]_{\tau_{V'V}} \\
*+[r]{\mathop{\underleftarrow{\lim}\,}\nolimits_{\text{$V:V\subseteq U$ open, $\bar V$ is compact}}\cE(V)} \ar[r]^(0.8){\pi_V} & *+[l]{\cE(V).\!} 
 }
\end{gathered}
\label{kh6eq5}
\e
Here $\tau_{V'V}$ is from Theorem \ref{kh2thm4}(c), and $\pi_V,\hat\pi_{V'}$ are the projections from the inverse limits. Thus, for each $V\subseteq U$ with $\bar V$ compact we get a morphism $\tau_{V'V}\ci\hat\pi_{V'}:{\underleftarrow{\lim}}_{V'}\hat\cE(V')\ra\cE(V)$. We can show using \eq{kh2eq42}--\eq{kh2eq43} that this is independent of the choice of $V'$ and compatible with the morphisms $\rho_{V_1V_2}:\cE(V_1)\ra\cE(V_2)$,  so by properties of the inverse limit ${\underleftarrow{\lim}}_V\cE(V)$ there is a unique morphism $\Xi$ such that \eq{kh6eq5} commutes for all such~$V,V'$.

From \eq{kh2eq42} with $V'$ in place of $U$ we have $\tau_{V'V}\ci\pi(V')=\rho_{V'V}$, where $\tau_{V'V},\pi(V')$ are used to define $\Xi,\Pi$ in \eq{kh6eq4}, and $\rho_{V'V}$ is used to define the inverse limit ${\underleftarrow{\lim}}_V\cE(V)$. Hence we see that $\Xi\ci\Pi=\id$. Similarly, from \eq{kh2eq42} we have $\pi(V)\ci\tau_{V'V}=\hat\rho_{V'V}$, which implies that $\Pi\ci\Xi=\id$, so $\Xi,\Pi$ are inverse. Therefore $\rm T$ in \eq{kh6eq4} is an isomorphism. This is the isomorphism \eq{kh2eq44} induced by the $\tau_{UV}:\hat\cE(U)\ra\cE(V)$, as we have to prove.

\subsubsection{Proof of Theorem \ref{kh2thm4}(f)}
\label{kh626}

The definition of $\cE$ c-soft in Definition \ref{kh2def8}(c) may be expressed like this: $\cE$ is c-soft if whenever $K\subseteq Y$ is compact, $V$ is an open neighbourhood of $K$ in $Y$, and $s\in\cE(V)$, then there exists an open neighbourhood $W$ of $K$ in $V$ and $s'\in\cE(Y)$ with $\rho_{VW}(s)=\rho_{YW}(s')$. This definition only needs $\cE$ to be a presheaf, not a sheaf, and so makes sense for both $\cE$ and~$\hat\cE$.

Suppose $\cE$ is c-soft. We will prove that $\hat\cE$ is c-soft. Let $K\subseteq Y$ be compact, $U$ be an open neighbourhood of $K$ in $Y$, and $\hat s\in\hat\cE(U)$. Choose an open neighbourhood $V$ of $K$ in $U$ such that the closure $\bar V$ of $V$ in $U$ is compact, which is possible as $K$ is compact and $Y$ is locally compact. Set $s=\tau_{UV}(\hat s)\in\cE(V)$, for $\tau_{UV}$ as in Theorem \ref{kh2thm4}(c), so that $\pi(V)(s)=\hat\rho_{UV}(\hat s)$. Since $\cE$ is c-soft, there exists an open neighbourhood $W$ of $K$ in $V$ and $s'\in\cE(Y)$ with $\rho_{VW}(s)=\rho_{YW}(s')$. Write $\hat s'=\pi(Y)(s')$. Then
\begin{align*}
\hat\rho_{UW}(\hat s)&=\hat\rho_{VW}\ci\hat\rho_{UV}(\hat s)=\hat\rho_{VW}\ci\pi(V)(s)=\pi(W)\ci\rho_{VW}(s)\\
&=\pi(W)\ci\rho_{YW}(s')=\hat\rho_{YW}\ci\pi(Y)(s')=\hat\rho_{YW}(\hat s'),
\end{align*}
using $\pi(V)(s)=\hat\rho_{UV}(\hat s)$, $\hat s'=\pi(Y)(s')$ and $\pi$ a presheaf morphism.

Thus, whenever $K\subseteq Y$ is compact, $U$ is an open neighbourhood of $K$ in $Y$, and $\hat s\in\hat\cE(U)$, there exists an open neighbourhood $W$ of $K$ in $V$ and $\hat s'\in\hat\cE(Y)$ with $\hat\rho_{VW}(\hat s)=\hat\rho_{YW}(\hat s')$, so $\hat\cE$ is c-soft, as we have to prove.

The last part of (f) follows from Theorems \ref{kh2thm3} and \ref{kh2thm4}(b),(d). Since $\hat\cE$ is a c-soft sheaf, as in Theorem \ref{kh2thm3} we may reconstruct $\hat\cE$ from a flabby cosheaf $\hat\ucE$ defined by $\hat\ucE(U)=\hat\cE_\cs(U)$ for open $U\subseteq Y$,
and $\hat\si_{VU}:\hat\ucE(V)\ra\hat\ucE(U)$ for $V\subseteq U\subseteq Y$ open is the inverse of the isomorphism
\e
\hat\rho_{UV}\vert_{\cdots}:\bigl\{\hat\al\in \hat\cE_\cs(U):\supp\hat\al\subseteq V\subseteq U\bigr\}\,{\buildrel\cong\over\longra}\,\hat\cE_\cs(V).
\label{kh6eq6}
\e
But Theorem \ref{kh2thm4}(d) implies that $\ucE_\cs(U)\cong\hat\ucE_\cs(U)$, and equation \eq{kh2eq41} in Theorem \ref{kh2thm4}(b) is lift of \eq{kh6eq6} from $\hat\ucE_\cs(U)$ to $\ucE_\cs(U)$ under the isomorphisms $\ucE_\cs(U)\cong\hat\ucE_\cs(U)$. Therefore we may define a flabby cosheaf $\ucE$ defined by $\ucE(U)=\cE_\cs(U)$ and $\si_{VU}:\ucE(V)\ra\ucE(U)$ for $V\subseteq U\subseteq Y$ open is the inverse of the isomorphism \eq{kh2eq41}, and we have an isomorphism $\ucE\cong\hat\ucE$, so $\hat\cE$ is the c-soft sheaf associated to $\ucE$ in Theorem \rm\ref{kh2thm3}(b). This completes the proof of Theorem~\ref{kh2thm4}.

\subsection{Proof of Theorem \ref{kh2thm6}}
\label{kh63}

Let $X$ be an effective orbifold of dimension $m$. Then $X$ has a natural locally closed stratification
\begin{equation*}
X=\coprod\nolimits_{\text{isomorphism classes of finite groups $\Ga$}}X_\Ga,
\end{equation*}
where $X_\Ga=\bigl\{x\in X:G_xX\cong\Ga\bigr\}$ is the {\it orbifold stratum\/} of $X$ with group $\Ga$. Then $X_{\{1\}}$ is open and dense in $X$ and a manifold, and each $X_\Ga$ for $\Ga\ne\{1\}$ is a disjoint union $X_\Ga=\coprod_{k=0}^{m-1}X_\Ga^k$ of manifolds $X_\Ga^k$ of dimension $k<m$, such that $X$ is locally modelled near each point of $X_\Ga^k$ on $\R^k\t(\R^{m-k}/\Ga)$ for some effective representation of $\Ga$ on $\R^{m-k}$. The closure $\,\ov{\!X}_\Ga$ of $X_\Ga$ in $X$ satisfies
\begin{equation*}
\,\ov{\!X}_\Ga\subseteq \coprod\nolimits_{\begin{subarray}{l}\text{isomorphism classes of finite groups $\De$:}\\ \text{$\Ga$ is isomorphic to a subgroup of $\De$}\end{subarray}}X_\De.
\end{equation*}

Let $N\gg 0$ (actually $N\ge 2m+1$ is sufficient), and let $f:X\ra\R^N$ be a generic smooth map. Then following well known arguments of Whitney on embeddings of manifolds in $\R^N$, we can show that:
\begin{itemize}
\setlength{\itemsep}{0pt}
\setlength{\parsep}{0pt}
\item[(i)] $f$ is injective, and a homeomorphism with its image $f(X)\subset\R^N$.
\item[(ii)] $f\vert_{X_\Ga^k}:X_\Ga^k\ra\R^N$ is an embedding of manifolds for all~$\Ga,k$. 
\item[(iii)] $f$ is nicely behaved in the normal directions to $X_\Ga^k$ in $X$ for all $\Ga,k$. That is, near each point of $X_\Ga^k$, $X$ looks like $\R^k\t(\R^{m-k}/\Ga)$, and up to local diffeomorphisms of $\R^N$, $f$ looks like $\id_{\R^k}\t C:\R^k\t(\R^{m-k}/\Ga)\ra\R^k\t\R^{N-k}$ for some quasi-homogeneous embedding $C:\R^{m-k}/\Ga\ra \R^{N-k}$ of $\R^{m-k}/\Ga$, roughly as a cone in~$\R^{N-k}$.
\item[(iv)] $f(X)$ is a Euclidean Neighbourhood Retract (ENR), that is, there exists an open neighbourhood $Y$ of $f(X)$ in $\R^N$ which retracts topologically onto $X$. Roughly speaking, $Y$ is a `tubular neighbourhood' of $f(X)$ in $\R^N$, except that the retraction $Y\ra X$ is not a disc fibration over the orbifold strata $X_\Ga$, but has more complicated contractible fibres.
\end{itemize}

Starting with the effective orbifold $X$, we have constructed a manifold $Y$ and a smooth map $f:X\ra Y$ such that $f:X\ra f(X)$ is a homeomorphism and $Y$ retracts onto $X$. Hence $f$ is a topological homotopy, so there exist a continuous homotopy inverse $g:Y\ra X$, and continuous maps $F:[0,1]\t X\ra X$ and $G:[0,1]\t Y\ra Y$ with $F(0,x)=g\ci f(x),$ $F(1,x)=x,$ $G(0,y)=f\ci g(y)$ and $G(1,y)=y$ for all $x\in X$ and~$y\in Y$.

The key point is to show that we can choose these $g,F,G$ to be smooth, so that $X$ and $Y$ are homotopic in $\Orbeff$, rather than just in topological spaces $\Top$. When $X$ is a manifold this is well known, and follows from the construction of tubular neighbourhoods of embedded submanifolds.

To make the homotopy inverse $g$ smooth, we choose $g$ near $f(X_\Ga^k)$ in $Y$ for all $\Ga$ by induction on increasing $k=0,1,\ldots,m$, that is, we choose $g$ on the orbifold strata of highest codimension first, with the property that $g:Y\ra X$ maps a small open neighbourhood of $\ov{f(X_\Ga^k)}$ in $Y$ to $\ov{X_\Ga^k}$ in $X$ for all $\Ga,k$. 

That is, near points $x$ in $X_\Ga^k$ and $f(x)$ in $Y$, we know that $X$ is locally modelled on $\R^k\t(\R^{m-k}/\Ga)$, and $Y$ locally modelled on $\R^k\t\R^{N-k}$, and $f$ on $\id_{\R^k}\t C:\R^k\t(\R^{m-k}/\Ga)\ra\R^k\t\R^{N-k}$ for a quasi-homogeneous embedding $C:\R^{m-k}/\Ga\ra\R^{N-k}$. If $x$ is not too close to a deeper orbifold stratum $X_\De^l$ for $l<k$, we can suppose $Y$ is locally modelled on $\R^k\t U$, for $U$ an open neighbourhood of $C(\R^{m-k}/\Ga)$ in $\R^{N-k}$ which retracts onto $C(\R^{m-k}/\Ga)$. Then we want $g:Y\ra X$ to be locally modelled on $\id_{\R^k}\t D:\R^k\t U\ra\R^k\t(\R^{m-k}/\Ga)$, where $D:U\ra\R^{m-k}/\Ga$ is a smooth homotopy inverse for $C:\R^{m-k}/\Ga\ra U\subseteq\R^{N-k}$. We want in particular that $D$ should map an open neighbourhood of $0=C(0/\Ga)$ in $U\subset\R^{N-k}$ to $0/\Ga$ in~$\R^{m-k}/\Ga$.

It is not difficult to see that, by mapping open neighbourhoods in $U$ of images under $C$ of orbifold strata in $\R^{m-k}/\Ga$ back to the same orbifold strata (or their closures), such smooth homotopy inverses $D:U\ra\R^{m-k}/\Ga$ exist. Then by an inductive process of choosing $g$ near orbifold strata of increasing dimension, and making $Y$ smaller if necessary, we find we can choose the homotopy inverse $g$ to be smooth, that is, a morphism $g:Y\ra X$ in~$\Orbeff$.

Having chosen $Y,f$ and $g$, then $f\ci g:Y\ra Y$ and $\id_Y:Y\ra Y$ are smooth maps of manifolds which are continuously homotopic, that is, there exists a continuous map $\ti G:[0,1]\t Y\ra Y$ with $\ti G(0,y)=f\ci g(y)$ and $\ti G(1,y)=y$ for all $y\in Y$. As smooth maps are dense in continuous maps, we can choose a small smooth perturbation $G:[0,1]\t Y\ra Y$ of $\ti G$ with $G(0,y)=f\ci g(y)$ and $G(1,y)=y$ for all $y\in Y$.

Similarly, there exists a continuous map $\ti F:[0,1]\t X\ra X$ with $\ti F(0,x)=g\ci f(x)$ and $\ti F(1,x)=x$ for all $x\in X$. For our definition of $\Orbeff$, by considering equivariant continuous and smooth maps we see that it is still true that smooth maps are dense in continuous maps. So again we can choose a small smooth perturbation $F:[0,1]\t X\ra X$ of $\ti F$ with $F(0,x)=g\ci f(x)$ and $F(1,x)=x$ for all $x\in X$. This completes the proof of Theorem~\ref{kh2thm6}.

Note that for the above proof to work, it is essential that morphisms in $\Orbeff$ from \S\ref{kh291} are continuous maps, rather than continuous maps plus extra data, as in the (2-)categories of orbifolds discussed in Remark~\ref{kh2rem5}(c).

\subsection{Proof of Theorem \ref{kh2thm7}}
\label{kh64}

For part (a), let $H_*(-;R),\ti H_*(-;R)$ be homology theories of effective orbifolds over $R$, in the sense of Definition \ref{kh2def15}. Restricting these to $\Man\subset\Orbeff$ gives homology theories of manifolds over $R$. Thus Theorem \ref{kh2thm1} gives canonical isomorphisms $I_{Y,Z}:H_*(Y,Z;R)\ra\ti H_*(Y,Z;R)$ for all {\it manifolds\/} $Y$ and open $Z\subseteq Y$ commuting with the given morphisms $f_*,\pd$ for manifolds and the given isomorphisms $H_0(*;R)\cong R\cong\ti H_0(*;R)$, and any other such morphisms $J_{Y,Z}:H_*(Y,Z;R)\ra\ti H_*(Y,Z;R)$ for manifolds $Y$ have~$J_{Y,Z}=I_{Y,Z}$.

Let $Y$ be an effective orbifold. Theorem \ref{kh2thm6} shows that we can construct a manifold $Y'$ and smooth maps $f:Y\ra Y'$, $g:Y'\ra Y,$ $F:[0,1]\t Y\ra Y$ and $G:[0,1]\t Y'\ra Y'$ with $F(0,y)=g\ci f(y)$, $F(1,y)=y,$ $G(0,y')=f\ci g(y')$ and $G(1,y')=y'$ for all $y\in Y$ and $y'\in Y'$. Thus for each $k\in\Z$ we have a diagram
\e
\begin{gathered}
\xymatrix@C=100pt@R=15pt{ *+[r]{H_k(Y;R)} \ar@<.5ex>[r]^{f_*} \ar@{.>}[d]^{I_{Y,\es}} & *+[l]{H_k(Y';R)} \ar@<.5ex>[l]^{g_*} \ar[d]_{I_{Y',\es}}^\cong \\
*+[r]{\ti H_k(Y;R)} \ar@<.5ex>[r]^{f_*} & *+[l]{\ti H_k(Y';R).\!} \ar@<.5ex>[l]^{g_*} }
\end{gathered}
\label{kh6eq7}
\e
Since $g\ci f:Y\ra Y$ and $\id_Y:Y\ra Y$ are homotopic by $F$, Axiom \ref{kh2ax1}(i),(iv) imply that $g_*\ci f_*=(g\ci f)_*=(\id_Y)_*=\id$ on both $H_k(Y;R)$ and $\ti H_k(Y;R)$. Similarly, using $G$ we have $f_*\ci g_*=\id$ on both $H_k(Y';R)$ and $\ti H_k(Y';R)$. So in the rows of \eq{kh6eq7}, the morphisms $f_*,g_*$ are inverse, and both are isomorphisms.

Therefore there exists a unique isomorphism $I_{Y,\es}:H_k(Y;R)\ra \ti H_k(Y;R)$ making \eq{kh6eq7} commute. If $\dot Y',\dot f,\dot g,\dot F,\dot G$ are alternative choices for $Y',f,g,F,G$, then by considering the diagram
\e
\begin{gathered}
\xymatrix@!0@C=90pt@R=15pt{ 
&& H_k(Y';R) \ar[ddd]_(0.3){I_{Y',\es}} \ar@<.5ex>[ddr]^{(\dot f\ci g)_*} \ar@<.5ex>[dll]^(0.4){g_*} \\
*+[r]{H_k(Y;R)} \ar[ddd]^{I_{Y,\es}} \ar@<.5ex>[urr]^{f_*} \ar@<.5ex>[drrr]^(0.53){\dot f_*} 
\\
&&& H_k(\dot Y';R) \ar[ddd]^{I_{\dot Y',\es}} \ar@<.5ex>[ulll]^{\dot g_*} \ar@<.5ex>[uul]^(0.7){(f\ci\dot g)_*}\\
&& \ti H_k(Y';R) \ar@<.5ex>[ddr]^{(\dot f\ci g)_*} \ar@<.5ex>[dll]^(0.4){g_*} \\
*+[r]{\ti H_k(Y;R)}  \ar@<.5ex>[urr]^{f_*} \ar@<.5ex>[drrr]^(0.53){\dot f_*}
\\
&&& \ti H_k(\dot Y';R),\! \ar@<.5ex>[ulll]^{\dot g_*} \ar@<.5ex>[uul]^(0.7){(f\ci\dot g)_*} }
\end{gathered}
\label{kh6eq8}
\e
we see that $I_{Y,\es}$ is independent of the choices of~$Y',f,g,F,G$.

Next let $Z\subseteq Y$ be an open suborbifold. From the proof of Theorem \ref{kh2thm6} in \S\ref{kh63}, we see that $f:Y\hookra Y'$ is a topological embedding, and we can choose an open submanifold $Z'\subseteq Y'$ such that $f\vert_Z:Z\ra Z'$ is a homotopy, with smooth homotopy inverse $\check g:Z'\ra Z$ and smooth maps $\check F:[0,1]\t Z\ra Z$ and $\check G:[0,1]\t Z'\ra Z'$ with $\check F(0,z)=\check g\ci f(z)$, $\check F(1,z)=z,$ $\check G(0,z')=f\ci\check g(z')$ and $\check G(1,z')=z'$ for all $z\in Z$ and~$z'\in Z'$. 

Note that we do {\it not\/} claim that we can choose $\check g=g\vert_{Z'}$, $\check F=F\vert_{[0,1]\t Z}$ and $\check G=G\vert_{[0,1]\t Z'}$, that is, we are not constructing a homotopy of pairs $f:(Y,Z)\ra (Y',Z')$. This is possible if $(Y,Z)$ is a {\it good\/} orbifold pair in the sense of \S\ref{kh61}, but may not be possible for general $(Y,Z)$. Instead, $f:(Y,Z)\ra (Y',Z')$ is a {\it pseudo-equivalence}, in the sense of \S\ref{kh61}.

Now consider the commutative diagram with exact rows
\begin{equation*}
\xymatrix@C=12pt@R=15pt{ \cdots \ar[r] & H_k(Z;R) \ar[d]_\cong^{(f\vert_Z)_*}\ar[r]_(0.47){i_*} & H_k(Y;R) \ar[d]_\cong^{f_*} \ar[r]_(0.45){j_*} & H_k(Y,Z;R) \ar[d]^{f_*} \ar[r]_(0.48)\pd & H_{k-1}(Z;R) \ar[d]_\cong^{(f\vert_Z)_*} \ar[r] & \cdots \\
\cdots \ar[r] & H_k(Z';R) \ar[r]^(0.47){i_*} & H_k(Y';R) \ar[r]^(0.45){j_*} & H_k(Y',Z';R) \ar[r]^(0.48)\pd & H_{k-1}(Z';R) \ar[r] & \cdots. }\!\!\!\!\!\!\!\!\!\!\!\!\!\!\!{}
\end{equation*}
The first, second and fourth columns are isomorphisms by the argument above, since $f:Y\ra Y'$ and $f\vert_Z:Z\ra Z'$ are homotopies in $\Orbeff$. The five lemma now implies that $f_*:H_k(Y,Z;R)\ra H_k(Y',Z';R)$ is an isomorphism, and similarly $f_*:\ti H_k(Y,Z;R)\ra\ti H_k(Y',Z';R)$ is an isomorphism. Define an isomorphism $I_{Y,Z}:H_k(Y,Z;R)\ra\ti H_k(Y,Z;R)$ by the commutative diagram
\begin{equation*}
\xymatrix@C=100pt@R=15pt{ *+[r]{H_k(Y,Z;R)} \ar[r]^{f_*}_\cong \ar@{.>}[d]^{I_{Y,Z}} & *+[l]{H_k(Y',Z';R)} \ar[d]_{I_{Y',Z'}}^\cong \\
*+[r]{\ti H_k(Y,Z;R)} \ar[r]^{f_*}_\cong & *+[l]{\ti H_k(Y',Z';R).\!} }
\end{equation*}
An argument similar to \eq{kh6eq8} shows that $I_{Y,Z}$ is independent of choices. 

By passing through the corresponding statements for manifolds using homotopies as above, we can show that these isomorphisms $I_{Y,Z}$ for all effective orbifolds $Y$ and open $Z\subseteq Y$ commute with the given morphisms $f_*,\pd$ in $H_*(-;R),\ti H_*(-;R)$ and the isomorphisms $H_0(*;R)\cong R\cong\ti H_0(*;R)$, and any other assignment of morphisms $J_{Y,Z}$ preserving all the structure satisfy $J_{Y,Z}=I_{Y,Z}$. This proves Theorem \ref{kh2thm7}(a). The proof of (b) is as for (a), except that we reverse directions of most of the morphisms $f_*,g_*,\pd,\ldots$ in the proof, to get~$f^*,g^*,\d,\ldots.$ 

\section{Proofs of results in \S\ref{kh4}}
\label{kh7}

\subsection{Proof of Proposition \ref{kh4prop1}}
\label{kh71}

Let $Y$ be a manifold. To show that $\pd:MC_k(Y;R)\ra MC_{k-1}(Y;R)$ in Definition \ref{kh4def1} is well defined, we have to prove that it takes relations Definition \ref{kh4def1}(i),(ii) in $MC_k(Y;R)$ to relations Definition \ref{kh4def1}(i),(ii) in $MC_{k-1}(Y;R)$. This is obvious for relation~(i).

Suppose $\sum_{i\in I}a_i\,[V_i,n,s_i,t_i]=0$ in $MC_k(Y;R)$ by relation (ii). Then there exists an open neighbourhood $X$ of $0$ in $\R^n$, such that $s_i:V_i\ra\R^n$ is proper over $X$ for all $i\in I$, and condition Definition \ref{kh4def1}(ii)$(*)$ holds. We will show that the conditions of Definition \ref{kh4def1}(ii) hold for $\sum_{i\in I}a_i\,[\pd V_i,n,s_i\ci i_{V_i},t_i\ci i_{V_i}]$ in $MC_{k-1}(Y;R)$ with the same open neighbourhood $X$, so that $\sum_{i\in I}a_i\,[\pd V_i,n,s_i\ci i_{V_i},t_i\ci i_{V_i}]=0$ by Definition \ref{kh4def1}(ii). This will prove that $\pd$ takes relation (ii) to relation (ii), and so is well defined.

The first part of (ii), that $s_i\ci i_{V_i}:\pd V_i\ra\R^n$ is proper over $X$ for all $i\in I$, holds as $i_{V_i}:\pd V_i\ra V_i$ is proper and $s_i:V_i\ra\R^n$ is proper over $X$. Suppose for a contradiction that condition $(*)$ does not hold for $\sum_{i\in I}a_i\,[\pd V_i,n,s_i\ci i_{V_i},t_i\ci i_{V_i}]$ in $MC_{k-1}(Y;R)$ with the given $X\subseteq\R^n$. Then there exists $(x,y)\in X\t Y$, such that for all $i\in I$ and $v'\in\pd V_i$ with $(s_i\ci i_{V_i},t_i\ci i_{V_i})(v')=(x,y)$, we have that $v'\in(\pd V_i)^\ci$ and 
\e
T_v(s_i\ci i_{V_i},t_i\ci i_{V_i}):T_v(\pd V_i)^\ci\longra T_xX\op T_yY
\label{kh7eq1}
\e
is injective, but there exists an oriented $(n+k-1)$-plane $P\subseteq T_xX\op T_yY$ such that \eq{kh4eq2} does not hold, that is,
\e
\begin{split}
&\sum_{\begin{subarray}{l} i\in I,\; v'\in(\pd V_i)^\ci:(s_i\ci i_{V_i},t_i\ci i_{V_i})(v')=(x,y),\;  T_{v'}(s_i\ci i_{V_i},t_i\ci i_{V_i})[T_{v'}(\pd V_i)^\ci]=P \\ \text{$T_{v'}(s_i\ci i_{V_i},t_i\ci i_{V_i}):T_{v'}(\pd V_i)^\ci\,{\buildrel\cong\over\longra}\,P$ is orientation-preserving}\end{subarray}} a_i\ne\\
&\sum_{\begin{subarray}{l} i\in I,\; v'\in(\pd V_i)^\ci:(s_i\ci i_{V_i},t_i\ci i_{V_i})(v')=(x,y),\;  T_{v'}(s_i\ci i_{V_i},t_i\ci i_{V_i})[T_{v'}(\pd V_i)^\ci]=P \\ \text{$T_{v'}(s_i\ci i_{V_i},t_i\ci i_{V_i}):T_{v'}(\pd V_i)^\ci\,{\buildrel\cong\over\longra}\,P$ is orientation-reversing}\end{subarray}\!\!\!\!\!\!\!} a_i\qquad\text{in $R$.}
\end{split}
\label{kh7eq2}
\e

As in Definition \ref{kh4def1}(ii), there are only finitely many points $v'$ in $\pd V_i$ for any $i\in I$ with $(s_i\ci i_{V_i},t_i\ci i_{V_i})(v')=(x,y)$. Write these points as $v'_1,\ldots,v'_N$ with $v_j'\in\pd V_{i_j}$ for $i_1,\ldots,i_N\in I$, where $v_1',\ldots,v_N'$ are distinct but $i_1,\ldots,i_N$ need not be distinct. Then $v_j'\in(\pd V_{i_j})^\ci$ and \eq{kh7eq1} is injective for $v_j'$, so $(s_{i_j}\ci i_{V_{i_j}},t_{i_j}\ci i_{V_{i_j}}):(\pd V_{i_j})^\ci\ra \R^n\t Y$ is an embedding of manifolds near $v_j'$ in $\pd V_j$. Suppose that $(x,y)\in X\t Y$ is chosen such that $N$ is minimal, such that \eq{kh7eq2} holds.

We may choose a small open neighbourhood $W$ of $(x,y)$ in $X\t Y\subseteq\R^n\t Y$ and open neighbourhoods $U_1,\ldots,U_N$ of $v_1',\ldots,v_N'$ in $(\pd V_{i_1})^\ci,\ldots,(\pd V_{i_N})^\ci$, such that $U_1,\ldots,U_N$ are disjoint, and $(s_{i_j}\ci i_{V_{i_j}},t_{i_j}\ci i_{V_{i_j}})\vert_{U_j}:U_j\ra\R^n\t Y$ is an embedding of manifolds, and if $i\in I$ and $v'\in\pd V_i$ then $(s_{i_j}\ci i_{V_{i_j}},t_{i_j}\ci i_{V_{i_j}})(v')$ lies in $W$ if and only if $i=i_j$ and $v'\in U_j$ for some unique~$j=1,\ldots,N$.

Write $\ti U_j=(s_{i_j}\ci i_{V_{i_j}},t_{i_j}\ci i_{V_{i_j}})[U_j]\subseteq W$. Then $\ti U_1,\ldots,\ti U_N$ are oriented embedded submanifolds of $W$ of dimension $n+k-1$, which are closed in $W$ as $(s_{i_j}\ci i_{V_{i_j}},t_{i_j}\ci i_{V_{i_j}})$ is proper over $W\subseteq X\t Y$, and $(x,y)\in\ti U_j$ for all $j$. Equation \eq{kh7eq2} may be rewritten
\e
\sum_{\begin{subarray}{l} \text{$j=1,\ldots,N$: $T_{(x,y)}\ti U_j=P$ and}\\
\text{$T_{(x,y)}\ti U_j,P$ have the same orientation} \end{subarray}\!\!\!\!\!\!\!} a_{i_j} \ne\sum_{\begin{subarray}{l} \text{$j=1,\ldots,N$: $T_{(x,y)}\ti U_j=P$ and}\\
\text{$T_{(x,y)}\ti U_j,P$ have the opposite orientation} \end{subarray}\!\!\!\!\!\!\!\!\!\!\!\!\!\!} a_{i_j}.
\label{kh7eq3}
\e

If the submanifolds $\ti U_1,\ldots,\ti U_N$ do not all coincide near $(x,y)$ in $W$, then we can find a point $(x',y')$ near $(x,y)$ in $W$ (which can be taken to be a generic point of $\ti U_j$ for some $j$) and an oriented $(n+k-1)$-plane $P'\subseteq T_{x'}X\op T_{y'}Y$ (which can be taken to be $T_{(x',y')}\ti U_j$) such that the analogue of \eq{kh7eq3} for $(x',y'),P'$ holds, but $(x',y')$ does not lie in all of $\ti U_1,\ldots,\ti U_N$. But then replacing $(x,y)$ by $(x',y')$ we can reduce $N$, contradicting the minimality of~$N$. 

Hence $\ti U_1,\ldots,\ti U_N$ all coincide near $(x,y)$, and making $W,U_j,\ti U_j$ smaller we can suppose that $\ti U_1=\cdots=\ti U_N$ with $P=T_{(x,y)}\ti U_j$, and also that the $\ti U_j$ are connected. Let the orientation of $\ti U_j$ be $\ep_j$ times the orientation of $\ti U_1$, where $\ep_1,\ldots,\ep_N\in\{\pm 1\}$ with $\ep_1=1$. Then \eq{kh7eq3} is equivalent to
\e
\sum_{j=1}^N\ep_ja_j\ne 0\qquad\text{in $R$.}
\label{kh7eq4}
\e

Making $W$ smaller if necessary, choose global coordinates $(w_1,\ldots,w_{k+n-1},\ab z_1,\ldots,z_l)$ on $W$, where $l=\dim Y-k+1$, such that $(x,y)=(0,0,\ldots,0)$, and
\begin{equation*}
\ti U_1=\cdots=\cdots=\ti U_N=\bigl\{(w_1,\ldots,w_{k+n-1},0,\ldots,0)\in W\bigr\},
\end{equation*} 
that is, $\ti U_j$ is the submanifold $z_1=\cdots=z_l=0$ in $W$, and suppose $\d w_1\w\cdots\w\d w_{k+n-1}$ is a positive form on the oriented manifold $\ti U_1$, so that $\ep_j\,\d w_1\w\cdots\w\d w_{k+n-1}$ is positive on $\ti U_j$ for $j=1,\ldots,N$.

For each $i\in I$, define $\ti V_i\subseteq V_i$ to be the open submanifold $(s_i,t_i)^{-1}(W)\subseteq V_i$, as an object of $\tManc$. Consider the morphisms in $\tManc$:
\ea
f_i:=(w_1,\ldots,w_{k+n-1})\ci (s_i,t_i)\vert_{\ti V_i}&:\ti V_i\longra\R^{n+k-1},
\label{kh7eq5}\\
g_i:=(w_1,\ldots,w_{k+n-1},z_1^2+\cdots+z_l^2)\ci (s_i,t_i)\vert_{\ti V_i}&:\ti V_i\longra\R^{n+k}.
\label{kh7eq6}
\ea
Since $\dim\ti V_i=n+k$, Assumption \ref{kh3ass7}(b) implies that there exists $S\subset\R^{n+k-1}$ with $\cH^{n+k-1}(S)=0$ such that \eq{kh7eq5} is a submersion on an open neighbourhood in $\ti V_i$ of the preimage of any $(w_1',\ldots,w_{n+k-1}')\in\R^{n+k-1}\sm S$ for all $i\in I$, and $\ti S\subset\R^{n+k}$ with $\cH^{n+k}(\ti S)=0$ with \eq{kh7eq6} a submersion on an open neighbourhood in $\ti V_i$ of the preimage of any $(w_1',\ldots,w_{n+k-1}',\de)\in\R^{n+k}\sm\ti S$ for all~$i\in I$.

Then $S\t\R\subset\R^{n+k}$ with $\cH^{n+k}(S\t\R)=0$, so $\cH^{n+k}\bigl((S\t\R)\t\ti S\bigr)=0$. Thus we may choose $(w_1',\ldots,w_{n+k-1}',\de)\in\R^{n+k}$ such that $w_1',\ldots,w_{n+k-1}',\de$ are very small, and $\de>0$, and $(w_1',\ldots,w_{n+k-1}')\notin S$, and $(w_1',\ldots,w_{n+k-1}',\de)\notin\ti S$. As $w_1',\ldots,w_{n+k-1}',\de$ are small we may assume that $(w_1',\ldots,w_{n+k-1}',z_1,\ldots,z_l)$ lies in $W$ for all $(z_1,\ldots,z_l)\in\R^l$ with~$z_1^2+\cdots+z_l^2\le\de$.

Define $h:(-\iy,0]\ra\R^{n+k}$ by
\begin{equation*}
h(x)=(w_1',\ldots,w_{n+k-1}',x+\de).
\end{equation*}
Then $h$ is a morphism in $\tManc$ by Assumption \ref{kh3ass4}(c). We claim that $g_i:\ti V_i\ra\R^{n+k}$ in  \eq{kh7eq6} and $h:(-\iy,0]\ra\R^{n+k}$ are transverse as morphisms in~$\tManc$. 

To prove this, suppose $v\in\ti V_i$ and $x\in(-\iy,0]$ with $g_i(v)=h(x)$. If $x=0$ then $g_i$ is a submersion in an open neighbourhood of $v$ as $g_i(v)=(w_1',\ldots,w_{n+k-1}',\de)$ in $\R^{n+k}\sm\ti S$, which implies that $g_i,h$ are transverse on open neighbourhoods of $v\in\ti V_i$ and $x\in(-\iy,0]$ by Assumption \ref{kh3ass5}(c). If $x\ne 0$, we split $\R^{n-k}=\R^{n-k-1}\t\R$ and write $g_i,h$ as direct products $g_i=(f_i,g_i')$ and $h=(h_1,h_2)$, where $f_i:\ti V_i\ra\R^{n+k-1}$ is a submersion near $v$ as $f_i(v)=(w_1',\ldots,w_{n+k-1}')$ in $\R^{n+k-1}\sm S$, and $h_2:(-\iy,0]\ra\R$, $h_2(x)=x+\de$ is a submersion on $(-\iy,0)\ni x$ by Assumption \ref{kh3ass5}(a). Hence $g_i,h$ are transverse on open neighbourhoods of $v\in\ti V_i$ and $x\in(-\iy,0]$ by Assumption \ref{kh3ass5}(f). This holds for all $v\in\ti V_i$ and $x\in(-\iy,0]$ with $g_i(v)=h(x)$, so $g_i,h$ are transverse by Assumption~\ref{kh3ass5}(e).

Therefore by Assumption \ref{kh3ass5}(c), the transverse fibre product
\e
T_i:=\ti V_i\t_{g_i,\R^{n+k},h}(-\iy,0]
\label{kh7eq7}
\e
exists in $\tManc$, with $\dim T_i=1$. Consider the diagram of topological spaces
\begin{equation*}
\xymatrix@C=120pt@R=15pt{
*+[r]{T_i} \ar[r]_(0.25){\pi_{\ti V_i}\t\pi_{(-\iy,0]}} & *+[l]{\bigl\{(v,x)\in\ti V_i\t(-\iy,0]:g_i(v)=h(x)\bigr\}} \ar[d]_{(s_i,t_i)\ci\pi_{\ti V_i}} \\
& *+[l]{\bigl\{(w_1',\ldots,w_{n+k-1}',z_1,\ldots,z_l):z_j\in\R,\;\> z_1^2+\cdots+z_l^2\le\de\bigr\}\subseteq W.} }
\end{equation*}
The top map is a homeomorphism by \eq{kh3eq8} in Assumption \ref{kh3ass5}(c). The right hand map is proper as $(s_i,t_i)$ is proper over $W\subseteq X\t Y$. The bottom space is a compact $l$-ball in $W$. Hence $T_i$ is compact.

Combining the given orientation on $\ti V_i\subseteq V_i$ with the standard orientations on $\R^{n+k},(-\iy,0]$ from Assumption \ref{kh3ass6}(k), by Assumption \ref{kh3ass6}(l) we have an orientation on $T_i$ in \eq{kh7eq7}. By Assumptions \ref{kh3ass6}(c), \ref{kh3ass6}(l),(m) and equations \eq{kh3eq2} and \eq{kh3eq9}, since $\pd(-\iy,0]=\{0\}$ in oriented objects in $\tManc$ we have
\e
\pd T_i\cong \bigl(\pd \ti V_i\t_{g_i\ci i_{\ti V_i},\R^{n+k},h}(-\iy,0]\bigr)\amalg \bigl(\ti V_i\t_{g_i,\R^{n+k},h\ci i_{(-\iy,0]}}\{0\}\bigr).
\label{kh7eq8}
\e

The first term on the r.h.s.\ of \eq{kh7eq8} is empty if $i\ne i_j$ for some $j=1,\ldots,N$, and if $i=i_j$ it is the transverse intersection of the hyperplane $\bigl\{(w_1,\ldots,w_{n+k-1},\ab 0):\ab w_j\in\R\bigr\}$ in $\R^{n+k}$, oriented so that $\ep_j\,\d w_1\w\cdots\w\d w_{k+n-1}$ is a positive form, and the ray $\bigl\{(w_1',\ldots,w_{n+k-1}',z):z\le\de\bigr\}$ in $\R^{n+k}$, oriented so that $\d z$ is a positive form. Thus, if $i=i_j$ then the first term is a single point $(w_1',\ldots,w_{n+k-1}',0)$, oriented with sign~$\ep_j$. 

The second term on the r.h.s.\ of \eq{kh7eq8} is the finite set of all points $v$ in $\ti V_i$ with $g_i(v)=(w_1',\ldots,w_{n+k-1}',\de)$, which is the same as the set of all points $v\in V_i$ with $(s_i,t_i)(v)=(w_1',\ldots,w_{n+k-1}',z_1,\ldots,z_l)\in W\subseteq X\t Y$ for some $(z_1,\ldots,z_l)\in\R^l$ with $z_1^2+\cdots+z_l^2=\de$. All such $v$ lie in $V_i^\ci$, with $\d g_i\vert_v:T_vV_i^\ci\ra\R^{n+k}$ an isomorphism. Now Assumption \ref{kh3ass6}(n) says that the number of points in \eq{kh7eq8}, counted with signs, is zero. Hence for each $i\in I$ we have
\ea
&\sum_{\begin{subarray}{l} (z_1,\ldots,z_l)\in\R^l: \\ z_1^2+\cdots+z_l^2=\de\end{subarray}}\,\,\sum_{\begin{subarray}{l} v\in V_i^\ci:(s_i,t_i)(v)= \\ (w_1',\ldots,w_{n+k-1}',z_1,\ldots,z_l)\end{subarray}\!\!\!\!\!\!\!\!\!\!\!\!\!\!\!}\begin{cases} 1, & \!\!\!\text{$\d g_i\vert_v:T_vV_i^\ci\!\ra\!\R^{n+k}$ orientation-reversing} \\
-1, & \!\!\!\text{$\d g_i\vert_v:T_vV_i^\ci\!\ra\!\R^{n+k}$ orientation-preserving} 
\end{cases}
\nonumber\\
&\qquad=\begin{cases} \ep_j, & i=i_j,\;\> j=1,\ldots,N, \\
0, & i\notin\{i_1,\ldots,i_N\}.
\end{cases}
\label{kh7eq9}
\ea

Multiply \eq{kh7eq9} by $a_i\in R$ and sum over all $i\in I$. The r.h.s.\ of the resulting equation is nonzero by \eq{kh7eq4}. But the l.h.s.\ is the sum over all points $(x,y)=(w_1',\ldots,w_{n+k-1}',z_1,\ldots,z_l)\in W\subseteq X\t Y$ for $(z_1,\ldots,z_l)\in\R^l$ with $z_1^2+\cdots+z_l^2=\de$, and over all $(n+k)$-planes $P\subseteq T_xX\op T_yY$ such that $\d(w_1,\ldots,w_{k+n-1},z_1^2+\cdots+z_l^2)\vert_{(x,y)}:T_xX\op T_yY\ra\R^{n+k}$ is an isomorphism, with $P$ oriented so that $\d(w_1,\ldots,w_{k+n-1},z_1^2+\cdots+z_l^2)\vert_{(x,y)}$ is orientation-reversing, of equation \eq{kh4eq2} for $\sum_{i\in I}a_i\,[V_i,n,s_i,t_i]$ at $(x,y),P$. Hence the l.h.s.\ of the resulting equation is zero, a contradiction.  

Therefore Definition \ref{kh4def1}(ii)$(*)$ holds for $\sum_{i\in I}a_i\,[\pd V_i,n,s_i\ci i_{V_i},t_i\ci i_{V_i}]$ in $MC_{k-1}(Y;R)$, so $\sum_{i\in I}a_i\,[\pd V_i,n,s_i\ci i_{V_i},t_i\ci i_{V_i}]=0$ by Definition \ref{kh4def1}(ii), and $\pd:MC_k(Y;R)\ra MC_{k-1}(Y;R)$ in Definition \ref{kh4def1} is well defined. This proves Proposition~\ref{kh4prop1}. 

\subsection{Proof of Theorem \ref{kh4thm1}}
\label{kh72}

We begin with some new notation.

\begin{dfn} Let $Y$ be a manifold and $T,U\subseteq Y$ be open, and write $i:T\hookra T\cup U$, $j:U\hookra T\cup U$ for the inclusions. Suppose $f:T\cup U\ra\R$ is smooth, with $\bigl\{y\in T\cup U:f(y)\ge 0\bigr\}\subseteq T$ and $\bigl\{y\in T\cup U:f(y)\le 0\bigr\}\subseteq U$. Such a function $f$ always exists; for instance, if $(\eta,1-\eta)$ is a partition of unity on $T\cup U$ subordinate to the open cover $(T,U)$, so that $\eta:T\cup U\ra\R$ is smooth with $\eta=0$ on $U\sm T$ and $\eta=1$ on $T\sm U$, then $f=\eta-\ha$ will do.

Define $R$-linear morphisms $\Pi_{T\cup U}^{T,f-},\Pi_{T\cup U}^{U,f+}$ for each $k\in\Z$ by
\e
\begin{split}
&\Pi_{T\cup U}^{T,f-}:MC_k(T\cup U;R)\longra MC_k(T;R),\\
&\Pi_{T\cup U}^{T,f-}:\bigl[V,n,(s_1,\ldots,s_n),t\bigr]\longmapsto \bigl[t^{-1}(T)\t(-\iy,0],n+1,\\
&\quad (s_1\ci\pi_{t^{-1}(T)},\ldots,s_n\ci\pi_{t^{-1}(T)},f\ci\pi_{t^{-1}(T)}+\pi_{(-\iy,0]}),t\ci\pi_{t^{-1}(T)}\bigr],\\
&\Pi_{T\cup U}^{U,f+}:MC_k(T\cup U;R)\longra MC_k(U;R),
\\
&\Pi_{T\cup U}^{U,f+}:\bigl[V,n,(s_1,\ldots,s_n),t\bigr]\longmapsto \bigl[t^{-1}(U)\t[0,\iy),n+1,\\
&\quad (s_1\ci\pi_{t^{-1}(U)},\ldots,s_n\ci\pi_{t^{-1}(U)},f\ci\pi_{t^{-1}(U)}+\pi_{[0,\iy)}),t\ci\pi_{t^{-1}(U)}\bigr],
\end{split}
\label{kh7eq10}
\e
for each generator $\bigl[V,n,(s_1,\ldots,s_n),t\bigr]$ in $MC_k(T\cup U;R)$.
\label{kh7def1}
\end{dfn}

\begin{prop}{\bf(i)} In the above, $\Pi_{T\cup U}^{T,f-}$ and\/ $\Pi_{T\cup U}^{U,f+}$ are well defined.
\smallskip

\noindent{\bf(ii)} We have $i_*\ci\Pi_{T\cup U}^{T,f-}+j_*\ci\Pi_{T\cup U}^{U,f+}=\id:MC_k(T\cup U;R)\ra MC_k(T\cup U;R)$.
\smallskip

\noindent{\bf(iii)} If\/ $[V,n,s,t]$ is a generator of\/ $MC_k(T;R)$ with\/ $f\ci t(v)>0$ for all\/ $v$ in $s^{-1}(0)\subseteq V$ then $\Pi_{T\cup U}^{T,f-}\ci i_*\bigl([V,n,s,t]\bigr)=[V,n,s,t]$.

Similarly, if\/ $[V,n,s,t]\in MC_k(U;R)$ with\/ $f\ci t(v)<0$ for all\/ $v$ in $s^{-1}(0)\subseteq V$ then $\Pi_{T\cup U}^{U,f+}\ci j_*\bigl([V,n,s,t]\bigr)=[V,n,s,t]$.
\label{kh7prop1}
\end{prop}

\begin{proof} For (i), observe that the conditions $\bigl\{y\in T\cup U:f(y)\ge 0\bigr\}\subseteq T$ and $\bigl\{y\in T\cup U:f(y)\le 0\bigr\}\subseteq U$ imply that we may define generators
\begin{equation*}
\bigl[T\t(-\iy,0],1,f\ci\pi_T+\pi_{(-\iy,0]},\pi_T\bigr],
\bigl[U\t[0,\iy),1,f\ci\pi_U+\pi_{[0,\iy)},\pi_U\bigr]
\end{equation*}
in $\cP MC^0(T\cup U;R)$, and then we can compare $\Pi_{T\cup U}^{T,f-},\Pi_{T\cup U}^{U,f+}$ in \eq{kh7eq10} with
\begin{align*}
-\cap\bigl[T\t(-\iy,0],1,f\!\ci\!\pi_T\!+\!\pi_{(-\iy,0]},\pi_T\bigr]&:MC_k(T\!\cup\! U;R)\!\ra\! MC_k(T\!\cup\! U;R),\\
-\cap\bigl[U\t[0,\iy),1,f\!\ci\!\pi_U\!+\!\pi_{[0,\iy)},\pi_U\bigr]&:MC_k(T\!\cup\! U;R)\!\ra\! MC_k(T\!\cup\! U;R),
\end{align*}
where the cap product
\begin{equation*}
\cap:MC_k(T\cup U;R)\t\cP MC^0(T\cup U;R)\longra MC_k(T\cup U;R)
\end{equation*}
is defined as in \S\ref{kh46}, except that we have reversed the order of $MC_k(T\cup U;R)$, $\cP MC^0(T\cup U;R)$. Then $\Pi_{T\cup U}^{T,f-},\Pi_{T\cup U}^{U,f+}$ are identical to these cap products, except that the targets are $MC_k(T;R),MC_k(U;R)$ rather than $MC_k(T\cup U;R)$. Hence the proof in \S\ref{kh46} and \S\ref{kh79} that $\cap$ is well defined also shows that $\Pi_{T\cup U}^{T,f-},\Pi_{T\cup U}^{U,f+}$ are well defined.

For (ii), if $\bigl[V,n,(s_1,\ldots,s_n),t\bigr]$ is a generator in $MC_k(T\cup U;R)$ then
\ea
\bigl(i_*&\ci\Pi_{T\cup U}^{T,f-}+j_*\ci\Pi_{T\cup U}^{U,f+}\bigr)\bigl(\bigl[V,n,(s_1,\ldots,s_n),t\bigr]\bigr)
\nonumber\\
&=\bigl[t^{-1}(T)\t(-\iy,0],n+1,(s_1\ci\pi_{t^{-1}(T)},\ldots,s_n\ci\pi_{t^{-1}(T)},
\nonumber\\
&\qquad\qquad\qquad f\ci\pi_{t^{-1}(T)}+\pi_{(-\iy,0]}),t\ci\pi_{t^{-1}(T)}\bigr]
\nonumber\\
&+\bigl[t^{-1}(U)\t[0,\iy),n+1,(s_1\ci\pi_{t^{-1}(U)},\ldots,s_n\ci\pi_{t^{-1}(U)},\nonumber\\
&\qquad\qquad\qquad f\ci\pi_{t^{-1}(U)}+\pi_{[0,\iy)}),t\ci\pi_{t^{-1}(U)}\bigr]
\nonumber\\
&=\bigl[V\t(-\iy,0],n+1,(s_1\ci\pi_V,\ldots,s_n\ci\pi_V,f\ci\pi_V+\pi_{(-\iy,0]}),t\ci\pi_V\bigr]
\nonumber\\
&+\bigl[V\t[0,\iy),n+1,(s_1\ci\pi_V,\ldots,s_n\ci\pi_V,f\ci\pi_V+\pi_{[0,\iy)}),t\ci\pi_V\bigr]
\nonumber\\
&=\bigl[V\t\R,n+1,(s_1\ci\pi_V,\ldots,s_n\ci\pi_V,f\ci\pi_V+\pi_\R),t\ci\pi_V\bigr]
\nonumber\\
&=\bigl[V\t\R,n+1,(s_1\ci\pi_V,\ldots,s_n\ci\pi_V,\pi_\R),t\ci\pi_V\bigr]
\nonumber\\&=\bigl[V,n,(s_1,\ldots,s_n),t\bigr].
\label{kh7eq11}
\ea
Here in the first step we use equations \eq{kh4eq7} and \eq{kh7eq10}. In the second we may enlarge $t^{-1}(T)\t(-\iy,0]$ to $V\t(-\iy,0]$ because the zeroes of $(s_1\ci\pi_V,\ldots,s_n\ci\pi_V,f\ci\pi_V+\pi_{[0,\iy)})$ are contained in $t^{-1}(T)\t(-\iy,0]\subseteq V\t(-\iy,0]$ as $\bigl\{y\in T\cup U:f(y)\ge 0\bigr\}\subseteq T$, and similarly for $t^{-1}(U)\t[0,\iy)$ and $V\t[0,\iy)$. In the third we combine $V\t(-\iy,0]$ and $V\t[0,\iy)$ into $V\t\R$ using relation Definition \ref{kh4def1}(ii). In the fourth we conjugate by the diffeomorphism $V\t\R\ra V\t\R$ mapping $(v,x)\mapsto (v,x-f(v))$, and in the fifth we use relation Definition \ref{kh4def1}(i). As \eq{kh7eq11} holds for all generators of $MC_k(T\cup U;R)$, part (ii) follows.

For (iii), suppose $[V,n,s,t]\in MC_k(T;R)$ with $f\ci t(v)>0$ for all $v$ in $s^{-1}(0)\subseteq V$, and write $s=(s_1,\ldots,s_n)$. Then as for \eq{kh7eq11} we have
\ea
\Pi_{T\cup U}^{T,f-}&\ci i_*\bigl([V,n,s,t]\bigr)\nonumber\\
&=\bigl[V\t(-\iy,0],n+1,(s_1\ci\pi_V,\ldots,s_n\ci\pi_V,f\ci\pi_V+\pi_{(-\iy,0]}),t\ci\pi_V\bigr]
\nonumber\\
&=\bigl[V\t\R,n+1,(s_1\ci\pi_V,\ldots,s_n\ci\pi_V,f\ci\pi_V+\pi_\R),t\ci\pi_V\bigr]
\nonumber\\
&=\bigl[V\t\R,n+1,(s_1\ci\pi_V,\ldots,s_n\ci\pi_V,\pi_\R),t\ci\pi_V\bigr]
\nonumber\\&=\bigl[V,n,(s_1,\ldots,s_n),t\bigr].
\label{kh7eq12}
\ea
Here in the first step we use \eq{kh4eq7} and \eq{kh7eq10}, noting that $t^{-1}(T)=V$ as $t$ maps $V\ra T\subseteq T\cup U$. In the second we use relation Definition \ref{kh4def1}(ii) in $MC_k(T;R)$, as the second and third lines of \eq{kh7eq12} coincide in a neighbourhood of $s^{-1}(0)$ in $V\t(-\iy,0]$ or $V\t\R$, since $f\ci t(v)>0$ for all $v$ in $s^{-1}(0)\subseteq V$. In the third we conjugate by the diffeomorphism $V\t\R\ra V\t\R$ mapping $(v,x)\mapsto (v,x-f(v))$, and in the fourth we use relation Definition \ref{kh4def1}(i). The second part of (iii) is proved in the same way. This completes the proposition.
\end{proof}

To prove Theorem \ref{kh4thm1}(a), suppose $T\subseteq U\subseteq Y$ are open, and write $i:T\hookra U$ for the inclusion. Suppose $\al\in MC_k(T;R)$ with $i_*(\al)=0$ in $MC_k(U;R)$. Write
\e
\al=\ts\sum_{a\in A}\al_a\,[V_a,n_a,s_a,t_a],
\label{kh7eq13}
\e
for $A$ a finite indexing set, $\al_a\in R$ and $[V_a,n_a,s_a,t_a]$ a generator of $MC_k(T;R)$. Then $\bigcup_{a\in A}t_a(s_a^{-1}(0))$ is a compact subset of $T$, so it is closed in $U$, and hence $\bigl(T, U\sm\bigcup_{i\in I}t_a(s_a^{-1}(0))\bigr)$ is an open cover of $U$. Therefore as in Definition \ref{kh7def1} we may choose smooth $f:U\ra\R$ with $\bigl\{y\in U:f(y)\ge 0\bigr\}\subseteq T$ and
\e
\bigl\{y\in U:f(y)\le 0\bigr\}\subseteq U\sm\ts\bigcup_{a\in A}t_a(s_a^{-1}(0)).
\label{kh7eq14}
\e
So Definition \ref{kh7def1} gives $\Pi_U^{T,f-}:MC_k(U;R)\ra MC_k(T;R)$, and Proposition \ref{kh7prop1}(iii) implies that for each $a\in A$ we have
\e
\Pi_U^{T,f-}\ci i_*\bigl([V_a,n_a,s_a,t_a]\bigr)=[V_a,n_a,s_a,t_a],
\label{kh7eq15}
\e
since \eq{kh7eq14} implies that $f\ci t_a(v)<0$ for all $v$ in $s_a^{-1}(0)\subseteq V_a$. Hence
\begin{align*}
\al&=\ts\sum_{a\in A}\al_a\,[V_a,n_a,s_a,t_a]=\sum_{a\in A}\al_a\,\Pi_U^{T,f-}\ci i_*\bigl([V_a,n_a,s_a,t_a]\bigr)\\
&=\Pi_U^{T,f-}\ci i_*(\al)=0,
\end{align*}
using \eq{kh7eq13}, \eq{kh7eq15} and $i_*(\al)=0$. Thus $i_*:MC_k(T;R)\ra MC_k(U;R)$ is injective, proving Theorem~\ref{kh4thm1}(a).

For (b), suppose $T,U\subseteq Y$ are open, and write $i:T\cap U\hookra T,$ $i':T\cap U\hookra U,$ $j:T\hookra T\cup U,$ $j':U\hookra T\cup U$ for the inclusions. Choose $f$ as in Definition \ref{kh7def1}, so that we have operators $\Pi_{T\cup U}^{T,f-},\Pi_{T\cup U}^{U,f+}$. Applying Definition \ref{kh7def1} with $\ti T=T$, $\ti U=T\cap U$, $\ti f=f\vert_T$ in place of $T,U,f$ gives operators 
\begin{align*}
\Pi_T^{T,f\vert_T-}&:MC_k(T;R)\longra MC_k(T;R),\\
\Pi_T^{T\cap U,f\vert_T+}&:MC_k(T;R)\longra MC_k(T\cap U;R),
\end{align*}
and applying Proposition \ref{kh7prop1}(ii) for these gives
\e
\Pi_T^{T,f\vert_T-}+i_*\ci\Pi_T^{T\cap U,f\vert_T+}=\id:MC_k(T;R)\longra MC_k(T;R).
\label{kh7eq16}
\e
Similarly we have operators
\ea
\Pi_U^{T\cap U,f\vert_U-}&:MC_k(U;R)\longra MC_k(T\cap U;R),
\nonumber\\
\Pi_U^{U,f\vert_U+}&:MC_k(U;R)\longra MC_k(U;R),\quad\text{with}
\nonumber\\
i'_*\ci\Pi_U^{T\cap U,f\vert_U-}&+\Pi_U^{U,f\vert_U+}=\id:MC_k(U;R)\longra MC_k(U;R).
\label{kh7eq17}
\ea

Comparing the actions of the two sides of each equation on generators $[V,n,s,t]$ using \eq{kh4eq7} and \eq{kh7eq10}, we see that
\ea
\Pi_T^{T,f\vert_T-}&=\Pi_{T\cup U}^{T,f-}\ci j_*:MC_k(T;R)\longra MC_k(T;R),
\label{kh7eq18}\\
\Pi_U^{U,f\vert_T+}&=\Pi_{T\cup U}^{U,f+}\ci j'_*:MC_k(U;R)\longra MC_k(U;R),
\label{kh7eq19}\\
i_*\ci\Pi_U^{T\cap U,f\vert_U-}&=\Pi_{T\cup U}^{T,f-}\ci j'_*:MC_k(U;R)\longra MC_k(T;R),
\label{kh7eq20}\\
i'_*\ci\Pi_T^{T\cap U,f\vert_T+}&=\Pi_{T\cup U}^{U,f+}\ci j_*:MC_k(T;R)\longra MC_k(U;R).
\label{kh7eq21}
\ea

As in \eq{kh4eq8}, we have to prove the following sequence is exact:
\e
\xymatrix@C=10.1pt{ 0 \ar[r] & MC_k(T\!\cap\! U;R) \ar[rr]^{i_*\op -i'_*} && {\begin{subarray}{l} \ts \; MC_k(T;R)\\ \ts\op MC_k(U;R) \end{subarray}} \ar[rr]^(0.45){j_*\op j'_*} && MC_k(T\!\cup\! U;R) \ar[r] & 0. }\!\!\!
\label{kh7eq22}
\e
Theorem \ref{kh4thm1}(a) implies that \eq{kh7eq22} is exact at the second term, since $i_*,i'_*$ are injective. Proposition \ref{kh7prop1}(ii) implies that $j_*\bigl(\Pi_{T\cup U}^{T,f-}(\al)\bigr)+j'_*\bigl(\Pi_{T\cup U}^{U,f+}(\al)\bigr)=\al$ for any $\al\in MC_k(T\cup U;R)$, so \eq{kh7eq22} is exact at the fourth term. 

To show \eq{kh7eq22} is exact at the third term, let $\be\!\in\! MC_k(T;R)$, $\ga\!\in\! MC_k(U;R)$ with $j_*(\be)+j'_*(\ga)=0$. Set $\al=\Pi_T^{T\cap U,f\vert_T+}(\be)-\Pi_U^{T\cap U,f\vert_U-}(\ga)$. Then
\begin{align*}
i_*(\al)&=i_*\ci\Pi_T^{T\cap U,f\vert_T+}(\be)-i_*\ci\Pi_U^{T\cap U,f\vert_U-}(\ga)\\
&=i_*\ci\Pi_T^{T\cap U,f\vert_T+}(\be)-\Pi_{T\cup U}^{T,f-}\ci j'_*(\ga)\\
&=i_*\ci\Pi_T^{T\cap U,f\vert_T+}(\be)+\Pi_{T\cup U}^{T,f-}\ci j_*(\be)\\
&=i_*\ci\Pi_T^{T\cap U,f\vert_T+}(\be)+\Pi_T^{T,f\vert_T-}(\be)=\be,
\end{align*}
using the definition of $\al$ in the first step, \eq{kh7eq20} in the second, $j_*(\be)+j'_*(\ga)=0$ in the third, \eq{kh7eq18} in the fourth, and \eq{kh7eq16} in the fifth. A similar proof using \eq{kh7eq21}, \eq{kh7eq19} and \eq{kh7eq17} yields $i'_*(\al)=-\ga$. Hence \eq{kh7eq22} is exact, proving Theorem~\ref{kh4thm1}(b).

For part (c), suppose $U_1\subseteq U_2\subseteq \cdots\subseteq Y$ are open with $U=\bigcup_{a=1}^\iy U_a$. Write $i_a:U_a\hookra U$ and $i_{a,b}:U_a\hookra U_b$, $a\le b$ for the inclusions. Then $i_b\ci i_{a,b}=i_a$, so
\begin{equation*}
(i_b)_*\ci(i_{a,b})_*=(i_a)_*:MC_k(U_a;R)\longra MC_k(U;R).
\end{equation*}
Thus the universal property of direct limits gives a unique morphism
\e
\pi:\underrightarrow{\lim}\,_{a=1}^\iy\, MC_k(U_a;R)\longra MC_k(U;R),
\label{kh7eq23}
\e
where the direct limit is over $(i_{a,b})_*:MC_k(U_a;R)\ra MC_k(U_b;R)$, such that
\e
(i_a)_*=\pi\ci\Pi_a:MC_k(U_a;R)\longra MC_k(U;R)
\label{kh7eq24}
\e
for all $a=1,2\ldots,$ where $\Pi_a:MC_k(U_a;R)\ra \underrightarrow{\lim}\,_{a=1}^\iy\, MC_k(U_a;R)$ is the natural morphism. We must show that $\pi$ in \eq{kh7eq23} is an isomorphism.

Suppose $\al\in \underrightarrow{\lim}\,_{a=1}^\iy\, MC_k(U_a;R)$ with $\pi(\al)=0$. Then $\al=\Pi_a(\al_a)$ for some $a\ge 1$ and $\al_a\in MC_k(U_a;R)$. But then $(i_a)_*(\al_a)=\pi\ci\Pi_a(\al_a)=\pi(\al)=0$ by \eq{kh7eq24}, so $\al_a=0$ as $(i_a)_*:MC_k(U_a;R)\ra MC_k(U;R)$ is injective by Theorem \ref{kh4thm1}(a), and thus $\al=\Pi_a(\al_a)=0$. Hence $\pi$ in \eq{kh7eq23} is injective.

Suppose $\al\in MC_k(U;R)$, and write 
\e
\al=\ts\sum_{c\in C}\al_c\,[V_c,n_c,s_c,t_c]
\label{kh7eq25}
\e
for $C$ a finite indexing set, $\al_c\in R$ and $[V_c,n_c,s_c,t_c]$ a generator of $MC_k(U;R)$. Then $\bigcup_{c\in C}t_c(s_c^{-1}(0))$ is a compact subset of $U$, and $\{U_1,U_2,\ldots\}$ is an open cover of $U$, so by taking a finite subcover $U_{a_1},\ldots,U_{a_n}$ for $\bigcup_{c\in C}t_c(s_c^{-1}(0))$ and setting $a=\max(a_1,\ldots,a_N)$ we see that there exists $a=1,2,\ldots$ such that $t_c(s_c^{-1}(0))\subseteq U_a\subseteq U$ for all $c\in C$. Then we have
\e
\begin{split}z
[V_c,n_c,s_c,t_c]&=\bigl[t_c^{-1}(U_a),n_c,s_c\vert_{t_c^{-1}(U_a)},t_c\vert_{t_c^{-1}(U_a)}\bigr]\\
&=(i_a)_*\bigl(\bigl[t_c^{-1}(U_a),n_c,s_c\vert_{t_c^{-1}(U_a)},t_c\vert_{t_c^{-1}(U_a)}\bigr]\bigr),
\end{split}
\label{kh7eq26}
\e
since $s_c^{-1}(0)\subseteq t_c^{-1}(U_a)\subseteq V_c$. So from \eq{kh7eq24}--\eq{kh7eq26} we see that
\begin{equation*}
\al=\pi\Bigl[\ts\sum_{c\in C}\al_c\,\Pi_a\bigl(\bigl[t_c^{-1}(U_a),n_c,s_c\vert_{t_c^{-1}(U_a)},t_c\vert_{t_c^{-1}(U_a)}\bigr]\bigr)\Bigr].
\end{equation*}
Thus $\pi$ in \eq{kh7eq23} is surjective, and so is an isomorphism. This proves Theorem \ref{kh4thm1}(c). The last two paragraphs of Theorem \ref{kh4thm1} are immediate.

\subsection{Proof of Proposition \ref{kh4prop2}}
\label{kh73}

Suppose $Y_1,Y_2$ are manifolds, $Z_1\subseteq Y_1,$ $Z_2\subseteq Y_2$ are open, and $g:Y_1\t[0,1]\ra Y_2$ is a smooth map of manifolds with $g(Z_1\t[0,1])\subseteq Z_2$. Define $f,f':Y_1\ra Y_2$ by $f(y)=g(y,0)$ and $f'(y)=g(y,1)$ for $y\in Y_1$. For all $k\in\Z$, define $G:MC_k(Y_1,Z_1;R)\ra MC_{k+1}(Y_2,Z_2;R)$ to be the unique $R$-linear map acting on generators by
\e
G:[V,n,s,t]\longmapsto (-1)^{\dim V}\bigl[V\t[0,1],n,s\ci\pi_V,g\ci(t\t\id_{[0,1]})\bigr],
\label{kh7eq27}
\e
where $V\t[0,1]$ has the product orientation of the given orientation on $V$ and the standard orientation on $[0,1]$. To show that $G$ is well-defined, extend $g$ to a smooth map $\ti g:\ti Y_1\ra Y_2$ for $\ti Y_1$ an open neighbourhood of $Y_1\t[0,1]$ in $Y_1\t\R$, and let $\ti Z_1$ be an open neighbourhood of $Z_1\t[0,1]$ in $Y_1\t\R$ with $\ti g(\ti Z_1)\subseteq Z_2$. Then $G$ is the composition of a morphism $\t[0,1]:MC_k(Y_1,Z_1;R)\ra MC_{k+1}(\ti Y_1,\ti Z_1;R)$ mapping $[V,n,s,t]\mapsto (-1)^{\dim V}\bigl[V\t[0,1],n+1,s\ci\pi_V,t\t\id_{[0,1]}\bigr]$, which is easy to see is well-defined, with the pushforward $\ti g_*:MC_{k+1}(\ti Y_1,\ti Z_1;R)\ra MC_{k+1}(Y_2,Z_2;R)$, which is well-defined as in Definitions \ref{kh4def2} and \ref{kh4def3}. We have
\begin{align*}
\pd&\ci G[V,n,s,t]\\
&=(-1)^{\dim V}\bigl[\pd\bigl(V\t[0,1]\bigr),n,s\ci\pi_V\ci i_{V\t[0,1]},g\ci (t\t\id_{[0,1]})\ci i_{V\t[0,1]}\bigr]\\
&=(-1)^{\dim V}\bigl[\pd V\t[0,1],n,s\ci i_V\ci\pi_{\pd V},g\ci((t\ci i_V)\t\id_{[0,1]})\bigr]\\
&\quad +\bigl[V\t(\pd[0,1]),n,s\ci\pi_V,g\ci(t\t i_{[0,1]})\bigr]\\
&=-G[\pd V,n,s\ci i_V,t\ci i_V]-[V,n,s,f\ci t]+[V,n,s,f'\ci t]\\
&=(-G\ci\pd-f_*+f'_*)[V,n,s,t],
\end{align*}
using \eq{kh4eq3} and \eq{kh7eq27} in the first step, Assumptions \ref{kh3ass3}(d) and \ref{kh3ass6}(h) and \eq{kh4eq5}--\eq{kh4eq6} in the second, $\pd[0,1]=-\{0\}\amalg\{1\}$, \eq{kh4eq5}--\eq{kh4eq6} and the definitions of $f,f',G$ in the third, and \eq{kh4eq3} and \eq{kh4eq7} in the fourth. As this holds for all generators $[V,n,s,t]$, we have
\begin{equation*}
\pd\ci G+G\ci\pd=f'_*-f_*:MC_k(Y_1,Z_1;R)\longra MC_k(Y_2,Z_2;R).
\end{equation*}
So $G$ is a chain homotopy from $f_*$ to $f_*'$ on M-chains, and thus $f_*=f_*'$ on M-homology. This proves Proposition~\ref{kh4prop2}.

\subsection{Proof of Theorem \ref{kh4thm2}}
\label{kh74}

Theorem \ref{kh4thm2} says that $MH_k(*;R)=0$ for $k\ne 0$ and $MH_0(*;R)\cong R$. Lemma \ref{kh4lem1} implies that $MH_k(*;R)=0$ for $k>0$. After introducing some notation and proving some auxiliary results in \S\ref{kh741}, we will show that $MH_k(*;R)=0$ for $k<0$ in \S\ref{kh742}, and that $MH_0(*;R)\cong R$ in~\S\ref{kh743}.

\subsubsection{Some auxiliary results}
\label{kh741}

We begin by introducing some notation. Definition \ref{kh7def2} and Propositions \ref{kh7prop2} and \ref{kh7prop3} actually work for $\widetilde{MC}_k(Y;R)$ for any manifold $Y$, assuming $k\le 0$ in Proposition \ref{kh7prop3}, but for simplicity we give them only for $Y=*$, the point, as this is all we need.

\begin{dfn} Let $*$ be the point, and $R$ a commutative ring. As in \S\ref{kh41}, $MC_k(*;R)$ is the $R$-module spanned by generators $[V,n,s,t]$ with $\dim V=n+k$ subject to relations Definition \ref{kh4def1}(i),(ii).  Any morphism $t:V\ra *$ is the projection $\pi:V\ra *$, so the $t$ in $[V,n,s,t]=[V,n,s,\pi]$ can basically be ignored. 

Write $\widetilde{MC}_k(*;R)$ for the $R$-module spanned by generators $[V,n,s,\pi]$ with $\dim V=n+k$ subject to relation Definition \ref{kh4def1}(ii) only. Then there is a surjective $R$-linear map
\begin{equation*}
\Pi:\widetilde{MC}_k(*;R)\longra MC_k(*;R),\qquad \Pi:[V,n,s,\pi]\longmapsto[V,n,s,\pi]
\end{equation*}
with kernel spanned by equation \eq{kh4eq1} from Definition~\ref{kh4def1}(i). 

Since the relation Definition \ref{kh4def1}(ii) in $\widetilde{MC}_k(*;R)$ involves only $[V_i,n,s_i,\pi]$ with $n$ fixed, we may write 
\begin{equation*}
\widetilde{MC}_k(*;R)=\ts\bigop_{n=0}^\iy \widetilde{MC}_k(*;R)^n
\end{equation*}
where $\widetilde{MC}_k(*;R)^n$ is spanned by generators $[V,n,s,\pi]$ with $n$ fixed, modulo relation Definition \ref{kh4def1}(ii) with $n$ fixed. 

We may define $\pd:\widetilde{MC}_k(*;R)\ra\widetilde{MC}_{k-1}(*;R)$ as in \S\ref{kh41}, and then $\Pi\ci\pd=\pd\ci\Pi$, and $\pd\ci\pd=0$, as the proof of this in Definition \ref{kh4def1} used only relation (ii). Also $\pd$ maps~$\widetilde{MC}_k(*;R)^n\ra\widetilde{MC}_{k-1}(*;R)^n$.
\label{kh7def2}
\end{dfn}

For the next two propositions we will consider the following situation. Let $R$ be a commutative ring, $k\in\Z$, and $\al\in MC_k(*;R)$ with $\pd\al=0$ in $MC_{k-1}(*;R)$. Choose a lift $\ti\al$ of $\al$ to $\widetilde{MC}_k(*;R)$, so that $\Pi(\ti\al)=\al$, and write
\e
\ti\al=\sum_{n=0}^N\sum_{i\in I^n}a_i^n\,\bigl[V_i^n,n,s_i^n,\pi\bigr]
\label{kh7eq28}
\e
where $N\in\N$, $I^0,I^1,\ldots,I^N$ are finite indexing sets, $a_i^n\in R$ and $[V_i^n,n,s_i^n,\pi]$ is a generator of $MC_k(*;R)$ for all $i,n$, as in Definition \ref{kh4def1}. Then $\Pi\ci\pd\ti\al=\pd\ci\Pi(\ti\al)=\pd\al=0$, so $\pd\ti\al$ lies in the subspace of $\widetilde{MC}_{k-1}(*;R)$ spanned by relations \eq{kh4eq1} from Definition \ref{kh4def1}(i). Hence by \eq{kh4eq3}, we may write
\ea
&\sum_{n=0}^N\sum_{i\in I^n}a_i^n\,\bigl[\pd V_i^n,n,s_i^n\ci i_{V_i^n},\pi\bigr]
\label{kh7eq29}\\
&=\sum_{n=0}^{N-1}\sum_{l=0}^n\sum_{j\in J^{n,l}}b_j^{n,l}\Bigl(\bigl[\dot V_j^{n,l},n,(\dot s_{j,1}^{n,l},\ldots,\dot s_{j,n}^{n,l}),\pi\bigr]-(-1)^{n-l}\bigl[\dot V_j^{n,l}\t\R,
\nonumber\\
&\qquad\qquad\qquad n+1,(\dot s_{j,1}^{n,l},\ldots,\dot s_{j,l}^{n,l},\pi_\R,\dot s_{j,l+1}^{n,l},\ldots,\dot s_{j,n}^{n,l}),\pi\bigr]\Bigr)
\nonumber
\ea
in $\widetilde{MC}_{k-1}(*;R)$. Here we increase $N\ge 0$ in \eq{kh7eq28} if necessary, setting $I^n=\es$ for any additional $n$, so that we can take the largest $n$ on the r.h.s.\ of \eq{kh7eq29} to be $N-1$. Also $J^{n,l}$ are finite indexing sets for $0\le l\le n<N$, and $b_j^{n,l}\in R$ and $[\dot V_j^{n,l},n,(\dot s_{j,1}^{n,l},\ldots,\dot s_{j,n}^{n,l}),\pi]$ are generators of $MC_{k-1}(*;R)$ for all $n,l,j$. The r.h.s.\ of \eq{kh7eq29} is applications of relation Definition \ref{kh4def1}(i), with $l$ in place of $i$. As a shorthand we write
\e
\dot s_j^{n,l}=(\dot s_{j,1}^{n,l},\ldots,\dot s_{j,n}^{n,l})\quad\text{and}\quad \acute s_j^{n,l}=(\dot s_{j,1}^{n,l},\ldots,\dot s_{j,l}^{n,l},\pi_\R,\dot s_{j,l+1}^{n,l},\ldots,\dot s_{j,n}^{n,l}).
\label{kh7eq30}
\e

Taking components of \eq{kh7eq29} in $\widetilde{MC}_{k-1}(*;R)^n$ for $n=0,\ldots,N$ shows that
\ea
\begin{split}
&\sum_{i\in I^n}a_i^n\,\bigl[\pd V_i^n,n,s_i^n\ci i_{V_i^n},\pi\bigr]=\sum_{l=0}^n\sum_{j\in J^{n,l}}b_j^{n,l}\,\bigl[\dot V_j^{n,l},n,\dot s_j^{n,l},\pi\bigr]\\
&-\sum_{l=0}^{n-1}\sum_{j\in J^{n-1,l}}(-1)^{n-1-l}b_j^{n-1,l}\,\bigl[\dot V_j^{n-1,l}\t\R,n,\acute s_j^{n-1,l},\pi\bigr],
\end{split}
\label{kh7eq31}
\ea
setting $J^{N,l}=\es$ for~$l=0,\ldots,N$.

We first show we can make the boundaries of the $\bigl[\dot V_j^{n,l},n,\dot s_j^{n,l},\pi\bigr]$ cancel.

\begin{prop} In the situation above, for all\/ $k\in\Z,$ taking $\al,\ti\al,N$ and the representation \eq{kh7eq28} for $\ti\al$ to be fixed, we may make alternative choices for the data $J^{n,l},b_j^{n,l},[\dot V_j^{n,l},n,\dot s_j^{n,l},\pi]$ on the r.h.s.\ of\/ \eq{kh7eq29}--\eq{kh7eq31} to ensure that for all\/ $n=0,\ldots,N-1$ and\/ $l=0,\ldots,n$ we have 
\e
\pd\raisebox{-4pt}{$\biggl[$}\sum_{j\in J^{n,l}}b_j^n\,\bigl[\dot V_j^{n,l},n,\dot s_j^{n,l},\pi\bigr]\raisebox{-4pt}{$\biggr]$}=0\qquad\text{in $\widetilde{MC}_{k-2}(*;R)^n$.}
\label{kh7eq32}
\e

\label{kh7prop2}
\end{prop}

\begin{proof} The proof is by a double induction, the outer induction being on decreasing $n'=N-1,N-2,\ldots,0$, and the inner induction being on $l'=0,\ldots,n'+1$. The inductive hypothesis is:
\begin{itemize}
\setlength{\itemsep}{0pt}
\setlength{\parsep}{0pt}
\item[$(\dag)_{n',l'}$] We can make alternative choices for the $J^{n,l},b_j^{n,l},[\dot V_j^{n,l},n,\dot s_j^{n,l},\pi]$ such that \eq{kh7eq32} holds for all $n=n'+1,\ldots,N-1$ and all $l=0,\ldots,n$, and also \eq{kh7eq32} holds for $n=n'$ and all~$l=0,\ldots,l'-1$.
\end{itemize}
The first case $(\dag)_{N-1,0}$ is vacuous, so the first step of the induction is trivial. Note that $(\dag)_{n',n'+1}$ and $(\dag)_{n'-1,0}$ are equivalent, as both say that \eq{kh7eq32} holds for all $n\ge n'$ and all $l=0,\ldots,n$. The final case $(\dag)_{0,1}$ says that \eq{kh7eq32} holds for all $n,l$, and so will prove the proposition.

In the inductive step, we suppose that $(\dag)_{n',l'}$ holds for some $n=0,\ldots,N-1$ and $l'=0,\ldots,n'$. We will show we can modify the $J^{n,l},b_j^{n,l},[\dot V_j^{n,l},n,\dot s_j^{n,l},\pi]$ to make \eq{kh7eq32} hold for $n=n'$, $l=l'$, whilst preserving the cases of \eq{kh7eq32} in $(\dag)_{n',l'}$. This proves $(\dag)_{n',l'+1}$. Hence by induction $(\dag)_{n',l'+1},(\dag)_{n',l'+2},(\dag)_{n',n'+1}$ hold, completing the inner induction. But $(\dag)_{n',n'+1}=(\dag)_{n'-1,0}$ if $n'>0$, so $(\dag)_{n'-1,0}$ holds, the inductive step in the outer induction. Continuing the induction we see that $(\dag)_{n'-1,1},(\dag)_{n'-1,2},\ldots,(\dag)_{n'-1,n'}$ hold, and so $(\dag)_{n'-2,0}$ holds if $n'>1$. Eventually we show $(\dag)_{0,1}$ holds, and finish.

So suppose $(\dag)_{n',l'}$ holds for some $n=0,\ldots,N-1$ and $l'=0,\ldots,n'$. Observe that \eq{kh7eq31} and \eq{kh7eq32} are equations in $\widetilde{MC}_*(*;R)^n$ which hold by an application of relation Definition \ref{kh4def1}(ii). Hence for each of these equations there exists an open neighbourhood $X$ of 0 in $\R^n$ such that for every generator $[V,n,s,\pi]$ occurring in the equation, $s:V\ra\R^n$ is proper over $X$, and condition Definition \ref{kh4def1}(ii)$(*)$ holds with this~$X$. 

Choose small $\ep>0$ such that $(-\ep,\ep)^{n'+1}\subseteq X\subseteq\R^{n'+1}$ for an allowed choice of open neighbourhood $X$ of 0 in $\R^{n'+1}$ for equation \eq{kh7eq31} with $n=n'+1$, and $(-\ep,\ep)^{n'+1}\subseteq X\subseteq\R^{n'+1}$ if $n'<N-1$ for an allowed choice of $0\in X\subseteq\R^{n'+1}$ for equation \eq{kh7eq32} with $n=n'+1$ and all $l=0,\ldots,n'+1$, and $(-\ep,\ep)^{n'}\subseteq X\subseteq\R^{n'}$ for an allowed choice of $0\in X\subseteq\R^{n'}$ for equation \eq{kh7eq32} with $n=n'$ and all $l=0,\ldots,l'-1$, where $(\dag)_{n',l'}$ guarantees \eq{kh7eq32} holds in these cases.

Suppose $[V,n'+1,(s_1,\ldots,s_{n'+1}),\pi]$ is one of the generators occurring in \eq{kh7eq31} for $n=n'+1$, so that $[V,n'+1,(s_1,\ldots,s_{n'+1}),\pi]$ is $\bigl[\pd V_i^{n'+1},n'+1,s_i^{n'+1}\ci i_{V_i^{n'+1}},\pi\bigr]$ or $\bigl[\dot V_j^{n'+1,l},n,\dot s_j^{n'+1,l},\pi\bigr]$ or $\bigl[\dot V_j^{n',l}\t\R,n'+1,\acute s_j^{n',l},\pi\bigr]$. Then $s_{l'+1}:V\ra\R$ is a morphism in $\tManc$. If $\dim V=k+n'+1\ge 1$ then Assumption \ref{kh3ass7}(b) implies that there exists $S\subset\R$ with $\cH^1(S)=0$ such that if $u\in\R\sm S$ then $s_{l'+1}:V\ra\R$ is a submersion in an open neighbourhood of $s_{l'+1}^{-1}(u)$. If $\dim V<1$ then $V$ is a 0-manifold or $V=\es$ by Assumption \ref{kh3ass4}(e), so $S=s_{l'+1}(V)$ has $\cH^1(S)=0$, and if $u\in\R\sm S$ then $s_{l'+1}^{-1}(u)=\es$, so trivially $s_{l'+1}:V\ra\R$ is a submersion in an open neighbourhood of~$s_{l'+1}^{-1}(u)$.

The complement of finitely many subsets $S\subset\R$ with $\cH^1(S)=0$ is dense in $\R$. Thus we may choose $u\in(-\ep,\ep)$ such that $s_{l'+1}:V\ra\R$ is a submersion in an open neighbourhood of $s_{l'+1}^{-1}(u)$ for all generators $[V,n'+1,(s_1,\ldots,s_{n'+1}),\pi]$ occurring in \eq{kh7eq31} for $n=n'+1$ (this implies the same thing holds for \eq{kh7eq32} for $n=n'+1$ and all $l=0,\ldots,n'+1$ with $s_{l'+1}:V\ra\R$, and for \eq{kh7eq32} for $n=n'$ and all $l=0,\ldots,l'-1$ for $s_{l'}:V\ra\R$), and also that $s^{n'+1}_{i,l'+1}:V_i^n\ra\R$ is a submersion in an open neighbourhood of $(s^{n'+1}_{i,l'+1})^{-1}(u)$ for all~$i\in I^{n'+1}$.

For each generator $[V,n'+1,(s_1,\ldots,s_{n'+1}),\pi]$ occurring in \eq{kh7eq31} for $n=n'+1$, and also for $[V,n'+1,(s_1,\ldots,s_{n'+1}),\pi]=[V_i^{n'+1},n'+1,s_i^{n'+1},\pi]$ for $i\in I^{n'+1}$, consider the morphisms $s_{l'+1}:V\ra\R$ and $u:*\ra\R$, $u:*\mapsto u$. Since $s_{l'+1}$ is a submersion near $s_{l'+1}^{-1}(u)$ in $V$, Assumption \ref{kh3ass5}(e) implies that $s_{l'+1},u$ are transverse, so a fibre product $V\t_{s_{l'+1},\R,u}*$ exists in $\tManc$ by Assumption \ref{kh3ass5}(c), and combining the given orientation on $V$ with the standard orientations on $\R,*$, Assumption \ref{kh3ass6}(l) gives an orientation on~$V\t_{s_{l'+1},\R,u}*$.

Thus as oriented objects of $\tManc$, we may define
\e
\begin{aligned}
W^{n'+1}_i&=V^{n'+1}_i\t_{s^{n'+1}_{i,l'+1},\R,u}*, && i\in I^{n'+1},\\
\dot W^{n'+1,l}_j&=\dot V^{n'+1,l}_j\t_{\dot s^{n'+1,l}_{j,l'+1},\R,u}*, && l=0,\ldots,n'+1,\;\> j\in J^{n'+1,l},\\
\dot W^{n',l}_j&=\dot V^{n',l}_j\t_{\dot s^{n',l}_{j,l'},\R,u}*, && l=0,\ldots,l'-1,\;\> j\in J^{n',l},\\
\dot W^{n',l}_j&=\dot V^{n',l}_j\t_{\dot s^{n',l}_{j,l'+1},\R,u}*, && l=l'+1,\ldots,n',\;\> j\in J^{n',l}.
\end{aligned}
\label{kh7eq33}
\e
Here the point of the final two lines of \eq{kh7eq33} is that starting with the generator $\bigl[\dot V_j^{n',l}\t\R,n'+1,\acute s_j^{n',l},\pi\bigr]$ in \eq{kh7eq31} for $n=n'+1$, we have
\begin{align*}
&(\dot V_j^{n',l}\t\R)\t_{\acute s_{j,l'+1}^{n',l},\R,u}*\cong \\
&\qquad\begin{cases}
(\dot V^{n',l}_j\t_{\dot s^{n',l}_{j,l'},\R,u}*)\t\R=\dot W^{n',l}_j\t\R, & l=0,\ldots,l'-1, \\
\dot V^{n',l'}_j, & l=l', \\
(\dot V^{n',l}_j\t_{\dot s^{n',l}_{j,l'+1},\R,u}*)\t\R=\dot W^{n',l}_j\t\R, & l=l'+1,\ldots,n',
\end{cases}
\end{align*}
using \eq{kh7eq30} and properties of (oriented) fibre products, including Assumption \ref{kh3ass6}(m) and \eq{kh3eq3}--\eq{kh3eq5}, which imply that $(\dot V_j^{n',l}\t\R)\t_\R *\cong (\dot V^{n',l}_j\t_\R *)\t\R$ holds in oriented objects in $\tManc$ without additional signs.

Now define morphisms in $\tManc$:
\begin{gather*}
\begin{aligned}
\hat s_i^{n'+1}&:W^{n'+1}_i\longra\R^{n'}, &  \check s_j^{n'+1,l}&:\dot W^{n'+1,l}_j\longra\R^{n'}, \\
\check s_j^{n',l}&:\dot W^{n',l}_j\longra\R^{n'-1},\;\> l\ne l', & \end{aligned}
\\
\begin{split}
\text{by}\quad \hat s_i^{n'+1}&:=(s^{n'+1}_{i,1},\ldots,s^{n'+1}_{i,l'},s^{n'+1}_{i,l'+2},\ldots, s^{n'+1}_{i,n'+1})\ci\pi_{V^{n'+1}_i},\\
\check s_j^{n'+1,l}&:=(\dot s^{n'+1,l}_{j,1},\ldots,\dot s^{n'+1,l}_{j,l'},\dot s^{n'+1,l}_{j,l'+2},\ldots,\dot s^{n'+1,l}_{j,n'+1})\ci\pi_{\dot V^{n'+1,l}_j},\\
\check s_j^{n',l}&:=(\dot s^{n',l}_{j,1},\ldots,\dot s^{n',l}_{j,l'-1},\dot s^{n',l}_{j,l'+1},\ldots,\dot s^{n',l}_{j,n'})\ci\pi_{\dot V^{n',l}_j},\quad l<l',\\
\check s_j^{n',l}&:=(\dot s^{n',l}_{j,1},\ldots,\dot s^{n',l}_{j,l'},\dot s^{n',l}_{j,l'+2},\ldots,\dot s^{n',l}_{j,n'})\ci\pi_{\dot V^{n',l}_j},\quad l>l'.\\
\end{split}
\end{gather*}
Here the $W^{n'+1}_i,\dot W^{n'+1,l}_j,\dot W^{n',l}_j$ are made by a fibre product over $\R$ of one of the coordinates in $\R^{n'+1},\R^{n'}$ from the morphisms $s^{n'+1}_i:V^{n'+1}_i\ra\R^{n'+1},\ldots,\dot s^{n',l}_j:\dot V^{n',l}_j\ra\R^{n'}$, and the morphisms $\hat s_i^{n'+1},\check s_j^{n'+1,l},\check s_j^{n',l}$ use the remaining coordinates in $\R^{n'+1},\R^{n'}$. We may now form generators
\begin{align*}
\bigl[W^{n'+1}_i,n',\hat s_i^{n'+1},\pi\bigr]&\!\in\! \widetilde{MC}_{k}(*;R)^{n'}, && i\in I^{n'+1},\\
\bigl[\dot W^{n'+1,l}_j,n',\check s_j^{n'+1,l},\pi\bigr]&\!\in\! \widetilde{MC}_{k-1}(*;R)^{n'}, && l\!=\!0,\ldots,n'\!+\!1,\; j\!\in\! J^{n'+1,l},\\
\bigl[\dot W^{n',l}_j,n'-1,\check s_j^{n',l},\pi\bigr]&\!\in\! \widetilde{MC}_{k-1}(*;R)^{n'-1}, && l\!=\!0,\ldots,n',\; l\!\ne\! l',\; j\!\in\! J^{n',l}.
\end{align*}
To see that these satisfy the properness conditions near 0 in $\R^{n'},\R^{n'-1}$ in Definition \ref{kh4def1}, note that by choice of $\ep$, the $V^{n'+1}_i,\dot V^{n'+1,l}_j,\dot V^{n',l}_j$ in the fibre products \eq{kh7eq33} are proper over $(-\ep,\ep)^{n'+1}\subset\R^{n'+1}$ or $(-\ep,\ep)^{n'}\subset\R^{n'}$, so as $u\in(-\ep,\ep)$, the corresponding $W^{n'+1}_i,\dot W^{n'+1,l}_j,\dot W^{n',l}_j$ are proper over $(-\ep,\ep)^{n'}\subset\R^{n'}$ or~$(-\ep,\ep)^{n'-1}\subset\R^{n'-1}$.

We now claim that \eq{kh7eq31} for $n=n'+1$ implies that
\ea
\begin{split}
&\sum_{i\in I^{n'+1}}a_i^{n'+1}\,\bigl[\pd W_i^{n'+1},n',\hat s_i^{n'+1}\ci i_{W_i^{n'+1}},\pi\bigr]\\
&=\sum_{l=0}^{n'+1}\sum_{j\in J^{n'+1,l}}b_j^{n'+1,l}\,\bigl[\dot W_j^{n'+1,l},n',\check s_j^{n'+1,l},\pi\bigr]\\
&\quad-\sum_{\begin{subarray}{l} l=0,\ldots,n',\\ l\ne l'\end{subarray}}
\sum_{j\in J^{n',l}}\begin{aligned}[t] &(-1)^{n'-l}b_j^{n',l}\,\bigl[\dot W_j^{n',l}\t\R,n',(\check s_{j,1}^{n',l},\ldots,\check s_{j,l}^{n',l},\\
&\quad\qquad \pi_\R,\check s_{j,l+1}^{n',l},\ldots,\check s_{j,n'}^{n',l})\ci\pi_{\dot W_j^{n',l}},\pi\bigr] 
\end{aligned}\\
&\quad-\sum_{j\in J^{n',l'}}(-1)^{n'-l'}b_j^{n',l'}\,\bigl[\dot V_j^{n',l'},n',\dot s_j^{n',l'},\pi\bigr]
\end{split}
\label{kh7eq34}
\ea
in $\widetilde{MC}_{k-1}(*;R)^{n'}$. To see this, note that \eq{kh7eq34} is essentially the result of restricting \eq{kh7eq31} for $n=n'+1$, thought of as living over $\R^{n'+1}$,  to the hyperplane $x_{l'+1}=u$ in $\R^{n'+1}$, writing $(x_1,\ldots,x_{n'+1})$ for the coordinates in $\R^{n'+1}$. Equation \eq{kh7eq31} is an application of relation Definition \ref{kh4def1}(ii) in $\widetilde{MC}_{k-1}(*;R)^{n'+1}$, where the condition Definition \ref{kh4def1}(ii)$(*)$ holds over $X\t *$ with $X=(-\ep,\ep)^{n'+1}\subset\R^{n'+1}$ by choice of $\ep$. From this we can deduce using Assumptions \ref{kh3ass5}, \ref{kh3ass6} that condition Definition \ref{kh4def1}(ii)$(*)$ holds for \eq{kh7eq34} over $X\t *$ with $X=(-\ep,\ep)^{n'}\subset\R^{n'}$, so \eq{kh7eq34} holds. 

Similarly, from \eq{kh7eq32} for $n=n'+1$ and $l=0,\ldots,n'+1$, and \eq{kh7eq32} for $n=n'$ and $l=0,\ldots,l'-1$, which hold by $(\dag)_{n',l'}$, we may deduce that
\ea
\pd\raisebox{-4pt}{$\biggl[$}\sum_{j\in J^{n'+1,l}}b_j^{n'+1}\,\bigl[\dot W^{n'+1,l}_j,n',\check s_j^{n'+1,l},\pi\bigr]
\raisebox{-4pt}{$\biggr]$}&=0,\quad l=0,\ldots,n'+1,
\label{kh7eq35}\\
\pd\raisebox{-4pt}{$\biggl[$}\sum_{j\in J^{n',l}}b_j^{n'}\,\bigl[\dot W^{n',l}_j,n'-1,\check s_j^{n',l},\pi\bigr]
\raisebox{-4pt}{$\biggr]$}&=0,\quad l=0,\ldots,l'-1.
\label{kh7eq36}
\ea

We will now define alternative choices $\ddot J^{n,l},\ddot b_j^{n,l},\bigl[\ddot V_j^{n,l},n,\ddot s_j^{n,l},\pi\bigr]$ for the data $J^{n,l},b_j^{n,l},\bigl[\dot V_j^{n,l},n,\dot s_j^{n,l},\pi\bigr]$ in \eq{kh7eq29}--\eq{kh7eq31}. For $n'\ne n$ and all $l,j$ we take
\e
\ddot J^{n,l}=J^{n,l},\quad \ddot b_j^{n,l}=b_j^{n,l},\quad \bigl[\ddot V_j^{n,l},n,\ddot s_j^{n,l},\pi\bigr]=\bigl[\dot V_j^{n,l},n,\dot s_j^{n,l},\pi\bigr],
\label{kh7eq37}
\e
that is, we make no change. When $n=n'$ and $l\ne l'$ we set
\ea
\ddot J^{n',l}&=J^{n',l}\amalg J^{n',l}=J^{n',l}_1\amalg J^{n',l}_2, \;\> l=0,\ldots,n',\;\> l\ne l', 
\label{kh7eq38}\\
\ddot b_j^{n',l}&=b_j^{n',l},\;\> \bigl[\ddot V_j^{n',l},n',\ddot s_j^{n',l},\pi\bigr]=\bigl[\dot V_j^{n',l},n',\dot s_j^{n',l},\pi\bigr],\;\> \text{all $l$, $j\in J^{n',l}_1$,}
\nonumber\\
\ddot b_j^{n',l}&=(-1)^{l'-l+1}b_j^{n',l},\qquad \bigl[\ddot V_j^{n',l},n',\ddot s_j^{n',l},\pi\bigr]=\bigl[\dot W_j^{n',l}\t\R,n',
\nonumber\\
(\check s_{j,1}^{n',l}&,\ldots,\check s_{j,l}^{n',l},\pi_\R,\check s_{j,l+1}^{n',l},\ldots,\check s_{j,n'}^{n',l})\!\ci\!\pi_{\dot W_j^{n',l}},\pi\bigr],\;\> \text{all $l$, $j\!\in\! J^{n',l}_2$,}
\nonumber
\ea
where we distinguish the two copies of $J^{n',l}$ in $\ddot J^{n',l}$ by writing them $J^{n',l}_1,J^{n',l}_2$, and when $n'=n$ and $l-l'$ we set
\e
\begin{split}
\ddot J^{n',l'}&=I^{n'+1}\amalg \ts\coprod_{l=0}^{n'+1}J^{n'+1,l}, \\
\ddot b_i^{n',l'}&=(-1)^{n'-l'+1}a_i^{n'+1},\qquad \bigl[\ddot V_i^{n',l'},n',\ddot s_i^{n',l'},\pi\bigr]=\\
&\quad\bigl[\pd W_i^{n'+1},n',\hat s_i^{n'+1}\ci i_{W_i^{n'+1}},\pi\bigr],\quad i\in I^{n'+1},\\
\ddot b_j^{n',l'}&=(-1)^{n'-l'}b_j^{n'+1,l},\qquad \bigl[\ddot V_j^{n',l'},n',\ddot s_j^{n',l'},\pi\bigr]=\\
&\quad\bigl[\dot W_j^{n'+1,l},n',\check s_j^{n'+1,l},\pi\bigr], \quad
j\in J^{n'+1,l},\; l=0,\ldots,n'+1.
\end{split}
\label{kh7eq39}
\e

We now claim that replacing $J^{n,l},b_j^{n,l},\bigl[\dot V_j^{n,l},n,\dot s_j^{n,l},\pi\bigr]$ by the alternative choices $\ddot J^{n,l},\ddot b_j^{n,l},\bigl[\ddot V_j^{n,l},n,\ddot s_j^{n,l},\pi\bigr]$, the analogue of \eq{kh7eq29} holds in $\widetilde{MC}_{k-1}(*;R)$. To see this, multiply \eq{kh7eq34} by $(-1)^{n'-l'}$, and use it to substitute for the first term $\sum_{j\in J^{n',l'}}b_j^{n',l'}\,\bigl[\dot V_j^{n',l'},n',\dot s_j^{n',l'},\pi\bigr]$ in the case $n=n'$, $l=l'$ in the sum on the r.h.s.\ of \eq{kh7eq29}. Also, modify \eq{kh7eq34} by replacing each generator $[V,n',(s_1,\ldots,s_{n'}),\pi]$ in \eq{kh7eq34} by $\bigl[V\t\R,n'+1,(s_1\ci\pi_V,\ldots,s_{l'}\ci\pi_V,\pi_\R,s_{l'+1}\ci\pi_V,\ldots,s_{n'}\ci\pi_V),\pi\bigr]$, and use this new equation to substitute for the second term $\sum_j(-1)^{n'-l'}b_j^{n',l'}\bigl[\dot V_j^{n',l'}\t\R,n'+1,(\dot s_{j,1}^{n',l'},\ldots,\dot s_{j',l'}^{n',l'},\pi_\R,\dot s_{j,l'+1}^{n',l'},\ldots,\dot s_{j,n'}^{n',l'}),\pi\bigr]$ in the case $n=n'$, $l=l'$ in the sum on the r.h.s.\ of~\eq{kh7eq29}.

In this way we get a modified version of \eq{kh7eq29}, in which the terms involving $\bigl[\dot V_j^{n,l},n,\dot s_j^{n,l},\pi\bigr]$ for $(n,l)\ne(n',l')$ remain unchanged, but the terms involving $\bigl[\dot V_j^{n',l'},n',\dot s_j^{n',l'},\pi\bigr]$ have been deleted, and replaced by the other terms in \eq{kh7eq34}. By \eq{kh7eq37}--\eq{kh7eq39}, this modified version is exactly the analogue of \eq{kh7eq29} with the $J^{n,l},b_j^{n,l},\bigl[\dot V_j^{n,l},n,\dot s_j^{n,l},\pi\bigr]$ replaced by the~$\ddot J^{n,l},\ddot b_j^{n,l},\bigl[\ddot V_j^{n,l},n,\ddot s_j^{n,l},\pi\bigr]$.

To see this, note that the terms $a_i^{n'+1}\,\bigl[\pd W_i^{n'+1},n',\hat s_i^{n'+1}\ci i_{W_i^{n'+1}},\pi\bigr]$ and $b_j^{n'+1,l}\,\bigl[\dot W_j^{n'+1,l},n',\check s_j^{n'+1,l},\pi\bigr]$ in \eq{kh7eq34} are transferred into terms of the form $\ddot b_j^{n',l'}\bigl[\ddot V_j^{n',l'},n',\ddot s_j^{n',l'},\pi\bigr]$ by \eq{kh7eq39}, with signs to compensate for the $-(-1)^{n'-l'}$ in the last line of \eq{kh7eq34}. The terms coming from $(-1)^{n'-l}b_j^{n',l}\,\bigl[\dot W_j^{n',l}\t\R,n',\ldots,\pi\bigr]$ in \eq{kh7eq34} for $l\ne l'$ are a little more subtle. When we replace $[V,n',s,\pi]$ in \eq{kh7eq34} by $\bigl[V\t\R,n'+1,\ldots,\pi\bigr]$ these terms yield 
$(-1)^{n'-l}b_j^{n',l}\,\bigl[\dot W_j^{n',l}\t\R\t\R,n'+1,\ldots,\pi\bigr]$, where the two $\R$ factors in $\dot W_j^{n',l}\t\R\t\R$ are inserted in positions $l$ and $l'$ in $\R^{n'+1}$. So we could regard these terms as contributing to either $\ddot J^{n',l},\ddot b_j^{n',l},\ldots$ or to $\ddot J_j^{n',l'},\ddot b_j^{n',l'},\ldots,$ and as in \eq{kh7eq38} we include them in~$\ddot J^{n',l},\ddot b_j^{n',l},\ldots.$

We also claim that our alternative choices $\ddot J^{n,l},\ddot b_j^{n,l},\bigl[\ddot V_j^{n,l},n,\ddot s_j^{n,l},\pi\bigr]$ satisfy $(\dag)_{n',l'+1}$. To see this,  observe that \eq{kh7eq32} for $\ddot J^{n,l},\ldots$ for $n=n'+1,\ldots,N-1$ and  $l=0,\ldots,n$ follows from \eq{kh7eq32} for $J^{n,l},\ldots$ for $n=n'+1,\ldots,N-1$ and  $l=0,\ldots,n$ by \eq{kh7eq37} (which holds by $(\dag)_{n',l'}$ for $J^{n,l},\ldots$). Also \eq{kh7eq32} for $\ddot J^{n,l},\ldots$ for $n=n'$ and $l=0,1,\ldots,l'-1$ follows from \eq{kh7eq38}, and from \eq{kh7eq32} for $J^{n,l},\ldots$ for $n=n'$ and $l=0,1,\ldots,l'-1$ (which holds by $(\dag)_{n',l'}$ for $J^{n,l},\ldots$), and from \eq{kh7eq36} for $l=0,1,\ldots,l'-1$. And \eq{kh7eq32} for $\ddot J^{n,l},\ldots$ for $n=n'$ and $l=l'$ follows from \eq{kh7eq39}, the fact that $\pd\bigl[\pd W_i^{n'+1},n',\hat s_i^{n'+1}\ci i_{W_i^{n'+1}},\pi\bigr]=0$ so that the terms from $I^{n'+1}\subseteq\ddot J^{n',l'}$ in \eq{kh7eq39} contribute zero to \eq{kh7eq32}, and equation \eq{kh7eq35} for $l=0,\ldots,n'+1$, which implies that the terms from $\coprod_{l=0}^{n'+1}J^{n'+1,l}\subseteq\ddot J^{n',l'}$ in \eq{kh7eq39} contribute zero to~\eq{kh7eq32}.

This proves $(\dag)_{n',l'+1}$, and completes the inductive step. So by induction, $(\dag)_{n',l'}$ holds for all $n'=N-1,N-2,\ldots,0$ and $l'=0,\ldots,n'+1$, so in particular $(\dag)_{0,1}$ holds. This completes the proof of Proposition~\ref{kh7prop2}.
\end{proof}

\begin{prop} In the situation above, by keeping $\al\in MC_k(*;R)$ with\/ $\pd\al=0$ fixed but changing the lift\/ $\ti\al$ of\/ $\al$ to $\widetilde{MC}_k(*;R)$ with\/ $\Pi(\ti\al)=\al,$ we may suppose that\/ $\pd\ti\al=0$ in $\widetilde{MC}_{k-1}(*;R)$.
\label{kh7prop3}
\end{prop}

\begin{proof} If $k>0$ the proposition is trivial since $MC_k(*;R)=\widetilde{MC}_k(*;R)=0$ as in Lemma \ref{kh4lem1}, so suppose $k\le 0$. We first choose an arbitrary lift $\ti\al$ of $\al$ to $\widetilde{MC}_k(*;R)$ of the form \eq{kh7eq28}. Then Proposition \ref{kh7prop2} shows that we may write $\pd\ti\al$ as in \eq{kh7eq29} in terms of data $J^{n,l},b_j^{n,l},[\dot V_j^{n,l},n,\dot s_j^{n,l},\pi]$ satisfying \eq{kh7eq32} for all $n=0,\ldots,N-1$ and $l=0,\ldots,n$. Since each such $\bigl[\dot V_j^{n,l},n,\dot s_j^{n,l},\pi\bigr]$ is a generator of $MC_{k-1}(*;R)$ we have $\dim\dot V_j^{n,l}=n+k-1<n$, as $k\le 0$. Consider the subset $\dot s_j^{n,l}(\dot V_j^{n,l})$ in $\R^n$. By differential geometry $\dot s_j^{n,l}[(\dot V_j^{n,l})^\ci]$ has Hausdorff dimension $\le n+k-1<n$, and $\cH^n[\dot s_j^{n,l}(\dot V_j^{n,l}\sm(\dot V_j^{n,l})^\ci)]=0$ by Assumption \ref{kh3ass7}(a) as $\dim\dot V_j^{n,l}\le n$. Hence~$\cH^n[\dot s_j^{n,l}(\dot V_j^{n,l})]=0$.

For each $n=0,1,\ldots,N-1$, choose an open neighbourhood $X^n$ of 0 in $\R^n$ such that $\dot s_j^{n,l}:\dot V_j^{n,l}\ra\R^n$ is proper over $X^n$ for all $l=0,\ldots,n$ and $j\in J^{n,l}$. This is possible by definition of generators in Definition \ref{kh4def1}. Also, equation \eq{kh7eq32} holds as an application of relation Definition \ref{kh4def1}(ii) in $\widetilde{MC}_{k-2}(*;R)^n$, which involves a condition $(*)$ upon points $(x,y)\in X\t *$ for $X$ an open neighbourhood of 0 in $\R^n$. Making $X^n$ smaller if necessary, we suppose  that Definition \ref{kh4def1}(ii)$(*)$ holds for \eq{kh7eq32} with~$X=X^n$.

Next, choose $x^n\in X^n$ such that 
\begin{itemize}
\setlength{\itemsep}{0pt}
\setlength{\parsep}{0pt}
\item[(a)] $\bigl\{\la x^n:\la\in[0,1]\bigr\}\subseteq X^n$, and
\item[(b)] $x^n\notin\dot s_j^{n,l}(\dot V_j^{n,l})$ for all $l=0,\ldots,n$ and $j\in J^{n,l}$.
\end{itemize}
Here (a) holds provided $x^n$ is small enough in $\R^n$, and (b) holds for generic $x^n$ as $\cH^n[\dot s_j^{n,l}(\dot V_j^{n,l})]=0$, so both are possible. Choose small $\de^n>0$ such that 
\begin{itemize}
\setlength{\itemsep}{0pt}
\setlength{\parsep}{0pt}
\item[(c)] $B_{\de^n}(\la x^n)\subseteq X^n$ for all $\la\in[0,1]$, with $B_{\de^n}(x^n)$ the open ball of radius $\de^n$ about $x$ in $\R^n$; and
\item[(d)] $B_{\de^n}(x^n)\cap \dot s_j^{n,l}(\dot V_j^{n,l})=\es$ for all $l=0,\ldots,n$ and $j\in J^{n,l}$.
\end{itemize}
Here (c),(d) are both possible if $\de^n$ is small enough, by (a),(b), since $X^n$ is open in $\R^n$, and $\dot s_j^{n,l}(\dot V_j^{n,l})\cap X^n$ is closed in $X^n$ as $\dot s_j^{n,l}$ is proper over~$X^n$.

For each $l=0,\ldots,n$ and $j\in J^{n,l}$, define $\check s_j^{n,l}:[0,1)\t\dot V_j^{n,l}\ra\R^n$ by
\begin{equation*}
\check s_j^{n,l}:(\la,v)\longmapsto \dot s_j^{n,l}(v)-\la x^n.
\end{equation*}
We claim that $\bigl[[0,1)\t\dot V_j^{n,l},n,\check s_j^{n,l},\pi\bigr]$ is a generator of $MC_k(*;R)$, where we give $[0,1)\t\dot V_j^{n,l}$ the product orientation of the standard orientation on $[0,1)$ and the given orientation on $\dot V_j^{n,l}$, as in Assumption~\ref{kh3ass6}(f),(k). 

The issue is to show that $\check s_j^{n,l}:[0,1)\t\dot V_j^{n,l}\ra\R^n$ is proper over an open neighbourhood of 0 in $\R^n$, noting that $[0,1)$ is noncompact. In fact (c),(d) above imply that $\check s_j^{n,l}$ is proper over $B_{\de^n}(0)$ in $\R^n$, since part (c), $\check s_j^{n,l}$ proper over $X^n$, and $[0,1]$ compact, imply that the obvious extension of $\check s_j^{n,l}$ to a map $[0,1]\t \dot V_j^{n,l}\ra\R^n$ is proper over $B_{\de^n}(0)$, and then part (d) implies that restricting to $\check s_j^{n,l}:[0,1)\t \dot V_j^{n,l}\ra\R^n$ is still proper over $B_{\de^n}(0)$, since $\{1\}\t\dot V_j^{n,l}$ does not map to $B_{\de^n}(0)$. Hence~$\bigl[[0,1)\t\dot V_j^{n,l},n,\check s_j^{n,l},\pi\bigr]\in MC_k(*;R)$.

Using equation \eq{kh4eq3}, Assumptions \ref{kh3ass3}(d) and \ref{kh3ass6}(h), and $\pd[0,1)=-\{0\}$ in oriented manifolds, we see that in $MC_{k-1}(*;R)$ or $\widetilde{MC}_{k-1}(*;R)$ we have
\e
\pd\bigl[[0,1)\!\t\!\dot V_j^{n,l},n,\check s_j^{n,l},\pi\bigr]\!=\!-\bigl[\dot V_j^{n,l},n,\dot s_j^{n,l},\pi\bigr]\!-\!\bigl[[0,1)\!\t\!\pd\dot V_j^{n,l},n,\hat s_j^{n,l},\pi\bigr],
\label{kh7eq40}
\e
where $\hat s_j^{n,l}:[0,1)\t\pd\dot V_j^{n,l}\ra\R^n$ maps~$\hat s_j^{n,l}:(\la,v')\mapsto \dot s_j^{n,l}\ci i_{\dot V_j^{n,l}}(v')-\la x^n$.

We may also define a generator $\bigl[[0,1)\t\dot V_j^{n,l}\t\R,n+1,\grave s_j^{n,l},\pi\bigr]$, where
\begin{align*}
\grave s_j^{n,l}=(&\hat s_{j,1}^{n,l}\ci\pi_{[0,1)\t\dot V_j^{n,l}},\ldots,\hat s_{j,l}^{n,l}\ci\pi_{[0,1)\t\dot V_j^{n,l}},\pi_\R,\\
&\hat s_{j,l+1}^{n,l}\ci\pi_{[0,1)\t\dot V_j^{n,l}},\ldots,\hat s_{j,n}^{n,l}\ci\pi_{[0,1)\t\dot V_j^{n,l}}),
\end{align*}
and then as for \eq{kh7eq40} we see that
\ea
&\pd\bigl[[0,1)\t\dot V_j^{n,l}\t\R,n+1,\grave s_j^{n,l},\pi\bigr]=-\bigl[\dot V_j^{n,l}\t\R,n+1,\acute s_j^{n,l},\pi\bigr]
\nonumber\\
&\qquad -\bigl[[0,1)\t\pd\dot V_j^{n,l}\t\R,n+1,\grave s_j^{n,l}\ci i_{[0,1)\t\dot V_j^{n,l}\t\R},\pi\bigr].
\label{kh7eq41}
\ea

Define a modification $\breve\al$ of $\ti\al$ in $\widetilde{MC}_k(*;R)$ by
\ea
\breve\al&=\sum_{n=0}^N\sum_{i\in I^n}a_i^n\,\bigl[V_i^n,n,s_i^n,\pi\bigr]
+\sum_{n=0}^{N-1}\sum_{l=0}^n\sum_{j\in J^{n,l}}b_j^{n,l}\,\bigl[[0,1)\t\dot V_j^{n,l},n,\check s_j^{n,l},\pi\bigr]
\nonumber\\
&-\sum_{n=0}^{N-1}\sum_{l=0}^n\sum_{j\in J^{n,l}}(-1)^{n-l}b_j^{n,l}\,\bigl[[0,1)\t\dot V_j^{n,l}\t\R,n+1,\grave s_j^{n,l},\pi\bigr].
\label{kh7eq42}
\ea
Then $\breve\al$ differs from $\ti\al$ by finitely many applications of relation Definition \ref{kh4def1}(i) in $MC_k(*;R)$, so $\breve\al$ is an alternative lift of $\al$ to $\widetilde{MC}_k(*;R)$, with $\Pi(\breve\al)=\al$. Applying $\pd$ to \eq{kh7eq42} and using equations \eq{kh7eq29}, \eq{kh7eq40} and \eq{kh7eq41} yields
\ea
&\pd\breve\al\!=\!-\!\sum_{n=0}^{N-1}\sum_{l=0}^n\sum_{j\in J^{n,l}}\!b_j^{n,l}
\bigl[[0,1)\!\t\!\pd\dot V_j^{n,l},n,\hat s_j^{n,l},\pi\bigr]\!+\!\sum_{n=0}^{N-1}\sum_{l=0}^n(-1)^{n-l}\cdot{}
\label{kh7eq43}\\
&\qquad\sum_{j\in J^{n,l}}b_j^{n,l}
\bigl[[0,1)\t\pd\dot V_j^{n,l}\!\t\!\R,n\!+\!1,\grave s_j^{n,l}\!\ci\! i_{[0,1)\t\dot V_j^{n,l}\t\R},\pi\bigr]
\nonumber
\ea
in~$\widetilde{MC}_{k-1}(*;R)$.

Now for all $n=0,\ldots,N-1$ and $l=0,\ldots,n$, equation \eq{kh7eq32} holds. For each fixed $n,l$, the two sums $\sum_{j\in J^{n,l}}b_j^{n,l}[\cdots]$ in \eq{kh7eq43} are each the result of applying an operation to the l.h.s.\ of \eq{kh7eq32}, where we replace each generator $[\pd V,n,s,\pi]$ in \eq{kh7eq32} by $\bigl[[0,1)\t\pd V,n,\hat s,\pi\bigr]$ or by $\bigl[[0,1)\t\pd V\t\R,n+1,\grave s,\pi\bigr]$. Hence \eq{kh7eq32} implies that these two sums $\sum_{j\in J^{n,l}}b_j^{n,l}[\cdots]$ in \eq{kh7eq43} are zero for all $n,l$, so $\pd\breve\al=0$ in $\widetilde{MC}_{k-1}(*;R)$. Replacing $\ti\al$ by $\breve\al$, Proposition \ref{kh7prop3} follows.
\end{proof}

\subsubsection{Proof that $MH_k(*;R)=0$ for $k<0$}
\label{kh742}

We can now prove the first part of Theorem \ref{kh4thm2}. We use the notation of \S\ref{kh741}. Let $R$ be a commutative ring, $k<0$, and $\al\in MC_k(*;R)$ with $\pd\al=0$ in $MC_{k-1}(*;R)$. Then Proposition \ref{kh7prop3} says that we may choose a lift $\ti\al$ of $\al$ to $\widetilde{MC}_k(*;R)$ with $\pd\ti\al=0$ in $\widetilde{MC}_{k-1}(*;R)$. Write $\ti\al$ as in \eq{kh7eq28}. Then for each $n=0,\ldots,N$ we have
\e
\sum_{i\in I^n}a_i^n\,\bigl[\pd V_i^n,n,s_i^n\ci i_{V_i^n},\pi\bigr]=0
\qquad\text{in $\widetilde{MC}_{k-1}(*;R)^n$.}
\label{kh7eq44}
\e

We follow the method of the proof of Proposition \ref{kh7prop3}, but applied to the $\bigl[V_i^n,n,s_i^n,\pi\bigr]$ for $i\in I^n$ rather than to the $\bigl[\dot V_j^{n,l},n,\dot s_j^{n,l},\pi\bigr]$ for $j\in J^{n,l}$. For each $n=0,1,\ldots,N$, we choose an open neighbourhood $X^n$ of 0 in $\R^n$ such that $s_i^n:V_i^n\ra\R^n$ is proper over $X^n$ for all $i\in I^n$. Equation \eq{kh7eq44} holds as an application of relation Definition \ref{kh4def1}(ii) in $\widetilde{MC}_{k-1}(*;R)^n$. Making $X^n$ smaller if necessary, we suppose that Definition \ref{kh4def1}(ii)$(*)$ holds for \eq{kh7eq44} with~$X=X^n$.

As in the proof of Proposition \ref{kh7prop3}, we choose $x^n\in X^n$ and $\de^n>0$ satisfying
\begin{itemize}
\setlength{\itemsep}{0pt}
\setlength{\parsep}{0pt}
\item[(a)] $\bigl\{\la x^n:\la\in[0,1]\bigr\}\subseteq X^n$;
\item[(b)] $x^n\notin s_i^n(V_i^n)$ for all $i\in I^n$;
\item[(c)] $B_{\de^n}(\la x^n)\subseteq X^n$ for all $\la\in[0,1]$; and
\item[(d)] $B_{\de^n}(x^n)\cap s_i^n(V_i^n)=\es$ for all $i\in I^n$.
\end{itemize}
Here (b) is possible as $\dim V^i_n=n+k<n$ since $k<0$, so $\cH^n[s_i^n(V_i^n)]=0$. For each $i\in I^n$, define $\check s_i^n:[0,1)\t V_i^n\ra\R^n$ by $\check s_i^n:(\la,v)\mapsto s_i^n(v)-\la x^n$. Then $\bigl[[0,1)\t V_i^n,n,\check s_i^n,\pi\bigr]$ is a generator of $MC_{k+1}(*;R)$, and as in \eq{kh7eq40}
\e
\pd\bigl[[0,1)\t V_i^n,n,\check s_i^n,\pi\bigr]=-\bigl[V_i^n,n,s_i^n,\pi\bigr]-\bigl[[0,1)\t\pd V_i^n,n,\hat s_i^n,\pi\bigr],
\label{kh7eq45}
\e
where $\hat s_i^n:[0,1)\t\pd V_i^n\ra\R^n$ is given by~$\hat s_i^n:(\la,v')\mapsto s_i^n\ci i_{V_i^n}(v')-\la x^n$.

Now define $\ti\be\in\widetilde{MC}_{k+1}(*;R)$ by
\e
\ti\be=-\sum_{n=0}^N\sum_{i\in I^n}a_i^n\,\bigl[[0,1)\t V_i^n,n,\check s_i^n,\pi\bigr].
\label{kh7eq46}
\e
Then combining equations \eq{kh7eq28}, \eq{kh7eq45} and \eq{kh7eq46} gives
\e
\pd\ti\be=\ti\al+\sum_{n=0}^N\sum_{i\in I^n}a_i^n\,\bigl[[0,1)\t\pd V_i^n,n,\hat s_i^n,\pi\bigr].
\label{kh7eq47}
\e
As for equations \eq{kh7eq32} and \eq{kh7eq43} in the proof of Proposition \ref{kh7prop3}, from \eq{kh7eq44} we see that for each $n=0,\ldots,N$ the sum $\sum_{i\in I^n}a_i^n[\cdots]$ in \eq{kh7eq47} is zero, so $\pd\ti\be=\ti\al$ in $\widetilde{MC}_k(*;R)$. Setting $\be=\Pi(\ti\be)$ in $MC_{k+1}(*;R)$, we have $\pd\be=\al$. Thus, for all $\al\in MC_k(*;R)$ with $\pd\al=0$ for $k<0$, we can find $\be\in MC_{k+1}(*;R)$ with $\pd\be=\al$. Hence $MH_k(*;R)=0$ for $k<0$, as we have to prove.

\subsubsection{Proof that $MH_0(*;R)\cong R$}
\label{kh743}

Since $MC_1(*;R)=0$ by Lemma \ref{kh4lem1}, we have
\e
MH_0(*;R)=\Ker\bigl(\pd:MC_0(*;R)\longra MC_{-1}(*;R)\bigr).
\label{kh7eq48}
\e
As $*$ has the standard orientation, Definition \ref{kh4def1} defines the fundamental cycle $[*]=[*,0,0,\id_*]\in MC_0(*;R)$, with $\pd[*]=0$. Define an $R$-linear map
\e
\io:R\longra \Ker\bigl(\pd:MC_0(*;R)\longra MC_{-1}(*;R)\bigr),\qquad \io:a\longmapsto a\,[*].
\label{kh7eq49}
\e
We will show that $\io$ is both surjective and injective, so that $MH_0(*;R)\cong R$, by an isomorphism identifying $[[*]]\in MH_0(*;R)$ with $1\in R$, as in Theorem~\ref{kh4thm2}.

Suppose $\al\in MC_0(*;R)$ with $\pd\al=0$ in $MC_{-1}(*;R)$. Using the notation $\widetilde{MC}_k(*;R),\widetilde{MC}_k(*;R)^n$ of \S\ref{kh741}, Proposition \ref{kh7prop3} gives a lift $\ti\al$ of $\al$ to $\widetilde{MC}_0(*;R)$ with $\Pi(\ti\al)=\al$ and $\pd\ti\al=0$ in $\widetilde{MC}_{-1}(*;R)$. As in \eq{kh7eq28}, write
\begin{equation*}
\ti\al=\sum_{n=0}^N\sum_{i\in I^n}a_i^n\,\bigl[V_i^n,n,s_i^n,\pi\bigr].
\end{equation*}
Then taking the component of $\pd\ti\al=0$ in $\widetilde{MC}_{-1}(*;R)^n$ for $n=0,\ldots,N$ yields
\e
\sum_{i\in I^n}a_i^n\,\bigl[\pd V_i^n,n,s_i^n\ci i_{V_i^n},\pi\bigr]=0.
\label{kh7eq50}
\e

By definition of generators $\bigl[V_i^n,n,s_i^n,\pi\bigr]\in MC_0(*;R)$ in \S\ref{kh41}, for each $n=0,\ldots,N$ we may choose an open neighbourhood $X^n$ of 0 in $\R^n$ such that $s_i^n:V_i^n\ra\R^n$ is proper over $X^n$ for all $i\in I^n$. Equation \eq{kh7eq50} holds as an application of relation Definition \ref{kh4def1}(ii) in $\widetilde{MC}_{-1}(*;R)$. Thus, making $X^n$ smaller we can suppose that Definition \ref{kh4def1}(ii)$(*)$ holds for \eq{kh7eq50} with $X=X^n$. Making $X^n$ still smaller, we can suppose it is a convex subset of~$\R^n$.

For each $n=0,\ldots,N$, define a subset $S^n\subseteq X^n\subseteq\R^n$ by
\e
\begin{split}
S^n=X^n\cap\Bigl[\ts\bigcup_{i\in I^n}s_i^n\bigl(\bigl\{v\in V_i^n:
\text{either $v\notin(V_i^n)^\ci$, or $v\in (V_i^n)^\ci$ and}\\
\text{$T_vs_i^n:T_v(V_i^n)^\ci\ra T_{s_i^n(v)}\R^n$ is not an isomorphism}\bigr\}\bigr)\Bigr].
\end{split}
\label{kh7eq51}
\e
As the subsets of $V_i^n$ in \eq{kh7eq51} are closed, and $s_i^n$ is proper over $X^n$, we see that $S^n$ is closed in $X^n$. Using Sard's Theorem and Assumption \ref{kh3ass7}(a) we can show that $\cH^n(S^n)=0$, which implies that~$X^n\sm S^n\ne\es$.

Define a function $\Phi^n:X^n\sm S^n\ra R$ by
\ea
&\Phi^n(x)=
\label{kh7eq52}\\
&\sum_{i\in I^n}\sum_{\begin{subarray}{l} v\in (V^n_i)^\ci: \\ s_i^n(v)=x \end{subarray}}\!a_i^n\cdot
\begin{cases} 1, & \!\!\!\text{$T_vs_i^n:T_v(V_i^n)^\ci\!\ra\! T_x\R^n$ is orientation-preserving,} \\ -1, & \!\!\!\text{$T_vs_i^n:T_v(V_i^n)^\ci\!\ra\! T_x\R^n$ is orientation-reversing.}
\end{cases}
\nonumber
\ea
This is well-defined as $s_i^n$ is a diffeomorphism near each $v$ in \eq{kh7eq52}, so $(s_i^n)^{-1}(x)$ has the discrete topology, and also $(s_i^n)^{-1}(x)$ is compact as $s_i^n$ is proper over $X^n\ne x$, so $(s_i^n)^{-1}(x)$ is finite. 

The proof of the next lemma is related to that of Proposition \ref{kh4prop1} in~\S\ref{kh71}.

\begin{lem} The function $\Phi^n:X^n\sm S^n\ra R$ is constant for $n=0,\ldots,N$.
\label{kh7lem1}
\end{lem}

\begin{proof} Since $s_i^n\vert_{(V_i^n)^\ci}:(V_i^n)^\ci\ra\R^n$ is proper and \'etale over $X^n\sm S^n$, we see that $\Phi^n$ is locally constant on $X^n\sm S^n$, and thus constant on each connected component of $X^n\sm S^n$. Let $x_0\ne x_1\in X^n\sm S^n$. We will show that $\Phi^n(x_0)=\Phi^n(x_1)$. Choose a linear isomorphism $\La:\R^n\ra\R^{n-1}\t\R$ with $\La(x_0)=(y,z_0)$ and $\La(x_1)=(y,z_1)$ for $y\in\R^{n-1}$ and $z_0\ne z_1\in\R$. Then $\pi_{\R^{n-1}}\ci\La\ci s_i^n:V_i^n\ra\R^{n-1}$ is a morphism in $\tManc$ for $i\in I^n$, where~$\dim V_i^n=n$. 

Assumption \ref{kh3ass7}(b) says that there is a subset $\ti S\subseteq\R^{n-1}$ with $\cH^{n-1}(\ti S)=0$ such that $\pi_{\R^{n-1}}\ci\La\ci s_i^n:V_i^n$ is a submersion near $(\pi_{\R^{n-1}}\ci\La\ci s_i^n)^{-1}(y')$ in $V_i^n$ for all $y'\in\R^{n-1}\sm \ti S$ and $i\in I^n$. As $\R^{n-1}\sm \ti S$ is dense in $\R^{n-1}$, this holds for $y'$ arbitrarily close to $y$. Thus we can find $x_0',x_1'$ close to $x_0,x_1$ in $X^n\sm S^n$ and in the same connected components of $X^n\sm S^n$ as $x_0,x_1$, so that $\Phi^n(x_0')=\Phi^n(x_0)$, $\Phi^n(x_1')=\Phi^n(x_1)$ as $\Phi^n$ is locally constant, with $\La(x_0')=(y',z_0')$ and $\La(x_1')=(y',z_1')$ for some $y'\in\R^{n-1}\sm \ti S$ close to $y$ and~$z_0'\ne z_1'\in\R$.

Now define $h:[0,1]\ra\R^n$ by $h(\la)=(1-\la)x_0'+\la x_1'$, so that $h(0)=x_0'$ and $h(1)=x_1'$. Then $h$ maps $[0,1]\ra X^n$ as $X^n$ is convex in $\R^n$. We claim that $s_i^n:V_i^n\ra\R^n$ and $h:[0,1]\ra\R^n$ are transverse morphisms in $\tManc$. To prove this, suppose $v\in V_i^n$ and $\la\in[0,1]$ with $s_i^n(v)=h(\la)$. If $\la=0$ or 1 then $h(\la)=x_0'$ or $x_1'$, so $h(\la)\in X^n\sm S^n$, and $s_i^n:V_i^n\ra\R^n$ is a submersion near $v$ in $V_i^n$ by \eq{kh7eq51}, which implies that $s_i^n,h$ are transverse on open neighbourhoods of $v\in V_i^n$ and $\la\in[0,1]$ by Assumption~\ref{kh3ass5}(c). 

If $\la\in(0,1)$ we identify $\R^n\cong\R^{n-1}\t\R$ using $\La$, and write $s_i^n:V_i^n\ra\R^{n-1}\t\R$ and $h:[0,1]\ra\R^{n-1}\t\R$ as direct products $(s_{i,1}^n,s_{i,2}^n)$ and $(h_1,h_2)$, where $s_{i,1}^n=\pi_{\R^{n-1}}\ci\La\ci s_i^n:V_i^n\ra\R^{n-1}$ is a submersion near $v\in (\pi_{\R^{n-1}}\ci\La\ci s_i^n)^{-1}(y')$ as $y'\in\R^{n-1}\sm \ti S$, and $h_2:[0,1]\ra\R$, $h_2(\la)=(1-\la)z_0'+\la z_1'$ is a submersion near $\la\in(0,1)$ as $z_0'=z_1'$. Hence $s_i^n,h$ are transverse on open neighbourhoods of $v\in V_i^n$ and $\la\in[0,1]$ by Assumption \ref{kh3ass5}(f). This holds for all $v\in V_i^n$ and $\la\in[0,1]$ with $s_i^n(v)=h(\la)$, so $s_i^n,h$ are transverse by Assumption~\ref{kh3ass5}(e).

Therefore by Assumption \ref{kh3ass5}(c), the transverse fibre product
\e
T_i^n:=V_i^n\t_{s_i^n,\R^n,h}[0,1]
\label{kh7eq53}
\e
exists in $\tManc$, with $\dim T_i^n=1$. Consider the diagram of topological spaces
\begin{equation*}
\xymatrix@C=50pt@R=15pt{
*+[r]{T_i^n} \ar[r]_(0.25){\pi_{V_i^n}\t\pi_{[0,1]}} & \bigl\{(v,\la)\in V_i^n\t[0,1]:s_i^n(v)=h(\la)\bigr\} \ar[r]_(0.7){\pi_{[0,1]}} \ar[d]^{\pi_{V_i^n}} & *+[l]{[0,1]} \ar[d]_h \\
& {\!\!\!\!\!\!\!\!\!\!\!\! V_i^n\subseteq (s_i^n)^{-1}(X^n)} \ar[r]^{s_i^n\vert_{\cdots}} & *+[l]{X^n.\!} }
\end{equation*}
The top left map is a homeomorphism by \eq{kh3eq8} in Assumption \ref{kh3ass5}(c). The bottom map is proper as $s_i^n$ is proper over $X^n$. Thus the top right morphism is proper, as the square is Cartesian. Hence as $[0,1]$ is compact, $T_i^n$ is compact.
 
Combining the given orientation on $V_i^n$ with the standard orientations on $\R^n,[0,1]$ from Assumption \ref{kh3ass6}(k), by Assumption \ref{kh3ass6}(l) we have an orientation on $T_i^n$ in \eq{kh7eq53}. By Assumptions \ref{kh3ass6}(c), \ref{kh3ass6}(l),(m) and equations \eq{kh3eq2} and \eq{kh3eq9}, since $\pd[0,1]=-\{0\}\amalg\{1\}$ in oriented objects in $\tManc$ and $h(0)=x_0'$, $h(1)=x_1'$, in oriented 0-manifolds we have
\e
\pd T_i^n\cong \bigl(\pd V_i^n\t_{s_i^n\ci i_{V_i^n},\R^n,h}[0,1]\bigr)\amalg -\bigl(V_i^n\t_{s_i^n,\R^n,x_0'}*\bigr)\amalg \bigl(V_i^n\t_{s_i^n,\R^n,x_1'}*\bigr).
\label{kh7eq54}
\e

Now Assumption \ref{kh3ass6}(n) says that the number of points in \eq{kh7eq54}, counted with signs, is zero. Hence we have
\ea
0&=\!\!\!\!\!\sum_{\begin{subarray}{l} \la\in(0,1),\\ v'\in (\pd V_i^n)^\ci:\\
s_i^n\ci i_{V_i^n}(v')=h(\la) \end{subarray}\!\!\!\!\!\!\!\!\!\!\!\!\!\!\!\!} \!\!{\begin{cases} 1, & \!\!\!\text{
$T_{v'}(s_i^n\!\ci\! i_{V_i^n})[T_{v'}(\pd V_i^n)^\ci]$ intersects $x_1'\!-\!x_0'$ positively in $\R^n$}\\
-1, & \!\!\!\text{
$T_{v'}(s_i^n\!\ci\! i_{V_i^n})[T_{v'}(\pd V_i^n)^\ci]$ intersects $x_1'\!-\!x_0'$ negatively in $\R^n$}\end{cases}}
\nonumber\\
&-\sum_{\begin{subarray}{l} v\in (V^n_i)^\ci: \\ s_i^n(v)=x_0' \end{subarray}}\!\cdot
\begin{cases} 1, & \!\!\!\text{$T_vs_i^n:T_v(V_i^n)^\ci\!\ra\! T_{x_0'}\R^n$ is orientation-preserving} \\ -1, & \!\!\!\text{$T_vs_i^n:T_v(V_i^n)^\ci\!\ra\! T_{x_0'}\R^n$ is orientation-reversing}
\end{cases}
\label{kh7eq55}\\
&+\sum_{\begin{subarray}{l} v\in (V^n_i)^\ci: \\ s_i^n(v)=x_1' \end{subarray}}\!\cdot
\begin{cases} 1, & \!\!\!\text{$T_vs_i^n:T_v(V_i^n)^\ci\!\ra\! T_{x_1'}\R^n$ is orientation-preserving} \\ -1, & \!\!\!\text{$T_vs_i^n:T_v(V_i^n)^\ci\!\ra\! T_{x_1'}\R^n$ is orientation-reversing.}
\end{cases}
\nonumber
\ea
Here the transverse fibre products in \eq{kh7eq54} are of dimension 0, so they are equal to their interiors by Assumption \ref{kh3ass4}(e), which are the fibre products of the interiors by Assumption \ref{kh3ass5}(d), so we may restrict to $v'\in (\pd V_i^n)^\ci$ and $\la\in(0,1)$ in the first sum in \eq{kh7eq55}, and to $v\in (V^n_i)^\ci$ in the second and third.

Multiply \eq{kh7eq55} by $a_i^n$ and sum over all $i\in I^n$. Using \eq{kh7eq52}, this yields
\ea
&\Phi^n(x_0')-\Phi^n(x_1')=\sum_{i\in I^n}a_i^n\cdot
\label{kh7eq56}\\
&\sum_{\begin{subarray}{l} \la\in(0,1),\\ v'\in (\pd V_i^n)^\ci:\\
s_i^n\ci i_{V_i^n}(v')=h(\la) \end{subarray}\!\!\!\!\!\!\!\!\!\!\!\!\!\!\!\!} \!\!{\begin{cases} 1, & \!\!\!\text{
$T_{v'}(s_i^n\!\ci\! i_{V_i^n})[T_{v'}(\pd V_i^n)^\ci]$ intersects $x_1'\!-\!x_0'$ positively in $\R^n$}\\
-1, & \!\!\!\text{
$T_{v'}(s_i^n\!\ci\! i_{V_i^n})[T_{v'}(\pd V_i^n)^\ci]$ intersects $x_1'\!-\!x_0'$ negatively in $\R^n$.}\end{cases}}
\nonumber
\ea
As above, condition Definition \ref{kh4def1}(ii)$(*)$ holds for \eq{kh7eq50} with $X=X^n$. Now the r.h.s.\ of \eq{kh7eq56} is the sum over all $\la\in(0,1)$ and all $(n-1)$-planes $P$ in $T_{h(\la)}\R^n$ transverse to $x_1'-x_0'$ and oriented so that $P$ intersects $x_1'-x_0'$ positively, of equation \eq{kh4eq2} in Definition \ref{kh4def1}(ii)$(*)$ for \eq{kh7eq50} at the point $h(\la)\in X^n$ and the oriented $(n-1)$-plane $P$. Therefore the r.h.s.\ of \eq{kh7eq56} is zero, and $\Phi^n(x_0')=\Phi^n(x_1')$. Since $\Phi^n(x_0')=\Phi^n(x_0)$, $\Phi^n(x_1')=\Phi^n(x_1)$, this gives $\Phi^n(x_0)=\Phi^n(x_1)$. Hence $\Phi^n:X^n\sm S^n\ra R$ is constant, proving the lemma.
\end{proof}

Lemma \ref{kh7lem1} implies that there exist $a^0,a^1,\ldots,a^N\in R$ with $\Phi^n(x)=a^n$ for all $n=0,\ldots,N$ and $x\in X^n\sm S^n$. We now claim that for $n=0,\ldots,N$ we have
\e
\sum_{i\in I^n}a_i^n\,\bigl[V_i^n,n,s_i^n,\pi\bigr]=a^n\bigl[\R^n,n,\id_{\R^n},\pi\bigr]\qquad\text{in $\widetilde{MC}_0(*;R)^n$.}
\label{kh7eq57}
\e
To see this, apply Definition \ref{kh4def1}(ii) to \eq{kh7eq57} with $X=X^n\subseteq\R^n$. The conditions on $(x,y)=(x,*)\in X\t Y=X\t *$ at the beginning of Definition \ref{kh4def1}(ii)$(*)$ are equivalent to $x\in X^n\sm S^n$, and then equation \eq{kh4eq2} at $x$ and $P=T_{(x,*)}(\R^n\t *)$ is equivalent to $\Phi^n(x)=a^n$, since $\Phi^n(x)$ in \eq{kh7eq52} is the contributions to \eq{kh4eq2} at $(x,*),P$ from $\sum_{i\in I^n}a_i^n\,\bigl[V_i^n,n,s_i^n,\pi\bigr]$, and $-a^n\bigl[\R^n,n,\id_{\R^n},\pi\bigr]$ contributes an additional $-a^n$ to \eq{kh4eq2} at $(x,*),P$. Hence \eq{kh4eq2} holds for all such $(x,*),P$, and Definition \ref{kh4def1}(ii) gives~\eq{kh7eq57}.

Applying $\Pi$ to project \eq{kh7eq57} to $MC_0(*;R)$ gives
\begin{align*}
\al&=\sum_{n=0}^N\sum_{i\in I^n}a_i^n\,\bigl[V_i^n,n,s_i^n,\pi\bigr]=\sum_{n=0}^Na^n\bigl[\R^n,n,\id_{\R^n},\pi\bigr]=\sum_{n=0}^Na^n\,[*]\\
&=\io(a^0+a^1+\cdots+a^N)\qquad\text{in $MC_0(*;R)$,}
\end{align*}
using \eq{kh7eq49} and that by Definition \ref{kh4def1}(i) in $MC_0(*;R)$ we have
\begin{equation*}
[*]=\bigl[\R^0,0,\id_{\R^0},\pi\bigr]=\bigl[\R,1,\id_{\R},\pi\bigr]=\bigl[\R^2,2,\id_{\R^2},\pi\bigr]=\cdots.
\end{equation*}
This proves that $\io$ in \eq{kh7eq49} is surjective.

To show that $\io$ is injective, suppose $r\in R$ with $\io(r)=0$. We will prove $r=0$. Regarding $r[*]$ as an element of $\widetilde{MC}_0(*;R)$ from \S\ref{kh741}, since $\Pi(r[*])=\io(r)=0$ in $MC_0(*;R)$, $r[*]$ lies in the kernel of $\Pi:\widetilde{MC}_0(*;R)\ra MC_0(*;R)$, which is spanned by equation \eq{kh4eq1} from relation Definition \ref{kh4def1}(i) in $MC_0(*;R)$. Thus as in \eq{kh7eq29}, in $\widetilde{MC}_0(*;R)$ we may write
\ea
r[*]&=\sum_{n=0}^N\sum_{l=0}^n\sum_{i\in I^{n,l}}a_i^{n,l}\Bigl(\bigl[V_i^{n,l},n,s_i^{n,l},\pi\bigr]-(-1)^{n-l}\bigl[V_i^{n,l}\t\R,
\nonumber\\
&\qquad n+1,(s_{i,1}^{n,l},\ldots,s_{i,l}^{n,l},\pi_\R,s_{i,l+1}^{n,l},\ldots,s_{i,n}^{n,l}),\pi\bigr]\Bigr).
\label{kh7eq58}
\ea
Taking components of \eq{kh7eq58} in $\widetilde{MC}_0(*;R)^n$ for $n=0,\ldots,N+1$ gives
\ea
r[*]&=\sum_{i\in I^{0,0}}a_i^{0,0}\bigl[V_i^{0,0},0,\pi_{\R^0},\pi\bigr],
\label{kh7eq59}\\
\begin{split}
0&=\sum_{l=0}^n\sum_{i\in I^{n,l}}a_i^{n,l}\bigl[V_i^{n,l},n,s_i^{n,l},\pi\bigr]\\
&-\sum_{l=0}^{n-1}\sum_{i\in I^{n-1,l}}
\begin{aligned}[t]&(-1)^{n-1-l}a_i^{n-1,l}\bigl[V_i^{n-1,l}\t\R,
n,(s_{i,1}^{n-1,l},\ldots,s_{i,l}^{n-1,l},\\
&\pi_\R,s_{i,l+1}^{n-1,l},\ldots,s_{i,n-1}^{n-1,l}),\pi\bigr],\quad n\!=\!1,\ldots,N,\!\!\!\!
\end{aligned}
\end{split}
\label{kh7eq60}\\
0&=-\sum_{l=0}^{N}\sum_{i\in I^{N,l}}
\begin{aligned}[t](-1)^{N-l}a_i^{N,l}\bigl[V_i^{N,l}\t\R,
N+1,(s_{i,1}^{N,l},\ldots,s_{i,l}^{N,l},\\
\pi_\R,s_{i,l+1}^{N,l},\ldots,s_{i,N}^{N,l}),\pi\bigr].
\end{aligned}
\label{kh7eq61}
\ea
By Proposition \ref{kh7prop2}, we may choose the representation \eq{kh7eq58} such that for all $n=0,\ldots,N$ and $l=0,\ldots,n$ we have 
\e
\sum_{i\in I^{n,l}}a_i^{n,l}\,\bigl[\pd V_i^{n,l},n,s_i^{n,l}\ci i_{V_i^{n,l}},\pi\bigr]=0\qquad\text{in $\widetilde{MC}_{-1}(*;R)^n$.}
\label{kh7eq62}
\e

The next part of the proof is very similar to the passage between equations \eq{kh7eq50} and \eq{kh7eq57} above. For each $n=0,\ldots,N$ and $l=0,\ldots,n$, choose an open neighbourhood $X^{n,l}$ of 0 in $\R^n$ such that $s_i^{n,l}:V_i^{n,l}\ra\R^n$ is proper over $X^{n,l}$ for all $i\in I^{n,l}$. Equation \eq{kh7eq62} holds as an application of relation Definition \ref{kh4def1}(ii) in $\widetilde{MC}_{-1}(*;R)$. Thus, making $X^{n,l}$ smaller we can suppose that Definition \ref{kh4def1}(ii)$(*)$ holds for \eq{kh7eq62} with $X=X^{n,l}$. Making $X^{n,l}$ still smaller, we can suppose it is a convex subset of~$\R^n$.

As in \eq{kh7eq51}, define a subset $S^{n,l}\subseteq X^{n,l}\subseteq\R^n$ by
\ea
&S^{n,l}=X^{n,l}\cap\Bigl[\ts\bigcup_{i\in I^{n,l}}s_i^{n,l}\bigl(\bigl\{v\in V_i^{n,l}:
\text{either $v\notin(V_i^{n,l})^\ci$, or $v\in (V_i^{n,l})^\ci$}
\nonumber\\
&\quad\text{and $T_vs_i^{n,l}:T_v(V_i^{n,l})^\ci\ra T_{s_i^{n,l}(v)}\R^n$ is not an isomorphism}\bigr\}\bigr)\Bigr].
\label{kh7eq63}
\ea
As the subsets of $V_i^{n,l}$ in \eq{kh7eq63} are closed, and $s_i^{n,l}$ is proper over $X^{n,l}$, we see that $S^{n,l}$ is closed in $X^{n,l}$. Using Sard's Theorem and Assumption \ref{kh3ass7}(a) we can show that $\cH^n(S^{n,l})=0$, which implies that~$X^{n,l}\sm S^{n,l}\ne\es$.

As in \eq{kh7eq52}, define a function $\Phi^{n,l}:X^{n,l}\sm S^{n,l}\ra R$ by
\ea
&\Phi^{n,l}(x)=
\label{kh7eq64}\\
&\sum_{i\in I^{n,l}}\sum_{\begin{subarray}{l} v\in (V^{n,l}_i)^\ci: \\ s_i^{n,l}(v)=x \end{subarray}}\!a_i^{n,l}\cdot
\begin{cases} 1, & \!\!\!\text{$T_vs_i^{n,l}:T_v(V_i^{n,l})^\ci\!\ra\! T_x\R^{n,l}$ is orientation-preserving,} \\ -1, & \!\!\!\text{$T_vs_i^{n,l}:T_v(V_i^{n,l})^\ci\!\ra\! T_x\R^n$ is orientation-reversing.}
\end{cases}
\nonumber
\ea
The proof of Lemma \ref{kh7lem1}, but using \eq{kh7eq62} in place of \eq{kh7eq50}, shows that $\Phi^{n,l}$ is constant. Hence there exist unique $a^{n,l}\in R$ with $\Phi^{n,l}(x)=a^{n,l}$ for all $n=0,\ldots,N$, $l=0,\ldots,n$ and $x\in X^{n,l}\sm S^{n,l}$.

For each $n=1,\ldots,N$, pick $(x_1,\ldots,x_n)\in\R^n$ satisfying the conditions:
\begin{itemize}
\setlength{\itemsep}{0pt}
\setlength{\parsep}{0pt}
\item[(a)] $(x_1,\ldots,x_n)\in X^{n,l}\sm S^{n,l}$ for all $l=0,\ldots,n$. 
\item[(b)] $(x_1,\ldots,x_l,x_{l+2},\ldots,x_n)\in X^{n-1,l}\sm S^{n-1,l}$ for all $l=0,\ldots,n-1$. 
\item[(c)] Equation \eq{kh7eq60} is an application of Definition \ref{kh4def1}(ii) in $\widetilde{MC}_0(*;R)$, which involves $0\in X\subseteq\R^n$ such that Definition \ref{kh4def1}(ii)$(*)$ holds for all suitable $(x,y)=(x,*)$ in $X\t Y=X\t *$. Then~$(x_1,\ldots,x_n)\in X$.
\end{itemize}
All these hold provided $(x_1,\ldots,x_n)\in\R^n$ is generic and small enough. Then (a)--(c) and \eq{kh7eq63} imply that $((x_1,\ldots,x_n),*)$ satisfies the conditions on $(x,y)$ in Definition \ref{kh4def1}(ii)$(*)$ for \eq{kh7eq60}, so equation \eq{kh4eq2} holds for \eq{kh7eq60} with $(x,y)=((x_1,\ldots,x_n),*)$ and~$P=T_x\R^n$. 

By \eq{kh7eq64}, the contributions to \eq{kh4eq2} from the terms $\sum_{i\in I^{n,l}}a_i^{n,l}\bigl[V_i^{n,l},\cdots\bigr]$ in \eq{kh7eq60} are $\Phi^{n,l}\bigl((x_1,\ldots,x_n)\bigr)=a^{n,l}$, which is defined by (a). Similarly, the contributions to \eq{kh4eq2} from the terms $\sum_{i\in I^{n-1,l}}(-1)^{n-1-l}a_i^{n-1,l}\bigl[V_i^{n-1,l}\t\R,\cdots\bigr]$ in \eq{kh7eq60} are $\Phi^{n-1,l}\bigl((x_1,\ldots,x_l,x_{l+2},\ldots,x_n)\bigr)=a^{n-1,l}$, which is defined by (b). Hence \eq{kh4eq2} implies that
\e
0=\ts\sum_{l=0}^na^{n,l}-\sum_{l=0}^{n-1}a^{n-1,l},\quad n=1,\ldots,N.
\label{kh7eq65}
\e

In the same way, equations \eq{kh7eq59} and \eq{kh7eq61} imply that
\ea
r&=a^{0,0},
\label{kh7eq66}\\
0&=-\ts\sum_{l=0}^Na^{N,l}.
\label{kh7eq67}
\ea
Taking the sum of equations \eq{kh7eq65} for $n=1,\ldots,N$ and \eq{kh7eq66}--\eq{kh7eq67}, the right hand sides cancel, yielding $r=0$. Therefore $\io$ is injective. We have now shown that $\io$ in \eq{kh7eq49} is surjective and injective, and thus an isomorphism, so by \eq{kh7eq48}--\eq{kh7eq49} we have $MH_0(*;R)\cong R$, completing the proof of Theorem~\ref{kh4thm2}.

\subsection{Proof of Proposition \ref{kh4prop4}}
\label{kh75}

We work in the situation of Definitions \ref{kh4def4} and \ref{kh4def5}. Let $[V',n,s',t']$ be as in \eq{kh4eq21}. Then we have a commutative diagram of topological spaces
\e
\begin{gathered}
\xymatrix@C=140pt@R=15pt{
*+[r]{V'} \ar@/^.7pc/[dr]^(0.75){(s',t')} \ar[d]^{(\pi_V,\pi_{Y_1})} \\
*+[r]{\bigl\{(v,y_1)\in V\t Y_1:t(v)=f(y_1)\bigr\}} \ar[r]_(0.7){s\t\id} \ar[d]^{\pi_V} 
& *+[l]{\R^n\t Y_1} \ar[d]_{\id\t f} \\
*+[r]{V} \ar[r]^{(s,t)} & *+[l]{\R^n\t Y_2.\!} }
\end{gathered}
\label{kh7eq68}
\e
Here the bottom square is Cartesian, and the bottom morphism $(s,t)$ is proper over an open neighbourhood $X$ of $\{0\}\t Y_2$ in $\R^n\t Y_2$, so the middle morphism $s\t\id$ is proper over $X':=(\id_{\R^n}\t f)^{-1}(X)$ in $\R^n\t Y_1$. Also the top left morphism $(\pi_V,\pi_{Y_1})$ is a homeomorphism by \eq{kh3eq8} in Assumption \ref{kh3ass5}(c). So from the top triangle in \eq{kh7eq68} we see that $(s',t'):V\ra\R^n\t Y_1$ is proper over $X'$, which is an open neighbourhood of $\{0\}\t Y_1$ in $\R^n\t Y_1$. Thus $[V',n,s',t']$ in \eq{kh4eq21} is a well-defined generator in $\cP MC^k(Y_1;R)$.

To see that $f^*$ is well-defined, we have to show that it maps relations Definition \ref{kh4def4}(i),(ii) in $\cP MC^k(Y_2;R)$ to relations (i),(ii) in $\cP MC^k(Y_1;R)$. For (i) this is obvious, since $(V\t\R)\t_{t\ci\pi_V,Y_2,f}Y_1\cong (V\t_{t,Y_2,f}Y_1)\t\R$. For (ii), suppose $\sum_{i\in I}a_i\,[V_i,n,s_i,t_i]=0$ in $\cP MC^k(Y_2;R)$ by relation (ii), using open $\{0\}\t Y_2\subseteq X_2\subseteq\R^n\t Y_2$. Define $X_1=(\id_{\R^n}\t f)^{-1}(X_2)\subseteq\R^n\t Y_1$. Then $X_1$ is an open neighbourhood of $\{0\}\t Y_1$ in $\R^n\t Y_1$. 

Set $[V_i',n,s_i',t_i']=f^*[V_i,n,s_i,t_i]$ as in \eq{kh4eq21} for $i\in I$. Then by Assumption \ref{kh3ass5}(c),(d) we have
\e
V_i^{\prime\ci}=V_i^\ci\t_{t\vert_{V_i^\ci},Y_2,f}Y_1\cong \bigl\{(v,y_1)\in V_i^\ci\t Y_1:t(v)=f(y_1)\bigr\}.
\label{kh7eq69}
\e
Suppose $(x,y_1)\in X_1$ such that for all $i\in I$ and $v'\in V_i'$ with $(s'_i,t'_i)(v')=(x,y_1)$ in $X_1$, we have $v'\in V_i^{\prime\ci}$ and
\e
T_{v'}(s'_i,t'_i):T_{v'}V_i^{\prime\ci}\longra T_x\R^n\op T_{y_1}Y_1
\label{kh7eq70}
\e
is injective. Set $y_2=f(y_1)$. For $v'$ as above, from \eq{kh7eq69} we have $v'=(v,y_1)$ for $v\in V_i^\ci$ with $t_i(v)=f(y_1)=y_2$ and $s_i(v)=s_i'(v')=x$, so that $(s_i,t_i)(v)=(x,y_2)$. Since \eq{kh7eq69} is a transverse fibre product of manifolds,
\e
\begin{gathered}
\xymatrix@C=100pt@R=15pt{ *+[r]{T_{v'}V_i^{\prime\ci}} \ar[d]^{T_{v'}t_i'} \ar[r]_{T_{v'}\pi_{V_i}} & *+[l]{T_vV_i^\ci} \ar[d]_{T_vt_i} \\
*+[r]{T_{y_1}Y_1} \ar[r]^{T_{y_1}f} & *+[l]{T_{y_2}Y_2} }
\end{gathered}
\label{kh7eq71}
\e
is Cartesian in $\mathop{\rm Vect}_\R$, and Assumption \ref{kh3ass6}(l) shows that the coorientations on $T_vt_i$ and $T_{v'}t_i'$ from $c_{t_i},c_{t_i'}$ are related as usual in~\eq{kh7eq71}. 

As \eq{kh7eq71} is Cartesian, we see that \eq{kh7eq70} is injective if and only if 
\begin{equation*}
T_v(s_i,t_i):T_vV_i\longra T_x\R^n\op T_{y_2}Y_2
\end{equation*} 
is injective. From all this we see that if $P'\subseteq T_x\R^n\op T_{y_1}Y_1$ is an $(m_1+n-k)$-plane with $\pi_{T_{y_1}Y_1}:P'\ra T_{y_1}Y_1$ cooriented, then either $P'$ lies in a Cartesian diagram
\e
\begin{gathered}
\xymatrix@C=90pt@R=15pt{ *+[r]{P'} \ar[d] \ar[r]_(0.3){\subseteq} & *+[l]{T_x\R^n\op T_{y_1}Y_1} \ar[d]_{\id\t T_{y_1}f} \\
*+[r]{P} \ar[r]^(0.3){\subseteq} & *+[l]{T_x\R^n\op T_{y_2}Y_2} }
\end{gathered}
\label{kh7eq72}
\e
in $\mathop{\rm Vect}_\R$ for some $(m_2+n-k)$-plane $P\subseteq T_x\R^n\op T_{y_2}Y_2$ with $\pi_{T_{y_2}Y_2}:P\ra T_{y_2}Y_2$ cooriented, in which case \eq{kh4eq17} for $\sum_{i\in I}a_i[V_i',n,s_i',t_i']$ at $(x,y_1),P'$ follows from \eq{kh4eq17} for $\sum_{i\in I}a_i[V_i,n,s_i,t_i]$ at $(x,y_2),P$, or else $P'$ lies in no such Cartesian diagram \eq{kh7eq72}, in which case \eq{kh4eq17} for $\sum_{i\in I}a_i[V_i',n,s_i',t_i']$ at $(x,y_1),P'$ is trivial as there are no $v'\in V_i'$ satisfying the conditions on either side.

This gives Definition \ref{kh4def4}(ii)$(*)$ for $\sum_{i\in I}a_i[V_i',n,s_i',t_i']$. Hence $\sum_{i\in I}a_i[V_i',\ab n,\ab s_i',\ab t_i']=0$ in $\cP MC^k(Y_1;R)$, and $f^*$ maps relation (ii) to relation (ii), so that $f^*$ is well-defined. This proves Proposition~\ref{kh4prop4}.

\subsection{Proof of Proposition \ref{kh4prop5}}
\label{kh76}

Here is the analogue for M-precochains of $\Pi_{T\cup U}^{T,f-},\Pi_{T\cup U}^{U,f+}$ in Definition \ref{kh7def1} in \S\ref{kh72}. Note that they map the opposite way: in \S\ref{kh72} we have $\Pi_{T\cup U}^{T,f-}:MC_k(T\cup U;R)\ra MC_k(T;R)$, but here we have $\Pi^{T\cup U}_{T,f-}:\cP MC^k(T;R)\ra \cP MC^k(T\cup U;R)$.

\begin{dfn} Let $Y$ be a manifold and $T,U\subseteq Y$ be open, and write $i:T\hookra T\cup U$, $j:U\hookra T\cup U$ for the inclusions. Suppose $f:T\cup U\ra\R$ is smooth, with $\bigl\{y\in T\cup U:f(y)\ge 0\bigr\}\subseteq T$ and $\bigl\{y\in T\cup U:f(y)\le 0\bigr\}\subseteq U$. As in Definition \ref{kh7def1}, such a function always exists.

Define $R$-linear morphisms $\Pi^{T\cup U}_{T,f-},\Pi^{T\cup U}_{U,f+}$ for each $k\in\Z$ by
\e
\begin{split}
&\Pi^{T\cup U}_{T,f-}:\cP MC^k(T;R)\longra \cP MC^k(T\cup U;R),\\
&\Pi^{T\cup U}_{T,f-}:\bigl[V,n,(s_1,\ldots,s_n),t\bigr]\longmapsto \bigl[V\t(-\iy,0],n+1,\\
&\quad (s_1\ci\pi_V,\ldots,s_n\ci\pi_V,f\ci\pi_V+\pi_{(-\iy,0]}),t\ci\pi_V\bigr],\\
&\Pi^{T\cup U}_{U,f+}:\cP MC^k(U;R)\longra \cP MC^k(T\cup U;R),
\\
&\Pi^{T\cup U}_{U,f+}:\bigl[V,n,(s_1,\ldots,s_n),t\bigr]\longmapsto \bigl[V\t[0,\iy),n+1,\\
&\quad (s_1\ci\pi_V,\ldots,s_n\ci\pi_V,f\ci\pi_V+\pi_{[0,\iy)}),t\ci\pi_V\bigr],
\end{split}
\label{kh7eq73}
\e
for each generator $\bigl[V,n,(s_1,\ldots,s_n),t\bigr]$ in $\cP MC^k(T;R)$ and $\cP MC^k(U;R)$.
\label{kh7def3}
\end{dfn}

Here is the analogue of Proposition \ref{kh7prop1}, with a very similar proof which we leave as an exercise. For (i), to show that $\bigl[V\t(-\iy,0],n+1,(s_1\ci\pi_V,\ldots,s_n\ci\pi_V,f\ci\pi_V+\pi_{(-\iy,0]}),t\ci\pi_V\bigr]$ in \eq{kh7eq73} is a generator of $\cP MC^k(T\cup U;R)$, note that $\bigl((s_1,\ldots,s_n),t\bigr):V\ra\R^n\t T$ proper near $\{0\}\t T$ in $\R^n\t T$ implies that 
\begin{equation*}
\bigl((s_1\ci\pi_V,\ldots,s_n\ci\pi_V,f\ci\pi_V\!+\!\pi_{(-\iy,0]}),t\ci\pi_V\bigr):V\!\t\!(-\iy,0]\!\longra\!\R^{n+1}\!\t \!(T\cup U)
\end{equation*}
is proper near $\{0\}\t (T\cup U)$ in $\R^{n+1}\t (T\cup U)$.

\begin{prop}{\bf(i)} In the above, $\Pi^{T\cup U}_{T,f-}$ and\/ $\Pi^{T\cup U}_{U,f+}$ are well defined.
\smallskip

\noindent{\bf(ii)} We have $\Pi^{T\cup U}_{T,f-}\ci i^*+\Pi^{T\cup U}_{U,f+}\ci j^*=\id:\cP MC^k(T\cup U;R)\ra \cP MC^k(T\cup U;R)$.

\smallskip

\noindent{\bf(iii)} If\/ $[V,n,s,t]$ is a generator of\/ $\cP MC^k(T;R)$ with\/ $f\ci t(v)>0$ for all\/ $v$ in $s^{-1}(0)\subseteq V$ then $i^*\ci\Pi^{T\cup U}_{T,f-}\bigl([V,n,s,t]\bigr)=[V,n,s,t]$.

Similarly, if\/ $[V,n,s,t]\!\in\! \cP MC^k(U;R)$ with\/ $f\ci t(v)\!<\!0$ for all\/ $v$ in $s^{-1}(0)\!\subseteq\! V$ then $j^*\ci\Pi^{T\cup U}_{U,f+}\bigl([V,n,s,t]\bigr)=[V,n,s,t]$.
\label{kh7prop4}
\end{prop}

To prove Proposition \ref{kh4prop5}(a), suppose $T,U\subseteq Y$ are open, and write $i:T\cap U\hookra T,$ $i':T\cap U\hookra U,$ $j:T\hookra T\cup U,$ $j':U\hookra T\cup U$ for the inclusions. Choose $f$ as in Definition \ref{kh7def3}, so that we have operators $\Pi^{T\cup U}_{T,f-},\Pi^{T\cup U}_{U,f+}$. Applying Definition \ref{kh7def3} with $\ti T=T$, $\ti U=T\cap U$, $\ti f=f\vert_T$ in place of $T,U,f$ gives operators 
\begin{align*}
&\Pi^T_{T,f\vert_T-}:\cP MC^k(T;R)\longra \cP MC^k(T;R),\\
&\Pi^T_{T\cap U,f\vert_T+}:\cP MC^k(T\cap U;R)\longra \cP MC^k(T;R),
\end{align*}
and applying Proposition \ref{kh7prop4}(ii) for these gives
\e
\Pi^T_{T,f\vert_T-}+\Pi^T_{T\cap U,f\vert_T+}\ci i^*=\id:\cP MC^k(T;R)\longra \cP MC^k(T;R).
\label{kh7eq74}
\e
Similarly we have operators
\ea
&\Pi^U_{T\cap U,f\vert_U-}:\cP MC^k(T\cap U;R)\longra \cP MC^k(U;R),
\nonumber\\
&\Pi^U_{U,f\vert_U+}:\cP MC^k(U;R)\longra \cP MC^k(U;R),\quad\text{with}
\nonumber\\
&\Pi^U_{T\cap U,f\vert_U-}\ci i'_*+\Pi^U_{U,f\vert_U+}=\id:\cP MC^k(U;R)\longra \cP MC^k(U;R).
\label{kh7eq75}
\ea

As in \eq{kh7eq18}--\eq{kh7eq21}, comparing the actions of the two sides of each equation on generators $[V,n,s,t]$ using \eq{kh4eq21} and \eq{kh7eq73}, we see that
\ea
\Pi^T_{T,f\vert_T-}&=j^*\ci\Pi^{T\cup U}_{T,f-}:\cP MC^k(T;R)\longra \cP MC^k(T;R),
\label{kh7eq76}\\
\Pi^U_{U,f\vert_T+}&=j^{\prime *}\ci\Pi_{T\cup U}^{U,f+}:\cP MC^k(U;R)\longra \cP MC^k(U;R),
\label{kh7eq77}\\
\Pi^U_{T\cap U,f\vert_U-}\ci i^*&=j^{\prime *}\ci\Pi^{T\cup U}_{T,f-}:\cP MC^k(U;R)\longra \cP MC^k(T;R),
\label{kh7eq78}\\
\Pi^T_{T\cap U,f\vert_T+}\ci i^{\prime *}&=j^*\ci\Pi^{T\cup U}_{U,f+}:\cP MC^k(T;R)\longra \cP MC^k(U;R).
\label{kh7eq79}
\ea

As in \eq{kh4eq22}, we have to prove the following sequence is exact:
\e
\xymatrix@C=9.5pt{ 0 \ar[r] & \cP MC^k(T\!\cup\! U;R) \ar[rr]^{j^*\op j^{\prime *}} && {\begin{subarray}{l} \ts \; \cP MC^k(T;R)\\ \ts\op\cP MC^k(U;R) \end{subarray}} \ar[rr]^(0.45){i^*\op -i^{\prime *}} && \cP MC^k(T\!\cap\! U;R). }
\label{kh7eq80}
\e
For any $\al\in MC^k(T\cup U;R)$, Proposition \ref{kh7prop4}(ii) implies that
\begin{equation*}
\Pi^{T\cup U}_{T,f-}\ci j^*(\al)+\Pi^{T\cup U}_{U,f+}\ci j^{\prime *}(\al)=\al,
\end{equation*}
so $j^*(\al)\op j^{\prime *}(\al)=0$ implies that $\al=0$, and \eq{kh7eq80} is exact at the second term. 

To show \eq{kh7eq80} is exact at the third term, suppose $\be\in \cP MC^k(T;R)$ and $\ga\in \cP MC^k(U;R)$ with $i^*(\be)=i^{\prime *}(\ga)$ in $\cP MC^k(T\cap U;R)$. Define $\al=\Pi^{T\cup U}_{T,f-}(\be)+\Pi^{T\cup U}_{U,f+}(\ga)$ in $\cP MC^k(T\cup U;R)$. Then
\begin{align*}
j^*(\al)&=j^*\ci\Pi^{T\cup U}_{T,f-}(\be)+j^*\ci\Pi^{T\cup U}_{U,f+}(\ga)=\Pi^T_{T,f\vert_T-}(\be)+\Pi^T_{T\cap U,f\vert_T+}\ci i^{\prime *}(\ga)\\
&=\Pi^T_{T,f\vert_T-}(\be)+\Pi^T_{T\cap U,f\vert_T+}\ci i^*(\be)=\be
\end{align*}
using the definition of $\al$ in the first step, \eq{kh7eq76} and \eq{kh7eq79} in the second, $i^*(\be)=i^{\prime *}(\ga)$ in the third, and \eq{kh7eq74} in the fourth. A similar proof using \eq{kh7eq77}, \eq{kh7eq78} and \eq{kh7eq75} shows that $j^{\prime *}(\al)=\ga$. Hence \eq{kh7eq80} is exact, proving Proposition~\ref{kh4prop5}(a).

For part (b), suppose $K\subseteq Y$ is closed, and $U$ is an open neighbourhood of $K$ in $Y$. Apply Definition \ref{kh7def3} with $T=Y\sm K$, so that $T\cup U=Y$, and we choose smooth $f:Y\ra\R$ with $K\subseteq f^{-1}\bigl((-\iy,0)\bigr)$ and $f^{-1}\bigl((-\iy,0]\bigr)\subseteq U$. Define $U'=f^{-1}\bigl((-\iy,0)\bigr)$, so that $U'$ is an open neighbourhood of $K$ in $U$ and $f>0$ on $U'$. Write $i:U'\hookra U$ and  $j:U'\hookra Y$ for the inclusions.

Then Definition \ref{kh7def3} gives $\Pi^Y_{U,f+}:\cP MC^k(U;R)\ra \cP MC^k(Y;R)$. Also Definition \ref{kh7def3} with $\ti T=\ti U=U'$ and $\ti f=f\vert_{U'}$ in place of $T,U,f$ gives
\begin{equation*}
\Pi^{U'}_{U',f\vert_{U'}+}:\cP MC^k(U';R)\longra \cP MC^k(U';R).
\end{equation*}
Comparing the actions of both sides on generators $[V,n,s,t]$ in $\cP MC^k(U;R)$ using \eq{kh4eq21} and \eq{kh7eq73}, we see that
\e
j^*\ci\Pi^Y_{U,f+}=\Pi^{U'}_{U',f\vert_{U'}+}\ci i^*:\cP MC^k(U;R)\longra \cP MC^k(U';R).
\label{kh7eq81}
\e
But as $f<0$ on $U'$, Proposition \ref{kh7prop4}(iii) implies that $\Pi^{U'}_{U',f\vert_{U'}+}=\id_{U'}^*\ci\Pi^{U'}_{U',f\vert_{U'}+}$ maps $[V,n,s,t]\mapsto [V,n,s,t]$, so that
\e
\Pi^{U'}_{U',f\vert_{U'}+}=\id:\cP MC^k(U';R)\longra \cP MC^k(U';R).
\label{kh7eq82}
\e
Now let $\al\in\cP MC^k(U;R)$, and set $\be=\Pi^Y_{U,f+}(\al)\in\cP MC^k(Y;R)$. Then
\begin{equation*}
i^*(\al)=\Pi^{U'}_{U',f\vert_{U'}+}\ci i^*(\al)=j^*\ci\Pi^Y_{U,f+}(\al)=
j^*(\be),
\end{equation*}
using \eq{kh7eq82} in the first step, \eq{kh7eq81} in the second, and the definition of $\be$ in the third. This proves Proposition~\ref{kh4prop5}(b).

\subsection{Proof of Proposition \ref{kh4prop6}}
\label{kh77}

First suppose that $Y_1,Y_2$ are manifolds and $g:Y_1\t[0,1]\ra Y_2$ is smooth. For all $k\in\Z$, define $G_{Y_1,Y_2}:\cP MC^k(Y_2;R)\ra\cP MC^{k-1}(Y_1;R)$ to be the unique $R$-linear map acting on generators by
\e
\begin{split}
G_{Y_1,Y_2}:[V,n&,s,t]\longmapsto (-1)^{\dim V}[V',n,s',t']:=\\
&(-1)^{\dim V}\bigl[V\t_{t,Y_2,g}(Y_1\t[0,1]),n,s\ci\pi_V,\pi_{Y_1}\ci\pi_{Y_1\t[0,1]}\bigr],
\end{split}
\label{kh7eq83}
\e
where $\pi_{Y_1}\ci\pi_{Y_1\t[0,1]}$ has the coorientation $c_{\pi_{Y_1}\ci\pi_{Y_1\t[0,1]}}=c_{\pi_{Y_1}}\ci c_{\pi_{Y_1\t[0,1]}}$ as in Assumption \ref{kh3ass6}(d), with $c_{\pi_{Y_1}}$ the coorientation on $\pi_{Y_1}:Y_1\t[0,1]\ra Y_1$ induced by the standard orientation on $[0,1]$ as in Assumption \ref{kh3ass6}(f),(k), and $c_{\pi_{Y_1\t[0,1]}}$ is the coorientation on $\pi_{\pi_{Y_1\t[0,1]}}:V'\ra Y_1\t[0,1]$ induced from the given coorientation $c_t$ on $t:V\ra Y_2$ by Assumption~\ref{kh3ass6}(l).

As $(s,t):V\ra\R^n\t Y_2$ is proper over an open neighbourhood $X$ of $\{0\}\t Y_2$ in $\R^n\t Y_2$, and fibre products in $\tManc$ map to fibre products in $\Top$, so $(s',\pi_{Y_1\t[0,1]}):V'\ra\R^n\t Y_1\t[0,1]$ is proper over $(\id_{\R^n}\t g)^{-1}(X)$, an open neighbourhood of $\{0\}\t Y_1\t[0,1]$ in $\R^n\t Y_1\t[0,1]$. By compactness of $[0,1]$ there exists an open neighbourhood $X'$ of $\{0\}\t Y_1$ in $\R^n\t Y_1$ with $X'\t[0,1]\subseteq(\id_{\R^n}\t g)^{-1}(X)$, and then $(s',t'):V'\ra\R^n\t Y_1$ is proper over $X'$, so $[V',n,s',t']$ is a generator of $\cP MC^{k-1}(Y_1;R)$.

To show that $G_{Y_1,Y_2}$ is well-defined, we must show that it maps relations Definition \ref{kh4def4}(i),(ii) in $\cP MC^k(Y_2;R)$ to relations (i),(ii) in $\cP MC^{k-1}(Y_1;R)$. This can be done by a proof similar to those for $f_*$ on M-chains in Definition \ref{kh4def2} and $f^*$ on M-precochains in Proposition \ref{kh4prop4}. In $\cP MC^k(Y_1;R)$ we have
\begin{align*}
\d&\ci G_{Y_1,Y_2}[V,n,s,t]=(-1)^{\dim V}[\pd V',n,s'\ci i_{V'},t'\ci i_{V'}]
\\
&=(-1)^{\dim V}\bigl[\pd V\t_{t\ci i_V,Y_2,g}(Y_1\t[0,1]),n,s\ci\pi_{\pd V},\pi_{Y_1}\ci\pi_{Y_1\t[0,1]}\bigr]\\
&\quad +\bigl[V\t_{t,Y_2,g\ci i_{Y_1\t [0,1]}}(\pd(Y_1\t[0,1])),n,s\ci\pi_V,\pi_{Y_1}\ci i_{Y_1\t[0,1]}\ci\pi_{\pd(Y_1\t[0,1])}\bigr]
\\
&=-G_{Y_1,Y_2}[\pd V,n,s\ci i_V,t\ci i_V]-[V\t_{t,Y_2,f}Y_1,n,s\ci\pi_V,\pi_{Y_1}]\\
&\quad +[V\t_{t,Y_2,f'}Y_1,n,s\ci\pi_V,\pi_{Y_1}]
\\
&=(-G_{Y_1,Y_2}\ci\d-f_*+f'_*)[V,n,s,t],
\end{align*}
using \eq{kh4eq20} and \eq{kh7eq83} in the first step, Assumptions \ref{kh3ass5}(c) and \ref{kh3ass6}(m) and \eq{kh4eq18}--\eq{kh4eq19} in the second, $\pd(Y_1\t[0,1])=-(Y_1\t\{0\})\amalg (Y_1\t\{1\})$ where the signs compare the natural coorientations for $\pi_{Y_1}\ci i_{Y_1\t[0,1]}:\pd(Y_1\t[0,1])\ra Y_1$ and for the obvious \'etale map $(Y_1\t\{0\})\amalg (Y_1\t\{1\})\ra Y_1$, equations \eq{kh4eq18}--\eq{kh4eq19}, and the definitions of $f,f',G_{Y_1,Y_2}$ in the third, and \eq{kh4eq20}--\eq{kh4eq21} in the fourth. As this holds for all generators $[V,n,s,t]$, we have
\begin{equation*}
\d\ci G_{Y_1,Y_2}+G_{Y_1,Y_2}\ci\d=f^{\prime *}-f^*:\cP MC^k(Y_2;R)\longra\cP MC^k(Y_1;R).
\end{equation*}

To promote $G_{Y_1,Y_2}$ from spaces $\cP MC^k(Y_i;R)$ to spaces $MC^k(Y_i;R)$, suppose $U_1\subseteq Y_1$ is open with $\bar U_1$ compact. Then $\bar U_1\t[0,1]$ is compact in $Y_1\t[0,1]$, so $g\bigl(\bar U_1\t[0,1]\bigr)$ is compact in $Y_2$, and we can choose an open neighbourhood $U_2$ of $g\bigl(\bar U_1\t[0,1]\bigr)$ in $Y_2$ with $\bar U_2$ compact. Then $g\vert_{U_1\t[0,1]}$ maps $U_1\t[0,1]\ra U_2$, and is a homotopy between $f\vert_{U_1}:U_1\ra U_2$ and $f'\vert_{U_1}:U_1\ra U_2$. Hence the definition above with $g\vert_{U_1\t[0,1]},f\vert_{U_1},f'\vert_{U_1}$ in place of $g,f,f'$ yields $G_{U_1,U_2}:\cP MC^k(U_2;R)\ra\cP MC^{k-1}(U_1;R)$ with
\e
\d\ci G_{U_1,U_2}+G_{U_1,U_2}\ci\d=f\vert_{U_1}^{\prime *}-f\vert_{U_1}^*:\cP MC^k(U_2;R)\ra\cP MC^k(U_1;R).
\label{kh7eq84}
\e

Comparing \eq{kh4eq21} and \eq{kh7eq83}, we see that such $G_{Y_1,Y_2},G_{U_1,U_2}$ are compatible with restrictions to open sets $U_1\subseteq Y_1$, $U_2\subseteq Y_2$. Thus as in the definition of $f^*$ in Definition \ref{kh4def6}, using the characterization \eq{kh4eq23} of $MC^{k-1}(Y_1;R)$ as an inverse limit, we deduce that there is a unique morphism $G_{Y_1,Y_2}:MC^k(Y_2;R)\ra MC^{k-1}(Y_1;R)$ such that the following commutes for all such $U_1,U_2$:
\begin{equation*}
\xymatrix@C=100pt@R=15pt{ *+[r]{MC^k(Y_2;R)} \ar[r]_{\tau_{Y_2U_2}} \ar[d]^{G_{Y_1,Y_2}} &
*+[l]{\cP MC^k(U_2;R)} \ar[d]_{G_{U_1,U_2}} \\
*+[r]{MC^{k-1}(Y_1;R)} \ar[r]^{\tau_{Y_1U_1}} & *+[l]{\cP MC^{k-1}(U_1;R).\!\!}  }
\end{equation*}
Here $\tau_{Y_aU_a}$ is as in Theorem \ref{kh2thm4}(c) for the strong presheaf $\cP\MC^*(Y_a;R)$ with sheafification $\MC^*(Y_a;R)$, which as in Theorem \ref{kh2thm4}(e) is the projection from the inverse limit \eq{kh4eq23} for $MC^*(Y_a;R)$. 

Now $\Pi\ci G_{Y_1,Y_2}=G_{Y_1,Y_2}\ci\Pi:\cP MC^k(Y_2;R)\ra MC^{k-1}(Y_1;R)$, so $G_{Y_1,Y_2}$ on generators $[V,n,s,t]$ in $MC^k(Y_2;R)$ is given by \eq{kh7eq83}. Equation \eq{kh7eq84} gives
\e
\d\ci G_{Y_1,Y_2}+G_{Y_1,Y_2}\ci\d=f^{\prime *}-f^*:MC^k(Y_2;R)\longra MC^k(Y_1;R).
\label{kh7eq85}
\e
So $G_{Y_1,Y_2}$ is a cochain homotopy from $f^*$ to $f^{\prime *}$ on M-cochains, and thus $f^*=f^{\prime *}:MH^k(Y_2;R)\ra MH^k(Y_1;R)$ on M-cohomology.

To extend this to relative M-cohomology, suppose as in the proposition that $Z_1\subseteq Y_1,$ $Z_2\subseteq Y_2$ are open with $g(Z_1\t[0,1])\subseteq Z_2$, so that $g\vert_{Z_1\t[0,1]}$ maps $Z_1\t[0,1]\ra Z_2$. Applying the above argument to $g\vert_{Z_1\t[0,1]}$ gives maps $G_{Z_1,Z_2}:MC^k(Z_2;R)\ra MC^{k-1}(Z_1;R)$ with
\e
\d\ci G_{Z_1,Z_2}+G_{Z_1,Z_2}\ci\d=f\vert_{Z_1}^{\prime *}-f\vert_{Z_1}^*:MC^k(Z_2;R)\longra MC^k(Z_1;R).
\label{kh7eq86}
\e
Compatibility with pullbacks implies that
\e
i_1^*\ci G_{Y_1,Y_2}=G_{Z_1,Z_2}\ci i_2^*:MC^k(Y_2;R)\longra MC^{k-1}(Z_1;R),
\label{kh7eq87}
\e
writing $i_a:Z_a\hookra Y_a$ for the inclusion for $a=1,2$. Define $G:MC^k(Y_2,Z_2;R)\ra MC^{k-1}(Y_1,Z_1;R)$ by $G:(\al,\be)\mapsto \bigl(G_{Y_1,Y_2}(\al),G_{Z_1,Z_2}(\be)\bigr)$. Then \eq{kh7eq85}--\eq{kh7eq87} imply that 
\begin{equation*}
\d\ci G+G\ci\d=f^{\prime *}-f^*:MC^k(Y_2,Z_2;R)\longra MC^k(Y_1,Z_1;R).
\end{equation*}
So $G$ is a cochain homotopy from $f^*$ to $f^{\prime *}$ on relative M-cochains, and thus $f^*=f^{\prime *}:MH^k(Y_2,Z_2;R)\ra MH^k(Y_1,Z_1;R)$ on relative M-cohomology. This proves  Proposition~\ref{kh4prop6}.

\subsection{Proof of Proposition \ref{kh4prop8}}
\label{kh78}

Let $Y$ be a manifold and $k\in\Z$. Then \S\ref{kh42} defines the $R$-modules $\cP MC^k(Y;R)$ and $MC^k(Y;R)$ and the projection $\Pi:\cP MC^k(Y;R)\ra MC^k(Y;R)$, and \S\ref{kh43} defines the $R$-submodules $\cP MC^k_\cs(Y;R)\subseteq\cP MC^k(Y;R)$ and $MC^k_\cs(Y;R)\subseteq MC^k(Y;R)$ of compactly-supported sections, where as in \eq{kh4eq41} the restriction
\e
\Pi\vert_{\cP MC^k_\cs(Y;R)}:\cP MC^k_\cs(Y;R)\longra MC^k_\cs(Y;R)
\label{kh7eq88}
\e
is an isomorphism. Note that as in Definition \ref{kh2def6}, as $\cP\MC^k(Y;R)$ is a presheaf, for $\al\in\cP MC^k(Y;R)$ to be compactly-supported means that there exist compact $K\subseteq Y$ with $i^*(\al)=0$ in $\cP MC^k(Y\sm K;R)$ for $i:Y\sm K\hookra Y$ the inclusion. This implies that $\supp\al$ is compact, but may not be equivalent to it.

Let $\widehat{MC}{}^k_\cs(Y;R)$ be the $R$-module spanned by compact generators $[V,n,s,t]$ in the sense of Definition \ref{kh4def10}, subject to relations Definition \ref{kh4def4}(i),(ii) applied to compact generators. Then we have a natural $R$-linear projection
\e
\Pi:\widehat{MC}{}^k_\cs(Y;R)\longra\cP MC^k_\cs(Y;R)
\label{kh7eq89}
\e
mapping $\Pi:[V,n,s,t]\mapsto[V,n,s,t]$ on (compact) generators, since every (compact) generator in $\widehat{MC}{}^k_\cs(Y;R)$ is also a generator in $\cP MC^k(Y;R)$ which is compactly-supported and so lies in $\cP MC^k_\cs(Y;R)\subseteq\cP MC^k(Y;R)$, and every relation in $\widehat{MC}{}^k_\cs(Y;R)$ is also a relation in~$\cP MC^k(Y;R)$.

We will show \eq{kh7eq89} is an isomorphism. As \eq{kh7eq88} is also an isomorphism, this will imply that $\widehat{MC}{}^k_\cs(Y;R)\cong MC^k_\cs(Y;R)$, proving Proposition \ref{kh4prop8}. Suppose $\al\in\cP MC^k_\cs(Y;R)$. Then there exists compact $K\subseteq Y$ with $i^*(\al)=0$ for $i:T\hookra Y$ the inclusion, where $T=Y\sm K$. Choose an open neighbourhood $U$ of $K$ in $Y$ with compact closure $\bar U$ in $Y$, and write $j:U\hookra Y$ for the inclusion. Then $Y=T\cup U$, so as in Definition \ref{kh7def3} we can choose smooth $f:Y\ra\R$ with $K\subseteq f^{-1}\bigl((-\iy,0)\bigr)$ and $f^{-1}\bigl((-\iy,0]\bigr)\subseteq U$, and define $R$-linear maps
\e
\begin{split}
&\Pi^Y_{T,f-}:\cP MC^k(T;R)\longra \cP MC^k(Y;R),\\
&\Pi^Y_{U,f+}:\cP MC^k(U;R)\longra \cP MC^k(Y;R),
\end{split}
\label{kh7eq90}
\e
which by Proposition \ref{kh7prop4}(ii) satisfy
\begin{equation*}
\Pi^Y_{T,f-}\ci i^*+\Pi^Y_{U,f+}\ci j^*=\id:\cP MC^k(Y;R)\longra \cP MC^k(Y;R).
\end{equation*}
As $i^*(\al)=0$, applying this to $\al$ gives $\al=\Pi^Y_{U,f+}\ci j^*(\al)$.

Suppose $[V,n,s,t]$ is a generator of $\cP MC^k(U;R)$, so that by \eq{kh7eq73} we have
\begin{align*}
\Pi^{T\cup U}_{U,f+}\bigl([V,n,s,t]\bigr)&=[\ti V,n+1,\ti s,\ti t]=\bigl[V\t[0,\iy),n+1,\\
&\quad (s_1\ci\pi_V,\ldots,s_n\ci\pi_V,f\ci\pi_V+\pi_{[0,\iy)}),t\ci\pi_V\bigr].
\end{align*}
Then $(\ti s,\ti t):\ti V\ra\R^{n+1}\t Y$ is proper over an open neighbourhood of $\{0\}\t Y$ in $\R^{n+1}\t Y$, as $[\ti V,n+1,\ti s,\ti t]$ is a generator of $\cP MC^k(Y;R)$. But $\ti t(\ti V)\subseteq U\subseteq\bar U\subseteq Y$, where $\bar U$ is compact, so projecting $\R^{n+1}\t Y\ra\R^{n+1}$ shows that $\ti s:\ti V\ra\R^{n+1}$ is proper near 0 in $\R^{n+1}$, so $[\ti V,n+1,\ti s,\ti t]$ is a compact generator of $\cP MC^k(Y;R)$, as in Definition \ref{kh4def10}, and lies in the image of $\Pi$ in \eq{kh7eq89}. Hence $\Pi^Y_{U,f+}$ in \eq{kh7eq90} maps into the image of \eq{kh7eq89}. Since $\al=\Pi^Y_{U,f+}\ci j^*(\al)$, it lies in the image of \eq{kh7eq89}, and thus $\Pi$ in \eq{kh7eq89} is surjective.

Next suppose $\hat\al\in\widehat{MC}{}^k_\cs(Y;R)$ with $\Pi(\hat\al)=0$ in $\cP MC^k_\cs(Y;R)$. Write
\e
\hat\al=\ts\sum_{a\in A}\hat\al_a\,[V_a,n_a,s_a,t_a],
\label{kh7eq91}
\e
where $A$ is a finite indexing set, $\al_a\in R$ and $[V_a,n_a,s_a,t_a]$ is a compact generator. Set $K=\bigcup_{a\in A}t_a(s_a^{-1}(0))\subseteq Y$, which is compact as each $t_a(s_a^{-1}(0))$ is compact. As above choose an open neighbourhood $U$ of $K$ in $Y$ with compact closure $\bar U$ in $Y$, and write $T=Y\sm K$ and $i:T\hookra Y$, $j:U\hookra Y$ for the inclusions. Choose smooth $f:Y\ra\R$ with $K\subseteq f^{-1}\bigl((-\iy,0)\bigr)$ and $f^{-1}\bigl((-\iy,0]\bigr)\subseteq U$, so that Definition \ref{kh7def3} defines $\Pi^Y_{T,f-},\Pi^Y_{U,f+}$ as in~\eq{kh7eq90}.

For functoriality under inclusions of open sets $j:U\hookra Y$, the spaces $\widehat{MC}{}^k_\cs(Y;R)$ behave like M-chains $MC_k(Y;R)$ in \S\ref{kh41}, or compactly-supported M-cochains $MC^k_\cs(Y;R)$ in \S\ref{kh43}, rather than like M-(pre)cochains $\cP MC^k(Y;R),\ab MC^k(Y;R)$ in \S\ref{kh42}. So as in \S\ref{kh42} and \S\ref{kh43} we may define the pushforward 
\begin{align*}
&j_*:\widehat{MC}{}^k_\cs(U;R)\longra \widehat{MC}{}^k_\cs(Y;R), \\
&j_*:[V,n,s,t]\longmapsto [V,n,s,t],
\end{align*}
and as in Definition \ref{kh7def1} for $MC_k(-;R)$ we may define
\begin{align*}
&\hat\Pi_Y^{U,f+}:\widehat{MC}{}^k_\cs(Y;R)\longra \widehat{MC}{}^k_\cs(U;R),
\\
&\hat\Pi_Y^{U,f+}:\bigl[V,n,(s_1,\ldots,s_n),t\bigr]\longmapsto \bigl[t^{-1}(U)\t[0,\iy),n+1,\\
&\quad (s_1\ci\pi_{t^{-1}(U)},\ldots,s_n\ci\pi_{t^{-1}(U)},f\ci\pi_{t^{-1}(U)}+\pi_{[0,\iy)}),t\ci\pi_{t^{-1}(U)}\bigr].
\end{align*}
Then comparing the actions on compact generators $[V,n,s,t]$ in $\widehat{MC}{}^k_\cs(Y;R)$, we see that the following rectangle commutes:
\e
\begin{gathered}
\xymatrix@C=70pt@R=17pt{ *+[r]{\widehat{MC}{}^k_\cs(Y;R)} \ar[d]^\Pi \ar[r]_(0.55){\hat\Pi_Y^{U,f+}} & \widehat{MC}{}^k_\cs(U;R) \ar[r]_(0.3){j_*} & *+[l]{\widehat{MC}{}^k_\cs(Y;R)} \ar[d]_\Pi \\
*+[r]{\cP MC^k(Y;R)} \ar[r]^(0.7){j^*} \ar@{.>}[urr]^(0.4)\Xi & \cP MC^k(U;R) \ar[r]^(0.45){\Pi^Y_{U,f+}} & *+[l]{\cP MC^k(Y;R).\!} }
\end{gathered}
\label{kh7eq92}
\e

In fact we can say more. As in the previous part, the operator $\Pi^Y_{U,f+}\ci j^*$ maps generators $[V,n,s,t]$ in $\cP MC^k(Y;R)$ to compact generators $[\ti V,n+1,\ti s,\ti t]$ in $\cP MC^k(Y;R)$, which thus lift to $\widehat{MC}{}^k_\cs(Y;R)$. The proofs in Propositions \ref{kh4prop4} and \ref{kh7prop4} that $j^*,\Pi^Y_{U,f+}$ are well defined show that they map relations (i),(ii) to relations (i),(ii), and checking the details we find that $\Pi^Y_{U,f+}\ci j^*$ maps relations (i),(ii) in possibly noncompact generators to relations (i),(ii) in compact generators, which lift to relations in $\widehat{MC}{}^k_\cs(Y;R)$. Therefore $\Pi^Y_{U,f+}\ci j^*$ factors as $\Pi\ci\Xi$ for $\Xi:\cP MC^k(Y;R)\ra \widehat{MC}{}^k_\cs(Y;R)$ making \eq{kh7eq92} commute.

Since $f>0$ on $t_a(s_a^{-1}(0))\subseteq Y$ by choice of $K,f$, the analogue of Proposition \ref{kh7prop1}(iii) for $\widehat{MC}{}^k_\cs(-;R),\hat\Pi_Y^{U,f+}$ shows that $j_*\ci\hat\Pi_Y^{U,f+}\bigl([V_a,n_a,s_a,t_a]\bigr)=[V_a,n_a,s_a,t_a]$ for $a\in A$, and so $j_*\ci\hat\Pi_Y^{U,f+}(\hat\al)=\hat\al$ by \eq{kh7eq91}, so \eq{kh7eq92} yields
\begin{equation*}
\hat\al=j_*\ci\hat\Pi_Y^{U,f+}(\hat\al)=\Xi\ci\Pi(\hat\al)=0.
\end{equation*}
Hence $\Pi$ in \eq{kh7eq89} is injective, and an isomorphism. This completes the proof of Proposition~\ref{kh4prop8}.

\subsection{Proof of Proposition \ref{kh4prop11}}
\label{kh79}

We work in the situation of Definition \ref{kh4def14}, and write $\dim Y=m$. First we show that the r.h.s.\ of \eq{kh4eq61} is a well defined generator of $\cP MC^{k+l}(Y;R)$. Here $\ti V=V\t_{t,Y,t'}V'$ is the fibre product in $\tManc$, which exists by Assumption \ref{kh3ass5}(c) as $t,t'$ are submersions and $Y$ is a manifold, with
\begin{equation*}
\dim\ti V\!=\!\dim V\!+\!\dim V'\!-\!\dim Y\!=\!(m\!+\!n\!-\!k)\!+\!(m\!+\!n'\!-\!l)\!-\!m\!=\!m+\!\ti n\!-\!(k\!+\!l).
\end{equation*}
The projections $\pi_V:\ti V\ra V$, $\pi_{V'}:\ti V\ra V'$ are submersions by Assumption \ref{kh3ass5}(c) as $t',t$ are, so $\ti t=t\ci\pi_V=t'\ci\pi_{V'}:\ti V\ra Y$ is a submersion by Assumption \ref{kh3ass5}(a). By Assumption \ref{kh3ass6}(l), the coorientation $c_{t'}$ on $t'$ induces a coorientation $c_{\pi_V}$ on $\pi_V:\ti V\ra V$. We then give $\ti t=t\ci\pi_V:\ti V\ra Y$ the product coorientation $c_{\ti t}=c_t\ci c_{\pi_V}$, as in Assumption~\ref{kh3ass6}(d).

Since $(s,t):V\ra\R^n\t Y$ and $(s',t'):V'\ra\R^{n'}\t Y$ are proper over open neighbourhoods $X,X'$ of $\{0\}\t Y$ in $\R^n\t Y$, $\R^{n'}\t Y$, we see that $(\ti s,\ti t):\ti V\ra\R^{\ti n}\t Y$ is proper over the open neighbourhood
\begin{align*}
X\t_YX'=\bigl\{&\bigl((x_1,\ldots,x_n,x_1',\ldots,x_{n'}'),y\bigr)\in \R^{\ti n}\t Y:\\
&\bigl((x_1,\ldots,x_n),y\bigr)\in X,\;\> 
\bigl((x_1',\ldots,x_{n'}'),y\bigr)\in X'\bigr\}
\end{align*}
of $\{0\}\t Y$ in $\R^{\ti n}\t Y$. This proves that $[\ti V,\ti n,\ti s,\ti t]$ in \eq{kh4eq61} is a generator of $\cP MC^{k+l}(Y;R)$ in the sense of Definition \ref{kh4def4}. To show that $\cup$ in \eq{kh4eq60} is well defined, we must show it takes relations Definition \ref{kh4def4}(i),(ii) in $\cP MC^k(Y;R),\cP MC^l(Y;R)$ to relations (i),(ii) in~$\cP MC^{k+l}(Y;R)$.

The sign $(-1)^{ln}$ in \eq{kh4eq61} is chosen to ensure that $\cup$ takes relation (i) in $\cP MC^k(Y;R),\cP MC^l(Y;R)$ to relation (i) in $\cP MC^{k+l}(Y;R)$. To see this, note that applying $-\cup[V',n',s',t']$ to \eq{kh4eq16} for $[V,n,s,t]$ in $\cP MC^k(Y;R)$ for some $i=0,\ldots,n$ gives
\e
\begin{split}
&(-1)^{ln}\bigl[V\t_{t,Y,t'}V',n+n',\ti s,\ti t\bigr]=\\
&(-1)^{l(n+1)}\cdot (-1)^{n-i}\bigl[(V\t \R)\t_{t\ci\pi_V,Y,t'}V',n+1+n',\check s,\check t\bigr],
\end{split}
\label{kh7eq93}
\e
where $\check s$ has the projection to $\R$ inserted in position $i$. Applying \eq{kh4eq16} for $[\ti V,\ti n,\ti s,\ti t]$ in $\cP MC^{k+l}(Y;R)$ for the same $i$ gives
\e
\bigl[V\t_{t,Y,t'}V',n+n',\ti s,\ti t\bigr]=
(-1)^{n+n'-i}\bigl[(V\t_{t,Y,t'}V')\t\R,n+n'+1,\check s,\check t\bigr].
\label{kh7eq94}
\e
Our orientation conventions imply that in cooriented manifolds over $Y$ we have
\begin{equation*}
(V\!\t\! \R)\!\t_Y\!V'\!=\!(-1)^{\dim\R(\dim Y+\dim V')}(V\!\t_Y\!V')\!\t\!\R\!=\!(-1)^{n'-l}(V\!\t_Y\!V')\!\t\!\R.
\end{equation*}
So comparing signs shows that \eq{kh7eq93} and \eq{kh7eq94} are equivalent, and $\cup$ maps relation (i) in $\cP MC^k(Y;R)$ to relation (i) in $\cP MC^{k+l}(Y;R)$.

Similarly, applying $[V,n,s,t]\cup -$ to \eq{kh4eq16} for $[V',n',s',t']$ in $\cP MC^l(Y;R)$ for some $i'=0,\ldots,n'$ gives
\e
\begin{split}
&(-1)^{ln}\bigl[V\t_YV',n+n',\ti s,\ti t\bigr]=\\
&(-1)^{ln}\cdot (-1)^{n'-i'}\bigl[V\t_{t,Y,t'\ci\pi_{V'}}(V'\t\R),n+n'+1,\hat s,\hat t\bigr],
\end{split}
\label{kh7eq95}
\e
where $\hat s$ has the projection to $\R$ inserted in position $n+i'$. Applying \eq{kh4eq16} for $[\ti V,\ti n,\ti s,\ti t]$ in $\cP MC^{k+l}(Y;R)$ for $i=n+i'$ gives
\e
\bigl[V\t_YV',n\!+\!n',\ti s,\ti t\bigr]\!=\!
(-1)^{n+n'-(n+i')}\bigl[(V\!\t_Y\!V')\!\t\!\R,n\!+\!n'\!+\!1,\hat s,\hat t\bigr].
\label{kh7eq96}
\e
Our orientation conventions imply that in cooriented manifolds over $Y$ we have $V\t_Y(V'\t\R)=(V\t_YV')\t\R$. So \eq{kh7eq95} and \eq{kh7eq96} are equivalent, and $\cup$ maps relation (i) in $\cP MC^l(Y;R)$ to relation (i) in~$\cP MC^{k+l}(Y;R)$.

To show that $\cup$ maps relation (ii) in $\cP MC^k(Y;R)$ to (ii) in $\cP MC^{k+l}(Y;R)$, suppose $\sum_{i\in I}a_i\,[V_i,n,s_i,t_i]$ satisfies Definition \ref{kh4def4}(ii)$(*)$ in $\cP MC^k(Y;R)$ for some open neighbourhood $X$ of $\{0\}\t Y$ in $\R^n\t Y$, so $\sum_{i\in I}a_i\,[V_i,n,s_i,t_i]=0$ in $\cP MC^k(Y;R)$. Apply $-\cup[V',n',s',t']$ for some $[V',n',s',t']\in\cP MC^l(Y;R)$. We will show that
\begin{equation*}
\ts(-1)^{ln}\sum_{i\in I}a_i\,[\ti V_i,\ti n,\ti s_i,\ti t_i]=0\qquad\text{in $\cP MC^{k+l}(Y;R)$,}
\end{equation*}
where $[V_i,n,s_i,t_i]\cup[V',n',s',t']=(-1)^{ln}[\ti V_i,\ti n,\ti s_i,\ti t_i]$ as in \eq{kh4eq61}. By Assumption \ref{kh3ass5}(c),(d) and \eq{kh4eq61} we have
\e
\ti V_i^\ci=V_i^\ci\t_{t\vert_{V_i^\ci},Y,t'\vert_{V^{\prime\ci}}}V^{\prime\ci}\cong\bigl\{(v,v')\in V_i^\ci\t V^{\prime\ci}:\text{$t_i(v)=t'(v')$ in $Y$}\bigr\}.
\label{kh7eq97}
\e

Since $[V',n',s',t']$ is a generator of $\cP MC^l(Y;R)$, $(s',t'):V'\ra \R^{n'}\t Y$ is a morphism in $\tManc$ proper over an open neighbourhood $X'$ of $\{0\}\t Y$ in $\R^{n'}\t Y$, with $t':V'\ra Y$ a submersion, and $\dim V'=m+n'-l$. Define
\begin{equation*}
\ti X=\bigl\{\bigl((x,x'),y\bigr):x\in\R^n,\; x'\in\R^{n'},\; y\in Y,\; (x,y)\in X,\; (x',y)\in X'\bigr\}.
\end{equation*}
Then $\ti X$ is an open neighbourhood of $\{0\}\t Y$ in $\R^{n+n'}\t Y$, as $X,X'$ are open neighbourhoods of $\{0\}\t Y$ in $\R^n\t Y$ and~$\R^{n'}\t Y$. 

Suppose $((x,x'),y)\in\ti X$ such that $\ti v\in\ti V_i^\ci$ and
\e
T_{\ti v}(\ti s_i,\ti t_i):T_{\ti v}\ti V_i^\ci\longra T_x\R^n\op T_{x'}\R^{n'}\op T_yY
\label{kh7eq98}
\e
is injective for all $i\in I$ and $\ti v\in\ti V_i$ with $(\ti s_i,\ti t_i)(\ti v)=((x,x'),y)$ in $\ti X$. By \eq{kh7eq97} we write $\ti v=(v,v')$ for $v\in V_i^\ci$ with $s_i(v)=x$ and $t_i(v)=y$ and $v'\in V^{\prime\ci}$ with $s'(v')=x'$ and $t'(v')=y$. As \eq{kh7eq97} is a transverse fibre product of manifolds,
\e
\begin{gathered}
\xymatrix@C=100pt@R=15pt{ *+[r]{T_{\ti v}\ti V_i^\ci} \ar[d]^{T_{\ti v}\pi_{V_i}} \ar[r]_{T_{\ti v}\pi_{V'}} & *+[l]{T_{v'}V^{\prime\ci}} \ar[d]_{T_{v'}t'} \\
*+[r]{T_vV_i} \ar[r]^{T_vt_i} & *+[l]{T_yY} }
\end{gathered}
\label{kh7eq99}
\e
is Cartesian in $\mathop{\rm Vect}_\R$, and Assumption \ref{kh3ass6}(l) shows that the coorientations on $T_vt_i,T_{\ti v}\pi_{V'}$ from $c_{t_i}$ and on $T_{v'}t',T_{\ti v}\pi_{V_i}$ from $c_{t'}$ are related as usual in~\eq{kh7eq99}. 

Using \eq{kh7eq99} Cartesian, we see that \eq{kh7eq98} is injective if and only if
\ea
T_v(s_i,t_i)&:T_vV_i^\ci\longra T_x\R^n\op T_yY,
\label{kh7eq100}\\
T_{v'}(s',t')&:T_{v'}V^{\prime\ci}\longra T_{x'}\R^{n'}\op T_yY,
\label{kh7eq101}
\ea
are injective. Therefore $((x,x'),y)\in\ti X$ satisfies the conditions in Definition \ref{kh4def4}(ii)$(*)$ for $\sum_{i\in I}a_i\,[\ti V_i,\ti n,\ti s_i,\ti t_i]$ in $\cP MC^{k+l}(Y;R)$ if and only if both
\begin{itemize}
\setlength{\itemsep}{0pt}
\setlength{\parsep}{0pt}
\item[(A)] $(x,y)\in X$ satisfies the conditions in Definition \ref{kh4def4}(ii)$(*)$ for $\sum_{i\in I}a_i\,[V_i,\ab n,\ab s_i,\ab t_i]$ in $\cP MC^k(Y;R)$, and 
\item[(B)] $v'\in V^{\prime\ci}$ with \eq{kh7eq101} injective for all $v'\in V'$ with $s'(v')=x'$, $t'(v')=y$.
\end{itemize}

Suppose $\ti P\subseteq T_x\R^n\op T_{x'}\R^{n'}\op T_yY$ is an $(m+\ti n-k-l)$-plane with $\pi_{T_yY}:\ti P\ra T_yY$ cooriented. Then either (i) $\ti P$ lies in a Cartesian diagram
\e
\begin{gathered}
\xymatrix@C=90pt@R=13pt{ *+[r]{\ti P} \ar[d] \ar[r] & *+[l]{P'} \ar[d] \\
*+[r]{P} \ar[r] & *+[l]{T_yY} }
\end{gathered}
\label{kh7eq102}
\e
which is a subdiagram of \eq{kh7eq101}, with $P$ an $(m+n-k)$-plane in $T_x\R^n\op T_yY$ with $\pi_{T_yY}:P\ra T_yY$ cooriented, and $P'$ an $(m+n'-l)$-plane in $T_{x'}\R^{n'}\op T_yY$ with $\pi_{T_yY}:P'\ra T_yY$ cooriented; or (ii) $\ti P$ lies in no such diagram~\eq{kh7eq102}.

In case (i), \eq{kh4eq17} for $\sum_{i\in I}a_i\,[\ti V_i,\ti n,\ti s_i,\ti t_i]$ at $((x,x'),y),\ti P$ follows by multiplying \eq{kh4eq17} for $\sum_{i\in I}a_i[V_i,n,s_i,t_i]$ at $(x,y),P$ by the number of $v'$ satisfying (B) above with $T_{v'}(s',t')[T_{v'}V^{\prime\ci}]=P'$, counted with signs according to how $T_{v'}(s',t'):T_{v'}V^{\prime\ci}\,{\buildrel\cong\over\longra}\,P'$ acts on coorientations over $T_yY$. In case (ii), \eq{kh4eq17} for $\sum_{i\in I}a_i\,[\ti V_i,\ti n,\ti s_i,\ti t_i]$ at $((x,x'),y),\ti P$ is trivial as there are no $v'\in V^{\prime\ci}$ satisfying the conditions on either side.
This gives $(*)$ in Definition \ref{kh4def4}(ii) for~$\sum_{i\in I}a_i\,[\ti V_i,\ti n,\ti s_i,\ti t_i]$. 

Hence $\sum_{i\in I}a_i\,[\ti V_i,\ti n,\ti s_i,\ti t_i]=0$ in $\cP MC^{k+l}(Y;R)$, so $-\cup[V',n',s',t']$ maps relation (ii) in $\cP MC^k(Y;R)$ to relation (ii) in $\cP MC^{k+l}(Y;R)$. Similarly, $[V,\ab n,\ab s,\ab t]\cup-$ maps relation (ii) in $\cP MC^l(Y;R)$ to relation (ii) in $\cP MC^{k+l}(Y;R)$. Therefore $\cup$ in \eq{kh4eq60}--\eq{kh4eq61} is well defined. This completes the proof.

\section{Proofs of results in \S\ref{kh5}}
\label{kh8}

\subsection{Proof of Theorem \ref{kh5thm4}}
\label{kh81}

We will follow the proof of Theorem \ref{kh4thm2} in \S\ref{kh74} closely. We must show that $MH_k^\dR(*;\R)=0$ for $k\ne 0$ and $MH_0^\dR(*;\R)\cong\R$. Lemma \ref{kh5lem1} implies that $MH_k^\dR(*;\R)=0$ for $k>0$. After introducing some notation and proving some auxiliary results in \S\ref{kh811}, we will show that $MH_k^\dR(*;\R)=0$ for $k<0$ in \S\ref{kh812}, and that $MH_0^\dR(*;\R)\cong\R$ in~\S\ref{kh813}.

\subsubsection{Some auxiliary results}
\label{kh811}

Following Definition \ref{kh7def2}, we define:

\begin{dfn} Let $*$ be the point. As in \S\ref{kh52}, $MC^\dR_k(*;\R)$ is the $\R$-vector space spanned by generators $[V,n,s,t,\om]$ with $\dim V=\deg\om+n+k$ subject to relations Definition \ref{kh5def4}(i)--(iii).  Any morphism $t:V\ra *$ is the projection $\pi:V\ra *$, so the $t$ in $[V,n,s,t,\om]=[V,n,s,\pi,\om]$ can basically be ignored. 

Write $\widetilde{MC}{}^\dR_k(*;\R)$ for the $\R$-vector space spanned by generators $[V,n,s,\pi,\om]$ with $\dim V=\deg\om+n+k$ subject to relation Definition \ref{kh5def4}(ii) only. Then there is a surjective $\R$-linear map
\begin{equation*}
\Pi:\widetilde{MC}{}^\dR_k(*;\R)\longra MC^\dR_k(*;\R),\qquad \Pi:[V,n,s,\pi,\om]\longmapsto[V,n,s,\pi,\om]
\end{equation*}
with kernel spanned by applications of Definition~\ref{kh5def4}(i),(iii).

As the relation Definition \ref{kh5def4}(ii) in $\widetilde{MC}{}^\dR_k(*;\R)$ involves only $[V_i,n,s_i,\pi,\om_i]$ with $n$ fixed, we may write 
\begin{equation*}
\widetilde{MC}{}^\dR_k(*;\R)=\ts\bigop_{n=0}^\iy \widetilde{MC}{}^\dR_k(*;\R)^n
\end{equation*}
where $\widetilde{MC}{}^\dR_k(*;\R)^n$ is spanned by generators $[V,n,s,\pi,\om]$ with $n$ fixed, modulo relation Definition \ref{kh5def4}(ii) with $n$ fixed. 

We may define $\pd:\widetilde{MC}{}^\dR_k(*;\R)\ra\widetilde{MC}{}^\dR_{k-1}(*;\R)$ as in \S\ref{kh52}, and then $\Pi\ci\pd=\pd\ci\Pi$, and $\pd\ci\pd=0$, as the proof of this in Definition \ref{kh5def5} used only relation (ii). Also $\pd$ maps~$\widetilde{MC}{}^\dR_k(*;\R)^n\ra\widetilde{MC}{}^\dR_{k-1}(*;\R)^n$.
\label{kh8def1}
\end{dfn}

For the next two propositions we will consider the following situation. Let $k\in\Z$, and $\al\in MC^\dR_k(*;\R)$ with $\pd\al=0$ in $MC^\dR_{k-1}(*;\R)$. Choose a lift $\ti\al$ of $\al$ to $\widetilde{MC}{}^\dR_k(*;\R)$, so that $\Pi(\ti\al)=\al$, and write
\e
\ti\al=\sum_{n=0}^N\sum_{i\in I^n}a_i^n\,\bigl[V_i^n,n,s_i^n,\pi,\om_i^n\bigr],
\label{kh8eq1}
\e
where $N\in\N$, $I^0,I^1,\ldots,I^N$ are finite indexing sets, $[V_i^n,n,s_i^n,\pi,\om_i^n]$ is a generator of $MC^\dR_k(*;\R)$ for all $i,n$, as in Definition \ref{kh5def4}, and $a_i^n\in\R$. Then $\Pi\ci\pd\ti\al=\pd\ci\Pi(\ti\al)=\pd\al=0$, so $\pd\ti\al$ lies in the subspace of $\widetilde{MC}{}^\dR_{k-1}(*;\R)$ spanned by relations Definition~\ref{kh5def4}(i),(iii). 

As in the definition of $\io:MC_k^\Q(*;\R)\ra MC_k(*;\R)$ in Definition \ref{kh5def2}, by replacing each $\bigl[V_i^n,n,s_i^n,\pi,\om_i^n\bigr]$ in \eq{kh8eq1} by
\begin{equation*}
\frac{1}{n!}\sum_{\si\in S_n}\sign(\si)\cdot\bigl[V_i^n,n,(s_{i,\si(1)}^n,s_{i,\si(2)}^n,\ldots,s_{i,\si(n)}^n),\pi,\om_i^n\bigr],
\end{equation*}
which changes the choice of $\ti\al$ but does not change $\al$, we may assume that the components of $\ti\al$ and $\pd\ti\al$ in $\widetilde{MC}{}^\dR_*(*;\R)^n$ are invariant under the action of $S_n$ on $\R^n$, and therefore that $\pd\ti\al$ lies in the subspace of $\widetilde{MC}{}^\dR_{k-1}(*;\R)$ spanned by relation Definition \ref{kh5def4}(i) only. Hence as for \eq{kh7eq29}, by \eq{kh5eq13} we may write
\ea
&\sum_{n=0}^N\sum_{i\in I^n}a_i^n\bigl(\bigl[\pd V_i^n,n,s_i^n\ci i_{V_i^n},\pi,\om_i^n\bigr]+(-1)^{n+k}\bigl[V_i^n,n,s_i^n,\pi,\d\om_i^n\bigr]\bigr)
\label{kh8eq2}\\
&=\sum_{n=0}^{N-1}\sum_{l=0}^n\sum_{j\in J^{n,l}}b_j^{n,l}\Bigl(\bigl[\dot V_j^{n,l},n,(\dot s_{j,1}^{n,l},\ldots,\dot s_{j,n}^{n,l}),\pi,\dot\om_j^{n,l}\bigr]-(-1)^{n-l}\bigl[\dot V_j^{n,l}\t\R,
\nonumber\\
&\qquad\qquad\qquad n+1,(\dot s_{j,1}^{n,l},\ldots,\dot s_{j,l}^{n,l},\pi_\R,\dot s_{j,l+1}^{n,l},\ldots,\dot s_{j,n}^{n,l}),\pi,\pi_{\dot V_j^{n,l}}^*(\dot\om_j^{n,l})\bigr]\Bigr)
\nonumber
\ea
in $\widetilde{MC}{}^\dR_{k-1}(*;\R)$. Define $\dot s_j^{n,l}$ and $\acute s_j^{n,l}$ as in \eq{kh7eq30}. Taking components of \eq{kh8eq2} in $\widetilde{MC}{}^\dR_{k-1}(*;\R)^n$ for $n=0,\ldots,N$ and setting $J^{N,l}=\es$ for $l=0,\ldots,N$ gives
\ea
\begin{split}
&\sum_{i\in I^n}a_i^n\bigl(\bigl[\pd V_i^n,n,s_i^n\ci i_{V_i^n},\pi,\om_i^n\bigr]+(-1)^{n+k}\bigl[V_i^n,n,s_i^n,\pi,\d\om_i^n\bigr]\bigr)\\
&=\sum_{l=0}^n\sum_{j\in J^{n,l}}b_j^{n,l}\,\bigl[\dot V_j^{n,l},n,\dot s_j^{n,l},\pi,\dot\om_j^{n,l}\bigr]\\
&-\sum_{l=0}^{n-1}\sum_{j\in J^{n-1,l}}(-1)^{n-1-l}b_j^{n-1,l}\,\bigl[\dot V_j^{n-1,l}\t\R,n,\acute s_j^{n-1,l},\pi,\pi_{\dot V_j^{n-1,l}}^*(\dot\om_j^{n-1,l})\bigr].
\end{split}
\label{kh8eq3}
\ea

Here is the analogue of Proposition \ref{kh7prop2}. It has essentially the same proof, including forms $\om_i^n,\dot\om_j^{n,l}$ throughout in a straightforward way.

\begin{prop} In the situation above, for all\/ $k\in\Z,$ taking $\al,\ti\al,N$ and the representation \eq{kh8eq1} for $\ti\al$ to be fixed, we may make alternative choices for the data $J^{n,l},b_j^{n,l},[\dot V_j^{n,l},n,\dot s_j^{n,l},\pi,\dot\om_j^{n,l}]$ on the r.h.s.\ of\/ \eq{kh8eq2}--\eq{kh8eq3} to ensure that for all\/ $n=0,\ldots,N-1$ and\/ $l=0,\ldots,n$ we have 
\e
\pd\raisebox{-4pt}{$\biggl[$}\sum_{j\in J^{n,l}}b_j^{n,l}\,\bigl[\dot V_j^{n,l},n,\dot s_j^{n,l},\pi,\dot\om_j^{n,l}\bigr]\raisebox{-4pt}{$\biggr]$}=0\qquad\text{in $\widetilde{MC}{}^\dR_{k-2}(*;\R)^n$.}
\label{kh8eq4}
\e

\label{kh8prop1}
\end{prop}

Here is the analogue of Proposition \ref{kh7prop3}. This time the proof needs nontrivial changes, since in \S\ref{kh741} the key point was that $\dim\dot V_j^{n,l}=n+k-1<n$ as $k\le 0$, but now $\dim\dot V_j^{n,l}=n+k-1+\deg \dot\om_j^{n,l}$, so we may not have~$\dim\dot V_j^{n,l}<n$.

\begin{prop} In the situation above, by keeping $\al\in MC_k^\dR(*;\R)$ with\/ $\pd\al=0$ fixed but changing the lift\/ $\ti\al$ of\/ $\al$ to $\widetilde{MC}{}^\dR_k(*;\R)$ with\/ $\Pi(\ti\al)=\al,$ we may suppose that\/ $\pd\ti\al=0$ in $\widetilde{MC}{}^\dR_{k-1}(*;\R)$.
\label{kh8prop2}
\end{prop}

\begin{proof} If $k>0$ the proposition is trivial since $MC^\dR_k(*;\R)=\widetilde{MC}{}^\dR_k(*;\R)=0$ as in Lemma \ref{kh5lem1}, so suppose $k\le 0$. We first choose an arbitrary lift $\ti\al$ of $\al$ to $\widetilde{MC}{}^\dR_k(*;\R)$ of the form \eq{kh8eq1}. Then Proposition \ref{kh8prop1} shows that we may write $\pd\ti\al$ as in \eq{kh8eq2} in terms of data $J^{n,l},b_j^{n,l},[\dot V_j^{n,l},n,\dot s_j^{n,l},\pi,\dot\om_j^{n,l}]$ satisfying \eq{kh8eq4} for all $n=0,\ldots,N-1$ and $l=0,\ldots,n$. 

For each $n=0,1,\ldots,N-1$, choose an open neighbourhood $X^n$ of 0 in $\R^n$ such that $\dot s_j^{n,l}\vert_{\supp\dot\om_j^{n,l}}:\supp\dot\om_j^{n,l}\ra\R^n$ is proper over $X^n$ for all $l=0,\ldots,n$ and $j\in J^{n,l}$. This is possible by definition of generators in Definition \ref{kh5def4}. Also, equation \eq{kh8eq4} holds as an application of relation Definition \ref{kh5def4}(ii) in $\widetilde{MC}{}^\dR_{k-2}(*;\R)^n$, which involves a condition $(*)$ including an open neighbourhood $X$ of 0 in $\R^n$. Making $X^n$ smaller if necessary, we suppose  that Definition \ref{kh5def4}(ii)$(*)$ holds for \eq{kh8eq4} with~$X=X^n$.

Next, choose $x^n\in X^n$ such that:
\begin{itemize}
\setlength{\itemsep}{0pt}
\setlength{\parsep}{0pt}
\item[(a)] $\bigl\{\la x^n:\la\in[0,1]\bigr\}\subseteq X^n$.
\item[(b)] $x^n\notin\dot s_j^{n,l}(\dot V_j^{n,l})$ for all $l=0,\ldots,n$ and $j\in J^{n,l}$ with $\dim\dot V_j^{n,l}<n$.
\item[(c)] $\dot s_j^{n,l}:\dot V_j^{n,l}\ra\R^n$ is a submersion in an open neighbourhood of $(\dot s_j^{n,l})^{-1}(x^n)$ in $\dot V_j^{n,l}$ for all $l=0,\ldots,n$ and $j\in J^{n,l}$ with $\dim\dot V_j^{n,l}\ge n$.
\end{itemize}
Here (a) holds provided $x^n$ is small enough in $\R^n$. Part (b) holds for generic $x^n$ as if $\dim\dot V_j^{n,l}<n$ then $\cH^n[\dot s_j^{n,l}(\dot V_j^{n,l})]=0$, since $\dot s_j^{n,l}[(\dot V_j^{n,l})^\ci]$ has Hausdorff dimension $<n$ by differential geometry, and $\cH^n[\dot s_j^{n,l}(\dot V_j^{n,l}\sm(\dot V_j^{n,l})^\ci)]=0$ by Assumption \ref{kh3ass7}(a). Part (c) holds for generic $x^n$ by Assumption~\ref{kh3ass7}(b).

Choose small $\de^n>0$ such that:
\begin{itemize}
\setlength{\itemsep}{0pt}
\setlength{\parsep}{0pt}
\item[(d)] $B_{\de^n}(\la x^n)\subseteq X^n$ for all $\la\in[0,1]$, with $B_{\de^n}(x^n)$ the open ball of radius $\de^n$ about $x$ in $\R^n$.
\item[(e)] $(\supp\dot\om_j^{n,l})\cap(\dot s_j^{n,l})^{-1}(B_{\de^n}(x^n))=\es$ for all $l=0,\ldots,n$ and $j\in J^{n,l}$ with~$\dim\dot V_j^{n,l}<n$.
\item[(f)] $\dot s_j^{n,l}:\dot V_j^{n,l}\ra\R^n$ is a submersion near $(\supp\dot\om_j^{n,l})\cap(\dot s_j^{n,l})^{-1}(B_{\de^n}(x^n))$ in $\dot V_j^{n,l}$ for all $l=0,\ldots,n$ and $j\in J^{n,l}$ with~$\dim\dot V_j^{n,l}\ge n$.
\end{itemize}
Here (d)--(f) are all possible if $\de^n$ is small enough, by (a)--(c), since $X^n$ is open in $\R^n$, and $\dot s_j^{n,l}\vert_{\supp\dot\om_j^{n,l}}:\supp\dot\om_j^{n,l}\ra\R^n$ is proper over~$X^n$.

For each $l=0,\ldots,n$ and $j\in J^{n,l}$, define $\check s_j^{n,l}:[0,1]\t\dot V_j^{n,l}\ra\R^n$ by
\begin{equation*}
\check s_j^{n,l}:(\la,v)\longmapsto \dot s_j^{n,l}(v)-\la x^n.
\end{equation*}
Then $\bigl[[0,1]\t\dot V_j^{n,l},n,\check s_j^{n,l},\pi,\pi_{\dot V_j^{n,l}}^*(\dot\om_j^{n,l})\bigr]$ is a generator of $MC^\dR_k(*;\R)$, where we give $[0,1]\t\dot V_j^{n,l}$ the product orientation of the standard orientation on $[0,1]$ and the given orientation on $\dot V_j^{n,l}$, as in Assumption~\ref{kh3ass6}(f),(k). 

Using equation \eq{kh5eq13}, Assumptions \ref{kh3ass3}(d), \ref{kh3ass6}(h), and $\pd[0,1]=-\{0\}\amalg\{1\}$ in oriented manifolds, we see that in $MC^\dR_{k-1}(*;\R)$ or $\widetilde{MC}{}^\dR_{k-1}(*;\R)$ we have
\ea
\pd&\bigl[[0,1]\t\dot V_j^{n,l},n,\check s_j^{n,l},\pi,\pi_{\dot V_j^{n,l}}^*(\dot\om_j^{n,l})\bigr]=-\bigl[\dot V_j^{n,l},n,\dot s_j^{n,l},\pi,\dot\om_j^{n,l}\bigr]
\nonumber\\
&+\bigl[\dot V_j^{n,l},n,\ddot s_j^{n,l},\pi,\dot\om_j^{n,l}\bigr]-\bigl[[0,1]\t\pd\dot V_j^{n,l},n,\hat s_j^{n,l},\pi,\pi_{\pd\dot V_j^{n,l}}^*\ci i_{\dot V_j^{n,l}}^*(\dot\om_j^{n,l})\bigr]
\nonumber\\
&+(-1)^{n+k}\bigl[[0,1]\t\dot V_j^{n,l},n,\check s_j^{n,l},\pi,\pi_{\dot V_j^{n,l}}^*(\d\dot\om_j^{n,l})\bigr],
\label{kh8eq5}
\ea
where $\hat s_j^{n,l}:[0,1]\t\pd\dot V_j^{n,l}\ra\R^n$ maps $\hat s_j^{n,l}:(\la,v')\mapsto \dot s_j^{n,l}\ci i_{\dot V_j^{n,l}}(v')-\la x^n$ and $\ddot s_j^{n,l}:\dot V_j^{n,l}\ra\R^n$ maps $v\mapsto\dot s_j^{n,l}(v)-x^n$.

Consider the term $\bigl[\dot V_j^{n,l},n,\ddot s_j^{n,l},\pi,\dot\om_j^{n,l}\bigr]$ in \eq{kh8eq5}. Divide into two cases (i) $\dim\dot V_j^{n,l}<n$, and (ii) $\dim\dot V_j^{n,l}\ge n$. In case (i), part (e) and the definition of $\ddot s_j^{n,l}$ imply that $\supp\dot\om_j^{n,l}\cap(\ddot s_j^{n,l})^{-1}(B_{\de^n}(0))=\es$, so that Definition \ref{kh5def4}(ii) gives
\e
\bigl[\dot V_j^{n,l},n,\ddot s_j^{n,l},\pi,\dot\om_j^{n,l}\bigr]=0\qquad \text{if $\dim\dot V_j^{n,l}<n$.}
\label{kh8eq6}
\e
In case (ii), choose an open neighbourhood $\dot W_j^{n,l}$ of $(\supp\dot\om_j^{n,l})\cap(\ddot s_j^{n,l})^{-1}(B_{\de^n}(0))$ in $(\ddot s_j^{n,l})^{-1}(B_{\de^n}(0))$ on which $\ddot s_j^{n,l}$ is a submersion, which is possible by (f), noting that $(\ddot s_j^{n,l})^{-1}(B_{\de^n}(0))=(\dot s_j^{n,l})^{-1}(B_{\de^n}(x^n))$.

Then $\ddot s_j^{n,l}\vert_{\dot W_j^{n,l}}:\dot W_j^{n,l}\ra B_{\de^n}(0)$ is a submersion, which is proper over $\supp(\dot\om_j^{n,l}\vert_{\dot W_j^{n,l}})$, and we give it the coorientation determined by the given orientation on $\dot V_j^{n,l}\supseteq \dot W_j^{n,l}$ and the standard orientation on $\R^n\supset B_{\de^n}(0)$. Thus Assumption  \ref{kh3ass8}(f) defines a pushforward form $(\ddot s_j^{n,l}\vert_{\dot W_j^{n,l}})_*(\dot\om_j^{n,l}\vert_{\dot W_j^{n,l}})$ on $B_{\de^n}(0)$, of degree $1-k$. Equation \eq{kh5eq12}, which is proved using only relation Definition \ref{kh5def4}(ii) and so also holds in $\widetilde{MC}{}^\dR_{k-1}(*;\R)$, implies that if $\dim\dot V_j^{n,l}\ge n$ then
\e
\begin{split}
\bigl[\dot V_j^{n,l},n&,\ddot s_j^{n,l},\pi,\dot\om_j^{n,l}\bigr]=\bigl[B_{\de^n}(0),n,\id,\pi,(\ddot s_j^{n,l}\vert_{\dot W_j^{n,l}})_*(\dot\om_j^{n,l}\vert_{\dot W_j^{n,l}})\bigr].
\end{split}
\label{kh8eq7}
\e

Define a $(1-k)$-form $\be^{n,l}$ on $B_{\de^n}(0)$ by
\e
\be^{n,l}=\sum_{j\in J^{n,l}:\dim\dot V_j^{n,l}\ge n}b_j^{n,l}
(\ddot s_j^{n,l}\vert_{\dot W_j^{n,l}})_*(\dot\om_j^{n,l}\vert_{\dot W_j^{n,l}}).
\label{kh8eq8}
\e
We now claim that for all $n=0,\ldots,N-1$ and $l=0,\ldots,n$ we have
\e
\begin{split}
&\pd\raisebox{-4pt}{$\biggl[$}\sum_{j\in J^{n,l}}b_j^{n,l}\,\bigl[[0,1]\t\dot V_j^{n,l},n,\check s_j^{n,l},\pi,\pi_{\dot V_j^{n,l}}^*(\dot\om_j^{n,l})\bigr]\raisebox{-4pt}{$\biggr]$}\\
&=-\sum_{j\in J^{n,l}}b_j^{n,l}\,\bigl[\dot V_j^{n,l},n,\dot s_j^{n,l},\pi,\dot\om_j^{n,l}\bigr]+
\bigl[B_{\de^n}(0),n,\id,\pi,\be^{n,l}\bigr]
\end{split}
\label{kh8eq9}
\e
in $\widetilde{MC}{}^\dR_{k-1}(*;\R)$. To see this, multiply \eq{kh8eq5} by $b_j^{n,l}$ and sum over $j\in J^{n,l}$. The l.h.s.\ of \eq{kh8eq5} yields the l.h.s.\ of \eq{kh8eq9}. The r.h.s.\ of \eq{kh8eq5} yields four sums. The first is the first sum on the r.h.s.\ of \eq{kh8eq9}. The second equals the second term on \eq{kh8eq9} by \eq{kh8eq6}--\eq{kh8eq8}. The sums coming from the third and fourth terms on the r.h.s.\ of \eq{kh8eq5} are the result of applying an operation to the l.h.s.\ of \eq{kh8eq4}, where we replace each generator $[V,n,s,\pi,\om]$ in \eq{kh8eq4} by $\bigl[[0,1]\t V,n,\hat s,\pi,\pi_{\pd V}^*(\om)\bigr]$. Hence \eq{kh8eq4} implies that the third and fourth sums add up to zero, proving \eq{kh8eq9}.

Applying $\pd$ to \eq{kh8eq9} and using $\pd\ci\pd=0$ and equations \eq{kh5eq13} and \eq{kh8eq4} gives
\e
\bigl[B_{\de^n}(0),n,\id,\pi,\d\be^{n,l}\bigr]=0\qquad\text{in $\widetilde{MC}{}^\dR_{k-2}(*;\R)$.}
\label{kh8eq10}
\e
This holds as an application of relation Definition \ref{kh5def4}(ii), so for some open neighbourhood $X$ of 0 in $B_{\de^n}(0)\subset\R^n$, Definition \ref{kh5def4}(ii)$(*)$ says that for all $(n+k-2)$-forms $\eta$ on $\R^n$ with $\supp\eta\subset X$ compact we have $\int_{B_{\de^n}(0)}\d\be^{n,l}\w\eta=0$. This implies that $\d\be^{n,l}\vert_X=0$. Making $\de^n>0$ smaller, we can suppose that $\d\be^{n,l}=0$ on $B_{\de^n}(0)$. Since $\be^{n,l}$ is a $(1-k)$-form for $k\le 0$, and $H^{1-k}_\dR\bigl(B_{\de^n}(0);\R\bigr)=0$, we see that $\be^{n,l}=\d\ga^{n,l}$ for some $(-k)$-form $\ga^{n,l}$ on $B_{\de^n}(0)$, so \eq{kh5eq13} gives
\e
\pd\bigl[B_{\de^n}(0),n,\id,\pi,\ga^{n,l}\bigr]=(-1)^{n+k}\bigl[B_{\de^n}(0),n,\id,\pi,\be^{n,l}\bigr]
\label{kh8eq11}
\e
in $\widetilde{MC}{}^\dR_{k-1}(*;\R)$.

Define a modification $\breve\al$ of $\ti\al$ in $\widetilde{MC}{}^\dR_k(*;\R)$ by
\ea
&\breve\al=\sum_{n=0}^N\sum_{i\in I^n}a_i^n\,\bigl[V_i^n,n,s_i^n,\pi,\dot\om_j^{n,l}\bigr]
-\sum_{n=0}^{N-1}\sum_{l=0}^n(-1)^{n+k}\bigl[B_{\de^n}(0),n,\id,\pi,\ga^{n,l}\bigr]
\nonumber\\
\begin{split}
&+\sum_{n=0}^{N-1}\sum_{l=0}^n(-1)^{k-l}\bigl[B_{\de^n}(0)\t\R,n+1,\id,\pi,\pi_{B_{\de^n}(0)}^*(\ga^{n,l})\bigr]\\
&+\sum_{n=0}^{N-1}\sum_{l=0}^n\sum_{j\in J^{n,l}}b_j^{n,l}\,\bigl[[0,1]\t\dot V_j^{n,l},n,\check s_j^{n,l},\pi,\pi_{\dot V_j^{n,l}}^*(\dot\om_j^{n,l})\bigr]
\end{split}
\label{kh8eq12}\\
&-\sum_{n=0}^{N-1}\sum_{l=0}^n\sum_{j\in J^{n,l}}(-1)^{n-l}b_j^{n,l}\,\bigl[[0,1]\!\t\!\dot V_j^{n,l}\!\t\!\R,n\!+\!1,\grave s_j^{n,l},\pi,\pi_{\dot V_j^{n,l}}^*(\dot\om_j^{n,l})\bigr],
\nonumber
\ea
where
\begin{align*}
\grave s_j^{n,l}=(&\hat s_{j,1}^{n,l}\ci\pi_{[0,1]\t\dot V_j^{n,l}},\ldots,\hat s_{j,l}^{n,l}\ci\pi_{[0,1]\t\dot V_j^{n,l}},\pi_\R,\\
&\hat s_{j,l+1}^{n,l}\ci\pi_{[0,1]\t\dot V_j^{n,l}},\ldots,\hat s_{j,n}^{n,l}\ci\pi_{[0,1]\t\dot V_j^{n,l}}).
\end{align*}
Then $\breve\al$ differs from $\ti\al$ by finitely many applications of relation Definition \ref{kh5def4}(i) in $MC^\dR_k(*;\R)$, so $\breve\al$ is an alternative lift of $\al$ to $\widetilde{MC}{}^\dR_k(*;\R)$, with $\Pi(\breve\al)=\al$. 

As in the end of the proof of Proposition \ref{kh7prop3}, applying $\pd$ to \eq{kh8eq12} and using equations \eq{kh8eq2}, \eq{kh8eq9}, \eq{kh8eq11} and analogues of \eq{kh8eq9}, \eq{kh8eq11} for the terms involving $[0,1]\t\dot V_j^{n,l}\t\R$ and $B_{\de^n}(0)\t\R$, we find that $\pd\breve\al=0$ in $\widetilde{MC}{}^\dR_{k-1}(*;\R)$. Replacing $\ti\al$ by $\breve\al$, Proposition \ref{kh8prop2} follows.
\end{proof}

\subsubsection{Proof that $MH_k^\dR(*;\R)=0$ for $k<0$}
\label{kh812}

We can now prove the first part of Theorem \ref{kh5thm4}. We use the notation of \S\ref{kh811}. Let $k<0$ and $\al\in MC^\dR_k(*;\R)$ with $\pd\al=0$ in $MC^\dR_{k-1}(*;\R)$. Then Proposition \ref{kh8prop2} says that we may choose a lift $\ti\al$ of $\al$ to $\widetilde{MC}{}^\dR_k(*;\R)$ with $\pd\ti\al=0$ in $\widetilde{MC}{}^\dR_{k-1}(*;\R)$. Write $\ti\al$ as in \eq{kh8eq1}. Then for each $n=0,\ldots,N$ we have
\e
\pd\sum_{i\in I^n}a_i^n\bigl[V_i^n,n,s_i^n,\pi,\om_i^n\bigr]=0
\qquad\text{in $\widetilde{MC}{}^\dR_{k-1}(*;\R)^n$.}
\label{kh8eq13}
\e

We follow the method of the proof of Proposition \ref{kh8prop2}, but applied to the $\bigl[V_i^n,n,s_i^n,\pi,\om_i^n\bigr]$ for $i\in I^n$ rather than to the $\bigl[\dot V_j^{n,l},n,\dot s_j^{n,l},\pi,\dot\om_j^{n,l}\bigr]$ for $j\in J^{n,l}$. For each $n=0,1,\ldots,N$, we choose an open neighbourhood $X^n$ of 0 in $\R^n$ such that $s_i^n\vert_{\supp\om_i^n}:\supp\om_i^n\ra\R^n$ is proper over $X^n$ for all $i\in I^n$. Equation \eq{kh8eq13} holds as an application of relation Definition \ref{kh5def4}(ii) in $\widetilde{MC}{}^\dR_{k-1}(*;\R)^n$. Making $X^n$ smaller if necessary, we suppose that Definition \ref{kh5def4}(ii)$(*)$ holds for \eq{kh8eq13} with~$X=X^n$.

As in the proof of Proposition \ref{kh8prop2}, we choose $x^n\in X^n$ and $\de^n>0$ satisfying
\begin{itemize}
\setlength{\itemsep}{0pt}
\setlength{\parsep}{0pt}
\item[(a)] $\bigl\{\la x^n:\la\in[0,1]\bigr\}\subseteq X^n$;
\item[(b)] $x^n\notin s_i^n(V_i^n)$ for all $i\in I^n$ with $\dim V_i^n<n$;
\item[(c)] $s_i^n:V_i^n\ra\R^n$ is a submersion in an open neighbourhood of $(s_i^n)^{-1}(x^n)$ in $V_i^n$ for all $i\in I^n$ with~$\dim V_i^n\ge n$;
\item[(d)] $B_{\de^n}(\la x^n)\subseteq X^n$ for all $\la\in[0,1]$; 
\item[(e)] $(\supp\om_i^n)\cap(s_i^n)^{-1}(B_{\de^n}(x^n))=\es$ for all $i\in I^n$ with $\dim V_i^n<n$; and
\item[(f)] $s_i^n:V_i^n\ra\R^n$ is a submersion near $(\supp\om_i^n)\cap(s_i^n)^{-1}(B_{\de^n}(x^n))$ in $V_i^n$ for all $i\in I^n$ with~$\dim V_i^n\ge n$.
\end{itemize}

For each $i\in I^n$, define $\check s_i^n:[0,1]\t V_i^n\ra\R^n$ by $\check s_i^n:(\la,v)\mapsto s_i^n(v)-\la x^n$. Then $\bigl[[0,1]\t V_i^n,n,\check s_i^n,\pi,\pi_{V_i^n}^*(\om_i^n)\bigr]$ is a generator of $MC^\dR_{k+1}(*;\R)$, and as in \eq{kh8eq5}, in $MC^\dR_k(*;\R)$ or $\widetilde{MC}{}^\dR_k(*;\R)$ we have
\begin{align*}
\pd\bigl[[0,1]\t V_i^n&,n,\check s_i^n,\pi,\pi_{V_i^n}^*(\om_i^n)\bigr]=-\bigl[V_i^n,n,s_i^n,\pi,\om_i^n\bigr]
\nonumber\\
&+\bigl[V_i^n,n,\ddot s_i^n,\pi,\om_i^n\bigr]-\bigl[[0,1]\t\pd V_i^n,n,\hat s_i^n,\pi,\pi_{\pd V_i^n}^*\ci i_{V_i^n}^*(\om_i^n)\bigr]
\nonumber\\
&+(-1)^{n+k+1}\bigl[[0,1]\t V_i^n,n,\check s_i^n,\pi,\pi_{V_i^n}^*(\d\om_i^n)\bigr],
\end{align*}
where $\hat s_i^n:[0,1]\t\pd V_i^n\ra\R^n$ maps $\hat s_i^n:(\la,v')\mapsto s_i^n\ci i_{V_i^n}(v')-\la x^n$ and $\ddot s_i^n:V_i^n\ra\R^n$ maps $v\mapsto s_i^n(v)-x^n$.

As for \eq{kh8eq6}, using (e) we see that $\bigl[V_i^n,n,\ddot s_i^n,\pi,\om_i^n\bigr]=0$ if $\dim V_i^n<n$. If $\dim V_i^n\ge n$, choose an open neighbourhood $W_i^n$ of $(\supp\om_i^n)\cap(\ddot s_i^n)^{-1}(B_{\de^n}(0))$ in $(\ddot s_i^n)^{-1}(B_{\de^n}(0))$ on which $\ddot s_i^n$ is a submersion, which is possible by (f). Then $\ddot s_i^n\vert_{W_i^n}:W_i^n\ra B_{\de^n}(0)$ is a submersion, which is proper over $\supp(\om_i^n\vert_{W_i^n})$, and we give it the coorientation determined by the given orientation on $V_i^n\supseteq W_i^n$ and the standard orientation on $\R^n\supset B_{\de^n}(0)$. Thus Assumption \ref{kh3ass8}(f) defines a pushforward form $(\ddot s_i^n\vert_{W_i^n})_*(\om_i^n\vert_{W_i^n})$ on $B_{\de^n}(0)$, of degree $-k$. As for \eq{kh8eq7}, equation \eq{kh5eq12} implies that
\begin{equation*}
\bigl[V_i^n,n,\ddot s_i^n,\pi,\om_i^n\bigr]=\bigl[B_{\de^n}(0),n,\id,\pi,(\ddot s_i^n\vert_{ W_i^n})_*(\om_i^n\vert_{ W_i^n})\bigr]\quad \text{if $\dim V_i^n\ge n$.}
\end{equation*}

As for \eq{kh8eq8}, define a $(-k)$-form $\be^n$ on $B_{\de^n}(0)$ by
\begin{equation*}
\be^n=\sum_{i\in I^n:\dim V_i^n\ge n}a_i^n(\ddot s_i^n\vert_{ W_i^n})_*(\om_i^n\vert_{ W_i^n}).
\end{equation*}
Then as for \eq{kh8eq9}, for all $n=0,\ldots,N$ we find that
\e
\begin{split}
&\pd\raisebox{-4pt}{$\biggl[$}\sum_{i\in I^n}a_i^n\,\bigl[[0,1]\t V_i^n,n,\check s_i^n,\pi,\pi_{ V_i^n}^*(\om_i^n)\bigr]\raisebox{-4pt}{$\biggr]$}\\
&=-\sum_{i\in I^n}a_i^n\,\bigl[ V_i^n,n, s_i^n,\pi,\om_i^n\bigr]+
\bigl[B_{\de^n}(0),n,\id,\pi,\be^n\bigr].
\end{split}
\label{kh8eq14}
\e
As for \eq{kh8eq10}, applying $\pd$ to \eq{kh8eq14} and using $\pd\ci\pd=0$, \eq{kh5eq13} and \eq{kh8eq13} gives
\e
\bigl[B_{\de^n}(0),n,\id,\pi,\d\be^n\bigr]=0\qquad\text{in $\widetilde{MC}{}^\dR_{k-1}(*;\R)$.}
\label{kh8eq15}
\e

As for the proof of \eq{kh8eq11}, making $\de^n>0$ smaller if necessary we can suppose that $\d\be^n=0$ on $B_{\de^n}(0)$, so as $\be^n$ is a $(-k)$-form on $B_{\de^n}(0)$ for $k<0$ we may write $\be^n=\d\ga^n$ for some $(-1-k)$-form $\ga^n$ on $B_{\de^n}(0)$, and \eq{kh5eq13} gives
\e
\pd\bigl[B_{\de^n}(0),n,\id,\pi,\ga^n\bigr]=(-1)^{n+k+1}\bigl[B_{\de^n}(0),n,\id,\pi,\be^n\bigr]
\label{kh8eq16}
\e
in $\widetilde{MC}{}^\dR_k(*;\R)$.

Now define $\ti\be\in\widetilde{MC}{}^\dR_{k+1}(*;\R)$ by
\begin{align*}
\ti\be=\,&-\sum_{n=0}^N\sum_{i\in I^n}a_i^n\,\bigl[[0,1]\t V_i^n,n,\check s_i^n,\pi,\pi_{V_i^n}^*(\om_i^n)\bigr]\\
&+\sum_{n=0}^N(-1)^{n+k+1}\bigl[B_{\de^n}(0),n,\id,\pi,\ga^n\bigr].
\end{align*}
Then combining equations \eq{kh8eq1}, \eq{kh8eq14} and \eq{kh8eq16} gives $\pd\ti\be=\ti\al$ in $\widetilde{MC}{}^\dR_k(*;\R)$. Setting $\be=\Pi(\ti\be)$ in $MC^\dR_{k+1}(*;\R)$, we have $\pd\be=\al$. Thus, for all $\al$ in $MC^\dR_k(*;\R)$ with $\pd\al=0$ for $k<0$, we can find $\be\in MC^\dR_{k+1}(*;\R)$ with $\pd\be=\al$. Hence $MH_k^\dR(*;\R)=0$ for $k<0$, as we have to prove.

\subsubsection{Proof that $MH_0^\dR(*;\R)\cong\R$}
\label{kh813}

Since $MC^\dR_1(*;\R)=0$ by Lemma \ref{kh5lem1}, we have
\e
MH_0^\dR(*;\R)=\Ker\bigl(\pd:MC^\dR_0(*;\R)\longra MC^\dR_{-1}(*;\R)\bigr).
\label{kh8eq17}
\e
As $*$ has the standard orientation, Definition \ref{kh5def5} defines the fundamental cycle $[*]=[*,0,0,\id_*,1_*]\in MC^\dR_0(*;\R)$, with $\pd[*]=0$. Define an $\R$-linear map
\e
\io:\R\longra \Ker\bigl(\pd:MC^\dR_0(*;\R)\longra MC^\dR_{-1}(*;\R)\bigr),\qquad \io:a\longmapsto a\,[*].
\label{kh8eq18}
\e
We will show that $\io$ is both surjective and injective, so that $MH_0^\dR(*;\R)\cong\R$, by an isomorphism identifying $[[*]]\in MH_0^\dR(*;\R)$ with $1\in\R$, as in Theorem~\ref{kh5thm4}.

Suppose $\al\in MC^\dR_0(*;\R)$ with $\pd\al=0$ in $MC^\dR_{-1}(*;\R)$. Using the notation $\widetilde{MC}{}^\dR_k(*;\R),\widetilde{MC}{}^\dR_k(*;\R)^n$ of \S\ref{kh811}, Proposition \ref{kh8prop2} gives a lift $\ti\al$ of $\al$ to $\widetilde{MC}{}^\dR_0(*;\R)$ with $\Pi(\ti\al)=\al$ and $\pd\ti\al=0$ in $\widetilde{MC}{}^\dR_{-1}(*;\R)$. Write $\ti\al$ as in \eq{kh8eq1}. 

In \S\ref{kh812}, the assumption that $k<0$ is first used between \eq{kh8eq15} and \eq{kh8eq16}, where $\be^n$ is a $(-k)$-form on $B_{\de^n}(0)$ with $\d\be^n=0$, so as $k<0$ we may write $\be^n=\d\ga^n$. Thus, when $k=0$ we may use the argument of \S\ref{kh812} as far as \eq{kh8eq15}. This gives constants $\de^n>0$ and functions (0-forms) $\be^n:B_{\de^n}(0)\ra\R$ for $n=0,\ldots,N$ with $\d\be^n=0$, so that $\be^n$ is constant, say $\be^n=b^n\cdot 1_{B_{\de^n}(0)}$ for $b^n\in\R$. Then in $\widetilde{MC}{}^\dR_0(*;\R)$ we have
\e
\begin{split}
\ti\al&=\ti\al+\pd\raisebox{-4pt}{$\biggl\{$}\sum_{n=0}^N\sum_{i\in I^n}a_i^n\,\bigl[[0,1]\t V_i^n,n,\check s_i^n,\pi,\pi_{ V_i^n}^*(\om_i^n)\bigr]\raisebox{-4pt}{$\biggr\}$}\\
&=\sum_{n=0}^Nb^n\,\bigl[B_{\de^n}(0),n,\id,\pi,1_{B_{\de^n}(0)}\bigr]=\sum_{n=0}^Nb^n\,\bigl[\R^n,n,\id,\pi,1_{\R^n}\bigr].
\end{split}
\label{kh8eq19}
\e
Here in the first step the term $\{\cdots\}$ in \eq{kh8eq19} is zero, since $MC^\dR_1(*;\R)=0$. In the second step we sum \eq{kh8eq14} from $n=0,\ldots,N$ and use \eq{kh8eq1} and $\be^n=b^n\cdot 1_{B_{\de^n}(0)}$, and in the third step we use Definition~\ref{kh5def4}(ii).

Now project \eq{kh8eq19} to $MC^\dR_0(*;\R)$ using $\Pi$. In $MC^\dR_0(*;\R)$ we have
\begin{equation*}
\bigl[\R^n,n,\id,\pi,1_{\R^n}\bigr]=[*,0,0,\id_*,1_*]=[*],
\end{equation*}
using  Definition \ref{kh5def4}(i) and induction on $n$. Hence \eq{kh8eq19} implies that $\al=\sum_{n=0}^Nb^n[*]=\io\bigl(\sum_{n=0}^Nb^n\bigr)$ in $MC^\dR_0(*;\R)$, so $\io$ in \eq{kh8eq18} is surjective.

To show that $\io$ is injective, suppose $r\in\R$ with $\io(r)=0$. We will prove $r=0$. Regarding $r[*]$ as an element of $\widetilde{MC}{}^\dR_0(*;\R)$ from \S\ref{kh811}, since $\Pi(r[*])=\io(r)=0$ in $MC^\dR_0(*;\R)$, $r[*]$ lies in the kernel of $\Pi:\widetilde{MC}{}^\dR_0(*;\R)\ra MC^\dR_0(*;\R)$, which is spanned by relations Definition \ref{kh5def4}(i),(iii) in $MC^\dR_0(*;\R)$. As in \S\ref{kh811}, by averaging over the action of the symmetric group $S_n$ permuting the coordinates of $\R^n$ in $\widetilde{MC}{}^\dR_0(*;\R)^n$, we can eliminate the applications of Definition \ref{kh5def4}(iii). Thus as in \eq{kh8eq2}, in $\widetilde{MC}{}^\dR_0(*;\R)$ we may write
\ea
r[*]&=\sum_{n=0}^N\sum_{l=0}^n\sum_{i\in I^{n,l}}a_i^{n,l}\Bigl(\bigl[V_i^{n,l},n,s_i^{n,l},\pi,\om_i^{n,l}\bigr]-(-1)^{n-l}\bigl[V_i^{n,l}\t\R,
\nonumber\\
&\qquad n+1,(s_{i,1}^{n,l},\ldots,s_{i,l}^{n,l},\pi_\R,s_{i,l+1}^{n,l},\ldots,s_{i,n}^{n,l}),\pi,\pi_{V_i^{n,l}}^*(\om_i^{n,l})\bigr]\Bigr).
\label{kh8eq20}
\ea

Taking components of \eq{kh8eq20} in $\widetilde{MC}{}^\dR_0(*;\R)^n$ for $n=0,\ldots,N+1$ gives
\ea
r[*]&=\sum_{i\in I^{0,0}}a_i^{0,0}\bigl[V_i^{0,0},0,\pi_{\R^0},\pi,\om_i^{0,0}\bigr],
\label{kh8eq21}\\
0&=\sum_{l=0}^n\sum_{i\in I^{n,l}}a_i^{n,l}\bigl[V_i^{n,l},n,s_i^{n,l},\pi,\om_i^{n,l}\bigr]
\label{kh8eq22}\\
&-\sum_{l=0}^{n-1}\sum_{i\in I^{n-1,l}}
\begin{aligned}[t]&(-1)^{n-1-l}a_i^{n-1,l}\bigl[V_i^{n-1,l}\t\R,
n,(s_{i,1}^{n-1,l},\ldots,s_{i,l}^{n-1,l},\pi_\R,\\
&s_{i,l+1}^{n-1,l},\ldots,s_{i,n-1}^{n-1,l}),\pi,\pi_{V_i^{n-1,l}}^*(\om_i^{n-1,l})\bigr],\quad n\!=\!1,\ldots,N,\!\!\!\!
\end{aligned}
\nonumber\\
0&=-\sum_{l=0}^{N}\sum_{i\in I^{N,l}}
\begin{aligned}[t](-1)^{N-l}a_i^{N,l}\bigl[V_i^{N,l}\t\R,
N+1,(s_{i,1}^{N,l},\ldots,s_{i,l}^{N,l},\\
\pi_\R,s_{i,l+1}^{N,l},\ldots,s_{i,N}^{N,l}),\pi,\pi_{V_i^{N,l}}^*(\om_i^{N,l})\bigr].
\end{aligned}
\label{kh8eq23}
\ea
By Proposition \ref{kh8prop1}, we may choose the representation \eq{kh8eq20} such that for all $n=0,\ldots,N$ and $l=0,\ldots,n$ we have 
\e
\pd\sum_{i\in I^{n,l}}a_i^{n,l}\,\bigl[V_i^{n,l},n,s_i^{n,l},\pi,\om_i^{n,l}\bigr]=0\qquad\text{in $\widetilde{MC}{}^\dR_{-1}(*;\R)^n$.}
\label{kh8eq24}
\e

By the argument used to prove \eq{kh8eq19}, starting with \eq{kh8eq24} rather than $\pd\ti\al=0$, we can show that there exist $b^{n,l}\in\R$ for $n=0,\ldots,N$ and $l=0,\ldots,n$ with
\e
\sum_{i\in I^{n,l}}a_i^{n,l}\,\bigl[V_i^{n,l},n,s_i^{n,l},\pi,\om_i^{n,l}\bigr]=b^{n,l}\,\bigl[\R^n,n,\id,\pi,1_{\R^n}\bigr]
\label{kh8eq25}
\e
in $\widetilde{MC}{}^\dR_0(*;\R)^n$. Inserting an $\R$ in the $(l+1)^{\rm th}$ coordinate in $\R^n$ throughout \eq{kh8eq25}, as in Definition \ref{kh5def4}(i), yields
\ea
&\sum_{i\in I^{n,l}\!\!}\!
(-\!1)^{n-l}a_i^{n,l}\bigl[V_i^{n,l}\!\t\!\R,
n\!+\!1,(s_{i,1}^{n,l},\ldots,s_{i,l}^{n,l},\pi_\R,s_{i,l+1}^{n,l},\ldots,s_{i,n}^{n,l}),\pi,\pi_{V_i^{n,l}}^*(\om_i^{n,l})\bigr]
\nonumber\\
&=b^{n,l}\,\bigl[\R^{n+1},n+1,\id,\pi,1_{\R^{n+1}}\bigr].
\label{kh8eq26}
\ea
Rewriting \eq{kh8eq21}--\eq{kh8eq23} using \eq{kh8eq25}--\eq{kh8eq26} now gives
\ea
r[*]&=b^{0,0}[*],
\label{kh8eq27}\\ 
0&=\bigl(\ts\sum_{l=0}^nb^{n,l}-\sum_{l=0}^{n-1}b^{n-1,l}\bigr)
\bigl[\R^n,n,\id,\pi,1_{\R^n}\bigr],\quad n\!=\!1,\ldots,N,
\label{kh8eq28}\\
0&=-\ts\sum_{l=0}^Nb^{N,l}\bigl[\R^{N+1},N+1,\id,\pi,1_{\R^{N+1}}\bigr].
\label{kh8eq29}
\ea

Equations \eq{kh8eq27}--\eq{kh8eq29} hold as applications of relation Definition \ref{kh5def4}(ii) in $\widetilde{MC}{}^\dR_0(*;\R)^n$. Thus for \eq{kh8eq28}, for example, there exists an open neighbourhood $X$ of 0 in $\R^n$ such that Definition \ref{kh5def4}(ii)$(*)$ holds, so that for all $n$-forms $\eta$ on $\R^n$ with $\supp\eta\subset X$ compact we have 
\begin{equation*}
\bigl(\ts\sum_{l=0}^nb^{n,l}-\sum_{l=0}^{n-1}b^{n-1,l}\bigr)\cdot \int_{\R^n}1_{\R^n}\w\eta=0.
\end{equation*}
Choosing such $\eta$ with $\int_{\R^n}\eta\ne 0$, we see that \eq{kh8eq27}--\eq{kh8eq29} are equivalent to
\ea
r&=b^{0,0},
\label{kh8eq30}\\ 
0&=\ts\sum_{l=0}^nb^{n,l}-\sum_{l=0}^{n-1}b^{n-1,l},\quad n\!=\!1,\ldots,N,
\label{kh8eq31}\\
0&=-\ts\sum_{l=0}^Nb^{N,l}.
\label{kh8eq32}
\ea
Taking the sum of equations \eq{kh8eq30}--\eq{kh8eq32} over all $n$ yields $r=0$. Therefore $\io$ is injective. We have now shown that $\io$ in \eq{kh8eq18} is surjective and injective, and thus an isomorphism, so by \eq{kh8eq17}--\eq{kh8eq18} we have $MH_0^\dR(*;\R)\cong\R$, completing the proof of Theorem~\ref{kh5thm4}.

\subsection{Proof of Proposition \ref{kh5prop3}}
\label{kh82}

Let $Y$ be a manifold and $k\in\Z$. To show that the maps $F_\Mh^\dRMh,F_\QMh^\dRMh$ in equation \eq{kh5eq20} of Example \ref{kh5ex1} are well defined, note that they clearly take relation (i) in $MC_k(Y;\R)$ to relation (i) in $MC_k^\dR(Y;\R)$, and relations (i),(iii) in $MC_k^\Q(Y;\R)$ to relations (i),(iii) in $MC_k^\dR(Y;\R)$. So we only need to show they take relation (ii) in $MC_k(Y;\R),MC_k^\Q(Y;\R)$ to relation (ii) in~$MC_k^\dR(Y;\R)$. 

Let $\sum_{i\in I}a_i\,[V_i,n,s_i,t_i]=0$ in $MC_k(Y;\R)$ or $MC_k^\Q(Y;\R)$ by relation (ii), so that Definition \ref{kh4def1}(ii)$(*)$ holds for some open $0\in X\subseteq\R^n$. Suppose for a contradiction that Definition \ref{kh5def4}(ii)$(*)$ with $\om_i=1_{V_i}$ for all $i\in I$ does not hold for the same $X$. Then there exists $\al\in C^\iy\bigl(\La^{n+k}T^*(\R^n\t Y)\bigr)$ with $\supp\al\subseteq K\t Y$ for some compact subset $K\subseteq X$ and
\e
\sum_{i\in I}a_i\int_{V_i}(s_i,t_i)^*(\al)\ne 0.
\label{kh8eq33}
\e

For each $(x,y)\in\supp\al$, choose an open neighbourhood $U_{x,y}$ of $(x,y)$ in $X\t Y$ diffeomorphic to $(-1,1)^m$ for $m=n+\dim Y$, with coordinates $(z_1,\ldots,z_m)$. Then $\bigl\{U_{x,y}:(x,y)\in\supp\al\bigr\}$ is an open cover of $\bigcup_{x,y}U_{x,y}\supseteq\supp\al$, so we can choose a subordinate locally finite smooth partition of unity $\bigl\{\eta_{x,y}:(x,y)\in\supp\al\bigr\}$. We have
\begin{equation*}
\sum_{i\in I}a_i\int_{V_i}(s_i,t_i)^*(\al)=\sum_{(x,y)\in\supp\al}\sum_{i\in I}a_i\int_{V_i}(s_i,t_i)^*(\eta_{x,y}\al),
\end{equation*}
where compact support of $(s_i,t_i)^*(\al)$ and local finiteness of $\{\eta_{x,y}\}$ imply that there are only finitely many nonzero terms in the sum over $(x,y)$. By \eq{kh8eq33}, for some $(x,y)\in\supp\al$ this term must be nonzero, so we can replace $\al$ by $\eta_{x,y}\al$, so that $(s_i,t_i)(V_i)\cap\supp\al\subseteq U_{x,y}\subseteq\R^n\t Y$ for all $i\in I$.

Using the coordinates $(z_1,\ldots,z_m)$ on $U_{x,y}$, write 
\begin{equation*}
\al\vert_{U_{x,y}}=\sum_{1\le j_1<\cdots<j_{n+k}\le m}\al_{j_1\cdots j_{n+k}}\d z_{j_1}\w\cdots\w\d z_{j_{n+k}}.
\end{equation*}
Splitting \eq{kh8eq33} into a sum over such $j_1,\ldots,j_{n+k}$, we see that we can choose $1\le j_1<\cdots<j_{n+k}\le m$ such that 
\e
\sum_{i\in I}a_i\int_{V_i}(s_i,t_i)^*(\al_{j_1\cdots j_{n+k}}\d z_{j_1}\w\cdots\w\d z_{j_{n+k}})\ne 0.
\label{kh8eq34}
\e
By reordering $z_1,\ldots,z_m$, we suppose for simplicity that $j_i=i$ for $i=1,\ldots,m$.

For each $i\in I$, define $\ti V_i^\ci\subseteq V_i^\ci$ to be the open subset of points $v$ in $V_i^\ci\cap (s_i,t_i)^{-1}(U_{x,y})$ such that $(s_i,t_i)^*(\d z_1)\vert_v,\ldots,(s_i,t_i)^*(\d z_{n+k})\vert_v$ form a basis of $T_v^*V_i^\ci$, and define $\ep_i:\ti V_i^\ci\ra\{1,-1\}$ by $\ep_i(v)=1$ if $(s_i,t_i)^*(\d z_1)\vert_v,\ab\ldots,\ab(s_i,t_i)^*(\d z_{n+k})\vert_v$ are an oriented basis, w.r.t.\ the orientation on $V_i$, and $\ep_i(v)=-1$ otherwise. Then we may write
\e
\begin{split}
0&\ne \sum_{i\in I}a_i\int_{V_i}(s_i,t_i)^*(\al_{1\cdots n+k}\,\d z_1\w\cdots\w\d z_{n+k})\\
&=\sum_{i\in I}a_i\int_{V_i^\ci}(s_i,t_i)^*(\al_{1\cdots n+k}\,\d z_1\w\cdots\w\d z_{n+k})\\
&=\sum_{i\in I}a_i\int_{\ti V_i^\ci}(s_i,t_i)^*(\al_{1\cdots n+k}\,\d z_1\w\cdots\w\d z_{n+k})\\
&=\int_{\begin{subarray}{l} (z_1,\ldots,z_{n+k})\\ \in(-1,1)^{n+k} \end{subarray}}
\sum_{\begin{subarray}{l} (z_{n+k+1},\ldots,z_m) \\ \in (-1,1)^{m-n-k}\end{subarray}}\al_{1\cdots n+k}(z_1,\ldots,z_m)\,\d z_1\w\cdots\w\d z_{n+k}\cdot{}\\
&\qquad
\biggl[\sum_{i\in I}\sum_{v\in\ti V_i^\ci:(s_i,t_i)(v)=(z_1,\ldots,z_m)}a_i\cdot \ep_i(v)\biggr]=0.
\end{split}
\label{kh8eq35}
\e
Here the first step of \eq{kh8eq35} is \eq{kh8eq34}, and the second is the definition Assumption \ref{kh3ass8}(e) of $\int_{V_i}\cdots$. For the third, if $v\in V_i^\ci\sm\ti V_i^\ci$ then either $v\notin (s_i,t_i)^{-1}(U_{x,y})$, so that $(s_i,t_i)^*(\al_{1\cdots n+k}\d z_1\w\cdots\w\d z_{n+k})\vert_v=0$ as $(s_i,t_i)(V_i)\cap\supp\al\subseteq U_{x,y}$, or $(s_i,t_i)^*(\d z_j)\vert_v$ for $j=1,\ldots,n+k$ are linearly dependent, so $(s_i,t_i)^*(\d z_1\w\cdots\w\d z_{n+k})\vert_v=0$. Thus in both cases the integrand in \eq{kh8eq35} is zero at $v\in V_i^\ci\sm\ti V_i^\ci$, and omitting such points $v$ does not change the integral.

Define $\Pi_i:\ti V_i^\ci\ra(-1,1)^{n+k}$ by composing $(s_i,t_i)\vert_{\ti V_i^\ci}:\ti V_i^\ci\ra U_{x,y}\cong (-1,1)^m$ with the projection $(-1,1)^m\ra(-1,1)^{n+k}$ to the first $n+k$ coordinates, $(z_1,\ldots,z_m)\mapsto(z_1,\ldots,z_{n+k})$. Then $\Pi_i:\ti V_i^\ci\ra(-1,1)^{n+k}$ is \'etale by definition of $\ti V_i^\ci$, and $\ep_i:\ti V_i^\ci\ra\{1,-1\}$ is 1 if $\Pi_i$ is orientation-preserving, and $-1$ otherwise. For the fourth step of \eq{kh8eq35}, we push forward the integral $\int_{\ti V_i^\ci}\cdots$ along $\Pi_i$ to an integral $\int_{(-1,1)^m}\cdots$. Since the integral $\int_{\ti V_i^\ci}\cdots$ is $L^1$ by Assumption \ref{kh3ass8}(e), the sums in the fifth step of \eq{kh8eq35} need not have only finitely many nonzero terms, but are absolutely convergent for almost all $(z_1,\ldots,z_{n+k})\in(-1,1)^{n+k}$, and the limit is an integrable function on $(-1,1)^{n+k}$ with integral equal to the fourth step of~\eq{kh8eq35}.

For the fifth step of \eq{kh8eq35}, note that the terms $[\cdots]$ on the last line, for fixed $(z_1,\ldots,z_m)\in U_{x,y}\subseteq X\t Y\subseteq\R^n\t Y$, may be interpreted as the sum over $(n+k)$-planes $P$ in $T_{(z_1,\ldots,z_m)}(\R^n\t Y)$ with isomorphic projection to $T_{(z_1,\ldots,z_{n+k})}(-1,1)^{n+k}$, of equation \eq{kh4eq2} in Definition \ref{kh4def1}(ii)$(*)$ at $(z_1,\ldots,z_m),\ab P$. So the terms $[\cdots]$ are zero, proving the fifth step, and completing~\eq{kh8eq35}. 

Equation \eq{kh8eq35} gives us the contradiction we need. Thus Definition \ref{kh5def4}(ii)$(*)$ with $\om_i=1_{V_i}$ does hold, so $\sum_{i\in I}a_i\,[V_i,n,s_i,t_i,1_{V_i}]=0$ in $MC_k^\dR(Y;\R)$ by Definition \ref{kh5def4}(ii). Hence the maps $F_\Mh^\dRMh,F_\QMh^\dRMh$ in \eq{kh5eq20} take relation (ii) in $MC_k(Y;\R),\ab MC_k^\Q(Y;\R)$ to relation (ii) in $MC_k^\dR(Y;\R)$, and are well defined.

\subsection{Proof of Proposition \ref{kh5prop5}}
\label{kh83}

The proof is similar in parts to that of Proposition \ref{kh4prop11} in \S\ref{kh79}. We work in the situation of Definition \ref{kh5def10}, and write $\dim Y=m$. First we show that the r.h.s.\ of \eq{kh5eq39} is a well defined generator of $\cP MC_\dR^{k+l}(Y;\R)$. Here $\ti V=V\t_{t,Y,t'}V'$ is the fibre product in $\tManc$, which exists by Assumption \ref{kh3ass5}(c) as $t,t'$ are submersions and $Y$ is a manifold, with
\begin{align*}
\dim\ti V=\dim V&+\dim V'-\dim Y=(\deg\om+m+n-k)\\
&+(\deg\om'+m+n'-l)-m=\deg\ti\om+m+\ti n-(k+l).
\end{align*}
The projections $\pi_V:\ti V\ra V$, $\pi_{V'}:\ti V\ra V'$ are submersions by Assumption \ref{kh3ass5}(c) as $t',t$ are, so $\ti t=t\ci\pi_V=t'\ci\pi_{V'}:\ti V\ra Y$ is a submersion by Assumption \ref{kh3ass5}(a). By Assumption \ref{kh3ass6}(l), the coorientation $c_{t'}$ on $t'$ induces a coorientation $c_{\pi_V}$ on $\pi_V:\ti V\ra V$. We then give $\ti t=t\ci\pi_V:\ti V\ra Y$ the product coorientation $c_{\ti t}=c_t\ci c_{\pi_V}$, as in Assumption~\ref{kh3ass6}(d).

Since $(s,t)\vert_{\supp\om}:\supp\om\ra\R^n\t Y$ and $(s',t')\vert_{\supp\om'}:\supp\om'\ra\R^{n'}\t Y$ are proper over open neighbourhoods $X,X'$ of $\{0\}\t Y$ in $\R^n\t Y$, $\R^{n'}\t Y$, we see that $(\ti s,\ti t)\vert_{\supp\ti\om}:\supp\ti\om\ra\R^{\ti n}\t Y$ is proper over the open neighbourhood
\begin{align*}
X\t_YX'=\bigl\{&\bigl((x_1,\ldots,x_n,x_1',\ldots,x_{n'}'),y\bigr)\in \R^{\ti n}\t Y:\\
&\bigl((x_1,\ldots,x_n),y\bigr)\in X,\;\> 
\bigl((x_1',\ldots,x_{n'}'),y\bigr)\in X'\bigr\}
\end{align*}
of $\{0\}\t Y$ in $\R^{\ti n}\t Y$. This proves that $[\ti V,\ti n,\ti s,\ti t,\ti\om]$ in \eq{kh5eq39} is a generator of $\cP MC_\dR^{k+l}(Y;\R)$ in the sense of Definition \ref{kh5def7}. To show that $\cup$ in \eq{kh5eq38}--\eq{kh5eq39} is well defined, we must show it takes relations Definition \ref{kh5def7}(i)--(iii) in $\cP MC_\dR^k(Y;\R),\cP MC_\dR^l(Y;\R)$ to relations (i)--(iii) in~$\cP MC_\dR^{k+l}(Y;\R)$.

We show that $\cup$ takes relation (i) in $\cP MC_\dR^k(Y;\R),\cP MC_\dR^l(Y;\R)$ to (i) in $\cP MC_\dR^{k+l}(Y;\R)$ as for the M-cohomology case in \S\ref{kh79}. It is also easy to see that $\cup$ in \eq{kh5eq38}--\eq{kh5eq39} takes relation (iii) in $\cP MC_\dR^k(Y;\R),\cP MC_\dR^l(Y;\R)$ to (iii) in $\cP MC_\dR^{k+l}(Y;\R)$. So we need to prove that $\cup$ maps relation (ii) to~(ii).

To show that $\cup$ maps (ii) in $\cP MC_\dR^k(Y;\R)$ to (ii) in $\cP MC_\dR^{k+l}(Y;\R)$, suppose $\sum_{i\in I}a_i\,[V_i,n,s_i,t_i,\om_i]$ satisfies Definition \ref{kh5def7}(ii)$(*)$ in $\cP MC_\dR^k(Y;\R)$ for some open neighbourhood $X$ of $\{0\}\t Y$ in $\R^n\t Y$, so that $\sum_{i\in I}a_i\,[V_i,n,s_i,t_i,\om_i]\ab =0$ in $\cP MC_\dR^k(Y;\R)$. Apply $-\cup[V',n',s',t',\om']$ for some $[V',n',s',t',\om']$ in $\cP MC_\dR^l(Y;\R)$. We will prove that in $\cP MC_\dR^{k+l}(Y;\R)$ we have
\e
\sum_{i\in I}(-1)^{nl+\deg\om_i(l+n')+\deg\om_i\deg\om'}a_i\,[\ti V_i,\ti n,\ti s_i,\ti t_i,\ti\om_i]=0,
\label{kh8eq36}
\e
where $[V_i,n,s_i,t_i,\om_i]\cup[V',n',s',t',\om']=(-1)^{nl+\deg\om_i(l+n')+\deg\om_i\deg\om'}[\ti V_i,\ab\ti n,\ab\ti s_i,\ab\ti t_i,\ab\ti\om_i]$ for $i\in I$ is as in equation~\eq{kh5eq39}.

As $[V',n',s',t',\om']$ is a generator of $\cP MC_\dR^l(Y;\R)$, $(s', t')\vert_{\supp\om'}:\supp\om'\ra \R^{n'}\t Y$ is proper over an open neighbourhood $X'$ of $\{0\}\t Y$ in $\R^{n'}\t Y$, with $t':V'\ra Y$ a submersion, and $\dim V'=\deg\om'+m+n'-l$. Define
\begin{equation*}
\ti X=\bigl\{\bigl((x,x'),y\bigr):x\in\R^n,\; x'\in\R^{n'},\; y\in Y,\; (x,y)\in X,\; (x',y)\in X'\bigr\}.
\end{equation*}
Then $\ti X$ is an open neighbourhood of $\{0\}\t Y$ in $\R^{n+n'}\t Y=\R^{\ti n}\t Y$, as $X,X'$ are open neighbourhoods of $\{0\}\t Y$ in $\R^n\t Y,\R^{n'}\t Y$. 

We will verify condition Definition \ref{kh5def7}(ii)$(*)$ for the expression \eq{kh8eq36} in $\cP MC_\dR^{k+l}(Y;\R)$. Suppose $\hat V$ is an oriented manifold of dimension $j$, $\hat t:\hat V\ra Y$ is smooth, and $\eta\in C^\iy\bigl(\La^{j+\ti n-k-l}T^*(\R^{\ti n}\t\hat V)\bigr)$ with $\supp\eta$ a compact subset of $(\id_{\R^{\ti n}}\t \hat t)^{-1}(\ti X)\subseteq\R^{\ti n}\t\hat V$. Then
\ea
&\sum_{i\in I}(-1)^{j\deg\ti\om_i+\deg\ti\om_i(\deg\ti\om_i-1)/2} \bigl[ (-1)^{nl+\deg\om_i(l+n')+\deg\om_i\deg\om'}\,a_i\bigr]\cdot{}
\nonumber\\
&\quad\int_{\ti V_i\t_{\ti t_i,Y,\hat t}\hat V}(\ti s_i\ci\pi_{\ti V_i},\pi_{\hat V})^*(\eta)\w\pi_{\ti V_i}^*(\ti\om_i)
\nonumber\\
&=\sum_{i\in I}\begin{aligned}[t]&
(-1)^{(j+l+n'+\deg\om')\deg\om_i+\deg\om_i(\deg\om_i-1)/2}\,a_i\cdot{}\\
&(-1)^{nl+j\deg\om'+\deg\om'(\deg\om'-1)/2+\deg\om_i\deg\om'}
\cdot{}\end{aligned}
\nonumber\\
&\int_{(V_i\t_{t_i,Y,t'}V')\t_{\pi_Y,Y,\hat t}\hat V}((s_i\ci\pi_{V_i},s'\ci\pi_{V'}),\pi_{\hat V})^*(\eta)\w(\pi_{V_i}^*(\om_i)\w\pi_{V'}^*(\om'))
\nonumber\\
&=\sum_{i\in I}\begin{aligned}[t]&(-1)^{(j+l+n'+\deg\om')\deg\om_i+\deg\om_i(\deg\om_i-1)/2}\,a_i\cdot{}\\
&(-1)^{nl+j\deg\om'+\deg\om'(\deg\om'-1)/2}\cdot{}\end{aligned}
\label{kh8eq37}\\
&\int_{V_i\t_{t_i,Y,\pi_Y}(V'\t_{t',Y,\hat t}\hat V)\!\!\!\!\!\!\!\!\!\!\!\!\!\!\!\!\!\!\!\!\!\!\!\!\!\!\!\!\!\!}(s_i\ci\pi_{V_i},\pi_{V'\t_Y\hat V})^*\bigl[(\id_{\R^n}\!\t\!(s'\ci\pi_{V'},\pi_{\hat V}))^*(\eta)\w\pi_{V'}^*(\om')\bigr]\w\pi_{V_i}^*(\om_i)
\nonumber\\
&=\sum_{i\in I}(-1)^{\check\jmath\deg\om_i+\deg\om_i(\deg\om_i-1)/2}a_i\int_{V_i\t_{t_i,Y,\check t}\check V\!\!}(s_i\ci\pi_{V_i},\pi_{\check V})^*(\check\eta)\w\pi_{V_i}^*(\om_i)=0.
\nonumber
\ea
 
Here in the first step we use \eq{kh5eq39} to substitute for $\ti V_i,\ti s_i,\ti t_i,\ti\om_i$. In the second we use the natural diffeomorphism $(V_i\t_YV')\t_Y\hat V\cong V_i\t_Y(V'\t_Y\hat V)$, by associativity of fibre products, and we rewrite the forms using properties of $\w$, acquiring an extra sign $(-1)^{\deg\om'\deg\om_i}$ from commuting $\pi_{V'}^*(\om')$ and $\pi_{V_i}^*(\om_i)$. Here $\bigl[(\id_{\R^n}\t(s'\ci\pi_{V'},\pi_{\hat V}))^*(\eta)\w\pi_{V'}^*(\om')\bigr]$ is a form on~$\R^n\t(V'\t_{t',Y,\hat t}\hat V)$.

In the third step of \eq{kh8eq37} we write $\check V=V'\t_{t',Y,\hat t}\hat V$, $\check t=\pi_Y:\check V\ra Y$, $\check\jmath=\dim\check V=j+n'-l+\deg\om'$, and set
\begin{equation*}
\check\eta=(-1)^{nl+j\deg\om'+\deg\om'(\deg\om'-1)/2}(\id_{\R^n}\t(s'\ci\pi_{V'},\pi_{\hat V}))^*(\eta)\w\pi_{V'}^*(\om')
\end{equation*}
in $C^\iy\bigl(\La^{\check\jmath+n-k}T^*(\R^{\ti n}\t\check V)\bigr)$. The fourth follows from Definition \ref{kh5def7}(ii)$(*)$ for $\sum_{i\in I}a_i\,[V_i,n,s_i,t_i,\om_i]$ in $\cP MC_\dR^k(Y;\R)$, with $\check V,\check\jmath,\check t$ in place of $V',l,t'$.

Equation \eq{kh8eq37} is \eq{kh5eq22} with \eq{kh8eq36} in place of $\sum_{i\in I}a_i\,[V_i,n,s_i,t_i,\om_i]$, and $\ti X,j,\hat V,\hat t$ in place of $X,l,V',t'$. Therefore \eq{kh8eq36} holds in $\cP MC_\dR^{k+l}(Y;\R)$, so $-\cup[V',n',s',t',\om']$ maps relation (ii) in $\cP MC_\dR^k(Y;\R)$ to relation (ii) in $\cP MC_\dR^{k+l}(Y;\R)$. Similarly, $[V,\ab n,\ab s,\ab t,\ab\om]\cup -$ maps relation (ii) in $\cP MC_\dR^l(Y;\R)$ to relation (ii) in $\cP MC_\dR^{k+l}(Y;\R)$. Hence $\cup$ in \eq{kh5eq38}--\eq{kh5eq39} is well defined. This completes the proof of Proposition~\ref{kh5prop5}.

\medskip

\noindent{\small\sc The Mathematical Institute, Radcliffe
Observatory Quarter, Woodstock Road, Oxford, OX2 6GG, U.K.

\noindent E-mail: {\tt joyce@maths.ox.ac.uk.}}

\end{document}